\documentclass[12pt,twoside,leqno,openany]{amsart}

\usepackage{amssymb,amsbsy,amsmath,amsfonts,amssymb,amscd,times,
graphics,color,xypic,footmisc,fancyhdr,multicol,fancybox,
graphicx,mathrsfs,rotating,ifthen,wasysym}
\usepackage[all]{xy}

\usepackage[T1]{fontenc}
\newcommand{\C}{\mathbb{C}}
\newcommand{\N}{\mathbb{N}}\newcommand{\R}{\mathbb{R}}
\newcommand{\U}{\mathbb{U}}\newcommand{\V}{\mathbb{V}}
\newcommand{\zero}[1]{\underline{#1}_{\circ}}
\let\mathcal\mathscr
\newtheorem{The}{Theorem}[section]
\newtheorem{Theorem}{Theorem}[section]
\newtheorem{Proposition}[The]{Proposition}
\newtheorem{Lemma}[The]{Lemma}
\newtheorem{Corollary}[The]{Corollary}
\theoremstyle{definition}
\newtheorem{Definition}[The]{Definition}
\newtheorem{Remark}[The]{Remark}

\newtheorem{Observation}[The]{Assertion}
\newtheorem{Scholium}[The]{Scholium}
\newtheorem{Openproblem}[The]{Open Problem}

\newcommand{\mathmotsf}[1]{\text{\footnotesize\sf #1}}
\newcommand{\smallmathmotsf}[1]{\text{\scriptsize\sf #1}}

\newcommand{\vf}{\vfill\end{document}}

\newcommand{\HEAD}[2]{%
\pagestyle{fancy}
\fancyhead[RO]{\scriptsize\sf\thepage}
\fancyhead[CO]{{\scriptsize\sf #1}}
\fancyhead[LE]{\scriptsize\sf\thepage}
\fancyhead[CE]{{\scriptsize\sf #2}}
\fancyfoot{}}


\begin{document}

\title[]{
Cartan equivalences for 
\\
5-dimensional CR-manifolds 
in $\mathbb C^4$
\\
belonging to General Class ${\sf III_1}$
}

\address{Department of Pure Mathematics,
University of Shahrekord, 88186-34141 Shahrekord, IRAN and School of
Mathematics, Institute for Research in Fundamental Sciences (IPM), P.
O. Box: 19395-5746, Tehran, IRAN}
\email{sabzevari@math.iut.ac.ir}

\address{D\'{e}partment de Math\'{e}matiques d'Orsay, B\^{a}timent 425,
Facult\'{e} des Sciences, Universit\'{e} Paris XI - Orsay, F-91405
Orsay Cedex, FRANCE}
\email{joel.merker@math.u-psud.fr}

\date{\number\year-\number\month-\number\day}

\maketitle

\begin{center}
Masoud {\sc Sabzevari}
and Jo\"el {\sc Merker}
\end{center}

\begin{abstract}
We reduce to various absolute parallelisms, namely to certain
$\{e\}$-structures on manifolds of dimensions 7, 6, 5, the biholomorphic
equivalence problem or the intrinsic CR equivalence problem for
generic submanifolds $M^5$ in $\mathbb C^4$ of CR dimension 1 and of codimension
3 that are maximally minimal and are geometry-preserving deformations
of one natural cubic model of Beloshapka, somewhere else called the
General Class $\sf III_1$ of 5-dimensional CR manifolds. Some inspiration
links exist with the treatment of the General Class $\sf II$ previously done
in 2007 by Beloshapka, Ezhov, Schmalz, and also with the
classification of nilpotent Lie algebras due to Goze, Khakimdjanov,
Remm.
\end{abstract}

\bigskip


\section{Introduction}
\label{introduction}
\HEAD{\ref{introduction}.~Introduction}{
Masoud {\sc Sabzevari} (Shahrekord) and Jo\"el {\sc Merker} (LM-Orsay)}

\medskip

For {\em systematic completeness} of our study of CR manifolds having:
\[
\text{\rm dimension}\,\,
\leqslant\,\,
\text{\bf 5},
\]
in the case of a
$3$-dimensional Levi nondegenerate hypersurface $M^3 \subset \C^2$, we
intentionally re-constructed in~\cite{ Merker-Sabzevari-M3-C2} an
explicit $\{ e\}$-structure on a certain natural
$8$-dimensional manifold $N^8 \longrightarrow M^3$, because our
current (wide) goal is similarly to perform completely explicit
constructions of $\{ e\}$-structures for CR equivalences in the six
General Classes:
\[
\text{\sf I},\ \ \ \ \
\text{\sf II},\ \ \ \ \
\text{\sf III}_{\text{\sf 1}},\ \ \ \ \
\text{\sf III}_{\text{\sf 2}},\ \ \ \ \
\text{\sf IV}_{\text{\sf 1}},\ \ \ \ \
\text{\sf IV}_{\text{\sf 2}},
\]
of all the possibly existing embedded
CR manifolds up to dimension $5$ that were
presented in~\cite{Merker-Pocchiola-Sabzevari-5-CR-II}.

The present memoir being specifically devoted to the:
\[
\text{\rm General Class}\,\,
\text{\sf III}_{\text{\sf 1}},
\]
let us review at first the General Class $\text{\sf I}$ in order to
appropriately recast the reader's thought in the right aspiration to
mathematical unification.

\subsection{Review of basic CR equivalences for Levi nondegenerate
hypersurfaces $M^3 \subset \C^2$.}
Let therefore $M^3 \subset \mathbb C^2$ be a (local) Levi-nondegenerate real
hypersurface of class at least $\mathcal{ C}^6$, graphed in
coordinates:
\[
(z,w)
=
\big(x+iy,\,u+iv\big)
\]
as:
\[
v
=
\varphi(x,y,u).
\]
An elementary normalization (\cite{ Merker-5-CR-III})
insures without loss of generality that:
\[
\varphi
=
x^2+y^2
+
{\rm O}(3).
\]

Starting with the (local) intrinsic generator for $T^{1, 0}M$:
\[
\mathcal{L}
:=
\frac{\partial}{\partial z}
-
\frac{\varphi_z}{i+\varphi_u}\,
\frac{\partial}{\partial u},
\]
having conjugate:
\[
\overline{\mathcal{L}}
=
\frac{\partial}{\partial\overline{z}}
-
\frac{\varphi_{\overline{z}}}{
-\,i+\varphi_u}\,
\frac{\partial}{\partial u},
\]
which in turn generates $T^{0, 1}M = \overline{ T^{ 1, 0} M}$, and
introducing:
\[
\mathcal{T}
:=
i\,\big[
\mathcal{L},\overline{\mathcal{L}}
\big],
\]
one gets by the Levi nondegeneracy assumption the natural frame:
\[
\big\{
\mathcal{T},
\overline{\mathcal{L}},
\mathcal{L}
\big\},
\]
for the complexified tangent bundle $\C \otimes_\R TM$.

\smallskip

We also introduce the dual coframe for $\C \otimes_\R T^*M$ consisting
of three $1$-forms denoted:
\[
\big\{
\rho_0,\,
\overline{\zeta}_0,\,
\zeta_0
\big\},
\]
namely satisfying by definition:
\[
\begin{array}{ccc}
\rho_0(\mathcal{T})=1
\ \ \ & \ \ \
\rho_0(\overline{\mathcal{L}})=0
\ \ \ & \ \ \
\rho_0(\mathcal{L})=0,
\\
\overline{\zeta}_0(\mathcal{T})=0
\ \ \ & \ \ \
\overline{\zeta_0}(\overline{\mathcal{L}})=1
\ \ \ & \ \ \
\overline{\zeta}_0\big(\mathcal{L}\big)=0,
\\
\zeta_0(\mathcal{T})=0
\ \ \ & \ \ \
\zeta_0(\overline{\mathcal{L}})=0
\ \ \ & \ \ \
\zeta_0\big(\mathcal{L}\big)=1.
\end{array}
\]

Abbreviating next:
\[
A
:=
-\,
\frac{\varphi_z}{i+\varphi_u},
\]
so that:
\[
\mathcal{T}
=
i\,
\bigg[
\frac{\partial}{\partial z}
+
A\,\frac{\partial}{\partial u},\,\,
\frac{\partial}{\partial\overline{z}}
+
\overline{A}\,
\frac{\partial}{\partial u}
\bigg],
\]
the {\sl Levi-factor} function:
\[
\ell
:=
i\,\big( \overline{A}_z + A\,\overline{A}_u -
A_{\overline{z}} - \overline{A}\,A_u \big),
\]
occurring in:
\[
\mathcal{T}
=
\ell\,\frac{\partial}{\partial u}
\]
is therefore {\em nowhere vanishing} by Levi nondegeneracy.

One also computes next:
\[
\aligned
\big[\mathcal{L},\,\mathcal{T}\big]
&
=
\bigg[
\frac{\partial}{\partial z}
+
A\,\frac{\partial}{\partial u},\,\,
\ell\,
\frac{\partial}{\partial u}
\bigg]
\\
&
=
\Big(
\ell_z
+
A\,\ell_u
-
\ell\,A_u
\Big)\,
\frac{\partial}{\partial u}
\\
&
=
\frac{\ell_z+A\,\ell_u-\ell\,A_u}{\ell}\,
\mathcal{T}.
\endaligned
\]
Then the appearing function:
\[
P
:=
\frac{\ell_z+A\,\ell_u-\ell\,A_u}{\ell},
\]
happens to be the single
one which enters the so-called {\sl initial Darboux structure:}
\[
\aligned
d\rho_0
&
=
P\,\rho_0\wedge\zeta_0
+
\overline{P}\,\rho_0\wedge\overline{\zeta}_0
+
i\,\zeta_0\wedge\overline{\zeta}_0,
\\
d\overline{\zeta}_0
&
=
0,
\\
d\zeta_0
&
=
0.
\endaligned
\]

\smallskip

As explained in~\cite{ Merker-5-CR-IV,
Merker-Sabzevari-M3-C2} and as is quite also very well known,
the initial ambiguity matrix group
for (local) biholomorphic
or CR equivalences of such hypersurfaces is:
\begin{eqnarray*}\footnotesize
\left\{
\left(\!
\begin{array}{ccc}
{\sf a}\overline{\sf a} & 0 & 0 \\
\overline{\sf b} & \overline{\sf a} & 0 \\
{\sf b} & 0 & {\sf a} \\
\end{array}
\!\right)
\,\in\,{\sf GL}_3(\C)
\colon
\ \ \
{\sf a}\in\C,\,\,
{\sf b}\in\C
\right\},
\end{eqnarray*}
just because through any
extrinsic local biholomorphism (or through any intrinsic
local CR-equivalence):
\[
h
\colon\ \ \
M
\,\longrightarrow\,
M',
\]
one has (\cite{ Merker-5-CR-IV}) for a certain coefficient-function $a$:
\[
\aligned
h_*'
\big(\mathcal{L}'\big)
&
=
a\,\mathcal{L},
\\
h_*'
\big(\overline{\mathcal{L}}'\big)
&
=
\overline{a}\,\overline{\mathcal{L}},
\endaligned
\]
whence:
\[
\aligned
h_*'\big(\mathcal{T}'\big)
&
=
h_*'\big(i\,\big[\mathcal{L}',\,\overline{\mathcal{L}}'\big]\big)
\\
&
=
i\,
\big[
h_*'\big(\mathcal{L}'\big),\,
h_*'\big(\overline{\mathcal{L}}'\big)
\big]
\\
&
=
i\,
\big[a\,\mathcal{L},\,\,
\overline{a}\,\overline{\mathcal{L}}\big]
\\
&
=
a\overline{a}\,
\underbrace{i\,\big[\mathcal{L},\,\overline{\mathcal{L}}\big]}_{
=
\,\mathcal{T}}
+
\underbrace{
i\,a\,\mathcal{L}\big(\overline{a}\big)}_{
=:\,\overline{b}}
\cdot
\overline{\mathcal{L}}
\underbrace{-
i\,\overline{a}\,\overline{\mathcal{L}}(a)}_{
=:\,b}
\cdot
\mathcal{L},
\endaligned
\]
so that setting:
\[
b
:=
-\,i\,\overline{a}\,\overline{\mathcal{L}}(a),
\]
forgetting how this further coefficient-function
is related to $a$, one obtains:
\[
h_*'\big(\mathcal{T}'\big)
=
a\overline{a}\,\mathcal{T}
+
\overline{b}\,\overline{\mathcal{L}}
+
b\,\mathcal{L}.
\]

\smallskip

Cartan's gist is to deal with the so-called {\sl lifted coframe:}
\[
\left(\!\!
\begin{array}{c}
\rho
\\
\overline{\zeta}
\\
\zeta
\end{array}
\!\!\right)
:=
\left(\!\!
\begin{array}{ccc}
{\sf a}\overline{\sf a} & 0 & 0
\\
\overline{\sf b} & \overline{\sf a} & 0
\\
{\sf b} & 0 & {\sf a}
\end{array}
\!\!\right)
\left(\!\!
\begin{array}{c}
\rho_0
\\
\overline{\zeta}_0
\\
\zeta_0
\end{array}
\!\!\right),
\]
in the space of $\big(x, y, u, {\sf a}, \overline{\sf a},
{\sf b}, \overline{\sf b}\big)$.

In~\cite{ Merker-Sabzevari-M3-C2}, after two absorbtions-normalizations
and after one prolongation, the desired equivalence
problem transforms to that of some\,\,---\,\,explicitly
computed\,\,---\,\,eight-dimensional coframe:
\[
\big\{
\rho,\zeta,\overline\zeta,\alpha,
\beta,\overline\alpha,\overline\beta,\delta
\big\}
\]
on a certain manifold $N^8 \longrightarrow M^3$
having $\{ e\}$-structure equations:
\begin{eqnarray*}
\begin{array}{ll}\footnotesize
d\rho=\alpha\wedge\rho+\overline\alpha\wedge\rho+i\,\zeta\wedge\overline\zeta,
&
\\
d\zeta=\beta\wedge\rho+\alpha\wedge\zeta, &
\\
d\overline\zeta=\overline\beta\wedge\rho+\overline\alpha\wedge\overline\zeta,
&
\\
d\alpha=\delta\wedge\rho+2\,i\,\zeta\wedge\overline{\beta}+i\,\overline{\zeta}\wedge\beta,

&
\\
d\beta=\delta\wedge\zeta+\beta\wedge\overline{\alpha}+{\frak
I}\,\overline\zeta\wedge\rho,
&
\\
d\overline\alpha=\delta\wedge\rho-2\,i\,\overline\zeta\wedge{\beta}-i\,{\zeta}\wedge\overline\beta,
&
\\
d\overline\beta=\delta\wedge\overline\zeta+\overline\beta\wedge{\alpha}+\overline{\frak
I}\,\zeta\wedge\rho,
&
\\
d\delta=\delta\wedge\alpha+\delta\wedge\overline\alpha+i\,\beta\wedge\overline\beta+{\frak
T}\,\rho\wedge\zeta+\overline{\frak T}\,\rho\wedge\overline\zeta, &
\end{array}
\end{eqnarray*}
with the single primary complex invariant:
\[
\aligned
\mathfrak{I}
&
:=
\frac{1}{6}\,
\frac{1}{{\sf a}\overline{\sf a}^3}\,
\Big(
-\,2\,
\overline{\mathcal{L}}\big(\mathcal{L}\big(
\overline{\mathcal{L}}\big(\overline{P}\big)\big)\big)
+
3\,\overline{\mathcal{L}}\big(
\overline{\mathcal{L}}\big(
\mathcal{L}\big(\overline{P}\big)\big)\big)
-\,
7\,\overline{P}\,\,
\overline{\mathcal{L}}\big(\mathcal{L}\big(\overline{P}\big)\big)
+
\\
&
\ \ \ \ \ \ \ \ \ \ \ \ \ \ \ \ \ \ \ \,
+
4\,\overline{P}\,
\mathcal{L}\big(\overline{\mathcal{L}}\big(
\overline{P}\big)\big)
-
\mathcal{L}\big(\overline{P}\big)\,
\overline{\mathcal{L}}\big(\overline{P}\big)
+
2\,\overline{P}\,\overline{P}\,\mathcal{L}
\big(\overline{P}\big)
\Big),
\endaligned
\]
and with one secondary invariant:
\[
\mathfrak{T}
=
\frac{1}{\overline{\sf a}}\,\bigg(\overline{\mathcal
L}(\overline{\frak I})-\overline P\,\overline{\frak
I}\bigg)-i\,\frac{\sf b}{{\sf a}\overline{\sf a}}\,\overline{\frak I}.
\]

\subsection{Explicitness obstacle.}
At the level of the function $P$, the completely explicit
formulas for $\mathfrak{ I}$, for $\mathfrak{ T}$
and for the $1$-forms constituting the $\{ e\}$-structure
remain writable on an article, {\em but not so anymore when one expresses
everything back in terms of the graphing function $\varphi$}.

Indeed, the real and imaginary parts $\Delta_1$ and $\Delta_2$ in:
\[
\mathfrak{I}
=
\frac{4}{{\sf a}\overline{\sf
a}^3}\big(\mathbf{\Delta}_1+i\,\mathbf{\Delta}_4\big)
\]
have numerators containing respectively
(\cite{ Merker-Sabzevari-CEJM}):
\[
{\bf 1\,553\,198}
\ \ \ \ \ \ \ \ \ \ \ \ \ \
\text{\rm and}
\ \ \ \ \ \ \ \ \ \ \ \ \ \
{\bf 1\,634\,457}
\]
monomials in the differential
ring in $\binom{ 6 + 3}{ 3} - 1 = 83$ variables:
\[
\mathbb{Z}\big[
\varphi_x,\varphi_y,
\varphi_{x^2},\varphi_{y^2},\varphi_{u^2},
\varphi_{xy},\varphi_{xu},\varphi_{yu},\dots\dots,
\varphi_{x^6},\varphi_{y^6},\varphi_{u^6},\dots
\big].
\]
Hence contrary to
the general case where $\varphi = \varphi ( x, y, u)$ does
depend upon the `{\sl CR-transversal}' variable $u$, in the so-called
{\sl rigid case} (often useful as a case of study-exploration) where
$\varphi = \varphi ( x, y)$ is independent of $u$ so that:
\[
P
=
\frac{\varphi_{zz\overline{z}}}{\varphi_{z\overline{z}}},
\]
one realizes that $\mathfrak{ I}$ is rather easily writable:
\[
\aligned
\mathfrak{I}
\bigg\vert_{
\substack{
\text{\sf rigid}\\
\text{\sf case}}}
&
=
\frac{1}{6}\,
\frac{1}{{\sf a}\overline{\sf a}^3}
\bigg(
\frac{\varphi_{z^2\overline{z}^4}}{\varphi_{z\overline{z}}}
-
6\,
\frac{\varphi_{z^2\overline{z}^3}\,\varphi_{z\overline{z}^2}}{
(\varphi_{z\overline{z}})^2}
-
\frac{\varphi_{z\overline{z}^4}\,\varphi_{z^2\overline{z}}}{
(\varphi_{z\overline{z}})^2}
-
4\,
\frac{\varphi_{z\overline{z}^3}\,\varphi_{z^2\overline{z}^2}}{
(\varphi_{z\overline{z}})^2}
+
\\
&
\ \ \ \ \ \ \ \ \ \ \ \ \ \ \ \ \ \ \
+
10\,
\frac{\varphi_{z\overline{z}^3}\,\varphi_{z^2\overline{z}}\,
\varphi_{z\overline{z}^2}}{
(\varphi_{z\overline{z}})^3}
+
15\,
\frac{(\varphi_{z\overline{z}^2})^2\,
\varphi_{z^2\overline{z}^2}}{
(\varphi_{z\overline{z}})^3}
-
15\,
\frac{(\varphi_{z\overline{z}^2})^3\,\varphi_{z^2\overline{z}}}{
(\varphi_{z\overline{z}})^4}
\bigg),
\endaligned
\]
and this therefore shows that
{\em there is a tremendous explosion of computational complexity
when one passes from the rigid case to the general case}.

Consequently, one expects an even much deeper computational
complexity when one addresses the question of passing to
CR manifolds of higher dimensions.

\subsection{Embedded CR submanifolds of CR dimension $1$ and nilpotent
Lie algebras.}
Consider now a general sufficiently smooth generic submanifold:
\[
M^{2+d}
\,\subset\,
\C^{1+d}
\]
having CR dimension $1$ and real codimension $d \geqslant 1$.
According to the background
article~\cite{Merker-Pocchiola-Sabzevari-5-CR-II},
the core bundle is:
\[
T^{1,0}M
:=
T^{1,0}\C^{1+n}
\cap
\big(
\C\otimes_\R TM
\big).
\]
To simplify (in fact, just a bit) the mathematical discussions,
we shall assume throughout that $M$ is {\em real analytic}.

The question is: to understand the possible initial geometries
of such CR manifolds, at least at a Zariski-generic point.
This question is not yet completely solved, because it opens
up infinitely many branches of classification.

Yet, one can present well known general considerations which
were already transparently explained in Sophus Lie's original
writings ({\em cf.} {\em e.g.}~\cite{ Engel-Lie-1888, Engel-Lie-1890,
Engel-Lie-1893}), though not targetly in a CR context.

Introduce the subdistributions of $\C \otimes_\R TM$:
\[
\aligned
D_\C^1M
&
:=
T^{1,0}M+T^{0,1}M,
\\
D_\C^2M
&
:=
{\sf Span}_{\mathcal{C}^\omega(M)}
\Big(
D_\C^1M
+
\big[T^{1,0}M,\,D_\C^1M\big]
+
\big[T^{0,1}M,\,D_\C^1M\big]
\Big),
\\
\\
D_\C^3M
&
:=
{\sf Span}_{\mathcal{C}^\omega(M)}
\Big(
D_\C^2M
+
\big[T^{1,0}M,\,D_\C^2M\big]
+
\big[T^{0,1}M,\,D_\C^2M\big]
\Big),
\endaligned
\]
and generally:
\[
D_\C^{k+1}M
:=
{\sf Span}_{\mathcal{C}^\omega(M)}
\Big(
D_\C^kM
+
\big[T^{1,0}M,\,D_\C^kM\big]
+
\big[T^{0,1}M,\,D_\C^kM\big]
\Big).
\]
By passing to some appropriate Zariski-open subset of $M$,
one may assume as is known that all these $D^k M$ become
true complex {\em vector subbundles} of $\C \otimes_\R TM$
having increasing ranks:
\[
2
=
r_1(M)
<
r_2(M)
<\cdots<
r_{k_M}(M)
=
r_{k_M+1}(M)
=
r_{k_M+2}(M)
=\cdots,
\]
until a first and final stabilization occurs. So we will admit
that the $M$ we consider enjoy constancies of such ranks,
and of several other invariants (finite in number) which
might happen to pop up later on.

As is known too, in the very special circumstance where:
\[
2
=
r_1(M)
=
r_2(M)
=
r_3(M)
=\cdots,
\]
namely when:
\[
\big[T^{1,0}M,\,T^{0,1}M\big]
\,\subset\,
T^{1,0}M
+
T^{0,1}M,
\]
the real analytic CR-generic submanifold $M^{ 2 + d} \subset \C^{ 1+d}$
is, locally in some small neighborhood of each of its points,
biholomorphic to $\C \times \R^d$, a degenerate case rapidly set away.
Sometimes, one also says that $M$ is {\sl Levi-flat},
or equivalently (just when the CR dimension equals $1$),
that $M$ is {\sl
holomorphically degenerate} (\cite{ Merker-2001}).

In fact, whenever the maximal possible rank:
\[
r_{k_M}(M)
=
r_{k_M+1}(M)
=
r_{k_M+2}(M)
=\cdots
\,<\,
2+d
=
\dim_\R\,M,
\]
is still smaller than the dimension of $M$ (not necessarily equal to
$2$), one realizes that $M$ is similarly, locally in a neighborhood of
a Zariski-generic point, biholomorphic to a CR-generic
submanifold of $\C^{ 1 + d}$ which is contained
in $2 + d - r_{ k_M} (M)$ straight real hyperplanes in
transverse intersection, a case which is also degenerate hence is
also disregarded, just because it essentially comes down to the case
of CR-generic submanifolds having smaller dimension that $2 + d$.

So the question is: to understand all the possible geometries
of such CR manifolds $M^{2 + d} \subset \C^{ 1 + d}$
of CR dimension $1$ that have the so-called
constant nonholonomic property that:
\[
D_\C^{k_M}M
=
\C\otimes_\R TM.
\]
In the Several
Complex Variables literature, such CR manifolds happen to be those called {\em
minimal in the sense of Tumanov}, or equivalently
{\sl of finite type in
the sense of Bloom-Graham}, and they happen to necessarily be also
simultaneously holomorphically nondegenerate too (an implication which
is true only in CR dimension $1$, as is
easily checked).

Beloshapka and his students, {\em e.g.} Shananina, Mamai, Kossovskiy
and others, have put some emphasis on the study of such a class of CR
manifolds, notably in the search for {\it nice} models which would
potentially reveal new Cartan geometries.

The truth is that this research field, like the one of hyperbolic
groups in the sense of Gromov, is {\em per se} opened to infinitely
many untamable branches of complexity, for one soon realizes after a
moment of reflection that the Lie algebras of infinitesimal CR
automorphisms (assuming for simplicity that everything is real
analytic) are deeply related to the classification of {\em nilpotent
Lie algebras}, an area which is known to be very rich and very infinite,
as Lie himself understood more than one century ago
({\em see} Chapter~28 in the English translation~\cite{Merker-2010}
of Volume~I of the {\sl Theorie der Transformationsgruppen}).
Section~\ref{nilpotent-5} here is devoted to review
the easiest part (only up
to dimension $5$) of the deep nilpotent
Lie algebra classification theorems
of Goze, Khakimdjanov, Remm up to dimension $8$
(\cite{Goze-web, Goze-2008, Goze-Khakimdjanov-1996}),
which already shows up an exploding ramification of very many branches.

In dimension $4$, there is a single irreducible nilpotent Lie algebra:
\[
\mathfrak{n}_4^1\colon
\ \ \ \ \ \ \ \ \ \
\left\{
\aligned
{}
[{\sf x}_1,{\sf x}_2]
&
=
{\sf x}_3,
\\
[{\sf x}_1,{\sf x}_3]
&
=
{\sf x}_4.
\endaligned\right.
\]
Correspondingly, as Beloshapka discovered in~1997 (\cite{
Beloshapka-1998}), a real analytic $4$-dimensional local CR-generic
submanifold $M^4 \subset \C^3$ of codimension $2$ whose complex
tangent bundle satisfies:
\[
\aligned
\C\otimes_\R TM
&
=
T^{1,0}M+T^{0,1}M
+
[T^{1,0}M,\,T^{0,1}M]
+
\big[T^{1,0}M,[T^{1,0}M,\,T^{0,1}M]\big]
+
\\
&
\ \ \ \ \ \ \ \ \ \ \ \ \ \ \ \ \ \ \ \ \ \ \ \ \ \ \ \ \ \ \ \ \ \ \ \ \
\ \ \ \ \ \ \ \ \ \ \ \ \ \ \ \ \ \ \ \ \ \ \ \ \ \
+
\big[T^{0,1}M,[T^{1,0}M,\,T^{0,1}M]\big]
\endaligned
\]
may always be represented, in suitable holomorphic coordinates:
\[
\big(z,w_1,w_2\big)
=
\big(
z,u_1+iv_1,u_2+iv_2
\big)
\]
by two complex defining equations of the
specific form:
\[
\aligned
v_1
&
=
z\overline{z}
+
{\rm O}_4(x,y,u_1,u_2),
\\
v_2
&
=
z\overline{z}\big(z+\overline{z}\big)
+
{\rm O}_4(x,y,u_1,u_2).
\endaligned
\]
Since then, such CR manifolds have been intensively studied,
by Beloshapka-Ezhov-Schmalz
who constructed a canonical Cartan connection (\cite{ BES-2008})
and who generalized Pinchuk-Vitushkin's germ extension phenomenon
(\cite{ BES-2006}),
by Gammel-Kossovskiy (\cite{ Gammel-Kossovskiy-2006}),
and by Beloshapka-Kossovskiy who provided a final complete
classification (\cite{ Beloshapka-Kossovskiy-2011}).

In~\cite{Merker-Pocchiola-Sabzevari-5-CR-II}, one refers
to the:
\[
\aligned
&
\boxed{\text{\sf General Class $\text{\sf II}$:}}
\\
&
\boxed{\,\,
\aligned
M^4\subset\C^3
\ \
&
\text{\rm with}\ \
\Big\{\mathcal{L},\,\overline{\mathcal{L}},\,\,
\big[\mathcal{L},\overline{\mathcal{L}}\big],\,\,
\big[\mathcal{L},\,\big[\mathcal{L},\overline{\mathcal{L}}\big]\big]
\Big\}\ \
\\
&
\text{\rm constituting a frame for}\ \
\C\otimes_\R TM.\,\,
\endaligned}
\endaligned
\]

The next natural General Class (\cite{ Merker-Pocchiola-Sabzevari-5-CR-II})
is the:
\[
\aligned
&
\boxed{\text{\sf General Class $\text{\sf III}_{\text{\sf 1}}$:}}
\\
&
\boxed{\,\,
\aligned
M^5\subset\C^4
\ \
&
\text{\rm with}\ \
\Big\{\mathcal{L},\,\overline{\mathcal{L}},\,\,
\big[\mathcal{L},\overline{\mathcal{L}}\big],\,\,
\big[\mathcal{L},\,
\big[\mathcal{L},\overline{\mathcal{L}}\big]\big],\,\,
\big[\overline{\mathcal{L}},\,
\big[\mathcal{L},\overline{\mathcal{L}}\big]\big]
\Big\}\,\,
\\
&
\text{\rm constituting a frame for}\ \
\C\otimes_\R TM,
\endaligned}
\endaligned
\]
and it is intrinsically related to the
irreducible $5$-dimensional nilpotent Lie algebra
(labeled here in the notation of Goze-Remm):
\[
\mathfrak{n}_5^4\colon
\ \ \ \ \ \ \ \ \ \
\left\{
\aligned
{}
[{\sf x}_1,{\sf x}_2]
&
=
{\sf x}_3,
\\
[{\sf x}_1,{\sf x}_3]
&
=
{\sf x}_4,
\\
[{\sf x}_2,{\sf x}_3]
&
=
{\sf x}_5.
\endaligned\right.
\]

Three years ago, we started to study CR equivalences of such CR
manifolds belonging
to the General Class $\text{\sf III}_{\text{\sf 1}}$, trying in the first months
to directly construct a Cartan connection as did
Beloshapka-Ezhov-Schmalz for the General Class $\text{\sf II}$. But
inspired by Chern's seminal 1939 paper on equivalences of third order
ordinary differential equations under contact transformations, we
realized that it would be better to perform at first a pure
exploration of the problem by employing the powerful tools of Cartan's
method of equivalence, in order to avoid
as much as possible those errors of understanding that are caused
by a too quick belief that certain features would somewhat easily generalize.

Because we have not
been since then aware of any other paper or preprint
or author having attacked the same problem, we decided to wait
until the study reached a point of maturity in which
everything could be presented in full computational details,
whatever complexity the theory has.

Of course in such a General
Class $\text{\sf III}_{\text{\sf 1}}$, it is known that the cubic
model $M_{\sf c}^5 \subset \C^4$ in coordinates:
\[
\big(z,w_1,w_2,w_3\big)
=
\big(z,\,u_1+iv_1,\,u_2+iv_2,\,u_3+iv_3\big)
\]
was also discovered by Beloshapka:
\[
\aligned
v_1
&
=
z\overline{z},
\\
v_2
&
=
z\overline{z}\,
\big(z+\overline{z}\big),
\\
v_3
&
=
z\overline{z}
\big(
-\,i\,z
+
i\,\overline{z}\big).
\endaligned
\]
But the Cartan invariants of the geometry-preserving deformations of
such a model have apparently never been studied, and such a study is
the main goal of the present memoir. Granted that the general class ${\sf IV_1}$ is already well studied since Chern-Moser (\cite{Chern-Moser}), a forthcoming paper by Samuel
Pocchiola will soon treat the General Class
$\text{\sf III}_{\text{\sf 2}}$ (as presented
in~\cite{ Merker-Pocchiola-Sabzevari-5-CR-II}), thus closing up the
study of CR equivalences of CR manifolds up to dimension $5$
(an overall systematic review
is planned to appear at the end).

With $J$ being
the standard complex structure of
$T\C^4$, one sets as usual
(\cite{ Merker-Pocchiola-Sabzevari-5-CR-II}):
\[
T^cM
:=
TM\cap J(TM),
\]
or equivalently:
\[
T^cM
:=
{\rm Re}\,T^{1,0}M.
\]
Our first elementary result appears in Section~\ref{Heisengerg-Beloshapka},
{\em cf.} also~\cite{ Merker-5-CR-III}.

\begin{Proposition}
Every real analytic 5-dimensional local CR-generic submanifold $M^5
\subset \C^4$ of codimension 3 which is maximally minimal, namely
which satisfies:
\[
\aligned
D^1M
&
=
T^cM
\ \ \ \ \
\text{has rank $2$},
\\
D^2M
&
=
T^cM+[T^cM,T^cM]
\ \ \ \ \
\text{has rank $3$},
\\
D^3M
&
=
T^cM+[T^cM,T^cM]
+
\big[T^cM,[T^cM,T^cM]\big]
\ \ \ \ \
\text{has maximal possible rank $5$},
\endaligned
\]
may be represented, in suitable holomorphic coordinates $(z, w_1, w_2,
w_3)$, by three complex defining equations of the specific form:
\[
\left[
\aligned
w_1-\overline{w}_1
&
=
2i\,z\overline{z}
+
\Pi_1\big(z,\overline{z},\overline{w}_1,\overline{w}_2,\overline{w}_3\big),
\\
w_2-\overline{w}_2
&
=
2i\,z\overline{z}(z+\overline{z})
+
\Pi_2\big(z,\overline{z},\overline{w}_1,\overline{w}_2,\overline{w}_3\big),
\\
w_3-\overline{w}_3
&
=
2\,z\overline{z}(z-\overline{z})
+
\Pi_3\big(z,\overline{z},\overline{w}_1,\overline{w}_2,\overline{w}_3\big),
\endaligned\right.
\]
where the three remainders $\Pi_1$, $\Pi_2$, $\Pi_3$ are all an
${\rm O} ( \vert z \vert^4) + z\overline{ z}\, {\rm O}( \vert w
\vert)$.

Conversely, for any choice of three such analytic functions enjoying
these conditions, the zero-locus of the three equations above represents a real analytic 5-dimensional
local CR-generic submanifold $M^5 \subset \C^4$ of codimension $3$
which is maximally minimal.
\qed
\end{Proposition}

Next, a general $(1, 0)$ holomorphic vector field:
\[
{\sf X}
=
Z(z,w)\,
\frac{\partial}{\partial z}
+
W^1(z,w)\,
\frac{\partial}{\partial w_2}
+
W^2(z,w)\,
\frac{\partial}{\partial w_2}
+
W^3(z,w)\,
\frac{\partial}{\partial w_3}
\]
is an infinitesimal CR automorphism of Beloshapka's cubic model $M_{\sf c}^5$
if by definition ${\sf X} + \overline{\sf X}$ is tangent to
$M_{\sf c}^5$. By analyzing in great details the
system of linear partial differential equations satisfied
by the unknown functions $Z$, $W^1$, $W^2$, $W^3$,
we obtain the second already known:

\begin{Proposition}
The Lie algebra $\mathfrak{ aut}_{ CR}
(M) = 2\, {\rm Re}\, \mathfrak{hol}(M)$ of the infinitesimal CR
automorphisms of the $5$-dimensional $3$-codimensional CR-generic
model cubic
$M_{\sf c}^5 \subset
\mathbb C^4$ represented by the three real
graphed equations:
\[
\left[
\aligned w_1-\overline{w}_1 & = 2i\,z\overline{z},
\\
w_2-\overline{w}_2 & = 2i\,z\overline{z}(z+\overline{z}),
\\
w_3-\overline{w}_3 & = 2\,z\overline{z}(z-\overline{z}),
\endaligned\right.
\]
is $7$-dimensional and it is generated by the $\mathbb{ R}$-linearly
independent real parts of the
following seven $(1, 0)$ holomorphic vector fields:
\[
\aligned T & :=
\partial_{w_1},
\\
S_1 & :=
\partial_{w_2},
\\
S_2 & :=
\partial_{w_3},
\\
 L_1 & :=
\partial_z
+ (2iz)\,\partial_{w_1} + (2iz^2+4w_1)\,\partial_{w_2} +
2z^2\,\partial_{w_3},
\\
L_2 & := i\,\partial_z + (2z)\,\partial_{w_1} +
(2z^2)\,\partial_{w_2} - (2iz^2-4w_1)\,\partial_{w_3},
\\
D & := z\,\partial_z + 2w_1\,\partial_{w_1} + 3w_2\,\partial_{w_2} +
3w_3\,\partial_{w_3},
\\
R & := iz\,\partial_z - w_3\,\partial_{w_2} + w_2\,\partial_{w_3},
\endaligned
\]
having Lie bracket commutator table:

\medskip
\begin{center}
\begin{tabular} [t] { c | c c c c c c c }
& $S_2$ & $S_1$ & $T$ & $L_2$ & $L_1$ & $D$ & $R$
\\
\hline $S_2$ & $0$ & $0$ & $0$ & $0$ & $0$ & $3S_2$ & $-S_1$
\\
$S_1$ & $*$ & $0$ & $0$ & $0$ & $0$ & $3S_1$ & $S_2$
\\
$T$ & $*$ & $*$ & $0$ & $4S_2$ & $4S_1$ & $2T$ & $0$
\\
$L_2$ & $*$ & $*$ & $*$ & $0$ & $-4T$ & $L_2$ & $-L_1$
\\
$L_1$ & $*$ & $*$ & $*$ & $*$ & $0$ & $L_1$ & $L_2$
\\
$D$ & $*$ & $*$ & $*$ & $*$ & $*$ & $0$ & $0$
\\
$R$ & $*$ & $*$ & $*$ & $*$ & $*$ & $*$ & $0$.
\end{tabular}
\end{center}
\end{Proposition}

One easily realizes that in the natural grading:
\[
\aligned
\mathfrak{g}_{-3}
&
:=
\mathmotsf{Span}_\R
\big<S_1,\,S_2\big>,
\\
\mathfrak{g}_{-2}
&
:=
\mathmotsf{Span}_\R
\big<T\big>,
\\
\mathfrak{g}_{-1}
&
:=
\mathmotsf{Span}_\R
\big<L_1,L_2 \big>,
\\
\mathfrak{g}_0
&
:=
\mathmotsf{Span}_\R
\big<D,R\big>,
\endaligned
\]
the above nilpotent Lie algebra $\mathfrak{ n}_5^4$
is isomorphic to:
\[
\mathfrak{g}_{-3}
\oplus
\mathfrak{g}_{-2}
\oplus
\mathfrak{g}_{-1},
\]
as is well known in Tanaka's theory. Of course,
a wealth of other correspondences between nilpotent Lie algebras
and CR manifolds of CR dimension $1$ certainly exist,
but we skip entering this already
much studied question in order to enter the core of a new
systematic effective development of Cartan's equivalence
method, an aspect which is, as we believe, not yet enough
developed in the mathematical literature as far
as hardest computations are concerned.

\smallskip

As it also appears independently in~\cite{ Merker-5-CR-IV},
we verify:

\begin{Proposition}
The initial ambiguity matrix associated to the local biholomorphic
equivalence problem between the cubic $5$-dimensional model CR-generic
submanifold $M^5_{\sf c}$ and any other maximally minimal CR-generic
$5$-dimensional submanifolds ${M'}^5 \subset \C^4$ under local
biholomorphic transformations is of the general form:
\[
\left(\!\!
\begin{array}{ccccc}
{\sf a}\overline{\sf a}\overline{\sf a} & 0 & \overline{\sf c} &
\overline{\sf e} & \overline{\sf d}
\\
0 & {\sf a}{\sf a}\overline{\sf a} & {\sf c} & {\sf d} & {\sf e}
\\
0 & 0 & {\sf a}\overline{\sf a} & \overline{\sf b} & {\sf b}
\\
0 & 0 & 0 & \overline{\sf a} & 0
\\
0 & 0 & 0 & 0 & {\sf a}
\end{array}
\!\!\right),
\]
where ${\sf a}$, ${\sf b}$, ${\sf c}$, ${\sf e}$, ${\sf d}$ are
complex numbers. Moreover, the collection of all these matrices makes
up a real $10$-dimensional matrix Lie subgroup of ${\sf GL}_5 ( \C)$.
\end{Proposition}

In Section~\ref{equivalence-model},
as a preliminary to higher level computations,
we perform the Cartan equivalence algorithm on the cubic
model $M_{\sf c}^5$ and we obtain an $\{ e\}$-structure
on a $7$-dimensional manifold of the form:
\[
\aligned d\sigma & =
\big(2\,\alpha+\overline{\alpha}\big) \wedge \sigma+
{\rho\wedge\zeta},
\\
d\overline{\sigma} & = \big(\alpha+2\,\overline{\alpha}\big) \wedge
\overline{\sigma}+ {\rho\wedge\overline{\zeta}},
\\
d\rho & = (\alpha+\overline{\alpha})\wedge\rho +
i\,\zeta\wedge\overline{\zeta},
\\
 d\zeta & = \alpha\wedge\zeta,
 \\
 d\overline{\zeta} & = \overline{\alpha}\wedge\overline{\zeta},
 \\
 d\alpha&=0,
 \\
 d\overline{\alpha}&=0,
\endaligned
\]
all structure functions being {\em constant},
and up to a very mild change
of basis, these are nothing but the Maurer-Cartan equations
on the $7$-dimensional Lie group associated to
the above $2$-dimensional (semidirect product)
extension of the nilpotent Lie algebra $\mathfrak{ n}_5^4$.
This performing of Cartan's method on the cubic
model $M_{\sf c}^5$ has the virtue of setting up
a kind of `{\sc gps}' for orientation in the
much deeper `computational jungle' of the general case.

\smallskip

In fact, it is in Section~\ref{initial-complex-frame},
that we start out the
main computations for general geometry-preserving
deformations $M^5 \subset \C^4$ of the cubic model
$M_{\sf c}^5 \subset \C^4$ with a first, already
computationally nontrivial, result.

\begin{Proposition}
For any local real analytic CR-generic submanifold $M^5 \subset
\C^4$ which is represented near the
origin as a graph:
\[
\aligned v_1&:= \varphi_1(x,y,u_1,u_2,u_3),
\\
v_2&:= \varphi_2(x,y,u_1,u_2,u_3),
\\
v_3&:= \varphi_3(x,y,u_1,u_2,u_3),
\endaligned
\]
in coordinates:
\[
\big(z,w_1,w_2,w_3\big) = \big(x+iy,u_1+iv_1,u_2+iv_2,u_3+iv_3\big),
\]
its fundamental intrinsic complex bundle $T^{0, 1}M$ is generated by:
\[
\boxed{
\aligned
\overline{\mathcal{L}} & =
\frac{\partial}{\partial\overline{z}}
+
{\sf A}_1\,
\frac{\partial}{\partial\overline{w}_1}
+
{\sf A}_2\,
\frac{\partial}{\partial\overline{w}_2}
+
{\sf A}_3\,
\frac{\partial}{\partial\overline{w}_3},
\endaligned}
\]
where:
\[
\footnotesize \aligned {\sf A}_1&=
\frac{\Lambda^1_1}{\Delta}+i\frac{\Lambda^1_2}{\Delta},
\\
{\sf A}_2&= \frac{\Lambda^2_1}{\Delta}+i\frac{\Lambda^2_2}{\Delta},
\\
{\sf A}_3&= \frac{\Lambda^3_1}{\Delta}+i\frac{\Lambda^3_2}{\Delta},
\endaligned
\]
where:
\[
\aligned \Delta&=\sigma^2+\tau^2,
\endaligned
\]
with:
\[
\aligned \sigma&=
\varphi_{3u_3}+\varphi_{1u_1}+\varphi_{2u_2}-\varphi_{1u_2}\varphi_{3u_1}\varphi_{2u_3}-
\varphi_{1u_3}\varphi_{2u_1}\varphi_{3u_2}+\varphi_{1u_2}\varphi_{2u_1}\varphi_{3u_3}-
\\
&-
\varphi_{1u_1}\varphi_{2u_2}\varphi_{3u_3}+\varphi_{1u_1}\varphi_{2u_3}\varphi_{3u_2}+\varphi_{1u_3}\varphi_{3u_1}\varphi_{2u_2},
\\
 \tau & =
-1+\varphi_{1u_1}\varphi_{2u_2}-\varphi_{2u_3}\varphi_{3u_2}-\varphi_{1u_3}\varphi_{3u_1}+
\varphi_{2u_2}\varphi_{3u_3}-\varphi_{1u_2}\varphi_{2u_1}+\varphi_{1u_1}\varphi_{3u_3},
\endaligned
\]
and where:
\[
\footnotesize \aligned \Lambda^1_1&=
\big(-\varphi_{3u_3}\varphi_{2x}\varphi_{1u_2}-\varphi_{1u_3}\varphi_{3y}+\varphi_{2u_2}\varphi_{1x}\varphi_{3u_3}+
\varphi_{3u_3}\varphi_{1y}-\varphi_{1x}-\varphi_{2y}\varphi_{1u_2}+
\\
&+\varphi_{2u_3}\varphi_{3x}\varphi_{1u_2}+\varphi_{2u_2}\varphi_{1y}-\varphi_{2u_3}\varphi_{3u_2}\varphi_{1x}-\varphi_{2u_2}\varphi_{1u_3}\varphi_{3x}+
\varphi_{2x}\varphi_{1u_3}\varphi_{3u_2}\big)\sigma +
\\
&+
\big(\varphi_{1u_3}\varphi_{3x}-\varphi_{1y}+\varphi_{2x}\varphi_{1u_2}+\varphi_{2u_3}\varphi_{1u_2}\varphi_{3y}-\varphi_{2u_2}\varphi_{1x}-
\varphi_{2u_3}\varphi_{3u_2}\varphi_{1y}-
\\
&-\varphi_{3u_3}\varphi_{1x}-\varphi_{2u_2}\varphi_{1u_3}\varphi_{3y}-\varphi_{3u_3}\varphi_{1u_2}\varphi_{2y}+\varphi_{1u_3}\varphi_{3u_2}\varphi_{2y}+
\varphi_{2u_2}\varphi_{3u_3}\varphi_{1y} \big)\tau,
\endaligned
\]
with similar formulas for $\Lambda_2^1$, $\Lambda_1^2$,
$\Lambda_2^2$, $\Lambda_1^3$, $\Lambda_2^3$.
\end{Proposition}

Next, introduce the third vector field:
\[
\mathcal T
:
=i\,\big[\mathcal L,\overline{\mathcal L}\big],
\]
which is {\em real:}
\[
\overline{\mathcal{T}}
=
\mathcal{T}.
\]
Direct computations provide the three numerators in:
\[
\boxed{ \mathcal T =
\frac{\Upsilon_1}{\Delta^3}\frac{\partial}{\partial u_1}+
\frac{\Upsilon_2}{\Delta^3}\frac{\partial}{\partial u_2} +
\frac{\Upsilon_3}{\Delta^3}\frac{\partial}{\partial u_3},}
\]
namely:
\begin{equation*}
\footnotesize\aligned \Upsilon_1 & =
 -\big(\Delta^2\Lambda^1_{2x}-\Delta\Delta_x\Lambda^1_2-
\Delta^2\Lambda^1_{1y}+\Delta\Delta_y\Lambda^1_1
+\Delta\Lambda^1_1\Lambda^1_{2u_1} -
\Delta\Lambda^1_2\Lambda^1_{1u_1}- \Delta\Lambda^2_2\Lambda^1_{1u_2}+
\\
&+ \Delta_{u_2}\Lambda^1_1\Lambda^2_2-
\Delta\Lambda^3_2\Lambda^1_{1u_3}+
\Delta_{u_3}\Lambda^3_2\Lambda^1_1+
\Delta\Lambda^2_1\Lambda^1_{2u_2}-
\Delta_{u_2}\Lambda^2_1\Lambda^1_{2}+
\Delta\Lambda^3_1\Lambda^1_{2u_3}-
\Delta_{u_3}\Lambda^3_1\Lambda^1_{2}\big),
\\
\Upsilon_2 & =
 - \big(\Delta^2\Lambda^2_{2x}- \Delta\Delta_x\Lambda^2_2+
\Delta\Lambda^1_{1}\Lambda^2_{2u_1}-
\Delta_{u_1}\Lambda^1_1\Lambda^2_2 -
\Delta^2\Lambda^2_{1y}+\Delta\Delta_y\Lambda^2_1-
\Delta\Lambda^1_2\Lambda^2_{1u_1}+
\\
&+
\Delta_{u_1}\Lambda^1_2\Lambda^2_1+\Delta\Lambda^2_1\Lambda^2_{2u_2}-
\Delta\Lambda^2_2\Lambda^2_{1u_2}+\Delta\Lambda^3_1\Lambda^2_{2u_3}-
\Delta_{u_3}\Lambda^3_1\Lambda^2_{2}-
\Delta\Lambda^3_2\Lambda^2_{1u_3}+
\Delta_{u_3}\Lambda^3_2\Lambda^2_{1}\big),
\\
\Upsilon_3 & = -
\big(\Delta^2\Lambda^3_{2x}-\Delta\Delta_x\Lambda^3_2+
\Delta\Lambda^1_{1}\Lambda^3_{2u_1} -
\Delta_{u_1}\Lambda^1_1\Lambda^3_2
-\Delta^2\Lambda^3_{1y}+\Delta\Delta_y\Lambda^3_1-
\Delta\Lambda^1_2\Lambda^3_{1u_1}+
\\
&+
\Delta_{u_1}\Lambda^1_2\Lambda^3_1-\Delta\Lambda^2_2\Lambda^3_{1u_2}+
\Delta_{u_2}\Lambda^2_2\Lambda^3_{1} +
\Delta\Lambda^3_1\Lambda^3_{2u_3}-
\Delta\Lambda^3_2\Lambda^3_{1u_3}+\Delta\Lambda^2_1\Lambda^3_{2u_2}-
\Delta_{u_2}\Lambda^2_1\Lambda^3_{2}\big).
\endaligned
\end{equation*}

Next, introduce the two Lie brackets of length three:
\[
\aligned
\mathcal S
&
:=
\big[
\mathcal L,\mathcal T],
\\
\overline{\mathcal S}
&
:=
\big[\,
\overline{\mathcal L},\mathcal T
\big].
\endaligned
\]
Again, direct computations provide the expressions:

\[
\aligned \mathcal S &=&
\frac{\Gamma^1_1-i\,\Gamma^1_2}{\Delta^5}\frac{\partial}{\partial
u_1}+\frac{\Gamma^2_1-i\,\Gamma^2_2}{\Delta^5}\frac{\partial}{\partial
u_2}+\frac{\Gamma^3_1-i\,\Gamma^3_2}{\Delta^5}\frac{\partial}{\partial
u_3},
\endaligned
\]
where, by allowing the two notational coincidences $x\equiv x_1$ and
$y\equiv x_2$, the numerators are (for $i=1,2$):
\begin{equation*}
\footnotesize\aligned \Gamma^1_i & = -2 \big({\textstyle\frac{1}{4}}
\Delta^2\Upsilon_{1{x_i}}- 3\Delta\Delta_{x_i}\Upsilon_1+
\Delta\Lambda^1_i\Upsilon_{1u_1}-
2\Delta_{u_1}\Lambda^1_i\Upsilon_{1}-
\Delta\Lambda^1_{iu_1}\Upsilon_{1}-
\Delta\Lambda^1_{iu_2}\Upsilon_{2}+
\\
&+ \Delta_{u_2}\Lambda^1_i\Upsilon_{2}-
\Delta\Lambda^1_{iu_3}\Upsilon_{3}+
\Delta_{u_3}\Lambda^1_i\Upsilon_{3}+
\Delta\Lambda^2_i\Upsilon_{1u_2}-
3\Delta_{u_2}\Lambda^2_i\Upsilon_{1}+
\Delta\Lambda^3_i\Upsilon_{1u_3}-
3\Delta_{u_3}\Lambda^3_i\Upsilon_{1}\big),
\\
\Gamma^2_i & = -2 \big( \Delta^2\Upsilon_{2{x_i}}-
3\Delta\Delta_{x_i}\Upsilon_2+ \Delta\Lambda^1_i\Upsilon_{2u_1}-
3\Delta_{u_1}\Lambda^1_i\Upsilon_{2}-
\Delta\Lambda^2_{iu_1}\Upsilon_{1}+
\Delta_{u_1}\Lambda^2_{i}\Upsilon_{1}+
\\
&+ \Delta\Lambda^2_i\Upsilon_{2u_2}-
2\Delta_{u_2}\Lambda^2_{i}\Upsilon_{2}-
\Delta\Lambda^2_{iu_2}\Upsilon_{2}-
\Delta\Lambda^2_{iu_3}\Upsilon_{3}+
\Delta_{u_3}\Lambda^2_i\Upsilon_{3}+
\Delta\Lambda^3_i\Upsilon_{2u_3}-
3\Delta_{u_3}\Lambda^3_i\Upsilon_{2}\big),
\\
\Gamma^3_i & = -2\big( \Delta^2\Upsilon_{3{x_i}}-
3\Delta\Delta_{x_i}\Upsilon_3+ \Delta\Lambda^1_i\Upsilon_{3u_1}-
3\Delta_{u_1}\Lambda^1_i\Upsilon_{3}+
\Delta\Lambda^2_{i}\Upsilon_{3u_2}-
3\Delta_{u_2}\Lambda^2_{i}\Upsilon_{3}-
\\
&- \Delta\Lambda^3_{iu_1}\Upsilon_{1}+
\Delta_{u_1}\Lambda^3_{i}\Upsilon_{1}-
\Delta\Lambda^3_{iu_2}\Upsilon_{2}+
\Delta_{u_2}\Lambda^3_{i}\Upsilon_{2}+
\Delta\Lambda^3_i\Upsilon_{3u_3}-
2\Delta_{u_3}\Lambda^3_i\Upsilon_{3}-
\Delta\Lambda^3_{iu_3}\Upsilon_{3}\big).
\endaligned
\end{equation*}

At this point, remind that quite straightforwardly from the definition,
one has:

\begin{Proposition}
A real analytic CR-generic submanifold $M^5 \subset \C^4$
belongs to the General Class $\text{\sf III}_{\text{\sf 1}}$
if and only if the collection of five vector fields:
\[
\Big\{ \overline{\mathcal{S}},\, \mathcal{S},\, \mathcal{T},\,
\overline{\mathcal{L}},\, \mathcal{L} \Big\},
\]
where:
\[
\aligned \mathcal{T} & :=
i\,\big[\mathcal{L},\overline{\mathcal{L}}\big],
\\
\mathcal{S} & := \big[\mathcal{L},\,\mathcal{T}\big],
\\
\overline{\mathcal{S}} & :=
\big[\overline{\mathcal{L}},\,\mathcal{T}\big],
\endaligned
\]
makes up a frame for $\C \otimes_\R TM^5$.\qed
\end{Proposition}

Having $5$ fields implies
that there are in sum $10$ pairwise
Lie brackets. Thus, in order to determine
the full Lie bracket structure, there
remain $7$ such brackets to be looked at.

The first group of four Lie brackets is:
\[
\big[\mathcal{L},\,\mathcal{S}\big], \ \ \ \ \ \ \ \ \
\big[\overline{\mathcal{L}},\,\mathcal{S}\big], \ \ \ \ \ \ \ \ \
\big[\mathcal{L},\,\overline{\mathcal{S}}\big], \ \ \ \ \ \ \ \ \
\big[\overline{\mathcal{L}},\,\overline{\mathcal{S}}\big].
\]

\begin{Lemma}
There are six complex-valued functions $P$,
$Q$, $R$, $A$, $B$, $C$ defined on $M$ such that:
\[
\aligned \big[\mathcal{L},\,\mathcal{S}\big] & = P\,\mathcal{T} +
Q\,\mathcal{S} + R\,\overline{\mathcal{S}},
\\
\big[\overline{\mathcal{L}},\,\mathcal{S}\big] & = A\,\mathcal{T} +
B\,\mathcal{S} + C\,\overline{\mathcal{S}},
\\
\big[\mathcal{L},\,\overline{\mathcal{S}}\big] & =
\overline{A}\,\mathcal{T} + \overline{C}\,\mathcal{S} +
\overline{B}\,\overline{\mathcal{S}},
\\
\big[ \overline{\mathcal{L}},\,\overline{\mathcal{S}}\big] & =
\overline{P}\,\mathcal{T} + \overline{R}\,\mathcal{S} +
\overline{Q}\,\overline{\mathcal{S}}.
\endaligned
\]
Moreover, the two brackets:
\[
\big[\overline{\mathcal{L}},\mathcal{S}\big]
=
\big[\mathcal{L},\overline{\mathcal{S}}\big],
\]
are real and equal. In particular, $A$ is a real-valued function and
$C=\overline{B}$.
\end{Lemma}

So we have:
\[
\boxed{ \aligned
\big[\mathcal{L},\,\mathcal{S}\big] & = P\,\mathcal{T} +
Q\,\mathcal{S} + R\,\overline{\mathcal{S}},
\\
\big[\overline{\mathcal{L}},\,\mathcal{S}\big] & = A\,\mathcal{T} +
B\,\mathcal{S} + \overline{B}\,\overline{\mathcal{S}},
\\
\big[\mathcal{L},\,\overline{\mathcal{S}}\big] & = A\,\mathcal{T} +
B\,\mathcal{S} + \overline{B}\,\overline{\mathcal{S}},
\\
\big[ \overline{\mathcal{L}},\,\overline{\mathcal{S}}\big] & =
\overline{P}\,\mathcal{T} + \overline{R}\,\mathcal{S} +
\overline{Q}\,\overline{\mathcal{S}}.
\endaligned}
\]

\smallskip

From Section~3 of \cite{ Merker-5-CR-V}, we know that the full expansion
of the above numerator differential polynomials $\Upsilon_1$,
$\Upsilon_2$, $\Upsilon_3$, $\Gamma_i^1$, $\Gamma_i^2$, $\Gamma_i^3$
involve dozens of millions of monomials
in the ${\bf 3} \cdot {\bf 55}$ partial derivatives:
\[
\Big(
\varphi_{1,x^jy^ku_1^{l_1}u_2^{l_2}u_3^{l_3}},\,\,
\varphi_{2,x^jy^ku_1^{l_1}u_2^{l_2}u_3^{l_3}},\,\,
\varphi_{3,x^jy^ku_1^{l_1}u_2^{l_2}u_3^{l_3}}
\Big)_{1\leqslant j+k+l_1+l_2+l_3\leqslant 3},
\]
of the three graphing functions $\varphi_1$,
$\varphi_2$, $\varphi_3$. {\em A fortiori}, the
numerators of the five functions:
\[
\boxed{\,\,
\aligned
P,\ \ \ \ \ \ \ \ \ \ \ \ \ \ \
&
Q,\ \ \ \ \ \ \ \ \ \ \ \ \ \ \
R,\,\,
\\
A,\ \ \ \ \ \ \ \ \ \ \ \ \ \ \
&
B,
\endaligned}
\]
would involve even much more monomials, without there being any
reasonable hope to deal with them systematically on any currently available
mostly powerful computer machine.

\medskip

{\em There is therefore some unavoidable practical
necessity of lowering the ambition of complete expliciteness by
raising up the level of calculations up to the so-denoted
five fundamental functions $P$, $Q$, $R$, $A$, $B$.}

\medskip

It yet remains to compute the $3$ among $10$ structure Lie brackets:
\[
\big[\mathcal{T},\,\mathcal{S}\big],
\ \ \ \ \ \ \ \ \ \ \ \ \
\big[\mathcal{T},\,\overline{\mathcal{S}}\big],
\ \ \ \ \ \ \ \ \ \ \ \ \
\big[\mathcal{S},\,\overline{\mathcal{S}}\big].
\]
A careful systematic
inspection of various Jacobi identities provides their expressions
in terms of the five fundamental functions $P$, $Q$, $R$, $A$, $B$.

\begin{Lemma}
The coefficients of the two Lie brackets:
\[
\aligned \big[\mathcal{T},\,\mathcal{S}\big] & = E\,\mathcal{T} +
F\,\mathcal{S} + G\,\overline{\mathcal{S}},
\\
\big[\mathcal{T},\,\overline{\mathcal{S}}\big] & =
\overline{E}\,\mathcal{T} + \overline{G}\,\mathcal{S} +
\overline{F}\,\overline{\mathcal{S}},
\endaligned
\]
are three complex-valued functions $E$, $F$, $G$ which can be
expressed as follows in terms of $P$, $Q$, $R$, $A$, $B$ and their
first-order frame derivatives:
\[
\aligned
E
&
=
-\,i\,\overline{\mathcal L}(P)
-
i\,A\,Q
-
i\,\overline{P}\,R
+
i\,\mathcal L(A)
+
i\,B\,P
+
i\,A\,\overline{B},
\\
F
&
=
-\,i\,\overline{\mathcal L}(Q)
-
i\,R\,\overline{R}
+
i\,A
+
i\,\mathcal{L}(B)
+
i\,B\,\overline{B},
\\
G
&
=
-\,i\,P
-
i\,\overline{B}\,Q
-
i\,R\,\overline{Q}
-
i\,\overline{\mathcal L}(R)
+
i\,B\,R
+
i\,\overline{B}\overline{B}
+
i\,\mathcal L(\overline{B}),
\endaligned
\]
while the coefficients of the last, tenth structure bracket:
\[
\aligned
\big[\mathcal{S},\,\overline{\mathcal{S}}\big]
&
=
i\,J\,\mathcal{T} + K\,\mathcal{S} -
\overline{K}\,\overline{\mathcal{S}},
\endaligned
\]
are one complex-valued function $K$ and one real-valued function
$J$ which can be expressed as follows in terms of $P$, $Q$, $R$,
$A$, $B$ and their frame derivatives up to order $2$:
\[
\!\!\!\!\!\!\!\!\!\!\!\!\!\!\!\!\!\!\!\!
\footnotesize
\aligned
-\,2\,J
&
=
-\,\overline{\mathcal{L}}\big(\overline{\mathcal{L}}(P)\big)
+
\overline{\mathcal{L}}\big(\mathcal{L}(A)\big)
+
\mathcal{L}\big(\overline{\mathcal{L}}(A)\big)
-
\mathcal{L}\big(\mathcal{L}(\overline{P})\big)
\\
&
\ \ \ \ \
-\,Q\,\overline{\mathcal{L}}(A)
-
2\,A\,\overline{\mathcal{L}}(Q)
-
R\,\overline{\mathcal{L}}(\overline{P})
-
2\,\overline{P}\,\overline{\mathcal{L}}(R)
-
2\,AR\overline{R}
-
2\,P\overline{P}
-
\overline{B}\overline{P}Q
-
\overline{P}\overline{Q}R
-
\\
&
\ \ \ \ \
-
\overline{R}\mathcal{L}(P)
-
2\,P\,\mathcal{L}(\overline{R})
-
\overline{Q}\,\mathcal{L}(A)
-
2\,A\mathcal{L}(\overline{Q})
-
PQ\overline{R}
-
BP\overline{Q}
+
\\
&
\ \ \ \ \
+
2\,P\overline{\mathcal{L}}(B)
+
B\overline{\mathcal{L}}(P)
+
2\,A\overline{\mathcal{L}}(\overline{B})
+
\overline{B}\,\overline{\mathcal{L}}(A)
+
2\,A\,\mathcal{L}(B)
+
2\,AA
+
2\,AB\overline{B}
+
2\,\overline{P}\,\mathcal{L}(\overline{B})
+
\\
&
\ \ \ \ \
+
B\overline{P}R
+
\overline{B}\overline{B}\overline{P}
+
B\,\mathcal{L}(A)
+
\overline{B}\,\mathcal{L}(\overline{P})
+
BBP
+
\overline{B}P\overline{R},
\endaligned
\]
\[
\!\!\!\!\!\!\!\!\!\!\!\!\!\!\!\!\!\!\!\!
\footnotesize
\aligned
2i\,K
&
=
-\,\overline{\mathcal{L}}\big(\overline{\mathcal{L}}(Q)\big)
+
\overline{\mathcal{L}}\big(\mathcal{L}(B)\big)
+
\mathcal{L}\big(\overline{\mathcal{L}}(B)\big)
-
\mathcal{L}\big(\mathcal{L}(\overline{R})\big)
-
\\
&
\ \ \ \ \
-
2\,\overline{R}\,\overline{\mathcal{L}}(R)
-
R\,\overline{\mathcal{L}}(\overline{R})
-
B\,\overline{\mathcal{L}}(Q)
-
BR\overline{R}
-
2\,P\overline{R}
-
\overline{Q}R\overline{R}
-
2\,\mathcal{L}(\overline{P})
-
\overline{R}\,\mathcal{L}(Q)
-
\\
&
\ \ \ \ \
-\,2\,Q\,\mathcal{L}(\overline{R})
-
\overline{Q}\,\mathcal{L}(B)
-
2\,B\,\mathcal{L}(\overline{Q})
-
A\overline{Q}
-
\overline{P}Q
-
QQ\overline{R}
-
BQ\overline{Q}
+
\\
&
\ \ \ \ \
+
2\,\overline{\mathcal{L}}(A)
+
\overline{B}\,\overline{\mathcal{L}}(B)
+
2\,B\,\overline{\mathcal{L}}(\overline{B})
+
3\,B\,\mathcal{L}(B)
+
3\,AB
+
BBQ
+
2\,BB\overline{B}
+
2\,\overline{R}\,\mathcal{L}(\overline{B})
+
\\
&
\ \ \ \ \
+
\overline{B}\overline{B}\overline{R}
+
\overline{B}\,\mathcal{L}(\overline{R})
+
\overline{B}\overline{P}
+
Q\,\overline{\mathcal{L}}(B).
\endaligned
\]
\end{Lemma}

Crucially too, since no existing computer machine is powerful enough
to compute everything in terms of the three graphing functions
$\varphi_1$, $\varphi_2$, $\varphi_3$, one must possess some means of
knowing as many as possible of the differential-algebraic relations
that these $5$ further coefficient-functions $E$, $F$, $G$, $J$, $K$
share in common with $P$, $Q$, $R$, $A$, $B$.
An inspection of higher, length-six, Jacobi
identities already explored in~\cite{ Merker-Sabzevari-CEJM} provides
$5$ such relations labeled $\overset{ 1}{=}$, $\overset{ 2}{=}$,
$\overset{ 3}{=}$, $\overset{ 4}{=}$, $\overset{ 5}{=}$, which will appear
to be very useful and important when performing Cartan's method later
on.
\[
\footnotesize
\aligned
0&\overset{1}{\equiv} 2\,\mathcal
L(\overline{\mathcal L}(P))-\mathcal L(\mathcal
L(A))-\overline{\mathcal L}(\mathcal L(P))-
\\
&2\,P\mathcal L(B)-B\mathcal L(P)-2\,A\mathcal L(\overline
B)-\overline B\mathcal L(A)+P\overline{\mathcal L}(Q)+A\mathcal L(Q)+
\\
&+2\,Q\mathcal L(A)-Q\overline{\mathcal L}(P)+A\overline{\mathcal
L}(R)+2\,R\mathcal L(\overline P)+\overline P\mathcal
L(R)-R\overline{\mathcal L}(A)-
\\
&-PB\overline B-A\overline B^2+PBQ+2\,AQ\overline
B-AQ^2-2\,ABR+2\,RP\overline R+2\,AR\overline Q-QR\overline
P-R\overline B\,\overline P,
\\
 0&\overset{2}{\equiv}2\,\mathcal L(\overline{\mathcal
L}(Q))-\mathcal L(\mathcal L(B))-\overline{\mathcal L}(\mathcal
L(Q))-
\\
&-2\,\mathcal L(A)-2\,B\mathcal L(\overline B)-\overline B\mathcal
L(B)+B\overline{\mathcal L}(R)+2\,R\mathcal L(\overline R)+\overline
R\mathcal L(R)-R\overline{\mathcal L}(B)+\overline{\mathcal L}(P)+
\\
&+2\,R\overline P+BQ\overline B-A\overline B-B\overline
B^2+AQ+QR\overline R+2\,BR\overline Q-2\,B^2R-R\overline B\,\overline
R,
\\
0&\overset{3}{\equiv}2\,\mathcal L(\overline{\mathcal L}(R))-\mathcal
L(\mathcal L(\overline B))-
\overline{\mathcal{L}}\big(\mathcal{L}(R)\big)-
\\
&-3\,\overline B\mathcal L(\overline B)+\overline B\mathcal
L(Q)+2\,Q\mathcal L(\overline B)-2\,R\mathcal L(B)-B\mathcal
L(R)+R\overline{\mathcal L}(Q)+\overline B\,\overline{\mathcal L}(R)+
\\
&+2\,R\mathcal L(\overline Q)+\overline Q\mathcal
L(R)-Q\overline{\mathcal L}(R)-\overline{\mathcal L}\mathcal
L(R)-R\overline{\mathcal L}(\overline B)+\mathcal L(P)+
\\
&+2\,Q\overline B^2-QP-Q^2\overline B-\overline B^3+P\overline
B-2\,AR-2\,BR\overline B+BQR+2\,R^2\overline R+R\overline
B\,\overline Q-QR\overline Q,
\\
0&\overset{4}{\equiv}-3\overline{\mathcal L}(\mathcal L(A))
-\mathcal L(\mathcal
L(\overline P))+3\,\mathcal
L(\overline{\mathcal L}(A))+\overline{\mathcal L}(\overline{\mathcal
L}(P))-
\\
&-2\,A\mathcal L(\overline Q)-\overline Q\mathcal L(A)+3\,B\mathcal
L(A)+3\,\overline B\mathcal L(\overline P)-3\,B\overline{\mathcal
L}(P)-3\,\overline B\,\overline{\mathcal L}(A)+2\,A\overline{\mathcal
L}(Q)-
\endaligned
\]
\[
\footnotesize\aligned
&+Q\overline{\mathcal L}(A)-2\,P\mathcal L(\overline R)-\overline
R\mathcal L(P)+2\,\overline P\,\overline{\mathcal
L}(R)+R\overline{\mathcal L}(\overline P)-
\\
&-BP\overline Q+3\,B^2P+2\,A\overline B\,\overline
Q-2\,BQA-3\,\overline B^2\overline P+Q\overline B\,\overline
P-PQ\overline R+3\,P\overline B\,\overline R-3\,BR\overline
P+R\overline P\,\overline Q,
\\
 0&\overset{5}{\equiv}-3\,\overline{\mathcal L}(\mathcal
L(B))+3\,\mathcal L(\overline{\mathcal L}(B))+\overline{\mathcal
L}(\overline{\mathcal L}(Q))-\mathcal L(\mathcal L(\overline R))+
\\
&+3\,B\mathcal L(B)-3\,\overline B\,\overline{\mathcal
L}(B)+Q\overline{\mathcal L}(B)-B\overline{\mathcal
L}(Q)-2\,Q\mathcal L(\overline R)-\overline R\mathcal L(Q)-
\\
&-2\,B\mathcal L(\overline Q)-\overline Q\mathcal L(B)+3\,\overline
B\mathcal L(\overline R)+2\,\overline R\,\overline{\mathcal
L}(R)+R\overline{\mathcal L}(\overline R)-2\,\mathcal L(\overline P)-
\\
&-Q\overline P-A\overline Q-BQ\overline Q+3\,AB+3\,\overline
B\,\overline P+2\,B\overline B\,\overline Q+B^2Q-Q^2\overline
R+4\,Q\overline B\,\overline R-
\\
&-3\,BR\overline R-3\,\overline B^2\overline R+R\overline
Q\,\overline R.
\endaligned
\]
Higher length $7$ or $8$ Jacobi relations should also be useful
if one wants to explore deeper the problem, but even when admitting
to compute everything only in terms of $P$, $Q$, $R$, $A$, $B$,
the formulas again show off a striking tendency to striking
swelling.

In any case, introducing the coframe of $1$-forms
generating the complexified cotangent bundle $\C \otimes_\R T^*M$:
\[
\big\{\overline{\sigma_0},\,\sigma_0,\,\rho_0,\,\overline{\zeta_0},\,\zeta_0\big\}
\ \ \ \text{\rm which is dual to the frame} \ \ \
\big\{\overline{\mathcal{S}},\,
\mathcal{S},\,\mathcal{T},\,\overline{\mathcal{L}},\,\mathcal{L}\big\},
\]
namely which satisfy by definition:
\[
\begin{array}{ccccc}
\overline{\sigma_0}(\overline{\mathcal{S}})=1 \ \ \ & \ \ \
\overline{\sigma_0}(\mathcal{S})=0 \ \ \ & \ \ \
\overline{\sigma_0}(\mathcal{T})=0 \ \ \ & \ \ \
\overline{\sigma_0}(\overline{\mathcal{L}})=0 \ \ \ & \ \ \
\overline{\sigma_0}\big(\mathcal{L}\big)=0,
\\
\sigma_0(\overline{\mathcal{S}})=0 \ \ \ & \ \ \
\sigma_0(\mathcal{S})=1 \ \ \ & \ \ \ \sigma_0(\mathcal{T})=0 \ \ \ &
\ \ \ \sigma_0(\overline{\mathcal{L}})=0 \ \ \ & \ \ \
\sigma_0\big(\mathcal{L}\big)=0,
\\
\rho_0(\overline{\mathcal{S}})=0 \ \ \ & \ \ \ \rho_0(\mathcal{S})=0
\ \ \ & \ \ \ \rho_0(\mathcal{T})=1 \ \ \ & \ \ \
\rho_0(\overline{\mathcal{L}})=0 \ \ \ & \ \ \ \rho_0(\mathcal{L})=0,
\\
\overline{\zeta_0}(\overline{\mathcal{S}})=0 \ \ \ & \ \ \
\overline{\zeta_0}(\mathcal{S})=0 \ \ \ & \ \ \
\overline{\zeta_0}(\mathcal{T})=0 \ \ \ & \ \ \
\overline{\zeta_0}(\overline{\mathcal{L}})=1 \ \ \ & \ \ \
\overline{\zeta_0}\big(\mathcal{L}\big)=0,
\\
\zeta_0(\overline{\mathcal{S}})=0 \ \ \ & \ \ \
\zeta_0(\mathcal{S})=0 \ \ \ & \ \ \ \zeta_0(\mathcal{T})=0 \ \ \ & \
\ \ \zeta_0(\overline{\mathcal{L}})=0 \ \ \ & \ \ \
\zeta_0\big(\mathcal{L}\big)=1,
\end{array}
\]
one determine its Darboux structure by reading vertically
a convenient auxiliary array:
\[
\footnotesize
\begin{array}{cccccccccccc}
& & \overline{\mathcal{S}} & & \mathcal{S} & & \mathcal{T} & &
\overline{\mathcal{L}} & & \mathcal{L}
\\
& & \boxed{d\overline{\sigma_0}} & & \boxed{d\sigma_0} & &
\boxed{d\rho_0} & & \boxed{d\overline{\zeta_0}} & & \boxed{d\zeta_0}
\\
\big[\overline{\mathcal{S}},\,\mathcal{S}\big] & = &
\overline{K}\cdot\overline{\mathcal{S}} & + &
-\,K\cdot\mathcal{S} & + & -i\,J\cdot\mathcal{T} & + & 0 &
+ & 0 & \boxed{\overline{\sigma_0}\wedge\sigma_0}
\\
\big[\overline{\mathcal{S}},\,\mathcal{T}\big] & = &
-\,\overline{F}\cdot\overline{\mathcal{S}} & + &
-\,\overline{G}\cdot\mathcal{S} & + &
-\,\overline{E}\cdot\mathcal{T} & + & 0 & + & 0 &
\boxed{\overline{\sigma_0}\wedge\rho_0}
\\
\big[\overline{\mathcal{S}},\,\overline{\mathcal{L}}\big] & = &
-\,\overline{Q}\cdot\overline{\mathcal{S}} & + &
-\,\overline{R}\cdot\mathcal{S} & + & -\,\overline{P}\cdot\mathcal{T}
& + & 0 & + & 0 & \boxed{\overline{\sigma_0}\wedge\overline{\zeta_0}}
\\
\big[\overline{\mathcal{S}},\,\mathcal{L}\big] & = &
-\,\overline{B}\cdot\overline{\mathcal{S}} & + &
-\,{B}\cdot\mathcal{S} & + & -A\cdot\mathcal{T} & + & 0 & + & 0 &
\boxed{\overline{\sigma_0}\wedge\zeta_0}
\\
\big[\mathcal{S},\,\mathcal{T}\big] & = &
-\,G\cdot\overline{\mathcal{S}} & + & -\,F\cdot \mathcal{S} & +
& -\,E\cdot\mathcal{T} & + & 0 & + & 0 &
\boxed{\sigma_0\wedge\rho_0}
\\
\big[\mathcal{S},\,\overline{\mathcal{L}}\big] & = &
-\,\overline{B}\cdot\overline{\mathcal{S}} & + & -\,B\cdot
\mathcal{S} & + & -\,A\cdot\mathcal{T} & + & 0 & + & 0 &
\boxed{\sigma_0\wedge\overline{\zeta_0}}
\\
\big[\mathcal{S},\,\mathcal{L}\big] & = &
-\,R\cdot\overline{\mathcal{S}} & + & -\,Q\cdot\mathcal{S} & + &
-\,P\cdot\mathcal{T} & + & 0 & + & 0 & \boxed{\sigma_0\wedge\zeta_0}
\\
\big[\mathcal{T},\,\overline{\mathcal{L}}\big] & = &
-\,\overline{\mathcal{S}} & + & 0 & + & 0 & + & 0 & + & 0 &
\boxed{\rho_0\wedge\overline{\zeta_0}}
\\
\big[\mathcal{T},\,\mathcal{L}\big] & = & 0 & + & -\,\mathcal{S} & +
& 0 & + & 0 & + & 0 & \boxed{\rho_0\wedge\zeta_0}
\\
\big[\overline{\mathcal{L}},\,\mathcal{L}\big] & = & 0 & + & 0 & + &
i\,\mathcal{T} & + & 0 & + & 0 &
\boxed{\overline{\zeta_0}\wedge\zeta_0}
\end{array}
\]
and this provides (minding an overall minus sign due to duality):
\[
\boxed{ \aligned d\overline{\sigma}_0 & =
-\,\overline{K} \cdot \overline{\sigma}_0\wedge\sigma_0 +
\overline{F}\cdot \overline{\sigma}_0\wedge\rho_0 +
\overline{Q}\cdot \overline{\sigma}_0\wedge\overline{\zeta}_0 +
\overline{B}\cdot \overline{\sigma}_0\wedge\zeta_0 +
\\
& \ \ \ \ \ + G\cdot \sigma_0\wedge\rho_0 + \overline{B}\cdot
\sigma_0\wedge\overline{\zeta}_0 + R\cdot \sigma_0\wedge\zeta_0 +
\rho_0\wedge\overline{\zeta}_0,
\\
d\sigma_0 & = K \cdot \overline{\sigma}_0\wedge\sigma_0 +
\overline{G}\cdot \overline{\sigma}_0\wedge\rho_0 +
\overline{R}\cdot \overline{\sigma}_0\wedge\overline{\zeta}_0 +
{B}\cdot \overline{\sigma}_0\wedge\zeta_0 +
\\
& \ \ \ \ \ + F\cdot \sigma_0\wedge\rho_0 + B\cdot
\sigma_0\wedge\overline{\zeta}_0 + Q\cdot \sigma_0\wedge\zeta_0 +
\rho_0\wedge\zeta_0,
\\
d\rho_0 & = i\,J\cdot \overline{\sigma}_0\wedge\sigma_0 +
\overline{E}\cdot \overline{\sigma}_0\wedge\rho_0 +
\overline{P}\cdot \overline{\sigma}_0\wedge\overline{\zeta}_0 +
A\cdot \overline{\sigma}_0\wedge\zeta_0 +
\\
& \ \ \ \ \ + E\cdot \sigma_0\wedge\rho_0 + A\cdot
\sigma_0\wedge\overline{\zeta}_0 + P\cdot \sigma_0\wedge\zeta_0 -
i\,\overline{\zeta}_0\wedge\zeta_0,
\\
d\overline{\zeta}_0 & = 0,
\\
d\zeta_0 & = 0.
\endaligned}
\]

Recall that exactly as for the model reviewed above,
the initial ambiguity matrix associated to the equivalence problem
under local biholomorphic transformations for maximally minimal
CR-generic $3$-codimensional submanifolds $M^5 \subset \C^4$ is of the
general form:
\[
\left(\!\!
\begin{array}{ccccc}
{\sf a}\overline{\sf a}\overline{\sf a} & 0 & \overline{\sf c} &
\overline{\sf e} & \overline{\sf d}
\\
0 & {\sf a}{\sf a}\overline{\sf a} & {\sf c} & {\sf d} & {\sf e}
\\
0 & 0 & {\sf a}\overline{\sf a} & \overline{\sf b} & {\sf b}
\\
0 & 0 & 0 & \overline{\sf a} & 0
\\
0 & 0 & 0 & 0 & {\sf a}
\end{array}
\!\!\right),
\]
where ${\sf a}$, ${\sf b}$, ${\sf c}$, ${\sf e}$, ${\sf d}$ are
complex numbers.
The so-called {\sl lifted coframe}
is then (one must transpose the above matrix):
\[
\aligned \left(\!\!
\begin{array}{c}
\overline{\sigma}
\\
\sigma
\\
\rho
\\
\overline{\zeta}
\\
\zeta
\end{array}
\!\!\right)
:=
\left(\!\!
\begin{array}{ccccc}
{\sf a}\overline{\sf a}\overline{\sf a} & 0 & 0 & 0 & 0
\\
0 & {\sf a}{\sf a}\overline{\sf a} & 0 & 0 & 0
\\
\overline{\sf c} & {\sf c} & {\sf a}\overline{\sf a} & 0 & 0
\\
\overline{\sf e} & {\sf d} & \overline{\sf b} & \overline{\sf a} & 0
\\
\overline{\sf d} & {\sf e} & {\sf b} & 0 & {\sf a}
\end{array}
\!\!\right)
\left(\!\!
\begin{array}{c}
\overline{\sigma_0}
\\
\sigma_0
\\
\rho_0
\\
\overline{\zeta_0}
\\
\zeta_0
\end{array}
\!\!\right),
\endaligned
\]
that is to say:
\[
\aligned \overline{\sigma} & = {\sf a}\overline{\sf a}\overline{\sf
a}\, \overline{\sigma_0},
\\
\sigma & = {\sf a}{\sf a}\overline{\sf a}\,\sigma_0,
\\
\rho & = \overline{\sf c}\,\overline{\sigma_0} + {\sf c}\,\sigma_0 +
{\sf a}\overline{\sf a}\,\rho_0,
\\
\overline{\zeta} & = \overline{\sf e}\,\overline{\sigma_0} + {\sf
d}\,\sigma_0 + \overline{b}\,\rho_0 + \overline{\sf
a}\,\overline{\zeta_0},
\\
\zeta & = \overline{\sf d}\,\overline{\sigma_0} + {\sf e}\,\sigma_0 +
{\sf b}\,\rho_0 + {\sf a}\,\zeta_0.
\endaligned
\]

To launch Cartan's method, with the Maurer-Cartan forms
(not writing their conjugates):
\[
\aligned {\alpha_1} & := \frac{d{\sf a}}{\sf a},
\\
{\alpha_2} &:= \frac{d{\sf c}}{{\sf a}^2\overline{\sf a}} -
\frac{{\sf c}\,d{\sf a}}{{\sf a}^3\overline{\sf a}} - \frac{{\sf
c}\,d\overline{\sf a}}{{\sf a}^2 \overline{\sf a}^2},
\\
{\alpha_3} & := -\,\frac{{\sf c}\,d{\sf b}}{{\sf a}^3\overline{\sf
a}^2} + \bigg( \frac{{\sf b}{\sf c}}{{\sf a}^4\overline{\sf a}^2} -
\frac{{\sf e}}{{\sf a}^3\overline{\sf a}} \bigg)\,d{\sf a} +
\frac{1}{{\sf a}^2\overline{\sf a}}\,d{\sf e},
\\
{\alpha_4} & := \frac{d\overline{\sf d}}{{\sf
a}\overline{\sf a}^2} - \frac{\overline{\sf c}\,d{\sf b}}{{\sf
a}^2\overline{\sf a}^3} + \bigg( \frac{{\sf b}\overline{\sf c}}{{\sf
a}^3\overline{\sf a}^3} - \frac{\overline{\sf d}}{{\sf
a}^2\overline{\sf a}^2} \bigg)\, d{\sf a},
\\
{\alpha_5} & := \frac{d{\sf b}}{{\sf a}\overline{\sf a}} - \frac{{\sf
b}\,d{\sf a}}{{\sf a}^2\overline{\sf a}},
\endaligned
\]
one computes the complete structure equations:
\begin{equation}
\aligned d\sigma & =
\big(2\,\alpha_1+\overline{\alpha}_1\big) \wedge \sigma+
\\
& \ \ \ \ \
+
U_1\,\sigma\wedge\overline{\sigma} +
U_2\,\sigma\wedge\rho + U_3\,\sigma\wedge\zeta +
U_4\,\sigma\wedge\overline{\zeta}
+
\\
&
\ \ \ \ \ \ \ \ \ \ \ \ \ \ \ \ \ \ \ \ \ \ \ \
+ U_5\,\overline{\sigma}\wedge\rho +
U_6\,\overline{\sigma}\wedge\zeta
+
U_7\,\overline{\sigma}\wedge\overline{\zeta}
+
\\
&
\ \ \ \ \ \ \ \ \ \ \ \ \ \ \ \ \ \ \ \ \ \
\ \ \ \ \ \ \ \ \ \ \ \ \ \ \ \ \ \ \ \ \
+ {\rho\wedge\zeta},
\endaligned
\end{equation}
\[
\aligned
d\rho & = \alpha_2\wedge\sigma +
\overline{\alpha}_2\wedge\overline{\sigma} + \alpha_1\wedge\rho +
\overline{\alpha}_1\wedge\rho +
\\
& \ \ \ \ \ +
V_1\,\sigma\wedge\overline{\sigma} +
V_2\,\sigma\wedge\rho + V_3\,\sigma\wedge\zeta +
V_4\,\sigma\wedge\overline{\zeta} +
\\
& \ \ \ \ \ \ \ \ \ \ \ \ \ \ \ \ \ \ \ \ \ \ \ \ +
\overline{V_2}\,\overline{\sigma}\wedge\rho +
\overline{V_4}\,\overline{\sigma}\wedge\zeta +
\overline{V_3}\,\overline{\sigma}\wedge\overline{\zeta} +
\\
& \ \ \ \ \ \ \ \ \ \ \ \ \ \ \ \ \ \ \ \ \ \ \ \ \ \ \ \ \ \ \ \ \ \
\ \ \ \ \ \ \ \ \ + V_8\,\rho\wedge\zeta +
\overline{V_8}\,\rho\wedge\overline{\zeta} +
\\
& \ \ \ \ \ \ \ \ \ \ \ \ \ \ \ \ \ \ \ \ \ \ \ \ \ \ \ \ \ \ \ \ \ \
\ \ \ \ \ \ \ \ \ \ \ \ \ \ \ \ \ \ \ \ \ \ \ \ \ \ \ \ \ \ +
i\,\zeta\wedge\overline{\zeta},
\endaligned
\]
\[
\aligned
d\zeta & = \alpha_3\wedge\sigma
+ \alpha_4\wedge\overline{\sigma} + \alpha_5\wedge\rho +
\alpha_1\wedge\zeta +
\\
& \ \ \ \ \ + W_1\,\sigma\wedge\overline{\sigma} +
W_2\,\sigma\wedge\rho + W_3\,\sigma\wedge\zeta +
W_4\,\sigma\wedge\overline{\zeta} +
\\
& \ \ \ \ \ \ \ \ \ \ \ \ \ \ \ \ \ \ \ \ \ \ \ \ +
W_5\,\overline{\sigma}\wedge\rho + W_6\,\overline{\sigma}\wedge\zeta
+ W_7\,\overline{\sigma}\wedge\overline{\zeta} +
\\
& \ \ \ \ \ \ \ \ \ \ \ \ \ \ \ \ \ \ \ \ \ \ \ \ \ \ \ \ \ \ \ \ \ \
\ \ \ \ \ \ \ \ \ + W_8\,\rho\wedge\zeta +
W_9\,\rho\wedge\overline{\zeta} +
\\
& \ \ \ \ \ \ \ \ \ \ \ \ \ \ \ \ \ \ \ \ \ \ \ \ \ \ \ \ \ \ \ \ \ \
\ \ \ \ \ \ \ \ \ \ \ \ \ \ \ \ \ \ \ \ \ \ \ \ \ \ \ \ \ +
W_{10}\,\zeta\wedge\overline{\zeta},
\endaligned
\]
expressed in terms of the lifted coframe; the explicit expressions
of the initial torsion coefficients $U_i$, $V_j$, $W_k$
are not reviewed in this introductory presentation, but they
appear in Section~\ref{abs-norm}.

\begin{Lemma}
In the first loop of Cartan's method, five essential
linear combinations of torsion coefficients appear:
\[
\aligned
U_5
&
=
\frac{1}{\overline{\sf a}^2}\,\overline{G}
-
\frac{\sf b}{{\sf a}\overline{\sf a}^2}\,B
-
\frac{\overline{\sf b}}{\overline{\sf a}^3}\,
\overline{R}
+
\frac{\overline{\sf d}}{{\sf a}\overline{\sf a}^2},
\\
U_6
&
=
\frac{1}{\overline{\sf a}}\,B
-
\frac{\overline{\sf c}}{{\sf a}\overline{\sf a}^2},
\\
U_7
&
=
\frac{{\sf a}}{\overline{\sf a}^2}\,\overline{R},
\\
U_3+\overline{U}_4-3\,V_8
&
=
\frac{1}{\sf a}\,Q
-
4\,\frac{\sf c}{{\sf a}^2\overline{\sf a}}
+
\frac{1}{\sf a}\,\overline{B}
-
3i\,
\frac{\overline{\sf b}}{{\sf a}\overline{\sf a}},
\\
\overline{U}_4-V_8-\overline{W}_{10}
&
=
\frac{1}{\sf a}\,\overline{B}
-
\frac{\sf c}{{\sf a}^2\overline{\sf a}}.
\endaligned
\]
\end{Lemma}

Precisely, this means that assigning the values:
\[
\aligned
0
&
=
U_5,
\\
0
&
=
U_6,
\\
1
&
=
U_7,
\\
0
&
=
U_3+\overline{U}_4-3\,V_8,
\\
0
&
=
\overline{U}_4
-
V_8
-
\overline{W}_{10},
\endaligned
\]
one can use these $5$ equations to {\sl normalize}
some of the group parameters.

Here, the third equation plays a special role, because when
the function $R \not\equiv 0$ is not identically zero,
one can then {\em normalize the special parameter ${\sf a}$
lying on the diagonal of the matrix group}.

In this case, it is rather easy to realize that one can
construct an absolute parallelism on the basis $M^5$,
which is $5$-dimensional (more will be said in a while).

Thus, we distinguish two branches in Cartan's method:

\smallskip\noindent$\square$\,
$R\equiv 0$,

\smallskip\noindent$\square$\,
$R \not\equiv 0$,

\smallskip\noindent
and we begin by exploring the first one.

\subsection{The branch $R = 0$} In this case, the above
$5$ normalizable expressions reduce to $3$ and we perform
the normalizations:
\[
\boxed{{\sf c}
:=
{\sf a}\overline{\sf a}\,\overline{B}.}
\]
\[
\boxed{{\sf b}
:=
{\sf a}\,
\big(
-i\,B
+
{\textstyle{\frac{i}{3}}}\,\overline{Q}
\big).
}
\]
\[
\boxed{
{\sf d}
=
\overline{\sf a}\,
\big(
-\,i\,\mathcal{L}\big(\overline{B}\big)
+
i\,P
+
{\textstyle{\frac{2i}{3}}}\,\overline{B}\,Q
\big).
}
\]

\begin{Lemma}
Within the branch $R \equiv 0$, after determining the
three group parameters ${\sf b}$, ${\sf c}$, ${\sf d}$
accordingly,
the further essential torsion functions become:
\[
\aligned
&
V_1,\ \ \ \ \ \ \ \ \ \ \ \ \
V_3,\ \ \ \ \ \ \ \ \ \ \ \ \
V_4,
\\
&
W_5,\ \ \ \ \ \ \ \ \ \ \ \ \
W_7,\ \ \ \ \ \ \ \ \ \ \ \ \
\\
&
X_1:=U_2+2\,W_8+\overline{W}_8,\ \ \ \ \ \ \ \ \ \ \ \ \
X_2:=\overline{U}_1+V_2+\overline{W}_6,
\endaligned
\]
and $V_4$ enables one to determine:
\[
\boxed{\,
{\sf e}
:=
{\sf a}\cdot
\big(
i\,\overline{\mathcal{L}}(\overline{B})
-
i\,A
-
2i\,B\overline{B}
+
{\textstyle{\frac{i}{3}}}\,
\overline{B}\,\overline{Q}
\big),\,
}
\]
while in fact:
\[
\boxed{V_3\equiv 0,}\ \ \ \ \ \ \ \ \
\boxed{W_7\equiv 0,}\ \ \ \ \ \ \ \ \
\boxed{X_2\equiv 0,}
\]
and while:
\[
\aligned
V_1
&
=
\frac{1}{{\sf a}^2\overline{\sf a}^2}\,
\Big(
\text{\footnotesize\sf long polynomial in the}\,\,
\big\{\mathcal{L},\overline{\mathcal{L}}\big\}
\text{\footnotesize\sf -derivatives of}\,\,
P,\,\,Q,\,\,R,\,\,A,\,\,B
\Big),
\\
W_5
&
=
\frac{1}{{\sf a}\overline{\sf a}^3}\,
\Big(
\text{\footnotesize\sf long polynomial in the}\,\,
\big\{\mathcal{L},\overline{\mathcal{L}}\big\}
\text{\footnotesize\sf -derivatives of}\,\,
P,\,\,Q,\,\,R,\,\,A,\,\,B
\Big),
\\
W_9
&
=
\frac{i}{9\overline{\sf a}^2}\bigg(18\,\overline{\mathcal
L}(B)-3\,\overline{\mathcal L}(\overline Q)-9\,\overline
P-12\,B\overline Q+9\,B^2+\overline Q^2\bigg),
\\
X_1
&
=
-\,\frac{i}{9{\sf a}\overline{\sf a}}\bigg(6\,\overline{\mathcal
L}(Q)+6\,\mathcal L(\overline Q)-18\,\mathcal
L(B)-18\,\overline{\mathcal L}(\overline B)+
\\
&
\ \ \ \ \ \ \ \ \ \ \ \ \ \ \ \ \ \ \
+
27\,B\overline
B-6\,BQ-6\,\overline B\,\overline Q+2\,Q\overline Q+9\,A\bigg).
\endaligned
\]
\end{Lemma}

In the next loop, new essential torsion coefficients appear:
\[
\aligned
&
W_1,\ \ \ \ \ \ \ \ \ \ \ \ \
W_2,\ \ \ \ \ \ \ \ \ \ \ \ \
W_4,
\\
&
Y:=V_2-W_3+\overline{W}_6,
\endaligned
\]
with:
\[
\aligned
W_1
&
=
\frac{1}{{\sf a}^2\overline{\sf a}^3}\,
\Big(
\text{\footnotesize\sf very long polynomial in the}\,\,
\big\{\mathcal{L},\overline{\mathcal{L}}\big\}
\text{\footnotesize\sf -derivatives of}\,\,
P,\,\,Q,\,\,R,\,\,A,\,\,B
\Big),
\\
W_2
&
=
\frac{1}{{\sf a}^2\overline{\sf a}^2}\,
\Big(
\text{\footnotesize\sf very long polynomial in the}\,\,
\big\{\mathcal{L},\overline{\mathcal{L}}\big\}
\text{\footnotesize\sf -derivatives of}\,\,
P,\,\,Q,\,\,R,\,\,A,\,\,B
\Big),
\endaligned
\]
but we easily show that:
\[
\boxed{\,Y\equiv 0,\,}
\]
while after some demanding computational explorations:
\[
W_4
=
-\,
\frac{6\,\overline B-2\,Q}{5\,\sf a}\,
W_9
+
\frac{\overline Q}{5\overline{\sf a}}\,
X_1
+
\frac{3}{5\,\sf a}
\mathcal{L}\big(W_9\big)
-
\frac{3}{5\overline{\sf a}}\overline{\mathcal L}\big(X_1\big).
\]

In quite summarized words, our first main
theorem is
as follows (more details about the ramification
of possible $\{e \}$-structures appear in the
last two sections of the memoir, and
more is also said about explicit curvatures).

\begin{Theorem}
Within the branch $R = 0$, the extrinsic biholomorphic equivalence
problem or
the intrinsic CR-equivalence problem for real analytic CR-generic
submanifolds $M^5 \subset \C^4$ that are maximally minimal in the
above sense, or equivalently that belong to the General Class
$\text{\sf III}_{\text{\sf 1}}$, reduces to various absolute parallelisms
namely to $\{ e\}$-structures on certain manifolds of dimension $6$,
or directly on the $5$-dimensional basis $M$, unless all existing
potentially normalizable torsion coefficients vanish identically, in
which case $M$ is (locally) biholomorphic to Beloshapka's cubic model
$M_{\sf c}^5$ with a characterization of such a condition being
explicit in terms of the five fundamental functions $P$, $Q$, $R$,
$A$, $B$.
\end{Theorem}

\subsection{The branch $R \not\equiv 0$}
In this branch, introducing a (locally defined,
after relocalization near points where $R \neq 0$ is
truly non-vanishing) real analytic function
${\bf A}_0$ satisfying:
\[
\frac{({\bf A}_0)^2}{\overline{\bf A}_0}
=
R,
\]
one can immediately also normalize ${\sf a}$, already during the first loop:
\begin{equation}
\boxed{\aligned
{\sf a}&:={\bf A}_0,
\\
{\sf b}&:={\bf A}_0\,\big(-i\,B+\textstyle{\frac{i}{3}}\,\overline
Q\big),
\\
{\sf c}&:={\bf A}_0\,{\bf \overline A}_0\,B,
\\
{\sf d}&:=\overline{\bf A}_0\,\Big(i\,\overline{\mathcal
L}(R)-i\,\mathcal L(\overline B)+\textstyle{\frac{4i}{3}}\,\overline
QR+iP+\textstyle{\frac{2i}{3}}\,\overline BQ-2i\,BR\Big).
\endaligned}
\end{equation}

Quite similarly to the branch $R \equiv 0$, the last group parameter
${\sf e}$ can also be normalized:
\[\boxed{\aligned
 {\sf e}=-\textstyle{\frac{i}{3}}\,{\bf A}_0\,\big(6\,B\overline
B-3\,\overline{\mathcal L}(\overline B)+3\,A-\overline B\,\overline
Q\big),
\endaligned}
\]
just because the
essential torsion coefficient $V_4$ stays essentially the same.

\smallskip

Our second and last main theorem therefore is as follows.

\begin{Theorem}
Within the branch $R \not\equiv 0$, the extrinsic
biholomorphic equivalence or the intrinsic CR-equivalence
problem for real analytic CR-generic submanifolds
$M^5 \subset \C^4$ that are maximally minimal in the
above sense, or else that belong to the
General Class $\text{\sf III}_{\text{\sf 1}}$,
always reduces to an absolute parallelisms on the $5$-dimensional
basis $M$.
\end{Theorem}

\bigskip

\begin{center}
\begin{minipage}[t]{11.75cm}
\baselineskip=0.35cm {\scriptsize

\centerline{\bf Table of contents}

\smallskip

{\bf \ref{introduction}. Introduction
\dotfill~\pageref{introduction}.}

{\bf \ref{standard-geometry}. Standard geometry of real affine planes
in $\C^N$ \dotfill~\pageref{standard-geometry}.}

{\bf \ref{Zariski-generic-CR}. Zariski-Generic CR behavior of real analytic
submanifolds in $\C^N$ \dotfill~\pageref{Zariski-generic-CR}.}

{\bf \ref{real-complex-defining}.
CR-Generic submanifolds $M\subset\mathbb{C}^{n+d}$:
real and complex
\dotfill~\pageref{real-complex-defining}.}

{\bf \ref{Heisengerg-Beloshapka}.
Heisenberg sphere in $\C^2$ and Beloshapka's higher dimensional models
\dotfill~\pageref{Heisengerg-Beloshapka}.}

{\bf \ref{nilpotent-5}. Symbol algebra and nilpotent Lie algebras up to dimension $5$
\dotfill~\pageref{nilpotent-5}.}

{\bf \ref{infinitesimal-CR}. Infinitesimal CR automorphisms:
$\mathfrak{ aut}_{CR} ( M) = {\rm Re} (\mathfrak{ hol}(M))$
\dotfill~\pageref{infinitesimal-CR}.}

{\bf \ref{geometric-invariants}. Geometric and analytic invariants of
CR-generic submanifolds $M \subset \C^{n+d}$
\dotfill~\pageref{geometric-invariants}.}

{\bf \ref{ineffective-access}. Ineffective access to the local Lie group
structure \dotfill~\pageref{ineffective-access}.}

{\bf \ref{infinitesimal-model}. Algebra of infinitesimal CR automorphisms
of the cubic model $M_{\sf c}^5 \subset \C^4$
\dotfill~\pageref{infinitesimal-model}.}

{\bf \ref{Tanaka-prolongations}. Tanaka prolongations
\dotfill~\pageref{Tanaka-prolongations}.}

{\bf \ref{equivalence-model}. Equivalence computations for the cubic model
\dotfill~\pageref{equivalence-model}.}

{\bf \ref{initial-complex-frame}. Initial complex frame
for geometry-preserving deformations of the model
\dotfill~\pageref{initial-complex-frame}.}

{\bf \ref{passage-coframe}. Passage to a dual coframe
and its Darboux-Cartan structure
\dotfill~\pageref{passage-coframe}.}

{\bf \ref{abs-norm}. Absorption and normalization
\dotfill~\pageref{abs-norm}.}

{\bf \ref{branch-R-0}. The branch $R\equiv 0$
\dotfill~\pageref{branch-R-0}.}

{\bf \ref{four-normalizations}. Four group parameter general normalizations
\dotfill~\pageref{four-normalizations}.}

{\bf \ref{branch-R-neq-0}. The branch $R\not\equiv 0$
\dotfill~\pageref{branch-R-neq-0}.}

}\end{minipage}
\end{center}

\section{Standard geometry of real affine planes in $\C^N$}
\label{standard-geometry}
\HEAD{\ref{standard-geometry}.~Standard geometry
of real affine planes in $\C^N$}{
Masoud {\sc Sabzevari} (Shahrekord) and Jo\"el {\sc Merker} (LM-Orsay)}

\subsection{Standard complex structure on $T\C^N$}
Let $N \geqslant 1$ be a positive integer and consider the complex
Euclidean space $\C^N$, equipped with the canonical coordinates
$(z_1, z_2, \dots, z_N)$, where the complex numbers:
\[
z_k = x_k + i\,y_k,
\]
for $k = 1, 2, \dots, N$, belong to $\C$, with:
\[
x_k = {\rm Re}\,z_k \ \ \ \ \ \ \ \ \ \ \ \ \text{\rm and} \ \ \ \ \
\ \ \ \ \ \ \ y_k = {\rm Im}\,z_k.
\]
The (real) tangent bundle $T\C^N$ is generated by the obvious frame
constituted by the $2N$ basic vector fields:
\[
{\textstyle{\frac{\partial}{\partial x_1}}},\ \
{\textstyle{\frac{\partial}{\partial y_1}}},\ \
{\textstyle{\frac{\partial}{\partial x_2}}},\ \
{\textstyle{\frac{\partial}{\partial y_2}}}, \ \ \cdots\cdots, \ \
{\textstyle{\frac{\partial}{\partial x_N}}},\ \
{\textstyle{\frac{\partial}{\partial y_N}}}.
\]
At an arbitrary fixed point $p \in \C^N$, a general vector tangent to
$\C^N$ therefore writes:
\[
\mathmotsf{vect}_p = {\sf x}_1\,{\textstyle{\frac{\partial}{\partial
x_1}}} \big\vert_p + {\sf y}_1\,{\textstyle{\frac{\partial}{\partial
y_1}}} \big\vert_p +\cdots\cdots+ {\sf
x}_N\,{\textstyle{\frac{\partial}{\partial x_N}}} \big\vert_p + {\sf
y}_N\,{\textstyle{\frac{\partial}{\partial y_N}}} \big\vert_p,
\]
where ${\sf x}_1, {\sf y}_1, \dots, {\sf x}_N, {\sf y}_N$ are free
real numbers.

Now, according to a standard definition, the {\sl complex structure}
$J$ is the endomorphism of $T\C^N$ which\,\,---\,\,in purely real
language\,\,---\,\,indicates how such a vector, when written using
complex notation:
\[
\mathmotsf{vect}_p = \big({\sf x}_1+i\,{\sf y}_1,\dots,{\sf
x}_N+i\,{\sf y}_N\big),
\]
is transformed after multiplication by $\sqrt{ -1}$, namely:
\[
i\,\mathmotsf{vect}_p = \big(-{\sf y}_1+i\,{\sf x}_1,\dots,-{\sf
y}_N+i\,{\sf x}_N\big),
\]
that is to say, in terms of the basic frame (without restricting the
considerations at a fixed point), $J$ is defined by:
\[
\aligned
J\big({\textstyle{\frac{\partial}{\partial x_k}}}\big)
&
:=
{\textstyle{\frac{\partial}{\partial y_k}}}
\ \ \ \ \ \ \ \ \ \ \ \ \ \ \ \
{\scriptstyle{(k\,=\,1\,\cdots\,N)}},
\\
J\big({\textstyle{\frac{\partial}{\partial y_k}}}\big)
&
:=
-{\textstyle{\frac{\partial}{\partial x_k}}}
\ \ \ \ \ \ \ \ \ \ \ \ \
{\scriptstyle{(k\,=\,1\,\cdots\,N)}}.
\endaligned
\]
Of course, one has:
\[
J^2= -\,\mathmotsf{Id}.
\]
In the case $N = 1$ of classical Complex Analysis, $J$ identifies
with $\frac{ \pi}{ 2}$-rotation of tangent vectors.

\subsection{Study of real affine subspaces of $\C^N$}
Next, given a real affine subspace $H_p \subset T_p \C^N$, a real
vector $\mathmotsf{vect}_p \in H_p$ will be a vector belonging to
some complex affine line\,\,---\,\,{\em i.e.} to some one-dimensional
$\C$-linear affine subspace of $\C^N$\,\,---\,\,if and only if
$J({\sf vect}_p)$ also belongs to $H_p$. Further, using $J^2 = - {\rm
Id}$, one easily convinces oneself that the intersection:
\[
H_p\cap J(H_p) =: H_p^c
\]
is the {\em largest} $J$-invariant real linear subspace of $H_p$,
hence also, it is the largest space in $H_p$ on which one can define
a true structure of a {\em complex} vector space. In fact, it
suffices to define the complex scalar multiplication of an arbitrary
vector simply by:
\[
(a+i\,b)\,\mathmotsf{vect}_p := a\,\mathmotsf{vect}_p +
b\,J(\mathmotsf{vect}_p),
\]
and one then verifies at once that this indeed equips $H_p^c$ with a
structure of vector space over $\C$. In particular, $H_p^c$ is
even-dimensional. The upper index $c$ reminds that $H_p^c$ is usually
called the {\em complex subspace} of $H_p$.

On the other hand and somehow `dually', the sum:
\[
H_p^{i_c} = H_p+J(H_p)
\]
is the smallest $J$-invariant real linear subspace of $T_p\C^N$ in
which $H_p$ is contained. Again, this space possesses a $\C$-linear
structure, hence it is even-dimensional too. The upper index $i_c$
reminds that $H_p^{ i_c}$ is usually called the {\em intrinsic
complexification} of $H_p$.

Thus in summary, we have the two inclusions:
\[
H_p^c \subset H_p \subset H_p^{i_c}
\]
inside $T_p \C^N$.

Concerning dimensions, we must introduce appropriate notations. We
will call $2n$ the real dimension $\dim_\R H_p^c$ and $2d$ the excess
dimension of $H_p^{ i_c}$ over $H_p^c$, namely:
\[
2n := \dim_\R H_p^c, \ \ \ \ \ \ \ \ \ \ 2n+2d := \dim_\R H_p^{i_c},
\]
with $n \geqslant 0$ and $d \geqslant 0$ in full generality. From the
definitions, it follows that there are exactly $d$ linearly
independent vectors ${\sf v}_1, \dots, {\sf v}_p$ in $H_p$ that are
also independent from $H_p^c$ such that the $d$ (linearly
independent) vectors $J ( {\sf v}_1)$, \dots, $J ( {\sf v}_d)$ are
{\em not} contained in $H_p$, from which it follows that:
\[
H_p^{i_c} = H_p \oplus \R J({\sf v}_1) \oplus\cdots\oplus \R J({\sf
v}_d).
\]
With these notations, we therefore have:
\[
\dim_\R H_p^c = 2n, \ \ \ \ \ \dim_\R H_p = 2n+d, \ \ \ \ \ \dim_\R
H_p^{i_c} = 2n+2d.
\]
Of course, $n + d \leqslant N$. Obviously, the (real) codimension of
$H_p$ in $T_p \C^N$ is then equal to:
\[
{\rm codim}_\R\,H_p = 2N-2n-d.
\]

Geometrically speaking, one should view our initial general real
vector subspace $H_p$ as sitting inside the complex vector subspace
$H_p^{ i_c} \simeq \C^{ n+d}$ and also as containing within itself
the complex vector subspace $H_p^c \simeq \C^n$. In fact and more
precisely, the true model geometric picture is:
\[
\C^n \subset \R^{2n+d} \subset \C^{n+d} \subset \C^N,
\]
as the following elementary lemma {\rm (\cite{Chirka-1991,
Boggess-1991})}, which we will reprove, confirms.

\begin{Lemma}
\label{H-z-w-t} With $H_p \subset T_p \C^N$ as above being an
arbitrary real affine subspace, there exist $N$ affine {\em complex}
coordinates:
\[
\big(z_1,\dots,z_n,w_1,\dots,w_d,t_1,\dots,t_{N-n-d}\big)
\]
centered at $p$ such that $H_p$ is represented by the following
combination of real and complex Cartesian linear equations:
\[
0 = {\rm Im}\,w_1 =\cdots= {\rm Im}\,w_d \ \ \ \ \ \text{\rm and} \ \
\ \ \ 0 = t_1 =\cdots= t_{N-n-d}.
\]
\end{Lemma}

\proof Firstly, one brings the complex-linear space $H_p^{ i_c}
\simeq \C^{ n+d }$ just to $\{ t_1 = \cdots = t_{ N - n - d} = 0\}$
using a complex-linear straightening. Secondly, inside this space
$\C^{ n + d}$, using another complex-linear straightening, one brings
$H_p^c \simeq \C^n$ to just $\{ w_1 = \cdots = w_d = 0 \}$ so that
$H_p^c$ is then spanned by the $z_1, \dots, z_n$ directions. Thirdly,
and lastly, again inside the $\C^{ n + d} = H_p^{ i_c}$, the affine
subspace $H_p$ under study is spanned by the $z_1, \dots, z_n$
directions and by yet $d$ real directions lying in $\{ 0 \} \times
\C^d$, hence $H_p$ is represented by $d$ linearly independent
Cartesian equations without any $z$ that are necessarily of the form:
\[
0 = \sum_{k=1}^d\,\alpha_{j,k}\,u_k + \sum_{k=1}^d\,\beta_{j,k}\,v_k
\ \ \ \ \ \ \ \ \ \ \ \ \ {\scriptstyle{(j\,=\,1\,\cdots\,d)}},
\]
where $w_k =: u_k + i \, u_k$, for certain real constants $\alpha_{
j, k}$ and $\beta_{ j,k}$. But such equations are visibly equivalent
to:
\[
0 = \sum_{k=1}^d\, {\rm Im} \big[
(\beta_{j,k}+i\,\alpha_{j,k})(u_k+i\,v_k) \big] \ \ \ \ \ \ \ \ \ \ \
\ \ {\scriptstyle{(j\,=\,1\,\cdots\,d)}},
\]
hence it suffices to make the $\C$-linear change of coordinates:
\[
w_j' = \sum_{k=1}^d\, \big(\beta_{j,k} + i\,\alpha_{ j,k}\big)\, w_k
\ \ \ \ \ \ \ \ \ \ \ \ \ {\scriptstyle{(j\,=\,1\,\cdots\,d)}}
\]
---\,\,the determinant
is again nonzero\,\,---\,\,in order to represent $H_p$ inside $H_p^{
i_c} \simeq \C^{ n + d}$ by just $0 = {\rm Im}\, w_1' = \cdots = {\rm
Im}\, w_d'$, as was to be proved.
\endproof

From this linear normalization lemma, we clearly see that the
quotient real vector space $T_p \C^N \big/ H_p$\,\,---\,\,which
somehow represents the `external', `normal'
space\,\,---\,\,decomposes as a direct sum of a complex vector space
$T_p \C^N \big/ H_p^{ i_c} \simeq \C^{ N - n - d}$ plus a real vector
space $H_p^{ i_c} / H_p \simeq \R^d$, while the original linear
subspace $H_p$ also decomposes as a direct sum of a complex vector
space $H_p^c \simeq \C^n$ plus a real vector space $H_p \big/ H_p^c
\simeq \R^d$, the two complex dimensions $N - n - d$ and $n$ being in
general distinct. We may therefore express a bit differently the view
that $H_p \simeq \C^n \times \R^d \times \{ 0\}^{ N - d - d}$.

\begin{Lemma}
\label{dimension-formula} To any arbitrary real affine subspace $H_p
\subset T_p \C^N$ are simultaneously associated its largest complex
subspace $H_p^c$ and its intrinsic complexification $H_p^{ i_c}$
satisfying:
\[
H_p^c = H_p\cap J(H_p) \subset H_p \subset H_p+J(H_p) = H_p^{i_c},
\]
where the two extra dimensions $\dim_\R \big( H_p / H_p^c)$ and
$\dim_\R \big( H_p^{ i_c} / H_p\big)$ coincide (underlined terms):
\begin{equation}
\label{extra-dimensions-coincide} \left\{ \aligned \dim_\R H_p & =
\dim_\R H_p^c + \underline{\big( \dim_\R H_p - \dim_\R H_p^c \big)}
\\
\dim_\R H_p^{i_c} & = \dim_\R H_p + \underline{\big( \dim_\R H_p -
\dim_\R H_p^c \big)}.
\endaligned\right.
\qed
\end{equation}
\end{Lemma}

In fact, this coincidence:
\[
\dim_\R H_p - \dim_\R H_p^c = \dim_\R H_p^{i_c} - \dim_\R H_p
\]
can also be seen by observing that the complex structure induces a
general isomorphism:
\[
J\colon H_p\big/H_p^c \overset{\simeq}{\longrightarrow}
H_p^{i_c}\big/H_p,
\]
because quotients match through: $J ( H_p^c) \subset H_p$, and
because $J$ is invertible, for $J^2 = - {\rm Id}$.

\smallskip

At present, let us say in advance that we shall mainly study {\em
generic} spaces, in the following sense:

\begin{Definition}
\label{totally-generic} The above arbitrary real vector subspace $H_p
\subset T_p \C^N$ is said to be:

\begin{itemize}

\smallskip\item[$\bullet$]
{\sl totally real} if $H_p^c = H_p \cap J H_p = \{ 0\}$, that is to
say, if $n = 0$;

\smallskip\item[$\bullet$]
{\sl generic} if $H_p^{ i_c} = H_p + J H_p = T_p \C^N$, that is to
say, if $n + d = N$;

\smallskip\item[$\bullet$]
{\sl maximally real} if it is both totally real and generic, that is
to say, if $n = 0$ and if $d = N$.

\end{itemize}\smallskip

\end{Definition}

We notice passim that some constraints on (co)dimensions exist.
Indeed, when $H_p$ is totally real, one has $n = 0$, whence $0 + d
\leqslant N$ and hence ${\rm codim}_\R H_p = 2N - 0 - d \geqslant N$.
But if $H_p$ is in addition maximally real, one has $0 + d = N$,
whence ${\rm codim}_\R H_p = N$ exactly. Finally, when $H_p$ is
generic, its real codimension:
\[
2N-2n-d = 2(n+d)-2n-d = d
\]
is simply equal to the dimension $d$ of its purely real part $H_p
\big/ H_p^c \simeq \R^d$, and one should remember this fact, which
can also be seen by reminding that the complex structure induces a
general isomorphism $ J\colon H_p\big/H_p^c \overset{ \simeq}{
\longrightarrow} H_p^{i_c} \big/H_p$, in which $H_p^{ i_c} \big/ H_p$
becomes $T_p\C^N \big/ H_p$ when $H_p$ is generic. In addition and
for later use, we specify explicitly the form of the defining
Cartesian equations of a generic affine subspace.

\begin{Corollary}
Let $H_p \subset T_p \C^{ n+d}$ as above be an arbitrary real affine
space which is generic and $d$-codimensional in $\C^N = \C^{ n+d}$.
Then there exist $n+d$ affine {\em complex} coordinates:
\[
\big(z_1,\dots,z_n,\,w_1,\dots,w_d\big) =
\big(x_1+iy_1,\dots,x_n+iy_n,\,u_1+iv_1,\dots,u_d+iv_d\big)
\]
centered at $p$ such that $H_p$ is represented by the following $d$
real equations:
\[
0 = {\rm Im}\,w_1 =\cdots= {\rm Im}\,w_d. \qed
\]

\end{Corollary}

In such a concrete coordinate representation in which:
\[
\!\!\!\!\!\!\!\!\!\!\!\!\!\!\!\!\!\!\!\!
\aligned H_p & = \R{\textstyle{\frac{\partial}{\partial
x_1}}}\big\vert_p \oplus \R{\textstyle{\frac{\partial}{\partial
y_1}}}\big\vert_p \oplus\cdots\oplus
\R{\textstyle{\frac{\partial}{\partial x_n}}}\big\vert_p \oplus
\R{\textstyle{\frac{\partial}{\partial y_n}}}\big\vert_p \oplus
\R{\textstyle{\frac{\partial}{\partial u_1}}}\big\vert_p
\oplus\cdots\oplus \R{\textstyle{\frac{\partial}{\partial
u_d}}}\big\vert_p,
\\
H_p^c & = \R{\textstyle{\frac{\partial}{\partial x_1}}}\big\vert_p
\oplus \R{\textstyle{\frac{\partial}{\partial y_1}}}\big\vert_p
\oplus\cdots\oplus \R{\textstyle{\frac{\partial}{\partial
x_n}}}\big\vert_p \oplus \R{\textstyle{\frac{\partial}{\partial
y_n}}}\big\vert_p,
\\
H_p\big/H_p^c & = \ \ \ \ \ \ \ \ \ \ \ \ \ \ \ \ \ \ \ \ \ \ \ \ \ \
\ \ \ \ \ \ \ \ \ \ \ \ \ \ \ \ \ \ \ \ \ \ \ \ \ \ \ \ \ \ \ \ \ \ \
\ \ \ \ \ \ \ \ \ \ \ \ \R{\textstyle{\frac{\partial}{\partial
u_1}}}\big\vert_p \oplus\cdots\oplus
\R{\textstyle{\frac{\partial}{\partial u_d}}}\big\vert_p,
\\
T_p\C^{n+d} \big/ H_p & = \ \ \ \ \ \ \ \ \ \ \ \ \ \ \ \ \ \ \ \ \ \
\ \ \ \ \ \ \ \ \ \ \ \ \ \ \ \ \ \ \ \ \ \ \ \ \ \ \ \ \ \ \ \ \ \ \
\ \ \ \ \ \ \ \ \ \ \ \ \ \ \ \
\R{\textstyle{\frac{\partial}{\partial v_1}}}\big\vert_p
\oplus\cdots\oplus \R{\textstyle{\frac{\partial}{\partial
v_d}}}\big\vert_p,
\endaligned
\]
one sees at once how $J$ induces an isomorphism $H_p \big/ H_p^c
\longrightarrow T_p \C^{ n+d} \big/ H_p$.

On the other hand, when $n + d \leqslant N - 1$ in full generality,
so that $H_p$ is {\em not} generic, the complex-codimensional part
$T_p \C^N \big/ H_p^{ i_c} \simeq \C^{ N - n - d}$ is nontrivial, but
one easily convinces oneself that $H_p$ becomes truly generic when it
is viewed {\em inside} its intrinsic complexification $H_p^{ i_c}
\simeq \C^{ n+d}$. Thus, from the point of view of understanding the
position of a real linear space within a complex linear space, one
may just drop the $\C^N / \C^{ N - n - d}$ and view directly $H_p$
sitting as a generic subspace of its intrinsic complexification $H_p^{
i_c} \simeq \C^{ n+d}$.

\subsection{Complexifications}
By tensoring with $\C$ the real tangent bundle $T \C^N$, we get the
complex vector bundle:
\[
\C \otimes_\R T\C^N = T\C^N \oplus i\,T\C^N,
\]
whose fiber $\simeq \C^N$ at an arbitrary point $p$ of $\C^N$
consists of all possible linear combinations:
\[
{\sf a}_1\,\frac{\partial}{\partial x_1}\Big\vert_p + {\sf
b}_1\,\frac{\partial}{\partial y_1}\Big\vert_p +\cdots+ {\sf
a}_N\,\frac{\partial}{\partial x_N}\Big\vert_p + {\sf
b}_N\,\frac{\partial}{\partial y_N}\Big\vert_p
\]
with free {\em complex} coefficients ${\sf a}_k$, ${\sf b}_k$.

Introduce also the following basic {\sl holomorphic} and {\sl
antiholomorphic} vector fields, for $k = 1, \dots, N$:
\[
\frac{\partial}{\partial z_k} := \frac{1}{2} \bigg(
\frac{\partial}{\partial x_k} - i\,\frac{\partial}{\partial y_k}
\bigg) \ \ \ \ \ \ \ \text{\rm and} \ \ \ \ \ \ \
\frac{\partial}{\partial\overline{z}_k} := \frac{1}{2} \bigg(
\frac{\partial}{\partial x_k} + i\,\frac{\partial}{\partial y_k}
\bigg),
\]
with of course inversely:
\[
\frac{\partial}{\partial x_k} = \frac{\partial}{\partial z_k} +
\frac{\partial}{\partial\overline{z}_k} \ \ \ \ \ \ \ \text{\rm and}
\ \ \ \ \ \ \ \frac{\partial}{\partial y_k} =
i\,\frac{\partial}{\partial z_k} -
i\,\frac{\partial}{\partial\overline{z}_k}.
\]
Then $\C \otimes_\R T\C^N$ happens to decompose as the direct sum of
two specific {\sl holomorphic} and {\sl antiholomorphic} bundles:
\[
\C \otimes_\R T\C^N =: T^{1,0}\C^N \oplus T^{0,1}\C^N,
\]
whose fibers at an arbitrary point $p \in \C^N$ are defined by:
\[
\aligned
T_p^{1,0}\C^N
&
:=
{\rm Span}_\C \bigg( \frac{\partial}{\partial
z_1}\Big\vert_p, \,\dots,\, \frac{\partial}{\partial z_N}\Big\vert_p
\bigg),
\\
T_p^{0,1}\C^N
&
:=
{\rm Span}_\C \bigg(
\frac{\partial}{\partial\overline{z}_1}\Big\vert_p, \,\dots,\,
\frac{\partial}{\partial\overline{z}_N}\Big\vert_p \bigg).
\endaligned
\]
Since $\frac{\partial}{ \partial z_k} = \frac{ \partial}{ \partial
x_k} - i J \big( \frac{ \partial}{ \partial x_k} \big)$ for $k = 1,
\dots, N$, one observes that one may also write:
\[
\aligned
T_p^{1,0}\C^N
&
=
{\rm Span}_\C \big\{ {\sf v}_p-iJ({\sf v}_p)\colon
{\sf v}_p\in T_p\C^N \big\},
\\
T_p^{0,1}\C^N
&
=
{\rm Span}_\C \big\{ {\sf v}_p+iJ({\sf v}_p)\colon
{\sf v}_p\in T_p\C^N \big\}.
\endaligned
\]

Now, let $H_p \subset T_p \C^N$ be an arbitrary vector space as
before. In accordance with what precedes, the tensored
complexification $\C \otimes_\R H_p^c$ of its maximal $J$-invariant
subspace $H^c = H_p \cap J ( H_p)$ decomposes as a direct sum:
\[
\C\otimes_\R H_p^c = H_p^{1,0}\oplus H_p^{0,1},
\]
where, quite similarly:
\begin{equation}
\label{definition-0-1}
\aligned
H_p^{1,0}
&
:=
\big\{ {\sf v}_p-iJ({\sf v}_p)
\colon {\sf v}_p\in H_p^c \big\},
\\
H_p^{0,1}
&
:=
\big\{ {\sf v}_p+iJ({\sf v}_p) \colon {\sf v}_p\in
H_p^c \big\}.
\endaligned
\end{equation}
On the other hand, in coordinates $(z_1, \dots, z_n, w_1, \dots, w_d,
t_1, \dots, t_{N - n - d})$ as in Lemma~\ref{H-z-w-t} above, we have
concretely:
\[
\aligned
H_p^{1,0}
&
=
{\rm Span}_\C \bigg( \frac{\partial}{\partial
z_1}\Big\vert_p, \,\dots,\, \frac{\partial}{\partial z_n}\Big\vert_p
\bigg),
\\
H_p^{0,1}
&
=
{\rm Span}_\C \bigg( \frac{\partial}{\partial\overline{z}_1}\Big\vert_p,
\,\dots,\, \frac{\partial}{\partial\overline{z}_n}\Big\vert_p \bigg).
\endaligned
\]

\section{Zariski-Generic CR behavior
\\
of real analytic submanifolds in $\C^N$} \label{Zariski-generic-CR}
\HEAD{\ref{Zariski-generic-CR}.~Zariski-Generic CR behavior
of real analytic submanifolds in $\C^N$}{
Masoud {\sc Sabzevari} (Shahrekord) and Jo\"el {\sc Merker} (LM-Orsay)}

\subsection{Admitted analyticity assumption}
Now, consider an arbitrary connected differentiable {\em submanifold}
$M$ of $\C^N$, not necessarily straight as were the affine spaces
$H_p$ above. In order to understand the interactions between the {\em
real differentiable structure} of $M$ and the {\em complex structure}
$J$ of $\C^N$, one should study the way how the real tangent planes
$T_p M$ behave with respect to $J$ when $p$ varies in $M$.

However, although $\dim_\R T_p M$ is constant\,\,---\,\,by definition
of a real manifold\,\,---, it is not at all true that in general, the
complex-tangent planes:
\[
T_p^cM := T_pM\cap J(T_pM)
\]
have constant dimensions as $p$ varies in $M$.

As in the ancient works of Sophus Lie and \'Elie Cartan and because
we shall mainly study Lie groups in CR geometry, we shall assume that
$M$ and all the subsequently appearing geometric objects are {\em
real analytic}, and there is a strong reason for this choice: this
will insure, among other things, that for every point $p \in M
\backslash \Sigma$ not lying in a certain closed {\em thin} subset
$\Sigma \subsetneqq M$, the dimension of $T_p M \cap J ( T_pM)$ will
be constant.

Precisely and concretely, the real analyticity
assumption\,\,---\,\,to be held throughout this memoir\,\,---\,\,will
thus be the following. Let $c$ be an integer with $0 \leqslant c
\leqslant 2N$ and use the coordinates $(z_1, \dots, z_N) \equiv (x_1,
y_1, \dots, x_N, y_N)$ on $\C^N$.

\begin{Definition}
A {\em real analytic submanifold} of $\C^N$ of codimension $c$ is a
closed subset $M \subset \C^N$ having the property that for every
point $p \in M$, there exists an open (small) neighborhood $\U_p$ of
$p$ in $\C^N$ and there exist $c$ {\em real analytic} functions
$\rho_1 ( x, y), \dots, \rho_c ( x, y)$ defined
and converging in $\U_p$
having independent real
differentials:
\[
c
=
{\rm rk}_\R
\left(
\begin{array}{cccccc}
\frac{\partial\rho_1}{\partial x_1} & \frac{\partial\rho_1}{\partial y_1}
& \cdots &
\frac{\partial\rho_1}{\partial x_N} & \frac{\partial\rho_1}{\partial y_N}
\\
\cdot\cdot & \cdot\cdot
& \cdots &
\cdot\cdot & \cdot\cdot
\\
\frac{\partial\rho_c}{\partial x_1} & \frac{\partial\rho_c}{\partial y_1}
& \cdots &
\frac{\partial\rho_c}{\partial x_N} & \frac{\partial\rho_c}{\partial y_N}
\end{array}
\right)
(x,y)
\]
at every point $(x, y) \in \U_p$, such that $M \cap \U_p$ consists
of exactly the points $(x, y)$ which satisfy the $c$ Cartesian
equations:
\[
0
=
\rho_1(x,y)
=\cdots=
\rho_c(x,y).
\]

\end{Definition}

As is well known, the rank assumption is precisely the one which
insures that the zero-set is {\em geometrically smooth}, {\em i.e.} is
a manifold. Of course, the neighborhoods $\U_p$
associated to points $p \in M$
may be assumed to be
plain small balls in which all the Taylor series:
\[
\rho_j(x,y)
=
\sum_{\alpha\in\N^N,\,\beta\in\N^N}\,
\rho_{j,\alpha,\beta}\,
x^\alpha\,y^\beta
\ \ \ \ \ \ \ \ \ \ \ \ \
{\scriptstyle{(j\,=\,1\,\cdots\,N)}}
\]
of the functions $\rho_j$ converge normally. Furthermore, as $M$ is
possibly global in $\C^N$, one must be able to compare two systems of
$c$ defining equations inside overlapping balls.

\begin{Lemma}
\label{overlapping-balls}
Whenever a point $p$ belongs to two such neighborhoods $\U_p'$ with
local defining functions $\rho_1' ( x, y), \dots, \rho_c' ( x, y)$ and
$\U_p''$ with local defining functions $\rho_1'' ( x, y), \dots,
\rho_c'' ( x, y)$, there exists a nonempty open subneighborhood $\V_p
\subset \U_p' \cap \U_p''$ and there exists an invertible $c \times c$
matrix $A = \big( a_{ jk} ( x, y) \big)_{ 1 \leqslant j \leqslant c}^{
1 \leqslant k \leqslant c}$ of analytic functions in $\V_p$ such
that:
\[
\rho_j''(x,y)
=
\sum_{k=1}^c\,
a_{j,k}(x,y)\,\rho_k'(x,y),
\]
or equivalently:
\[
\rho_j'(x,y)
=
\sum_{k=1}^c\,
a_{j,k}^{-1}(x,y)\,\rho_k''(x,y),
\]
for every $(x, y) \in \V_p$.
\end{Lemma}

\proof
Leaving out the technical details, the main reason why this is true is the
following. After a straightening, $M \cap \U_p'$ is defined by:
\[
0
=
s_1
=\cdots=
s_c,
\]
in some real coordinates $(s_1, \dots, s_c, s_{
c+1}, \dots, s_{ 2N})$ on $\C^N \cong \R^{ 2N}$ with
$p$ being the origin. Then any other local system of
defining equations $0 = \rho_1 (s) = \cdots = \rho_c (s)$ for $M$ in a
possibly smaller subset $\V_p \subset \U_p'$ must be so that:
\[
0
\equiv
\rho_j(0,\dots,0,s_{c+1},\dots,s_{2N})
\ \ \ \ \ \ \ \ \ \ \ \ \ {\scriptstyle{(j\,=\,1\,\cdots\,d)}},
\]
whence it immediately follows that the power series of
each $\rho_j$ writes under the form:
\[
\rho_j(s)
=
s_1\,a_{j,1}(s)
+\cdots+
s_c\,a_{j,c}(s)
\ \ \ \ \ \ \ \ \ \ \ \ \ {\scriptstyle{(j\,=\,1\,\cdots\,d)}},
\]
for some remainder power series $a_{ j, 1} ( s)$, \dots, $a_{ j, c} (
s)$. But the assumption that the full Jacobian matrix $\big( \frac{
\partial \rho_j}{ \partial s_k} (s) \big)_{ 1\leqslant k \leqslant
2N}^{ 1 \leqslant j \leqslant c}$ has rank $c$ at every $s \in \U_p'$
and the condition $T_p M = \{ 0 = s_1 = \cdots = s_c\}$ imply that in
fact, already the leftmost $c \times c$ minor is nonzero: $\det \big(
\frac{ \partial \rho_j}{ \partial s_k} (0) \big)_{ 1\leqslant k
\leqslant c}^{ 1 \leqslant j \leqslant c} \neq 0$. It therefore
follows by applying the operators $\frac{ \partial}{ \partial s_k}
\big\vert_{ s = 0}$ to the $\rho_j ( s)$ written above
that $\det \big( a_{ j,k} ( 0) \big)_{
1\leqslant j \leqslant c}^{ 1 \leqslant k \leqslant c} \neq 0$,
whence by continuity the same determinant $\det \big( a_{ j,k}
( s) \big)_{ 1\leqslant j \leqslant c}^{ 1 \leqslant k \leqslant c}$
does not vanish for $s$ near the origin, and this completes the
essence of the argument.
\endproof

\subsection{Generic constancy of complex tangential data}
As we said by anticipation, outside
some thin real analytic subset, the geometric
behavior is always fine.

\begin{Proposition}
\label{generic-CR-behavior}
Let $M$ be an arbitrary {\rm connected} real analytic submanifold of
$\C^N$.  Then there exist two integers $n_M \geqslant 0$ and $d_M
\geqslant 0$ with $n_M + d_M \leqslant N$, and there exists a proper
real analytic subset $\Sigma \subsetneqq M$ such that, for every point
$p \in M \backslash \Sigma$ not lying in $\Sigma$, the space:
\[
T_p^cM
:=
T_pM\cap J(T_pM)
\]
has constant real dimension $2n_M$, and such that in addition, the space:
\[
T_p^{i_c}M
:=
T_pM+J(T_pM)
\]
also has constant real dimension, equal to $2n_M + 2d_M$.
\end{Proposition}

Usually, the first space $T_p^c M$ is called the {\sl complex tangent
plane} to $M$ at $p$, while no special name is given, in the
literature, to the intrinsic complexification $T_p^{ i_c} M$ of
$T_pM$. Thus, for all $p$ not in $\Sigma$, the tangent plane $T_p M$
behaves constantly with respect to $J$, like an affine space as studied
previously. Now, before entering the proof, what is a real analytic
subset, and why is it thin?

\subsection{Basic structure of real analytic subsets of $\C^N$}
By passing
to the charts of an atlas, it suffices to consider the case where $M$ is
some real Euclidean space $\R^e$.

\begin{Definition}
A {\sl real analytic subset} of $\R^e$, $e \geqslant 1$, equipped with
coordinates $(s_1, \dots, s_e)$ is a {\em closed} subset $\Sigma \subset \R^e$
having the property that for every point $p \in \Sigma$, there exists
an open (small) neighborhood $\U_p$ of $p$ in $\R^e$ and there exists
a finite number of real analytic functions:
\[
\rho_1( s_1,\dots,s_e),
\dots\dots,
\rho_c(s_1,\dots,s_e)
\]
defined in $\U_p$ such that $\Sigma \cap \U_p$ consists of exactly the
points $(x, y)$ which satisfy the $c$ Cartesian equations:
\[
0
=
\rho_1(s_1,\dots,s_e)
=\cdots=
\rho_c(s_1,\dots,s_e).
\]
\end{Definition}

Notice that no rank condition is required on the differentials of the
$\rho_j$, so that $\Sigma$ is allowed to have
completely arbitrary singularities.
The closedness condition must be emphasized (\cite{ Le-Massey-2007,
Trotman-2007}).

Thus, to define a real analytic subset of a connected real analytic
(abstract)
manifold $M$, one sets up the same definition, intrinsically to $M \cong \R^{
\dim M}$, the quantities $(s_1, \dots, s_{\dim M})$ being any system
of local real analytic coordinates on $M$. However, in the context we
will be dealing with in this memoir, only {\em local} analytic geometric
objects will be studied, sitting inside some fixed small ball of
$\C^N$, with all concerned power series converging normally in such a
small ball.

A real analytic subset $\Sigma \subset M$ is said to be {\sl
proper} if it is not equal to the whole of $M$. As $M$ was assumed
to be connected, it happens that any such proper $\Sigma \subset M$
may be shown to be closed and nowhere dense, so that $M \backslash
\Sigma$ is open with $\overline{ M \backslash \Sigma} = M$. In fact,
more is true, because as is well known, every real analytic
subset may be {\sl stratified}.

\begin{Definition}
(\cite{ Le-Massey-2007, Trotman-2007})
A {\sl stratification} of a real analytic subset $\Sigma \subset M$ of
some real analytic manifold $M$ is a collection of geometrically
smooth real analytic submanifolds $S_\alpha$ of $M$, for $\alpha$
running in some index set $A$, which constitutes a {\em partition} of
$M$:
\[
\bigcup_{\alpha\in A}\,
S_\alpha
=
\Sigma,
\ \ \ \ \ \ \ \ \ \
\text{\rm with}
\ \ \ \ \ \ \ \ \ \
\emptyset
=
S_\alpha
\cap
S_\beta
\ \ \ \ \
\text{\rm for}\ \
\alpha\neq\beta,
\]
which is in addition {\sl locally finite}, namely
satisfies for every compact subset $K \Subset M$:
\[
\mathmotsf{Card}\,
\big\{
\alpha\in A\colon\
S_\alpha\cap K
\neq\emptyset
\big\}
<
\infty,
\]
and lastly, which satisfies the so-called {\sl frontier condition}:
\[
\aligned
\text{\rm if $S_{\alpha} \neq S_\beta$ and $S_\alpha \cap
\overline{S_\beta} \neq \emptyset$, then $S_\alpha \subset
\overline{ S_\beta}$ and $\dim S_\alpha \leqslant \dim S_\beta -
1$}. \endaligned
\]
\end{Definition}

We shall admit without proof the next classical result of
stratifiability, {\em see}~\cite{ Trotman-2007} and the references
therein. Importantly, it implies that the complement $M \backslash
\Sigma$ of a proper real analytic subset is open, for $\Sigma$ then is
a locally finite union of submanifolds of $M$ of dimensions equal to
$1$, $2$, \dots up to at most $n - 1$.

\begin{Theorem}
Any real analytic subset of a real analytic manifold admits
a stratification. \qed
\end{Theorem}

\begin{Definition}
A {\sl Zariski-generic} pointwise property on a real analytic
manifold $M$ is meant a property that holds true at every
point $p \in M \backslash \Sigma$ outside some {\em proper} real
analytic subset $\Sigma \subset M$.
\end{Definition}

The precise terminology {\sl Zariski-generic} is chosen in
order to avoid confusion with the notion of {\sl CR-generic} submanifold
of $\C^N$ ({\em see} below)\,\,---\,\,some authors used the term {\sl generating} in
the
past.

\subsection{Ranks and generic ranks of matrices and mappings}
As a second preliminary before entering the proof of
Proposition~\ref{generic-CR-behavior}, we now study the (generic)
ranks of real analytic matrices and of real analytic mappings, some
two useful model cases in which some exceptional real analytic subset
naturally appear.

Let $e \geqslant 1$, let $(s_1, \dots, s_e)$ be the canonical
coordinates on $\R^e$, let $a \geqslant 1$, let
$b \geqslant 1$ and consider an $a \times b$ matrix:
\[
\Psi(s)
=
\big(\psi_j^k(s)\big)_{1\leqslant j\leqslant a}^{
1\leqslant k\leqslant b}
\]
of functions $\psi_j^k ( s)$ that are real analytic in some small
neighborhood of the origin in $\R^e$. For every integer $r$ such that
$1 \leqslant r \leqslant \min ( a, b)$, one may form the collection of
all $r \times r$ determinants (minors) that are extracted from this
matrix:
\[
\Psi_{j_1,\dots,j_r}^{k_1,\dots,k_r}(s)
:=
\left\vert
\begin{array}{ccc}
\psi_{j_1}^{k_1}(s) & \cdots & \psi_{j_1}^{k_r}(s)
\\
\cdot\cdot & \cdots & \cdot\cdot
\\
\psi_{j_r}^{k_1}(s) & \cdots & \psi_{j_r}^{k_r}(s)
\end{array}
\right\vert
\ \ \ \ \ \ \ \ \ \ \ \ \
\aligned
&
{\scriptstyle{(1\,\leqslant k_1\,<\,\cdots\,<\,k_r\,\leqslant\,b)}},
\\
&
{\scriptstyle{(1\,\leqslant j_1\,<\,\cdots\,<\,j_r\,\leqslant\,a)}}.
\endaligned
\]
Starting from the largest possible size $r := \min ( a, b)$, if all
these determinants are identically zero (as functions of the variables
$s_1, \dots, s_e$), then one passes from the size $r$ to the lower
size $r - 1$, one forms all the $( r - 1) \times ( r - 1)$ minors, and
one tests again whether they all vanish identically or not, and so on.

Then by definition,
the {\sl generic rank} $r^*$ of the matrix-valued function $\Psi
(s)$ is the largest integer $r$ having the property that at least one
$r \times r$ minor is not identically zero, while all higher minors
are identically zero. One has $r^* = 0$ if and only if all entry
functions $\psi_j^k ( s)$ are identically zero (uninteresting
case), and otherwise, one has in full generality $1 \leqslant r^*
\leqslant \min ( a, b)$. Most importantly, if one introduces the
locus:
\[
\Sigma
:=
\Big\{
s\in\R^e:\,
\Psi_{j_1,\dots,j_{r^*}}^{k_1,\dots,k_{r^*}}
(s)
=
0
\colon
\forall\,
j_1,\dots,j_{r^*},\
\forall\,
k_1,\dots,k_{r^*}
\Big\}
\]
of the points $s$ at which all the $r^* \times r^*$ minors vanish,
then this locus clearly is a {\em proper} real analytic subset of
$\R^e$, for at least one function
$\Psi_{ j_1, \dots, j_{r^*}}^{ k_1,\dots, k_{r^*}}$ is not
identically zero.

Furthermore and by construction, at every point $s \in \R^e \backslash
\Sigma$, at least one $r^* \times r^*$ minor is nonzero, and because
all minors of higher size vanish identically by definition of $r^*$,
we deduce the following remarkable property: {\em at every point $s
\in \R^e \backslash \Sigma$ near the origin, the rank of the
matrix-valued function $\Psi ( s)$ is maximal, equal to its generic
rank $r^*$}. A particular case is when $r^* = a = b$, so that the
square
matrix $\Psi ( s)$ is invertible at every point $s \in \R^e \backslash
\Sigma$ near the origin.

Conceptionally speaking, the {\sl generic
rank} of a matrix-valued function is equal to its rank {\em at a generic
point}. What matters for us is that the exceptional real analytic set
of `bad' points is explicitly described as the zero-set of a
collection of minors, which are real analytic functions concretely
known in terms of the initial data $\psi_j^k ( s)$.

These considerations apply directly to the study of the (generic) rank
of any local real analytic map:
\[
(s_1,\dots,s_e)
\longmapsto
\big(
\Phi_1(s_1,\dots,s_e),
\,\dots,\,
\Phi_b(s_1,\dots,s_e)
\big),
\]
for the rank of this map at any point $s$ is equal to the rank
at that point $s$ of its associated Jacobian matrix:
\[
\left(
\begin{array}{ccc}
\frac{\partial\Phi_1}{\partial s_1}
& \cdots &
\frac{\partial\Phi_1}{\partial s_e}
\\
\cdot\cdot
& \cdots &
\cdot\cdot
\\
\frac{\partial\Phi_b}{\partial s_1}
& \cdots &
\frac{\partial\Phi_b}{\partial s_e}
\end{array}
\right)
(s).
\]
If the real analytic objects are globally defined, one verifies that
the exceptional real analytic subsets defined in two coordinate charts
match together and we get the following basic useful observation.

\begin{Lemma}
Given any matrix of real analytic functions defined on a real analytic
manifold, or given any real analytic mapping between two real analytic
manifolds, the set of points where its rank is maximal, equal to the
generic rank, is the complement of a certain {\em proper} real
analytic subset, which may be empty, and in any case, which may be
explicitly described in terms of the matrix, or in terms of the
mapping.\qed
\end{Lemma}

\subsection{Proof of Proposition~\ref{generic-CR-behavior}}
We can now prove the proposition left above.
In the first part of the proof, we work locally, and in the second
part, we show how to glue the local reasonings.

Near an arbitrary point $p \in M$, the real analytic submanifold $M
\subset \C^N$
is represented by $c \leqslant N$ real analytic Cartesian equations of
the form:
\[
0
=
\rho_1(x,y)
=\cdots=
\rho_c(x,y),
\]
in coordinates $(z_1, \dots, z_N) = (x_1 + iy_1, \dots, x_N + iy_N)$
vanishing at $p$, and the local geometric smoothness (`manifoldness') of
$M$ amounts to the assumption that the $c \times N$ Jacobian
matrix:
\[
\left(
\begin{array}{ccccccc}
\frac{\partial\rho_1}{\partial x_1} & \frac{\partial\rho_1}{\partial y_1}
& \cdots &
\frac{\partial\rho_1}{\partial x_N} & \frac{\partial\rho_1}{\partial y_N}
\\
\cdot\cdot & \cdot\cdot
& \cdots &
\cdot\cdot & \cdot\cdot
\\
\frac{\partial\rho_c}{\partial x_1} & \frac{\partial\rho_c}{\partial y_1}
& \cdots &
\frac{\partial\rho_c}{\partial x_N} & \frac{\partial\rho_c}{\partial y_N}
\end{array}
\right)
(x,y)
\]
has rank $c$ everywhere near the origin. Then a vector based at any
point of coordinates $(x, y)$ lying close to the origin:
\[
{\sf v}
\big\vert_{(x,y)}
=
{\sf x}_1\,
{\textstyle{\frac{\partial}{\partial x_1}}}
+
{\sf y}_1\,
{\textstyle{\frac{\partial}{\partial y_1}}}
+\cdots+
{\sf x}_N\,
{\textstyle{\frac{\partial}{\partial x_N}}}
+
{\sf y}_N\,
{\textstyle{\frac{\partial}{\partial y_N}}}
\big\vert_{(x,y)}
\]
belongs to the tangent space $T_{ (x, y)} M$ if and only
the column vector ${}^\tau ( {\sf x}_1, {\sf y}_1,
\dots, {\sf x}_N, {\sf y}_N)$ belongs to the
{\em kernel} of this Jacobian matrix.

On the other hand, the $J$-rotated vector:
\[
J\big({\sf v}\vert_{(x,y)}\big)
=
-{\sf y}_1\,
{\textstyle{\frac{\partial}{\partial x_1}}}
+
{\sf x}_1\,
{\textstyle{\frac{\partial}{\partial y_1}}}
+\cdots
-{\sf y}_N\,
{\textstyle{\frac{\partial}{\partial x_N}}}
+
{\sf x}_N\,
{\textstyle{\frac{\partial}{\partial y_N}}}
\big\vert_{(x,y)}
\]
stays tangent to $M$ at $(x, y)$ if and only if it also
belongs to the same kernel. Equivalently, the initial
vector ${\sf v} \vert_{ x, y}$ belongs to the
kernel of the associated $c \times N$ auxiliary matrix:
\[
\left(
\begin{array}{ccccccc}
\frac{\partial\rho_1}{\partial y_1} & -\frac{\partial\rho_1}{\partial x_1}
& \cdots &
\frac{\partial\rho_1}{\partial y_N} & -\frac{\partial\rho_1}{\partial x_N}
\\
\cdot\cdot & \cdot\cdot
& \cdots &
\cdot\cdot & \cdot\cdot
\\
\frac{\partial\rho_c}{\partial y_1} & -\frac{\partial\rho_c}{\partial x_1}
& \cdots &
\frac{\partial\rho_c}{\partial y_N} & -\frac{\partial\rho_c}{\partial x_N}
\end{array}
\right)
(x,y).
\]
In sum, such a general vector ${\sf v} \vert_{ (x, y)}$ belongs to the
complex tangent plane:
\[
T_{(x,y)}^c M
=
T_{(x,y)}M
\cap
J(T_{(x, y)}M)
\]
if and only if the column vector ${}^\tau ( {\sf x}_1,
{\sf y}_1, \dots, {\sf x}_N, {\sf y}_N)$ belongs to the kernel of the
$2c \times 2N$ matrix (we now use index notation to denote partial
derivatives):
\[
\left(
\begin{array}{cccccc}
\rho_{1,x_1} & \rho_{1,y_1}
& \cdots &
\rho_{1,x_N} & \rho_{1,y_N}
\\
\rho_{1,y_1} & -\rho_{1,x_1}
& \cdots &
\rho_{1,y_N} & -\rho_{1,x_N}
\\
\cdot\cdot & \cdot\cdot
& \cdots &
\cdot\cdot & \cdot\cdot
\\
\rho_{c,x_1} & \rho_{c,y_1}
& \cdots &
\rho_{c,x_N} & \rho_{c,y_N}
\\
\rho_{c,y_1} & -\rho_{c,x_1}
& \cdots &
\rho_{c,y_N} & -\rho_{c,x_N}
\end{array}
\right)
(x,y).
\]
The kernel of this matrix is necessarily even-dimensional, because
$T_{ (x, y)}^c M$ itself\,\,---\,\,a complex
vector space\,\,---\,\,is even-dimensional, as we already know. Thus
in the process of forming minors, we may restrict attention to minors
of size $(2j) \times (2j)$, from $j = \max ( 2c, 2N) = 2N$, downwards
to $j = 2$, and skip identically vanishing minors until some
minor which is
not identically zero is found. Denote then by $2N - 2n_M$ the largest even
integer $2j$ such that there exists a $(2j) \times (2j)$ minor which is
not identically zero, and define $\Sigma$ to be the
real analytic subset which is the zero-set of all
minors of size $(2N - 2n_M) \times
(2N - 2n_M)$, namely:
\[
\Sigma\colon
\ \ \ \ \
0
=
\Delta_1(x,y)
=\cdots=
\Delta_{K}(x,y),
\]
where the number $K$ of such minors
is just equal to the
binomial product $\binom{ 2N}{ 2N - 2n_M} \binom{ 2c}{ 2N - 2n_M}$.

Then by construction, a point $(x, y)$ does not belong to $\Sigma$ if
and only if $T_{ (x, y)}^cM$ has maximal
possible dimension $2n_M$:
\[
\dim_\R T_{(x,y)}^cM
=
2\,n_M,
\ \ \ \ \ \ \ \ \ \
\forall\,(x,y)\not\in\Sigma,
\]
this is the first
property claimed by Proposition~\ref{generic-CR-behavior}.

Next, if we set:
\[
\aligned
d_M
:=
&\,
\dim_\R T_{(x,y)}M
-
\dim_\R\,T_{(x,y)}^cM
\\
=
&\,
2N-c-2\,n_M,
\endaligned
\]
we deduce from an application of the second formula in
equation~\thetag{ \ref{extra-dimensions-coincide}} above that
the complex dimension:
\[
\aligned
\dim_\C T_{(x,y)}^{i_c}M
&
=
{\textstyle{\frac{1}{2}}}\,
\dim_\R T_{(x,y)}^{i_c}M
\\
&
=
{\textstyle{\frac{1}{2}}}\,
\big[
\dim_\R T_pM
+
(\dim_\R T_pM-\dim_\R T_p^cM)
\big]
\\
&
=
{\textstyle{\frac{1}{2}}}\,
\big[
2N-c
+
(2N-c-2n_M)\big)
\big]
\\
&
=
2N-c-n_M
\\
&
=
n_M+d_M
\endaligned
\]
is also constant, as was claimed by the second property of
Proposition~\ref{generic-CR-behavior}.

\smallskip

In order to glue these local reasonings, we observe that if $M$ is
represented by two systems of equations $\rho_j' = 0$,
$j = 1, \dots, c$ and $\rho_j '' =
0$, $j = 1, \dots, c$,
so that, according to Lemma~\ref{overlapping-balls}, one has
$\rho'' = A \, \rho'$ for some invertible matrix $A$, then the Jacobian
matrix of $\rho''$, after restriction to $M$, is equal to $A$ times
the Jacobian matrix of $\rho'$. A similar relation holds between the
two associated auxiliary matrices as above, and it follows from the
theory of matrices that the two zero-sets of minors coincide.
\qed.

\begin{Scholium}
In the proof of
Proposition~\thetag{ \ref{generic-CR-behavior}}, $\Sigma$ is
exactly the set of points $p \in M$ at which $\dim T_p^cM \neq n_M$,
and in fact, this dimension can only increase:
\[
\Sigma
=
\big\{
p\in M\colon\,\,
\dim T_p^cM\geqslant 1+n_M
\big\},
\]
whence the open subset $M \backslash \Sigma$ gathers exactly all
(generic) points of $M$ at which a constant $J$-tangential behavior holds.
\qed
\end{Scholium}

\subsection{Reduction of real analytic local
CR submanifolds to CR-generic submanifolds}
The preceding considerations showed that
it is justified to delineate the following (classical) concepts.

\begin{Definition}
A real analytic submanifold $M \subset \C^N$ is said to be:

\begin{itemize}

\smallskip\item[$\bullet$]
{\sl totally real} if $T_p^cM = T_p M \cap J(T_pM) = \{ 0\}$ is null,
at {\em every} point $p \in M$;

\smallskip\item[$\bullet$]
{\sl holomorphic} if $T_p M = J ( T_pM)$ is fully complex, at {\em every} point
$p \in M$;

\smallskip\item[$\bullet$]
{\sl CR-generic} if $T_pM + J ( T_pM) = T_p \C^N$ generates the whole
ambient tangent space, at {\em every} point $p \in M$;

\smallskip\item[$\bullet$]
{\sl Cauchy-Riemann}\,\,---\,\,{\sl CR} for short\,\,---\,\,if the complex
dimension of $T_p^cM$ is {\em constant}, as $p$ varies in $M$, namely equal
to a certain fixed integer $n_M$.

\end{itemize}\smallskip

\end{Definition}

Obviously, a totally real or holomorphic manifold $M$ is
Cauchy-Riemann. Also a CR-generic $M$ is CR too, because the dimension
formula:
\[
\dim_\R(H+G)
=
\dim_\R H
+
\dim_\R G
-
\dim_\R(H\cap G)
\]
for vector subspaces applied to $H := T_p M$ and to $G := J( T_pM)$
having the same dimension yields if one assumes $T_p M + J ( T_p M) =
T_p \C^N$:
\[
\dim_\R
\big(
T_pM\cap J(T_pM)
\big)
=
2\,\dim_\R T_pM
-
2N,
\]
which is indeed constant independently of the base point (we already saw
this argument in Lemma~\ref{dimension-formula} and after
Definition~\ref{totally-generic}).

The concept of CR submanifold of $\C^N$ embraces that of totally real,
holomorphic and CR-generic submanifolds, but the next proposition ({\em
see}~\cite{ Boggess-1991} for a proof), shows that after possibly
passing to a smaller $\C^N$, every local real analytic CR submanifold
becomes in fact {\em CR-generic}.

\begin{Proposition}
Every connected real analytic CR submanifold of $M \subset \C^N$ of
any CR dimension:
\[
n_M
=
\mathmotsf{rank}\,
\big(TM\cap J(TM)\big)
\]
and of any intrinsic real codimension:
\[
d_M
:=
\dim_\R M
-
2\,n_M,
\]
is locally contained in a certain uniquely defined smallest `germ' of
holomorphic submanifold $M^{ i_c}$ spread along $M$:
\[
M\subset
M^{i_c}
\subset\C^N,
\]
called its {\sl intrinsic complexification} which has complex dimension
equal to:
\[
\dim_\C
M^{i_c}
=
n_M+d_M,
\]
and in addition, $M$ is CR-generic within its intrinsic complexification:
\[
T_pM+J(T_pM)
=
T_pM^{i_c}
\ \ \ \ \ \ \ \ \ \ \ \ \
{\scriptstyle{(p\,\in\,M)}}.
\]

Furthermore, at every point $p \in M$, there exist centered affine holomorphic
coordinates:
\[
\aligned
\big(
&
z_1,\dots,z_{n_M},
w_1,\dots,w_{d_M},t_1,\dots,t_{N-n_M-d_M}
\big)
\in
\C^{n_M}
\times
\C^{d_M}
\times
\C^{N-n_M-d_M},
\\
&
z_1=x_1+i\,y_1,\dots\dots,
z_{n_M}=x_{n_M}+i\,y_{n_M},
\\
&
w_1=u_1+i\,v_1,\dots\dots,
w_{d_M}=u_{d_M}+i\,v_{d_M},
\endaligned
\]
vanishing at that point in which $M$ is locally represented as
the zero-locus of $d_M + 2(N-n_M-d_M)$ real Cartesian equations of
the form:
\[
\left[
\aligned
v_1
&
=
\varphi_1(x,y,u),
\\
\cdots
&
\cdots\cdots\cdots\cdots,
\\
v_{n_M}
&
=
\varphi_{n_M}(x,y,u),
\endaligned\right.
\ \ \ \ \ \ \ \ \ \ \ \ \ \ \ \
\left[
\aligned
t_1
&
=
\Psi_1(z,w),
\\
\cdots\cdots\cdots
&
\cdots\cdots\cdots\cdots\cdots,
\\
t_{N-n_M-d_M}
&
=
\Psi_{N-n_M-d_M}(z,w),
\endaligned\right.
\]
---\,\,implicitly, one takes the real and the imaginary parts
of each one of the $(N-n_M-d_M)$ complex
equations of the second group\,\,---,
in which the $\varphi_j (x, y, u)$ are local real analytic
functions while the $\Psi_k (z, w)$ are local {\em
holomorphic} functions. In such a representation, the holomorphic
submanifold
$M^{ i_c}$ is locally represented as the zero-locus
of the second group of {\em holomorphic} equations.

Lastly, by performing the natural {\em holomorphic} change of coordinates:
\[
t_1'
:=
t_1-\Psi_1(z,w),
\dots\dots\dots,
t_{N-n_M-d_M}'
:=
t_{N-n_M-d_M}
-
\Psi_{N-n_M-d_M}(z,w),
\]
one straightens out the intrinsic complexification $M^{ i_c}$ locally to
become the complex $(n_M+d_M)$-dimensional complex Euclidean space:
\[
\big\{
0
=
t_1'
=\cdots=
t_{N-n_M-d_M}'
\big\},
\]
so that the original CR submanifold $M \subset \C^N$ may be viewed as
sitting in the new complex Euclidean space $\C^{ n_M + d_M}$ of
smaller dimension and as being {\em CR-generic} there.
\qed
\end{Proposition}

These facts then justifies that the equivalence problem under
biholomorphic mappings for arbitrary real analytic submanifolds of
$\C^N$\,\,---\,\,understood mainly at Zariski-generic points similarly
as was the case in Sophus Lie's and \'Elie Cartan's works and
as is usual in its
contemporary prolongations as well\,\,---\,\,comes down to studying
{\em CR-generic submanifolds} in some appropriate $\C^N$.

Of course, when $M$ is of null codimension, $M \equiv \C^N$ locally,
and no equivalence problem exists. Also, in the case where the CR
dimension $n_M$ of $M$ is null, only one local model exists.

\begin{Proposition}
Every CR-generic real analytic submanifold $M \subset \C^N$ which is
totally real, namely of CR dimension $n_M = 0$, is locally
biholomorphically equivalent to a real $N$-dimensional hyperplane,
{\em e.g.} to $\R^N$ sitting in $\R^N + i\, \R^N$.
\qed
\end{Proposition}

Thus, we have now fully justified the (known) fact that only the
consideration of CR-generic submanifolds that have {\em positive}
codimension and {\em positive} CR dimension opens up mathematical
problems.

In this memoir, we will mainly consider the case of CR dimension $1$,
and the equivalence problem will appear to be not at all completely settled
when the codimension is high.

\section{CR-Generic real analytic
submanifolds $M\subset\mathbb{C}^{n+d}$:
\\
real and complex}
\label{real-complex-defining}
\HEAD{\ref{real-complex-defining}.~CR-Generic
real analytic submanifolds $M\subset\mathbb{C}^{n+d}$:
real and complex}{
Masoud {\sc Sabzevari} (Shahrekord) and Jo\"el {\sc Merker} (LM-Orsay)}

\subsection{Real and complex local equations for CR-generic
submanifolds} Consider therefore a local {\em real analytic} submanifold
$M\subset \C^N$ of positive (real) codimension $d \geqslant 1$ which
is {\sl CR-generic} in the sense that its tangent planes:
\[
T_pM+J(T_pM)
=
T_p\C^N
\ \ \ \ \ \ \ \ \ \ \ \ \
{\scriptstyle{(p\,\in\,M)}}
\]
generate the whole ambient tangent plane $T_p\C^N$ over complex
numbers, and which has {\em positive} CR dimension:
\[
n
:=
\dim_\R
\big(
T_pM
\cap
J(T_pM)
\big).
\]
Thus, $N = n+d$, and the letter $N$ will not be used anymore.

In any system of affine holomorphic coordinates:
\[
\aligned
(z,w)
=
\big(
&
z_1,\dots,z_n,\,
w_1,\dots,w_d
\big),
\\
&
z_1=x_1+i\,y_1,\dots\dots\dots,
z_n=x_n+i\,y_n,
\\
&
w_1
=u_1+i\,v_1,\dots\dots\dots,
w_d=u_d+i\,y_d,
\endaligned
\]
centered at one reference point $p_0
\in M$\,\,---\,\,the associated {\em origin}
$0$\,\,---\,\, and for which:
\[
T_{p_0}M
=
\big\{
{\rm Im}\,w_j=0\colon\,
j=1,\dots,d\big\},
\]
the CR-generic submanifold $M \subset \C^{ n+d}$ is locally represented
by $d$ {\em real analytic} equations of the form:
\begin{equation}
\label{real-equations-M}
v_j
=
\varphi_j(x,y,u)
\ \ \ \ \ \ \ \ \ \ \ \ \ {\scriptstyle{(j\,=\,1\,\cdots\,d)}}.
\end{equation}
Viewed geometrically, $M$ is
a graph over its $d$-codimensional plane $T_{p_0}M \subset
T_{p_0} \C^{ n+d}$, with of course the
property that the first order jet of each graphing function
$\varphi_j$ vanishes at the origin:
\[
0 = \varphi_j(0) =
\partial_{x_k}\varphi_j(0)
=
\partial_{y_k}\varphi_j(0)
=
\partial_{u_{j'}}\varphi_j(0)
\ \ \ \ \ \ \ \ \ \ \ \ \ {\scriptstyle{(k\,=\,1\,\cdots\,n\,;\,\,\,
j,\,\,j'\,=\,1\,\cdots\,d)}}.
\]

Let us rewrite these $d$ Cartesian real equations as:
\[
{\textstyle{\frac{w_j-\overline{w}_j}{2\,i}}}
=
\varphi_j \big( {\textstyle{\frac{z+\overline{z}}{2}}},\,
{\textstyle{\frac{z-\overline{z}}{2\,i}}},\,
{\textstyle{\frac{w+\overline{w}}{2}}} \big)
\ \ \ \ \ \ \ \ \ \ \ \ \
{\scriptstyle{(j\,=\,1\,\cdots\,d)}}.
\]
Since the right-hand sides $\varphi_j$ are all an ${\rm O} ( 2)$, we
can then apply the analytic implicit function theorem in order to
solve these equations for the $d$ variables $w_j$, $j = 1, \dots,
d$. Performing this, we obtain a collection of $d$ equations of the
shape:
\[
w_j
=
\Theta_j
\big(z,\,\overline{z},\,\overline{w}\big)
\ \ \ \ \ \ \ \ \ \ \ \ \
{\scriptstyle{(j\,=\,1\,\cdots\,d)}},
\]
whose right-hand side power series converge of course near the origin:
\[
(0,0,0)
\in
\C^n\times\C^n\times\C^d.
\]
Since $d\varphi ( 0) =
0$, one has $\Theta = - \overline{ w} + {\sf order}\,2\,{\sf terms}$.
In fact, the functions $\Theta_j$ are {\em analytic} with respect to
their variables $(z, \overline{ z}, \overline{ w})$, hence they expand
in convergent Taylor series, say under the form:
\[
\Theta_j\big(z,\,\overline{z},\,\overline{w}\big)
=
\sum_{\alpha\in\mathbb{N}^n,\,\beta\in\mathbb{N}^n,\,\gamma\,\in\,
\mathbb{N}^d\atop
\vert\alpha\vert+\vert\beta\vert+\vert\gamma\vert\geqslant 1}\,
\Theta_{j,\alpha,\beta,\gamma}\,
z^\alpha\,\overline{z}^\beta\,\overline{w}^\gamma \in
\C\big\{z,\,\overline{z},\,\overline{w}\big\},
\]
the coefficients $\Theta_{ j, \alpha, \beta, \gamma} \in \C^d$ being
in general {\em non-real complex numbers}, because
of the presence of $i = \sqrt{ -1}$ in the
rewritten Cartesian equations. These sorts of complex equations will appear
to be more convenient to deal with in the sequel,
{\em cf.}~\cite{ Merker-2001, Merker-2005a, Merker-2005b, Merker-2011,
Merker-Hachtroudi}, hence let us explain in which precise,
rigorous sense they are equivalent to the original ones
$v_j = \varphi_j ( x, y, u)$.

Initially, our functions $\varphi_j ( x, y, u)$ were in fact all
{\em real-valued}, namely, in their Taylor series:
\[
\varphi_j(x,y,u)
=
\sum_{
\alpha\in\N^n,\,\,\beta\in\N^n,\,\,\gamma\in\N^d
\atop
\vert\alpha\vert+\vert\beta\vert+\vert\gamma\vert
\geqslant 2}\,
\varphi_{j,\alpha,\beta,\gamma}\,
x^\alpha\,y^\beta\,u^\gamma
\ \ \ \ \ \ \ \ \ \ \ \ \
{\scriptstyle{(j\,=\,1\,\cdots\,d)}},
\]
the appearing Taylor coefficients were all {\em real}:
\[
\varphi_{j,\alpha,\beta,\gamma}
\in
\R.
\]
This feature may be expressed under the form:
\[
\small
\aligned
\overline{\varphi_j(x,y,u)}
\equiv
\sum_{\vert\alpha\vert+\vert\beta\vert+\vert\gamma\vert\geqslant 1}\,
\overline{\varphi_{j,\alpha,\beta,\gamma}}\,\,
\overline{x}^\alpha\,\overline{y}^\beta\,\overline{u}^\gamma
\equiv
\sum_{\vert\alpha\vert+\vert\beta\vert+\vert\gamma\vert\geqslant 1}\,
\varphi_{j,\alpha,\beta,\gamma}\,x^\alpha y^\beta u^\gamma
\equiv
\varphi_j(x,y,u)
\endaligned
\]
of equations that are identically satisfied in the ring $\C\{ x, y, u\}$.

\medskip\noindent{\bf Natural principle for the conjugation of Taylor
series.} {\em With $t = (t_1, \dots, t_c) \in \C^c$ being
complex-valued variables, given an arbitrary complex Taylor series:
\[
\Phi(t)
=
\Phi(t_1,\dots,t_c)
=
\sum_{(\gamma_1,\dots,\gamma_c)\in\N^c}\,
\Phi_{\gamma_1,\dots,\gamma_c}\,
(t_1)^{\gamma_1}
\cdots
(t_c)^{\gamma_c}
=
\sum_{\gamma\in\N^c}\,
\Phi_\gamma\,t^\gamma,
\]
convergent or not, and having complex coefficients:
\[
\Phi_{\gamma_1,\dots,\gamma_c}
\in
\C,
\]
one defines the new Taylor series:
\[
\overline{\Phi}
(t)
:=
\sum_{\gamma\in\N^c}\,
\overline{\Phi_\gamma}\,
t^\gamma
\]
by conjugating {\em only} its complex coefficients, so that the
conjugation operator (overline) can be applied independently and
separately over functions and over variables as shown by the
functional identity:
\[
\overline{\Phi(t_1,\dots,t_c)}
\equiv
\overline{\Phi}(\overline{t_1},\dots,\overline{t_c}).
\]
}\medskip

But now, coming back to the $d$ {\em complex} equations:
\[
w_j
=
\Theta_j\big(z,\overline{z},\overline{w}\big)
\ \ \ \ \ \ \ \ \ \ \ \ \
{\scriptstyle{(j\,=\,1\,\cdots\,d)}},
\]
which we obtained through the implicit function theorem,
how can they represent a {\em real} $d$-codimensional
submanifold of $\C^{ n+d}$? For in principle, they
provide not $d$, but $2\,d$ real equations:
\[
\aligned
0
&
=
{\rm Re}
\big[
w_j-\Theta_j\big(z,\overline{z},\overline{w}\big)
\big]
\ \ \ \ \ \ \ \ \ \ \ \ \
{\scriptstyle{(j\,=\,1\,\cdots\,d)}},
\\
0
&
=
{\rm Im}
\big[
w_j-\Theta_j\big(z,\overline{z},\overline{w}\big)
\big]
\ \ \ \ \ \ \ \ \ \ \ \ \
{\scriptstyle{(j\,=\,1\,\cdots\,d)}},
\endaligned
\]
which is twice what is appropriate. Fortunately\,\,---\,\,{\em
cf.}~\cite{ Merker-2005a}, \S~2.1.13\,\,---, the complex power
series $\Theta_j \big( z, \overline{ z}, \overline{ w} \big)$ are {\em
not} arbitrary, they keep a track of reality.

\begin{Theorem}
The $d$ complex analytic power series:
\[
\Theta_j\big(z,\overline{z},\overline{w}\big)
=
\sum_{\vert\alpha\vert+\vert\beta\vert+\vert\gamma\vert\geqslant 1}\,
\Theta_{j,\alpha,\beta,\gamma}\,
z^\alpha\,
\overline{z}^\beta\,
\overline{w}^\gamma
\]
together with their respective complex conjugate series:
\[
\overline{\Theta}_j
=
\overline{\Theta}_j\big(\overline{z},z,w)
=
\sum_{\vert\alpha\vert+\vert\beta\vert+\vert\gamma\vert\geqslant 1}\,
\overline{\Theta_{j,\alpha,\beta,\gamma}}\,\,
\overline{z}^\alpha\,z^\beta\,w^\gamma \in
\C\big\{\overline{z},\,z,\,w\big\}^d
\]
satisfy the two\,\,---\,\,equivalent by conjugation\,\,---\,\,collections
of $d$
functional equations:
\begin{equation}
\label{reality-Theta} \aligned \overline{w}_j \equiv & \
\overline{\Theta}_j
\big(\overline{z},z,\Theta(z,\overline{z},\overline{w})\big) \ \ \ \
\ \ \ \ \ \ \ \ \ {\scriptstyle{(j\,=\,1\,\cdots\,d)}},
\\
w_j \equiv & \ \Theta_j\big(z,\overline{z},
\overline{\Theta}(\overline{z},z,w)\big) \ \ \ \ \ \ \ \ \ \ \ \ \
{\scriptstyle{(j\,=\,1\,\cdots\,d)}},
\endaligned
\end{equation}
identically in
$\C\big\{ z, \overline{ z}, \overline{ w} \big\}$ and in $\C\big\{ \overline{ z}, z,
w\big\}$, respectively.

Conversely, given any collection
of $d$ local analytic power series:
\[
\Theta_j\big(z,\overline{z},\overline{w}\big)
=
\sum_{\vert\alpha\vert+\vert\beta\vert+\vert\gamma\vert\geqslant 1}\,
\Theta_{j,\alpha,\beta,\gamma}\,
z^\alpha\,
\overline{z}^\beta\,
\overline{w}^\gamma
\ \ \ \ \ \ \ \ \ \ \ \ \
{\scriptstyle{(j\,=\,1\,\cdots\,d)}}
\]
having complex coefficients $\Theta_{ j, \alpha, \beta, \gamma}
\in \C$ and satisfying:
\[
\Theta_j
=
-\,\overline{w}_j
+
\mathmotsf{second order terms},
\]
which, in conjunction with their conjugates $\overline{ \Theta_j}
\big( \overline{ z}, z, w \big)$, satisfy this pair of (equivalent)
functional equations, then the two zero-sets:
\[
\big\{ 0=-\,w+\Theta\big(z,\,\overline{z},\,\overline{w}\big) \big\}
\ \ \ \ \ \ \ \ \ \ \text{\em and} \ \ \ \ \ \ \ \ \ \
\big\{0=-\,\overline{w}+ \overline{\Theta}\big(\overline{z},\,z,\,w\big)
\big\}
\]
coincide and define a local CR-generic $d$-codimensional real analytic
submanifold passing through the origin in $\C^{ n + d}$.\qed
\end{Theorem}

In fact, one may also
show (\cite{ Merker-2005a, Merker-2005b}) that there is an
invertible $d \times d$ matrix $a ( z, w, \overline{ z}, \overline{
w})$ of analytic functions defined near the origin such that one has:
\[
w-\Theta(z,\overline{z},\overline{w}) \equiv
a(z,w,\overline{z},\overline{w})\, \big[
\overline{w}-\overline{\Theta}(\overline{z},z,w) \big],
\]
identically in $\C \{ z, w, \overline{ z}, \overline{ w} \}^d$, whence
the coincidence of the two zero-sets immediately follows,
but we will not need this.

\subsection{Rigid CR-generic submanifolds}
Sometimes, it is advisable to restrict attention
to those CR-generic submanifolds, usually called {\sl rigid},
for which the right-hand side graphing functions are
{\em all independent of the variables $(u_1, \dots, u_d)$}:
\[
v_j
=
\varphi_j(x,y)
\ \ \ \ \ \ \ \ \ \ \ \ \
{\scriptstyle{(j\,=\,1\,\cdots\,d)}}.
\]
In this case, the associated complex defining equations
are most simply computed:
\[
w_j
=
\overline{w}_j
+
2i\,\Phi_j\big(z,\overline{z}\big)
\ \ \ \ \ \ \ \ \ \ \ \ \
{\scriptstyle{(j\,=\,1\,\cdots\,d)}},
\]
with:
\[
\Phi_j\big(z,\overline{z}\big)
:=
\varphi_j(x,y)
\ \ \ \ \ \ \ \ \ \ \ \ \
{\scriptstyle{(j\,=\,1\,\cdots\,d)}}.
\]

\subsection{Existence of normal coordinates}
Up to now, we have made a distinction between writing {\em complex
analytic} functions like the $\Theta_j \big( z, \overline{ z},
\overline{ w} \big)$, and writing {\em real analytic} functions like
the $\varphi_j(x, y, u)$ using the real and imaginary parts $x$ and $y$
of $z$. But since $x = \frac{ z + \overline{ z}}{ 2}$ and $y = \frac{
z - \overline{ z}}{ 2\, i}$, we can also consider that the latter
functions:
\[
\varphi_j\big(
{\textstyle{\frac{z+\overline{z}}{2}}},\,
{\textstyle{\frac{z+\overline{z}}{2}}},\,u
\big)
\]
after expansion in convergent Taylor series and reorganization
of its monomials, depend on $\big(z, \overline{ z}, u\big)$.
(Observe {\em passim}
that this rewriting would fail if $\varphi_j$ were only smooth.)
Hence by a slight abuse of notation, we will sometimes
accept to also write $\varphi_j (z, \overline{ z}, u)$
instead of $\varphi_j ( x, y, u)$.

\begin{Theorem}
\label{normal-coordinates}
Let $M \subset \C^{ n+d}$ be a local real analytic CR-generic
submanifold, let $p_0 \in M$ be one of its points, and assume it to be
represented, in coordinates $(z, w) \in \C^n \times \C^d$
centered at $p_0$,
simultaneously by $d$ real defining equations and by $d$ complex
defining equations of the form:
\[
\aligned
v_j
&
=
\varphi_j\big(z,\overline{z},u\big)
\ \ \ \ \ \ \ \ \ \ \ \ \
{\scriptstyle{(j\,=\,1\,\cdots\,d)}},
\\
w_j
&
=
\Theta_j\big(z,\overline{z},\overline{w}\big)
\ \ \ \ \ \ \ \ \ \ \ \ \
{\scriptstyle{(j\,=\,1\,\cdots\,d)}}.
\endaligned
\]
Then there exists a local biholomorphic change of coordinates
$h \colon\, (z, w) \longmapsto (z', w')$ fixing the origin and having
the specific property of leaving unchanged the $z$-coordinates:
\[
z'
=
z,
\ \ \ \ \ \ \ \ \ \ \
w'
=
g(z,w),
\]
such that the image $M' := h(M)$\,\,---\,\,again a CR-generic
submanifold of $\C^{ n + d}$ passing through the origin\,\,---\,\,has
new real and complex defining equations:
\[
\aligned
v_j'
&
=
\varphi_j'\big(z',\overline{z}',u'\big)
\ \ \ \ \ \ \ \ \ \ \ \ \
{\scriptstyle{(j\,=\,1\,\cdots\,d)}},
\\
w_j'
&
=
\Theta_j'\big(z',\overline{z}',\overline{w}'\big)
\ \ \ \ \ \ \ \ \ \ \ \ \
{\scriptstyle{(j\,=\,1\,\cdots\,d)}}
\endaligned
\]
with right-hand side graphing functions becoming identically
zero whenever one of its arguments $z'$ or $\overline{ z}'$ is null:
\[
\aligned
0
&
\equiv
\varphi_j'\big(0,\overline{z}',u'\big)
\equiv
\varphi_j'\big(z',0,u'\big)
\ \ \ \ \ \ \ \ \ \ \ \ \
{\scriptstyle{(j\,=\,1\,\cdots\,d)}},
\\
0
&
\equiv
\Theta_j'\big(0,\overline{z}',\overline{w}'\big)
\equiv
\Theta_j'\big(z',0,\overline{w}'\big)
\ \ \ \ \ \ \ \ \ \ \ \ \
{\scriptstyle{(j\,=\,1\,\cdots\,d)}}.
\endaligned
\]
\end{Theorem}

\subsection{Fundamental $(1, 0)$ and $(0, 1)$ fields in terms of real
defining equations}
As above, let $M \subset \C^{ n+d}$ be a real analytic CR-generic
submanifold of positive codimension $d \geqslant 1$ and of positive CR
dimension $n \geqslant 1$. Since the real dimension of $T_p^c M = T_p
M \cap J ( T_p M)$ is constantly equal to $n$ at every point of $M$,
the collection of complex-tangential subplanes $\big( T_p^cM \big)_{ p
\in M}$ organizes coherently as a real vector subbundle of $TM$
having rank $2\,n$.

Furthermore, to keep track on $M$ of the basic
holomorphic and antiholomorphic
tangent vectors which exist on $\C^{ n+d}$:
\[
\aligned
&
\frac{\partial}{\partial z_1},
\dots,
\frac{\partial}{\partial z_n},\ \
\frac{\partial}{\partial w_1},
\dots,
\frac{\partial}{\partial w_d},
\\
&
\frac{\partial}{\partial\overline{z}_1},
\dots,
\frac{\partial}{\partial\overline{z}_n},\ \
\frac{\partial}{\partial\overline{w}_1},
\dots,
\frac{\partial}{\partial\overline{w}_d},
\endaligned
\]
it is natural, in view of the
definitions~\thetag{ \ref{definition-0-1}} of
$T_p^{ 1, 0} M$ and of $T_p^{ 0, 1}M$, to define,
at every point $p \in M$,
two {\em complex} vector subspaces:
\[
\aligned
T_p^{1,0}M
&
:=
\mathmotsf{Span}_\C\,
\big\{
{\sf v}_p-iJ({\sf v}_p)\colon\,
{\sf v}_p\in T_p^cM
\big\}
\subset
T_p^{1,0}M,
\\
T_p^{0,1}M
&
:=
\mathmotsf{Span}_\C\,
\big\{
{\sf v}_p+iJ({\sf v}_p)\colon\,
{\sf v}_p\in T_p^cM
\big\}
\subset
T_p^{0,1}M
\endaligned
\]
that are of course conjugate to each other:
\[
T_p^{0,1}M
=
\overline{T_p^{1,0}M}.
\]
One then easily convinces oneself that, as $p$ varies on $M$, these
two collections of spaces organize coherently as two complex vector
bundles on $M$ of rank $n$, and that one may also define them as
being:
\[
\aligned
T^{1,0}M
&
:=
T^{1,0}\C^{n+d}
\cap
\big(\C\otimes_\R TM\big),
\\
T^{0,1}M
&
:=
T^{0,1}\C^{n+d}
\cap
\big(\C\otimes_\R TM\big).
\endaligned
\]
Visibly also, both of them are complex vector {\em subbundles} of:
\[
\C\otimes_\R T\C^{n+d}.
\]

In all what follows, working with a given local real analytic CR-generic
submanifold represented as above, it will be necessary to express
explicitly two (conjugate) bases for these two fundamental vector
bundles.

Thus, with some usual Cartesian equations $v_j = \varphi_j ( x, y, u)$
in which $T_{p_0} M = \big\{ {\rm Im}\, w = 0 \big\}$, we have of course at
the origin:
\[
T_{p_0}^{1,0}M
=
\mathmotsf{Span}_\C
\bigg(
\frac{\partial}{\partial z_1}
\bigg\vert_0,\,
\dots,\,
\frac{\partial}{\partial z_n}
\bigg\vert_0
\bigg).
\]
It follows geometrically that a local basis of $(1, 0)$-vector
fields tangent to $M$, namely a {\sl local frame} for
the antiholomorphic tangent bundle
$T^{ 0, 1} M$, will necessarily be of the form:
\[
\overline{\mathcal{L}_k}
:=
\frac{\partial}{\partial\overline{z}_k}
+
\sum_{l=1}^d\,
{\sf A}_{k,l}(x,y,u)\,
\frac{ \partial}{\partial\overline{w}_l}
\ \ \ \ \ \ \ \ \ \ \ \ \
{\scriptstyle{(k\,=\,1\,\cdots\,n)}},
\]
for certain uniquely defined real analytic functions
${\sf A}_{ k, l} ( x, y, u)$ that may be computed elementarily.
Indeed, the condition that these $\overline{ \mathcal{L}_k}$ be tangent
to $M$, namely tangent to the zero-locus of the
$d$ graphed equations:
\[
{\textstyle{\frac{w_j-\overline{w}_j}{2\,i}}}
=
\varphi_j
\big(x,y,{\textstyle{\frac{w+\overline{w}}{2}}}\big)
\ \ \ \ \ \ \ \ \ \ \ \ \
{\scriptstyle{(j\,=\,1\,\cdots\,d)}},
\]
writes down as:
\[
0
=
\overline{\mathcal{L}_k}
\Big[
-\,{\textstyle{\frac{w_j}{2\,i}}}
+
{\textstyle{\frac{\overline{w}_j}{2\,i}}}
+
\varphi_j
\big(x,y,{\textstyle{\frac{w+\overline{w}}{2}}}\big)
\Big]
\ \ \ \ \ \ \ \ \ \ \ \ \
{\scriptstyle{(j\,=\,1\,\cdots\,d)}}.
\]
Equivalently,
this condition amounts to requiring that, for
every fixed $k = 1, \dots, n$, the following $d$
affine-linear equations are identically
satisfied on $M$ by the $d$ unknowns ${\sf A}_{ k, 1}, \dots,
{\sf A}_{ k, d}$:
\[
0
=
{\textstyle{\frac{1}{2\,i}}}\,{\sf A}_{k,j}
+
\varphi_{j,\overline{z}_k}
+
{\textstyle{\frac{1}{2}}}\,
{\sf A}_{k,1}\,\varphi_{j,u_1}
+\cdots+
{\textstyle{\frac{1}{2}}}\,
{\sf A}_{k,d}\,\varphi_{j,u_d}
\ \ \ \ \ \ \ \ \ \ \ \ \
{\scriptstyle{(j\,=\,1\,\cdots\,d)}}.
\]
Thanks to the fact that by assumption all the $\varphi_{j, \overline{
z}_k}$ and all the $\varphi_{ j, u_{ j'}}$ vanish at the
origin\,\,---\,\,the central point\,\,---, a unique local solution
exists for each $k$, which, in abbreviated matrix notation writes,
shortly:
\[
{\sf A}_k
=
2\,\big(i\,I_{d\times d}-\varphi_u\big)^{-1}
\cdot
\varphi_{\overline{z}_k},
\]
where the $d \times d$ matrix:
\[
\varphi_u
=
\big(
\varphi_{j,u_{j'}}\big
)_{
1\leqslant j\leqslant d}^{
1\leqslant j'\leqslant d}
\]
has row index $j$, where $\varphi_{ \overline{ z}_k}$ is the $d \times
1$ matrix $(\varphi_{ j, \overline{ z}_k})_{ 1 \leqslant j \leqslant
d}$\,\,---\,\,a column vector\,\,---, and where:
\[
{\sf A}_k
=
\big({\sf A}_{k,j}\big)_{
1\leqslant j\leqslant d}.
\]
Of course, this solution is real analytic in a (possibly shrunk)
neighborhood of the origin.

We notice here that these vector fields $\overline{ \mathcal{ L}_k}$
are considered {\em extrinsically}, namely they involve the extra
vector fields $\frac{ \partial}{ \partial v_1}, \dots, \frac{
\partial}{ \partial v_d}$ living in $\C^{ n + d}$ that are {\em not}
tangent to $M$. In order to get {\em intrinsic} sections of $T^{ 0,
1} M$, since $M$ is naturally equipped with the coordinates $(x, y,
u)$, we must naturally drop the $\frac{ \partial}{ \partial v_j}$ and
we obtain:
\[
\overline{\mathcal{L}_k}
\big\vert_M
=
\frac{\partial}{\partial\overline{z_k}}
+
\sum_{l=1}^d\,
{\textstyle{\frac{1}{2}}}\,{\sf A}_{k,l}\,
\frac{\partial}{\partial u_l}
\ \ \ \ \ \ \ \ \ \ \ \ \
{\scriptstyle{(k\,=\,1\,\cdots\,n)}}.
\]
Of course, a local frame for $T^{1, 0}M$ in a neighborhood of the
origin is obtained by plain complex conjugation:
\[
\mathcal{L}_k
\big\vert_M
=
\frac{\partial}{\partial z_k}
+
\sum_{l=1}^d\,
{\textstyle{\frac{1}{2}}}\,\overline{{\sf A}_{k,l}}\,
\frac{\partial}{\partial u_l}
\ \ \ \ \ \ \ \ \ \ \ \ \
{\scriptstyle{(k\,=\,1\,\cdots\,n)}}.
\]

\subsection{Holomorphic and antiholomorphic tangent vector fields}
When computing the coefficients ${\sf A}_{ k,l}$, a somehow unpleasant
matrix inversion was needed at the moment, and this happens to cause
some differential algebra swelling troubles as soon as one enters a
more-in-depth study of the equivalence problem, {\em cf.} what will
follow. On the other hand, when dealing with the (equivalent) complex
defining equations:
\[
\overline{w}_j
=
\overline{\Theta}_j\big(\overline{z},z,w\big)
\ \ \ \ \ \ \ \ \ \ \ \ \
{\scriptstyle{(j\,=\,1\,\cdots\,d)}},
\]
it is clear that the conditions of tangency:
\[
\aligned
0
&
=
\overline{\mathcal{L}_k}
\Big[
-\,\overline{w}_j
+
\overline{\Theta}_j\big(\overline{z},z,w\big)
\Big]
\\
&
=
-\,{\sf A}_{k,j}
+
\frac{\partial\overline{\Theta}_j}{\partial\overline{z}_k}
\big(\overline{z},z,w\big)
\ \ \ \ \ \ \ \ \ \ \ \ \
{\scriptstyle{(j\,=\,1\,\cdots\,d)}}
\endaligned
\]
solves straightforwardly, and we deduce that in such a representation,
the bundle $T^{ 1,0}M$ and its
conjugate $T^{ 0, 1} M$ are generated,
respectively, by the two collections of mutually independent $(1, 0)$
and $(0, 1)$ vector fields:
\[
\aligned \mathcal{L}_k & = \frac{\partial}{\partial z_k} +
\sum_{j=1}^d\, \frac{\partial\Theta_j}{\partial
z_k}\big(z,\overline{z},\overline{w}\big)\, \frac{\partial}{\partial w_j} \ \
\ \ \ \ \ \ \ \ \ \ \ {\scriptstyle{(k\,=\,1\,\cdots\,n)}},
\\
\overline{\mathcal{L}}_k & = \frac{\partial}{\partial\overline{z}_k}
+ \sum_{j=1}^d\,
\frac{\partial\overline{\Theta}_j}{\partial\overline{z}_k}
\big(\overline{z},z,w\big)\,
\frac{\partial}{\partial\overline{w}_j} \ \ \ \
\ \ \ \ \ \ \ \ \ {\scriptstyle{(k\,=\,1\,\cdots\,n)}}.
\endaligned
\]

Then with such a pair of frames, it becomes immediately easier to
compute somewhat explicitly some iterated brackets like for instance:
\[
\big[\mathcal{ L}_{k_1},\,
\overline{\mathcal{L}}_{ k_2}\big],
\ \ \ \ \ \ \ \
\big[\mathcal{L}_{k_1},\,
[\mathcal{L}_{k_2},\,\overline{\mathcal{L}}_{k_3}]\big],
\]
while this task is much more difficult when using the real
representation with the coefficients:
\[
{\sf A}_k
=
2\,\big(i\,I_{d\times d}-\varphi_u\big)^{-1}
\cdot
\varphi_{\overline{z}_k}.
\]
But before pushing further the general theory, it
is now great time to exhibit some paradigmatic examples.

\section{Heisenberg sphere in $\C^2$
\\
and Beloshapka's higher dimensional models}
\label{Heisengerg-Beloshapka}
\HEAD{\ref{Heisengerg-Beloshapka}.~Heisenberg
sphere in $\C^2$ and Beloshapka's higher dimensional
models}{
Masoud {\sc Sabzevari} (Shahrekord) and Jo\"el {\sc Merker} (LM-Orsay)}

\subsection{Heisenberg sphere in $\C^2$ and its deformations}
With $n = 1$, $d = 1$ and $(z, w) \in \C^2$, it is known that the complex
equation:
\[
w
=
\Theta\big(z,\overline{z},\overline{w}\big)
\]
of any real analytic hypersurface $M^3 \subset \C^2$ which
satisfies:
\[
T^{1,0}M
+
T^{0,1}M
+
\big[T^{1,0}M,\,T^{0,1}M\big]
=
\C\otimes_\R TM
\]
---\,\,such are usually called {\sl Levi nondegenerate}\,\,---, may be
brought to the form:
\[
w-\overline{w}
=
2i\,z\overline{z}
+
\mathmotsf{terms of order $\geqslant 3$}.
\]
Moreover, for {\em any} remainder of order $\geqslant 3$, the obtained
hypersurface is Levi nondegenerate.

A concrete proof consists in examining the first-order terms
in the Taylor series of the graphing function:
\[
{\textstyle{\frac{w-\overline{w}}{2\,i}}}
=
\varphi(x,y,u)
=
\alpha\,z^2
+
c\,z\overline{z}
+
\overline{\alpha}\,\overline{z}^2
+
{\rm O}_3\big(z,\overline{z}\big)
+
u\,
{\rm O}_1\big(x,y,u),
\]
with $\alpha \in \C$ and $c \in \R$. If the coordinates are already
normal in the sense of Theorem~\ref{normal-coordinates}, one has
$\alpha = 0$, otherwise, a plain replacement of $w$ by $w' := w -
2i\,\alpha\,z^2$ makes $\alpha = 0$. The constant $c$ happens
to be unremovable
by means of local biholomorphic changes (invariance
of the Levi form) of variables, and one
makes $c = \pm 1$ by substituting $z' := c^{- 1/2}\, z$, and
lastly $c = 1$ by replacing $w' := \pm w$.

\smallskip

Next, the related '{\sl model}' is the one for which the
remainder is identically, namely the so-called {\sl Heisenberg sphere}
$\mathbb{ H}^3 \subset \C^2$ having quadratic equation:
\[
\label{Heisenberg-sphere}
w-\overline{w}
=
2i\,z\overline{z}.
\]
On the other hand, it is known that the unit real $3$-sphere
$S^3$ in $\C^2$ having equation:
\[
1
=
z'\overline{z}'
+
w'\overline{w}'
\]
plays the remarkable r\^ole, in CR geometry, of being
the universal model in $\C^2$. But in fact, $S^3 \setminus \{ p_\infty \}$
with $p_\infty := (0,-1)$, is biholomorphic, through the so-called
{\sl Cayley transform}:
\[
(z,w)\longmapsto
\big(
{\textstyle{\frac{-\,4z}{4i+w}}},\,\,
{\textstyle{\frac{4-iw}{4i+w}}}
\big)
=:
(z',w')
\]
to the above {\sl Heisenberg sphere}.

\subsection{Beloshapka's cubic fourfold in $\C^3$}
Assuming that the CR dimension $n = 1$ is smallest possible, the next
example of a `{\sl universal}' model was introduced in codimension $d
= 2$ by Beloshapka in~\cite{ Beloshapka-1998}. In coordinates $(z,
w_1, w_2)$, it is the cubic:
\[
\label{beloshapka-cubic}
\left[
\aligned
w_1-\overline{w}_1
&
=
2i\,z\overline{z},
\\
w_2-\overline{w}_2
&
=
2i\,z\overline{z}(z+\overline{z}).
\endaligned\right.
\]
Here is
the way it may be introduced. Consider the two
graphed equations of a generic $M^4 \subset \C^3$:
\[
\left[
\aligned
v_1
&
=
\varphi_1(x,y,u_1,u_2),
\\
v_2
&
=
\varphi_2(x,y,u_1,u_2),
\endaligned\right.
\]
assume that the coordinates are normal, and find the `simplest possible'
model. The first equation is brought to the Heisenberg
form. Next, one looks at the first terms in the Taylor
series of the second equation:
\[
\left[
\aligned
v_1
&
=
z\overline{z},
\\
v_2
&
=
a\,z\overline{z}
+
\beta\,z^2\overline{z}
+
\overline{\beta}\,\overline{z}^2z
+
\underbrace{
z\overline{z}\,
\Big[
{\rm O}_2\big(z,\overline{z}\big)
+
u_1\cdot\mathmotsf{remainder}
+
u_2\cdot\mathmotsf{remainder}
\Big]}_{\smallmathmotsf{all monomials are of weighted
order $\geqslant 4$}},
\endaligned\right.
\]
with $a \in \R$ and $\beta \in \C$; notice that, since the coordinates
are normal, namely since $0 \equiv \varphi_2 \big( 0, \overline{ z},
u_1, u_2) \equiv \varphi_2 \big( z, 0, u_1, u_2\big)$, and since $z
\in \C$ is a {\em single} complex variable, all monomials in the
Taylor series must be divisible by $z\overline{ z}$. If one then
assigns natural weights to the variables:
\[
\mathmotsf{weight}(z)
:=
1,
\ \ \ \ \ \ \ \ \ \ \ \ \ \
\mathmotsf{weight}(w_1)
:=
2,
\ \ \ \ \ \ \ \ \ \ \ \ \ \
\mathmotsf{weight}(w_2)
:=
3,
\]
all the remainder terms
are of weighted order $\geqslant 4$, hence they may be dropped if one
just seeks a simple model with at most {\em cubic} terms:
\[
\left[
\aligned
v_1
&
=
z\overline{z},
\\
v_2
&
=
a\,z\overline{z}
+
\beta\,z^2\overline{z}
+
\overline{\beta}\,\overline{z}^2z.
\endaligned\right.
\]

Here of course, a plain subtraction $w_2 ' := w_2 - a\, w_1$ makes $a
= 0$. Next, it is natural to assume that $\beta \neq 0$ (otherwise,
there is degeneration), and replacing $z$ by $\lambda\,z$ with a
$\lambda \in \C$ satisfying $\beta\, \lambda^2 \overline{ \lambda} =
1$, one arrives at the so-called {\sl Beloshapka cubic}:
\[
\left[
\aligned
v_1
&
=
z\overline{z},
\\
v_2
&
=
z^2\overline{z}
+
\overline{z}^2z.
\endaligned\right.
\]
Its geometry-preserving deformations may be introduced in a
coordinate-invariant manner as follows (we skip the proof, because
full details of a more substantial case will be provided
in a while).

\begin{Proposition}
\label{deformations-Beloshapka}
A real analytic $4$-dimensional local CR-generic submanifold $M^4
\subset \C^3$ of codimension $2$ whose complex tangent bundle
satisfies the two equivalent conditions:
\[
\aligned
TM
&
=
T^cM
+
[T^cM,\,T^cM]
+
\big[T^cM,[T^cM,\,T^cM]\big],
\\
\C\otimes_\R TM
&
=
T^{1,0}M+T^{0,1}M
+
[T^{1,0}M,\,T^{0,1}M]
+
\big[T^{1,0}M,[T^{1,0}M,\,T^{0,1}M]\big]
+
\\
&
\ \ \ \ \ \ \ \ \ \ \ \ \ \ \ \ \ \ \ \ \ \ \ \ \ \ \ \ \ \ \ \ \ \ \ \ \
\ \ \ \ \ \ \ \ \ \ \ \ \ \ \ \ \ \ \ \ \ \ \ \ \ \
+
\big[T^{0,1}M,[T^{1,0}M,\,T^{0,1}M]\big]
\endaligned
\]
may always be represented, in suitable holomorphic coordinates
$(z, w_1, w_2)$ by two complex defining equations of the
specific form:
\[
\left[
\aligned
w_1-\overline{w}_1
&
=
2i\,z\overline{z}
+
{\rm O}_{\smallmathmotsf{weighted}}(3),
\\
w_2-\overline{w}_2
&
=
2i\,z\overline{z}\big(z+\overline{z}\big)
+
{\rm O}_{\smallmathmotsf{weighted}}(4).
\qed
\endaligned\right.
\]
\end{Proposition}

\subsection{Coordinatewise introduction of a cubic model $M^5 \subset \C^4$}

Consider now a real analytic, five-dimensional local
real analytic CR submanifold $M^5
\subset \mathbb{ C}^4$ which is CR-generic, hence of CR dimension $1$,
and let $p_0 \in M^5$ be one of its points. There are holomorphic
coordinates:
\[
\big(z,w_1,w_2,w_3\big)
=
\big(
x+iy,\,u_1+iv_1,\,u_2+iv_2,\,u_3+iv_3\big)
\]
vanishing at $p_0$ in which $M^5$ can be represented as a graph
of the form:
\[
\left[
\aligned
{\textstyle{\frac{w_1-\overline{w}_1}{2i}}}
=
v_1
&
=
\varphi_1(x,y,u_1,u_2,u_3)
=
\psi_1(z,\overline{z})
+
{\rm O}(u_1)
+
{\rm O}(u_2)
+
{\rm O}(u_3),
\\
{\textstyle{\frac{w_2-\overline{w}_2}{2i}}}
=
v_2
&
=
\varphi_2(x,y,u_1,u_2,u_3)
=
\psi_2(z,\overline{z})
+
{\rm O}(u_1)
+
{\rm O}(u_2)
+
{\rm O}(u_3),
\\
{\textstyle{\frac{w_3-\overline{w}_3}{2i}}}
=
v_3
&
=
\varphi_3(x,y,u_1,u_2,u_3)
=
\psi_3(z,\overline{z})
+
{\rm O}(u_1)
+
{\rm O}(u_2)
+
{\rm O}(u_3),
\endaligned\right.
\]
over its tangent plane:
\[
T_{p_0} M^5
=
\big\{v_1=v_2=v_3=0\big\}
\]
by means of three real analytic graphing functions $\varphi_1$,
$\varphi_2$, $\varphi_3$
which vanish, together with all their first order
derivatives at the origin. Here by reality of right-hand sides, the
pure $(z, \overline{ z})$ functions must satisfy:
\[
\overline{\psi}_j(\overline{z},z)
\equiv
\psi_j(z,\overline{z})
\ \ \ \ \ \ \ \ \ \ \ \ \ {\scriptstyle{(j\,=\,1,\,2,\,3)}},
\]
identically in $\mathbb{ C}\{ z, \overline{ z}\}$. Replacing if
necessary $w_j$ by the new holomorphic variable $w_j - 2i\, \psi_j (
z, 0)$, we may assume that $\psi_j(z,\overline{z}) = {\rm
O}(z\overline{z})$, whence there are constants $c_j \in \mathbb{ R}$
and $\alpha_j \in \mathbb{ C}$ such that the terms of order $\leqslant
3$ look like:
\[
\psi_j(z,\overline{z})
=
c_j\,z\overline{z}
+
\alpha_j\,z^2\overline{z}
+
\overline{\alpha}_j\,\overline{z}^2z
+
{\rm O}_4(z,\overline{z})
\ \ \ \ \ \ \ \ \ \ \ \ \ {\scriptstyle{(j\,=\,1,\,2,\,3)}}.
\]
In fact, since the coordinates may freely be assumed
to be normal in the sense of Theorem~\ref{normal-coordinates},
we can even assume that the remainders ${\rm O}_4 \big( z, \overline{ z}
\big) = z\overline{ z} \, {\rm O}_2 \big( z, \overline{ z} \big)$
are divisible by $z \overline{ z}$.

Now, we shall make the following first (among two) nondegeneracy
assumption:

\medskip\noindent
{\sf Hypothesis 1:}
At least one of the above three
real constants $c_1$, $c_2$, $c_3$ is nonzero, say $c_1 \neq 0$.

\medskip
Under this assumption, by replacing $w_1$ by $\frac{ 1}{ c_1} \, w_1$,
we arrange $c_1 = 1$, and then, by replacing $w_2$ by $w_2 - c_2\,
w_1$ and $w_3$ by $w_3 - c_3\, w_1$, we come to $c_2 = c_3 = 0$. Let
us therefore rewrite the three equations as follows, using the same
letters $\alpha_j$ which may have changed in the process:
\[
\left[
\aligned
v_1
&
=
z\overline{z}
+
\alpha_1\,z^2\overline{z}
+
\overline{\alpha}_1\,\overline{z}^2z
+
{\rm O}_4(z,\overline{z})
+
{\rm O}(u_1)
+
{\rm O}(u_2)
+
{\rm O}(u_3),
\\
v_2
&
=
\ \ \ \ \ \ \ \ \
\alpha_2\,z^2\overline{z}
+
\overline{\alpha}_2\,\overline{z}^2z
+
{\rm O}_4(z,\overline{z})
+
{\rm O}(u_1)
+
{\rm O}(u_2)
+
{\rm O}(u_3),
\\
v_3
&
=
\ \ \ \ \ \ \ \ \
\alpha_3\,z^2\overline{z}
+
\overline{\alpha}_3\,\overline{z}^2z
+
{\rm O}_4(z,\overline{z})
+
{\rm O}(u_1)
+
{\rm O}(u_2)
+
{\rm O}(u_3).
\endaligned\right.
\]

\medskip\noindent
{\sf Hypothesis 2:} The two complex numbers
$\alpha_2, \alpha_3 \in \mathbb{ C}$ above are $\mathbb{ R}$-linearly
independent.

\medskip
Hence firstly, $\alpha_2 \neq 0$ is nonzero and replacing $z$ by
$\lambda\, z$ with a $\lambda \in \mathbb{ C}$ satisfying $1 =
\alpha_2\, \lambda^2\, \overline{ \lambda}$, we arrange $\alpha_2 = 1$
(to keep $v_1 = z \overline{ z} + {\rm O} (3)$, it suffices to
simultaneously replace $w_1$ by $\frac{ 1}{ \lambda \overline{
\lambda}} \, w_1$).

Secondly, writing $\alpha_3 = \alpha_3 ' + i \,
\alpha_3''$, we may replace $w_3$ by $w_3 - \alpha_3' \, w_2$ to
arrange that $\alpha_3 = i \alpha_3''$ becomes purely imaginary. Then
$\alpha_3 '' \neq 0$ too by Hypothesis 2, and replacing $w_3$
by $\frac{ 1}{ - \alpha_3''}\, w_3$, we arrive at $\alpha_3 = -i$.

Thirdly and Lastly, writing $\alpha_1 = \alpha_1' + i \alpha_1''$ and replacing
$w_1$ by $w_1 - \alpha_1' \, w_2 + \alpha_1'' \, w_3$, we come to
$\alpha_1 = 0$, whence the three equations of $M^5$ have been reduced to
the following initial general form:
\[
\left[
\aligned
{\textstyle{\frac{w_1-\overline{w}_1}{2i}}}
&
=
z\overline{z}
+
\ \ \ \ \ \ \ \ \ \ \ \ \ \ \ \ \ \
\ \ \ \ \ \ \,
+
{\rm O}_4(z,\overline{z})
+
{\rm O}(u_1)
+
{\rm O}(u_2)
+
{\rm O}(u_3),
\\
{\textstyle{\frac{w_2-\overline{w}_2}{2i}}}
&
=
\ \ \ \ \ \ \ \ \
z^2\overline{z}+\overline{z}^2z
\ \ \ \ \ \ \
+
{\rm O}_4(z,\overline{z})
+
{\rm O}(u_1)
+
{\rm O}(u_2)
+
{\rm O}(u_3),
\\
{\textstyle{\frac{w_3-\overline{w}_3}{2i}}}
&
=
\ \ \ \ \ \ \ \ \
-i\,z^2\overline{z}+i\,\overline{z}^2z
+
{\rm O}_4(z,\overline{z})
+
{\rm O}(u_1)
+
{\rm O}(u_2)
+
{\rm O}(u_3).
\endaligned\right.
\]

If we assume that all remainders vanish, we get the
following model of cubic $5$-dimensional
real algebraic CR-generic submanifold of $\C^4$:\[
\left[
\aligned
{\textstyle{\frac{w_1-\overline{w}_1}{2i}}}
&
=
z\overline{z},
\\
{\textstyle{\frac{w_2-\overline{w}_2}{2i}}}
&
=
z^2\overline{z}+\overline{z}^2z,
\\
{\textstyle{\frac{w_3-\overline{w}_3}{2i}}}
&
=
-i\,z^2\overline{z}+i\,\overline{z}^2z.
\endaligned\right.
\]
It generalizes both the Heisenberg sphere $\mathbb{ H}^3
\subset \C^2$ of CR dimension $1$
having defining equation:
\[
w-\overline{w}
=
2i\,z\overline{z},
\]
and
Beloshapka's four-dimensional cubic $\mathbb{ B}^4 \subset \C^3$ of CR
dimension $1$ having the two defining equations:
\[
\left[
\aligned
w_1-\overline{w}_1
&
=
2i\,z\overline{z}
\\
w_2-\overline{w}_2
&
=
2i\,z\overline{z}(z+\overline{z}).
\endaligned\right.
\]

\subsection{Invariant introduction of the cubic model $M^5 \subset \C^4$}

The two hypotheses made above about the CR-generic submanifold $M^5
\subset \C^4$ can be reformulated in a way which shows well that it is
completely invariant and independent of any choice of coordinates.

Indeed, for a general CR-generic $n$-dimensional and $d$-codimensional $M \subset
\C^{ n+d}$, one may look at all the iterated brackets between the
local sections of the complex tangent bundle:
\[
T^cM
=
TM\cap JTM,
\]
which is here of real rank $2$, since $M$ is of CR dimension $1$.
More precisely, one introduces the subsequent subdistributions:
\[
\aligned
D^1M
:=
T^cM,
\ \ \ \ \
D^2M
&
:=
{\rm Span}_{\mathcal{C}^\omega(M)}
\big(
D^1M
+
\big[T^cM,\,D^1M\big]
\big),
\\
D^3M
&
:=
{\rm Span}_{\mathcal{C}^\omega(M)}
\big(
D^2M
+
\big[T^cM,\,D^2M\big]
\big),
\ \ \ \ \ \
\dots
\endaligned
\]
of $TM$ that are linearly generated, over the algebra $\mathcal{
C}^\omega ( M)$ of real analytic functions on $M$, by all possible
iterated Lie brackets between the local sections of $T^cM$, with of
course $D^1 M \subset D^2 M \subset D^3 M \subset \cdots$.

Both in the ancient Lie-Cartan theory and in the more recent field of
subRiemannian geometry ({\em see} {\em e.g.} the survey
\cite{Bellaiche-1996}), it is usual to assume {\sl strong uniformity},
namely that for each $i = 1, 2, 3, \dots$, the dimensions of the
$D_p^iM$ are all locally constant as $p$ runs in $M$, hence are fully
constant if $M$ is thought of as being localized around one of its
points (and connected too, as will always be assumed implicitly). So,
all $D^i M$ are subbundles of $TM$. Furthermore, it is natural to
assume in addition that $M$ is {\em minimal} ({\em see} \cite{
Tumanov-1988, Merker-2005b, Merker-Porten-2006}), namely that:
\[
D^iM
=
TM
\ \ \ \ \
\text{\rm for all}\ \
i\geqslant i^*\ \
\text{\rm large enough}.
\]
Lastly, as a first step in the study of such differential structures,
it is also natural to assume that the ranks $r_1(M)$, $r_2(M)$,
$r_3(M)$, \dots, of the subbundles $D^1M$, $D^2M$, $D^3M$, \dots,
increase as much as possible. We propose to say simply that a CR-generic
$M$ for which the $r_i ( M)$ increase maximally is {\sl maximally
minimal}. There is a general problem, that we will not touch here,
of describing the structure of all possible maximally minimal
CR manifolds, and the concept of free Lie algebra
(\cite{Gershkovich-1988, Gershkovich-Vershik-1994,
Reutenauer-1993}) is concerned.

In our case of a CR-generic 3-codimensional $M^5 \subset \C^4$, because
$T^c M^5$ has rank:
\[
r_1(M^5)
=
2\,{\rm CRdim}\,M^5
=
2,
\]
and because the Lie bracket
is skew-symmetric,
$D^2M^5$ can at most be of rank $3$. Bracketing then
$D^2M^5$ with some two linearly independent local
sections of $T^cM^5$ can
at most yield two more independent vector fields,
so $r_3 ( M^5)$ can at
most be equal to $5 = \dim M^5$. Thus on the agreement that the $r_i (
M^5)$ are maximal possible, it suffices in fact to jump only up to level
$i = 3$ to reach {\sl minimality}, namely:
\[
D^3M^5
=
TM^5.
\]

We will
therefore study the class of real
analytic 5-dimensional CR-generic 3-codimensional submanifolds
$M^5 \subset \C^4$ for which:
\[
\aligned
r_1(M^5)
&
=
2
=
2\,
{\rm CRdim}\,M^5,
\\
r_2(M^5)
&
=
3,
\\
r_3(M^5)
&
=
5
=
\dim M^5.
\endaligned
\]
In other words and using the terminology
introduced a moment ago, we will study {\em maximally minimal} real
analytic CR-generic 3-codimensional submanifolds $M^5 \subset \C^4$.

For completeness and briefly, let us observe that a real analytic
hypersurface $M^3 \subset \C^2$ is maximally minimal if and only if it
is Levi nondegenerate (at every point). Furthermore, the so-called
{\sl Engel CR manifolds} $M^4 \subset \C^3$ of codimension $2$ that
are deformations of Beloshapka's cubic as stated in
Proposition~\ref{deformations-Beloshapka} are maximally minimal too,
with:
\[
\aligned
r_1(M^4)
&
=
2
=
2\,{\rm CRdim}\,M^4,
\\
r_2(M^4)
&
=
3,
\\
r_3(M^4)
&
=
4
=
\dim M^4.
\endaligned
\]
However, in this latter case, a specific phenomenon occurs, because
the maximal freedom for iterated Lie brackets between two linearly
independent sections $\xi_1$ and $\xi_2$ of $T^cM^4$ may (as already
seen) yield in general {\em five} linearly independent vector fields
({\em see} also \S2.1.3 p.~9 in~\cite{Merker-Porten-2006}):
\[
\xi_1,\ \ \ \ \
\xi_2,\ \ \ \ \
[\xi_1,\,\xi_2],\ \ \ \ \
\big[\xi_1,\,[\xi_1,\,\xi_2]\big],\ \ \ \ \
\big[\xi_2,\,[\xi_1,\,\xi_2]\big],\ \ \ \ \
\]
while the dimension of $M^4 \subset \C^3$ is equal to $4 < 5$, which
imposes a constraint. This is the true reason why, in~\cite{
BES-2007}, there is a {\em distinguished} complex-tangential
direction field $\xi_0$, namely a section of $T^cM^4$ such that $\big[
\xi_0, \, [ \xi_1, \, \xi_2] \big]$ is {\em not} linearly
independent of $\xi_1$, $\xi_2$ and $[ \xi_1, \, \xi_2]$. But in our
case of a maximally minimal $M^5 \subset \C^4$, the 5-dimensionality
of $M$ drops such a dimensional constraint and hence, there is no
distinguished direction field, which gives a situation somewhat
analogous to the paradigmatic case of a Levi nondegenerate $M^3
\subset \C^2$.

Thus, let $M^5 \subset \C^4$ be a local 5-dimensional real analytic
CR-generic submanifold having codimension 3. In suitable holomorphic
coordinates $(z, w_1, w_2, w_3)$ centered at some point of $M^5$, the
three {\sl complex defining equations} of $M^5$ are of the form:
\begin{equation}
\label{three-complex-defining-equations}
\left[
\aligned
w_1-\overline{w}_1
&
=
\Xi_1\big(z,\overline{z},\overline{w}_1,\overline{w}_2,\overline{w}_3\big),
\\
w_2-\overline{w}_2
&
=
\Xi_2\big(z,\overline{z},\overline{w}_1,\overline{w}_2,\overline{w}_3\big),
\\
w_3-\overline{w}_3
&
=
\Xi_3\big(z,\overline{z},\overline{w}_1,\overline{w}_2,\overline{w}_3\big),
\endaligned\right.
\end{equation}
where $\Xi_1$, $\Xi_2$ and $\Xi_3$ are analytic functions defined in a
neighborhood of the origin in $\C^5$ such that the
functions:
\[
\Theta_j\big(z,\overline{z},\overline{w}_1,\overline{w}_2,\overline{w}_3\big)
=
\overline{w}_j
+
\Xi_j\big(z,\overline{z},\overline{w}_1,\overline{w}_2,\overline{w}_3\big)
\ \ \ \ \ \ \ \ \ \ \ \ \
{\scriptstyle{(j\,=\,1,\,2,\,3)}}
\]
are subjected to the
reality condition~\thetag{ \ref{reality-Theta}}.

\subsection{Initial normalization of defining equations}

Now, we want to interpret the condition of maximal minimality in terms
of the three complex defining equations~\thetag{
\ref{three-complex-defining-equations}}. As said above, $M^5$ is local,
real analytic,
connected and it passes through the origin. Clearly, a (single) vector
field generating the bundle $T^{ 1, 0} M^5$ of $(1, 0)$-fields tangent
to $M^5$
is:
\[
\mathcal{L}
:=
\frac{\partial}{\partial z}
+
\Xi_{1,z}\,
\frac{\partial}{\partial w_1}
+
\Xi_{2,z}\,
\frac{\partial}{\partial w_2}
+
\Xi_{3,z}\,
\frac{\partial}{\partial w_3}.
\]
Let $\overline{ \mathcal{ L}}$ denote the complex conjugate to
$\mathcal{ L}$:
\[
\overline{\mathcal{L}}
=
\frac{\partial}{\partial\overline{z}}
+
\overline{\Xi}_{1,\overline{z}}\,
\frac{\partial}{\partial\overline{w}_1}
+
\overline{\Xi}_{2,\overline{z}}\,
\frac{\partial}{\partial\overline{w}_2}
+
\overline{\Xi}_{3,\overline{z}}\,
\frac{\partial}{\partial\overline{w}_3},
\]
which is also tangent to $M^5$, since it annihilates the three equations
conjugate to~\thetag{ \ref{three-complex-defining-equations}} that
are known to be equivalent to them, whence $\overline{ \mathcal{ L}}$
generates $T^{ 0, 1} M^5$.

It is often easier and more natural to work with the {\sl extrinsic
complexification} $M^{ e_c}$ of $M^5$, having the same
equations~\thetag{ \ref{three-complex-defining-equations}}, but with
$\overline{ z}$, $\overline{ w}_1$, $\overline{ w}_2$ and $\overline{
w}_3$ being considered as {\em new independent} complex variables
that we will denote:
\[
\big(
\underline{z},\,\underline{w}_1,\,\underline{w}_2,\,\underline{w}_3
\big)
\]
as in Subsection~\ref{extrinsic-complexification} below.
Then according to this subsection, two equivalent
collections of $d$ Cartesian equations of
this extrinsic complexification\,\,---\,\,which is
now a true holomorphic submanifold
of $\C^4 \times \C^4 = \C^8$ equipped with
coordinates $\big( z, w_1, w_2, w_3, \underline{ z},
\underline{ w}_1, \underline{ w}_2, \underline{ w}_3 \big)$\,\,---\,\,are:
\[
\aligned
w_j-\underline{w}_j
&
=
\Xi_j\big(z,\underline{z},\underline{w}_1,\underline{w}_2,
\underline{w}_3\big)
\ \ \ \ \ \ \ \ \ \ \ \ \
{\scriptstyle{(j\,=\,1,\,2,\,3)}},
\\
\underline{w}_j-w_j
&
=
\overline{\Xi}_j
\big(\underline{z},z,w_1,w_2,w_3\big)
\ \ \ \ \ \ \ \ \ \ \ \ \
{\scriptstyle{(j\,=\,1,\,2,\,3)}}.
\endaligned
\]
Furthermore, the extrinsic complexifications of the $(1, 0)$ and
of the $(0, 1)$ vector fields are:
\[
\aligned
\mathcal{L}
&
:=
\frac{\partial}{\partial z}
+
\sum_{l=1}^3\,
\Xi_{l,z}\big(z,\underline{z},\underline{w}_1,
\underline{w}_2,\underline{w}_3\big)\,
\frac{\partial}{\partial w_l},
\\
\underline{\mathcal{L}}
&
:=
\frac{\partial}{\partial\underline{z}}
+
\sum_{l=1}^3\,
\overline{\Xi}_{l,\underline{z}}
\big(\underline{z},z,w_1,w_2,w_3\big)\,
\frac{\partial}{\partial\underline{w}_l}.
\endaligned
\]

Lastly, since $M^{ e_c}$ is a $5$-dimensional holomorphic submanifold
of $\C^8$ and since it has two equivalent $3$-tuples of Cartesian
equations, two equivalent local frames for its holomorphic tangent
bundle near the origin are:
\[
\Big(
\mathcal{L},\ \
\underline{\mathcal{L}},\ \
{\textstyle{\frac{\partial}{\partial\underline{w}_1}}},\ \
{\textstyle{\frac{\partial}{\partial\underline{w}_2}}},\ \
{\textstyle{\frac{\partial}{\partial\underline{w}_3}}}
\Big)
\ \ \ \ \ \ \ \ \ \ \ \ \ \
\text{\rm and}
\ \ \ \ \ \ \ \ \ \ \ \ \ \
\Big(
\underline{\mathcal{L}},\ \
\mathcal{L},\ \
{\textstyle{\frac{\partial}{\partial w_1}}},\ \
{\textstyle{\frac{\partial}{\partial w_2}}},\ \
{\textstyle{\frac{\partial}{\partial w_3}}}
\Big).
\]

\begin{Lemma}
The local CR-generic $3$-codimensional
real analytic submanifold $M^5 \subset \C^4$ is
maximally minimal at the origin if and only if the five vector fields:
\[
\mathcal{L},\ \ \ \ \
\underline{\mathcal{L}},\ \ \ \ \
\big[\mathcal{L},\,\underline{\mathcal{L}}\big],\ \ \ \
\big[\mathcal{L},\,
\big[\mathcal{L},\,\underline{\mathcal{L}}\big]\big],\ \ \ \ \
\big[\underline{\mathcal{L}},\,
\big[\mathcal{L},\,\underline{\mathcal{L}}\big]\big]
\]
are $\C$-linearly independent at the origin,
hence at all points of $M^{ e_c}$ (which, again, is small and
localized around the origin).
\end{Lemma}

\proof
Indeed, to check this claim, one may consider:
\[
L_1
:=
{\rm Re}\,\mathcal{L}
\ \ \ \ \
\text{\rm and}
\ \ \ \ \
L_2
:=
{\rm Im}\,\mathcal{L}
\]
as two sections generating $T^c M^5 = {\rm Re}\, T^{ 1, 0}
M^5$ and compare real and complex linear spans.
\endproof

At first, such an independency condition
requires at least that the first three
vector fields $\mathcal{ L}$,
$\underline{ \mathcal{ L}}$ and $\big[ \mathcal{ L}, \, \underline{
\mathcal{ L}} \big]$ are $\C$-linearly independent at the origin. But
a direct computation yields the expression of the bracket:
\begin{equation}
\label{L-barL}
\small
\aligned
\big[
\mathcal{L},\,\underline{\mathcal{L}}
\big]
&
=
\bigg[
\frac{\partial}{\partial z}
+
\sum_{l=1}^3\,
\Xi_{l,z}\,
\frac{\partial}{\partial w_l},
\ \
\frac{\partial}{\partial\underline{z}}
+
\sum_{j=1}^3\,
\underline{\Xi}_{j,\underline{z}}\,
\frac{\partial}{\partial\underline{w}_j}
\bigg]
\\
&
=
\sum_{j=1}^3\,
\bigg(
\underline{\Xi}_{j,\underline{z}z}
+
\sum_{l=1}^3\,\Xi_{l,z}\,
\underline{\Xi}_{j,\underline{z}w_l}
\bigg)
\frac{\partial}{\partial\underline{w}_j}
-
\sum_{j=1}^3\,
\bigg(
\Xi_{j,z\underline{z}}
+
\sum_{l=1}^3\,
\underline{\Xi}_{l,\underline{z}}\,
\Xi_{j,z\underline{w}_l}
\bigg)
\frac{\partial}{\partial w_j},
\endaligned
\end{equation}
and since $0 = \Theta_{ l, z} ( 0) = \Theta_{ l, \underline{ z}} ( 0)$
for $l = 1, 2, 3$, we realize that these three vectors based at the
origin have the following simple values:
\begin{equation}
\label{three-L-L}
\aligned
\mathcal{L}\big\vert_0
&
=
{\textstyle{\frac{\partial}{\partial z}}}
\big\vert_0,
\\
\underline{\mathcal{L}}\big\vert_0
&
=
{\textstyle{\frac{\partial}{\partial\underline{z}}}}
\big\vert_0,
\\
\big[\mathcal{L},\,\underline{\mathcal{L}}\big]
\big\vert_0
&
=
\underline{\Xi}_{1,\underline{z}z}(0)\,
{\textstyle{\frac{\partial}{\partial\underline{w}_1}}}
\big\vert_0
+
\underline{\Xi}_{2,\underline{z}z}(0)\,
{\textstyle{\frac{\partial}{\partial\underline{w}_2}}}
\big\vert_0
+
\underline{\Xi}_{3,\underline{z}z}(0)\,
{\textstyle{\frac{\partial}{\partial\underline{w}_3}}}
\big\vert_0
-
\\
&
\ \ \ \ \
-
\Xi_{1,z\underline{z}}(0)\,
{\textstyle{\frac{\partial}{\partial w_1}}}
\big\vert_0
-
\Xi_{2,z\underline{z}}(0)\,
{\textstyle{\frac{\partial}{\partial w_2}}}
\big\vert_0
-
\Xi_{3,z\underline{z}}(0)\,
{\textstyle{\frac{\partial}{\partial w_3}}}
\big\vert_0.
\endaligned
\end{equation}
Without loss of generality, and in order to simplify a bit our next
computations, we will assume that the coordinates $(z, w_1, w_2,
w_3)$ are normal from the beginning,
namely that:
\begin{equation}
\label{normality-conditions}
0
\equiv
\Xi_j\big(0,\underline{z},\underline{w}_1,\underline{w}_2,\underline{w}_3\big)
\equiv
\Xi_j\big(z,0,\underline{w}_1,\underline{w}_2,\underline{w}_3\big)
\ \ \ \ \ \ \ \ \ \ \ \ \ {\scriptstyle{(j\,=\,1,\,\,2,\,\,3)}}.
\end{equation}
It follows
that each $\Xi_j$ is a multiple of the product $z\underline{ z}$,
namely of the form $\Xi_j = z\underline{ z}\, \Xi_j^\sim$, and in
particular, is an ${\rm O} ( z\overline{ z})$; obviously, the same
 also holds true of each $\underline{ \Xi}_j$. So necessarily,
there can be only one second-order term in each $\Xi_j$ and it is of the
form $\lambda_j \, z\underline{ z}$ for some $\lambda_j \in \C$,
because only the constant in $\Xi_j^\sim$ is concerned. But the
reality condition~\thetag{ \ref{reality-Theta}} implies that each
$\lambda_j$ belongs to $i\, \R$, and hence we can write
$\lambda_j = 2i\, c_j \, z \underline{ z}$ for some $c_j \in
\R$. We therefore get:
\[
\aligned
w_1-\underline{w}_1
&
=
2ic_1\,z\underline{z}
+
{\rm O}(3),
\\
w_2-\underline{w}_2
&
=
2ic_2\,z\underline{z}
+
{\rm O}(3),
\\
w_3-\underline{w}_3
&
=
2ic_3\,z\underline{z}
+
{\rm O}(3).
\endaligned
\]

After these preparations, looking at the three expressions~\thetag{
\ref{three-L-L}}, we deduce that for ${\rm Span}_\C \big( \mathcal{ L}
\big\vert_0, \, \, \underline{ \mathcal{ L}} \big\vert_0,\,\, \big[
\mathcal{ L}, \, \underline{ \mathcal{ L}} \big] \big\vert_0 \big)$ to
be 3-dimensional, it is necessary and sufficient that at least one
$c_j$ be nonzero, let us say: $c_1 \neq 0$, after
permuting the coordinates, if necessary. Replacing then $z$ by $\sqrt{
c_1} \, z$ (for some complex square root of $c_1$), we get $c_1 = 1$.
But then in addition,
replacing $w_2$ by $w_2 - c_2 w_1$ and $w_3$ by $w_3 -
c_3 w_1$, we make $c_2 = c_3 = 0$. In summary, the equations
of $M^{ e_c}$ receive the form:
\[
\aligned
w_1-\underline{w}_1
&
=
2i\,z\underline{z}
+
{\rm O}(3)
=
\Xi_1,
\\
w_2-\underline{w}_2
&
=
\ \ \ \ \ \ \ \ \ \ \ \ \
{\rm O}(3)
=
\Xi_2,
\\
w_3-\underline{w}_3
&
=
\ \ \ \ \ \ \ \ \ \ \ \ \
{\rm O}(3)
=
\Xi_3,
\endaligned
\]
and the coordinates are still normal. Now, we must examine
the ${\rm O} ( 3)$ terms.

The first length-three bracket may be computed completely and the
result appears to be a long four-lines expression:
\begin{equation}
\!\!\!\!\!\!\!\!\!\!\!\!\!\!\!\!\!\!\!\!
\!\!\!\!\!\!\!\!\!\!\!\!\!\!\!\!\!\!\!\!
\label{L-L-barL}
\scriptsize
\aligned
&
\big[\mathcal{L},\,\,
\big[\mathcal{L},\,\underline{\mathcal{L}}\big]
\big]
=
\\
&
=
\bigg[
\frac{\partial}{\partial z}
+
\sum_{k=1}^3\,
\Xi_{k,z}\,
\frac{\partial}{\partial w_k},
\ \ \
\sum_{j=1}^3\,
\bigg(
\underline{\Xi}_{j,\underline{z}z}
+
\sum_{l=1}^3\,\Xi_{l,z}\,
\underline{\Xi}_{j,\underline{z}w_l}
\bigg)
\frac{\partial}{\partial\underline{w}_j}
-
\sum_{j=1}^3\,
\bigg(
\Xi_{j,z\underline{z}}
+
\sum_{l=1}^3\,
\underline{\Xi}_{l,\underline{z}}\,
\Xi_{j,z\underline{w}_l}
\bigg)
\frac{\partial}{\partial w_j}
\bigg]
=
\\
&
=
\sum_{j=1}^3\,
\bigg(
\underline{\Xi}_{j,\underline{z}zz}
+
\sum_{l=1}^3\,
\big(
\Xi_{l,zz}\,\underline{\Xi}_{j,\underline{z}w_l}
+
\Xi_{l,z}\,\underline{\Xi}_{j,\underline{z}w_lz}
\big)
+
\sum_{k=1}^3\,
\Xi_{k,z}\,
\underline{\Xi}_{j,\underline{z}zw_k}
+
\sum_{k=1}^3\,\sum_{l=1}^3\,
\Xi_{k,z}
\Big(
\Xi_{l,zw_k}\,\underline{\Xi}_{j,\underline{z}w_l}
+
\Xi_{l,z}\,\underline{\Xi}_{j,\underline{z}w_lw_k}
\Big)
\bigg)
\frac{\partial}{\partial\underline{w}_j}
-
\\
&
-
\sum_{j=1}^3\,
\bigg(
\Xi_{j,z\underline{z}z}
+
\sum_{l=1}^3\,
\big(
\underline{\Xi}_{l,\underline{z}z}\,\Xi_{j,z\underline{w}_l}
+
\underline{\Xi}_{l,\underline{z}}\,\Xi_{j,z\underline{w}_lz}
\big)
+
\sum_{k=1}^3\,
\Xi_{k,z}\,\Xi_{j,z\underline{z}w_k}
+
\sum_{k=1}^3\,\sum_{l=1}^3\,
\Xi_{k,z}\,
\Big(
\underline{\Xi}_{l,\underline{z}w_k}\,
\Xi_{j,z\underline{w}_l}
+
\underline{\Xi}_{l,\underline{z}}\,
\Xi_{j,z\underline{w}_lw_k}
\Big)
\bigg)
\frac{\partial}{\partial w_j}
-
\\
&
-
\sum_{k=1}^3\,
\bigg(
\sum_{j=1}^3\,
\Big(
\underline{\Xi}_{j,\underline{z}z}
+
\sum_{l=1}^3\,\Xi_{l,z}\,
\underline{\Xi}_{j,\underline{z}w_l}
\Big)\,
\Xi_{k,z\underline{w}_j}
\bigg)
\frac{\partial}{\partial w_k}
+
\\
&
+
\sum_{k=1}^3\,
\bigg(
\sum_{j=1}^3\,
\Big(
\Xi_{j,z\underline{z}}
+
\sum_{l=1}^3\,
\underline{\Xi}_{l,\underline{z}}\,
\Xi_{j,z\underline{w}_l}
\Big)\,
\Xi_{k,zw_j}
\bigg)
\frac{\partial}{\partial w_k}.
\endaligned
\end{equation}
But the normality conditions~\thetag{ \ref{normality-conditions}}
entail that:
\[
0
=
\Xi_{j,z}(0)
=
\Xi_{l,zz}(0)
=
\Xi_{k,z\underline{w}_j}(0)
=
\Xi_{k,zw_j}(0),
\]
and consequently, the value at zero of this length-three
bracket is just equal to:
\[
\big[
\mathcal{L},\,\,
\big[\mathcal{L},\,\underline{\mathcal{L}}\big]
\big]
\big\vert_0
=
\sum_{j=1}^3\,
\underline{\Xi}_{j,\underline{z}zz}(0)\,
\frac{\partial}{\partial\underline{w}_j}
-
\sum_{j=1}^3\,\Xi_{j,z\underline{z}z}(0)\,
\frac{\partial}{\partial w_j}.
\]
In a completely similar way, we also get:
\[
\big[
\underline{\mathcal{L}},\,\,
\big[\mathcal{L},\,\underline{\mathcal{L}}\big]
\big]
\big\vert_0
=
\sum_{j=1}^3\,
\underline{\Xi}_{j,\underline{z}z\underline{z}}(0)\,
\frac{\partial}{\partial\underline{w}_j}
-
\sum_{j=1}^3\,
\Xi_{j,z\underline{z}\underline{z}}(0)\,
\frac{\partial}{\partial w_j}.
\]

Now, looking at those third-order terms in the three equations
of $M^5$ which only involve $z$ and $\underline{ z}$, we
may write:
\[
\aligned
w_1-\underline{w}_1
&
=
2i\,z\underline{z}
+
\alpha\,z^2\underline{z}
+
\widetilde{\alpha}\,\underline{z}^2z
+
{\rm O}(\vert z\vert^4)
+
z\underline{z}\,{\rm O}(\vert w\vert)
\\
w_2-\underline{w}_2
&
=
\ \ \ \ \ \ \ \ \ \ \ \ \
\beta\,z^2\underline{z}
+
\widetilde{\beta}\,\underline{z}^2z
+
{\rm O}(\vert z\vert^4)
+
z\underline{z}\,{\rm O}(\vert w\vert)
\\
w_3-\underline{w}_3
&
=
\ \ \ \ \ \ \ \ \ \ \ \ \
\gamma\,z^2\underline{z}
+
\widetilde{\gamma}\,\underline{z}^2z
+
{\rm O}(\vert z\vert^4)
+
z\underline{z}\,{\rm O}(\vert w\vert),
\endaligned
\]
for some complex constants $\alpha, \widetilde{ \alpha}, \beta,
\widetilde{ \beta}, \gamma, \widetilde{ \gamma} \in \C$. But we claim
that $\widetilde{ \alpha} = -\overline{ \alpha}$, that $\widetilde{
\beta} = -\overline{ \beta}$ and that $\widetilde{ \gamma} =
-\overline{ \gamma}$, in fact. Indeed, the reality condition may be
inspected by conjugating-complexifying the above equations, which yields:
\[
\aligned
\underline{w}_1-w_1
&
=
-2i\,\underline{z}z
+
\overline{\alpha}\,\underline{z}^2z
+
\overline{\widetilde{\alpha}}\,z^2\underline{z}
+
{\rm O}(\vert z\vert^4)
+
z\underline{z}\,{\rm O}(\vert w\vert)
\\
\underline{w}_2-w_2
&
=
\ \ \ \ \ \ \ \ \ \ \ \ \ \ \ \ \
\overline{\beta}\,\underline{z}^2z
+
\overline{\widetilde{\beta}}\,z^2\underline{z}
+
{\rm O}(\vert z\vert^4)
+
z\underline{z}\,{\rm O}(\vert w\vert)
\\
\underline{w}_3-w_3
&
=
\ \ \ \ \ \ \ \ \ \ \ \ \ \ \ \ \
\overline{\gamma}\,\underline{z}^2z
+
\overline{\widetilde{\gamma}}\,z^2\underline{z}
+
{\rm O}(\vert z\vert^4)
+
z\underline{z}\,{\rm O}(\vert w\vert),
\endaligned
\]
which should yield equations that are {\em equivalent} the previous
ones, and since the two remainders ${\rm O} (\vert z\vert^4)$ and
${\rm O} ( \vert w \vert)$ visibly do not interfere at all,
it follows by identifications of the monomials
$z^2 \underline{ z}$ and $\underline{ z}^2 z$ that $\widetilde{ \alpha} =
-\overline{ \alpha}$, that $\widetilde{ \beta} =
-\overline{ \beta}$ and that $\widetilde{ \gamma} =
-\overline{ \gamma}$, as was claimed.
As a result, the equations of $M^5$ are:
\[
\aligned
w_1-\underline{w}_1
&
=
2i\,z\underline{z}
+
\alpha\,z^2\underline{z}
-
\overline{\alpha}\,\underline{z}^2z
+
{\rm O}(\vert z\vert^4)
+
z\underline{z}\,{\rm O}(\vert w\vert),
\\
w_2-\underline{w}_2
&
=
\ \ \ \ \ \ \ \ \ \ \ \ \
\beta\,z^2\underline{z}
-
\overline{\beta}\,\underline{z}^2z
+
{\rm O}(\vert z\vert^4)
+
z\underline{z}\,{\rm O}(\vert w\vert),
\\
w_3-\underline{w}_3
&
=
\ \ \ \ \ \ \ \ \ \ \ \ \
\gamma\,z^2\underline{z}
-
\overline{\gamma}\,\underline{z}^2z
+
{\rm O}(\vert z\vert^4)
+
z\underline{z}\,{\rm O}(\vert w\vert).
\endaligned
\]

On the intrinsic complexification $M^{ e_c}$, we now
choose the five complex coordinates $z$, $\underline{ z}$,
$\underline{ w}_1$, $\underline{ w}_2$ and $\underline{ w}_3$.
But what has been already seen, we know that,
in these five intrinsic coordinates of $M^{ e_c}$\,\,---\,\,which means that
we drop $\frac{ \partial}{ \partial w_1}$,
$\frac{ \partial}{ \partial w_2}$ and $\frac{ \partial}{ \partial
w_3}$\,\,---, we have:
\[
\mathcal{L}\big\vert_0
=
{\textstyle{\frac{\partial}{\partial z}}}\big\vert_0,
\ \ \ \ \ \ \ \ \ \
\underline{\mathcal{L}}\big\vert_0
=
{\textstyle{\frac{\partial}{\partial\underline{z}}}}\big\vert_0,
\ \ \ \ \ \ \ \ \ \
\big[\mathcal{L},\,\underline{\mathcal{L}}\big]\big\vert_0
=
{\textstyle{\frac{\partial}{\partial\underline{w}_1}}}\big\vert_0.
\]
On the other hand, also in the coordinates of $M^{ e_c}$, we know that:
\[
\aligned
\big[\mathcal{L},\,\,
\big[\mathcal{L},\,\underline{\mathcal{L}}\big]
\big]
\big\vert_0
&
=
\underline{\Xi}_{1,z\underline{z}z}(0)\,
{\textstyle{\frac{\partial}{\partial\underline{w}_1}}}\big\vert_0
+
\underline{\Xi}_{2,z\underline{z}z}(0)\,
{\textstyle{\frac{\partial}{\partial\underline{w}_2}}}\big\vert_0
+
\underline{\Xi}_{3,z\underline{z}z}(0)\,
{\textstyle{\frac{\partial}{\partial\underline{w}_3}}}\big\vert_0
\\
&
=
\alpha\,
{\textstyle{\frac{\partial}{\partial\underline{w}_1}}}\big\vert_0
+
\beta\,
{\textstyle{\frac{\partial}{\partial\underline{w}_2}}}\big\vert_0
+
\gamma
{\textstyle{\frac{\partial}{\partial\underline{w}_3}}}\big\vert_0,
\endaligned
\]
and quite similarly, that:
\[
\aligned
\big[\underline{\mathcal{L}},\,\,
\big[\mathcal{L},\,\underline{\mathcal{L}}\big]
\big]
\big\vert_0
&
=
\underline{\Xi}_{1,\underline{z}z\underline{z}}(0)\,
{\textstyle{\frac{\partial}{\partial\underline{w}_1}}}\big\vert_0
+
\underline{\Xi}_{2,\underline{z}z\underline{z}}(0)\,
{\textstyle{\frac{\partial}{\partial\underline{w}_2}}}\big\vert_0
+
\underline{\Xi}_{3,\underline{z}z\underline{z}}(0)\,
{\textstyle{\frac{\partial}{\partial\underline{w}_3}}}\big\vert_0
\\
&
=
-\overline{\alpha}\,
{\textstyle{\frac{\partial}{\partial\underline{w}_1}}}\big\vert_0
-
\overline{\beta}\,
{\textstyle{\frac{\partial}{\partial\underline{w}_2}}}\big\vert_0
-
\overline{\gamma}
{\textstyle{\frac{\partial}{\partial\underline{w}_3}}}\big\vert_0.
\endaligned
\]
Consequently, for these five complex vectors based
at the origin to generate:
\[
T_0M^{e_c}
=
\C
{\textstyle{\frac{\partial}{\partial z}}}
\oplus
\C
{\textstyle{\frac{\partial}{\partial\underline{z}}}}
\oplus
\C
{\textstyle{\frac{\partial}{\partial\underline{w}_1}}}
\oplus
\C
{\textstyle{\frac{\partial}{\partial\underline{w}_2}}}
\oplus
\C
{\textstyle{\frac{\partial}{\partial\underline{w}_3}}},
\]
it is necessary and sufficient that the $2 \times 2$ determinant
(we drop minus signs):
\[
\left\vert\!\!
\begin{array}{cc}
\beta & \gamma
\\
\overline{\beta} & \overline{\gamma}
\end{array}
\!\!\right\vert
\neq
0
\]
be nonzero. But things are not completely finished, for some pleasant
simplifications are still available. At least, this last condition
requires $\beta \neq 0$. Then we can make $\beta = 2i$. Indeed, if we
write $\beta = \vert \beta \vert \, e^{ i \arg \beta }$, this can
simply be done by just dilating $z$ to $\lambda \, z$ with $\lambda :=
i \, \vert \frac{ 2}{ \beta} \vert^{ 1/3}\, e^{ -i \arg \beta}$, and
at the same time, we replace $w_1$ by $\frac{ 1}{ \lambda \overline{
\lambda}} \, w_1$ so as to keep unchanged the term $2i \, z \underline{
z}$ of the first equation; in the process, $\alpha$ and $\gamma$
change a bit, but we do not introduce a new name for these modified
constants. Again, the $2 \times 2$ determinant must be nonzero. Since
$\beta = 2i$ is purely imaginary, this means that $\gamma = \gamma ' +
i \, \gamma ''$ is {\rm not} purely imaginary, namely that $\gamma '
\neq 0$. Replacing $w_3$ by $w_3 - \frac{ \gamma''}{ 2} \, w_2$, one
makes $\gamma'' = 0$, {\em i.e.} $\gamma = \gamma '$ with $\gamma'
\neq 0$ real. But then a final dilation $z \mapsto \mu \, z$ for some
appropriate $\mu \in \C$ makes $\gamma'' = 2$, while a simultaneous
dilation of $w_1$ and of $w_2$ keeps unchanged the terms already
simplified.

\smallskip

What we have seen so far can be summarized in the following
basic statement.

\begin{Proposition}
\label{Proposition-initial-form}
Every real analytic 5-dimensional local CR-generic submanifold $M^5
\subset \C^4$ of codimension 3 which is maximally minimal, namely satisfies:
\[
\aligned
D^1M
&
=
T^cM
\ \ \ \ \
\text{has rank $2$},
\\
D^2M
&
=
T^cM+[T^cM,T^cM]
\ \ \ \ \
\text{has rank $3$},
\\
D^3M
&
=
T^cM+[T^cM,T^cM]
+
\big[T^cM,[T^cM,T^cM]\big]
\ \ \ \ \
\text{has maximal possible rank $5$},
\endaligned
\]
may be represented, in suitable holomorphic coordinates $(z, w_1, w_2,
w_3)$, by three complex defining equations of the specific form:
\begin{equation}
\label{initial-form}
\left[
\aligned
w_1-\overline{w}_1
&
=
2i\,z\overline{z}
+
\Pi_1\big(z,\overline{z},\overline{w}_1,\overline{w}_2,\overline{w}_3\big),
\\
w_2-\overline{w}_2
&
=
2i\,z\overline{z}(z+\overline{z})
+
\Pi_2\big(z,\overline{z},\overline{w}_1,\overline{w}_2,\overline{w}_3\big),
\\
w_3-\overline{w}_3
&
=
2\,z\overline{z}(z-\overline{z})
+
\Pi_3\big(z,\overline{z},\overline{w}_1,\overline{w}_2,\overline{w}_3\big),
\endaligned\right.
\end{equation}
where the three remainders $\Pi_1$, $\Pi_2$ and $\Pi_3$ are all an
${\rm O} ( \vert z \vert^4) + z\overline{ z}\, {\rm O}( \vert w
\vert)$ and satisfy both the reality condition~\thetag{
\ref{reality-Theta}} and the two normality conditions:
\[
0
\equiv
\Pi_j\big(0,\overline{z},\overline{w}_1,\overline{w}_2,\overline{w}_3\big)
\equiv
\Pi_j\big(z,0,\overline{w}_1,\overline{w}_2,\overline{w}_3\big)
\ \ \ \ \ \ \ \ \ \ \ \ \ {\scriptstyle{(j\,=\,1,\,\,2,\,\,3)}}.
\]

Conversely, for any choice of three such analytic functions enjoying
these conditions, the zero-locus of the three equations~\thetag{
\ref{initial-form}} above represents a real analytic 5-dimensional
local CR-generic submanifold $M^5 \subset \C^4$ of codimension $3$
which is maximally minimal.
\qed
\end{Proposition}

Since it is visibly natural to assign the weights:
\[
\mathmotsf{weight}(z)
:=
1,
\ \ \ \ \ \ \ \ \ \
\mathmotsf{weight}(w_1)
:=
2,
\ \ \ \ \ \ \ \ \ \
\mathmotsf{weight}(w_2)
:=
3,
\ \ \ \ \ \ \ \ \ \
\mathmotsf{weight}(w_3)
:=
3,
\]
one can write more-in-brief that these remainders are:
\[
\Pi_1
=
{\rm O}_{\smallmathmotsf{weighted}}(4),
\ \ \ \ \ \ \ \ \ \
\Pi_2
=
{\rm O}_{\smallmathmotsf{weighted}}(4),
\ \ \ \ \ \ \ \ \ \
\Pi_3
=
{\rm O}_{\smallmathmotsf{weighted}}(4).
\]

\subsection{Shananina's and Mamai's models}
Still in CR dimension $n = 1$, what happens for higher codimensions $d
\geqslant 4$? The above class of maximally minimal CR-generic $M^5
\subset \C^4$ already appears among the first members of Shananina's
(\cite{ Shananina-2000}) and Mamai's (\cite{ Mamai-2009}) lists. Since
the principal goal of the present memoir is to apply Cartan's method of
equivalence to these $M^5 \subset \C^4$ after having set up firm,
conceptual and computational grounds, let us briefly review some of
these results.

In codimension $d = 4$, after reductions and simplifications that are
similar to the one explained above, one adds either an equation of the
form:
\[
v_4
=
z^3\overline{z}
+
\overline{z}^3z
+
c\,z^2\overline{z}^2
+
{\rm O}_{\smallmathmotsf{weighted}}(5),
\]
where the constant $c \in \R$ happens to be a {\em nonremovable} parameter,
or an equation of the form:
\[
v_4
=
z^2\overline{z}^2
+
{\rm O}_{\smallmathmotsf{weighted}}(5).
\]

\begin{Openproblem}
For these classes of $4$-codimensional nondegenerate real analytic
CR-generic submanifolds $M^6 \subset \C^5$ of CR dimension $1$,
perform Cartan's equivalence method in order to construct an absolute
parallelism on certain appropriate principal bundles, and compute
explicitly the
related biholomorphic invariants.
\end{Openproblem}

We hope to come to that in a future publication, the interest being
that it would be the first instance\,\,---\,\,in CR
geometry\,\,---\,\,of a study of geometry-preserving deformations of
models in which a non-removable (real) parameter exists (\cite{
Kossovskiy-2012}).

In codimension $d = 5$, after reductions and simplifications
that are similar to the ones explained above, one comes
(\cite{ Shananina-2000}) either to:
\[
\left[
\aligned
v_4
&
=
z^3\overline{z}
+
\overline{z}^3z
+
b\,z^2\overline{z}^2
+
{\rm O}_{\smallmathmotsf{weighted}}(5),
\\
v_5
&
=
-i\,z^3\overline{z}
+
i\,\overline{z}^3z
+
c\,z^2\overline{z}^2
+
{\rm O}_{\smallmathmotsf{weighted}}(5),
\endaligned\right.
\]
for some two nonremovable real constants $b, c \in \R$, or to:
\[
\left[
\aligned
v_4
&
=
z^3\overline{z}
+
\overline{z}^3z
+
{\rm O}_{\smallmathmotsf{weighted}}(5),
\\
v_5
&
=
z^2\overline{z}^2
+
{\rm O}_{\smallmathmotsf{weighted}}(5).
\endaligned\right.
\]

In codimension $d = 6$, since the real vector space of real
quartic polynomials in $(z, \overline{z})$ that
are divisible by the product $z \overline{z}$
is $3$-dimensional,
generated by ${\rm Re}\, z^3\overline{ z}$,
${\rm Im}\, z^3\overline{ z}$ and
$z^2 \overline{z}^2$,
easy $\R$-linear transformations
in the $(w_4, w_5, w_6)$-space yield that the natural
polynomial
CR-generic model
$M^7 \subset \C^6$ has equations:
\[
\left[
\aligned
v_1
&
=
z\overline{z},
\\
v_2
&
=
z^2\overline{z}+\overline{z}^2z,
\\
v_3
&
=
-\,i\,z^2\overline{z}+i\,\overline{z}^2z,
\endaligned\right.
\ \ \ \ \ \ \ \ \ \ \ \ \ \ \ \
\left[
\aligned
v_4
&
=
z^3\overline{z}+\overline{z}^3z,
\\
v_5
&
=
-\,i\,z^3\overline{z}+i\,\overline{z}^3z,
\\
v_6
&
=
z^2\overline{z}^2.
\endaligned\right.
\]

Similar, more refined and nonrigid models, exists in the higher
codimensions $7 \leqslant d \leqslant 12$, {\em see}~\cite{
Shananina-2000, Mamai-2009} where the Lie algebra of infinitesimal
CR automorphisms of the corresponding models is also presented. Of
course, performing effectively Cartan's method for all these
geometries would be `fantastic'.

\section{
Symbol algebra
\\
and nilpotent Lie algebras up to dimension $5$} \label{nilpotent-5}
\HEAD{\ref{nilpotent-5}.~Symbol algebra and nilpotent
Lie algebras up to dimension $5$}{
Masoud {\sc Sabzevari} (Shahrekord) and Jo\"el {\sc Merker} (LM-Orsay)}

\subsection{Zariski-generic invariants of completely non-holonomic
complex-tangential distributions}
Now, coming back to a general real analytic CR-generic
$M \subset \C^{ n+d}$ which is connected,
for every open subset $U \subset M$, denote by:
\[
\Gamma\big(U,T^cM\big)
\]
the $\mathcal{ C}^\omega ( U)$-module of sections
of $T^c M$ over $U$, namely complex tangent vector
fields on $U$. Also, let $\Gamma \big( T^c M \big)$ denote
the sheaf of (local) sections of $T^c M$.

Next, set $D^1 M := \Gamma \big( T^c M \big)$, set:
\[
D^2M
:=
D^1M
+
\big[D^1M,\,D^1M\big]
\]
(usual Lie brackets between vector fields) and
inductively, define for every integer $k \geqslant 2$:
\[
D^{k+1}M
:=
D^kM
+
\big[D^1M,\,D^kM\big].
\]
In this way, one obtains a nested family of sheaves of local sections
of $TM$:
\[
\Gamma\big(T^cM\big)
=
D^1M
\subset
D^2M
\subset
\cdots
\subset
\Gamma\big(D^{k-1}M\big)
\subset
\Gamma\big(D^kM\big)
\subset
\cdots
\subset
\Gamma\big(TM\big).
\]

We shall assume that the complex-tangential bundle $T^cM$ is
{\sl completely nonholonomic} in the sense that at every point
$p \in M$, there exists an integer $k (p) \geqslant 1$
and there exists an open neighborhood $U_p$ of $p$ in $M$
such that:
\[
\Gamma\big(U_p,\,D^{k(p)}M\big)
\equiv
\Gamma\big(U_p,TM\big);
\]
this condition means that the collection of all possible Lie brackets between
complex-tangential fields up to a sufficiently high
order generate the whole complex tangent bundle to $M$
near $p$.

Another way of expressing this is as follows. At an
arbitrary fixed point $p \in M$,
introduce the vector subspaces of $T_p M$:
\[
\aligned
D^1M(p)
&
:=
\big\{
X_1(p)\colon\,
X_1\,\,\text{\rm is a local section of $D^1M$ near $p$}
\big\}
=
T_p^cM,
\\
\cdots\cdots
&
\cdots\cdots\cdots\cdots\cdots\cdots\cdots\cdots\cdots\cdots\cdots\cdots
\cdots\cdots\cdots\cdots
\\
D^kM(p)
&
:=
\big\{
X_k(p)\colon\,
X_k\,\,\text{\rm is a local section of $D^kM$ near $p$}
\big\},
\endaligned
\]
which are of course nested:
\[
D^1M(p)
\subset
D^2M(p)
\subset\cdots\subset
D^{k-1}M(p)
\subset
D^kM(p)
\subset
\cdots.
\]
Then this sequence grows with $k$, hence if one introduces a
notation for their dimensions:
\[
{\bf n}_k(p)
:=
D^kM(p),
\]
one has the inequalities:
\[
2n
=
{\bf n}_1(p)
\leqslant
{\bf n}_2(p)
\leqslant
\cdots
\leqslant
{\bf n}_k(p)
\leqslant
{\bf n}_{k+1}(p)
\leqslant
\cdots.
\]
Further, the assumption that $T^c M$ is completely nonholonomic
reformulates as the existence, at every point $p$, of
an integer ${\sf h}(p) \geqslant 2$ such that:
\[
{\bf n}_{{\sf h}(p)}(p)
=
2n+d
=
\dim_\R\,T_pM.
\]
By technically analyzing the ranks in a system of local charts
on $M$ in terms of matrices written in
coordinates, one may establish:

\begin{Proposition}
Under the assumption that the CR-generic submanifold $M \subset \C^{
n+d}$ is connected and real analytic, there exists a proper
exceptional real
analytic subset of $M$, an integer ${\sf h} \geqslant 2$
and integers:
\[
2n={\bf n}_1
<
{\bf n}_2
<
\cdots
<
{\bf n}_{\sf h}
=
2n+d
\]
with the properties that the degree of non-holonomy:
\[
{\sf h}(p)
=
{\sf h}
\]
is constant at every point outside
the exceptional set, and that moreover the gained dimensions:
\[
{\bf n}_1(p)
=
n_1,
\ \ \ \ \
{\bf n}_2(p)
=
{\bf n}_2,
\ \ \
\dots\dots,
\ \ \
{\bf n}_{{\sf h}(p)}(p)
=
{\bf n}_{\sf h},
\]
are also constant at every point outside the
exceptional set.
\qed
\end{Proposition}

In other words, the sheaves $D^k M$ are true {\em real vector
subbundles} of $TM$. Denote by $D^kM (p) \subset T_pM$
their fibers at points $p \in M$.

By double induction on $k\geqslant 1$ and on $l \geqslant 1$,
one verifies
using the Jacobi identity (\cite{ Tanaka-1970, Yamaguchi-1993}) that:
\[
\big[D^kM,\,D^lM\big]
\subset
D^{k+l}M.
\]
Also, if on an open subset
$U \subset M \backslash \Sigma$, one has for a certain integer
$k$ the one-step stabilization:
\[
\Gamma\big(U,D^{k+1}M\big)
=
\Gamma\big(U,D^kM\big),
\]
then this stabilization is inherited by higher order bundles:
\[
\Gamma\big(U,D^{k+l}M\big)
=
\Gamma\big(U,D^kM\big),
\]
for every $l \geqslant 1$.

\subsection{The Tanaka symbol Lie algebra}
Now, at any point $p \in M \backslash \Sigma$, introduce
for any integer $k \geqslant 2$ the quotient spaces\,\,---\,\,mind
the negative lower indices\,\,---:
\[
\mathfrak{g}_{-k}(p)
:=
D^k(p)
\big/
D^{k-1}(p),
\]
together with the canonical projections:
\[
\mathmotsf{proj}_k(p)
\colon\ \ \
D^k(p)
\longrightarrow
D^k(p)
\big/
D^{k-1}(p).
\]

\begin{Definition}
At a Zariski-generic point $p \in M \backslash \Sigma$,
the {\sl Tanaka symbol Lie algebra} is the graded sum of
the quotient vector spaces:
\[
\mathfrak{m}(p)
:=
\bigoplus_{k=1}^{k={\sf h}_M}\,
\mathfrak{g}_{-k}(p)
\]
associated to the filtration of the vector subbundles over $M
\backslash \Sigma$:
\[
D^1M
\subset
D^2M
\subset\cdots\subset
D^{k-1}M
\subset
D^kM
\subset\cdots.
\]
\end{Definition}

Thanks to taking quotients,
a natural Lie bracket $[\cdot, \cdot]$ exists on
$\mathfrak{ m} (p)$, which justifies
the name `{\sl algebra}', and it is defined as follows.
Take two elements:
\[
{\sf x}
\in
\mathfrak{g}_{-k}(p)
\ \ \ \ \ \ \ \ \ \ \ \ \ \ \
\text{\rm and}
\ \ \ \ \ \ \ \ \ \ \ \ \ \ \
{\sf y}
\in
\mathfrak{g}_{-l}(p)
\]
in two different quotients.
Both have representatives:
\[
\widetilde{\sf x}
\in
D^kM(p)
\ \ \ \ \ \ \ \ \ \ \ \ \ \ \
\text{\rm and}
\ \ \ \ \ \ \ \ \ \ \ \ \ \ \
\widetilde{\sf y}
\in
D^lM(p),
\]
which are plain tangent vectors in $T_pM$.
Take any two local vector fields defined in some
open neighborhood $U_p$ of $p$ in $M$:
\[
\widetilde{X}
\in
\Gamma\big(U_p,D^kM(p)\big)
\ \ \ \ \ \ \ \ \ \ \ \ \ \ \
\text{\rm and}
\ \ \ \ \ \ \ \ \ \ \ \ \ \ \
\widetilde{Y}
\in
\Gamma\big(U_p,D^lM(p)\big)
\]
which `extend' these two fixed vectors in the sense that:
\[
\widetilde{X}\big\vert_p
=
\widetilde{\sf x}
\ \ \ \ \ \ \ \ \ \ \ \ \ \ \
\text{\rm and}
\ \ \ \ \ \ \ \ \ \ \ \ \ \ \
\widetilde{Y}\big\vert_p
=
\widetilde{\sf y}.
\]
Next, compute the usual Lie bracket between these
two vector fields, and
then by definition, the Lie bracket between the
two original elements is its projection:
\[
\aligned
{}
[{\sf x},{\sf y}]
&
:=
\mathmotsf{proj}_{k+l}(p)
\Big(
\underbrace{
\big[
\widetilde{X},\,
\widetilde{Y}
\big]
\big\vert_p}_{
\in D^{k+l}M(p)}
\Big)
\\
&\
\in
\mathfrak{g}_{-k-l}(p).
\endaligned
\]
One verifies (\cite{ Tanaka-1970}, Lemma~1.1;
\cite{ Yamaguchi-1993}, pp.~420--421) that the result
is independent of the choices made and that:

\begin{Proposition}
Endowed with this bracket operation, the Tanaka symbol Lie algebra
$\mathfrak{ m}(p)$ for $p \in M \backslash \Sigma$
becomes a {\sl nilpotent} graded Lie algebra with:
\[
\dim_\R\,
\mathfrak{m}(p)
=
\dim_\R\,M
\]
which in addition is generated by $\mathfrak{ g}_{ -1} ( p)$:
\[
\mathfrak{g}_{-k}(p)
=
\big[
\mathfrak{g}_{-1}(p),\,\,
\mathfrak{g}_{-k+1}(p)
\big].
\qed
\]
\end{Proposition}

Conversely (\cite{ Tanaka-1970, Yamaguchi-1993}),
it is elementary to verify that:

\begin{Theorem}
Let a nilpotent Lie algebra:
\[
\mathfrak{m}
=
\bigoplus_{k=1}^{k=\mu}\,
\mathfrak{g}_{-k}
\]
be graded:
\[
\big[
\mathfrak{g}_{-k},\,\mathfrak{g}_{-l}
\big]
\subset
\mathfrak{g}_{-k-l}
\ \ \ \ \ \ \ \ \ \ \ \ \
{\scriptstyle{(\mathfrak{g}_{-\nu}\,=\,0\,\,\,
\text{\rm for}\,\,\,\nu\,\geqslant\,\mu\,+\,1)}},
\]
and satisfy the generating condition:
\[
\mathfrak{g}_{-k-1}
=
\big[
\mathfrak{g}_{-1},\,\mathfrak{g}_{-k}
\big].
\]
Then in the unique simply connected Lie group $G ( \mathfrak{ m})$
associated with $\mathfrak{ m}$, the distribution $D^1 \subset T G(
\mathfrak{ m})$ corresponding to $\mathfrak{ g}_{ -1}$ is completely
non-holonomic and its derived distributions have Lie bracket
structure isomorphic to that of $\mathfrak{ m}$.
\qed
\end{Theorem}

Lastly, without writing out a proof which could require extended
tedious technical considerations, we would like to emphasize that, on
an arbitrary real analytic CR-generic submanifold $M \subset \C^{
n+d}$, the behavior of iterated Lie brackets between sections of the
complex-tangential distribution happens to be constant at a
Zariski-generic point.  We give a precise name $\Sigma_{ CR}$ to the
appearing exceptional set, for a second one
$\Sigma_{\smallmathmotsf{Segre}}$ will be introduced in
Theorem~\ref{Sigma-Segre}, and from the point of view of studying the
biholomorphic equivalence problem at Zariski-generic points of
CR-generic real analytic submanifolds, one naturally has to avoid the
{\em union}:
\[
\Sigma_{CR}
\cup
\Sigma_{\smallmathmotsf{Segre}}.
\]

\begin{Theorem}
\label{Sigma-CR}
Under the assumption that the CR-generic submanifold $M \subset \C^{
n+d}$ is connected real analytic and that its complex-tangential
distribution $T^c M$ is completely non-holonomic at a Zariski-generic
point, there exists a certain proper real analytic subset:
\[
\Sigma_{\smallmathmotsf{CR}}
\subset
M
\]
such that
the Tanaka symbol Lie algebras $\mathfrak{ m}(p)$ of $T^cM$
are all {\em isomorphic}\,\,---\,\,hence
have same dimensional growths\,\,---\,\,
at every point:
\[
p\in
M
\backslash
\Sigma_{CR}.
\qed
\]
\end{Theorem}

For the production of model generic submanifolds which would
complement Beloshapka's approach, one could classify in small
dimensions those Lie algebras that are possible Tanaka symbol Lie
algebras for the distribution of complex-tangential planes in a
CR-generic $M \subset \C^{ n+d}$.

A more general problem concerns the classification, up to
isomorphisms, of nilpotent Lie algebras, without the generating
condition. They have been classified up to dimensions
$7$ and $8$ and we now review the concepts
and the classification up to dimension $5$.

\subsection{Isomorphic finite-dimensional Lie algebras}
Let $\mathfrak{ g}$ be a real or complex abstract Lie algebra of
finite dimension $r < \infty$, equipped with a Lie bracket operator
denoted as usual by $[ \cdot, \, \cdot ]$, or sometimes by $[ \cdot,
\, \cdot ]_{ \mathfrak{ g}}$ when a precision is needed. Two Lie
algebras $\mathfrak{ g}$ and $\widetilde{ \mathfrak{ g}}$ of the same
dimension $r = \widetilde{ r}$ are said to be {\sl isomorphic} if
there is a linear isomorphism $\phi \colon \mathfrak{ g} \to
\widetilde{ \mathfrak{ g}}$ which transfers properly the Lie bracket
structure, namely which satisfies:
\[
\phi\big([{\sf x},\,{\sf y}]_{\mathfrak{g}}\big)
=
[\phi({\sf x}),\,\phi({\sf y})]_{\widetilde{\mathfrak{g}}}
\]
for any two elements ${\sf x}, {\sf y} \in \mathfrak{ g}$.

Following Lie (Chap.~17 in~\cite{ Engel-Lie-1888}), this abstract
definition of isomorphism can be made more effective by introducing
some two linear bases ${\sf x}_1$, ${\sf x}_2$, \dots, ${\sf x}_r$ and
$\widetilde{\sf x}_1, \widetilde{\sf x}_2, \dots, \widetilde{\sf x}_r$
of $\mathfrak{ g}$ and of $\widetilde{ \mathfrak{ g}}$, so that one
can write:
\[
\mathfrak{g}
=
\C{\sf x}_1
\oplus
\C{\sf x}_2
\oplus\cdots\oplus
\C{\sf x}_r
\ \ \ \ \
\text{\rm and}
\ \ \ \ \
\widetilde{\mathfrak{g}}
=
\C\widetilde{\sf x}_1
\oplus
\C\widetilde{\sf x}_2
\oplus\cdots\oplus
\C\widetilde{\sf x}_r,
\]
as plain vector spaces. Then the datum of such two Lie algebras
$\mathfrak{ g}$ and $\widetilde{ \mathfrak{ g}}$ in terms of such
kinds of specific bases is clearly equivalent to the datum of the
so-called {\sl structure constants} $c_{ jk}^s$ and $\widetilde{ c}_{
jk}^s$ which appear in all possible brackets:
\[
[{\sf x}_j,\,{\sf x}_k]
=
\sum_{s=1}^r\,c_{jk}^s\,{\sf x}_s
\ \ \ \ \
\text{\rm and}
\ \ \ \ \
[\widetilde{\sf x}_j,\,\widetilde{\sf x}_k]
=
\sum_{s=1}^r\,c_{jk}^s\,\widetilde{\sf x}_s
\ \ \ \ \ \ \ \ \ \ \ \ \ {\scriptstyle{(j,\,\,k\,=\,1\,\cdots\,r)}}.
\]
With these notations, the two (arbitrary) $r$-dimensional Lie algebras
$\mathfrak{ g}$ and $\widetilde{ \mathfrak{ g}}$ happen to be
isomorphic if and only if one has:
\[
\phi\big([{\sf x}_j,{\sf x}_k]_{\mathfrak{g}}\big)
=
\big([\phi({\sf x}_j),\,\phi({\sf x}_k)]_{
\widetilde{\mathfrak{g}}}\big),
\]
for any two integers $j, k = 1, 2, \dots, r$. Of course,
there are constants $\phi_{ jr}$ so that
$\phi ( {\sf x}_j) = \sum_{ s=1}^r\, \phi_{ jr}\,
\widetilde{ \sf x}_r$.

We evaluate firstly
the left-hand side:
\[
\aligned
\phi\big([{\sf x}_j,{\sf x}_k]_{\mathfrak{g}}\big)
&
=
\phi\bigg(
\sum_{t=1}^r\,c_{jk}^t\,{\sf x}_t
\bigg)
=
\sum_{t=1}^r\,c_{jk}^t\,\phi({\sf x}_t)
=
\sum_{s=1}^r\,\bigg(
\underline{\sum_{t=1}^r\,\phi_{ts}\,c_{jk}^t}
\bigg)\,\widetilde{x}_s,
\endaligned
\]
and secondly, we do the same for the right-hand side:
\[
\aligned
\big([\phi({\sf x}_j),\,\phi({\sf x}_k)]_{
\widetilde{\mathfrak{g}}}\big)
&
=
\bigg[
\sum_{l=1}^r\,\phi_{jl}\,\widetilde{\sf x}_l,\,\,
\sum_{m=1}^r\,\phi_{km}\,\widetilde{\sf x}_m
\bigg]_{\widetilde{\mathfrak{g}}}
\\
&
=
\sum_{l=1}^r\,\sum_{m=1}^r\,\phi_{jl}\,\phi_{km}\,
[\widetilde{\sf x}_l,\,\widetilde{\sf x}_m]_{
\widetilde{\mathfrak{g}}}
\\
&
=
\sum_{s=1}^r\,
\bigg(
\underline{\sum_{l=1}^r\,\sum_{m=1}^r\,
\phi_{jl}\,\phi_{km}\,\widetilde{c}_{lm}^{\,s}}
\bigg)\,\widetilde{\sf x}_s.
\endaligned
\]
The two terms that have been underlined therefore identify for any $j,
k, s = 1, 2, \dots, r$, and we thus get the family of relations:
\begin{equation}
\label{c-c-non-finished}
\sum_{t=1}^r\,\phi_{ts}\,c_{jk}^t
=
\sum_{l=1}^r\,\sum_{m=1}^r\,
\phi_{jl}\,\phi_{km}\,\widetilde{c}_{lm}^{\,s}
\ \ \ \ \ \ \ \ \ \ \ \ \
{\scriptstyle{(j,\,\,k,\,\,s\,=\,1\,\cdots\,r)}}.
\end{equation}
In order to finish the computation, introduce the inverse matrix
$\phi^{ -1} ( \widetilde{ \sf x}_j) = \sum_{ l=1}^r\, \phi_{ jl}^{
-1}\, {\sf x}_l$, with hence the basic, defining properties that:
\[
\sum_{l=1}^r\,\phi_{jl}\,\phi_{lk}^{-1}
=
\delta_{jk}
\ \ \ \ \ \ \
\text{\rm and inversely:}
\ \ \ \ \ \ \
\sum_{l=1}^r\,\phi_{jl}^{-1}\,\phi_{lk}
=
\delta_{jk},
\]
for any $j, k = 1, 2, \dots, r$. Thus, we may multiply
\thetag{ \ref{c-c-non-finished}} by $\phi_{ su}^{ -1}$,
where $u$ is arbitrary between $1$ and $r$, and
sum with respect to $s$, which yields:
\[
\aligned
c_{jk}^u
=
\sum_{t=1}^r\,\delta_{tu}\,c_{jk}^t
&
=
\sum_{s=1}^r\,\sum_{t=1}^r\,
\phi_{ts}\,\phi_{su}^{-1}\,c_{jk}^t
=
\sum_{l=1}^r\,\sum_{m=1}^r\,\sum_{s=1}^r\,
\phi_{jl}\,\phi_{km}\,\phi_{su}^{-1}\,
\widetilde{c}_{lm}^{\,s}
\\
&
\ \ \ \ \ \ \ \ \ \ \ \ \ \ \
{\scriptstyle{(j,\,\,k,\,\,u\,=\,1\,\cdots\,r)}}.
\endaligned
\]
This is the way how the two collection of structure constants
must be related when there exists a Lie algebra isomorphism
between $\mathfrak{ g}$ and $\widetilde{ \mathfrak{ g}}$
and we summarize as follows the result gained.

\begin{Proposition}
{\rm (\cite{Engel-Lie-1888, Merker-Engel-Lie}, Chap.~17)}
Consider two arbitrary real or complex Lie algebras
having the same dimension $r$ which are described by generators:
\[
\mathfrak{g}
=
{\rm Span}\big({\sf x}_1,{\sf x}_2,\dots,{\sf x}_r\big)
\ \ \ \ \ \ \
\text{\rm and}
\ \ \ \ \ \ \
\mathfrak{g}'
=
{\rm Span}\big(\widetilde{\sf x}_1,\widetilde{\sf x}_2,
\dots,\widetilde{\sf x}_r\big)
\]
and whose Lie bracket structures, in terms of their respective
generators ${\sf x}_k$ and $\widetilde{ \sf x}_k$, read:
\[
[{\sf x}_j,{\sf x}_k]
=
\sum_{s=1}^r\,c_{jk}^s\,{\sf x}_s
\ \ \ \ \ \ \
\text{\rm and}
\ \ \ \ \ \ \
[\widetilde{\sf x}_j,\widetilde{\sf x}_k]
=
\sum_{s=1}^r\,\widetilde{c}_{jk}^{\,s}\,\widetilde{\sf x}_s,
\]
where the $c_{ jk}^s$ and the $\widetilde{ c}_{ jk}^{\,s}$
are constants subjected to skew symmetry and to Jacobi identity:
\[
\!\!\!\!\!\!\!\!\!\!\!\!\!\!\!\!\!\!\!\!
\left\{
\aligned
0
&
=
c_{jk}^s+c_{kj}^s,
\\
0
&
=
{\textstyle{\sum_{s=1}^r}}
\big(
c_{kl}^s\,c_{js}^m
+
c_{jk}^s\,c_{ls}^m
+
c_{lj}^s\,c_{ks}^m
\big)
\endaligned\right.
\ \ \ \ \ \ \
\text{\rm and}
\ \ \ \ \ \ \
\left\{
\aligned
0
&
=
\widetilde{c}_{jk}^{\,s}+\widetilde{c}_{kj}^{\,s},
\\
0
&
=
{\textstyle{\sum_{s=1}^r}}
\big(
\widetilde{c}_{kl}^{\,s}\,\widetilde{c}_{js}^{\,m}
+
\widetilde{c}_{jk}^{\,s}\,\widetilde{c}_{ls}^{\,m}
+
\widetilde{c}_{lj}^{\,s}\,\widetilde{c}_{ks}^{\,m}
\big).
\endaligned\right.
\]
Then $\mathfrak{ g}$ and $\widetilde{ \mathfrak{ g}}$ are isomorphic
as Lie algebras if and only it is possible to find a system of $r^2$
real or complex numbers $\phi_{ lm}$ with nonzero $\det ( \phi_{
lm})_{ 1 \leqslant l \leqslant r}^{ 1 \leqslant m \leqslant r} \neq
0$ such that the two collections of constant structures exchange
through the following formulas:
\[
c_{jk}^s
=
\sum_{l=1}^r\,\sum_{m=1}^r\,\sum_{t=1}^r\,
\phi_{jl}\,\phi_{km}\,\phi_{ts}^{-1}\,
\widetilde{c}_{lm}^{\,t},
\]
where $j$, $k$, $s$ are arbitrary integers between $1$ and $r$,
or equivalently, through the inverse formulas:
\[
\widetilde{c}_{jk}^s
=
\sum_{l=1}^r\,\sum_{m=1}^r\,\sum_{t=1}^r\,
\phi_{jl}^{-1}\,\phi_{km}^{-1}\,\phi_{ts}\,
c_{lm}^t.
\]
\end{Proposition}

\subsection{Decreasing weak derived sequence}
Again,
let $\mathfrak{ g}$
be a real or complex abstract Lie algebra of finite dimension $r <
\infty$, equipped with a Lie bracket operator denoted as usual by $[
\cdot, \, \cdot ]$. Introduce the following sequence of subspaces of
$\mathfrak{ g}$:
\[
\left\{
\aligned
\mathcal{N}^{-1}(\mathfrak{g})
&
:=
\mathfrak{g},
\\
\mathcal{N}^{-2}(\mathfrak{g})
&
:=
[\mathfrak{g},\,\mathfrak{g}]
=
\big[\mathfrak{g},\,\mathcal{N}^{-1}(\mathfrak{g})\big]
\\
\mathcal{N}^{-k-1}(\mathfrak{g})
&
:=
\big[\mathfrak{g},\,\mathcal{N}^{-k}(\mathfrak{g})\big]
\ \ \ \ \ \
\text{\rm for any}\ \
k\geqslant 2.
\endaligned\right.
\]
Then the arising {\sl weak derived sequence}\footnote{\,
It should be distinguished from the {\sl strong derived sequence} of
iterated commutators starting also with $\mathcal{ R}^{ -2} (
\mathfrak{ g}) := [ \mathfrak{ g}, \mathfrak{ g}]$ but continuing with
$\mathcal{ R}^{ - k - 1} ( \mathfrak{ g}) := \big[ \mathcal{ R}^{ - k}
( \mathfrak{ g}), \, \mathcal{ R}^{ - k} (\mathfrak{ g}) \big]$, so
that $\mathcal{ R}^{ - k} (\mathfrak{ g}) \subset \mathcal{ N}^{ - k}
(\mathfrak{ g})$ for any $k \geqslant -1$.
}: 

\[
\cdots
\subset
\mathcal{N}^{-k-1}(\mathfrak{g})
\subset
\mathcal{N}^{-k}(\mathfrak{g})
\subset\cdots\subset
\mathcal{N}^{-2}(\mathfrak{g})
\subset
\mathcal{N}^{-1}(\mathfrak{g})
=
\mathfrak{g}
\]
is constituted only of {\em ideals} of $\mathfrak{ g}$, namely of
subalgebras $\mathfrak{ n}$ of $\mathfrak{ g}$ satisfying $[
\mathfrak{ n}, \, \mathfrak{ g}] \subset \mathfrak{ n}$, as may be
verified by induction using the Jacobi identity. Also, because
$\mathfrak{ g}$ is finite-dimensional, this sequence must stabilize
after some steps, that is to say, there must exist an integer $k
\geqslant 1$ such that $\mathcal{ N}^{ - k - 1} ( \mathfrak{ g}) =
\mathcal{ N}^{ - k} ( \mathfrak{ g})$.

\begin{Definition}
A real or complex Lie algebra $\mathfrak{ g}$ is said to be {\sl nilpotent}
if its weak derived sequences ends up to zero, namely if there is an
integer $-k$ such that:
\[
\mathcal{N}^{-k}(\mathfrak{g})
=
\{0\}.
\]
When this occurs, the smallest integer $\mu + 1$ such that $\mathcal{ N}^{
- \mu - 1} ( \mathfrak{ g}) = \{ 0\}$ is called the {\sl nilindex} of
$\mathfrak{ g}$, while the integer $\mu$, the largest such that
$\mathcal{ N}^{ - \mu} (\mathfrak{ g} )\neq 0$, is called the
{\sl kind} of $\mathfrak{ g}$ ({\em cf.}~\cite{ Tanaka-1970}).
\end{Definition}

The nilindex of a nilpotent Lie algebra $\mathfrak{ g}$ is
of course $\leqslant \dim \mathfrak{ g} + 1$.

As a basic example, a Lie algebra is Abelian, namely $[ \mathfrak{ g},
\mathfrak{ g}] = 0$, if and only if it is nilpotent with nilindex
smallest possible, equal to $2$. Another example, in dimension $r =
3$, is the Heisenberg Lie algebra ${\rm Span}_\C ( {\sf x}_1, {\sf
x}_2, {\sf x}_3)$, having structure:
\[
[{\sf x}_1,{\sf x}_2]
=
{\sf x}_3,
\ \ \ \ \ \ \
[{\sf x}_1,{\sf x}_3]
=
0,
\ \ \ \ \ \ \
[{\sf x}_2,{\sf x}_3]
=
0.
\]

A fundamental tool for the classification of nilpotent Lie algebras is
the so-called classical {\sl theorem of Engel} (\cite{ Knapp-2002}),
which states that a finite-dimensional $\mathfrak{ g}$ is
nilpotent if and only if, for every ${\sf x} \in \mathfrak{ g}$, the
associated adjoint endomorphism:
\[
\aligned
{\rm ad}({\sf x})\colon\ \ \
\mathfrak{g}
&
\longrightarrow\mathfrak{g}
\\
{\sf y}
&
\longmapsto
\big[{\sf x},{\sf y}]
\endaligned
\]
is nilpotent in ${\rm End} ( \mathfrak{ g})$, namely ${\rm ad} ( {\sf
x})^{ \circ s} \equiv 0$ for all $s$ large enough. It is well known,
then, that in fact ${\rm ad} ( {\sf x})^{ \circ \dim \mathfrak{ g}}
\equiv 0$, uniformly for all ${\sf x}$ belonging to the nilpotent
$\mathfrak{ g}$. But since ${\rm ad} ( {\sf x}) ( {\sf x}) = [ {\sf
x}, {\sf x}] = 0$, which shows that ${\sf x}$ is an eigenvector with
zero eigenvalue, in the Jordan bloc decomposition by invariant linear
subspaces, the largest dimension of an invariant subspace on which
${\rm ad} ( {\sf x})$ is not identically zero is in any case
$\leqslant \dim \mathfrak{ g} - 1$. It follows that:
\[
{\rm ad}({\sf x})^{\circ(\dim\mathfrak{g}-1)}
\equiv
0
\ \ \ \ \ \ \ \ \ \ \ \ \
{\scriptstyle{({\sf x}\,\,\in\,\,\mathfrak{g})}}.
\]

\begin{Definition}
If $\mathfrak{ g}$ is nilpotent, the smallest integer $s
\leqslant \dim \mathfrak{ g} - 1$ such that
${\rm ad} ( {\sf x})^{ \circ s} \equiv 0$ for all ${\sf x} \in
\mathfrak{ g}$ is called the {\sl nilpotency order} of $\mathfrak{
g}$. A nilpotent complex Lie algebra of dimension $r$ is said to be
{\sl filiform} if its nilpotency order equals $\dim \mathfrak{ g} - 1$.
\end{Definition}

For instance, the four-dimensional Lie algebra ${\rm Span}_\C ( {\sf
x}_1, {\sf x}_2, {\sf x}_3, {\sf x}_4)$ having structure:
\[
[{\sf x}_1,{\sf x}_2]={\sf x}_3,
\ \ \ \ \ \ \ \ \ \ \
[{\sf x}_1,{\sf x}_3]={\sf x}_4
\]
is filiform: look at ${\rm ad} ( {\sf x}_1)$.

\begin{Definition}
Let $\mathfrak{ g}$ be a nilpotent complex Lie algebra of dimension
$n$. For every ${\sf x} \in \mathfrak{ g}$, let $c ( {\sf x})$ be the
decreasing sequence of the dimensions of Jordan blocs of the nilpotent
endomorphism ${\rm ad} ( {\sf x})$. The {\sl characteristic sequence}
of $\mathfrak{ g}$ (Goze's invariant) is the sequence:
\[
c(\mathfrak{g})
:=
\max
\big\{
c({\sf x})\colon\,
{\sf x}\in\mathfrak{g}\backslash[\mathfrak{g},\mathfrak{g}]
\big\}.
\]
\end{Definition}

The characteristic sequence then obviously constitutes a {\sl
partition} of $r$ and will be written as $c ( \mathfrak{ g}) =
(\ell_1, \ell_2, \dots, \ell_r)$ with $\ell_1 \geqslant \ell_2
\geqslant \cdots \geqslant \ell_r \geqslant 0$ and $\ell_1 + \ell_2 +
\cdots + \ell_r = r$. A nilpotent Lie algebra $\mathfrak{ g}$ is
filiform if and only if $c ( \mathfrak{ g}) = (r-1, 1, 0, \dots, 0)$.

\subsection{Classification of complex nilpotent Lie algebras up to
dimension five}
In the next paragraphs, we follow Goze's classification results~\cite{
Goze-web, Goze-2008} closely. We now consider an arbitrary
$r$-dimensional nilpotent {\em complex} Lie algebra $\mathfrak{ g}$ of
dimension $r \leqslant 5$ with generators ${\sf x}_1$, ${\sf x}_2$,
\dots, ${\sf x}_r$.

\smallskip\noindent{\footnotesize\sf Dimension 1:}
There exists only one nilpotent Lie algebra of dimension $1$: the
Abelian Lie algebra, and we will denote it by: $\mathfrak{ a}_1$.

\smallskip\noindent{\footnotesize\sf Dimension 2:}
In dimension $2$, there is no indecomposable complex nilpotent Lie
algebra, only $\mathfrak{ a}_2 := \mathfrak{ a}_1 \oplus \mathfrak{
a}_1$ exists.

\smallskip\noindent{\footnotesize\sf Dimension 3:}
In dimension $3$, leaving aside $\mathfrak{ a}_3 := \mathfrak{ a}_1
\oplus \mathfrak{ a}_1 \oplus \mathfrak{ a}_1$, there exists a single
indecomposable complex nilpotent Lie algebra, the Heisenberg Lie
algebra:
\[
\boxed{
\mathfrak{n}_3^1\colon
\ \ \ \ \ \ \
[{\sf x}_1,\,{\sf x}_2]
=
{\sf x}_3}\,.
\]
By convention here, brackets that are not written are implicitly
assumed to be zero. We shall observe later that $\mathfrak{ n}_1^3$
is the Tanaka symbol algebra of any Levi nondegenerate hypersurface
$M^3 \subset \C^2$.

\smallskip\noindent{\footnotesize\sf Dimension 4:}
In dimension $4$, leaving aside the two decomposable
nilpotent complex Lie algebras,
\[
\mathfrak{a}_4:=\mathfrak{a}_1^{\oplus 4}
\ \ \ \ \
\text{\rm and}
\ \ \ \ \
\mathfrak{a}_1\oplus\mathfrak{n}_1^3,
\]
there again exists only a single indecomposable complex
nilpotent Lie algebra, whose structure is:
\[
\boxed{
\mathfrak{n}_4^1\colon
\ \ \ \ \ \ \ \ \ \
\left\{
\aligned
{}
[{\sf x}_1,{\sf x}_2]
&
=
{\sf x}_3
\\
[{\sf x}_1,{\sf x}_3]
&
=
{\sf x}_4.
\endaligned\right.}
\]

\smallskip\noindent{\footnotesize\sf Dimension 5:}
Next, in dimension $r = 5$,
leaving aside the
three decomposable nilpotent complex Lie algebras:
\[
\mathfrak{a}_5:=\mathfrak{a}_1^{\oplus 5},
\ \ \ \ \ \ \
\mathfrak{n}_3^1\oplus\mathfrak{a}_2
\ \ \ \ \ \ \
\text{\rm and}
\ \ \ \ \ \ \
\mathfrak{n}_4^1\oplus\mathfrak{a}_1,
\]
there exist six mutually nonisomorphic nilpotent complex Lie algebras
that are gathered as follows according to their respective Goze
invariants.

\smallskip$\square$
\underline{\small\sf $c(\mathfrak{g})=(4,1)$ (filiform case):}
\[
\boxed{
\mathfrak{n}_5^1\colon
\ \ \ \ \ \ \ \ \ \
\left\{
\aligned
{}
[{\sf x}_1,{\sf x}_2]
&
=
{\sf x}_3
\\
[{\sf x}_1,{\sf x}_3]
&
=
{\sf x}_4
\\
[{\sf x}_1,{\sf x}_4]
&
=
{\sf x}_5
\endaligned\right.}
\ \ \ \ \ \ \ \ \ \ \ \ \ \ \ \ \
\boxed{
\mathfrak{n}_5^2\colon
\ \ \ \ \ \ \ \ \ \
\left\{
\aligned
{}
[{\sf x}_1,{\sf x}_2]
&
=
{\sf x}_3
\\
[{\sf x}_1,{\sf x}_3]
&
=
{\sf x}_4
\\
[{\sf x}_1,{\sf x}_4]
&
=
{\sf x}_5
\\
[{\sf x}_2,{\sf x}_3]
&
=
{\sf x}_5
\endaligned\right.}
\]

\smallskip$\square$
\label{n-5-4}
\underline{\small\sf $c(\mathfrak{g})=(3,1,1)$:}
\[
\boxed{
\mathfrak{n}_5^3\colon
\ \ \ \ \ \ \ \ \ \
\left\{
\aligned
{}
[{\sf x}_1,{\sf x}_2]
&
=
{\sf x}_3
\\
[{\sf x}_1,{\sf x}_3]
&
=
{\sf x}_4
\\
[{\sf x}_2,{\sf x}_5]
&
=
{\sf x}_4
\endaligned\right.}
\ \ \ \ \ \ \ \ \ \ \ \ \ \ \ \ \
\boxed{
\mathfrak{n}_5^4\colon
\ \ \ \ \ \ \ \ \ \
\left\{
\aligned
{}
[{\sf x}_1,{\sf x}_2]
&
=
{\sf x}_3
\\
[{\sf x}_1,{\sf x}_3]
&
=
{\sf x}_4
\\
[{\sf x}_2,{\sf x}_3]
&
=
{\sf x}_5
\endaligned\right.}
\]

\smallskip$\square$
\underline{\small\sf $c(\mathfrak{g})=(2,2,1)$:}
\[
\boxed{
\mathfrak{n}_5^5\colon
\ \ \ \ \ \ \ \ \ \
\left\{
\aligned
{}
[{\sf x}_1,{\sf x}_2]
&
=
{\sf x}_3
\\
[{\sf x}_1,{\sf x}_4]
&
=
{\sf x}_5
\endaligned\right.}
\]

\smallskip$\square$
\underline{\small\sf $c(\mathfrak{g})=(2,1,1,1)$:}
\[
\boxed{
\mathfrak{n}_5^5\colon
\ \ \ \ \ \ \ \ \ \
\left\{
\aligned
{}
[{\sf x}_1,{\sf x}_2]
&
=
{\sf x}_3
\\
[{\sf x}_4,{\sf x}_5]
&
=
{\sf x}_3
\endaligned\right.}
\]

\subsection{Graded nilpotent Lie algebras}
On p.~11 of~\cite{Tanaka-1970},
one finds up to dimension $5$ the possible dimensional growths of
graded nilpotent Lie algebras which come from the Tanaka symbol of a
$2$-dimensional distribution, but without
the Lie bracket structure. In fact, the corresponding structures
may be read off from the above list:

\medskip\noindent$\square$
in dimension $3$, with dimensional growth $(2, 1)$, one finds
$\mathfrak{ n}_3^1$;

\medskip\noindent$\square$
in dimension $4$, with dimensional growth $(2, 1, 1)$, one finds
$\mathfrak{ n}_4^1$;

\medskip\noindent$\square$
in dimension $5$, with dimensional growth $(2,1, 2)$, one finds
$\mathfrak{ n}_5^4$ which corresponds to our cubic $M_{\sf c}^5
\subset \C^4$;

\medskip\noindent$\square$
in dimension $5$, with growth vector $(2, 1, 1, 1)$, one
finds $\mathfrak{ n}_5^1$, and this would correspond
to a yet unstudied class of CR-generic submanifolds
$M^5 \subset \C^4$.

\begin{Openproblem}
Classify, in small dimensions, real Lie algebras:
\[
\mathfrak{g}_{-\mu}
\oplus
\cdots
\oplus
\mathfrak{g}_{-2}
\oplus
\mathfrak{g}_{-1}
\]
that are graded, nilpotent and satisfy the generating
condition:
\[
\mathfrak{g}_{-k-1}
=
\big[\mathfrak{g}_{-1},\,\mathfrak{g}_{-k}\big].
\]
Apply the gained classifications
to coordinate-independent production of model CR-generic
submanifolds of complex Euclidean spaces.
\end{Openproblem}

\section{Infinitesimal CR automorphisms:
$\mathfrak{ aut}_{CR} ( M) = {\rm Re} (\mathfrak{ hol}(M))$}
\label{infinitesimal-CR}
\HEAD{\ref{infinitesimal-CR}.~Infinitesimal CR automorphisms:
$\mathfrak{ aut}_{CR} ( M) = {\rm Re} (\mathfrak{ hol}(M))$}{
Masoud {\sc Sabzevari} (Shahrekord) and Jo\"el {\sc Merker} (LM-Orsay)}

\subsection{Extrinsic holomorphic definition}
According to a standard, important definition (\cite{ Stanton-1996,
Beloshapka-1998, BES-2007}), a {\sl (local) infinitesimal
CR-automorphism of $M$} is a $(1, 0)$ vector field
having {\em holomorphic} coefficients:
\begin{equation}
\label{X-infinitesimal-CR}
{\sf X}
=
\sum_{k=1}^n\,Z^k(z,w)
\frac{\partial}{\partial z_k}
+
\sum_{j=1}^d\,W^j(z,w)
\frac{\partial}{\partial w_j}
\end{equation}
the real part of which:
\[
{\rm Re}\,X
=
{\textstyle{\frac{1}{2}}}
\big(
{\sf X}+\overline{\sf X}
\big)
\]
which is tangent to $M$. Importantly, one should notice here that, contrary
to the $(1, 0)$ generators $\mathcal{ L}_k$ of $T^{ 1, 0} M$, such an
${\sf X}$ is supposed to have purely holomorphic coefficients, whereas
the coefficients $\frac{ \partial \Theta_j}{ \partial z_k} ( z,
\overline{ z}, \overline{ w})$ of the $\mathcal{ L}_k$
are\,\,---\,\,most often\,\,---\,\,neither purely holomorphic, nor
purely antiholomorphic, but only real analytic.

This condition of tangency,
much studied by Beloshapka and his school, will be
explored in depth below because knowing all such ${\sf X}$ is the same
as knowing the {\em CR symmetries} of $M$,
and this knowledge lies in the heart of
the problem of classifying local analytic CR manifolds up to
biholomorphisms.

By integration, the {\em real} flow:
\[
(t,z,w)
\longmapsto
\exp(t\,{\sf X})(z,w)
\ \ \ \ \ \ \ \ \ \ \ \ \
{\scriptstyle{(t\,\in\,\R\,\,\,
\text{\rm small})}}
\]
constitutes a local one-parameter group of local biholomorphisms of
$\C^n$, and because ${\sf X}$ is tangent to $M$, this flow leaves $M$
invariant, that is to say: through this flow, points of $M$ are
transferred to points of $M$ (more details may be found in~\cite{
Stanton-1996}). We note {\em passim} that this
real flow coincides with restricting the consideration of the {\em
complex} (holomorphic) flow:
\[
(\tau,z,w)
\longmapsto
\exp(\tau\,{\sf X})(z,w)
\ \ \ \ \ \ \ \ \ \ \ \ \
{\scriptstyle{(\tau\,\in\,\C\,\,\,
\text{\rm small})}}
\]
to a {\em real} time parameter $\tau := t \in \R$.
Conversely, one may show:

\begin{Lemma}
If $M \subset \C^{n+d}$ is a CR-generic submanifold and if $(z, w)
\longmapsto \phi_t ( z, w)$ is a local one-parameter group of {\em
holomorphic} self-transformations of $\C^{n+d}$ which stabilizes $M$
locally, then the vector field:
\[
{\textstyle{\frac{d}{dt}}}\big\vert_0
\big(
\phi_t(z,w)
\big)
\]
has holomorphic coefficients and its real part is tangent to $M$.
\qed
\end{Lemma}

From fundamental facts of of Lie theory, if $\mathfrak{ hol} (M)$ is
finite-dimensional,
then necessarily, it constitutes a {\em real} Lie algebra, namely if
${\sf X}_1, \dots, {\sf X}_r$ denote any basis of $\mathfrak{ hol}(
M)$, there are {\em real} structure constants $c_{ jk}^s \in \R$ such
that:
\begin{equation}
\label{real-structure-constants}
\big[{\sf X}_j,\,{\sf X}_k\big]
=
\sum_{s=1}^r\,c_{jk}^s\,{\sf X}_s.
\end{equation}
For an explicitly given $M \subset \C^{ n+d}$, determining a basis of
the Lie algebra $\mathfrak{ hol} ( M)$ is a natural problem for which
some systematic procedures exist ({\em see} below). The
groundbreaking works of Sophus Lie and his collaborators,
Friedrich Engel, Georg Scheffers and others showed that the
most fundamental question in concern here is to draw lists of possible
Lie algebras $\mathfrak{ hol} ( M)$ which would classify
possible $M$'s according to their CR symmetries,
cf.~\cite{Engel-Lie-1888, Engel-Lie-1893, Merker-Engel-Lie,
Merker-Hermann}.

\subsection{Intrinsic CR definition}
On the other hand, if one prefers to view the CR-generic manifold $M$ in a
purely intrinsic way, one may consider the local group ${\rm Aut}_{
CR} ( M)$ of automorphisms of the CR structure, namely of local
$\mathcal{ C}^\infty$ diffeomorphisms $g \colon M \to M$ (close to the
identity mapping) which satisfy:
\[
dg_p(T_p^cM)
=
T_{g(p)}^cM
\ \ \ \ \
\text{\rm and}
\ \ \ \ \
dg_p\big(J({\sf v}_p)\big)
=
J_{g(p)}\big(dg_p({\sf v}_p)\big)
\]
at any point $p \in M$ and for any vector ${\sf v}_p \in T_p^c M$. In other
words, $g$ belongs to ${\rm Aut}_{ CR} ( M)$ if and only if it is a
(local) CR-diffeomorphism of $M$, namely a diffeomorphism which
respects the (intrinsic) CR structure of $M$.

As did Lie most of the
time in his original theory (\cite{ Engel-Lie-1888, Engel-Lie-1893}), we
shall consider only a neighborhood of the identity mapping, hence all
our groups will be {\sl local Lie groups}; the reader
is referred to~\cite{ Olver-1986, Merker-Engel-Lie} for fundamentals about
local Lie groups in general, especially concerning the fact that it is
essentially useless to point out open sets and domains in which
mappings and transformations are defined, some superfluous
details we shall dispense ourselves with.

Accordingly, let:
\[
\mathfrak{aut}_{CR}(M)
\]
denote the collection of all (real) vector fields ${\sf Y}$ on $M$ the
flow of which $(t, p) \mapsto \exp ( t {\sf Y}) ( p)$ becomes a local
CR diffeomorphism of $M$. When ${\rm Aut}_{ CR} (M)$ is a
finite-dimensional Lie group, $\mathfrak{ aut}_{ CR} ( M)$ is just its
Lie algebra. The principles and the proof of the following assertion
date back to Sophus Lie's monographs.

\begin{Lemma}
{\rm (\cite{Engel-Lie-1888}, Chap.~8)} A local real
analytic vector field ${\sf Y}$ on $M$ belongs to $\mathfrak{ aut}_{
CR} ( M)$, if and only if for every local section $L$ of the complex
tangent bundle $T^cM$, the Lie bracket $[ {\sf Y}, \, L]$ is again a
section of $T^cM$.
\qed
\end{Lemma}

\subsection{Coincidence between extrinsic and intrinsic CR automorphisms}
In all cases which are of interest, namely when $M$ is nondegenerate
in a sense that we will make precise just later, such
real analytic flows $(t, p) \mapsto \exp ( t {\sf Y}) ( p)$ happen
to extend as local {\em biholomorphic} maps from a neighborhood of
$M$ in $\C^{n+d}$. In all these cases which cover a broad universe of
yet unstudied CR structures, one has the fundamental relation:
\[
\boxed{
\mathfrak{aut}_{CR}(M)
=
{\rm Re}
\big(
\mathfrak{hol}(M)
\big)}\,,
\]
where both sides are finite-dimensional, spanned by vector fields
whose coefficients are expandable in converging power series. Thus,
one may work exclusively with the {\em holomorphic} vector fields
generating $\mathfrak{ hol} ( M)$, as we will do from now on. And in
any case, there will be no confusion to call an {\em infinitesimal CR
automorphism} either the holomorphic vector field ${\sf X} \in
\mathfrak{ hol} ( M)$ or its real part $\frac{ 1}{ 2} ( {\sf X} +
\overline{\sf X}) \in \mathfrak{ aut}_{CR} ( M)$.

Since holomorphic vector fields obviously commute with antiholomorphic
vector fields, we deduce from~\thetag{ \ref{real-structure-constants}}
that when $\mathfrak{ hol} ( M) = \R {\sf X}_1 \oplus \cdots \oplus \R
{\sf X}_r$ is $r$-dimensional, the real parts of the ${\sf X}_j$ which
generate $\mathfrak{ aut}_{ CR} ( M)$ simply have the same (real)
structure constants:
\begin{equation}
\aligned
\big[{\sf X}_j+\overline{\sf X}_j,\,
{\sf X}_k+\overline{\sf X}_k\big]
&
=
\big[{\sf X}_j,\,{\sf X}_k\big]
+
\big[\overline{\sf X}_j,\,\overline{\sf X}_k\big]
\\
&
=
\sum_{s=1}^r\,c_{jk}^s\,
\big({\sf X}_s+\overline{\sf X}_s\big).
\endaligned
\end{equation}

\subsection{Isotropy Lie subalgebras} At a fixed point $p \in M$,
one may consider the Lie subalgebras $\mathfrak{ hol}( M, p)$ of
$\mathfrak{ hol} ( M)$ and $\mathfrak{ aut}_{ CR}( M, p)$ of
$\mathfrak{ aut}_{ CR} ( M)$ consisting of vector fields whose values
vanish at $p$. These are the Lie algebra of the subgroups ${\rm Hol}
( M, p)$ of ${\rm Hol} ( M)$ and ${\rm Aut}_{ CR} ( M, p)$ of ${\rm
Aut}_{ CR} ( M)$ consisting of maps that fix the point $p$. One has
$\mathfrak{ aut}_{ CR} ( M, p) = {\rm Re} \big( {\sf hol} ( M, p)
\big)$.

\subsection{Extrinsic complexification}
\label{extrinsic-complexification}
As is known in local CR geometry, it is natural to introduce new
independent complex variables $(\underline{ z}, \underline{ w}) \in
\C^n \times \C^d$\,\,---\,\,underlining
here should {\em not} be confused with complex
conjugating\,\,---\,\,and to define the so-called {\sl extrinsic
complexification} $M^{e_c}$ of $M$ as being the {\em holomorphic}
$d$-codimensional submanifold of $\C^{n+d}\times \C^{n+d}$ equipped
with the $2n + 2d$ coordinates $(z, w, \underline{ z}, \underline{
w})$ which is defined by the $d$ holomorphic equations:
\[
w_j = \Theta_j\big(z,\underline{z},\underline{w}\big)
\ \ \ \ \ \ \ \ \ \ \ \ \
{\scriptstyle{(j\,=\,1\,\cdots\,d)}}.
\]
Notice that the replacement of $( \overline{ z}, \overline{ w})$ by $(
\underline{ z}, \underline{ w})$ in the Taylor series of $\Theta$ is
really meaningful:
\[
\Theta_j(z,\underline{z},\underline{w}) :=
\sum_{\alpha\in\mathbb{N}^n,\,\beta\in\mathbb{N}^n,\,\gamma\in\mathbb{N}^d}\,
\Theta_{j,\alpha,\beta,\gamma}\,
z^\alpha\,\underline{z}^\beta\,\underline{w}^\gamma,
\]
thanks to the fact that the series in question converges locally.

Equivalently, $M^{ e_c}$ is defined by the $d$ equations:
\[
\underline{w}_j
=
\overline{\Theta}_j
\big(\underline{z},z,w\big)
\ \ \ \ \ \ \ \ \ \ \ \ \
{\scriptstyle{(j\,=\,1\,\cdots\,d)}}.
\]
Then $M$ is recovered from $M^{ e_c}$ by just replacing these
independent variables $( \underline{ z}, \underline{ w})$ by the
original conjugates $( \overline{ z}, \overline{ w})$.

Of course, the extrinsic complexifications of the $(1, 0)$
and of the $(0, 1)$ tangent vector fields are:
\[
\aligned \mathcal{L}_k^{e_c} & = \frac{\partial}{\partial z_k} +
\sum_{j=1}^d\, \frac{\partial\Theta_j}{\partial z_k}
(z,\underline{z},\underline{w})\, \frac{\partial}{\partial w_j} \ \ \
\ \ \ \ \ \ \ \ \ \ {\scriptstyle{(k\,=\,1\,\cdots\,n)}},
\\
\underline{\mathcal{L}}_k^{e_c} & =
\frac{\partial}{\partial\underline{z}_k} + \sum_{j=1}^d\,
\frac{\partial\overline{\Theta}_j}{\partial\underline{z}_k}
(\underline{z},z,w)\, \frac{\partial}{\partial\underline{w}_j} \ \ \
\ \ \ \ \ \ \ \ \ \ {\scriptstyle{(k\,=\,1\,\cdots\,n)}}.
\endaligned
\]

Lastly, we shall constantly use the following standard
uniqueness principle.

\begin{Lemma}
\label{uniqueness-lemma}
With a CR-generic real analytic $M \subset \C^{ n+d}$ as above,
consider a complex-valued converging power series:
\[
\Phi
=
\Phi\big(z,w,\overline{z},\overline{w}\big)
=
\sum_{\alpha\in\mathbb{N}^n,\,\beta\in\mathbb{N}^d,\,
\gamma\in\mathbb{N}^n,\,\delta\in\mathbb{N}^d}\,
\Phi_{\alpha,\beta,\gamma,\delta}\,
z^\alpha\,w^\beta\,\overline{z}^\gamma\,\overline{w}^\delta
\]
in $\C \{ z, w, \overline{ z}, \overline{ w} \}$ having complex
coefficients $\Phi_{ \alpha, \beta, \gamma, \delta} \in \C$. Then the
following four properties are equivalent:

\begin{itemize}

\smallskip\item[$\bullet$]
$\Phi$ takes only the value zero when the point $(z, w)$ varies
on $M \subset \C^n$;

\smallskip\item[$\bullet$]
the extrinsic complexification of $\Phi$:
\[
\Phi^{e_c}
=
\Phi^{e_c}\big(z,w,\underline{z},\underline{ w}\big)
:=
\sum_{\alpha,\,\beta,\,
\gamma,\,\delta}\,
\Phi_{\alpha,\beta,\gamma,\delta}\,
z^\alpha\,w^\beta\,\underline{z}^\gamma\,\underline{w}^\delta
\]
takes only the value zero when the point $(z, w, \underline{ z},
\underline{ w})$ varies on $M^{ e_c} \subset \C^{ 2n+2d}$;

\smallskip\item[$\bullet$]
after replacing $\underline{ w}$ by $\overline{ \Theta} ( \underline{
z}, z, w)$ in the extrinsic complexification $\Phi^{ e_c}$ of $\Phi$,
the result is an identically zero series in $\C\{ \underline{ z}, z,
w \}$:
\[
0 \equiv
\sum_{\alpha,\,\beta,\,\gamma,\,\delta}
\Phi_{\alpha,\beta,\gamma,\delta}\,
z^\alpha\,w^\beta\,\underline{z}^\gamma\,
\big[\overline{\Theta}(\underline{z},z,w)\big]^\delta;
\]

\item[$\bullet$] after replacing $w$ by $\Theta \big( z, \underline{ z},
\underline{ w} \big)$ in the extrinsic complexification $\Phi^{ e_c}$, the
result is an identically zero power series in $\C\{ z, \underline{
z}, \underline{ w} \}^d$, namely:
\[
0 \equiv
\sum_{\alpha,\,\beta,\,\gamma,\,\delta}\,
\Phi_{\alpha,\beta,\gamma,\delta}\, z^\alpha\,
\big[\Theta(z,\underline{z},\underline{w})\big]^\beta\,
\underline{z}^\gamma\,\underline{w}^\delta.
\qed
\]

\end{itemize}\smallskip

\end{Lemma}

\subsection{Tangency equations for the determination of
$\mathfrak{ aut}_{ CR}(M)$} In order to compute $\mathfrak{ hol} ( M)$
for an explicitly given CR-generic submanifold $M \subset \C^{n+d}$,
it is most convenient to work with the extrinsic complexification of
its complex defining equations:
\begin{equation}
\label{complex-equations-M}
w_j
=
\Theta_j(z,\underline{z},\underline{w})
\ \ \ \ \ \ \ \ \ \ \ \ \
{\scriptstyle{(j\,=\,1\,\cdots\,d)}}.
\end{equation}
We will assume that $M$ is {\sl rigid}, in the
sense that its real defining equations:
\[
v_j
=
\varphi_j(x,y)
\ \ \ \ \ \ \ \ \ \ \ \ \
{\scriptstyle{(j\,=\,1\,\cdots\,d)}}
\]
do {\em not} depend upon the variables $u = (u_1, \dots, u_d)$.  Two
justifications of this simplification are: 1) the explicit
computations presented in Section~\ref{infinitesimal-model} concern
our cubic model $M_{\sf c}^5 \subset \C^4$, and this model is rigid;
2) the presentation of general formulas for the determination of
infinitesimal CR automorphisms with non necessarily rigid CR-generic
real analytic
$M \subset \C^{ n+d}$ has already been made in~\cite{ AMS-2011}.

So let $M$ be rigid and write
the defining equations of its extrinsic complexification $M^{ e_c}$ as:
\[
w_j
=
\underline{w}_j
+
2i\,\Phi_j\big(z,\underline{z}\big)
\ \ \ \ \ \ \ \ \ \ \ \ \
{\scriptstyle{(j\,=\,1\,\cdots\,d)}},
\]
with the slight change of notation:
\[
\Phi_j\big(z,\overline{z}\big)
\equiv
\varphi_j(x,y).
\]
The extrinsic complexification:
\[
\!\!\!\!\!\!\!\!\!\!\!\!\!\!\!\!\!\!\!\!
{\sf X}
+
\underline{\sf X}
=
\sum_{k=1}^n\,
Z^k(z,w)\,
\frac{\partial}{\partial z_k}
+
\sum_{j=1}^d\,W^j(z,w)\,
\frac{\partial}{\partial w_j}
+
\sum_{k=1}^n\,
\overline{Z}^k\big(\underline{z},\underline{w}\big)\,
\frac{\partial}{\partial\underline{z}_k}
+
\sum_{l=1}^d\,
\overline{W}^l\big(\underline{z},\underline{w}\big)\,
\frac{\partial}{\partial\underline{w}_l}
\]
of (twice) the real part $X + \overline{ X}$ of a $(1, 0)$ vector
field having holomorphic coefficients:
\[
{\sf X}
=
\sum_{k=1}^n\,
Z^k(z,w)\,
\frac{\partial}{\partial z_k}
+
\sum_{j=1}^d\,W^j(z,w)\,
\frac{\partial}{\partial w_j}
\]
is tangent to $M^{ e_c}$ if and only if
it annihilates its equations {\em identically on restriction to them}
(an application of the uniqueness
Lemma~\ref{uniqueness-lemma} is required to pass
from $M$ to $M^{ e_c}$):
\[
0
\equiv
\big(
{\sf X}+\underline{\sf X}
\big)
\Big[
w_j
-
\underline{w}_j
-
2i\,\Phi_j\big(z,\underline{z}\big)
\Big]
\Big\vert_{M^{e_c}}
\ \ \ \ \ \ \ \ \ \ \ \
{\scriptstyle{(j\,=\,1\,\cdots\,d)}}.
\]
Since restricting to $M^{ e_c}$ simply means replacing
$w$ by $\underline{ w} + 2i\, \Phi ( z, \underline{ z})$, these
equations write out in greater length:
\[
\aligned
0
&
\equiv
\bigg[
W^j(z,w)
-
2i\,\sum_{k=1}^n\,
Z^k(z,w)\,
\frac{\partial\Phi_j}{\partial z_k}
\big(z,\underline{z}\big)
-
\\
&
\ \ \ \ \ \ \ \ \ \
-
\overline{W}^j\big(\underline{z},\underline{w}\big)
-
2i\,\sum_{k=1}^n\,
\overline{Z}^k\big(\underline{z},\underline{w}\big)\,
\frac{\partial\Phi_j}{\partial\underline{z}_k}
\big(z,\underline{z}\big)
\bigg]
\bigg\vert_{w=\underline{w}+2i\,\Phi(z,\underline{z})}
\ \ \ \ \ \ \ \ \ \ \ \
{\scriptstyle{(j\,=\,1\,\cdots\,d)}},
\endaligned
\]
that is to say, after really performing the
mentioned replacement of $w$:
\[
\aligned
0
&
\equiv
W^j\big(z,\underline{w}+2i\,\Phi(z,\underline{z})\big)
-
2i\,\sum_{k=1}^n\,
Z^k\big(z,\underline{w}+2i\,\Phi(z,\underline{z})\big)\,
\frac{\partial\Phi_j}{\partial z_k}
\big(z,\underline{z}\big)
-
\\
&
\ \ \ \ \
-\,
\overline{W}^j\big(\underline{z},\underline{w}\big)
-
2i\,
\sum_{k=1}^n\,
\overline{Z}^k\big(\underline{z},\underline{w}\big)\,
\frac{\partial\Phi_j}{\partial\underline{z}_k}
\big(z,\underline{z}\big)
\ \ \ \ \ \ \ \ \ \ \ \ \ \ \ \ \ \ \ \ \ \ \ \ \
{\scriptstyle{(j\,=\,1\,\cdots\,d)}}.
\endaligned
\]
Now, if we expand in partial Taylor series with respect to
$z$ and to $\underline{z}$ all the coefficients of
${\sf X}$ and $\underline{\sf X}$:
\[
\left[
\aligned
Z^k(z,w)
&
=
\sum_{\alpha\in\N^n}\,
z^\alpha\,Z^{k,\alpha}(w),
\\
W^j(z,w)
&
=
\sum_{\alpha\in\N^n}\,
z^\alpha\,W^{j,\alpha}(w),
\endaligned\right.
\ \ \ \ \ \ \ \ \ \ \ \ \ \ \ \ \ \ \ \ \ \ \ \ \ \ \ \ \ \ \ \
\left[
\aligned
\overline{Z}^k\big(\underline{z},\underline{w}\big)
&
=
\sum_{\beta\in\N^n}\,
\underline{z}^\beta\,
\overline{Z}^{k,\beta}\big(\underline{w}\big),
\\
\overline{W}^j\big(\underline{z},\underline{w}\big)
&
=
\sum_{\beta\in\N^n}\,
\underline{z}^\beta\,
\overline{W}^{j,\beta}
\big(\underline{w}\big),
\endaligned\right.
\]
the obtained equations rewrite under the form:
\[
\aligned
0
&
\equiv
\sum_{\alpha\in\N^n}\,z^\alpha\,
W^{j,\alpha}
\big(\underline{w}+2i\,\Phi(z,\underline{z})\big)
-
\\
&
\ \ \ \ \
-\,2i\,
\sum_{k=1}^n\,\sum_{\alpha\in\N^n}\,
z^\alpha\,
Z^{k,\alpha}
\big(\underline{w}+2i\,\Phi(z,\underline{z})\big)\,
\frac{\partial\Phi_j}{\partial z_k}
\big(z,\underline{z}\big)
-
\\
&
\ \ \ \ \
-\,
\sum_{\beta\in\N^n}\,
\underline{z}^\beta\,
\overline{W}^{j,\beta}\big(\underline{w}\big)
-
\\
&
\ \ \ \ \
-\,2i\,
\sum_{k=1}^n\,\sum_{\beta\in\N^n}\,
\underline{z}^\beta\,
\overline{Z}^{k,\beta}
\big(\underline{w}\big)\,
\frac{\partial\Phi_j}{\partial\underline{z}_k}
\big(z,\underline{z}\big)
\ \ \ \ \ \ \ \ \ \ \ \ \
{\scriptstyle{(j\,=\,1\,\cdots\,d)}}.
\endaligned
\]
Next, applying the general Taylor expansion:
\[
\Lambda\big(\underline{w}+2i\,\Phi(z,\underline{z})\big)
=
\sum_{\gamma\in\N^d}\,
\frac{\partial^{\vert\gamma\vert}\Lambda}{
\partial\underline{w}^\gamma}\,
\frac{1}{\gamma!}\,
\big(2i\,\Phi(z,\underline{z})\big)^\gamma,
\]
we continue to expand as follows our $d$ equations:
\[
\aligned
0
&
\equiv
\sum_{\alpha\in\N^n}\,\sum_{\gamma\in\N^d}\,
\frac{1}{\gamma!}\,z^\alpha\,
\big(2i\,\Phi(z,\underline{z})\big)^\gamma\,
\frac{\partial^{\vert\gamma\vert}W^{j,\alpha}}{
\partial\underline{w}^\gamma}
\big(\underline{w}\big)
-
\\
&
\ \ \ \ \
-\,2i\,
\sum_{k=1}^n\,\sum_{\alpha\in\N^n}\,\sum_{\gamma\in\N^d}\,
\frac{1}{\gamma!}\,z^\alpha\,
\big(2i\,\Phi(z,\underline{z})\big)^\gamma\,
\frac{\partial\Phi_j}{\partial z_k}
\big(z,\underline{z}\big)\,
\frac{\partial^{\vert\gamma\vert}Z^{k,\alpha}}{
\partial\underline{w}^\gamma}
\big(\underline{w}\big)
-
\\
&
\ \ \ \ \
-\,
\sum_{\beta\in\N^n}\,
\underline{z}^\beta\,
\overline{W}^{j,\beta}
\big(\underline{w}\big)
-
\\
&
\ \ \ \ \
-\,2i\,
\sum_{k=1}^n\,\sum_{\beta\in\N^n}\,
\underline{z}^\beta\,
\frac{\partial\Phi_j}{\partial\underline{z}_k}
\big(z,\underline{z}\big)\,
\overline{Z}^{k,\beta}\big(\underline{w}\big)
\ \ \ \ \ \ \ \ \ \ \ \ \ \ \
{\scriptstyle{(j\,=\,1\,\cdots\,d)}}.
\endaligned
\]
After performing
a complete expansion in powers of $\big( z, \overline{ z}\big)$,
we obtain a certain family of linear expressions in the
partial derivatives of the unknown sub-coefficients
$Z^{k,\alpha} (w)$, $W^{j, \alpha} (w)$
of the
infinitesimal CR automorphism ${\sf X}$\,\,---\,\,and
of their complexifications
as well
$\overline{ Z}^{k, \beta} \big(
\underline{ w} \big)$,
$\overline{ W}^{j, \beta} \big( \underline{ w} \big)$\,\,---:
\[
\aligned
0
&
\equiv
\sum_{\alpha\in\N^n}\,\sum_{\beta\in\N^n}\,
z^\alpha\underline{z}^\beta
\cdot
\\
&
\cdot
\mathmotsf{Linear-Expression}_{j,\alpha,\beta}
\bigg(
\frac{\partial^{\vert\gamma'\vert}W^{j',\alpha'}}{
\partial\underline{w}^{\gamma'}}
\big(\underline{w}\big),\ \
\frac{\partial^{\vert\gamma'\vert}Z^{k',\alpha'}}{
\partial\underline{w}^{\gamma'}}
\big(\underline{w}\big),\ \
\overline{W}^{j',\alpha'}\big(\underline{w}\big),\ \
\overline{Z}^{k',\beta'}\big(\underline{w}\big)
\bigg)
\\
&
\ \ \ \ \ \ \ \ \ \ \ \ \ \ \ \ \ \ \ \ \ \ \ \ \ \ \ \ \ \ \ \
\ \ \ \ \ \ \ \ \ \ \ \ \ \ \ \ \ \ \ \ \ \ \ \ \ \ \ \ \ \ \ \
{\scriptstyle{(j\,=\,1\,\cdots\,d)}},
\endaligned
\]
the coefficients of these linear expressions being certain functions
of $\underline{ w}$ which depend exclusively upon the defining
equation of $M^{ e_c}$; more advanced information on the latter
coefficients is provided in~\cite{ AMS-2011} in the general (non-rigid)
case.

In summary, by identifying to zero the coefficients of all the
monomials $z^\alpha \overline{ z}^\beta$ in the above identically
satisfied $d$ equations, we have proved that ${\sf X}$ is an
infinitesimal CR automorphism of a rigid real analytic CR-generic
submanifold if and only some (infinite) system
of linear partial differential equations is satisfied:
\[
\aligned
0
&
=
\mathmotsf{Linear-Expression}_{j,\alpha,\beta}
\bigg(
\frac{\partial^{\vert\gamma'\vert}W^{j',\alpha'}}{
\partial\underline{w}^{\gamma'}}
\big(\underline{w}\big),\ \
\frac{\partial^{\vert\gamma'\vert}Z^{k',\alpha'}}{
\partial\underline{w}^{\gamma'}}
\big(\underline{w}\big),\ \
\overline{W}^{j',\alpha'}\big(\underline{w}\big),\ \
\overline{Z}^{k',\beta'}\big(\underline{w}\big)
\bigg)
\\
&
\ \ \ \ \ \ \ \ \ \ \ \ \ \ \ \ \ \ \ \ \ \ \ \ \ \ \ \ \ \ \ \
\ \ \ \ \ \ \ \ \ \ \ \ \ \ \
{\scriptstyle{(j\,=\,1\,\cdots\,d\,;\,\,\,
\alpha\,\in\,\N^n\,;\,\,\,
\beta\,\in\,\N^n)}}.
\endaligned
\]
by an infinite number of unknown functions of $\underline{ w}$.
In Section~\ref{infinitesimal-model},
we shall illustrate these general considerations by providing the full
details of the computation of $\mathfrak{ aut}_{ CR} \big( M_{\sf c}^5
\big)$ for our model $5$-cubic.

It is worth noting that, jointed with Amir Hashemi and Benyamin M.-Alizadeh, recently we have prepared two articles concerning the computations of the Lie algebras of infinitesimal CR-automorphisms. In \cite{SHAM}, we have employed applicable techniques of {\it differential algebra} to provide an effective algorithm to treat systematically solving the {\sc pde} systems arising among the computations. Moreover in \cite{SHAM-2}, it is provided a powerful and fast algorithm to perform computations in the case of (parametric and non-parametric) homogeneous and weighted homogeneous CR-manifolds\,\,---\,\,as Beloshapka's models are. This algorithm employs just simple techniques of linear algebra instead of constructing and solving the mentioned {\sc pde} systems ({\it see} section \ref{infinitesimal-model}). We also have implemented the designed algorithm by means of the effective tools of the new and modern concept {\it comprehensive Gr\"obner systems} for considering the parametric cases.

\section{
Geometric and analytic invariants of
\\
CR-generic submanifolds $M \subset \C^{n+d}$}
\label{geometric-invariants}
\HEAD{\ref{geometric-invariants}.~
Geometric and analytic invariants of
CR-generic submanifolds $M \subset \C^{n+d}$}{
Masoud {\sc Sabzevari} (Shahrekord) and Jo\"el {\sc Merker} (LM-Orsay)}

\subsection{Essential holomorphic dimension and Levi multitype}
Assume again that the CR-generic real analytic
submanifold $M \subset \C^{ n+d}$ is
connected and let as before its extrinsic complexification $M^{ e_c}$
be represented by $d$ holomorphic defining equations of the form:
\[
w_j
=
\Theta_j\big(z,\underline{z},\underline{w}\big)
\ \ \ \ \ \ \ \ \ \ \ \ \
{\scriptstyle{(j\,=\,1\,\cdots\,d)}}.
\]
For any integer $\kappa \in \N$, let us introduce, the
{\sl morphism of
$\kappa$-th jets} of the holomorphic functions:
\[
z
\longmapsto
\Theta_j\big(z,\underline{z},\underline{w})
\]
with respect to only its $n$ first
variables $(z_1, \dots, z_n)$, that is to say
precisely, the map:
\[
\mathmotsf{Segre-jet}_\kappa
\colon
\big(z,\underline{z},\underline{w}\big)
\longmapsto
\bigg(
z,\,
\Big(
{\textstyle{\frac{1}{\beta!}}}\,
\partial_{z}^{\beta}
\Theta_{j}
(z,\underline{z},\underline{w})
\Big)_{1\leqslant j\leqslant d,
\
\beta\in \N^n,\
\vert\beta\vert\leqslant
\kappa}
\bigg)
\in
\C^{n+d\,\binom{n+\kappa}{n}},
\]
whose target components just collect all the partial derivatives of
order $\leqslant \kappa$ of the $\Theta_j$ with respect to the $z_k$,
adding the point $z$ as its first $n$ components.

It is established in~\cite{ Merker-2002, Merker-2005a, Merker-2005b} that
the rank properties of this map\,\,---\,\,called there the morphism of
$\kappa$-jets of complexified Segre varieties\,\,---\,\,are independent of the
choice of the coordinates $(z, w)$, hence also independent of the
complex graphing functions $\Theta_j$, and for this reason, this map
gives access to {\em invariant properties} of the CR-generic
submanifold $M$, namely to properties that are invariant under
biholomorphisms of $\C^{ n+d}$.
Thus, we will not dwell on the invariant character
of this map and just admit it here.

Clearly, the ranks at the origin $0$
and the generic ranks increase with $\kappa$:
\[
\aligned
\mathmotsf{rank}_0
\big(\mathmotsf{Segre-jet}_{\kappa+1}\big)
&
\geqslant
\mathmotsf{rank}_0
\big(\mathmotsf{Segre-jet}_\kappa\big)
\\
\mathmotsf{genrank}
\big(\mathmotsf{Segre-jet}_{\kappa+1}
&
\geqslant
\mathmotsf{genrank}
\big(\mathmotsf{Segre-jet}_\kappa\big).
\endaligned
\]
Concretely, generic ranks are tested by examining all minors of the
Jacobian matrix of this map. But before entering
examination of minors, we
mention the following elementary fact.

\begin{Lemma}
{\rm (\cite{ Merker-2002, Merker-2005a, Merker-2005b})}
If this generic rank does not increase stepwise at a certain jet level
$\kappa^*$:
\[
\mathmotsf{genrank}
\big(\text{\small\sf Segre-jet}_{\kappa^*+1}\big)
=
\mathmotsf{genrank}
\big(\text{\small\sf Segre-jet}_{\kappa^*}\big),
\]
then it remains constantly stabilized for any
jet level $\kappa \geqslant \kappa^*$:
\[
\mathmotsf{genrank}
\big(\text{\small\sf Segre-jet}_\kappa\big)
=
\mathmotsf{genrank}
\big(\text{\small\sf Segre-jet}_{\kappa^*}\big).
\qed
\]
\end{Lemma}

Furthermore, the connectedness of $M$
and the invariance of the Segre jet map imply
({\em ibidem}) that
$\kappa^*$ is the same at every point of $M$ and
in every system of coordinates.

\begin{Definition}
Accordingly, denote now by $\kappa_M$ the smallest integer $\kappa$
such that the generic rank of
the Segre jet map does not
increase after $\kappa$, namely:
\[
\mathmotsf{genrank}
\big(\text{\small\sf Segre-jet}_{\kappa_M-1}\big)
<
\mathmotsf{genrank}
\big(\text{\small\sf Segre-jet}_{\kappa_M}\big)
=
\mathmotsf{genrank}
\big(\text{\small\sf Segre-jet}_{\kappa_M+1}\big)
=
\cdots.
\]
\end{Definition}

Since generic ranks increase (strictly) from $\kappa = 0$, $1$,
$2$, \dots up to $\kappa = \kappa_M$, this integer
is always bounded above by (a better bound follows
below):
\[
\kappa_M
\leqslant
2n+d.
\]

Before explaining the geometric meaning of this maximal generic rank,
let us make the following simple observation. Since the map:
\[
\C^d\ni
\underline{w}
\longmapsto
\Theta\big(z,\underline{z},\underline{w}\big)
\in\C^d
\]
is already of rank $d$ at every point $w \in \C^d$ near the origin
because of~\thetag{ \ref{reality-Theta}}, it follows immediately that
the (generic) rank of the zero-th order Segre-jet map satisfies
already:
\[
\aligned
n+d
&
=
\mathmotsf{genrk}
\big[
(z,\underline{z},\underline{w})
\longmapsto
\big(z,\Theta(z,\underline{z},\underline{w})\big)
\big]
\\
&
=
\mathmotsf{genrk}
\big[
\mathmotsf{Segre-jet}_0
\big],
\endaligned
\]
and it follows {\em passim} that:
\[
\kappa_M
\leqslant
n.
\]

\begin{Definition}
Accordingly, decompose as:
\[
n+d+n_M
\]
the maximal possible generic rank of the
Segre jet maps, with a certain (nonnegative) integer:
\[
n_M
\leqslant
n.
\]
that is to say, the generic rank of $\text{\rm
Segre-jet}_{ \kappa_M}$.
\end{Definition}

Most importantly, the following crucial
statement shows that it is natural to call the integer:
\[
n_M+d
\]
the {\sl essential holomorphic dimension} of $M$. Indeed, let $\Delta
:= \{ \zeta \in \C \colon \vert \zeta \vert < 1 \}$ denote the unit
disc in $\C$, considered as a (small) open piece of it.

\begin{Theorem}
{\rm (\cite{ Merker-2002, Merker-2005a, Merker-2005b})}
Locally in a neighborhood of a Zariski-generic point $p \in M$, the
CR-generic submanifold $M \subset \C^{ n+d}$ is biholomorphically
equivalent to the product:
\[
\underline{M}_p
\times
\Delta^{n-n_M}
\]
of a certain CR-generic submanifold $\underline{ M }_p$ of
codimension $d$ in $\C^{ n_M + d}$ by a complex polydisc $\Delta^{ n -
n_M}$. In addition, no such $\underline{ M}_p$ is biholomorphic to a
product of $\Delta$ by a CR-generic submanifold in a complex Euclidean
space of smaller dimension.
\qed
\end{Theorem}

So in this precise sense, the integer:
\[
n_M+d
\leqslant
n+d
\]
is the
smallest possible integer such that, locally around a generic point,
$M$ looks like a CR-generic submanifold of $\C^{ n_M + d}$, modulo an
innocuous factor $\C^{ n - n_M}$. Even more precisely, this means that
{\em the defining functions of $\underline{ M}_p$ are completely
independent of the variables of the $\C^{ n - n_M}$}, so that
$\underline{ M}_p$ essentially lives in the space of smaller dimension
$\C^{ n_M + d}$ (in some appropriately chosen holomorphic coordinates).

\smallskip

To finish, we review a more precise combinatorics of the generic ranks of the
Segre-jet maps.

\begin{Theorem}
\label{Sigma-Segre}
{\rm (\cite{ Merker-2002, Merker-2005a, Merker-2005b})}
Let $M \subset \C^{ n + d }$ be a {\rm connected} real analytic
CR-generic submanifold of codimension $d \geqslant 1$ and of CR
dimension $n \geqslant 1$. Then there exist well defined integers:
\[
\ell_M
\geqslant 0,
\ \ \ \ \
\lambda_M^0
\geqslant 1,
\ \ \ \ \
\lambda_M^1
\geqslant 1, \ \
\dots\dots,\ \
\lambda_M^{\ell_M}
\geqslant 1
\]
and there exists a proper real analytic subset:
\[
\Sigma_{\smallmathmotsf{Segre}}
\]
of $M$ such
that for every point $p \in M \backslash
\Sigma_{\smallmathmotsf{Segre}}$ and for every system of
coordinates $(z,w)$ vanishing at $p$ in which $M$ is represented by
defining equations in the standard form:
\[
w_j
=
\Theta_j
\big(z,\overline{z},\overline{w}\big)
\ \ \ \ \ \ \ \ \ \ \ \ \
{\scriptstyle{(j\,=\,1\,\cdots\,d)}},
\]
then the following three properties hold:

\begin{itemize}

\smallskip\item[$\bullet$]
$\ell_{ M} \leqslant n_M$;

\smallskip\item[$\bullet$]
$\lambda_M^0 = d$;

\smallskip\item[$\bullet$]
for every $\kappa = 0, 1, \dots, \ell_{ M}$, the mapping of $k$-th
order jets of the $\kappa$-th Segre jet map:
\[
\big(z,\underline{z},\underline{w}\big)
\longmapsto
\bigg(
z,\,
\Big(
{\textstyle{\frac{1}{\beta!}}}\,
\partial_{z}^{\beta}
\Theta_{j}
(z,\underline{z},\underline{w})
\Big)_{1\leqslant j\leqslant d,
\
\beta\in \N^n,\
\vert\beta\vert\leqslant
\kappa}
\bigg),
\]
is equal to:
\[
n
+
\lambda_M^0
+
\lambda_M^1
+\cdots+
\lambda_M^\kappa
\]
at the origin $(0,0)$.

\smallskip\item[$\bullet$]
the essential holomorphic dimension $n_M$ of $M$ is equal to:
\[
n_M
=
d
+
\lambda_M^1
+\cdots+
\lambda_M^{\ell_{M}};
\]

\end{itemize}

\end{Theorem}

\subsection{Generic constancy of CR-geometric invariants}
Recalling that our main objects of study are {\em completely arbitrary}
connected real analytic $d$-codimensional CR-generic submanifolds $M
\subset \C^{ n+d}$, with $d \geqslant 1$ and $n \geqslant 1$,
the two fundamental, preliminary, Theorems~\ref{Sigma-CR}
and~\ref{Sigma-Segre} have shown that, by avoiding certain two
exceptional proper real analytic subsets\,\,---\,\,which might be
complicated\,\,---:
\[
\Sigma_{CR}
\ \ \ \ \ \ \ \ \ \
\text{\rm and}
\ \ \ \ \ \ \ \ \ \
\Sigma_{\smallmathmotsf{Segre}},
\]
namely by re-localizing the study only at points:
\[
p
\in
M
\big\backslash
\big(
\Sigma_{CR}
\cup
\Sigma_{\smallmathmotsf{Segre}}
\big),
\]
one comes to:

\medskip\noindent
{\small\sf I: constancy and maximality
of the Tanaka symbol Lie bracket structure;}

\medskip\noindent
{\small\sf II: constancy and maximality
of the stepwise ranks of the Segre-jet at a
Zariski-generic point,}

\medskip\noindent
where by `maximality', we of course mean
complete non-holonomy and maximal possible essential holomorphic
dimension $n_M = n$.

\subsection{The specificity of CR dimension $1$}
Remarkably, there is only {\em one} variable $z$ in CR dimension $n =
1$. Then we leave as an exercise to the reader to verify that the
fulfillement of Condition~II with $n_M = n$ is just equivalent to the
requirement that in some $d$ real equations for $M$ centered at a
point $p_0 \in M \backslash \Sigma_{\smallmathmotsf{Segre}}$:
\[
v_j
=
\varphi_j(x,y,u)
\ \ \ \ \ \ \ \ \ \ \ \ \
{\scriptstyle{(j\,=\,1\,\cdots\,d)}},
\]
for at least one index $j_0$, there is a {\em nonzero} quadratic
monomial $c_{j_0}\, z\overline{z}$ in the Taylor series of $\varphi_{
j_0}$.  Hence in CR dimension $n = 1$, Condition~II is almost
automatically
satisfied, while a truly infinite algebraic complexity happens to come
from the first Condition~I.

\section{Ineffective access to the local Lie group structure}
\label{ineffective-access}
\HEAD{\ref{ineffective-access}.~Ineffective access
to the local Lie group structure}{
Masoud {\sc Sabzevari} (Shahrekord) and Jo\"el {\sc Merker} (LM-Orsay)}

A fundamental theorem was established in~\cite{ Gaussier-Merker-2004},
Theorem~4.1, but because it does not solve either the biholomorphic
equivalence problem or the effective the description of CR
automorphism groups, we want to briefly comment on its defects after
restating it.  A detailed definition of the concept of {\sl local Lie
group} based on Sophus Lie' original presentation appears in~\cite{
Merker-Engel-Lie, Merker-Hermann}.

\begin{Theorem}
Let $M \subset \C^{ n + d}$ be a connected real analytic CR-generic
submanifold of positive codimension $d \geqslant 1$ and of positive CR
dimension $d \geqslant 1$ and assume it to be maximally non-holonomic
and that its essential holomorphic dimension is the ambient one $n +
d$. Then at every point:
\[
p
\in
M
\big\backslash
\big(
\Sigma_{CR}
\cup
\Sigma_{\smallmathmotsf{Segre}}
\big)
\]
in a neighborhood of which all basic CR-geometric invariants behave
constantly, the abstract group of local biholomorphic transformations:
\[
\mathmotsf{Hol}(M)
\]
fixing $M$\,\,---\,\,but not necessarily $p$\,\,---,
is a finite-dimensional Lie group of dimension bounded above by:
\[
n\,\frac{(2nd+5n)!}{n!\,\,(2nd+4n)!}.
\qed
\]
\end{Theorem}

Such a bound is directly related to jet determination
of local biholomorphisms fixing $M$, and then
at least two defects exist.

\medskip\noindent$\square$\,
This bound on the jet order, or quite
similarly, some other more recent refined bounds,
are usually much above what is truly required\,\,---\,\,a price
to pay for generality\,\,---, as the finer study of specific
CR structures shows, {\em confer exempli gratia} the
geometry-preserving deformations of the above $5$-dimensional
cubic model $M^5 \subset \C^4$.

\medskip\noindent$\square$\,
Most generally, among jets less than say an optimal jet order bound
for the determination of local biholomorphisms fixing $M$, not all
jets below the bound are free, but many are dependent in terms of a
specific set of independent jets, {\em confer} Section~1 in~\cite{
Merker-2008}, and {\em confer} also well known features of
(differential) Gr\"obner bases. Hence determination results
by only one jet order, and the accompanying
bound about the dimension of ${\sf Hol} (M)$, appear to be a rather
rough approach of the reality.

\section{Algebra of infinitesimal CR automorphisms
\\
of the cubic model $M_{\sf c}^5 \subset \C^4$}
\label{infinitesimal-model}
\HEAD{\ref{infinitesimal-model}.~Algebra of infinitesimal CR automorphisms
of the cubic model $M_{\sf c}^5 \subset \C^4$}{
Masoud {\sc Sabzevari} (Shahrekord) and Jo\"el {\sc Merker} (LM-Orsay)}

\subsection{Tangency equations}
\label{Tangency}
Consider the extrinsic complexification of the cubic five-dimensional
model CR-generic submanifold $M_{\sf c}^5 \subset \mathbb{ C}^4$
defined in coordinates $(z, w_1, w_2, w_3, \underline{ z}, \underline{
w}_1, \underline{ w}_2, \underline{ w}_3)$ by the three holomorphic
equations:
\[
\left[
\aligned
w_1-\underline{w}_1 & = 2iz\underline{z},
\\
w_2-\underline{w}_2 & = 2iz\underline{z}(z+\underline{z}),
\\
w_3-\underline{w}_3 & = 2z\underline{z}(z-\underline{z}).
\endaligned\right.
\]
A general $(1, 0)$ holomorphic vector field:
\[
{\sf X}
=
Z(z,w)\,
\frac{\partial}{\partial z}
+
W^1(z,w)\,
\frac{\partial}{\partial w_2}
+
W^2(z,w)\,
\frac{\partial}{\partial w_2}
+
W^3(z,w)\,
\frac{\partial}{\partial w_3}
\]
is an infinitesimal CR automorphism of the cubic model
if and only if
its local holomorphic
coefficients $Z$, $W^1$,
$W^2$, $W^3$ and their conjugates
$\overline{ Z}$, $\overline{ W}^1$,
$\overline{ W}^2$, $\overline{ W}^3$
satisfy the following three equations:
\begin{eqnarray}
\label{first-tangency} && 0 \equiv \big[ W^1 - \overline{W}^1 -
2i\underline{z}Z - 2iz\overline{Z}
\big]_{w=\underline{w}+2i\,\Phi(z,\underline{z})},
\\
\label{second-tangency} && 0 \equiv \big[W^2 - \overline{W}^2 -
4iz\underline{z}Z - 2i\underline{z}^2Z - 2iz^2\overline{Z} -
4iz\underline{z}\overline{Z}
\big]_{w=\underline{w}+2i\,\Phi(z,\underline{z})},
\\
\label{third-tangency} && 0 \equiv \big[W^3 - \overline{W}^3 -
4z\underline{z}Z - 2z^2\overline{Z} + 2\underline{z}^2Z +
4z\underline{z}\overline{Z}
\big]_{w=\underline{w}+2i\,\Phi(z,\underline{z})},
\end{eqnarray}
identically in $\mathbb{ C} \{ z, \underline{ z}, \underline{ w} \}$.
Since these coefficient functions are
analytic, we may expand them with respect to the powers of $z
\in \C$:
\[
Z(z,w) = \sum_{k\in\mathbb{N}}\,z^k\,Z_k(w_1,w_2,w_3) \ \ \ \ \ \ \ \
\text{\rm and} \ \ \ \ \ \ \ \ W^i(z,w) =
\sum_{k\in\mathbb{N}}\,z^k\,W_k^i(w_1,w_2,w_3).
\]
Our current aim is to find closed polynomial expressions for these
holomorphic functions $Z(z,w)$,
$W^1 ( z, w)$, $W^2 ( z, w )$, $W^3 ( z, w)$ by analyzing
this system of three identically satisfied equations.

\begin{Lemma}
\label{lem-mth-paert}
The Taylor expansions with respect to $z$ are relatively polynomial:
\[
\aligned
Z(z,w) & = Z_0(w)+z\,Z_1(w)+z^2\,Z_2(w)+z^3\,Z_3(w),
\\
W^1(z,w) & = W_0^1(w)+z\,W_1^1(w),
\\
W^2(z,w) & = W_0^2(w)+z\,W_1^2(w)+z^2\,W_2^2(w),
\\
W^3(z,w) & = W_0^3(w)+z\,W_1^3(w)+z^2\,W_2^3(w).
\endaligned
\]
\end{Lemma}
\proof After expansion, the first equation reads:
\[
\aligned 0 \equiv & \sum_{k\in\mathbb{N}}\, z^k \Big[ W_k^1 \big(
\underline{w}_1+2iz\underline{z},\,
\underline{w}_2+2iz^2\underline{z}+2iz\underline{z}^2,\,
\underline{w}_3+2z^2\underline{z}-2z\underline{z}^2 \big) -
\\
& \ \ \ \ \ \ \ \ \ \ - 2i\underline{z}\,Z_k \big(
\underline{w}_1+2iz\underline{z},\,
\underline{w}_2+2iz^2\underline{z}+2iz\underline{z}^2,\,
\underline{w}_3+2z^2\underline{z}-2z\underline{z}^2 \big) \Big] +
\\
& + \sum_{k\in\mathbb{N}}\, \underline{z}^k \Big[ -\overline{W}_k^1
\big( \underline{w}_1,\underline{w}_2,\underline{w}_3 \big) - 2iz\,
\overline{Z}_k \big( \underline{w}_1,\underline{w}_2,\underline{w}_3
\big) \Big].
\endaligned
\]
We may expand further each $W_k^1$ and each $Z_k$ using simply
the Taylor series formula for a general holomorphic $\Lambda = \Lambda (
\underline{ w}_1, \underline{ w}_2, \underline{ w}_3)$:
\begin{equation}
\label{taylor} \aligned & \ \ \ \ \ \
\Lambda\big(
\underline{w}_1+2iz\underline{z},\,
\underline{w}_2+2iz^2\underline{z}+2iz\underline{z}^2,\,
\underline{w}_3+2z^2\underline{z}-2z\underline{z}^2 \big) =
\\
 &
\sum_{l_1,\,l_2,\,l_3\in\mathbb{N}}\,
\Lambda_{\underline{w}_1^{l_1}\underline{w}_2^{l_2}\underline{w}_3^{l_3}}\big(
\underline{w}_1,\underline{w}_2,\underline{w}_3 \big)
\frac{(2iz\underline{z})^{l_1}}{l_1!}\,
\frac{(2iz^2\underline{z}+2iz\underline{z}^2)^{l_2}}{l_2!}\,
\frac{(2z^2\underline{z}-2z\underline{z}^2)^{l_3}}{l_3!}\,.
\endaligned
\end{equation}
Chasing then the coefficient of $\underline{ z}^k$ for every $k
\geqslant 2$ in the equation obtained after such an (unwritten)
expansion, we see that the first two lines give absolutely no
contribution, and that from the third line, it only comes: $0 \equiv
\overline{ W}_k^1 ( \underline{ w})$, which is what was claimed about
the $W_k^1$: all of them vanish identically for every $k \geqslant 2$.

Next, chasing the coefficient of $z \underline{ z}^{ k'}$ for every
$k' \geqslant 4$, we get $0 \equiv \overline{ Z}_{ k'} ( \underline{
  w})$, which is what was claimed about the $Z_k ( w)$.

The two remaining families of vanishing equations $0 \equiv
\overline{W}_k^2 ( \underline{ w}) \equiv \overline{ W}_k^3 (
\underline{ w})$ for $k \geqslant 3$ are obtained in a completely
similar way by looking at the second tangency equation~\thetag{
\ref{second-tangency}} and as well at the third~\thetag{
\ref{third-tangency}}.
\endproof

Granted this remarkable relative polynomialness, our next aim is to determine
the expressions of the twelve holomorphic functions:
\[
\aligned
&
Z_0(w), \ \ \ \ \ Z_1(w),\ \ \ \ \ Z_2(w),\ \ \ \ \ Z_3(w),
\\
&
W_0^1(w),\ \ \ \ \ W_1^1(w),
\\
&
W_0^2(w),\ \ \ \ \ W_1^2(w),\ \ \ \ \ W_2^2(w),
\\
&
W_0^3(w),\ \ \ \ \ W_1^3(w),\ \ \ \ \ W_2^3(w)
\endaligned
\]
of only the variables $(w_1, w_2, w_3)$.

To begin with, let us replace the just obtained expressions of the four
functions $Z$, $W^1$ ,$W^2$, $W^3$ in the fundamental equations
\thetag{\ref{first-tangency}}, \thetag{\ref{second-tangency}} and
\thetag{\ref{third-tangency}}:
\begin{equation}
\label{newfirst-tangency}
\footnotesize\aligned 0  \equiv & W_0^1
\big( \underline{w}_1+2iz\underline{z},\,\,
\underline{w}_2+2iz^2\underline{z}+2iz\underline{z}^2,\,\,
\underline{w}_3+2z^2\underline{z}-2z\underline{z}^2 \big) +
\\
&+ z\,W_1^1 \big( \underline{w}_1+2iz\underline{z},\,\,
\underline{w}_2+2iz^2\underline{z}+2iz\underline{z}^2,\,\,
\underline{w}_3+2z^2\underline{z}-2z\underline{z}^2 \big) -
\\
& - \overline{W}_0^1
(\underline{w}_1,\underline{w}_2,\underline{w}_3) - \underline{z}\,
\overline{W}_1^1 (\underline{w}_1,\underline{w}_2,\underline{w}_3) -
\\
& -2i\underline{z}\,Z_0 \big(
\underline{w}_1+2iz\underline{z},\,\,
\underline{w}_2+2iz^2\underline{z}+2iz\underline{z}^2,\,\,
\underline{w}_3+2z^2\underline{z}-2z\underline{z}^2 \big) -
\\
& -2i\underline{z}z\,Z_1 \big(
\underline{w}_1+2iz\underline{z},\,\,
\underline{w}_2+2iz^2\underline{z}+2iz\underline{z}^2,\,\,
\underline{w}_3+2z^2\underline{z}-2z\underline{z}^2 \big) -
\\
& -2i\underline{z}z^2\,Z_2 \big(
\underline{w}_1+2iz\underline{z},\,\,
\underline{w}_2+2iz^2\underline{z}+2iz\underline{z}^2,\,\,
\underline{w}_3+2z^2\underline{z}-2z\underline{z}^2 \big) -
\\
& -2i\underline{z}z^3 Z_3 \big(
\underline{w}_1+2iz\underline{z},\,\,
\underline{w}_2+2iz^2\underline{z}+2iz\underline{z}^2,\,\,
\underline{w}_3+2z^2\underline{z}-2z\underline{z}^2 \big)-
\\
& -2iz\, \overline{Z}_0
(\underline{w}_1,\underline{w}_2,\underline{w}_3) - 2iz\underline{z}\,
\overline{Z}_1 (\underline{w}_1,\underline{w}_2,\underline{w}_3) -
2iz\underline{z}^2\, \overline{Z}_2
(\underline{w}_1,\underline{w}_2,\underline{w}_3)-
\\
& -2iz\underline{z}^3\, \overline{Z}_3
(\underline{w}_1,\underline{w}_2,\underline{w}_3),
\endaligned
\end{equation}
\begin{equation}
\label{newsecond-tangency}
\footnotesize\aligned 0 \equiv & W_0^2 \big( \underline{w}_1+2iz\underline{z},\,\,
\underline{w}_2+2iz^2\underline{z}+2iz\underline{z}^2,\,\,
\underline{w}_3+2z^2\underline{z}-2z\underline{z}^2 \big) +
\\
& + z\,W_1^2 \big( \underline{w}_1+2iz\underline{z},\,\,
\underline{w}_2+2iz^2\underline{z}+2iz\underline{z}^2,\,\,
\underline{w}_3+2z^2\underline{z}-2z\underline{z}^2 \big) +
\\
& + z^2W_2^2\big( \underline{w}_1+2iz\underline{z},\,\,
\underline{w}_2+2iz^2\underline{z}+2iz\underline{z}^2,\,\,
\underline{w}_3+2z^2\underline{z}-2z\underline{z}^2 \big) -
\\
& - \overline{W}_0^2
(\underline{w}_1,\underline{w}_2,\underline{w}_3) - \underline{z}\,
\overline{W}_1^2 (\underline{w}_1,\underline{w}_2,\underline{w}_3) -
\underline{z}^2\overline{W}_2^2(\underline{w}_1,\underline{w}_2,\underline{w}_3)-
\\
& -4i\underline{z}z\,Z_0 \big(
\underline{w}_1+2iz\underline{z},\,\,
\underline{w}_2+2iz^2\underline{z}+2iz\underline{z}^2,\,\,
\underline{w}_3+2z^2\underline{z}-2z\underline{z}^2 \big) -
\\
& -4i\underline{z}z^2\,Z_1 \big(
\underline{w}_1+2iz\underline{z},\,\,
\underline{w}_2+2iz^2\underline{z}+2iz\underline{z}^2,\,\,
\underline{w}_3+2z^2\underline{z}-2z\underline{z}^2 \big) -
\\
& -4i\underline{z}z^3\,Z_2 \big(
\underline{w}_1+2iz\underline{z},\,\,
\underline{w}_2+2iz^2\underline{z}+2iz\underline{z}^2,\,\,
\underline{w}_3+2z^2\underline{z}-2z\underline{z}^2 \big) -
\\
&-4i\underline{z}z^4\,Z_3 \big(
\underline{w}_1+2iz\underline{z},\,\,
\underline{w}_2+2iz^2\underline{z}+2iz\underline{z}^2,\,\,
\underline{w}_3+2z^2\underline{z}-2z\underline{z}^2 \big) -
\\
 & -2i\underline{z}^2\,Z_0 \big( \underline{w}_1+2iz\underline{z},\,\,
\underline{w}_2+2iz^2\underline{z}+2iz\underline{z}^2,\,\,
\underline{w}_3+2z^2\underline{z}-2z\underline{z}^2 \big) -
\\
&  -2i\underline{z}^2z\,Z_1 \big(
\underline{w}_1+2iz\underline{z},\,\,
\underline{w}_2+2iz^2\underline{z}+2iz\underline{z}^2,\,\,
\underline{w}_3+2z^2\underline{z}-2z\underline{z}^2 \big) -
\\
& -2i\underline{z}^2z^2\,Z_2 \big(
\underline{w}_1+2iz\underline{z},\,\,
\underline{w}_2+2iz^2\underline{z}+2iz\underline{z}^2,\,\,
\underline{w}_3+2z^2\underline{z}-2z\underline{z}^2 \big) -
\\
& -2i\underline{z}^2z^3\,Z_3 \big(
\underline{w}_1+2iz\underline{z},\,\,
\underline{w}_2+2iz^2\underline{z}+2iz\underline{z}^2,\,\,
\underline{w}_3+2z^2\underline{z}-2z\underline{z}^2 \big) -
\\
 &  -2iz^2\, \overline{Z}_0
(\underline{w}_1,\underline{w}_2,\underline{w}_3) - 2iz^2\underline{z}\,
\overline{Z}_1 (\underline{w}_1,\underline{w}_2,\underline{w}_3) -
2iz^2\underline{z}^2\, \overline{Z}_2
(\underline{w}_1,\underline{w}_2,\underline{w}_3) -
\\
& -2iz^2\underline{z}^3\, \overline{Z}_3
(\underline{w}_1,\underline{w}_2,\underline{w}_3)-4iz\underline{z}\,
\overline{Z}_0 (\underline{w}_1,\underline{w}_2,\underline{w}_3) -
4iz\underline{z}^2\, \overline{Z}_1
(\underline{w}_1,\underline{w}_2,\underline{w}_3) -
\\
& -4iz\underline{z}^3\, \overline{Z}_2
(\underline{w}_1,\underline{w}_2,\underline{w}_3)-4iz\underline{z}^4\,
\overline{Z}_3 (\underline{w}_1,\underline{w}_2,\underline{w}_3),
\endaligned
\end{equation}
\begin{equation}
\label{newthird-tangency}
\footnotesize \aligned 0 \equiv & W_0^3
\big( \underline{w}_1+2iz\underline{z},\,\,
\underline{w}_2+2iz^2\underline{z}+2iz\underline{z}^2,\,\,
\underline{w}_3+2z^2\underline{z}-2z\underline{z}^2 \big) +
\\
& + z\,W_1^3 \big( \underline{w}_1+2iz\underline{z},\,\,
\underline{w}_2+2iz^2\underline{z}+2iz\underline{z}^2,\,\,
\underline{w}_3+2z^2\underline{z}-2z\underline{z}^2 \big) +
\\
& + z^2W_2^3\big( \underline{w}_1+2iz\underline{z},\,\,
\underline{w}_2+2iz^2\underline{z}+2iz\underline{z}^2,\,\,
\underline{w}_3+2z^2\underline{z}-2z\underline{z}^2 \big) -
\\
& - \overline{W}_0^3
(\underline{w}_1,\underline{w}_2,\underline{w}_3) - \underline{z}\,
\overline{W}_1^3 (\underline{w}_1,\underline{w}_2,\underline{w}_3) -
\underline{z}^2\overline{W}_2^3(\underline{w}_1,\underline{w}_2,\underline{w}_3)-
\\
& -4\underline{z}z\,Z_0 \big(
\underline{w}_1+2iz\underline{z},\,\,
\underline{w}_2+2iz^2\underline{z}+2iz\underline{z}^2,\,\,
\underline{w}_3+2z^2\underline{z}-2z\underline{z}^2 \big) -
\\
& -4\underline{z}z^2\,Z_1 \big(
\underline{w}_1+2iz\underline{z},\,\,
\underline{w}_2+2iz^2\underline{z}+2iz\underline{z}^2,\,\,
\underline{w}_3+2z^2\underline{z}-2z\underline{z}^2 \big) -
\\
& -4\underline{z}z^3\,Z_2 \big(
\underline{w}_1+2iz\underline{z},\,\,
\underline{w}_2+2iz^2\underline{z}+2iz\underline{z}^2,\,\,
\underline{w}_3+2z^2\underline{z}-2z\underline{z}^2 \big) -
\endaligned
\end{equation}
\begin{equation*}
\footnotesize\aligned
& -4\underline{z}z^4\,Z_3 \big(
\underline{w}_1+2iz\underline{z},\,\,
\underline{w}_2+2iz^2\underline{z}+2iz\underline{z}^2,\,\,
\underline{w}_3+2z^2\underline{z}-2z\underline{z}^2 \big) -
\\
& -2z^2\,\overline{Z}_0
(\underline{w}_1,\underline{w}_2,\underline{w}_3) -
2\underline{z}z^2\,\overline{Z}_1
(\underline{w}_1,\underline{w}_2,\underline{w}_3) -
2\underline{z}^2z^2\,\overline{Z}_2
(\underline{w}_1,\underline{w}_2,\underline{w}_3) -
\\
& - 2\underline{z}^3z^2\,\overline{Z}_3
(\underline{w}_1,\underline{w}_2,\underline{w}_3)+
\\
& +2\underline{z}^2\,Z_0 \big(
\underline{w}_1+2iz\underline{z},\,\,
\underline{w}_2+2iz^2\underline{z}+2iz\underline{z}^2,\,\,
\underline{w}_3+2z^2\underline{z}-2z\underline{z}^2 \big) +
\\
 & +2\underline{z}^2z\,Z_1\big(
\underline{w}_1+2iz\underline{z},\,\,
\underline{w}_2+2iz^2\underline{z}+2iz\underline{z}^2,\,\,
\underline{w}_3+2z^2\underline{z}-2z\underline{z}^2 \big) +
\\
& +2\underline{z}^2z^2\,Z_2\big(
\underline{w}_1+2iz\underline{z},\,\,
\underline{w}_2+2iz^2\underline{z}+2iz\underline{z}^2,\,\,
\underline{w}_3+2z^2\underline{z}-2z\underline{z}^2 \big) +
\\
&  +2\underline{z}^2z^3\,Z_3\big(
\underline{w}_1+2iz\underline{z},\,\,
\underline{w}_2+2iz^2\underline{z}+2iz\underline{z}^2,\,\,
\underline{w}_3+2z^2\underline{z}-2z\underline{z}^2 \big) +
\\
&  +4z\underline{z}\, \overline{Z}_0
(\underline{w}_1,\underline{w}_2,\underline{w}_3) + 4z\underline{z}^2\,
\overline{Z}_1 (\underline{w}_1,\underline{w}_2,\underline{w}_3) +
4z\underline{z}^3\, \overline{Z}_2
(\underline{w}_1,\underline{w}_2,\underline{w}_3)+
\\
& + 4z\underline{z}^4\, \overline{Z}_3
(\underline{w}_1,\underline{w}_2,\underline{w}_3).
\endaligned
\end{equation*}
Next, applying~\thetag{\ref{taylor}}, we insert the corresponding
Taylor series of the holomorphic functions $Z(z,w)$, $W^1(z,w)$,
$W^2(z,w)$, $W^3(z,w)$ in the above expressions.  Then we extract the
coefficients of the monomials $z^\mu \underline{ z}^\nu$ for small
values of $\mu$ and $\nu$, and these coefficients must all vanish
identically.

First of all, for $(\mu,\nu)=(0,0)$ we get the following expressions
in the cases of~\thetag{\ref{newfirst-tangency}},
\thetag{\ref{newsecond-tangency}},
\thetag{\ref{newthird-tangency}}, respectively:
\begin{eqnarray}
\label{00g} 0 &\equiv&
W^1_0(\underline{w})-\overline{W}^1_0(\underline{w}),
 \\
\label{00h}
 0
 &\equiv& W^2_0(\underline{w})-\overline{W}^2_0(\underline{w}),
 \\
 \label{00k}
 0
 &\equiv& W^3_0(\underline{w})-\overline{W}^3_0(\underline{w}).
\end{eqnarray}
This means that the functions $W^1_0(w),W^2_0(w)$, $W^3_0(w)$ are all
real, {\em i.e.} have real Taylor coefficients.

For $(\mu,\nu)=(1,0)$, equation~\thetag{\ref{newfirst-tangency}} gives the
equality:
\begin{eqnarray}
\label{10g}
0
\equiv W^1_1(\underline{w})-2i\overline{Z}_0(\underline{w}),
\end{eqnarray}
while we have the following quite advantageous vanishing result from
considering the same exponents $(\mu,\nu)=(1,0)$
in~\thetag{\ref{newsecond-tangency}} and
in~\thetag{\ref{newthird-tangency}}:

\begin{Lemma}
\label{W^2_1vanish} The holomorphic functions $W^2_1(w)$ and
$W^3_1(w)$ vanish identically.
\end{Lemma}

Next, for $(\mu,\nu) = (2,0)$ no new useful information comes
from~\thetag{\ref{newfirst-tangency}}, while
from~\thetag{\ref{newsecond-tangency}} and~\thetag{\ref{newthird-tangency}}
we obtain:
\begin{eqnarray}
\label{20h}
0
&\equiv& W^2_2(\underline{w})-2i\overline{Z}_0(\underline{w}),
\\
\label{20k}
0
&\equiv& W^3_2(\underline{w})-2\overline{Z}_0(\underline{w}).
\end{eqnarray}
\noindent
Next, for $(\mu,\nu) = (1,1)$ we obtained from the
three mentioned equations:
\begin{eqnarray}
\label{11g} 0 &\equiv&
W^1_{0\underline{w}_1}(\underline{w})-Z_1(\underline{w})-\overline{Z}_1(\underline{w}),
\\
\label{11h} 0 &\equiv&
W^2_{0\underline{w}_1}(\underline{w})-2Z_0(\underline{w})-2\overline{Z}_0(\underline{w}),
\\
\label{11k} 0 &\equiv&
iW^3_{0\underline{w}_1}(\underline{w})-2Z_0(\underline{w})+2\overline{Z}_0(\underline{w}),
\end{eqnarray}
and for $(\mu,\nu) = (2,1)$ we get:
\begin{eqnarray}
\label{21g} 0 &\equiv&
iW^1_{0\underline{w}_2}(\underline{w})
+
W^1_{0\underline{w}_3}(\underline{w})+iW^1_{1\underline{w}_1}(\underline{w})-iZ_2(\underline{w}),
\\
\label{21h}
0
&\equiv&
iW^2_{0\underline{w}_2}(\underline{w})
+
W^2_{0\underline{w}_3}(\underline{w})-2iZ_1(\underline{w})-i\overline{Z}_1(\underline{w}),
\\
\label{21k}
0
&\equiv&
iW^3_{0\underline{w}_2}(\underline{w})
+
W^3_{0\underline{w}_3}(\underline{w})
-
2Z_1(\underline{w})-\overline{Z}_1(\underline{w}).
\end{eqnarray}

Now let us continue the process for $(\mu,\nu) = (3,1)$. In this case,
\thetag{\ref{newfirst-tangency}}
and~\thetag{\ref{newthird-tangency}} do not
provide any useful information, while
the following equality can be obtained by
inspecting~\thetag{\ref{newsecond-tangency}}:
\begin{eqnarray}
\label{31h} 0
&\equiv&
W^2_{2\underline{w}_1}(\underline{w})-2Z_2(\underline{w}).
\end{eqnarray}

\begin{Lemma}
\label{z2vanish} The holomorphic function $Z_2(w)$ vanishes identically.
\end{Lemma}

\proof For $(\mu,\nu)=(2,2)$, we get two equations
from~\thetag{\ref{newsecond-tangency}} and~\thetag{\ref{newthird-tangency}}:
\begin{eqnarray}
\label{22h}
0
&\equiv&
-W^2_{0{\underline{w}_1}^2}(\underline{w})
+
4Z_{0\underline{w}_1}(\underline{w})-iZ_2(\underline{w})-i\overline{Z}_2(\underline{w}),
\\
\label{22k}
0
&\equiv&
-W^3_{0{\underline{w}_1}^2}(\underline{w})
-
4iZ_{0\underline{w}_1}(\underline{w})
-
\overline{Z}_2(\underline{w})+Z_2(\underline{w}).
\end{eqnarray}
By differentiating once equation~\thetag{\ref{11h}}
with respect to $\underline{w}_1$ and then
replacing the value of $W^2_{0{\underline{w}_1}^2}(\underline{w})$
in~\thetag{\ref{22h}}, we obtain:
\begin{eqnarray}
\label{2211h}
2\big(
Z_{0\underline{w}_1}(\underline{w})
-
\overline{Z}_{0\underline{w}_1}(\underline{w})\big)
\equiv
i\big(Z_2(\underline{w})+\overline{Z}_2(\underline{w})\big).
\end{eqnarray}
One can apply the same line of reasoning to the equations~\thetag{
\ref{22k}} and~\thetag{\ref{11k}}, and obtain, respectively:
\begin{eqnarray}
\label{2211k}
2\big(Z_{0\underline{w}_1}(\underline{w}\big)
+
\overline{Z}_{0\underline{w}_1}(\underline{w}))
\equiv
i\big(\overline{Z}_2(\underline{w})-Z_2(\underline{w})\big).
\end{eqnarray}
Now comparison of~\thetag{\ref{2211h}} and~\thetag{\ref{2211k}}
yields that:
\begin{eqnarray}
\label{2211hk} 0\equiv
i\overline{Z}_2(\underline{w})-2Z_{0\underline{w}_1}(\underline{w}).
\end{eqnarray}
On the other hand, according to~\thetag{\ref{20h}} one can replace
$W^2_{2\underline{w}_1}(\underline{w})$ by
$2i\overline{Z}_{0\underline{w}_1}(\underline{w})$ in~\thetag{\ref{31h}}.
Thus we have:
\begin{eqnarray}
\label{3120h} 0\equiv
Z_2(\underline{w})-i\overline{Z}_{0\underline{w}_1}(\underline{w}).
\end{eqnarray}
Now, comparing~\thetag{\ref{2211hk}} and~\thetag{\ref{3120h}}
immediately yield
$Z_2(w)\equiv0$, as desired.
\endproof

Now, equations~\thetag{\ref{3120h}}, \thetag{\ref{10g}},
\thetag{\ref{20h}}, \thetag{\ref{20k}}, \thetag{\ref{11h}} and
\thetag{\ref{11k}} imply the
identical vanishing of the the following six functions:
\begin{eqnarray}
\label{firsres}
0&\equiv&
Z_{0\underline{w}_1}(\underline{w})
\equiv
W^1_{1\underline{w}_1}(\underline{w})
\equiv
W^2_{2\underline{w}_1}(\underline{w}) \equiv
\\
\nonumber
&\equiv&
W^3_{2\underline{w}_1}(\underline{w})
\equiv
W^2_{0{\underline{w}_1}^2}(\underline{w}) \equiv
W^3_{0{\underline{w}_1}^2}(\underline{w}).
\end{eqnarray}

In particular, the identical vanishing of the holomorphic function
$Z_{0\underline{w_1}}$ implies a third advantageous fact.

\begin{Lemma}
\label{Z3}
The holomorphic function $Z_3(w)$ vanishes, identically.
\end{Lemma}

\proof
It suffices to look at the case $(\mu,\nu) = (1,4)$ in
\thetag{\ref{newsecond-tangency}}:
\[
0\equiv
\zero{Z_{0\underline{w}_1}(\underline{w})}+i\overline{Z}_3(\underline{w}).
\]
\endproof

\smallskip

Taking account of the three
vanishing Lemmas~\ref{W^2_1vanish}, \ref{z2vanish} and
\ref{Z3}, the initial forms in~\ref{lem-mth-paert}
of three of our four functions simplify:
\begin{eqnarray}
\label{impz} Z(z,w)&=&Z_0(w)+zZ_1(w),
\\
\label{w2im2} W^2(z,w) &=& W^2_0(w)+z^2W^2_2(w),
\\
\label{w3im2} W^3(z,w) &=& W^3_0(w)+z^2W^3_2(w).
\end{eqnarray}
Moreover, the equation~\thetag{\ref{21g}} changes into the following
form:
\begin{eqnarray}
\label{new21g} 0\equiv
iW^1_{0\underline{w}_2}(\underline{w})+W^1_{0\underline{w}_3}(\underline{w}).
\end{eqnarray}
Since the function $W^1_{0}(w)$ is real
by~\thetag{\ref{00g}}, this last expression~\thetag{\ref{new21g}}
together with its conjugate
yield that $W^1_{0\underline{w}_2}(\underline{w})$ and
$W^1_{0\underline{w}_3}(\underline{w})$ vanish together, {\it i.e.}
$W^1_0(w)$ is independent of the variables $w_2$ and $w_3$. Using
this fact and differentiating~\thetag{\ref{11g}} once with respect
to $\underline{w}_2$ and $\underline{w}_3$ implies that
$Z_{1\underline{ w}_i}(\underline{ w}) + \overline{
Z}_{1\underline{ w}_i}(\underline{ w})$
vanishes for $i = 2,3$. In other words, $Z_{1\underline{ w}_i}$ is a real
function for $i = 2,3$.

On the other hand, differentiating
equations~\thetag{\ref{21h}} and~\thetag{\ref{21k}} with respect to
$\underline{w}_2$ and to~$\underline{w}_3$
yields that:
\begin{eqnarray}
\label{new21h}
 0&\equiv&
iW^2_{0\underline{w}_2\underline{w}_i}+W^2_{0\underline{w}_3\underline{w}_i}-iZ_{1\underline{w}_i}(\underline{w}),
\\
\label{new21k}
0&\equiv&
iW^3_{0\underline{w}_2\underline{w}_i}+W^3_{0\underline{w}_3\underline{w}_i}-Z_{1\underline{w}_i}(\underline{w})
\ \ \ \ \ \ \ \ \ \ \ \
{\scriptstyle{(i\,=\ 2, \,3)}}.
\end{eqnarray}
But as~\thetag{\ref{00h}}
and~\thetag{\ref{00k}} meant, the two function $W^2_0(\underline{w})$ and
$W^3_0(\underline{w})$ are real. Then according to
equation~\thetag{\ref{new21h}}, we have
$W^2_{0\underline{w}_3\underline{w}_i}(\underline{w})\equiv0$ for
$i=2,3$.

Now, using again~\thetag{\ref{new21h}} for $i=3$ immediately implies
that $Z_{1\underline{w}_3}(\underline{w})\equiv 0$. In a similar way,
equation~\thetag{\ref{new21k}} yields the vanishing of the two
differentiated functions
$W^3_{0\underline{w}_2\underline{w}_i}(\underline{w})$ and
$Z_{1\underline{w}_2}(\underline{w})$. It follows from this fact
together with the vanishing of $Z_{1\underline{w}_i}(\underline{w})$,
$W^2_{0\underline{w}_3\underline{w}_i}(\underline{w})$ and
$W^3_{0\underline{w}_2\underline{w}_i}(\underline{w})$ for
$i=2,3$ that\,\,---\,\,{\em see}~\thetag{\ref{new21h}}
and~\thetag{\ref{new21k}} again\,\,---:
\begin{eqnarray}
\label{secres}
0
&\equiv&
Z_{1\underline{w}_i}(\underline{w})
\equiv
W^k_{0{\underline{w}_i}{\underline{w}_j}}(\underline{w})
\ \ \ \ \ \ \ \ \ \ \ \
{\scriptstyle{(i,\,j,\,k\,=\,2,\,3)}}.
\end{eqnarray}
Furthermore, equation~\thetag{\ref{newfirst-tangency}} gives the
following equality after inspecting $(\mu,\nu) = (1,3)$:
\begin{eqnarray}
\label{13g}
0\equiv Z_{0\underline{w}_2}(\underline{w})
+
iZ_{0\underline{w}_3}(\underline{w}).
\end{eqnarray}

\begin{Lemma}
\label{z0constant}
The two holomorphic functions $Z_0(w)$ and $Z_1(w)$
are constant.
\end{Lemma}

\proof
Inspection of the fundamental equation
\thetag{\ref{newsecond-tangency}} for $(\mu,\nu)=(3,2)$
and for $(\mu, \nu) = (2,3)$ respectively gives:
\begin{eqnarray}
\label{32h}
0
&\equiv&
-2W^2_{0\underline{w}_1\underline{w}_2}(\underline{w})
+2iW^2_{0\underline{w}_1\underline{w}_3}(\underline{w})
+
iW^2_{2\underline{w}_2}-
\\
\nonumber
&-&W^2_{2\underline{w}_3}(\underline{w})+4Z_{0\underline{w}_2}(\underline{w})
-4iZ_{0\underline{w}_3}(\underline{w})
+
4Z_{1\underline{w}_1}(\underline{w}),
\\
\label{23h}
0
&\equiv&
-2W^2_{0\underline{w}_1\underline{w}_2}(\underline{w})-2iW^2_{0\underline{w}_1\underline{w}_3}(\underline{w})
+6Z_{0\underline{w}_2}(\underline{w})+
\\
\nonumber &+&
2iZ_{0\underline{w}_3}(\underline{w})+2Z_{1\underline{w}_1}(\underline{w}).
\end{eqnarray}
Using~\thetag{\ref{20h}} and~\thetag{\ref{11h}}, we can replace
$W^2_{0\underline{w}_1\underline{w}_i}(\underline{w})$ by
$2Z_{0\underline{w}_i}(\underline{w})+2\overline{
Z}_{0\underline{w}_i}(\underline{w})$
and $W^2_{2\underline{w}_i}(\underline{w})$ by
$2i\overline{Z}_{0\underline{w}_i}(\underline{w})$ for $i=2,3$.
Moreover, according to~\thetag{\ref{13g}}, we can replace
$\overline{Z}_{0\underline{w}_3}$ by
$-i\overline{Z}_{0\underline{w}_2}$. These substitutions change
\thetag{\ref{32h}} and~\thetag{\ref{23h}} into:
\begin{eqnarray}
\label{new32h} 0
&\equiv&
\overline{Z}_{0\underline{w}_2}(\underline{w})
-
Z_{1\underline{w}_1}(\underline{w}),
\\
\label{new23h} 0 &\equiv&
2Z_{0\underline{w}_2}(\underline{w})
-
4\overline{Z}_{0\underline{w}_2}(\underline{w})
+
Z_{1\underline{w}_1}(\underline{w}).
\end{eqnarray}
Eliminating the function $Z_{1\underline{w}_1}(\underline{w})$
from these expressions implies that:
\begin{eqnarray}
0\equiv
2Z_{0\underline{w}_2}(\underline{w})-3\overline{Z}_{0\underline{w}_2}(\underline{w}),
\end{eqnarray}
and together with its conjugate, this last
equation yields that the holomorphic function
$Z_{0\underline{w}_2}(\underline{w}_2)$ vanishes identically.

Thanks to this,
\thetag{\ref{13g}} immediately implies that
$Z_{0\underline{w}_3}(w)\equiv 0$. Furthermore, we know also
from~\thetag{\ref{firsres}} that $Z_{0\underline{w}_1}(w)\equiv 0$ and
hence the holomorphic function $Z_0(w)$ is constant.
On the other hand, according to~\thetag{\ref{new32h}} we have
$Z_{1\underline{w}_1}(\underline{w})\equiv 0$, which, together
with~\thetag{\ref{secres}} yields that the holomorphic function $Z_1(w)$
is constant too.
\endproof

According to the above lemma we have the following forms for the
holomorphic functions $Z_0(w)$ and $Z_1(w)$:
\begin{eqnarray}
\label{f1type}
Z_1(w)&:=&{\sf d}+i{\sf r},
\\
\label{f0type} Z_0(w)&:=&{\sf l}_1+i{\sf l}_2,
\end{eqnarray}
in terms of four complex numbers ${\sf d,r},{\sf l}_1$ and ${\sf l}_2$. Now,
\thetag{\ref{20h}} and~\thetag{\ref{20k}} immediately imply that:
\begin{eqnarray}
\label{h2type}
W^2_2(w)&=& 2{\sf l}_2+2i{\sf l}_1,
\\
\label{k2type}
W^3_2(w)&=& 2{\sf l}_1-2i{\sf l}_2.
\end{eqnarray}
Moreover, differentiating once~\thetag{\ref{11h}}
and~\thetag{\ref{11k}} with respect to $\underline{w}_i, i=1,2,3$, yields
that the holomorphic functions $W^k_{0\underline{w}_1\underline{w}_i}$
vanish for $k=2,3$. Hence taking account of~\thetag{\ref{secres}}, we
have:
\begin{eqnarray}
0\equiv W^k_{0\underline{w}_i\underline{w}_j}
\ \ \ \ \
{\scriptstyle{(i,j=\ 1, \ 2, \,3, \, \ \ k=\ 2, \ 3)}}.
\end{eqnarray}
More precisely, the degrees of the functions $W^2_0(w)$ and
$W^3_0(w)$ with respect to the variables $w_i, i=1,2,3$ are
equal to $1$. Hence, according to~\thetag{\ref{11h}},
\thetag{\ref{11k}}, \thetag{\ref{21h}} and~\thetag{\ref{21k}}, we
can write:
\begin{eqnarray}
\label{h0type} W^2_0(w) &:=&
4{\sf l}_1w_1+3{\sf d}w_2-{\sf r}w_3+{\sf s}_1,
\\
\label{k0type}
W^3_0(w)
&:=&
4{\sf l}_2w_1+{\sf r}w_2+3{\sf d}w_3+{\sf s}_2,
\end{eqnarray}
with two complex numbers ${\sf s}_1$ and ${\sf s}_1$. Moreover,
\thetag{\ref{10g}} and~\thetag{\ref{f0type}} help us to realize the
expression of $W^1_1(w)$ as follow:
\begin{eqnarray}
\label{g1type} W^1_1(w)
=
2{\sf l}_2+2i{\sf l}_1.
\end{eqnarray}
Additionally as we saw, the degree of the real holomorphic function
$W^1_0(w)$ with respect to the variable $w_1$ is one
({\em compare}~\thetag{\ref{11g}} with~\thetag{\ref{f1type}}) and also this
function is independent of the variables $w_2$ and $w_3$. Hence we
have the following expression for $W^1_0(w)$:
\begin{eqnarray}
\label{g0type} W^1_0(w):=2{\sf d}w_1+{\sf t},
\end{eqnarray}
for an arbitrary complex number ${\sf t}$.

Now, the above process has determined explicitly the expressions of
the holomorphic functions $Z(z,w)$ and $W^i(z,w), i=1,2,3$. Indeed,
according to the obtained results we have
exactly found {\em seven} real numerical constants:
\[
{\sf d},
\ \ \
{\sf r},
\ \ \
{\sf l}_1,
\ \ \
{\sf l}_2,
\ \ \
{\sf t},
\ \ \
{\sf s}_1,
\ \ \
{\sf s}_2,
\]
which give us seven $\mathbb R$-linearly independent infinitesimal
automorphisms of our model. More precisely, by verification of the
results we could find the expressions of the desired holomorphic
functions as follows:
\begin{eqnarray*}
Z(z,w) &=& Z_0(w) + Z_1(w)z = {\sf l}_1+ i{\sf l}_2 + ({\sf d}+i{\sf
r})z,
\\
W^1(z,w) &=&
W^1_0(w)
+
W^1_1(w)z
=
2{\sf d}w_1
+
{\sf t}
+
2
({\sf l}_2+i{\sf l}_1)z,
\\
W^2(z,w)
&=&
W^2_0(w)
+
W^2_2(w)z^2
=
4{\sf l}_1w_1
+
3{\sf d}w_2
-
{\sf r}w_3
+
{\sf s}_1
+
2({\sf l}_2+i{\sf l}_1)z^2
\\
W^3(z,w)
&=&
W^3_0(w)
+
W^3_2(w)z^2
=
4{\sf l}_2w_1
+
{\sf r}w_2
+
3{\sf d}w_3
+
{\sf s}_2
+
2
({\sf l}_1-i{\sf l}_2)
z^2.
\end{eqnarray*}

Hence we have the following {\em detailed} confirmation of one of Shananina's
computations (\cite{ Shananina-2000}):

\begin{Proposition}
The Lie algebra $\mathfrak{ aut}_{ CR}
(M) = 2\, {\rm Re}\, \mathfrak{hol}(M)$ of the infinitesimal CR
automorphisms of the $5$-dimensional $3$-codimensional CR-generic
model cubic
$M_{\sf c}^5 \subset
\mathbb C^4$ represented by the three real
graphed equations:
\[
\left[
\aligned w_1-\overline{w}_1 & = 2i\,z\overline{z},
\\
w_2-\overline{w}_2 & = 2i\,z\overline{z}(z+\overline{z}),
\\
w_3-\overline{w}_3 & = 2\,z\overline{z}(z-\overline{z}),
\endaligned\right.
\]
is $7$-dimensional and it is generated by the $\mathbb{ R}$-linearly
independent real parts of the
following seven $(1, 0)$ holomorphic vector fields:
\[
\aligned T & :=
\partial_{w_1},
\\
S_1 & :=
\partial_{w_2},
\\
S_2 & :=
\partial_{w_3},
\\
 L_1 & :=
\partial_z
+ (2iz)\,\partial_{w_1} + (2iz^2+4w_1)\,\partial_{w_2} +
2z^2\,\partial_{w_3},
\\
L_2 & := i\,\partial_z + (2z)\,\partial_{w_1} +
(2z^2)\,\partial_{w_2} - (2iz^2-4w_1)\,\partial_{w_3},
\\
D & := z\,\partial_z + 2w_1\,\partial_{w_1} + 3w_2\,\partial_{w_2} +
3w_3\,\partial_{w_3},
\\
R & := iz\,\partial_z - w_3\,\partial_{w_2} + w_2\,\partial_{w_3}.
\endaligned
\]
\end{Proposition}

Computing all commutators between any two
of these seven generators
of $\frak {aut}_{ CR}(M)$ gives the following complete
table of Lie brackets:

\medskip
\begin{center}
\begin{tabular} [t] { c | c c c c c c c }
& $S_2$ & $S_1$ & $T$ & $L_2$ & $L_1$ & $D$ & $R$
\\
\hline $S_2$ & $0$ & $0$ & $0$ & $0$ & $0$ & $3S_2$ & $-S_1$
\\
$S_1$ & $*$ & $0$ & $0$ & $0$ & $0$ & $3S_1$ & $S_2$
\\
$T$ & $*$ & $*$ & $0$ & $4S_2$ & $4S_1$ & $2T$ & $0$
\\
$L_2$ & $*$ & $*$ & $*$ & $0$ & $-4T$ & $L_2$ & $-L_1$
\\
$L_1$ & $*$ & $*$ & $*$ & $*$ & $0$ & $L_1$ & $L_2$
\\
$D$ & $*$ & $*$ & $*$ & $*$ & $*$ & $0$ & $0$
\\
$R$ & $*$ & $*$ & $*$ & $*$ & $*$ & $*$ & $0$.
\end{tabular}
\end{center}

\noindent
Moreover, we would
like to observe that $\frak{aut}_{ CR}(M)$ is a $3$-graded Lie
algebra with nilpotent negative part, in the sense of Tanaka. More
precisely, with the notation:
\[
\mathfrak{g}
:=
\mathfrak{aut}_{ CR}(M),
\]
if we further set:
\[
\aligned
\mathfrak{g}_{-3}
&
:=
\mathmotsf{Span}_\R
\big<S_1,\,S_2\big>,
\\
\mathfrak{g}_{-2}
&
:=
\mathmotsf{Span}_\R
\big<T\big>,
\\
\mathfrak{g}_{-1}
&
:=
\mathmotsf{Span}_\R
\big<L_1,L_2 \big>,
\\
\mathfrak{g}_0
&
:=
\mathmotsf{Span}_\R
\big<D,R\big>,
\endaligned
\]
then one readily checks:
\[
\frak g=\frak g_{-3}\oplus\frak g_{-2}\oplus\frak g_{-1}\oplus\frak
g_{0}.
\]
Furthermore, with the convention that $\frak g_k = \{0\}$ for either
$k\leqslant -4$ or $k \geqslant 2$,
one may then verify the property that:
\[
[\frak g_{k_1},\frak g_{k_2}]\subseteq\frak g_{k_1+k_2},
\]
for any two integers $k_1,k_2\in\mathbb Z$.
Accordingly, the Lie subalgebra:
\[
\frak g_-:=\frak g_{-3}\oplus\frak g_{-2}\oplus\frak g_{-1}
\]
is called the {\sl Levi-Tanaka algebra} of the CR-manifold $M$. This
subalgebra is isomorphic to the Lie algebra $\frak n^4_5$ in Goze's
classification presented on p.~\pageref{n-5-4} above.

\section{Tanaka prolongations}
\label{Tanaka-prolongations}
\HEAD{\ref{Tanaka-prolongations}.~Tanaka Prolongations}{
Masoud {\sc Sabzevari} (Shahrekord) and Jo\"el {\sc Merker} (LM-Orsay)}

In the former section, we computed the Levi-Tanaka algebra associated
to our model cubic CR-manifold $M_{\sf c}^5 \subset\mathbb
C^4$. Generally, Tanaka's prolongation ({\it see} \cite{Tanaka-1970,
Yamaguchi-1993}) of a graded Lie algebra:
\[
\frak g_-
:=
\frak g_{-\mu}\oplus\ldots\oplus\frak g_{-2}\oplus\frak
g_{-1}
\]
is an algebraic procedure to generate a graded Lie algebra $\frak
g=\frak g_-\oplus\frak g_0\oplus\frak g_1\oplus\ldots$ that
determines to a large extent the geometric properties of the CR
structure ({\it see} \cite{AMS-2011} for a computational example).
We want to show that the prolongation of our model Lie algebra
$\mathfrak{ g}_{ -3} \oplus \mathfrak{ g}_{ -2}
\oplus \mathfrak{ g}_{ -1}$ regives
$\mathfrak{ g}_0 = \left< D, R \right>$.

Consider therefore a finite-dimensional graded real Lie algebra indexed by
negative integers:
\[
\mathfrak{g}_-
=
\mathfrak{g}_{-\mu} \oplus \cdots \oplus \mathfrak{g}_{-2} \oplus
\mathfrak{g}_{-1},
\]
satisfying $[ \mathfrak{ g}_{ - l_1}, \, \mathfrak{ g}_{ -l_2}
] \subset \mathfrak{ g}_{ -l_1 - l_2}$ with the convention that
$\mathfrak{ g}_k = 0$ for $k \leqslant - \mu - 1$. Following~\cite{
Tanaka-1970}, $\mathfrak{ g}_-$ will be said to be {\sl of $\mu$-th
kind}. Assume that there is a complex structure $J \colon \mathfrak{
g}_{ -1} \to \mathfrak{ g}_{ -1}$ such that $J^2 = - \mathmotsf{Id}$, whence
$\mathfrak{ g}_{ -1}$ is even-dimensional and bears a natural
structure of a complex vector space. Tanaka's prolongation of
$\mathfrak{ g}_-$ is an algebraic procedure which generates a certain
larger graded Lie algebra:
\[
\mathfrak{g} = \mathfrak{g}_- \oplus \mathfrak{g}_0 \oplus \mathfrak{g}_1 \oplus
\mathfrak{g}_2
\oplus \cdots
\]
in the following way.

By definition, the order-zero component $\mathfrak{ g}_0$ consists of
all linear endomorphisms ${\sf d} \colon \mathfrak{ g}_- \to
\mathfrak{ g}_-$ which preserve gradation: ${\sf d} ( \mathfrak{ g}_k)
\subset \mathfrak{ g}_k$, which respect the complex structure: ${\sf
d} ( J \, {\sf x}) = J {\sf d} ( {\sf x})$ for all ${\sf x} \in
\mathfrak{ g}_{ -1}$ and which are {\em derivations}, namely satisfy
${\sf d} ( [ {\sf y}, \, {\sf z}]) = [ {\sf d} ( {\sf x}), \, {\sf y}]
+ [{\sf x}, \, {\sf d} ( {\sf y}) ]$ for every ${\sf y}, {\sf z} \in
\mathfrak{ g}_-$. Then the bracket between a ${\sf d} \in \mathfrak{
g}_0$ and an ${\sf x} \in \mathfrak{ g}_-$ is simply defined by $[
{\sf d}, \, {\sf x}] := {\sf d} ( {\sf x})$, while the bracket between
{\em two} elements ${\sf d}', {\sf d}'' \in \mathfrak{ g}_0$ is
defined to be the commutator ${\sf d}'\circ {\sf d}'' - {\sf d}''
\circ {\sf d}'$ between endomorphisms. One checks at once that Jacobi
relations hold, hence $\mathfrak{ g}_- \oplus \mathfrak{ g}_0$ becomes
a true Lie algebra.

By contrast, for any $l \geqslant 1$, no constraint with respect to
$J$ is required. Assuming that the components $\mathfrak{ g}_{ l'}$
are already constructed for any $l' \leqslant l - 1$, the $l$-th
component $\mathfrak{ g}_l$ of the prolongation consists of
$l$-shifted graded linear morphisms $\mathfrak{ g}_- \to \mathfrak{
g}_- \oplus \mathfrak{ g}_0 \oplus \mathfrak{ g}_1 \oplus \cdots
\oplus \mathfrak{ g}_{ l-1}$ that are derivations, namely:
\begin{equation}
\label{Tanakapro}
\mathfrak{g}_l
=
\Big\{ {\sf d} \in
\bigoplus_{k\leqslant-1}\, \mathmotsf{
Lin}(\mathfrak{g}_k,\,\mathfrak{g}_{k+l})
\colon
{\sf d}([{\sf y},\,{\sf z}]) =
[{\sf d}({\sf
y}),\,{\sf z}] + [{\sf y},\,{\sf d}({\sf z})],
\ \ \ \ \
\forall\, {\sf y},\,{\sf
z}\in\mathfrak{g}_- \Big\}.
\end{equation}
Now, for ${\sf d} \in \mathfrak{ g}_k$ and ${\sf e} \in \mathfrak{
g}_l$, by induction on the integer $k + l \geqslant 0$, one defines
the bracket $[ {\sf d}, \, {\sf e} ] \in \mathfrak{ g}_{ k
+ l} \otimes \mathfrak{ g}_-^*$ by:
\begin{equation}
\label{d-e-x} [{\sf d},\,{\sf e}]({\sf x}) = \big[[{\sf d},\,{\sf x}],\,{\sf e}\big]
+ \big[{\sf
d},\,[{\sf e},\,{\sf x}]\big] \ \ \ \ \ \ \ \ \ \ \text{\rm for}\ \ {\sf
x}\in\mathfrak{g}_-.
\end{equation}
One notes that, for $k = l = 0$, this definition coincides with the
above one for $[ \mathfrak{ g}_0, \, \mathfrak{ g}_0]$. It follows by
induction (\cite{ Tanaka-1970, Yamaguchi-1993}) that $[ {\sf d}, \, {\sf e} ]
\in \mathfrak{ g}_{ k+l}$ and that with this bracket, the sum
$\mathfrak{ g}_- \bigoplus_{ k \geqslant 1}\, \mathfrak{ g}_k$ becomes
a graded Lie algebra, because the general Jacobi identity:
\[
0
=
\big[[{\sf d},\,{\sf e}],\,{\sf f}\big]
+
\big[[{\sf f},\,{\sf d}],\,{\sf
e}\big] + \big[[{\sf
e},\,{\sf f}],\,{\sf d}\big]
\]
for ${\sf d} \in \mathfrak{ g}_k$, ${\sf e} \in \mathfrak{ g}_l$ and
${\sf f} \in \mathfrak{ g}_m$ follows by definition when one of $k$,
$l$, $m$ is negative, and can be shown by induction on the integer $k
+ l + m \geqslant 0$ when all of $k$, $l$, $m$ are nonnegative.

\subsection{Tanaka prolongation of the Levi-Tanaka algebra $\frak g_-$}
Now let us find the Tanaka prolongation $\frak g$ of the Lie algebra
$\frak g_-=\frak g_{-3}\oplus\frak g_{-2}\oplus\frak g_{-1}$ with
$\frak g_{-3}=<{\sf s}_1,{\sf s}_2>$, with $\frak g_{-2}=<{\sf
t}>$, with $\frak g_{-1}=<{\sf l}_1,{\sf l}_2>$ and with the same
table of commutators as presented in the former section.\footnote{\,
For reasons of coherence\,\,---\,\,as will be realized at the end of
this section\,\,---, our notations are quite similar to the
ones of the preceding section.
} 
According to the definition, the zero-order component $\frak g_0$ of
this Lie algebra is the subalgebra containing all derivations ${\sf
d}=({\sf d}_1,{\sf d}_2,{\sf d}_3)$ with ${\sf d}_{i}\in\mathmotsf{
End}(\frak g_{-i},\frak g_{-i})$, $i=1,2,3$. Assume the value of the
components of ${\sf d}$ on the basis elements as follows:
\[\aligned
{\sf d}_1({\sf l}_1) &=
 r_1{\sf l}_1
 +
 r_2{\sf l}_2,
 \ \ \ \
 {\sf d}_1({\sf l}_2)
 =
r_3{\sf l}_1 +
 r_4{\sf l}_2, \ \ \ \
 {\sf d}_2({\sf t})
 =
 k{\sf t}
 \\
{\sf d}_3({\sf s}_1) &=
 m_1{\sf s}_1
 +
 m_2{\sf s}_2,
 \ \ \ \ \
 {\sf d}_3({\sf s}_2)
 =
 m_3{\sf s}_1
 +
 m_4{\sf s}_2,
 \endaligned
\]
for some nine real unknown constants. Preserving the complex
structure $J$ by ${\sf d}$ implies that (notice that $J({\sf
l}_1)={\sf l}_2$):
\[
 r_1=r_4, \ \ \ \ \ r_3=-r_2.
 \]
Furthermore, since ${\sf d}$ is a derivation we can obtain some other
relations within the coefficients $r_i,k,m_i$, $i=1,\ldots,4$. At
first, this property gives the following equality:
\[\aligned
{\sf d}([{\sf l}_1,{\sf l}_2]) &=[{\sf d}({\sf l}_1),{\sf l}_2]-
[{\sf d}({\sf l}_2),{\sf l}_1]
\\
&= [r_1{\sf l}_1+r_2{\sf l}_2,{\sf l}_2] - [r_3{\sf l}_1+r_4{\sf
l}_2,{\sf l}_1],
\endaligned
 \]
 which can be read as:
 \[
 k{\sf t} = r_1{\sf t} + r_4{\sf t}
\]
or equivalently:
\[
k = r_1 + r_4.
\]
Applying the same to the value ${\sf d}([{\sf l}_1,{\sf t}])$ gives:
\[
m_3
 =
 r_3,
 \ \ \ \
 m_4
 =
 r_4+k.
 \]
Other values of ${\sf d}$ do not have new result and can be
disregarded. It follows from these relations that $\frak g_0$ is
two-dimensional and is generated by two derivations:
\[
{\sf d}: \ \ \ {\sf l}_1\mapsto {\sf l}_1, \ \ \ {\sf l}_2\mapsto
{\sf l}_2, \ \ \ {\sf t}\mapsto 2{\sf t}, \ \ \ {\sf s}_1\mapsto
3{\sf s}_1, \ \ \ {\sf s}_2\mapsto 3{\sf s}_2,
\]
\[
{\sf r}: \ \ \ {\sf l}_1\mapsto -{\sf l}_2, \ \ \ {\sf l}_2\mapsto
{\sf l}_1, \ \ \ {\sf t}\mapsto 0, \ \ \ {\sf s}_1\mapsto -{\sf
s}_2, \ \ \ {\sf s}_2\mapsto {\sf s}_1.
\]
Next component $\frak g_1$ of $\frak g$ is the set of all liner maps
${\sf d}=({\sf d}_1,{\sf d}_2,{\sf d}_3)$ with ${\sf d}_i\in\mathmotsf{
Lin}\big(\frak g_{-i},\frak g_{1-i}\big)$, $i=1,2,3$ satisfying the
fundamental equation introduced in (\ref{Tanakapro}). One can write
the values of the map $\sf d$ as follows:
\[\aligned
{\sf d}_1({\sf l}_1)&=r_1{\sf d}+r_2{\sf r}, \ \ \ {\sf d}_1({\sf
l}_2)=r_3{\sf d}+r_4{\sf r}, \ \ \ {\sf d}_2({\sf t})=k_1{\sf
l}_1+k_2{\sf l}_2,
\\
{\sf d}_3({\sf s}_1)&=m_1{\sf t}, \ \ \ {\sf d}_3({\sf
s}_2)=m_2{\sf t}.
\endaligned
\]
 Applying the equality
(\ref{Tanakapro}) for ${\sf y}={\sf l}_1$ and ${\sf z}={\sf l}_2$
gives:
\[
{\sf d}([{\sf l}_1,{\sf l}_2]) = [r_1{\sf d}+r_2{\sf r},{\sf l}_2] -
[r_3{\sf d}+r_4{\sf r},{\sf l}_1],
\]
which can be read as:
\[
 k_1{\sf l}_1
 +
 k_2{\sf l}_2
 =
 r_1{\sf l}_2
 +
 r_2{\sf l}_1
 -
 r_3{\sf l}_1
 +
 r_4{\sf l}_2
 \]
 which immediately implies that:
\[
 k_1
 =
 r_2-r_3,
 \ \ \ \
 k_2
 =
 r_1+r_4.
\]
Similar inspecting of the other values ${\sf d}\big([{\sf y},{\sf
z}]\big)$ for ${\sf y,z}={\sf l}_1,{\sf l}_2,{\sf t},{\sf s}_1$
and ${\sf s}_2$ gives, moreover, the following relations between the
coefficients (here we present only the useful equalities):
\[
 k_1 = r_2 - r_3, \ \ k_2 = r_1 + r_4,
 \ \
m_1 = k_2 + 2r _ 1, \ \ r_2 = 0,
 \ \
m_1 + 3r_1 = 0,
 \ \
\]
\[
r_1 = 0,
 \ \
m_2 + r_2 = 0, \ \ m_2 = 2r_3 - k_1,
 \ \
r_3 = 0,
 \ \
m_1 - r_4 = 0, \ \ \ldots.
\]
It follows from the above relations that $r_i=k_j=m_t=0$ for
$i=1,\ldots 4$ and $j,t=1,2$ which means that the subalgebra $\frak
g_1$ is trivial, {\it i.e.} $\frak g_1=0$. Accordingly, the transitivity of Tanaka algebras (\cite{Tanaka-1970, Yamaguchi-1993}) and also fundamentality of the Levi-Tanaka subalgebra $\frak g_-$ (\cite{Beloshapka-2004}) imply that all components $\frak g_l$, $l>0$ vanish identically ({\it see} \cite{SHAM-2}, Proposition 4.6). So we
have:
\[
\frak g=\frak g_{-3}\oplus\frak g_{-2}\oplus\frak g_{-1}\oplus\frak
g_0.
\]
This prolonged Lie algebra is seven-dimensional with the basis
elements ${\sf d},{\sf r},{\sf l}_1,{\sf l}_2,{\sf t},{\sf s}_1$ and
${\sf s}_2$ and according to the definition it has the following
table of commutators:
\medskip
\begin{center}
\begin{tabular} [t] { c | c c c c c c c }
& ${\sf s}_2$ & ${\sf s}_1$ & ${\sf t}$ & ${\sf l}_2$ & ${\sf l}_1$
& ${\sf d}$ & ${\sf r}$
\\
\hline ${\sf s}_2$ & $0$ & $0$ & $0$ & $0$ & $0$ & $3{\sf s}_2$ &
$-{\sf s}_1$
\\
${\sf s}_1$ & $*$ & $0$ & $0$ & $0$ & $0$ & $3{\sf s}_1$ & ${\sf
s}_2$
\\
${\sf t}$ & $*$ & $*$ & $0$ & $4{\sf s}_2$ & $4{\sf s}_1$ & $2{\sf
t}$ & $0$
\\
${\sf l}_2$ & $*$ & $*$ & $*$ & $0$ & $-4{\sf t}$ & ${\sf l}_2$ &
$-{\sf l}_1$
\\
${\sf l}_1$ & $*$ & $*$ & $*$ & $*$ & $0$ & ${\sf l}_1$ & ${\sf l}_2$
\\
${\sf d}$ & $*$ & $*$ & $*$ & $*$ & $*$ & $0$ & $0$
\\
${\sf r}$ & $*$ & $*$ & $*$ & $*$ & $*$ & $*$ & $0$ %
\end{tabular}
\end{center}
By the correspondence $X\mapsto {\sf x}$, for $X=D,R,\ldots,S_2$, one
easily verifies that this Lie algebra is isomorphic to the algebra of
infinitesimal automorphisms $\frak {aut}(M)$ computed in the former
subsection\footnote{\,
Indeed, it is proved ({\it see} \cite{Beloshapka-2004, Gammel-Kossovskiy-2006, Yamaguchi-1993}) that
the Tanaka prolongation of the Levi-Tanaka algebra associated to $M$
{\it contains} the algebra of infinitesimal automorphisms
$\frak{aut}(M)$.
}. 
That is why we chose the same letters for generators of the Lie
algebra as was done in
Section~\ref{infinitesimal-model}.

\section{Equivalence computations for the cubic model}
\label{equivalence-model}
\HEAD{\ref{equivalence-model}.~Equivalence computations
for the cubic model}{
Masoud {\sc Sabzevari} (Shahrekord) and Jo\"el {\sc Merker} (LM-Orsay)}

\subsection{Initial Lie bracket structure}
For the 3-codimensional {\em cubic model} CR-generic submanifold $M^5_{\sf c}
\subset \C^4$ represented as a graph:
\begin{equation}
\label{def-eq}
\left[
\aligned v_1 & = 2\,i\,z\,\overline{z},
\\
v_2 & = 2\,i\,\big(z^2\,\overline{z}+z\,\overline{z}^2\big),
\\
v_3 & = 2\,\big(z^2\,\overline{z}-z\,\overline{z}^2\big),
\endaligned\right.
\end{equation}
the $(0, 1)$-complex tangent bundle $T^{0,1}{ M}^5_{\sf c}$ is
spanned by the single $(0,1)$-vector field:
\[
\overline{\mathcal{L}} = \frac{\partial}{\partial\underline{ z}} -
2\,i\,z\,\frac{\partial}{\partial \overline{w}_1}
-\big(2\,i\,z^2+4\,i\,z\,\overline{z}\big)\,\frac{\partial}{\partial
\overline{w}_2}
-\big(2\,z^2-4\,z\,\overline{z}\big)\,\frac{\partial}{\partial\overline{
w}_3}.
\]
Here, the expression of $\overline{\mathcal L}$ is presented as a
vector field which lives in a neighborhood of $M^5_{\sf c}$ in
$\mathbb C^4$, while $M^5_{\sf c}$ itself, is a real five-dimensional
hypersurface equipped with the five real coordinates
$x,y,u_1,u_2,u_3$. But, in order to express $\overline{\mathcal L}$
{\it intrinsically}, one must drop $\frac{\partial}{\partial
v_1},\,\frac{\partial}{\partial v_2}$ and $\frac{\partial}{\partial
v_3}$ and also simultaneously replace the $v_i$ by their expressions in
\thetag{\ref{def-eq}} for $i=1,2,3$ in its coefficients.
Then, after expanding $\overline{\mathcal L}$ in real and
imaginary parts:
\[
\aligned \overline{\mathcal L}\big\vert_M
&
= \frac{\partial}{\partial\overline{
z}} -\,i\,z\,\frac{\partial}{\partial {u}_1}
-\big(i\,z^2+2\,i\,z\,\overline{z}\big)\,\frac{\partial}{\partial {u}_2}
-\big(z^2-2\,z\,\overline{z}\big)\,\frac{\partial}{\partial{u}_3},
\endaligned
\]
one gains a result that can now be summarized as follows.

\begin{Proposition}
For the cubic 5-dimensional real algebraic CR-generic
model submanifold $M^5_{\sf c} \subset \C^4$ represented near the origin
as a graph:
\[
\left[
\aligned v_1&:= 2\,i\,z\,\overline{z},
\\
v_2&:= 2\,i\,\big(z^2\,\overline{z}+z\,\overline{z}^2\big),
\\
v_3&:= 2\,\big(z^2\,\overline{z}-z\,\overline{z}^2\big),
\endaligned\right.
\]
in coordinates:
\[
\big(z,w_1,w_2,w_3\big)
=
\big(x+iy,u_1+iv_1,u_2+iv_2,u_3+iv_3\big),
\]
its complex bundles $T^{0,1}M^5_{\sf c}$ and $T^{1,0}M^5_{\sf c}$
are generated by:
\[
\boxed{ \aligned \overline{\mathcal L} & =
\frac{\partial}{\partial\underline{ z}}
-\,i\,z\,\frac{\partial}{\partial {u}_1}
-\big(i\,z^2+2\,i\,z\,\overline{z}\big)\,\frac{\partial}{\partial {u}_2}
-\big(z^2-2\,z\,\overline{z}\big)\,\frac{\partial}{\partial{u}_3}
\ \ \ \ \ \ \ \ \text{\rm and}\
\\
{\mathcal L} & = \frac{\partial}{\partial{ z}}
+\,i\,\overline{z}\,\frac{\partial}{\partial {u}_1}
+\big(i\,\overline{z}^2+2\,i\,z\,\overline{z}\big)\,\frac{\partial}{\partial
{u}_2}
-\big(\overline{z}^2-2\,z\,\overline{z}\big)\,\frac{\partial}{\partial{u}_3}.
\endaligned}
\]
\end{Proposition}

\subsection{Length-two Lie bracket}
Between these two complex vector fields $\mathcal L$ and
$\overline{\mathcal L}$, there is of course only one Lie bracket
$\big[ \mathcal L,\overline{\mathcal L} \big]$ of length two. This
vector field is in fact {\it imaginary}. In order to get a {\it real}
vector field, we multiply it by $i$:
\[
\mathcal T
:=
i\,\big[\mathcal L,\overline{\mathcal L}\big].
\]
A direct computation yields its expression:
\[
\boxed{
\mathcal T = 2\,\frac{\partial}{\partial u_1}
+
4\,\big(z+\overline{z}\big)\frac{\partial}{\partial u_2}
-
4\,i\,\big(z-\overline{z}\big)\frac{\partial}{\partial u_3}.}
\]

\subsection{Length-three Lie brackets}
In this length, we have two Lie brackets:
\[
\mathcal S
:=
\big[
\mathcal L,\mathcal T],
\ \ \ \ \
\overline{\mathcal S}
:=
\big[\,
\overline{\mathcal L},\mathcal T
\big].
\]
Again, direct easy computations provide the following expressions for
them:
\[\boxed{
\aligned \mathcal S &= 4\,\frac{\partial}{\partial
u_2}-4\,i\,\frac{\partial}{\partial u_3}
\ \ \ \ \ \ \ \ \text{\rm and}\
\\
\overline{\mathcal S} &= 4\,\frac{\partial}{\partial
u_2}+4\,i\,\frac{\partial}{\partial u_3}.
\endaligned}
\]

\begin{Lemma}
The five vector fields $\overline{\mathcal L},{\mathcal L},\mathcal
T,\overline{\mathcal S},\mathcal S$ constitute a (complex) frame for
$TM^5_{\sf c}\otimes_{\mathbb R}\mathbb C$.
\end{Lemma}

\proof
It is sufficient to see from their expressions that these
vector fields are linearly independent and hence they constitute a
frame.
\endproof

\subsection{Other iterated Lie brackets}
We saw that the collection of five vector fields:
\[
\Big\{ \overline{\mathcal{S}},\, \mathcal{S},\, \mathcal{T},\,
\overline{\mathcal{L}},\, \mathcal{L} \Big\},
\]
where:
\[
\aligned \mathcal{T} & :=
i\big[\mathcal{L},\overline{\mathcal{L}}\big],
\\
\mathcal{S} & := \big[\mathcal{L},\,\mathcal{T}\big],
\\
\overline{\mathcal{S}} & :=
\big[\overline{\mathcal{L}},\,\mathcal{T}\big],
\endaligned
\]
makes up a frame for $\C \otimes_\R TM_{\sf c}^5$. Having five vector
fields implies that there are in sum $\binom{5}{2} = 10$ Lie brackets
between them. Thus, there remain seven such brackets to be looked
at. However, simple computations show that all of these remaining
vector fields vanish identically, namely we have:
\[
\aligned \big[\mathcal{L},\,\mathcal{S}\big] & = 0, \ \ \ \
\big[\overline{\mathcal{L}},\,\mathcal{S}\big] =0, \ \ \ \
\big[\mathcal{L},\,\overline{\mathcal{S}}\big] =0,
\\
\big[ \overline{\mathcal{L}},\,\overline{\mathcal{S}}\big] & =0, \ \
\ \ \big[\mathcal{T},\,\mathcal{S}\big]=0, \ \ \ \
\big[\mathcal{T},\,\overline{\mathcal{S}}\big]=0,
\\
\big[\mathcal{S},\,\overline{\mathcal{S}}\big]&=0.
\endaligned
\]
When we will study arbitrary geometry-preserving deformations
of this cubic model, the corresponding seven supplementary Lie
brackets will be highly more complicated.

\subsection{Passage to a dual coframe
and its Darboux-Cartan structure.}
On the natural agreement that the coframe:
\[
\big\{du_3,du_2,du_1,\,dz,\,d\overline{z}\big\}
\]
is
dual to the frame:
\[
\big\{ {\textstyle{\frac{\partial}{\partial
u_3}}},\, {\textstyle{\frac{\partial}{\partial u_2}}},\,
{\textstyle{\frac{\partial}{\partial u_1}}},\,
{\textstyle{\frac{\partial}{\partial z}}},\,
{\textstyle{\frac{\partial}{\partial\overline{z}}}} \big\},
\]
let us introduce the coframe:
\[
\big\{\overline{\sigma_0},\,\sigma_0,\,\rho_0,\,\overline{\zeta_0},\,\zeta_0\big\}
\ \ \ \text{\rm which is dual to the frame} \ \ \
\big\{\overline{\mathcal{S}},\,
\mathcal{S},\,\mathcal{T},\,\overline{\mathcal{L}},\,\mathcal{L}\big\},
\]
that is to say which satisfies by definition:
\[
\begin{array}{ccccc}
\overline{\sigma_0}(\overline{\mathcal{S}})=1 \ \ \ & \ \ \
\overline{\sigma_0}(\mathcal{S})=0 \ \ \ & \ \ \
\overline{\sigma_0}(\mathcal{T})=0 \ \ \ & \ \ \
\overline{\sigma_0}(\overline{\mathcal{L}})=0 \ \ \ & \ \ \
\overline{\sigma_0}\big(\mathcal{L}\big)=0,
\\
\sigma_0(\overline{\mathcal{S}})=0 \ \ \ & \ \ \
\sigma_0(\mathcal{S})=1 \ \ \ & \ \ \ \sigma_0(\mathcal{T})=0 \ \ \ &
\ \ \ \sigma_0(\overline{\mathcal{L}})=0 \ \ \ & \ \ \
\sigma_0\big(\mathcal{L}\big)=0,
\\
\rho_0(\overline{\mathcal{S}})=0 \ \ \ & \ \ \ \rho_0(\mathcal{S})=0
\ \ \ & \ \ \ \rho_0(\mathcal{T})=1 \ \ \ & \ \ \
\rho_0(\overline{\mathcal{L}})=0 \ \ \ & \ \ \ \rho_0(\mathcal{L})=0,
\\
\overline{\zeta_0}(\overline{\mathcal{S}})=0 \ \ \ & \ \ \
\overline{\zeta_0}(\mathcal{S})=0 \ \ \ & \ \ \
\overline{\zeta_0}(\mathcal{T})=0 \ \ \ & \ \ \
\overline{\zeta_0}(\overline{\mathcal{L}})=1 \ \ \ & \ \ \
\overline{\zeta_0}\big(\mathcal{L}\big)=0,
\\
\zeta_0(\overline{\mathcal{S}})=0 \ \ \ & \ \ \
\zeta_0(\mathcal{S})=0 \ \ \ & \ \ \ \zeta_0(\mathcal{T})=0 \ \ \ & \
\ \ \zeta_0(\overline{\mathcal{L}})=0 \ \ \ & \ \ \
\zeta_0\big(\mathcal{L}\big)=1.
\end{array}
\]
Since neither $\mathcal{ T}$, nor $\mathcal{ S}$, nor $\overline{
\mathcal{ S}}$ incorporates any $\frac{ \partial}{ \partial u_j}$, $j
= 1, 2, 3$, we have:
\[
\zeta_0 = dz \ \ \ \ \ \ \ \ \ \ \ \ \ \text{\rm and} \ \ \ \ \ \ \ \
\ \ \ \ \ \overline{\zeta_0} = d\overline{z}.
\]
In order to launch the Cartan algorithm of
equivalence for the cubic model, initially we need the
expressions of the five 2-forms
$d\overline{\sigma}_0,d\sigma_0,d\rho_0,d\overline{\zeta}_0,
d\zeta_0$ in terms of the wedge products of
$\overline{\sigma}_0,\sigma_0,\rho_0,\overline{\zeta}_0, \zeta_0$. To
find them, we use the following well known
duality correspondence.

\begin{Lemma}
\label{lem-d}
Given a frame $\big\{ \mathcal{ L}_1, \dots, \mathcal{
L}_n\big\}$ on an open subset of $\R^n$ enjoying the Lie structure:
\[
\big[\mathcal{L}_{i_1},\,\mathcal{L}_{i_2}\big] =
\sum_{k=1}^n\,a_{i_1,i_2}^k\,\mathcal{L}_k \ \ \ \ \ \ \ \ \ \ \ \ \
{\scriptstyle{(1\,\leqslant\,i_1\,<\,i_2\,\leqslant\,n)}},
\]
where the $a_{ i_1, i_2}^k$ are certain functions on $\R^n$, the dual
coframe $\{ \omega^1, \dots, \omega^n \}$ satisfying by definition:
\[
\omega^k
\big(\mathcal{L}_i\big)
=
\delta_i^k
\]
enjoys a quite similar Darboux-Cartan structure, up to an overall minus sign:
\[
d\omega^k = - \sum_{1\leqslant i_1<i_2\leqslant n}\,
a_{i_1,i_2}^k\,\omega^{i_1}\wedge\omega^{i_2}
\ \ \ \ \ \ \ \ \ \ \ \ \ {\scriptstyle{(k\,=\,1\,\cdots\,n)}}.
\]
\end{Lemma}

\proof
Just apply the so-called Cartan formula $d\omega ( \mathcal{
X}, \mathcal{ Y}) = \mathcal{ X} \big( \omega ( \mathcal{ Y})\big) -
\mathcal{ Y} \big( \omega ( \mathcal{ X})\big) - \omega\big( [
\mathcal{ X}, \mathcal{ Y}]\big)$.
\endproof

Thank to this Lemma, minding the overall minus sign, we can readily
find the expressions of the exterior derivatives of our five $1$-forms
that provide the associated Darboux-Cartan structure:
\begin{equation}
\label{d0-model}
\boxed{\aligned d\overline{\sigma}_0 & =
\rho_0\wedge\overline{\zeta}_0,
\ \ \ \ \ \ \ \ \ \ \ \ \ \ \ \ \ \
d\sigma_0
=
\rho_0\wedge\zeta_0,\
\\
d\rho_0 & = i\,\zeta_0\wedge\overline{\zeta_0},
\\
d\overline{\zeta_0} & = 0,
\ \ \ \ \ \ \ \ \ \ \ \ \ \
\ \ \ \ \ \ \ \ \ \ \ \ \ \
d\zeta_0= 0.
\endaligned}
\end{equation}

\subsection{Ambiguity matrix}
Our next goal is to set up, in a coordinate-free manner, the {\sl
Cartan ambiguity matrix} associated to the {\sl problem of local
biholomorphic equivalences}:
\[
h\colon\ \ \ (z,w)
\longmapsto
\big(f(z,w),\,g(z,w)\big)
=:
(z',w')
\]
between our cubic 5-dimensional CR-generic cubic model $M^5_{\sf c}$
and another arbitrary local real analytic CR-generic maximally minimal real
submanifold ${M'}^5 \subset \mathbb{C}^4$ in coordinates $\big( z', w_1',
w_2', w_3' \big)$.
Naturally, we assume that
${M'}^5$ is also equipped with a collection of five vector fields:
\[
\Big\{
\mathcal{L}',\,\,\overline{\mathcal{L}}',\,\, \mathcal T', \,\,
\mathcal S', \,\, \overline{\mathcal S}'
\Big\}
\]
where $\mathcal{L}'$ is a local generator of $T^{1,0} M'$ and where:
\[
\mathcal{T}'
:=
i\,\big[\mathcal{L}',\,\overline{\mathcal{L}}'\big],
\ \ \ \ \ \ \ \
\mathcal{S'}
:=
\big[\mathcal{L}',\,\mathcal{T}'\big],
\ \ \ \ \ \ \ \
\overline{\mathcal{S}}'
:=
\big[\overline{\mathcal{L}}',\,\mathcal{T}'\big],
\]
which so that they also make up a frame for its rank-five complexified
tangent bundle $\C \otimes_\R T{M'}^5$, in the sense that, at every
point $p' \in {M'}^5$ one has:
\[
\aligned
\C\otimes T_{p'}{M'}^5
&=
\C\,\mathcal{L}'\big\vert_{p'} \oplus
\C\,\overline{\mathcal{L}}'\big\vert_{p'}
\oplus\C\,\mathcal{T}'\big\vert_{p'}\oplus\C\,
\mathcal{S}'\big\vert_{p'}\oplus\C\,\overline{\mathcal{L}}'\big\vert_{p'}.
\endaligned
\]
The fact that ${M'}^5$ is maximally minimal in the sense of
Proposition~\ref{Proposition-initial-form} if and only these fields
make up such a frame is essentially tautological, but it will be made
more explicit in Proposition~\ref{frame-5-general} below.

\smallskip

Thus, suppose that two given hypersurfaces $M^5_{\sf c}$ and ${M'}^5$
are CR-equivalent under some (possibly unknown) local
equivalence:
\[
h\colon\ \ \
M^5_{\sf c}
\longrightarrow
{M'}^5
\]
which is a biholomorphism from some neighborhood of $M^5_{\sf c}$ onto
some neighborhood of ${M'}^5$. Then, the associated differential of
$h$:
\[
h_\ast
\colon\ \ \
TM_{\sf c}^5
\longrightarrow
T{M'}^5
\]
induces a push-forward complexified map, still
denoted with the same symbol:
\[
h_\ast
\colon\ \ \
\mathbb{C}\otimes
TM_{\sf c}^5\longrightarrow
\mathbb{C}
\otimes
T{M'}^5,
\]
which naturally defined by ({\it see} \cite{Boggess-1991}, Subsection
3.1):
\[
\aligned
h_\ast\big({\sf z}\otimes_\R \mathcal X\big)
:=
{\sf z}\otimes_\R
h_\ast(\mathcal X), \ \ \ \ \ \ \ {\footnotesize {\sf z}\in \mathbb
C, \ \ \ \mathcal X\in T_pM^5_{\sf c}},\ \ \
p\in M_{\sf c}^5.
\endaligned
\]

\begin{Proposition}
\label{10-ambiguity-matrix}
The initial ambiguity matrix associated to the local biholomorphic
equivalence problem between the cubic $5$-dimensional model CR-generic
submanifold $M^5_{\sf c}$ and any other maximally minimal CR-generic
$5$-dimensional submanifolds ${M'}^5 \subset \C^4$ under local
biholomorphic transformations is of the general form:
\[
\left(\!\!
\begin{array}{ccccc}
{\sf a}\overline{\sf a}\overline{\sf a} & 0 & \overline{\sf c} &
\overline{\sf e} & \overline{\sf d}
\\
0 & {\sf a}{\sf a}\overline{\sf a} & {\sf c} & {\sf d} & {\sf e}
\\
0 & 0 & {\sf a}\overline{\sf a} & \overline{\sf b} & {\sf b}
\\
0 & 0 & 0 & \overline{\sf a} & 0
\\
0 & 0 & 0 & 0 & {\sf a}
\end{array}
\!\!\right),
\]
where ${\sf a}$, ${\sf b}$, ${\sf c}$, ${\sf e}$, ${\sf d}$ are
complex numbers. Moreover, the collection of all these matrices makes
up a real $10$-dimensional matrix Lie subgroup of ${\sf GL}_5 ( \C)$.
\end{Proposition}

\proof
As a CR mapping, $h \colon M_{\sf c}^5
\longrightarrow {M'}^5$ respects the complex structure ({\it see}
\cite{Boggess-1991}, page 149, Definition 1), and hence we necessarily
have:
\[\aligned
h_{\ast}(\mathcal L) & := {\sf a}'\mathcal L',
\endaligned
\]
for some nonzero complex function ${\sf a}'$ defined on ${M'}^5$. By
conjugation, it obviously comes that one also has:
\[
\aligned h_{\ast}\big(\overline{\mathcal L}\big) & := \overline{{\sf
a}}'\,\overline{\mathcal L}'.
\endaligned
\]
Next, let us look at what happens with Lie brackets. Since the
differential commutes with brackets, we have:
\[
\aligned h_{\ast}(\mathcal T)
= h_{\ast}\big(i[\mathcal
L,\overline{\mathcal L}]\big)=
i\,h_\ast\big([\mathcal L,\overline{\mathcal
L}]\big)
=i\,
\big[h_\ast(\mathcal L),h_\ast(\overline{\mathcal L})\big]
&
=i\big[{\sf
a}'\mathcal L',\overline{\sf a}'\overline{\mathcal L}'\big],
\endaligned
\]
and if we expand this last bracket, we obtain:
\[
\aligned
h_{\ast}(\mathcal T)
&
=
{\sf a}'\overline{\sf a}'\cdot i\,
\big[\mathcal{L}',\,\overline{\mathcal{L}'}\big]
\underbrace{-
i\,\overline{\sf a}'\overline{\mathcal{L}}'({\sf a}')}_{
=:\,{\sf b}'}
\cdot
\mathcal{L}'
+
i\,{\sf a}'\,\mathcal{L}'\big(\overline{\sf a}'\big)
\cdot
\overline{\mathcal{L}}'
\\
&=: {\sf a}'\overline{\sf a}'\mathcal T'+{\sf b}'\mathcal
L'+\overline{\sf b}'\overline{\mathcal L}',
\endaligned
\]
by introducing\,\,---\,\,in accordance with the general principles of
Cartan's approach\,\,---\,\,a name ${\sf b}'$ for a certain
complicated function that might remain unknown as long the problem of
equivalence is not settled.

Now, a quite similar computation for the next Lie bracket:
\[
\aligned h_{\ast}(\mathcal S)&= h_{\ast} \big([\mathcal L,\mathcal
T]\big)=\big[h_{\ast}(\mathcal L),h_{\ast}(\mathcal T)\big]=
\\
&=
\big[{\sf a}'\mathcal L',{\sf a}'\overline{\sf a}'\mathcal T'+{\sf
b}'\mathcal L'+\overline{\sf b}'\overline{\mathcal L}'\big]=
\\
&={\sf a'a'\overline{a}'}\mathcal S'+0\mathcal S'+\underbrace{\big({\sf
a}'\mathcal L'\big(\sf a'\overline{a}'\big)-{\it i}\,{\sf
a}'\overline{b}'\big)}_{=:{\sf c}'}\mathcal T'+
\\
& \ \ \ \ \ +\underbrace{\big(-{\sf a'\overline{a}'}\mathcal T'({\sf
a}')+{\sf a}'\overline{\mathcal L}'({\sf b}')-{\sf b}'\mathcal
L'({\sf a}')-\overline{\sf b}'\overline{\mathcal L}'({\sf
a}')\big)}_{=:{\sf e}'}\mathcal L'+
\\
& \ \ \ \ \ +\underbrace{{\sf a}'\mathcal L'\big(\overline{{\sf
b}}'\big)}_{=:{\sf d}'}\overline{\mathcal L}',
\endaligned
\]
shows us that:
\[
h_{\ast}(\mathcal S)= {\sf a'a'\overline{a}'}\mathcal S'+0\mathcal
S'+{\sf c}'\mathcal T'+{\sf e}'\mathcal L'+{\sf d}'\overline{\mathcal
L}',
\]
for some three complex-valued functions ${\sf c}'$, ${\sf d}'$ and
${\sf e}'$ defined on ${M'}^5$. This means that the initial ambiguity
matrix associated to the equivalence problem under local
biholomorphic maps for maximally minimal CR-generic submanifolds ${M'}^5
\subset \C^4$ is of the general form\,\,---\,\,we drop the primes in the group
variables ${\sf a}', {\sf b}', {\sf c}', {\sf d}', {\sf e}'$\,\,---:
\begin{eqnarray*}
\left(%
\begin{array}{c}
\overline{\mathcal S} \\
\mathcal S \\
\mathcal T \\
\overline{\mathcal L} \\
\mathcal L \\
\end{array}%
\right) = \left(\!\!
\begin{array}{ccccc}
{\sf a}\overline{\sf a}\overline{\sf a} & 0 & \overline{\sf c} &
\overline{\sf e} & \overline{\sf d}
\\
0 & {\sf a}{\sf a}\overline{\sf a} & {\sf c} & {\sf d} & {\sf e}
\\
0 & 0 & {\sf a}\overline{\sf a} & \overline{\sf b} & {\sf b}
\\
0 & 0 & 0 & \overline{\sf a} & 0
\\
0 & 0 & 0 & 0 & {\sf a}
\end{array}
\!\!\right)
\left(%
\begin{array}{c}
\overline{\mathcal S}' \\
\mathcal S' \\
\mathcal T' \\
\overline{\mathcal L}' \\
\mathcal L' \\
\end{array}%
\right),
\end{eqnarray*}
as announced. The group property follows by verifying that the
product of any two such matrices is again a matrix of this general
form.
\endproof

\subsection{Setting up the equivalence problem}
According to the general principles (\cite{ Gardner-1989, Olver-1995}),
the so-called {\sl lifted coframe} in terms
of the dual basis of $1$-form then becomes,
after a plain matrix transposition:
\begin{equation}
\label{lift-coframe-model}
\aligned \left(\!\!
\begin{array}{c}
\overline{\sigma}
\\
\sigma
\\
\rho
\\
\overline{\zeta}
\\
\zeta
\end{array}
\!\!\right) := \underbrace{\left(\!\!
\begin{array}{ccccc}
{\sf a}\overline{\sf a}\overline{\sf a} & 0 & 0 & 0 & 0
\\
0 & {\sf a}{\sf a}\overline{\sf a} & 0 & 0 & 0
\\
\overline{\sf c} & {\sf c} & {\sf a}\overline{\sf a} & 0 & 0
\\
\overline{\sf e} & {\sf d} & \overline{\sf b} & \overline{\sf a} & 0
\\
\overline{\sf d} & {\sf e} & {\sf b} & 0 & {\sf a}
\end{array}
\!\!\right)}_{=:g} \left(\!\!
\begin{array}{c}
\overline{\sigma_0}
\\
\sigma_0
\\
\rho_0
\\
\overline{\zeta_0}
\\
\zeta_0
\end{array}
\!\!\right),
\endaligned
\end{equation}
that is to say:
\[
\aligned \overline{\sigma} & = {\sf a}\overline{\sf a}\overline{\sf
a}\, \overline{\sigma_0},
\\
\sigma & = {\sf a}{\sf a}\overline{\sf a}\,\sigma_0,
\\
\rho & = \overline{\sf c}\,\overline{\sigma_0} + {\sf c}\,\sigma_0 +
{\sf a}\overline{\sf a}\,\rho_0,
\\
\overline{\zeta} & = \overline{\sf e}\,\overline{\sigma_0} + {\sf
d}\,\sigma_0 + \overline{b}\,\rho_0 + \overline{\sf
a}\,\overline{\zeta_0},
\\
\zeta & = \overline{\sf d}\,\overline{\sigma_0} + {\sf e}\,\sigma_0 +
{\sf b}\,\rho_0 + {\sf a}\,\zeta_0.
\endaligned
\]
Of course, the $1$-form $\rho$ is real and the $1$-forms
$\overline{ \sigma}$ and $\overline{ \zeta}$ are the conjugate of
$\sigma$ and $\zeta$.

The main aim of the next sections will be to construct one (or more)
absolute parallelism(s) on a certain principal bundle(s) over a
general maximally minimal $M^5 \subset \C^4$ by performing the Cartan
equivalence method with initial data a coframe $\big\{ \overline{
\sigma}, \sigma, \rho, \overline{ \zeta}, \zeta\big\}$ related to
$M^5$\,\,---\,\,to be computed explicitly later\,\,---\,\,and
the structure group:
\[\footnotesize\aligned
G:=\left\{g=\left(\!\!
\begin{array}{ccccc}
{\sf a}\overline{\sf a}\overline{\sf a} & 0 & 0 & 0 & 0
\\
0 & {\sf a}{\sf a}\overline{\sf a} & 0 & 0 & 0
\\
\overline{\sf c} & {\sf c} & {\sf a}\overline{\sf a} & 0 & 0
\\
\overline{\sf e} & {\sf d} & \overline{\sf b} & \overline{\sf a} & 0
\\
\overline{\sf d} & {\sf e} & {\sf b} & 0 & {\sf a}
\end{array}
\!\!\right), \ \ \ \ \ {\sf a,b,c,d,e}\in\mathbb C\right\}.
\endaligned
\]
In this way, we will encounter invariants of the desired holomorphic
equivalence. {\em But our tricky strategy of approach is to perform
this algorithm beforehand in the case where the five $1$-form $\big\{
\overline{ \sigma}, \sigma, \rho, \overline{ \zeta}, \zeta\big\}$
come from the model cubic $M_{\sf c}^5 \subset \C^4$, because the
much simpler computations will serve as a guide before treating the
more delicate general case of an arbitrary maximally minimal $M^5
\subset \C^4$.} We essentially admit that the reader knows the so-called {\sl
Cartan algorithm} which consists of three major parts {\sl
absorbtion, normalization} and {\sl prolongation}, but we briefly
provide a summary
({\em cf.}~\cite{Olver-1995}, pp.~305--310).

\subsection{Absorption and normalization: general features}

Suppose that on an $n$-dimensional (local) manifold
$M$, one has an initial coframe:
\[
\theta_0
=
\big(\theta_0^1,\dots,\theta_0^n \big)^t
\]
of $1$-forms written as a column vector. Suppose that a
certain closed matrix structure
subgroup $G\subset {\sf GL}(n,\mathbb R)$ specifies a certain geometric
equivalence problem, and
introduce the lifted coframe:
\[
\theta=g\cdot\theta_0.
\]
Differentiating both sides
of this equality gives in concise notation:
\begin{equation}
\label{d-g-theta}
d\theta=dg\wedge\theta_0
+
g\cdot d\theta_0.
\end{equation}
Assume that the expressions of $d\theta^1_0,\ldots,d\theta^n_0$ in
terms of the initial 2-forms $\theta^j_0\wedge\theta^k_0$ with $1\leqslant
j<k\leqslant n$, are at hand thanks to some preliminary computations
related to the specific geometric features of the problem under
study, say of the form:
\[
d\theta_0^i
=
\sum_{1\leqslant j<k\leqslant n}\,
T_{0jk}^i
\cdot
\theta_0^j\wedge\theta_0^k
\ \ \ \ \ \ \ \ \ \ \ \ \
{\scriptstyle{(i\,=\,1\,\cdots\,n)}},
\]
for some explicitly known functions $T_{ 0\bullet\bullet}^\bullet$.
For $i = 1, \dots, n$, let $g^{(i)}$ denote the $i$-th row of the
matrix $g$. Passing through the inverted lifting:
\[
\theta_0
=
g^{-1}\cdot\theta,
\]
it is generally
possible\,\,---\,\,after computations that are almost always delicate
(even on a computer) in
nontrivial applications\,\,---\,\,to re-express each scalar expression
$g^{(i)} \cdot d\theta_0$ in terms
of the lifted basis $\theta^j \wedge \theta^k$ of
$2$-forms, and one gets expressions of the quite similar general form:
\[
g^{(i)}\cdot d\theta_0
=
\sum_{1\leqslant j<k\leqslant n}\,
T^i_{jk}\cdot
\theta^j\wedge\theta^k,
\]
in which there appear certain functions $T^\bullet_{\bullet\bullet}$,
called {\sl torsion coefficients}, which express explicitly in terms
of the initial structure functions $T_{0jk}^i$ and in terms of the
group parameters. This treats the second term of the right-hand side
of~\thetag{ \ref{d-g-theta}}.

For the first term, one rewrites:
\[
dg\wedge\theta_0=\underbrace{dg\cdot g^{-1}}_{\omega_{\sf MC}}
\,\wedge\,
\underbrace{g\cdot\theta_0}_{\theta},
\]
and there naturally appears an $n\times n$ matrix:
\[
\aligned
\omega_{\sf MC}
=
&\,
\big(
(\omega_{\sf MC})_j^i
\big)_{1\leqslant j\leqslant n}^{1\leqslant i\leqslant n}
\\
:=
&\,
\sum_{k=1}^n\,
dg_k^i\,
\big(g^{-1}\big)_j^k
\endaligned
\]
(where $i$ is the index of rows)
which is called the {\sl Maurer-Cartan matrix} associated to $G$.

As a result, the expression of $d\theta$ rewritten in terms of the
lifted coframe of $2$-forms $\theta^j \wedge \theta^k$ can be read as
follows:
\begin{equation}
\aligned
\label{equiv-1-general}
d\theta^i
=
\sum_{j=1}^n\,
\big(\omega_{\sf MC}\big)_j^i
\wedge\theta^j
+
\sum_{1\leqslant j<k\leqslant n}\,
T_{jk}^i\cdot
\theta^j\wedge\theta^k
\ \ \ \ \ \ \ \ \ \ \ \ \
{\scriptstyle{(i\,=\,1\,\cdots\,n)}}.
\endaligned
\end{equation}
The obtained equations are called {\sl structure equations}.

To pursue further, it is necessary to express each entry of the
Maurer-Cartan matrix in terms of a basis of $1$-forms living on the
abstract Lie group corresponding to $G$. Thus, with $r := \dim_\R G$,
and with a basis of left-invariant $1$-forms on the group in question,
say $\alpha^1, \dots, \alpha^r$, one can decompose:
\[
\aligned
\big(\omega_{\sf MC}\big)_j^i
=
\sum_{s=1}^r\,
a^i_{js}\,\alpha^s
\ \ \ \ \ \ \ \ \ \ \ \ \
{\scriptstyle{(i,\,\,j\,=\,1\,\dots\,n)}},
\endaligned
\]
in terms of certain {\em constants} $a^\bullet_{\bullet\bullet}$.
Consequently, the equations~\thetag{\ref{equiv-1-general}} can be
brought into the form:
\begin{equation}
\label{equiv-1-general-entries}
\aligned
d\theta^i
=
\sum_{k=1}^n\bigg(\sum_{s=1}^r\,a_{ks}^i\,\alpha^s
+
\sum_{j=1}^{k-1}T^i_{jk}\,\theta^j\bigg)\wedge\theta^k
\ \ \ \ \ \ \ \ \ \ \ \ \
{\scriptstyle{(i\,=\,1\,\cdots\,n)}}.
\endaligned
\end{equation}
Since the constant coefficients $a^\bullet_{ \bullet \bullet}$ depend
merely on the structure Lie group $G$, for another similar lifted
coframe $\widetilde{ \theta} = \widetilde{g} \cdot \widetilde{ \theta
}_0$ with the same group $\widetilde{ G} = G$ but on another manifold
$\widetilde{ M}$ with another initial coframe $\big( \widetilde{
\theta}_0^1, \dots, \widetilde{ \theta}_0^n \big)$, one has in a
completely similar way:
\begin{equation}
\label{equiv-1-tilde} \aligned
d\widetilde{\theta}^i
=
\sum_{k=1}^n\bigg(\sum_{s=1}^r\,a_{ks}^i\,\widetilde{\alpha}^s
+
\sum_{j=1}^{k-1}\widetilde{T}^i_{jk}\,
\widetilde{\theta}^j\bigg)\wedge\widetilde{\theta}^k
\ \ \ \ \ \ \ \ \ \ \ \ \
{\scriptstyle{(i\,=\,1\,\cdots\,n)}},
\endaligned
\end{equation}
with {\em unchanged} constants $a_{ks}^i$. So, if an equivalence
holds
between the two initial coframes,
according to the fundamental result of the
theory (\cite{ Gardner-1989}), this equivalence lifts up and it provides
an equality\,\,---\,\,through an unwritten pull-back\,\,---\,\,
between the two lifted coframes:
\[
\theta=\widetilde{\theta}.
\]
Applying exterior differentiation, it follows at once that:
\[
d\theta
=
d\widetilde{\theta}.
\]
But then we may subtract the two
representation~\thetag{\ref{equiv-1-general-entries}} and
\thetag{\ref{equiv-1-tilde}}
of the $d\theta^i$, and this gives us:
\begin{equation}
\label{theta=tilde}
\aligned
0
&
\equiv
\sum_{k=1}^n
\bigg(
\underbrace{
\sum_{s=1}^r\,
a_{ks}^i\,\big(
\widetilde{\alpha}^s-\alpha^s\big)
+
\sum_{j=1}^{k-1}\,
\big(
\widetilde{T}^i_{jk}-T^i_{jk}\big)\theta^j
}_{=:\eta_k^i}
\bigg)
\wedge\theta^k
\ \ \ \ \ \ \ \ \ \ \ \ \
{\scriptstyle{(i\,=\,1\,\cdots\,n)}}.
\endaligned
\end{equation}

\begin{Lemma}
{\sc (Cartan)}
\label{Cartan-lemma}
Let $\big\{ \vartheta^1,\ldots,\vartheta^n \big\}$ be a set of locally
defined linearly independent 1-forms. Then some $n$ arbitrary
$1$-forms $\eta_1,\ldots,\eta_n$ satisfy $\sum_{k=1}^n\,
\eta_k\wedge\vartheta^k=0$ if and only if they write as $\eta_k =
\sum_{l=1}^n\, A_{kl}\, \vartheta^l$ for some symmetric matrix of
functions $A_{kl} = A_{lk}$.
\qed
\end{Lemma}

The Cartan's Lemma plays an essential role in the theory of
equivalence problems; here, by applying it to the
equality~\thetag{\ref{theta=tilde}}, we obtain that, for each $i = 1,
\ldots, n$, there exist functions $A^i_{kl}$ with $A^i_{ kl} =
A^i_{ lk}$ such that:
\begin{equation}
\label{cartan-lemma-1}
\aligned
\sum_{s=1}^r\,
a_{ks}^i\,
\big(
\widetilde{\alpha}^s-\alpha^s\big)
&
+
\sum_{j=1}^{k-1}\,
\big(
\widetilde{T}^i_{jk}-T^i_{jk}\big)\,
\theta^j
=
\sum_{l=1}^n\,A^i_{kl}\,\theta^l
\\
&
\ \ \ \ \ \ \
{\scriptstyle{(i,\,\,k\,=\,1\,\cdots\,n)}}.
\endaligned
\end{equation}
Next, multiplying both sides of each of these $n$ equations by
$\theta^1 \wedge \cdots \wedge \theta^n$ brings the equalities:
\[
\aligned
\sum_{s=1}^r\,a_{ks}^i\,
&
\big(\widetilde{\alpha}^s-\alpha^s\big)
\wedge\theta^1\wedge\cdots\wedge\theta^n=0
\\
&
\ \ \ \ \ \ \ \ \ \
{\scriptstyle{(i\,=\,1\,\cdots\,n)}}.
\endaligned
\]
If we apply once again Cartan's Lemma to each of these $n$
relations, we conclude that for every $s=1,\ldots,r$, there exist $n$
functions $z^s_1, \dots, z_n^s$ defined on the base manifold $M$ such
that:
\begin{equation*}
\widetilde{\alpha}^s
=
\alpha^s
+
\sum_{j=1}^n\,z^s_j\,\theta^j
\ \ \ \ \ \ \ \ \ \ \ \ \
{\scriptstyle{(s\,=\,1\,\cdots\,r)}}.
\end{equation*}
Substituting this expression into~\thetag{\ref{theta=tilde}} gives:
\[
0\equiv\sum_{k=1}^n
\bigg(
\sum_{s=1}^r\,\sum_{j=1}^n\,
a_{ks}^i\,z_j^s\,\theta^j
+
\sum_{j=1}^{k-1}\,\big(
\widetilde{T}_{jk}^i
-
T_{jk}^i
\big)\,
\theta^j\bigg)
\wedge
\theta^k.
\]
Then all the $\frac{ n ( n-1)}{ 2}$ coefficients of the basis
$2$-forms $\theta^j \wedge \theta^k$ for $1\leqslant j< k\leqslant n$
must vanish. Extracting these coefficients and equating them to zero
yields:
\begin{equation*}
\label{t=tield}
\widetilde{T}_{jk}^i
=
T_{jk}^i
+
\sum_{s=1}^r\,
\big(
a_{js}^i\,z_k^s
-
a_{ks}^i\,z_j^s
\big)
\ \ \ \ \ \ \ \ \ \ \ \ \
{\scriptstyle{(i\,=\,1\,\cdots\,n\,;\,\,\,
1\,\leqslant\,j\,<\,k\,\leqslant\,n)}}.
\end{equation*}

\begin{Proposition}
\label{prop-changes}
In the structure equations
\thetag{\ref{equiv-1-general-entries}}{\em :}
\begin{equation*}
\aligned
d\theta^i
=
\sum_{k=1}^n\bigg(\sum_{s=1}^r\,a_{ks}^i\,\alpha^s
+
\sum_{j=1}^{k-1}T^i_{jk}\,\theta^j\bigg)\wedge\theta^k
\ \ \ \ \ \ \ \ \ \ \ \ \
{\scriptstyle{(i\,=\,1\,\cdots\,n)}},
\endaligned
\end{equation*}
one can replace each Maurer-Cartan form $\alpha^s$ and each torsion
coefficient $T^i_{jk}$ with:
\begin{equation}
\label{possible-changes}
\boxed{\aligned
\alpha^s
&
\longmapsto
\alpha^s+\sum_{j=1}^n\,z^s_j\,\theta^j
\ \ \ \ \ \ \ \ \ \ \ \ \ \ \ \ \ \ \ \ \ \ \ \ \ \ \ \ \
{\scriptstyle{(s\,=\,1\,\cdots\,r)}},
\\
T^i_{jk}&\longmapsto T^i_{jk}
+
\sum_{s=1}^r\,
\big(
a_{js}^i\,z_k^s
-
a_{ks}^i\,z_j^s
\big) \ \ \ \ \ \ \ \ \ \ \
{\scriptstyle{(i\,=\,1\,\cdots\,n\,;\,\,\,
1\,\leqslant\,j\,<\,k\,\leqslant\,n)}},\,
\endaligned}
\end{equation}
for some arbitrary functions $z^\bullet_\bullet$
on the base manifold $M$.
\qed
\end{Proposition}

\subsection{Absorbtion and normalization for the model}
\label{g-minus-1}

Now,
differentiating both sides of~\thetag{\ref{lift-coframe-model}} gives:
\begin{eqnarray*}
d\,\left(%
\begin{array}{c}
\overline{\sigma} \\
\sigma \\
\rho \\
\overline{\zeta} \\
\zeta \\
\end{array}%
\right) =
dg\,
\wedge\,\left(%
\begin{array}{c}
\overline{\sigma}_0 \\
\sigma_0 \\
\rho_0 \\
\overline{\zeta}_0 \\
\zeta_0 \\
\end{array}%
\right) + g\cdot\left(%
\begin{array}{c}
d\overline{\sigma}_0 \\
d\sigma_0 \\
d\rho_0 \\
d\overline{\zeta}_0 \\
d\zeta_0 \\
\end{array}%
\right).
\end{eqnarray*}
Also,
differentiating the entries of the matrix $g$ and putting the
expressions of $d \overline{ \sigma}_0$, $d\sigma_0$, $d\rho_0$,
$d\overline{\zeta}_0$, $d\zeta_0$ in the second term of the above
equation gives the structure equations as follows:
\begin{equation}
\small\aligned
\label{struc-eq-1-model}
d\,\left(%
\begin{array}{c}
\overline{\sigma} \\
\sigma \\
\rho \\
\overline{\zeta} \\
\zeta \\
\end{array}%
\right)&=\underbrace{
\left(%
\begin{array}{ccccc}
\overline{\sf a}^2\,d{\sf a}+2 {\sf a}\overline{\sf a}\,
d\overline{\sf a} & 0 & 0
& 0 & 0 \\
0 & 2{\sf a}\overline{\sf a}\,d{\sf a}
+
{\sf a}^2\overline{\sf a}\,d\overline{\sf
a} & 0 & 0 & 0 \\
 d\overline{\sf c} & d{\sf c} & {\sf a}\,
d\overline{\sf a}+\overline{\sf a}\,d{\sf
a} & 0 & 0 \\
d\overline{\sf e} & d{\sf d} & d\overline{\sf b} & d\overline{\sf a} & 0 \\
d\overline{\sf d} & d{\sf e} & d{\sf b} & 0 & d{\sf a} \\
\end{array}%
\right)}_{dg} \wedge
\left(%
\begin{array}{c}
\overline{\sigma}_0 \\
\sigma_0 \\
\rho_0 \\
\overline{\zeta}_0 \\
\zeta_0 \\
\end{array}%
\right)
+
\\
& \ \ \ \ \
+ \left(\!\!
\begin{array}{c}
{\sf a}\overline{\sf a}\overline{\sf a}\,d\overline{\sigma}_0
\\
{\sf a}{\sf a}\overline{\sf a}\,d\sigma_0
\\
\overline{\sf c}\,d\overline{\sigma}_0 + {\sf c}\,d\sigma_0 + {\sf
a}\overline{\sf a}\,d\rho_0
\\
\overline{\sf e}\,d\overline{\sigma}_0+ {\sf
d}\,d\sigma_0+\overline{\sf b}\,d\rho_0+ \overline{\sf
a}\,d\overline{\zeta}_0
\\
\overline{\sf d}\,d\overline{\sigma}_0+{\sf e}\,d\sigma_0+ {\sf
b}\,d\rho_0+ {\sf a}\,d\zeta_0
\end{array}
\!\!\right).
\endaligned
\end{equation}
Moreover, since
the determinant of the matrix $g$ is $({\sf a}\overline{\sf
a})^5$ and since the variable $\sf a$, lying on the diagonal of the
matrix group, must necessarily be nonzero, $g$ is invertible with
inverse:
\begin{eqnarray*}
g^{-1}=\left(%
\begin{array}{ccccc}
\frac{1}{{\sf a}\overline{\sf a}^2} & 0 & 0 & 0 & 0 \\
0 & \frac{1}{{\sf a}^2\overline{\sf a}} & 0 & 0 & 0 \\
-\frac{\overline{\sf c}}{{\sf a}^2\overline{\sf a}^3} & -\frac{\sf c}{{\sf
a}^3\overline{\sf a}^2} &
\frac{1}{{\sf a}\overline{\sf a}} & 0 & 0 \\
\frac{\overline{\sf b}\,\overline{\sf c}-\overline{\sf e}{\sf a}\overline{\sf
a}}{{\sf a}^2\overline{\sf a}^4} &
\frac{\overline{\sf b}\sf c-{\sf a}\overline{a}\sf d}{{\sf a}^3\overline{\sf a}^3}
& -\frac{\overline{\sf b}}{{\sf a}\overline{\sf a}^2} &
\frac{1}{\overline{\sf a}} & 0 \\
\frac{{\sf b}\overline{\sf c}-{\sf a}\overline{\sf a}\overline{\sf d}}{{\sf
a}^3\overline{\sf a}^3} &
\frac{{\sf bc}-{\sf ea}\overline{\sf a}}{{\sf a}^4\overline{\sf a}^2} & -\frac{\sf
b}{{\sf a}^2\overline{\sf a}} &
0 & \frac{1}{{\sf a}} \\
\end{array}%
\right).
\end{eqnarray*}
Multiplying this matrix by the transpose of
$\big( \overline{\sigma}, \sigma, \rho, \overline{ \zeta}, \zeta
\big)$ gives the inverse expressions:
\begin{equation}
\aligned \sigma_0 & = \frac{1}{{\sf
a}^2\overline{\sf a}}\,\sigma,
\\
\rho_0 & = -\,\frac{\overline{\sf c}}{{\sf a}^2\overline{\sf a}^3}\,
\overline{\sigma} - \frac{\sf c}{{\sf a}^3\overline{\sf a}^2}\,\sigma
+ \frac{1}{{\sf a}\overline{\sf a}}\,\rho,
\\
\zeta_0 & = \frac{{\sf b}\overline{\sf c}-{\sf a}\overline{\sf
a}\overline{\sf d}}{ {\sf a}^3\overline{\sf a}^3}\,\overline{\sigma}
+ \frac{{\sf b}{\sf c}-{\sf a}\overline{\sf a}{\sf e}}{ {\sf
a}^4\overline{\sf a}^2}\,\sigma - \frac{{\sf b}}{{\sf
a}^2\overline{\sf a}}\,\rho + \frac{1}{{\sf a}}\,\zeta,
\endaligned
\end{equation}
with plain conjugations to obtain $\overline{ \sigma}_0$ and
$\overline{ \zeta}_0$\,\,---\,\,remembering that both $\rho_0$ and
$\rho$ are real.

Thus, one can then replace the first term $dg \wedge
(\overline{\sigma}_0, \sigma_0, \rho_0, \overline{ \zeta}_0,
\zeta_0)^t$ of~\thetag{\ref{struc-eq-1-model}} by:
\[
(dg\cdot g^{-1})
\wedge
\big[
g\cdot\big(\overline{\sigma}_0,\sigma_0,\rho_0,\overline{
\zeta}_0,\zeta_0\big)^t\big],
\]
and this, according to
\thetag{\ref{lift-coframe-model}}, gives:
\begin{equation*}
\aligned
dg\wedge\left(
\begin{array}{c}
\overline{\sigma}_0 \\
\sigma_0 \\
\rho_0 \\
\overline{\zeta}_0 \\
\zeta_0 \\
\end{array}
\right)
=
\underbrace{\left(
\begin{array}{ccccc}
2\,\overline{\alpha}_1+\alpha_1 & 0 & 0 & 0 & 0 \\
0 & 2\,\alpha_1+\overline{\alpha}_1 & 0 & 0 & 0 \\
\overline{\alpha}_2 & {\alpha}_2 & \alpha_1+\overline{\alpha}_1 & 0 & 0 \\
\overline{\alpha}_3 & \overline{\alpha}_4 & \overline{\alpha}_5 &
\overline{\alpha}_1 & 0 \\
{\alpha}_4 & {\alpha}_3 & {\alpha}_5 & 0 & {\alpha}_1 \\
\end{array}%
\right)}_{\omega_{\sf MC}:=dg\cdot g^{-1}}\,\,\,
\wedge\!\!\!\!\!
\underbrace{\left(%
\begin{array}{c}
\overline{\sigma} \\
\sigma \\
\rho \\
\overline{\zeta} \\
\zeta \\
\end{array}%
\right),}_{
g\cdot(\overline{\sigma}_0,\sigma_0,\rho_0,\overline{\zeta}_0,\zeta_0)^t}
\endaligned
\end{equation*}
where the {\sl Maurer-Cartan 1-forms}\,\,---\,\,minding the typographical
distinction between the group variable
${\sf d}$ and the standard symbol $d$ for exterior
differentiation\,\,---\,\,are:
\[
\aligned
{\alpha_1} & := \frac{d{\sf a}}{\sf a},
\\
{\alpha_2} &:= \frac{d{\sf c}}{{\sf a}^2\overline{\sf a}} -
\frac{{\sf c}\,d{\sf a}}{{\sf a}^3\overline{\sf a}} - \frac{{\sf
c}\,d\overline{\sf a}}{{\sf a}^2 \overline{\sf a}^2},
\\
{\alpha_3} & := -\,\frac{{\sf c}\,d{\sf b}}{{\sf a}^3\overline{\sf
a}^2} + \bigg( \frac{{\sf b}{\sf c}}{{\sf a}^4\overline{\sf a}^2} -
\frac{{\sf e}}{{\sf a}^3\overline{\sf a}} \bigg)\,d{\sf a} +
\frac{1}{{\sf a}^2\overline{\sf a}}\,d{\sf e},
\\
{\alpha_4} & := \frac{d\overline{\sf d}}{{\sf a}\overline{\sf a}^2} -
\frac{\overline{\sf c}\,d{\sf b}}{{\sf a}^2\overline{\sf a}^3} +
\bigg( \frac{{\sf b}\overline{\sf c}}{{\sf a}^3\overline{\sf a}^3} -
\frac{\overline{\sf d}}{{\sf a}^2\overline{\sf a}^2} \bigg)\, d{\sf
a},
\\
{\alpha_5} & := \frac{d{\sf b}}{{\sf a}\overline{\sf a}} - \frac{{\sf
b}\,d{\sf a}}{{\sf a}^2\overline{\sf a}}.
\endaligned
\]
Here, the so-defined $5\times 5$ matrix $\omega_{\sf MC}$ of $1$-forms
is of course the {\sl Maurer-Cartan form} associated to our
$10$-dimensional structure group $G$.

Next, the above expressions of $\sigma_0$, $\rho_0$, $\zeta_0$ in
terms of the 1-forms $\sigma$, $\rho$, $\zeta$ and the conjugate
expressions as well enable us to express the five 2-forms $d\overline{
\sigma}_0$, $d \sigma_0$, $d\rho_0$, $d\overline{ \zeta}_0$,
$d\zeta_0$ of~\thetag{\ref{d0-model}} in terms of exterior products by
pairs of the $1$-forms $\overline{ \sigma}$, $\sigma$, $\rho$,
$\overline{ \zeta}$, $d \zeta$. After non-painful
computations, we obtain:
\[
\aligned
d\sigma_0
&
=
\sigma\wedge\overline{\sigma}\,
\bigg(
\frac{{\sf c}\overline{\sf d}}{
{\sf a}^5\overline{\sf a}^4}
-
\frac{{\sf e}\overline{\sf c}}{
{\sf a}^5\overline{\sf a}^4}
\bigg)
+
\sigma\wedge\rho\,
\bigg(
\frac{{\sf e}}{
{\sf a}^4\overline{\sf a}^4}
\bigg)
+
\sigma\wedge\zeta\,
\bigg(
\frac{{\sf c}}{{\sf a}^4\overline{\sf a}^2}
\bigg)
+
\\
&
\ \ \ \ \
+
\overline{\sigma}\wedge\rho\,
\bigg(
\frac{\overline{\sf d}}{{\sf a}^3\overline{\sf a}^3}
\bigg)
+
\overline{\sigma}\wedge\zeta\,
\bigg(
-
\frac{\overline{\sf c}}{{\sf a}^3\overline{\sf a}^3}
\bigg)
+
\rho\wedge\zeta\,
\bigg(
\frac{1}{{\sf a}^2\overline{\sf a}}
\bigg).
\endaligned
\]
and:
\[
\aligned
d\rho_0
&
=
\sigma\wedge\overline{\sigma}\,
\bigg(
-
i\,
\frac{{\sf b}{\sf c}\overline{\sf e}}{
{\sf a}^5\overline{\sf a}^5}
-
i\,\frac{\overline{\sf b}\overline{\sf c}{\sf e}}{
{\sf a}^5\overline{\sf a}^5}
+
i\,
\frac{{\sf e}\overline{\sf e}}{
{\sf a}^4\overline{\sf a}^4}
+
i\,
\frac{{\sf b}\overline{\sf c}{\sf d}}{
{\sf a}^5\overline{\sf a}^5}
+
i\,
\frac{\overline{\sf b}{\sf c}\overline{\sf d}}{
{\sf a}^5\overline{\sf a}^5}
-
i\,
\frac{{\sf d}\overline{\sf d}}{
{\sf a}^4\overline{\sf a}^4}
\bigg)
+
\\
&
\ \ \ \ \
+
\sigma\wedge\rho
\bigg(
i\,
\frac{\overline{\sf b}{\sf e}}{{\sf a}^4\overline{\sf a}^3}
-
i\,
\frac{{\sf b}{\sf d}}{{\sf a}^4\overline{\sf a}^3}
\bigg)
+
\sigma\wedge\zeta
\bigg(
-
i\,
\frac{\overline{\sf b}{\sf c}}{{\sf a}^4\overline{\sf a}^3}
+
i\,
\frac{{\sf d}}{{\sf a}^3\overline{\sf a}^2}
\bigg)
\\
&
+
\sigma\wedge\overline{\zeta}\,
\bigg(
i\,
\frac{{\sf b}{\sf c}}{
{\sf a}^4\overline{\sf a}^3}
-
i\,\frac{{\sf e}}{{\sf a}^3\overline{\sf a}^2}
\bigg)
+
\overline{\sigma}\wedge\rho\,
\bigg(
-
i\,
\frac{{\sf b}\overline{\sf e}}{
{\sf a}^3\overline{\sf a}^4}
+
i\,
\frac{\overline{\sf b}\overline{\sf d}}{
{\sf a}^3\overline{\sf a}^4}
\bigg)
+
\endaligned
\]
\[
\aligned
&
+
\overline{\sigma}\wedge\zeta\,
\bigg(
-
i\,\frac{\overline{\sf b}\overline{\sf c}}{
{\sf a}^3\overline{\sf a}^4}
+
i\,
\frac{\overline{\sf e}}{
{\sf a}^2\overline{\sf a}^3}
\bigg)
+
\overline{\sigma}\wedge\overline{\zeta}\,
\bigg(
i\,
\frac{{\sf b}\overline{\sf c}}{{\sf a}^3\overline{\sf a}^4}
-
i\,
\frac{\overline{\sf d}}{{\sf a}^2\overline{\sf a}^3}
\bigg)
+
\\
&
+
\rho\wedge\zeta\,
\bigg(
i\,\frac{\overline{\sf b}}{{\sf a}^2\overline{\sf a}^2}
\bigg)
+
\rho\wedge\overline{\zeta}\,
\bigg(
-\,i\,
\frac{{\sf b}}{{\sf a}^2\overline{\sf a}^2}
\bigg)
+
\frac{i}{{\sf a}\overline{\sf a}}\,
\zeta\wedge\overline{\sf a},
\endaligned
\]
while trivially:
\[
d\zeta_0
=
0.
\]
After that, we can express $d\overline{ \sigma}$, $d \sigma$, $d
\rho$, $d \overline{ \zeta}$ in the structure equation
\thetag{ \ref{struc-eq-1-model}} in terms of the one forms $\overline{
\sigma}$, $\sigma$, $\rho$, $\overline{ \zeta}$ instead of
$\overline{ \sigma}_0$, $\sigma_0$, $\rho_0$, $\overline{ \zeta}_0$
and we obtain:
\[
\aligned
d\sigma &
=
\big(2\,\alpha_1+\overline{\alpha}_1\big) \wedge
\sigma+
\\
& \ \ \ \ \ + \bigg( \frac{{\sf c}\overline{\sf d}}{{\sf
a}^3\overline{\sf a}^3}
-\frac{\overline{\sf c}{\sf e}}{{\sf
a}^3\overline{\sf a}^3}
\bigg)\, \sigma\wedge\overline{\sigma} +
\bigg( \frac{{\sf e}}{{\sf a}^2\overline{\sf a}} \bigg)\,
\sigma\wedge\rho + \bigg( -\,\frac{{\sf c}}{{\sf a}^2\overline{\sf
a}} \bigg)\, \sigma\wedge\zeta + 0 +
\\
& \ \ \ \ \ \ \ \ \ \ \ \ \ \ \ \ \ \ \ \ \ \ \ \ \ \ \ \ \ \ \ \ \ \
\ \ \ \ \ \ \ \ \ \ \ \ \ \,
+ \bigg( \frac{\overline{\sf d}}{{\sf
a}\overline{\sf a}^2} \bigg)\, \overline{\sigma}\wedge\rho + \bigg(
-\,\frac{\overline{\sf c}}{{\sf a}\overline{\sf a}^2} \bigg)\,
\overline{\sigma}\wedge\zeta + 0 +
\\
& \ \ \ \ \ \ \ \ \ \ \ \ \ \ \ \ \ \ \ \ \ \ \ \ \ \ \ \ \ \ \ \ \ \
\ \ \ \ \ \ \ \ \ \ \ \ \ \ \ \ \ \ \ \ \ \ \ \ \ \ \ \ \ \ \ \ \ \ \
\ \ \ \ \ \ \ \ \ \ \ \
 \ \ \ \ \ \ \ \ \ \ \ \ \
+ \rho\wedge\zeta + 0 +
\\
& \ \ \ \ \ \ \ \ \ \ \ \ \ \ \ \ \ \ \ \ \ \ \ \ \ \ \ \ \ \ \ \ \ \
\ \ \ \ \ \ \ \ \ \ \ \ \ \ \ \ \ \ \ \ \ \ \ \ \ \ \ \ \ \ \ \ \ \ \
\ \ \ \ \ \ \ \ \ \ \ \ \ \ \ \ \ \ \ \ \ \ \ \ \ \ \ \ \ \ \ \ \ \ \
\ \ \ \,
+0.
\endaligned
\]
For the order between the $10$ two-forms, we choose:
\[
\aligned
\sigma\wedge\overline{\sigma},\ \ \ \ \
\sigma\wedge\rho,\ \ \ \ \
\sigma\wedge\zeta,\ \ \ \ \
\sigma\wedge\overline{\zeta},
\\
\overline{\sigma}\wedge\rho,\ \ \ \ \
\overline{\sigma}\wedge\zeta,\ \ \ \ \
\overline{\sigma}\wedge\overline{\zeta},
\\
\rho\wedge\zeta,\ \ \ \ \
\rho\wedge\overline{\zeta},
\\
\zeta\wedge\overline{\zeta}.
\endaligned
\]
Let us abbreviate this first structure equation by
introducing specific numbered letters $U_\bullet$ for
the torsion coefficients:
\begin{equation}
\label{dsigma}
\aligned d\sigma & =
\big(2\,\alpha_1+\overline{\alpha}_1\big) \wedge \sigma+
\\
& \ \ \ \ \
+
U_1\,\sigma\wedge\overline{\sigma} +
U_2\,\sigma\wedge\rho + U_3\,\sigma\wedge\zeta +
\\
& \ \ \ \ \ + U_5\,\overline{\sigma}\wedge\rho +
U_6\,\overline{\sigma}\wedge\zeta + {\rho\wedge\zeta}.
\endaligned
\end{equation}
Similarly, one may compute elementarily the torsion
coefficients $V_1$, $V_2$, $V_3$, $V_4$, $V_8$ which appear
in the finalized expression of:
\begin{equation}
\label{drho}
 \aligned d\rho & = \alpha_2\wedge\sigma +
\overline{\alpha}_2\wedge\overline{\sigma} + \alpha_1\wedge\rho +
\overline{\alpha}_1\wedge\rho +
\\
& \ \ \ \ \ + V_1\,\sigma\wedge\overline{\sigma} +
V_2\,\sigma\wedge\rho + V_3\,\sigma\wedge\zeta +
V_4\,\sigma\wedge\overline{\zeta} +
\\
& \ \ \ \ \ \ \ \ \ \ \ \ \ \ \ \ \ \ \ \ \ \ \ \ +
\overline{V_2}\,\overline{\sigma}\wedge\rho +
\overline{V_4}\,\overline{\sigma}\wedge\zeta +
\overline{V_3}\,\overline{\sigma}\wedge\overline{\zeta} +
\\
& \ \ \ \ \ \ \ \ \ \ \ \ \ \ \ \ \ \ \ \ \ \ \ \ \ \ \ \ \ \ \ \ \ \
\ \ \ \ \ \ \ \ \ + V_8\,\rho\wedge\zeta +
\overline{V_8}\,\rho\wedge\overline{\zeta} +
\\
& \ \ \ \ \ \ \ \ \ \ \ \ \ \ \ \ \ \ \ \ \ \ \ \ \ \ \ \ \ \ \ \ \ \
\ \ \ \ \ \ \ \ \ \ \ \ \ \ \ \ \ \ \ \ \ \ \ \ \ \ \ \ \ \ +
i\,\zeta\wedge\overline{\zeta},
\endaligned
\end{equation}
and their explicit expressions are:
\[
\aligned
V_1
&
=
\frac{{\sf c}{\sf d}}{{\sf a}^3\overline{\sf a}^3}
+
\frac{{\sf c}{\sf c}\overline{\sf d}}{
{\sf a}^5\overline{\sf a}^4}
-
\frac{{\sf e}{\sf c}\overline{\sf c}}{
{\sf a}^5\overline{\sf a}^4}
-
\frac{\overline{\sf c}\overline{\sf c}{\sf d}}{
{\sf a}^4\overline{\sf a}^5}
+
\frac{{\sf c}\overline{\sf c}\overline{\sf e}}{
{\sf a}^4\overline{\sf a}^5}
-
\\
&
\ \ \ \ \
-\,i\,
\frac{{\sf b}{\sf c}\overline{\sf e}}{
{\sf a}^4\overline{\sf a}^4}
-
i\,
\frac{\overline{\sf b}\overline{\sf c}{\sf e}}{
{\sf a}^4\overline{\sf a}^4}
+
i\,
\frac{{\sf e}\overline{\sf e}}{
{\sf a}^3\overline{\sf a}^3}
+
i\,
\frac{{\sf b}\overline{\sf c}{\sf d}}{
{\sf a}^4\overline{\sf a}^4}
+
i\,
\frac{\overline{\sf b}{\sf c}\overline{\sf d}}{
{\sf a}^4\overline{\sf a}^4}
-
i\,
\frac{{\sf d}\overline{\sf d}}{
{\sf a}^3\overline{\sf a}^3},
\\
V_2
&
=
\frac{{\sf c}{\sf e}}{
{\sf a}^4\overline{\sf a}^2}
+
\frac{\overline{\sf c}{\sf d}}{
{\sf a}^3\overline{\sf a}^3}
+
i\,
\frac{\overline{\sf b}{\sf e}}{
{\sf a}^3\overline{\sf a}^2}
-
i\,
\frac{{\sf b}{\sf d}}{
{\sf a}^3\overline{\sf a}^2},
\endaligned
\]
\[
\aligned
V_3
&
=
-\,\frac{{\sf c}{\sf c}}{
{\sf a}^4\overline{\sf a}^4}
-
i\,
\frac{\overline{\sf b}{\sf c}}{
{\sf a}^3\overline{\sf a}^2}
+
i\,
\frac{{\sf d}}{
{\sf a}^2\overline{\sf a}},
\\
V_4
&
=
-\,
\frac{{\sf c}\overline{\sf c}}{
{\sf a}^3\overline{\sf a}^3}
+
i\,
\frac{{\sf b}{\sf c}}{
{\sf a}^3\overline{\sf a}^2}
-
i\,
\frac{{\sf e}}{
{\sf a}^2\overline{\sf a}},
\\
V_8
&
=
\frac{{\sf c}}{
{\sf a}^2\overline{\sf a}}
+
i\,
\frac{\overline{\sf b}}{
{\sf a}\overline{\sf a}}.
\endaligned
\]
Lastly, the $10$ torsion coefficients $W_\bullet$ in:
\begin{equation}
\label{dzeta}
\aligned d\zeta & = \alpha_3\wedge\sigma
+ \alpha_4\wedge\overline{\sigma} + \alpha_5\wedge\rho +
\alpha_1\wedge\zeta +
\\
& \ \ \ \ \ + W_1\,\sigma\wedge\overline{\sigma} +
W_2\,\sigma\wedge\rho + W_3\,\sigma\wedge\zeta +
W_4\,\sigma\wedge\overline{\zeta} +
\\
& \ \ \ \ \ \ \ \ \ \ \ \ \ \ \ \ \ \ \ \ \ \ \ \ +
W_5\,\overline{\sigma}\wedge\rho + W_6\,\overline{\sigma}\wedge\zeta
+ W_7\,\overline{\sigma}\wedge\overline{\zeta} +
\\
& \ \ \ \ \ \ \ \ \ \ \ \ \ \ \ \ \ \ \ \ \ \ \ \ \ \ \ \ \ \ \ \ \ \
\ \ \ \ \ \ \ \ \ + W_8\,\rho\wedge\zeta +
W_9\,\rho\wedge\overline{\zeta} +
\\
& \ \ \ \ \ \ \ \ \ \ \ \ \ \ \ \ \ \ \ \ \ \ \ \ \ \ \ \ \ \ \ \ \ \
\ \ \ \ \ \ \ \ \ \ \ \ \ \ \ \ \ \ \ \ \ \ \ \ \ \ \ \ \ +
W_{10}\,\zeta\wedge\overline{\zeta},
\endaligned
\end{equation}
are equal to:
\[
\aligned
W_1
&
=
\frac{{\sf d}{\sf e}}{
{\sf a}^3\overline{\sf a}^4}
+
\frac{{\sf c}\overline{\sf d}{\sf e}}{
{\sf a}^5\overline{\sf a}^4}
-
\frac{\overline{\sf c}{\sf e}{\sf e}}{
{\sf a}^5\overline{\sf a}^4}
-
\frac{\overline{\sf c}{\sf d}\overline{\sf d}}{
{\sf a}^4\overline{\sf a}^5}
+
\frac{{\sf c}\overline{\sf d}\overline{\sf e}}{
{\sf a}^4\overline{\sf a}^5}
-
\\
&
\ \ \ \ \
-
i\,\frac{{\sf b}{\sf b}{\sf c}\overline{\sf e}}{
{\sf a}^5\overline{\sf a}^5}
-
i\,
\frac{{\sf b}\overline{\sf b}\overline{\sf c}{\sf e}}{
{\sf a}^5\overline{\sf a}^5}
+
i\,\frac{{\sf b}{\sf e}\overline{\sf e}}{
{\sf a}^4\overline{\sf a}^4}
+
i\,\frac{{\sf b}{\sf b}\overline{\sf c}{\sf d}}{
{\sf a}^5\overline{\sf a}^5}
+
i\,\frac{{\sf b}\overline{\sf b}{\sf c}\overline{\sf d}}{
{\sf a}^5\overline{\sf a}^5}
-
i\,
\frac{{\sf b}{\sf d}\overline{\sf d}}{
{\sf a}^4\overline{\sf a}^4},
\\
W_2
&
=
\frac{{\sf e}{\sf e}}{
{\sf a}^4\overline{\sf a}^4}
+
\frac{{\sf d}\overline{\sf d}}{
{\sf a}^3\overline{\sf a}^3}
+
i\,\frac{{\sf b}\overline{\sf b}{\sf e}}{
{\sf a}^4\overline{\sf a}^3}
-
i\,
\frac{{\sf b}{\sf b}{\sf d}}{
{\sf a}^4\overline{\sf a}^3},
\\
W_3
&
=
-\,
\frac{{\sf c}{\sf e}}{
{\sf a}^4\overline{\sf a}^2}
-i\,
\frac{{\sf b}\overline{\sf b}{\sf c}}{
{\sf a}^4\overline{\sf a}^3}
+
i\,
\frac{{\sf b}{\sf d}}{
{\sf a}^3\overline{\sf a}^2},
\\
W_4
&
=
-\frac{{\sf c}\overline{\sf d}}{
{\sf a}^3\overline{\sf a}^3}
+
i\,\frac{{\sf b}{\sf b}{\sf c}}{
{\sf a}^4\overline{\sf a}^4}
-
i\,
\frac{{\sf b}{\sf e}}{{\sf a}^3\overline{\sf a}^2},
\\
W_5
&
=
\frac{\overline{\sf d}{\sf e}}{
{\sf a}^3\overline{\sf a}^3}
+
\frac{\overline{\sf e}\overline{\sf d}}{
{\sf a}^2\overline{\sf a}^4}
-
i\,\frac{{\sf b}{\sf b}\overline{\sf e}}{
{\sf a}^3\overline{\sf a}^4}
+
i\,
\frac{{\sf b}\overline{\sf b}\overline{\sf d}}{
{\sf a}^3\overline{\sf a}^4},
\\
W_6
&
=
-\,\frac{\overline{\sf c}{\sf e}}{
{\sf a}^3\overline{\sf a}^3}
-
i\,
\frac{{\sf b}\overline{\sf b}\overline{\sf c}}{
{\sf a}^3\overline{\sf a}^4}
+
i\,\frac{{\sf b}\overline{\sf e}}{
{\sf a}^2\overline{\sf a}^3},
\\
W_7
&
=
-\,\frac{\overline{\sf c}\overline{\sf d}}{
{\sf a}^2\overline{\sf a}^4}
+
i\,
\frac{{\sf b}{\sf b}\overline{\sf c}}{
{\sf a}^3\overline{\sf a}^4}
-
i\,
\frac{{\sf b}\overline{\sf d}}{\
{\sf a}^2\overline{\sf a}^3},
\\
W_8
&
=
\frac{{\sf e}}{
{\sf a}^2\overline{\sf a}}
+
i\,
\frac{{\sf b}}{
{\sf a}^2\overline{\sf a}^2},
\\
W_9
&
=
\frac{\overline{\sf d}}{
{\sf a}\overline{\sf a}^2}
-
i\,
\frac{{\sf b}{\sf b}}{
{\sf a}^2\overline{\sf a}^2},
\\
W_{10}
&
=
i\,
\frac{{\sf b}}{
{\sf a}\overline{\sf a}}.
\endaligned
\]

\subsection{First-loop absorbtion}
\label{first-loop-model}
Now, we are ready to realize which of the above torsion coefficients
are {\it normalizable}. According to Proposition~\ref{prop-changes}, we
must modify the five $1$-forms $\alpha_1$, $\alpha_2$, $\alpha_3$,
$\alpha_4$, $\alpha_5$ by adding to them general
linear combinations of the
$1$-forms $\sigma$, $\overline{ \sigma}$, $\rho$, $\zeta$,
$\overline{ \zeta}$:
\[
\aligned
\alpha_1
&
\longmapsto
\alpha_1
+
p_1\,\sigma
+
q_1\,\overline{\sigma}
+
r_1\,\rho
+
s_1\,\zeta
+
t_1\,\overline{\zeta},
\\
\alpha_2
&
\longmapsto
\alpha_2
+
p_2\,\sigma
+
q_2\,\overline{\sigma}
+
r_2\,\rho
+
s_2\,\zeta
+
t_2\,\overline{\zeta},
\\
\alpha_3
&
\longmapsto
\alpha_3
+
p_3\,\sigma
+
q_3\,\overline{\sigma}
+
r_3\,\rho
+
s_3\,\zeta
+
t_3\,\overline{\zeta},
\\
\alpha_4
&
\longmapsto
\alpha_4
+
p_4\,\sigma
+
q_4\,\overline{\sigma}
+
r_4\,\rho
+
s_4\,\zeta
+
t_4\,\overline{\zeta},
\\
\alpha_5
&
\longmapsto
\alpha_5
+
p_5\,\sigma
+
q_5\,\overline{\sigma}
+
r_5\,\rho
+
s_5\,\zeta
+
t_5\,\overline{\zeta},
\endaligned
\]
with $25$ arbitrary real analytic functions $p_i$, $q_i$, $r_i$,
$s_i$, $t_i$, hence by conjugation, we
also have the two useful replacements:
\[
\aligned
\overline{\alpha}_1
&
\longmapsto
\overline{\alpha}_1
+
\overline{q}_1\,\sigma
+
\overline{p}_1\,\overline{\sigma}
+
\overline{r}_1\,\rho
+
\overline{t}_1\,\zeta
+
t_1\,\overline{\zeta},
\\
\overline{\alpha}_2
&
\longmapsto
\overline{\alpha}_2
+
\overline{q}_2\,\sigma
+
\overline{p}_2\,\overline{\sigma}
+
\overline{r}_2\,\rho
+
\overline{t}_2\,\zeta
+
t_2\,\overline{\zeta}.
\endaligned
\]
Performing the replacement with these modified Maurer-Cartan forms,
we obtain the three new structure equations
which change the expressions of $d\sigma,d\rho$ and $d\zeta$
in~\thetag{\ref{dsigma}}, \thetag{\ref{drho}} and
\thetag{\ref{dzeta}} as follows:
\[
\!\!\!\!\!\!\!\!\!\!\!\!\!\!\!\!\!\!\!\!
\small
\aligned
d\sigma
&
=
\big(2\alpha_1+\overline{\alpha}_1\big)
\wedge
\sigma
+
\\
&
+
\sigma\wedge\overline{\sigma}
\big[U_1-2q_1-\overline{p}_1\big]
+
\sigma\wedge\rho
\big[
U_2-2r_1-\overline{r}_1
\big]
+
\sigma\wedge\zeta
\big[
U_3-2s_1-\overline{t}_1
\big]
+
\sigma\wedge\overline{\zeta}
\big[
0-2t_1-\overline{s}_1
\big]
\\
&
\ \ \ \ \ \ \ \ \ \ \ \ \ \ \ \ \ \ \ \ \ \ \ \ \ \ \ \ \ \ \ \ \ \ \
\ \ \ \ \ \ \,
+
\overline{\sigma}\wedge\rho
\big[U_5\big]
+
\overline{\sigma}\wedge\zeta
\big[U_6\big]
+
0
+
\\
&
\ \ \ \ \ \ \ \ \ \ \ \ \ \ \ \ \ \ \ \ \ \ \ \ \ \ \ \ \ \ \ \ \ \ \
\ \ \ \ \ \ \ \ \ \ \ \ \ \ \ \ \ \ \ \ \ \ \ \ \ \ \ \ \,
+
\rho\wedge\zeta,
\endaligned
\]
\[
\!\!\!\!\!\!\!\!\!\!\!\!\!\!\!\!\!\!\!\!
\aligned
d\rho
&
=
\alpha_2\wedge\sigma
+
\overline{\alpha}_2\wedge\overline{\sigma}
+
\alpha_1\wedge\rho
+
\overline{\alpha}_1\wedge\rho
+
\\
&
\ \ \ \ \
+
\sigma\wedge\overline{\sigma}
\big[
V_1-q_2-\overline{q}_2
\big]
+
\sigma\wedge\rho
\big[
V_2-r_2+p_1+\overline{q}_1
\big]
+
\sigma\wedge\zeta
\big[
V_3-s_2
\big]
+
\sigma\wedge\overline{\zeta}
\big[
V_4-t_2
\big]
\\
&
\ \ \ \ \ \ \ \ \ \ \ \ \ \ \ \ \ \ \ \ \ \ \ \ \ \ \ \ \ \ \ \ \ \ \
+
\overline{\sigma}\wedge\rho
\big[
\overline{V}_2-\overline{r}_2+q_1+\overline{p}_1
\big]
+
\overline{\sigma}\wedge\zeta
\big[
\overline{V}_4-\overline{t}_2
\big]
+
\overline{\sigma}\wedge\overline{\zeta}
\big[
\overline{V}_3-\overline{s}_2
\big]
+
\\
&
\ \ \ \ \ \ \ \ \ \ \ \ \ \ \ \ \ \ \ \ \ \ \ \ \ \ \ \ \ \ \ \ \ \ \
\ \ \ \ \ \ \ \ \ \ \ \ \ \ \ \ \ \ \ \ \ \ \ \ \ \ \ \ \ \ \ \ \ \ \
+
\rho\wedge\zeta
\big[
V_8-s_1-\overline{t}_1
\big]
+
\rho\wedge\overline{\zeta}
\big[
\overline{V}_8-t_1-\overline{s}_1
\big]
+
\\
&
\ \ \ \ \ \ \ \ \ \ \ \ \ \ \ \ \ \ \ \ \ \ \ \ \ \ \ \ \ \ \ \ \ \ \
\ \ \ \ \ \ \ \ \ \ \ \ \ \ \ \ \ \ \ \ \ \ \ \ \ \ \ \ \ \ \ \ \ \ \
\ \ \ \ \ \ \ \ \ \ \ \ \ \ \ \ \ \ \ \ \ \ \ \ \ \ \ \ \ \ \ \ \ \ \
\ \
+
i\,\zeta\wedge\overline{\zeta},
\endaligned
\]
\[
\!\!\!\!\!\!\!\!\!\!\!\!\!\!\!\!\!\!\!\!
\aligned
d\zeta
&
=
\alpha_3\wedge\sigma
+
\alpha_4\wedge\overline{\sigma}
+
\alpha_5\wedge\rho
+
\alpha_1\wedge\zeta
+
\\
&
\ \ \ \ \
+
\sigma\wedge\overline{\sigma}
\big[
W_1-q_3+p_4
\big]
+
\sigma\wedge\rho
\big[
W_2-r_3+p_5
\big]
+
\sigma\wedge\zeta
\big[
W_3-s_3+p_1
\big]
+
\sigma\wedge\overline{\zeta}
\big[
W_4-t_3
\big]
\\
&
\ \ \ \ \ \ \ \ \ \ \ \ \ \ \ \ \ \ \ \ \ \ \ \ \ \ \ \ \ \ \ \ \ \ \
+
\overline{\sigma}\wedge\rho
\big[
W_5-r_4+q_5
\big]
+
\overline{\sigma}\wedge\zeta
\big[
W_6-s_4+q_1
\big]
+
\overline{\sigma}\wedge\overline{\zeta}
\big[
W_7-t_4
\big]
+
\\
&
\ \ \ \ \ \ \ \ \ \ \ \ \ \ \ \ \ \ \ \ \ \ \ \ \ \ \ \ \ \ \ \ \ \ \
\ \ \ \ \ \ \ \ \ \ \ \ \ \ \ \ \ \ \ \ \ \ \ \ \ \ \ \ \ \ \ \ \ \ \
+
\rho\wedge\zeta
\big[
W_8-s_5+r_1
\big]
+
\rho\wedge\overline{\zeta}
\big[
W_9-t_5
\big]
+
\\
&
\ \ \ \ \ \ \ \ \ \ \ \ \ \ \ \ \ \ \ \ \ \ \ \ \ \ \ \ \ \ \ \ \ \ \
\ \ \ \ \ \ \ \ \ \ \ \ \ \ \ \ \ \ \ \ \ \ \ \ \ \ \ \ \ \ \ \ \ \ \
\ \ \ \ \ \ \ \ \ \ \ \ \ \ \ \ \ \ \ \ \ \ \ \ \ \ \ \ \ \ \ \ \ \ \
\ \ \
+
\zeta\wedge\overline{\zeta}
\big[
W_{10}-t_1
\big].
\endaligned
\]

\subsection{First loop normalization}

Now, according to the general
principles (\cite{ Olver-1995}, Chapter~10), in order to know what
are the precise linear combinations of the $20$ torsion coefficients:
\[
\aligned
&
U_1,\ \ \ \ \ U_2,\ \ \ \ \ U_3,\ \ \ \ \ U_5, \ \ \ \ \ U_6,
\\
&
V_1,\ \ \ \ \ V_2,\ \ \ \ \ V_3,\ \ \ \ \ V_4,\ \ \ \ \ V_8,
\\
&
W_1,\ \ \ \ \ W_2,\ \ \ \ \ W_3,\ \ \ \ \ W_4,\ \ \ \ \
W_5,\ \ \ \ \ W_6,\ \ \ \ \ W_7,\ \ \ \ \ W_8, \ \ \ \ \
W_9,\ \ \ \ \ W_{10}
\endaligned
\]
that are {\em necessarily normalizable},
one must determine all possible linear
combinations of the following $6 + 5 + 10 = 21$
equations\,\,---\,\,including their (unwritten) conjugates\,\,---:
\begin{equation}
\label{1-norm-eqs}
\left[
\aligned
U_1
&
=
2q_1+\overline{p}_1,
\\
U_2
&
=
2r_1+\overline{r}_1,
\\
U_3
&
=
2s_1+\overline{t}_1,
\\
0
&
=
2t_1+\overline{s}_1,
\\
U_5
&
=
0,
\\
U_6
&
=
0,
\endaligned\right.
\ \ \ \ \ \ \ \ \ \ \ \ \ \ \ \ \ \ \ \ \
\left[
\aligned
V_1
&
=
q_2-\overline{q}_2,
\\
V_2
&
=
r_2-p_1-\overline{q}_1,
\\
V_3
&
=
s_2,
\\
V_4
&
=
t_2,
\\
V_8
&
=
s_1+\overline{t}_1,
\endaligned\right.
\ \ \ \ \ \ \ \ \ \ \ \ \ \ \ \ \ \ \ \ \
\left[
\aligned
W_1
&
=
q_3-p_4,
\\
W_2
&
=
r_3-p_5,
\\
W_3
&
=
s_3-p_1,
\\
W_4
&
=
t_3,
\\
W_5
&
=
r_4-q_5,
\\
W_6
&
=
s_4-q_1,
\\
W_7
&
=
t_4,
\\
W_8
&
=
s_5-r_1,
\\
W_9
&
=
t_5,
\\
W_{10}
&
=
t_1.
\endaligned\right.
\end{equation}
so as to obtain null right-hand sides, though without exchanging any
left-hand side term with any right-hand side term.

One then easily convinces oneself just visually that some complete
appropriate linear combinations are:
\[
\aligned
0
&
=
U_5,
\\
0
&
=
U_6,
\\
0
&
=
U_3-3\,V_8,
\\
0
&
=
V_8+\overline{W}_{10},
\endaligned
\]
and that is all. This means that these four right
linear combinations are normalizable\,\,---\,\,minding that
this synoptic way of finding normalizable linear combinations
with the notation:
\[
`0
=
\mathmotsf{combination of torsion coefficients'}
\]
does not
necessarily mean that one will
assign the value $0$ to them, and in fact, some normalizable
right-hand sides could well be assigned
the value $1$ in certain circumstances.

Now, if one just replaces the appearing torsion coefficients
with the values computed a while ago, one plainly obtains:
\[
\aligned
U_5
&
=
\frac{\overline{\sf d}}{{\sf a}\overline{\sf a}^2},
\\
U_6
&
=
-\,
\frac{\overline{\sf c}}{
{\sf a}\overline{\sf a}^2},
\\
U_3-3\,V_8
&
=
-\,
\frac{{\sf c}}{{\sf a}^2\overline{\sf a}}
-
3\,
\frac{{\sf c}}{
{\sf a}^2\overline{\sf a}}
-
3i\,
\frac{\overline{\sf b}}{
{\sf a}\overline{\sf a}},
\\
V_8
+
\overline{W}_{10}
&
=
\frac{{\sf c}}{{\sf a}^2\overline{\sf a}}.
\endaligned
\]
Since none of the $3$ group parameters ${\sf b}$, ${\sf c}$, ${\sf d}$
appearing in denominator place here does belong to the
diagonal of our initial lower triangular
matrix subgroup of ${\sf GL}_5 ( \C)$, it is clear
that we can normalize all of them to zero, namely
we can set:
\[
\boxed{
\aligned
{\sf b}:=0,
\ \ \ \ \ \ \ \ \ \
{\sf c}:=0,
\ \ \ \ \ \ \ \ \ \
{\sf d}:=0.
\endaligned}
\]

\subsection{Second-loop absorbtion and normalization}
This assignment then considerably simplifies a lot of the
expressions of the three differentials
$d\sigma$, $d\rho$, $d\zeta$
in~\thetag{\ref{dsigma}}, \thetag{\ref{drho}}, \thetag{\ref{dzeta}}).

Also, the Maurer-Cartan
1-forms $\alpha_2$, $\alpha_4$, $\alpha_5$
then vanish identically and
the new Maurer-Cartan matrix $\omega_{\sf MC}$ changes into the
form\,\,---\,\,we employ a new letter
$\beta$ instead of $\alpha$ in this second loop and we re-number
them\,\,---:
\[\footnotesize
\aligned
\omega_{\sf MC}=\left(%
\begin{array}{ccccc}
2\,\overline{\beta}_1+\beta_1 & 0 & 0 & 0 & 0 \\
0 & 2\,\beta_1+\overline{\beta}_1 & 0 & 0 & 0 \\
0 & 0 & \beta_1+\overline{\beta}_1 & 0 & 0 \\
\overline{\beta}_2 & 0 & 0 & \overline{\beta}_1 & 0 \\
0 & {\beta}_2 & 0 & 0 & {\beta}_1 \\
\end{array}
\right),
\endaligned
\]
with two non-vanishing 1-forms:
\[
\aligned
\beta_1
&
:=
\frac{d{\sf a}}{\sf a},
\\
\beta_2
&
:=
-\,\frac{{\sf e}}{{\sf a}^3\overline{\sf a}}
d{\sf a}
+
\frac{1}{{\sf a}^2\overline{\sf a}}\,d{\sf e}.
\endaligned
\]
Moreover, among the $20$ torsion coefficients $U_\bullet$,
$V_\bullet$, $W_\bullet$, only the following five simplified ones
remain non-vanishing:
\[
\aligned
U_2&=\frac{\sf e}{{\sf a}^2\overline{\sf a}},
\\
V_1&=\frac{i\,{\sf e}\overline{\sf e}}{{\sf a}^3\overline{\sf a}^3},
\\
V_4&=-\frac{i\,{\sf e}}{{\sf a}^2\overline{\sf a}},
\\
W_2&=\frac{{\sf e}^2}{{\sf a}^4\overline{\sf a}^2},
\\
W_8&=\frac{\sf e}{{\sf a}^2\overline{\sf a}}.
\endaligned
\]
Hence, one has the following reduced expressions for $d\sigma$,
$d\rho$, $d\zeta$:
\begin{equation}
\label{d2}
\aligned d\sigma & =
\big(2\,\beta_1+\overline{\beta}_1\big) \wedge \sigma+
U_2\,\sigma\wedge\rho + {\rho\wedge\zeta},
\\
d\rho & = (\beta_1+\overline{\beta}_1)\wedge\rho
+
V_1\,\sigma\wedge\overline{\sigma} +
V_4\,\sigma\wedge\overline{\zeta} +
\overline{V}_4\,\overline{\sigma}\wedge\zeta +
i\,\zeta\wedge\overline{\zeta},
\\
d\zeta & = \beta_2\wedge\sigma +
\beta_1\wedge\zeta + W_2\,\sigma\wedge\rho + W_8\,\rho\wedge\zeta.
\endaligned
\end{equation}

Again we must modify the two $1$-forms $\beta_1$,
$\beta_2$ by adding to them general linear combinations
of the $1$-forms $\sigma$, $\overline{ \sigma}$, $\rho$,
$\zeta$, $\overline{ \zeta}$:
\[
\aligned
\beta_1
&
\longmapsto
\beta_1
+
p_1\sigma
+
q_1\overline{\sigma}
+
r_1\rho
+
s_1\zeta
+
t_1\overline{\zeta},
\\
\beta_2
&
\longmapsto
\beta_2
+
p_2\sigma
+
q_2\overline{\sigma}
+
r_2\rho
+
s_2\zeta
+
t_2\overline{\zeta},
\endaligned
\]
with $10$ arbitrary real analytic functions $p_i$, $q_i$,
$r_i$, $s_i$, $t_i$. To get the new absorption equations,
it just suffices to set ${\sf b} := 0$,
${\sf c} := 0$, ${\sf d} := 0$ in the ones obtained
a moment ago, remembering that some torsion coefficients
vanish (are normalized) and that there is a change of
numbering due to $\beta_2 := \alpha_3$:
\[
\left[
\aligned
U_1
&
=
2q_1+\overline{p}_1,
\\
U_2
&
=
2r_1+\overline{r}_1,
\\
0
&
=
2s_1+\overline{t}_1,
\\
0
&
=
2t_1+\overline{s}_1,
\\
0
&
=
0,
\\
0
&
=
0,
\endaligned\right.
\ \ \ \ \ \ \ \ \ \ \ \ \ \ \ \ \ \ \ \ \
\left[
\aligned
V_1
&
=
0,
\\
0
&
=
-\,p_1-\overline{q}_1,
\\
0
&
=
0,
\\
V_4
&
=
0,
\\
0
&
=
s_1+\overline{t}_1,
\endaligned\right.
\ \ \ \ \ \ \ \ \ \ \ \ \ \ \ \ \ \ \ \ \
\left[
\aligned
0
&
=
q_2,
\\
W_2
&
=
r_2,
\\
0
&
=
s_2,
\\
0
&
=
t_2,
\\
0
&
=
0,
\\
0
&
=
-\,q_1,
\\
0
&
=
0,
\\
W_8
&
=
-\,r_1,
\\
0
&
=
0,
\\
0
&
=
t_1.
\endaligned\right.
\]
Visually, one sees that three linear combinations of torsion
coefficients are normalizable at this step, which are represented as
the equations:
\[
\aligned
V_1
&
=
0,
\\
V_4
&
=
0,
\\
2\,W_8+\overline{W}_8+U_2
&
=0.
\endaligned
\]
Thus, we can normalize:
\[
\aligned
V_1
&
=
\frac{i\,{\sf e}\overline{\sf e}}{{\sf
a}^3\overline{\sf a}^3},
\\
V_4
&
=
-\,\frac{i\,{\sf e}}{{\sf a}^2\overline{\sf a}},
\\
2\,W_8+\overline{W}_8+U_2
&
=
\frac{3\,{\sf e}}{{\sf a}^2\overline{\sf
a}}+\frac{\overline{\sf e}}{{\sf a}\overline{\sf a}^2},
\endaligned
\]
and immediately, this amounts to just annihilating:
\[
\boxed{
{\sf e}:=0}\,.
\]

After really setting this group parameter null, all the torsion
coefficients vanish identically, the Maurer-Cartan 1-form $\beta_2$
also annihilates, and the expression of $\beta_1$ changes into the
simpler form:
\[
\beta_1=\frac{d\sf a}{\sf a}.
\]
Finally, one has the following greatly simplified expressions for
$d\sigma$, $d\rho$, $d\zeta$:
\begin{equation}
\label{d3}
\aligned d\sigma & =
\big(2\,\beta_1+\overline{\beta}_1\big)
\wedge
\sigma
+
{\rho\wedge\zeta},
\\
d\rho & =
\big(\beta_1+\overline{\beta}_1\big)
\wedge\rho
+
i\,\zeta\wedge\overline{\zeta},
\\
d\zeta & = \beta_1\wedge\zeta,
\endaligned
\end{equation}
in which no more nonconstant torsion coefficient appears.
Therefore, executing once more an absorbtion-normalization loop is not
valuable and we should start the {\sl prolongation step}.
But beforehand,
let us present an important application of Cartan's Lemma.

\begin{Lemma}
\label{unique-alpha} The Maurer-Cartan form $\beta_1=\frac{d\sf
a}{\sf a}$ is the only 1-form which enjoys the structure equations
\thetag{\ref{d3}}.
\end{Lemma}

\proof
Let $\beta_1'$ be another 1-form, enjoying the structure equations
\thetag{\ref{d3}}. Then subtracting by pairs the expressions of
$d\sigma$, $d\rho$, $d\zeta$ with $\beta_1$ and with $\beta_1'$
immediately gives:
\[
\aligned
0
&
\equiv
\big(
2\,\beta_1+\overline{\beta}_1
-
2\,\beta_1'-\overline{\beta}_1'
\big)
\wedge\sigma,
\\
0
&
\equiv
\big(
\beta_1+\overline{\beta}_1
-
\beta_1'-\overline{\beta}_1'
\big)
\wedge\rho,
\\
0
&
\equiv
\big(\beta_1-\beta_1')\wedge\zeta.
\endaligned
\]
Applying the Cartan's Lemma~\ref{Cartan-lemma} on the second and
on the third equations yields that one has:
\begin{equation}
\label{AB}
\aligned
\beta_1-\beta_1'+\overline{\beta}_1-\overline{\beta}_1'
&
=
A\,\rho,
\\
\beta_1-\beta_1'
&
=
B\,\zeta,
\endaligned
\end{equation}
for some certain functions $A$, $B$. Now, substituting the
expression $\beta'_1=\beta_1-B\,\zeta$ obtained from the second
equality into the first one gives:
\[
B\,\zeta+\overline{B}\,\overline{\zeta}-A\,\rho=0,
\]
and consequently, we have $A=B=0$, due to the fact that $\sigma$,
$\rho$, $\zeta$ are linearly independent. Now, the second equation of
\thetag{\ref{AB}} immediately implies that:
\[
\beta_1'=\beta_1,
\]
and hence the Maurer-Cartan 1-form $\beta_1$ is unique, as was
claimed.
\endproof

\subsection{Prolongation: checking firstly non-involutiveness}
\label{involutive}

At this stage, we shall be interested in testing whether the obtained
coframe $\big\{ \sigma, \overline{\sigma}, \rho, \zeta,
\overline{\zeta} \big\}$ is involutive. If this holds, the structure
group is infinite-dimensional, else we have to prolong the lifted
coframe. We stick to Chapter~11 of Olvers' book~\cite{ Olver-1995}.

Assume that after performing
all loops of the absorbtion-normalization
procedure, the structure equations for the
lifted coframe have the form:
\[
\aligned
d\theta^i
=
\sum_{k=1}^n\bigg(\sum_{s=1}^r\,a_{ks}^i\,\alpha^s
+
\sum_{j=1}^{k-1}T^i_{jk}\,\theta^j\bigg)\wedge\theta^k
\ \ \ \ \ \ \ \ \ \ \ \ \
{\scriptstyle{(i\,=\,1\,\cdots\,n)}},
\endaligned
\]
where $\alpha^1, \dots, \alpha^r$ are a basis of Maurer-Cartan 1-forms
on the group, that is to say, assume that none of the
essential\,\,---\,\,{\em i.e.}  unabsorbable\,\,---\,\,torsion
coefficients depend explicitly on the group parameters, since
otherwise, some parameters would be normalizable.  Modifying each
Maurer-Cartan form by adding to it a linear combination of coframe
elements:
\[
\alpha^s
\longmapsto
\alpha^s+\sum_{j=1}^n\,z^s_j\,\theta^j
\ \ \ \ \ \ \ \ \ \ \ \ \ \ \ \ \ \ \ \ \ \ \ \ \ \ \ \ \
{\scriptstyle{(s\,=\,1\,\cdots\,r)}},
\]
results in the {\sl linear absorption equations}:
\[
\sum_{s=1}^r\,
\big(
a_{ks}^i\,z_j^s
-
a_{js}^i\,z_k^s
\big)
=
T_{jk}^i
\ \ \ \ \ \ \ \ \ \ \
{\scriptstyle{(i\,=\,1\,\cdots\,n\,;\,\,\,
1\,\leqslant\,j\,<\,k\,\leqslant\,n)}}.
\]

\begin{Definition}
The {\sl degree of indeterminacy}
of a lifted coframe is the number of free variables in the associated
homogeneous linear equations:
\[
\sum_{s=1}^r\,
\big(
a_{ks}^i\,z_j^s
-
a_{js}^i\,z_k^s
\big)
=
0
\ \ \ \ \ \ \ \ \ \ \
{\scriptstyle{(i\,=\,1\,\cdots\,n\,;\,\,\,
1\,\leqslant\,j\,<\,k\,\leqslant\,n)}}.
\]
\end{Definition}

Next, for a given vector $v = (v^1, \dots, v^n) \in \R^n$, introduce
$I(v)$ to be the $n \times r$ matrix with the entries:
\[
I(v)^i_s
:=
\sum_{l=1}^n\,a^i_{sl}\,v^l
\ \ \ \ \ \ \ \ \ \ \ \ \
{\scriptstyle{(i\,=\,1\,\cdots\,n\,;\,\,\,
s\,=\,1\,\cdots\,r)}}.
\]

\begin{Definition}
The {\sl first $n-1$ reduced characters} $s'_1, s_2',\ldots,s'_{n-1}$ are
defined, for $k = 1, 2, \dots, n-1$, recursively by:
\[\footnotesize
\aligned s'_1+s'_2+\ldots+s'_k
=
\mathmotsf{max}\,\left\{ \mathmotsf{rank}\,\left(%
\begin{array}{c}
I(v_1) \\
I(v_2) \\
\vdots \\
I(v_k) \\
\end{array}%
\right)\colon
\ \ \ \ v_1,\ldots,v_k\in\mathbb R^n \right\},
\endaligned
\]
the ranks being always $\leqslant r$. Moreover, the final $n$-th
reduced character is defined by:
\[
s'_1+s'_2+\ldots+s'_{n-1}+s'_n=r.
\]
\end{Definition}

\begin{Definition}
Lastly, the coframe $\theta$ is said to be
{\sl involutive} if the value of the sum:
\[
\sum_{k=1}^n\,k\,s'_k
\]
is equal to the degree of indeterminancy.
\end{Definition}

Now, Lemma~\ref{unique-alpha}
showed that if we execute again the absorbtion procedure
on the structure equation~\thetag{\ref{d3}} by replacing $\beta_1$
with:
\[
\beta_1
+p_1\sigma+q_1\overline{\sigma}+r_1\rho+s_1\zeta+t_1\overline{\zeta},
\]
then in order to annihilate
all the coefficients in the new expressions of $d\sigma$,
$d\rho$,
$d\zeta$, the only solution is:
\[
p_1=q_1=r_1=s_1=t_1:=0.
\]
In other words, the number of free variables is null in our case.

On the other hand, we claim that the reduced characters cannot all be
null, and hence the lifted coframe is certainly non-involutive. Let us
check this fact by constructing the $5\times 2$ matrix:
\[\footnotesize\aligned
\alpha \ \ \ \ \ \ \ \ \overline{\alpha}
\ \ \ \ \ \ \ \ \ \ \ \ \ \ \ \ \ \ \ \ \
\ \ \\
I(v):=\begin{array}{cc}
\left(%
\begin{array}{cc}
2\,v^\sigma & v^\sigma \\
v^{\overline{\sigma}} & 2\,v^{\overline{\sigma}} \\
v^\rho & v^\rho \\
v^\zeta & 0 \\
0 & v^{\overline{\zeta}} \\
\end{array}%
\right)
& \begin{array}{c}
d\sigma \\
d\overline{\sigma} \\
d\rho \\
d\zeta \\
d\overline{\zeta} \\
\end{array}
\end{array}
\endaligned
\]
according to the structure equations~\thetag{\ref{d3}} and their
conjugates. Then, by the above definition we have:
\[
s'_1
=
\mathmotsf{rank}
\big(I(v)\big)
=
2,
\]
which is the maximum rank of the above matrix for all arbitrary
vectors $\big( v^\sigma, v^{\overline{\sigma}}, v^\rho, v^\zeta,
v^{\overline{\zeta}}) \in \mathbb{ R}^5$.  Consequently, the sum
$\sum_{k=1}^5\,k\,s'_k \geqslant 2$ is {\em larger} than the number,
$0$, of free variables.

In conclusion, the associated lifted coframe $\sigma$,
$\overline{\sigma}$, $\rho$, $\zeta$, $\overline{\zeta}$ is certainly
non-involutive and we have to start the {\it prolongation procedure}.

\subsection{Prolongation}
Now, we are ready to prolong the structure equation~\thetag{\ref{d3}}
of the non-involutive coframe
$\big\{\sigma,\overline{\sigma},\rho,\zeta,
\overline{\zeta} \big\}$\,\,---\,\,now, we rename
$\beta_1$ as $\alpha$\,\,---:
\begin{equation*}
\aligned d\sigma & =
\big(2\,\alpha+\overline{\alpha}\big) \wedge \sigma+
{\rho\wedge\zeta},
\\
d\rho & =
\big(\alpha+\overline{\alpha}\big)\wedge\rho +
i\,\zeta\wedge\overline{\zeta},
\\
 d\zeta & = \alpha\wedge\zeta.
\endaligned
\end{equation*}
We prolong the base manifold $M^5_{\sf c}$ to the
prolonged space $M^5_{\sf c}\times G^{\sf red}$,
where $G^{\sf red}$ is the reduced structure group:
\[\footnotesize\aligned
G^{\sf red}:=\left\{g=\left(\!\!
\begin{array}{ccccc}
{\sf a}\overline{\sf a}\,\overline{\sf a} & 0 & 0 & 0 & 0
\\
0 & {\sf a}{\sf a}\overline{\sf a} & 0 & 0 & 0
\\
0 & 0 & {\sf a}\overline{\sf a} & 0 & 0
\\
0 & 0 & 0 & \overline{\sf a} & 0
\\
0 & 0 & 0 & 0 & {\sf a}
\end{array}
\!\!\right)\colon
\ \ \ \ \ {\sf a}\in\mathbb C\right\}.
\endaligned
\]
The prolonged space $M^5_{\sf c}\times G^{\sf red}$
is in fact a submanifold of
the complex space $\mathbb C^5:=\mathbb C_{(z,w_1,w_2,w_3,{\sf
a})}$. The idea behind the prolongation procedure is based upon the
following fundamental proposition ({\it see} Proposition 12.1 page
375 of \cite{Olver-1995} for the general assertion):

\begin{Proposition}
\label{lifted-prop}
Let $\theta$ and $\theta'$ be two lifted coframes on
two manifolds $M$ and $M'$ having the same structure group $G$,
let $\alpha$ and $\alpha'$ be the modified Maurer-Cartan
forms obtained by solving the absorbtion equations and
assume that neither
group-dependent essential torsion coefficients exist
nor free absorption variables remain. Then there exists
a diffeomorphism $\Phi:M\rightarrow M'$ mapping $\theta$ to $\theta'$
for some choice of group parameters if and only if there is a
diffeomorphism $\Psi:M\times G\longrightarrow M'\times G$ mapping the
prolonged coframe $\{\theta,\alpha\}$ to $\{\theta',\alpha'\}$.
\qed
\end{Proposition}

Accordingly, one transforms the equivalence problem for
the 5-dimensional
cubic $M^5_{\sf c}$ into an equivalence problem on the prolonged
submanifold $M^5_{\sf c}\times G^{\sf red}$ equipped with the coframe
$\big\{ \sigma, \overline{\sigma}, \rho, \zeta, \overline{\zeta},
\alpha, \overline{\alpha} \big\}$, enlarged with the additional
$1$-form $\alpha=\frac{d\sf a}{\sf a}$. Since $d\alpha =
d\log{\sf a} =0$, we obtain gratuitously the following\,\,---\,\,fully
expressed\,\,---\,\,structure equations on this space:
\begin{equation}
\label{d4}
\aligned d\sigma & =
\big(2\,\alpha+\overline{\alpha}\big) \wedge \sigma+
{\rho\wedge\zeta},
\\
d\rho & = (\alpha+\overline{\alpha})\wedge\rho +
i\,\zeta\wedge\overline{\zeta},
\\
 d\zeta & = \alpha\wedge\zeta,
 \\
 d\alpha&=0.
\endaligned
\end{equation}
and this is nothing but the final,
desired $\{ e\}$-{\it structure}.

\smallskip

Now, recall that if
$\{{\sf v}_1,\ldots, {\sf v}_r\}$ is a basis of a real Lie algebra
$\frak g$ of left-invariant vector fields on
an $r$-dimensional Lie group $G$ which have Lie brackets:
\[
\big[{\sf v}_i,{\sf v}_j\big]
=
\sum_{k=1}^r\,c^k_{ij}\,{\sf v}_k
\ \ \ \ \ \ \ \ \ \ \ \ \
{\scriptstyle{(1\,\leqslant\,i\,<\,j\,\leqslant\,r)}},
\]
with certain {\sl structure constants}
$c^\bullet_{\bullet\bullet}$,
and if $\alpha^1,\ldots,\alpha^r$
is the dual Maurer-Cartan basis of left-invariant
1-forms, then their structure equations:
\[
d\alpha^k
=
-\sum_{1\leqslant i<j\leqslant r}\,
c^k_{ij}\,\alpha^i\wedge\alpha^j
\ \ \ \ \ \ \ \ \ \ \ \ \
{\scriptstyle{(k\,=\,1\,\cdots\,r)}}
\]
have the same structure coefficients $c^\bullet_{\bullet\bullet}$
up to an overall minus sign.

\smallskip

As in the case of a Lie group,
in our case~\thetag{ \ref{d4}}, the structure functions are constant,
hence the associated coframe is of
{\sl rank zero} ({\it see} \cite{Olver-1995}, pages 266--268).
To proceed the problem, now we need the following result.

\begin{Theorem}
\label{Olver-2}
(see \cite{Olver-1995}, page 268, Theorem 8.16). Let $\theta$ be a rank zero coframe on an $m$-dimensional manifold $M$, with constant structure functions $T^k_{ij}=-c^k_{ij}$. Let $G$ be an $m$-dimensional Lie group whose Lie algebra $\frak g$ has structure constants $c^k_{ij}$ relative to a basis $\{{\sf v}_1,\ldots,{\sf v}_m\}$, and let $\{\alpha^1,\ldots,\alpha^m\}$ denote the dual basis of Maurer-Cartan forms. Then, there exists a local diffeomorphism $\Phi:M\rightarrow G$ mapping the given coframe to the Maurer-Cartan coframe on $G$.
\end{Theorem}

Accordingly, every
generic submanifold ${ M'}^5 \subset \C^4$
which is equivalent to
$M^5_{\sf c}$ corresponds, after prolongation,
to a $7$-dimensional Lie group $G$ whose Lie
algebra $\frak g$ has the same structure constants as the structure
coefficients of~\thetag{\ref{d4}}. Let us examine this Lie algebra.

\smallskip

According to the reminder, taking account of the
overall minus sign,
we find the following complete table of commutators:

\medskip
\begin{center}
\label{table-5-restricted}
\begin{tabular}[t]{ l | c c c c c c c }
& ${\sf v}^{\sigma}$ & ${\sf v}^{\overline{\sigma}}$ & ${\sf
v}^{\rho}$ & ${\sf v}^{\zeta}$ & ${\sf v}^{\overline{\zeta}}$ & ${\sf
v}^{\alpha}$ & ${\sf v}^{\overline{\alpha}}$
\\
\hline ${\sf v}^{{\sigma}}$ & $0$ & $0$ & $0$ & $0$ & $0$ & $2\,{\sf
v}^{{\sigma}}$ & ${\sf v}^{{\sigma}}$
\\
${\sf v}^{\overline{\sigma}}$ & $*$ & $0$ & $0$ & $0$ & $0$ & ${\sf
v}^{\overline{\sigma}}$ & $2\,{\sf v}^{\overline{\sigma}}$
\\
${\sf v}^{\rho}$ & $*$ & $*$ & $0$ & $-{\sf v}^{{\sigma}}$ & $-{\sf
v}^{\overline{\sigma}}$ & ${\sf v}^{\rho}$ & ${\sf v}^{\rho}$
\\
${\sf v}^{\zeta}$ & $*$ & $*$ & $*$ & $0$ & $-i\,{\sf v}^{\rho}$ &
${\sf v}^{\zeta}$ & $0$
\\
${\sf v}^{\overline{\zeta}}$ & $*$ & $*$ & $*$ & $*$ & $0$ & $0$ &
${\sf v}^{\overline{\zeta}}$
\\
${\sf v}^{\alpha}$ & $*$ & $*$ & $*$ & $*$ & $*$ & $0$ & $0$
\\
${\sf v}^{\overline{\alpha}}$ & $*$ & $*$ & $*$ & $*$ & $*$ & $*$ & $0$.
\end{tabular}
\end{center}

\noindent
A visual inspection of this table shows that
this Lie algebra $\frak g$ is a $3$-{\sl graded algebra} in the sense of
Tanaka, of the form:
\[
\aligned
\frak g:=
&\,
\frak g_{-3}\oplus\frak g_{-2}\oplus\frak g_{-1}\oplus\frak
g_{0},
\\
&\,
[\frak g_i,\frak g_j]\subset\frak g_{i+j}
\ \ \ \ \ \ \ \ \ \ \ \ \
{\scriptstyle (i,\,j\,=\,-3,\,-2\,-1,\,0)},
\endaligned
\]
where:
\[
\aligned \frak g_{-3}&:=\langle {\sf v}^\sigma, {\sf
v}^{\overline{\sigma}}\rangle,
\\
\frak g_{-2}&:=\langle {\sf v}^\rho\rangle,
\\
\frak g_{-1}&:=\langle {\sf v}^\zeta, {\sf
v}^{\overline{\zeta}}\rangle,
\\
\mathfrak{g}_0&:=\langle {\sf v}^\alpha, {\sf
v}^{\overline{\alpha}}\rangle.
\endaligned
\]
But remember that
in Section~\ref{infinitesimal-model},
we computed the Lie algebra
$\frak{aut}_{CR}(M^5_{\sf c})$ of infinitesimal CR-automorphisms of
our cubic model $M^5_{\sf c}$ and there,
we showed that it is 7-dimensional of
the form:
\[
\frak {aut}_{CR}\big(M^5_{\sf c}\big)
:=
\frak a_{-3}\oplus\frak
a_{-2}\oplus\frak a_{-1}\oplus\frak a_{0},
\]
where:
\[
\aligned
\frak a_{-3}
&
:=
\langle S_1,S_2\rangle,
\\
\frak a_{-2}
&
:=
\langle T\rangle,
\\
\frak a_{-1}
&
:=
\langle L_1,L_2\rangle,
\\
\frak
a_{0}
&
:=
\langle D, R\rangle,
\endaligned
\]
and that it is equipped with the following table of
commutators:

\medskip
\begin{center}
\begin{tabular} [t] { c | c c c c c c c }
& $S_2$ & $S_1$ & $T$ & $L_2$ & $L_1$ & $D$ & $R$
\\
\hline $S_2$ & $0$ & $0$ & $0$ & $0$ & $0$ & $3S_2$ & $-S_1$
\\
$S_1$ & $*$ & $0$ & $0$ & $0$ & $0$ & $3S_1$ & $S_2$
\\
$T$ & $*$ & $*$ & $0$ & $4S_2$ & $4S_1$ & $2T$ & $0$
\\
$L_2$ & $*$ & $*$ & $*$ & $0$ & $-4T$ & $L_2$ & $-L_1$
\\
$L_1$ & $*$ & $*$ & $*$ & $*$ & $0$ & $L_1$ & $L_2$
\\
$D$ & $*$ & $*$ & $*$ & $*$ & $*$ & $0$ & $0$
\\
$R$ & $*$ & $*$ & $*$ & $*$ & $*$ & $*$ & $0$.
\end{tabular}
\end{center}

Introduce now the linear map:
\[
\Psi\colon\ \ \
\frak{aut}_{CR}(M^5_{\sf c})
\longrightarrow
\frak g,
\]
having the following values on the basis elements of $\frak
{aut}_{CR}(M^5_{\sf c})$:

\medskip
\begin{center}
\begin{tabular} [t] { l | l }
$\Psi$ & \ \ \ \ \ \ \ $\longmapsto$
\\
\hline $S_1$ &
$-\frac{i}{2}\,v^{\sigma}-\frac{i}{4}\,v^{\overline{\sigma}}$
\\
 $S_2$ &
$\frac{1}{2}\,v^{\sigma}-\frac{1}{4}\,v^{\overline{\sigma}}$
\\
$T$ & $v^{\rho}$
\\
$L_1$ & $2\,i\,v^{\zeta}+i\,v^{\overline{\zeta}}$
\\
$L_2$ & $-2\,v^{\zeta}+v^{\overline{\zeta}}$
\\
$D$ & $v^{\alpha}+v^{\overline{\alpha}}$
\\
$R$ & $i\,v^{\alpha}-i\,v^{\overline{\alpha}}$.
\end{tabular}
\end{center}
One checks that this linear map $\Psi$ is an isomorphism.
Consequently, the Lie algebra $\frak{aut}_{CR}(M^5_{\sf c})$ of
our cubic model $M^5_{\sf c}$ is isomorphic to the Lie algebra $\frak
g$ obtained after
prolongation.

\begin{Theorem}
\label{theorem-cubic}
A 5-dimensional maximally minimal real analytic CR-generic local
submanifold ${M'}^5\subset\mathbb C^4$ is equivalent, through some
local biholomorphic transformation, to the cubic model $M^5_{\sf c}
\subset \C^4$, if and only if the structure equations of the lifted
coframe $\{\sigma,\overline\sigma,\rho,\zeta,\overline\zeta,\alpha,\overline\alpha\}$
associated to the 7-dimensional prolonged space $M^5\times G^{\sf red}$ are in the form \thetag{\ref{d4}}, if and only if they
have structure constants of a Lie algebra
isomorphic to $\mathfrak{ aut}_{ CR} \big( M_{\sf c}^5 \big)$
 where $G^{\sf red}$ is the 2-dimensional matrix group:
 \[\footnotesize\aligned
G^{\sf red}:=\left\{g=\left(\!\!
\begin{array}{ccccc}
{\sf a}\overline{\sf a}\,\overline{\sf a} & 0 & 0 & 0 & 0
\\
0 & {\sf a}{\sf a}\overline{\sf a} & 0 & 0 & 0
\\
0 & 0 & {\sf a}\overline{\sf a} & 0 & 0
\\
0 & 0 & 0 & \overline{\sf a} & 0
\\
0 & 0 & 0 & 0 & {\sf a}
\end{array}
\!\!\right)\colon
\ \ \ \ \ {\sf a}\in\mathbb C\right\},
\endaligned
\]
 and, if and only if the prolonged space $M^5\times G^{\sf red}$ is equivalent, through an $\{e\}$-structure to ${\sf Aut}_{CR}(M^5_{\sf c})$, the symmetry group of the model $M^5_{\sf c}$ with the associated Lie algebra $\mathfrak{ aut}_{ CR} \big( M_{\sf c}^5 \big)$.
\end{Theorem}

\proof
As we already observed, the structure constants of the structure equations \thetag{\ref{d4}} are precisely equal to those of the Lie algebra $\mathfrak{ aut}_{ CR} \big( M_{\sf c}^5 \big)$. Furthermore, the Lie algebra $\mathfrak{ aut}_{ CR} \big( M_{\sf c}^5 \big)$ associates to the Lie group ${\sf Aut}_{ CR} \big( M_{\sf c}^5 \big)$. Now, the result comes immediately from Theorem \ref{Olver-2}.
\endproof

One should notice that actually the problem of equivalency to our cubic model $M^5_{\sf c}$ is not fully-completed, yet. Being more precise and according to the first part of the above theorem, to realize the equivalency of an arbitrary 5-dimensional CR-manifold ${M'}^5$ to this model, we need to have the structure equations of the lifted coframe, associated to it. In Corollary \ref{complete-model} we will finalize the solution of the current equivalence problem.

\begin{Remark}
Let us denot by $G_-$ the Lie subgroup of ${\sf Aut}_{CR}(M^5_{\sf c})$ associated to the Levi-Tanaka subalgebra $\frak g_-$ of $\mathfrak{ aut}_{ CR} \big( M_{\sf c}^5 \big)$. It is known ({\it see} \cite{Beloshapka-2004}, Proposition 3) that this group is diffeomorphic to the model $M^5_{\sf c}$. Then, the above theorem asserts that the {\it prolonged} space $M^5_{\sf c}\times G^{\sf red}$ is equivalent to the Tanaka {\it prolongation} of the Lie subgroup $G_-$, namely ${\sf Aut}_{CR}(M^5_{\sf c})$. This may show the close coherency between two concepts of prolongation in the senses of Cartan and Tanaka.
\end{Remark}

\begin{Corollary}
\label{corollary-cubic}
Let $M^5$ be an under consideration 5-dimensional CR-manifold which is equivalent, through some local biholomorphism to the model $M^5_{\sf c}$. Then, the symmetry group of the corresponding 7-dimensional prolonged space $M^5\times G^{\sf red}$ is ${\sf Aut}_{CR}(M^5_{\sf c})$.
\end{Corollary}

\proof
The assertion is a direct consequence of Corollary 14.20 page 435 of \cite{Olver-1995}.
\endproof

\section{Initial complex frame
\\
for geometry-preserving deformations of the model}
\label{initial-complex-frame}
\HEAD{\ref{initial-complex-frame}.~
Initial complex frame
for geometry-preserving deformations of the model}{
Masoud {\sc Sabzevari} (Shahrekord) and Jo\"el {\sc Merker} (LM-Orsay)}

\subsection{Initial Lie bracket structure}
Given a 3-codimensional CR-generic submanifold $M^5 \subset \C^4$
represented as a graph:
\[
\aligned v_1 & = \varphi_1\big(x,y,u_1,u_2,u_3\big),
\\
v_2 & = \varphi_2\big(x,y,u_1,u_2,u_3\big),
\\
v_3 & = \varphi_3\big(x,y,u_1,u_2,u_3\big),
\endaligned
\]
with $\varphi_j ( 0) = 0$ and $d \varphi_j ( 0) = 0$, the $(0,
1)$-complex tangent bundle $T^{0,1}  M^5$ is spanned by a
certain single $(0,1)$-vector field of the form:
\[
\overline{\mathcal{L}} = \frac{\partial}{\partial\overline{ z}} +
{\sf A}_1\,\frac{\partial}{\partial \overline{w}_1} + {\sf
A}_2\,\frac{\partial}{\partial \overline{w}_2} + {\sf
A}_3\,\frac{\partial}{\partial\overline{ w}_3},
\]
where the coefficient-functions ${\sf A}_1$, ${\sf A}_2$, ${\sf A}_3$
may be expressed explicitly in terms of the first-order jets $J_{ x,
y, u_1, u_2, u_3}^1 \varphi_j$ of the functions $\varphi_j$. Writing
the tangency equations gives that its three (unknown) coefficients
${\sf A}_1$, ${\sf A}_2$, ${\sf A}_3$ should satisfy the three
equations:
\[
\footnotesize\aligned \frac{i}{2}{\sf A}_1 & =
\varphi_{1\overline{z}}+ \frac{1}{2}{\sf
A}_1\varphi_{1u_1}+\frac{1}{2}{\sf A}_2\varphi_{1u_2}+\frac{1}{2}{\sf
A}_3\varphi_{1u_3},
\\
\frac{i}{2}{\sf A}_2 & = \varphi_{2\overline{z}}+ \frac{1}{2}{\sf
A}_1\varphi_{2u_1}+\frac{1}{2}{\sf A}_2\varphi_{2u_2}+\frac{1}{2}{\sf
A}_3\varphi_{2u_3},
\\
\frac{i}{2}{\sf A}_3 & = \varphi_{3\overline{z}}+ \frac{1}{2}{\sf
A}_1\varphi_{3u_1}+\frac{1}{2}{\sf A}_2\varphi_{3u_2}+\frac{1}{2}{\sf
A}_3\varphi_{3u_3}.
\endaligned
\]
Equivalently, these equations can be read as a non-homogeneous system
with right hand-side vanishing at the origin:
\[
\footnotesize\aligned -\frac{1}{2i}{\sf A}_1&=
\frac{1}{2}\varphi_{1x}-\frac{1}{2i}\varphi_{1y}+\frac{1}{2}{\sf
A}_1\varphi_{1u_1}+\frac{1}{2}{\sf A}_2\varphi_{1u_2}+\frac{1}{2}{\sf
A}_3\varphi_{1u_3},
\\
-\frac{1}{2i}{\sf A}_2&=
\frac{1}{2}\varphi_{2x}-\frac{1}{2i}\varphi_{2y}+\frac{1}{2}{\sf
A}_1\varphi_{2u_1}+\frac{1}{2}{\sf A}_2\varphi_{2u_2}+\frac{1}{2}{\sf
A}_3\varphi_{2u_3},
\\
-\frac{1}{2i}{\sf A}_3&=
\frac{1}{2}\varphi_{3x}-\frac{1}{2i}\varphi_{2y}+\frac{1}{2}{\sf
A}_1\varphi_{3u_1}+\frac{1}{2}{\sf A}_2\varphi_{3u_2}+\frac{1}{2}{\sf
A}_3\varphi_{3u_3}.
\endaligned
\]
Solving this linear system of three equations and three unknowns
${\sf A_1},{\sf A}_2$, ${\sf A}_3$ and expanding the solutions
according to their real and imaginary parts provides the rational
solutions:
\[
\footnotesize \aligned {\sf A}_1&=
\frac{\Lambda^1_1}{\Delta}+i\frac{\Lambda^1_2}{\Delta},
\\
{\sf A}_2&= \frac{\Lambda^2_1}{\Delta}+i\frac{\Lambda^2_2}{\Delta},
\\
{\sf A}_3&= \frac{\Lambda^3_1}{\Delta}+i\frac{\Lambda^3_2}{\Delta},
\endaligned
\]
in which the denominator $\Delta$, in terms of the first-order
derivatives of $\varphi_1$, $\varphi_2$, $\varphi_3$, is:
\[
\aligned \Delta&=\sigma^2+\tau^2,
\endaligned
\]
with the squared functions of:
\[
\footnotesize \aligned \sigma&=
\varphi_{3u_3}+\varphi_{1u_1}+\varphi_{2u_2}-\varphi_{1u_2}\varphi_{3u_1}\varphi_{2u_3}-
\varphi_{1u_3}\varphi_{2u_1}\varphi_{3u_2}+\varphi_{1u_2}\varphi_{2u_1}\varphi_{3u_3}-
\\
&-
\varphi_{1u_1}\varphi_{2u_2}\varphi_{3u_3}+\varphi_{1u_1}\varphi_{2u_3}\varphi_{3u_2}+\varphi_{1u_3}\varphi_{3u_1}\varphi_{2u_2},
\\
 \tau & =
-1+\varphi_{1u_1}\varphi_{2u_2}-\varphi_{2u_3}\varphi_{3u_2}-\varphi_{1u_3}\varphi_{3u_1}+
\varphi_{2u_2}\varphi_{3u_3}-\varphi_{1u_2}\varphi_{2u_1}+\varphi_{1u_1}\varphi_{3u_3},
\endaligned
\]
and in which the numerators $\Lambda^\bullet_\bullet$ of ${\sf A}_1,
{\sf A}_2$ and ${\sf A}_3$ are in the (longer) forms:
\[
\footnotesize \aligned \Lambda^1_1&=
\big(-\varphi_{3u_3}\varphi_{2x}\varphi_{1u_2}-\varphi_{1u_3}\varphi_{3y}+\varphi_{2u_2}\varphi_{1x}\varphi_{3u_3}+
\varphi_{3u_3}\varphi_{1y}-\varphi_{1x}-\varphi_{2y}\varphi_{1u_2}+
\\
&+\varphi_{2u_3}\varphi_{3x}\varphi_{1u_2}+\varphi_{2u_2}\varphi_{1y}-\varphi_{2u_3}\varphi_{3u_2}\varphi_{1x}-\varphi_{2u_2}\varphi_{1u_3}\varphi_{3x}+
\varphi_{2x}\varphi_{1u_3}\varphi_{3u_2}\big)\sigma +
\\
&+
\big(\varphi_{1u_3}\varphi_{3x}-\varphi_{1y}+\varphi_{2x}\varphi_{1u_2}+\varphi_{2u_3}\varphi_{1u_2}\varphi_{3y}-\varphi_{2u_2}\varphi_{1x}-
\varphi_{2u_3}\varphi_{3u_2}\varphi_{1y}-
\\
&-\varphi_{3u_3}\varphi_{1x}-\varphi_{2u_2}\varphi_{1u_3}\varphi_{3y}-\varphi_{3u_3}\varphi_{1u_2}\varphi_{2y}+\varphi_{1u_3}\varphi_{3u_2}\varphi_{2y}+
\varphi_{2u_2}\varphi_{3u_3}\varphi_{1y} \big)\tau,
\endaligned
\]
\[
\footnotesize \aligned \Lambda^1_2&=
\big(\varphi_{1u_3}\varphi_{3x}-\varphi_{1y}+\varphi_{2x}\varphi_{1u_2}+\varphi_{2u_3}\varphi_{1u_2}\varphi_{3y}-\varphi_{2u_2}\varphi_{1x}-
\varphi_{2u_3}\varphi_{3u_2}\varphi_{1y}-
\\
&-\varphi_{3u_3}\varphi_{1x}-\varphi_{2u_2}\varphi_{1u_3}\varphi_{3y}-
\varphi_{3u_3}\varphi_{1u_2}\varphi_{2y}+
\varphi_{1u_3}\varphi_{3u_2}\varphi_{2y}+\varphi_{2u_2}\varphi_{3u_3}\varphi_{1y}\big)\sigma
-
\\
&-
\big(-\varphi_{3u_3}\varphi_{2x}\varphi_{1u_2}-\varphi_{1u_3}\varphi_{3y}+\varphi_{2u_2}\varphi_{1x}\varphi_{3u_3}+\varphi_{3u_3}\varphi_{1y}-
\varphi_{1x}-\varphi_{2y}\varphi_{1u_2}+
\\
&+\varphi_{2u_3}\varphi_{3x}\varphi_{1u_2}+\varphi_{2u_2}\varphi_{1y}-
\varphi_{2u_3}\varphi_{3u_2}\varphi_{1x}-
\varphi_{2u_2}\varphi_{1u_3}\varphi_{3x}+\varphi_{2x}\varphi_{1u_3}\varphi_{3u_2}\big)\tau,
\endaligned
\]
\[
\footnotesize \aligned \Lambda^2_1&=
\big(-\varphi_{2x}+\varphi_{3u_3}\varphi_{2y}+\varphi_{1u_3}\varphi_{2u_1}\varphi_{3x}-\varphi_{2u_3}\varphi_{3y}-
\varphi_{1u_3}\varphi_{3u_1}\varphi_{2x}- \varphi_{2u_1}\varphi_{1y}-
\\
&-\varphi_{2u_1}\varphi_{3u_3}\varphi_{1x}+\varphi_{1u_1}\varphi_{2y}-\varphi_{1u_1}\varphi_{2u_3}\varphi_{3x}+
\varphi_{3u_1}\varphi_{2u_3}\varphi_{1x}+\varphi_{1u_1}\varphi_{3u_3}\varphi_{2x}\big)\sigma+
\\
&
+\big(-\varphi_{1u_1}\varphi_{2u_3}\varphi_{3y}+\varphi_{1u_3}\varphi_{2u_1}\varphi_{3y}-\varphi_{1u_3}\varphi_{3u_1}\varphi_{2y}+
\varphi_{3u_1}\varphi_{2u_3}\varphi_{1y}+\varphi_{2u_3}\varphi_{3x}-
\\
&-\varphi_{2u_1}\varphi_{3u_3}\varphi_{1y}-\varphi_{3u_3}\varphi_{2x}+
\varphi_{1u_1}\varphi_{3u_3}\varphi_{2y}+\varphi_{2u_1}\varphi_{1x}-\varphi_{1u_1}\varphi_{2x}-\varphi_{2y}\big)\tau,
\endaligned
\]
\[
\footnotesize \aligned \Lambda^2_2&=
\big(-\varphi_{1u_1}\varphi_{2u_3}\varphi_{3y}+\varphi_{1u_3}\varphi_{2u_1}\varphi_{3y}-\varphi_{1u_3}\varphi_{3u_1}\varphi_{2y}+
\varphi_{3u_1}\varphi_{2u_3}\varphi_{1y}+\varphi_{2u_3}\varphi_{3x}-
\\
&-\varphi_{2u_1}\varphi_{3u_3}\varphi_{1y}-\varphi_{3u_3}\varphi_{2x}+
\varphi_{1u_1}\varphi_{3u_3}\varphi_{2y}+\varphi_{2u_1}\varphi_{1x}-\varphi_{1u_1}\varphi_{2x}-\varphi_{2y}\big)\sigma-
\\
&-
\big(-\varphi_{2x}+\varphi_{3u_3}\varphi_{2y}+\varphi_{1u_3}\varphi_{2u_1}\varphi_{3x}-\varphi_{2u_3}\varphi_{3y}-
\varphi_{1u_3}\varphi_{3u_1}\varphi_{2x}-\varphi_{2u_1}\varphi_{1y}-
\\
&-
\varphi_{2u_1}\varphi_{3u_3}\varphi_{1x}+\varphi_{1u_1}\varphi_{2y}-
\varphi_{1u_1}\varphi_{2u_3}\varphi_{3x}+\varphi_{3u_1}\varphi_{2u_3}\varphi_{1x}+\varphi_{1u_1}\varphi_{3u_3}\varphi_{2x}\big)\tau,
\endaligned
\]
\[
\footnotesize \aligned \Lambda^3_1&=
\big(-\varphi_{2u_1}\varphi_{1u_2}\varphi_{3x}-\varphi_{3u_1}\varphi_{1y}+\varphi_{2u_1}\varphi_{3u_2}\varphi_{1x}+
\varphi_{1u_1}\varphi_{3y}-\varphi_{3u_1}\varphi_{2u_2}\varphi_{1x}-
\\
&-\varphi_{3u_2}\varphi_{2y}+\varphi_{3u_1}\varphi_{1u_2}\varphi_{2x}+
\varphi_{2u_2}\varphi_{3y}-\varphi_{3u_2}\varphi_{1u_1}\varphi_{2x}+\varphi_{1u_1}\varphi_{2u_2}\varphi_{3x}-\varphi_{3x}\big)\sigma+
\\
&+
\big(-\varphi_{3u_2}\varphi_{1u_1}\varphi_{2y}+\varphi_{2u_1}\varphi_{3u_2}\varphi_{1y}+\varphi_{3u_1}\varphi_{1x}+
\varphi_{3u_1}\varphi_{1u_2}\varphi_{2y}-\varphi_{2u_1}\varphi_{1u_2}\varphi_{3y}+
\\
&+\varphi_{3u_2}\varphi_{2x}-\varphi_{1u_1}\varphi_{3x}-
\varphi_{3u_1}\varphi_{2u_2}\varphi_{1y}+\varphi_{1u_1}\varphi_{2u_2}\varphi_{3y}-\varphi_{3y}-\varphi_{2u_2}\varphi_{3x}\big)\tau,
\endaligned
\]
\[
\footnotesize \aligned \Lambda^3_2&=
\big(-\varphi_{3u_2}\varphi_{1u_1}\varphi_{2y}+\varphi_{2u_1}\varphi_{3u_2}\varphi_{1y}+\varphi_{3u_1}\varphi_{1x}+
\varphi_{3u_1}\varphi_{1u_2}\varphi_{2y}-\varphi_{2u_1}\varphi_{1u_2}\varphi_{3y}+
\\
&+ \varphi_{3u_2}\varphi_{2x}-\varphi_{1u_1}\varphi_{3x}-
\varphi_{3u_1}\varphi_{2u_2}\varphi_{1y}+\varphi_{1u_1}\varphi_{2u_2}\varphi_{3y}-\varphi_{3y}-\varphi_{2u_2}\varphi_{3x}\big)\sigma-
\\
&-
\big(-\varphi_{2u_1}\varphi_{1u_2}\varphi_{3x}-\varphi_{3u_1}\varphi_{1y}+\varphi_{2u_1}\varphi_{3u_2}\varphi_{1x}+
\varphi_{1u_1}\varphi_{3y}-\varphi_{3u_1}\varphi_{2u_2}\varphi_{1x}-
\\
&-\varphi_{3u_2}\varphi_{2y}+\varphi_{3u_1}\varphi_{1u_2}\varphi_{2x}+
\varphi_{2u_2}\varphi_{3y}-\varphi_{3u_2}\varphi_{1u_1}\varphi_{2x}+\varphi_{1u_1}\varphi_{2u_2}\varphi_{3x}-\varphi_{3x}\big)\tau,
\endaligned
\]
again in terms of the first-order derivatives of $\varphi_1$,
$\varphi_2$, $\varphi_3$.

Here, the expression of $\overline{\mathcal L}$ is presented as a
vector field which lives in a neighborhood of $M^5$ in $\mathbb C^4$,
while $M^5$ itself, is a real five-dimensional hypersurface equipped
with the five real coordinates $x,y,u_1,u_2,u_3$. But, in order to
express $\overline{\mathcal L}$ {\it intrinsically}, one must drop
$\frac{\partial}{\partial v_1},\,\frac{\partial}{\partial v_2}$ and
$\frac{\partial}{\partial v_3}$ and also replace $v_i$ by
$\varphi_i(x,y,u_1,u_2,u_3)$ for $i=1,2,3$, simultaneously in its
expression. Then, after expanding $\overline{\mathcal L}$ in real and
imaginary parts:
\[
\aligned \overline{\mathcal L} & = \frac{\partial}{\partial\overline{
z}}+ \frac{{\sf A}_1}{2}\frac{\partial}{\partial u_1}+ \frac{{\sf
A}_2}{2}\frac{\partial}{\partial u_2}+ \frac{{\sf
A}_3}{2}\frac{\partial}{\partial u_3},
\endaligned
\]
one gains a result that can now be summarized as follows.

\begin{Proposition}
For any local real analytic CR-generic submanifold $M^5 \subset
\C^4$ which is represented near the
origin as a graph:
\[
\left[
\aligned v_1&:= \varphi_1(x,y,u_1,u_2,u_3),
\\
v_2&:= \varphi_2(x,y,u_1,u_2,u_3),
\\
v_3&:= \varphi_3(x,y,u_1,u_2,u_3),
 \endaligned\right.
\]
in coordinates:
\[
\big(z,w_1,w_2,w_3\big) = \big(x+iy,u_1+iv_1,u_2+iv_2,u_3+iv_3\big),
\]
its complex bundle $T^{1,0}M$
is generated by:
\[
\boxed{ \aligned \overline{\mathcal L} & = \frac{\partial}{\partial
\overline{z}}+ \frac{{\sf A}_1}{2}\frac{\partial}{\partial u_1}+
\frac{{\sf A}_2}{2}\frac{\partial}{\partial u_2}+ \frac{{\sf
A}_3}{2}\frac{\partial}{\partial u_3},
\endaligned}
\]
whose numerators and denominators have the explicit expressions shown
above in terms of the first order jets
$J_{x,y,u_1,u_2,u_3}^1\varphi_j, \, j=1,2,3$.
\end{Proposition}

In particular, for the cubic model $M^5_{\sf c}\subset\mathbb C^4$,
represented as the graph:
\[
\left[
\aligned v_1&:= x^2+y^2,
\\
v_2&:= 2\,x^3+2\,xy^2,
\\
v_3&:= 2\,x^2y+2\,y^3,
\endaligned\right.
\]
 we have:
\[\footnotesize\aligned
\overline{\mathcal L}_{\sf
c}&=\frac{\partial}{\partial\overline{z}}+\big(y-ix\big)\,\frac{\partial}{\partial
u_1} + \big(2xy-i(3x^2+y^2)\big)\,\frac{\partial}{\partial u_2} +
\big(x^2+3y^2-2ixy\big)\,\frac{\partial}{\partial u_3}
\\
{\mathcal L}_{\sf
c}&=\frac{\partial}{\partial{z}}+\big(y+ix\big)\,\frac{\partial}{\partial
u_1} + \big(2xy+i(3x^2+y^2)\big)\,\frac{\partial}{\partial u_2} +
\big(x^2+3y^2+2ixy\big)\,\frac{\partial}{\partial u_3}.
\endaligned
\]

\subsection{Length-two Lie bracket}
Between the two already presented complex vector fields $\mathcal L$
and $\overline{\mathcal L}$, of course there is only one Lie bracket
$[\mathcal L,\overline{\mathcal L}]$ of length two. This vector field
is in fact {\it imaginary}. Being more precise, expressing
$\overline{\mathcal L}:=L_1+iL_2$ in real and imaginary parts, then
the commutator:
\begin{eqnarray*}
\big[\mathcal L,\overline{\mathcal L}\big]&=&
\big[L_1-iL_2,L_1+iL_2\big],
\\
&=&2i\,\big[L_1,L_2\big]
\end{eqnarray*}
is imaginary. Hence let us denote the third
({\it real}) vector field by:
\[
\mathcal T
:
=i\,\big[\mathcal L,\overline{\mathcal L}\big].
\]
Performing direct but painful computations provides the expression
of:
\[
\boxed{ \mathcal T =
\frac{\Upsilon_1}{\Delta^3}\frac{\partial}{\partial u_1}+
\frac{\Upsilon_2}{\Delta^3}\frac{\partial}{\partial u_2} +
\frac{\Upsilon_3}{\Delta^3}\frac{\partial}{\partial u_3},}
\]
in which the three new numerators $\Upsilon_i$ are given as follows:
\begin{equation*}
\footnotesize\aligned \Upsilon_1 & =
 -\big(\Delta^2\Lambda^1_{2x}-\Delta\Delta_x\Lambda^1_2-
\Delta^2\Lambda^1_{1y}+\Delta\Delta_y\Lambda^1_1
+\Delta\Lambda^1_1\Lambda^1_{2u_1} -
\Delta\Lambda^1_2\Lambda^1_{1u_1}- \Delta\Lambda^2_2\Lambda^1_{1u_2}+
\\
&+ \Delta_{u_2}\Lambda^1_1\Lambda^2_2-
\Delta\Lambda^3_2\Lambda^1_{1u_3}+
\Delta_{u_3}\Lambda^3_2\Lambda^1_1+
\Delta\Lambda^2_1\Lambda^1_{2u_2}-
\Delta_{u_2}\Lambda^2_1\Lambda^1_{2}+
\Delta\Lambda^3_1\Lambda^1_{2u_3}-
\Delta_{u_3}\Lambda^3_1\Lambda^1_{2}\big),
\\
\Upsilon_2 & =
 - \big(\Delta^2\Lambda^2_{2x}- \Delta\Delta_x\Lambda^2_2+
\Delta\Lambda^1_{1}\Lambda^2_{2u_1}-
\Delta_{u_1}\Lambda^1_1\Lambda^2_2 -
\Delta^2\Lambda^2_{1y}+\Delta\Delta_y\Lambda^2_1-
\Delta\Lambda^1_2\Lambda^2_{1u_1}+
\\
&+
\Delta_{u_1}\Lambda^1_2\Lambda^2_1+\Delta\Lambda^2_1\Lambda^2_{2u_2}-
\Delta\Lambda^2_2\Lambda^2_{1u_2}+\Delta\Lambda^3_1\Lambda^2_{2u_3}-
\Delta_{u_3}\Lambda^3_1\Lambda^2_{2}-
\Delta\Lambda^3_2\Lambda^2_{1u_3}+
\Delta_{u_3}\Lambda^3_2\Lambda^2_{1}\big),
\\
\Upsilon_3 & = -
\big(\Delta^2\Lambda^3_{2x}-\Delta\Delta_x\Lambda^3_2+
\Delta\Lambda^1_{1}\Lambda^3_{2u_1} -
\Delta_{u_1}\Lambda^1_1\Lambda^3_2
-\Delta^2\Lambda^3_{1y}+\Delta\Delta_y\Lambda^3_1-
\Delta\Lambda^1_2\Lambda^3_{1u_1}+
\\
&+
\Delta_{u_1}\Lambda^1_2\Lambda^3_1-\Delta\Lambda^2_2\Lambda^3_{1u_2}+
\Delta_{u_2}\Lambda^2_2\Lambda^3_{1} +
\Delta\Lambda^3_1\Lambda^3_{2u_3}-
\Delta\Lambda^3_2\Lambda^3_{1u_3}+\Delta\Lambda^2_1\Lambda^3_{2u_2}-
\Delta_{u_2}\Lambda^2_1\Lambda^3_{2}\big).
\endaligned
\end{equation*}

In particular, for the cubic model $M_{\sf c}^5 \subset \C^3$, we
have:
\[
\mathcal T_{\sf c} = 4\,\frac{\partial}{\partial u_1} +
16x\,\frac{\partial}{\partial u_2} +16y\,\frac{\partial}{\partial
u_3}.
\]

\subsection{Length-three Lie brackets}
In this length, we have two Lie brackets:
\[
\mathcal S:=[\mathcal L,\mathcal T]
\ \ \ \ \ \ \ \ \ \ \
\text{\rm and}
\ \ \ \ \ \ \ \ \ \ \
\overline{\mathcal
S}:=[\overline{\mathcal L},\mathcal T].
\]
According to the performed computations, we have the explicit
expressions of these two {\it complex} vector fields in terms of the
defining functions $\varphi_1,\varphi_2,\varphi_3$ as:
\[
\aligned \mathcal S &=&
\frac{\Gamma^1_1-i\Gamma^1_2}{\Delta^5}\frac{\partial}{\partial
u_1}+\frac{\Gamma^2_1-i\Gamma^2_2}{\Delta^5}\frac{\partial}{\partial
u_2}+\frac{\Gamma^3_1-i\Gamma^3_2}{\Delta^5}\frac{\partial}{\partial
u_3},
\endaligned
\]
where, by allowing the two notational coincidences $x\equiv x_1$ and
$y\equiv x_2$, the numerators are (for $i=1,2$):
\begin{equation*}
\footnotesize\aligned \Gamma^1_i & = -2 \big({\textstyle\frac{1}{4}}
\Delta^2\Upsilon_{1{x_i}}- 3\Delta\Delta_{x_i}\Upsilon_1+
\Delta\Lambda^1_i\Upsilon_{1u_1}-
2\Delta_{u_1}\Lambda^1_i\Upsilon_{1}-
\Delta\Lambda^1_{iu_1}\Upsilon_{1}-
\Delta\Lambda^1_{iu_2}\Upsilon_{2}+
\\
&+ \Delta_{u_2}\Lambda^1_i\Upsilon_{2}-
\Delta\Lambda^1_{iu_3}\Upsilon_{3}+
\Delta_{u_3}\Lambda^1_i\Upsilon_{3}+
\Delta\Lambda^2_i\Upsilon_{1u_2}-
3\Delta_{u_2}\Lambda^2_i\Upsilon_{1}+
\Delta\Lambda^3_i\Upsilon_{1u_3}-
3\Delta_{u_3}\Lambda^3_i\Upsilon_{1}\big),
\\
\Gamma^2_i & = -2 \big( \Delta^2\Upsilon_{2{x_i}}-
3\Delta\Delta_{x_i}\Upsilon_2+ \Delta\Lambda^1_i\Upsilon_{2u_1}-
3\Delta_{u_1}\Lambda^1_i\Upsilon_{2}-
\Delta\Lambda^2_{iu_1}\Upsilon_{1}+
\Delta_{u_1}\Lambda^2_{i}\Upsilon_{1}+
\\
&+ \Delta\Lambda^2_i\Upsilon_{2u_2}-
2\Delta_{u_2}\Lambda^2_{i}\Upsilon_{2}-
\Delta\Lambda^2_{iu_2}\Upsilon_{2}-
\Delta\Lambda^2_{iu_3}\Upsilon_{3}+
\Delta_{u_3}\Lambda^2_i\Upsilon_{3}+
\Delta\Lambda^3_i\Upsilon_{2u_3}-
3\Delta_{u_3}\Lambda^3_i\Upsilon_{2}\big),
\\
\Gamma^3_i & = -2\big( \Delta^2\Upsilon_{3{x_i}}-
3\Delta\Delta_{x_i}\Upsilon_3+ \Delta\Lambda^1_i\Upsilon_{3u_1}-
3\Delta_{u_1}\Lambda^1_i\Upsilon_{3}+
\Delta\Lambda^2_{i}\Upsilon_{3u_2}-
3\Delta_{u_2}\Lambda^2_{i}\Upsilon_{3}-
\\
&- \Delta\Lambda^3_{iu_1}\Upsilon_{1}+
\Delta_{u_1}\Lambda^3_{i}\Upsilon_{1}-
\Delta\Lambda^3_{iu_2}\Upsilon_{2}+
\Delta_{u_2}\Lambda^3_{i}\Upsilon_{2}+
\Delta\Lambda^3_i\Upsilon_{3u_3}-
2\Delta_{u_3}\Lambda^3_i\Upsilon_{3}-
\Delta\Lambda^3_{iu_3}\Upsilon_{3}\big).
\endaligned
\end{equation*}

In particular for the cubic model $M^5_{\sf c}\subset\mathbb C^4$,
the above expressions give:
\[
\aligned
\mathcal S_{\sf c}&=8\frac{\partial}{\partial
u_2}-8i\frac{\partial}{\partial u_3},
\\
\overline{\mathcal S}_{\sf c}&=8\frac{\partial}{\partial
u_2}+8i\frac{\partial}{\partial u_3}.
\endaligned
\]

\begin{Proposition}
\label{frame-5-general}
The five vector fields $\overline{\mathcal L},{\mathcal L},\mathcal
T,\overline{\mathcal S},\mathcal S$ constitute a (complex) frame for
$TM^5\otimes_{\mathbb R}\mathbb C$.
\end{Proposition}

\proof Thanks to assumption that remainders in the graphing functions
$\varphi_1$, $\varphi_2$ and $\varphi_3$ are all an ${\sf OW} ( 4)$,
the values of these vector fields at the origin are the same as for
the corresponding ones of the model $M_{\sf c}^5 \subset \C^4$,
namely:
\[\footnotesize
\aligned \overline{\mathcal L}\big\vert_0 & = \overline{\mathcal
L}_{{\sf c}}\big\vert_0 = \frac{\partial}{\partial\overline{z}}, \ \
\ \ \ \ \ \mathcal L\big\vert_0 = \mathcal L_{{\sf c}}\big\vert_0 =
\frac{\partial}{\partial z},
\\
\mathcal T\big\vert_0 & = \mathcal T_{{\sf c}}\big\vert_0 =
4\,\frac{\partial}{\partial u_1},
\\
\overline{\mathcal S}|_0 & = \overline{\mathcal S}_{{\sf
c}}\big\vert_0 = \,8\frac{\partial}{\partial
u_2}+8i\frac{\partial}{\partial u_3}, \ \ \ \ \ \ \ \mathcal
S\big\vert_0 = \mathcal S_{{\sf c}}\big\vert_0 =
\,8\frac{\partial}{\partial u_2}-8i\frac{\partial}{\partial u_3},
\endaligned
\]
whence the five vectors $\overline{\mathcal L},{\mathcal L},\mathcal
T,\overline{\mathcal S},\mathcal S$ are linearly independent in a
neighborhood of the origin and they constitute a local frame for
$M^5$.
\endproof

\subsection{Other iterated Lie brackets}
We saw that the collection of five vector fields:
\[
\Big\{ \overline{\mathcal{S}},\, \mathcal{S},\, \mathcal{T},\,
\overline{\mathcal{L}},\, \mathcal{L} \Big\},
\]
where:
\[
\aligned \mathcal{T} & :=
i\big[\mathcal{L},\overline{\mathcal{L}}\big],
\\
\mathcal{S} & := \big[\mathcal{L},\,\mathcal{T}\big],
\\
\overline{\mathcal{S}} & :=
\big[\overline{\mathcal{L}},\,\mathcal{T}\big],
\endaligned
\]
makes up a frame for $\C \otimes TM^5$. Having five fields implies
that there are in sum ten Lie brackets between them. Thus, there
remain seven such brackets to be looked at.

Let us start with the following group of four Lie brackets:
\[
\big[\mathcal{L},\,\mathcal{S}\big], \ \ \ \ \ \ \ \ \
\big[\overline{\mathcal{L}},\,\mathcal{S}\big], \ \ \ \ \ \ \ \ \
\big[\mathcal{L},\,\overline{\mathcal{S}}\big], \ \ \ \ \ \ \ \ \
\big[\overline{\mathcal{L}},\,\overline{\mathcal{S}}\big].
\]
Because of the coefficient $\frac{ \partial}{
\partial z}$ in $\mathcal{ L}$ is 1, we observe that each one of these four Lie
brackets is a linear combination of just $\frac{ \partial}{ \partial
u_1}$, $\frac{ \partial}{ \partial u_2}$, $\frac{ \partial}{ \partial
u_3}$. For the same reason, we point out that $\mathcal{ T}$,
$\mathcal{ S}$ and $\overline{ \mathcal{S}}$ were also already a
linear combination of just $\frac{ \partial}{ \partial u_1}$, $\frac{
\partial}{ \partial u_2}$, $\frac{ \partial}{ \partial u_3}$. Thus,
we are sure that there are six complex-valued functions $P$, $Q$,
$R$, $A$, $B$, $C$ defined on $M^5$ such that:
\[
\aligned \big[\mathcal{L},\,\mathcal{S}\big] & = P\,\mathcal{T} +
Q\,\mathcal{S} + R\,\overline{\mathcal{S}},
\\
\big[\overline{\mathcal{L}},\,\mathcal{S}\big] & = A\,\mathcal{T} +
B\,\mathcal{S} + C\,\overline{\mathcal{S}},
\\
\big[\mathcal{L},\,\overline{\mathcal{S}}\big] & =
\overline{A}\,\mathcal{T} + \overline{C}\,\mathcal{S} +
\overline{B}\,\overline{\mathcal{S}},
\\
\big[ \overline{\mathcal{L}},\,\overline{\mathcal{S}}\big] & =
\overline{P}\,\mathcal{T} + \overline{R}\,\mathcal{S} +
\overline{Q}\,\overline{\mathcal{S}}.
\endaligned
\]

\begin{Lemma}
In fact, the above vector fields $[\overline{\mathcal L},\mathcal S]$
and $[\mathcal L,\overline{\mathcal S}]$ are real and equal. In
particular, $A$ is a real-valued function and $C=\overline{B}$.
\end{Lemma}

\proof Indeed, given any two real or complex vector fields $H_1$ and
$H_2$ on any manifold, one always has the following consequence of
the Jacobi identity ({\em cf.}~\cite{Merker-Sabzevari-CEJM}, eq.~\thetag{15}, p.~1817):
\[
\big[H_2,\,\big[H_1,\,[H_1,\,H_2]\big]\big] =
\big[H_1,\,\big[H_2,\,[H_1,\,H_2]\big]\big],
\]
an identity which is true in free Lie algebras. Applied to $H_1 :=
\mathcal{ L}$ and $H_2 := \overline{ \mathcal{ L}}$, this identity
gives that the following two vector fields which are visibly
conjugate to each other:
\[
\aligned \big[\overline{\mathcal{L}},\, \mathcal{S}\big] =
\big[\overline{\mathcal{L}},\,\big[\mathcal{L},\,\mathcal{T}\big]\big]
& = \big[\overline{\mathcal{L}},\,\big[\mathcal{L},\,i[\mathcal{L},\,
\overline{\mathcal{L}}]\big]\big] =
i\big[\overline{\mathcal{L}},\,\big[\mathcal{L},\,[\mathcal{L},\,
\overline{\mathcal{L}}]\big]\big],
\\
\big[\mathcal{L},\, \overline{\mathcal{S}}\big] =
\big[\mathcal{L},\,\big[\overline{\mathcal{L}},\,\mathcal{T}\big]\big]
& = \big[\mathcal{L},\,\big[\overline{\mathcal{L}},\,i[\mathcal{L},\,
\overline{\mathcal{L}}]\big]\big] =
i\big[\mathcal{L},\,\big[\overline{\mathcal{L}},\,[\mathcal{L},\,
\overline{\mathcal{L}}]\big]\big],
\endaligned
\]
are also in fact {\em equal}. This immediately implies that $A$ is a
real function and $C=\overline{B}$.
\endproof

Consequently, the expression of the third of the above four Lie
brackets simplifies as:
\begin{equation}
\label{lenght-two} \boxed{ \aligned
\big[\mathcal{L},\,\mathcal{S}\big] & = P\,\mathcal{T} +
Q\,\mathcal{S} + R\,\overline{\mathcal{S}},
\\
\big[\overline{\mathcal{L}},\,\mathcal{S}\big] & = A\,\mathcal{T} +
B\,\mathcal{S} + \overline{B}\,\overline{\mathcal{S}},
\\
\big[\mathcal{L},\,\overline{\mathcal{S}}\big] & = A\,\mathcal{T} +
B\,\mathcal{S} + \overline{B}\,\overline{\mathcal{S}},
\\
\big[ \overline{\mathcal{L}},\,\overline{\mathcal{S}}\big] & =
\overline{P}\,\mathcal{T} + \overline{R}\,\mathcal{S} +
\overline{Q}\,\overline{\mathcal{S}}.
\endaligned}
\end{equation}

Expressing these five functions $A, B, P, Q$ and $R$ explicitly in
terms of the 4-th order jets $J_{ x, y, u_1, u_2, u_3}^4 \varphi_j$
of the three graphing functions $\varphi_j$ brings formulas that are
of an already quite impressive size. Nonetheless, we aim at
presenting them, in semi-expanded form. At first, we need to have the
expressions of the three {\it coordinate fields}
$\frac{\partial}{\partial u_i}, i=1,2,3$ in terms of the five complex
vector fields $\overline{\mathcal L},\mathcal L, \mathcal T,
\overline{\mathcal S}, \mathcal S$. According to our computations we
have:
\begin{equation}
\label{inversion-d-du-T-S} \aligned \frac{\partial}{\partial u_1} &
=(\Pi^1_{t})\mathcal T+(\Pi^1_{s_1}-i\Pi^1_{s_2})\mathcal
S+(\Pi^1_{s_1}+i\Pi^1_{s_2})\overline{\mathcal S},
\\
\frac{\partial}{\partial u_2} & = (\Pi^2_{t})\mathcal
T+(\Pi^2_{s_1}-i\Pi^2_{s_2})\mathcal
S+(\Pi^2_{s_1}+i\Pi^2_{s_2})\overline{\mathcal S},
\\
\frac{\partial}{\partial u_3} & = (\Pi^3_{t})\mathcal
T+(\Pi^3_{s_1}-i\Pi^3_{s_2})\mathcal
S+(\Pi^3_{s_1}+i\Pi^3_{s_2})\overline{\mathcal S},
\endaligned
\end{equation}
where the coefficients are the real functions:
\[
\footnotesize\aligned \Pi^1_t
&=\Delta^3\big(\Gamma^2_2\Gamma^3_1-\Gamma^3_2\Gamma^2_1\big)/4\Pi, \
\ \ \ \ \ \ \ \ \ \ \ \ \ \ \ \ \ \ \ \ \Pi^1_{s_1}=
\Delta^5\big(\Upsilon_3 \Gamma^2_2-\Gamma^3_2 \Upsilon_2 \big)/16\Pi,
\\
\Pi^1_{s_2}&=\Delta^5\big(\Upsilon_2 \Gamma^3_1-\Upsilon_3 \Gamma^2_1
\big)/16\Pi, \ \ \ \ \ \ \ \ \ \ \ \ \ \ \ \ \ \ \ \ \
\Pi^2_t=\Delta^3\big(\Gamma^3_1 \Gamma^1_2-\Gamma^1_1 \Gamma^3_2
\big)/4\Pi,
\\
\Pi^2_{s_1}&=\Delta^5\big(\Upsilon_3\Gamma^1_2-\Upsilon_1\Gamma^3_2
\big)/16\Pi, \ \ \ \ \ \ \ \ \ \ \ \ \ \ \ \ \ \ \ \
\Pi^2_{s_2}=\Delta^5\big(\Upsilon_1 \Gamma^3_1-\Gamma^1_1
\Upsilon_3\big)/16\Pi,
\\
\Pi^3_t&=\Delta^3\big(\Gamma^2_1 \Gamma^1_2-\Gamma^1_1
\Gamma^2_2\big)/16\Pi, \ \ \ \ \ \ \ \ \ \ \ \ \ \ \ \ \ \ \ \ \
\Pi^3_{s_1}=\Delta^5\big(\Upsilon_2 \Gamma^1_2-\Upsilon_1 \Gamma^2_2
\big)/16\Pi,
\\
\Pi^3_{s_3}&=\Delta^5\big(\Upsilon_1 \Gamma^2_1-\Upsilon_2
\Gamma^1_1\big)/16\Pi,
\endaligned
\]
and where the denominator $\Pi$ explicitly is:
\[
\Pi=-\Upsilon_1 \Gamma^2_2 \Gamma^3_1- \Upsilon_2 \Gamma^1_1
\Gamma^3_2+ \Upsilon_2 \Gamma^1_2 \Gamma^3_1+\Upsilon_1 \Gamma^3_2
\Gamma^2_1+ \Upsilon_3 \Gamma^1_1 \Gamma^2_2-\Upsilon_3 \Gamma^1_2
\Gamma^2_1.
\]
Computing directly the Lie brackets $[\mathcal L,\mathcal S]$ and
$[\overline{\mathcal L},\mathcal S]$ and putting the above
expressions of $\frac{\partial}{\partial u_i}, i=1,2,3$, the
following expressions bring for $A,B,C,P,Q$ and $R$ (as mentioned
before, the expressions of these functions are too much
extensive. That is
why we divide them into several sub-terms):
\begin{equation}
\label{PQRABC} \aligned
P&=\Phi^1_3-\Phi^1_1+2\,i\,\Phi^1_2,
\\
\nonumber
Q&=
\textstyle{\frac{1}{4}}\big(\Phi^2_1-\Phi^2_3-2\,\Phi^3_2\big)+\frac{i}{4}\big(\Phi^3_3-\Phi^3_1-2\,\Phi^2_2\big),
\\
\nonumber R&=
\textstyle{\frac{1}{4}}\big(\Phi^2_1-\Phi^2_3+2\,\Phi^3_2\big)+\frac{i}{4}\big(\Phi^3_1-\Phi^3_3-2\,\Phi^2_2\big),
\\
\nonumber A&=-\Phi^1_1-\Phi^1_3,
\\
\nonumber B&=
\textstyle{\frac{1}{4}}\big(\Phi^2_1+\Phi^2_3\big)+\frac{i}{4}\big(\Phi^3_1+\Phi^3_3\big),
\endaligned
\end{equation}
where the terms $\Phi^\bullet_\bullet$ are:
\[
\footnotesize\aligned \Phi^1_1&=
\frac{-\Gamma^1_1\Gamma^3_2\Omega^2_1-
\Gamma^1_2\Gamma^2_1\Omega^3_1+ \Gamma^3_1\Gamma^1_2\Omega^2_1+
\Gamma^1_1\Omega^3_1\Gamma^2_2+ \Omega^1_1\Gamma^2_1\Gamma^3_2-
\Omega^1_1\Gamma^3_1\Gamma^2_2}{\Delta^2\Sigma},
\\
\Phi^2_1&= -\frac{\Upsilon_3\Gamma^1_2\Omega^2_1-
\Upsilon_3\Omega^1_1\Gamma^2_2- \Omega^2_1\Upsilon_1\Gamma^3_2+
\Upsilon_2\Omega^1_1\Gamma^3_2- \Omega^3_1\Gamma^1_2\Upsilon_2+
\Omega^3_1\Upsilon_1\Gamma^2_2}{\Sigma},
\\
\Phi^3_1&= \frac{\Gamma^2_1\Omega^3_1\Upsilon_1-
\Gamma^2_1\Upsilon_3\Omega^1_1+ \Upsilon_3\Gamma^1_1\Omega^2_1-
\Omega^3_1\Gamma^1_1\Upsilon_2- \Omega^2_1\Gamma^3_1\Upsilon_1+
\Gamma^3_1\Upsilon_2\Omega^1_1}{\Sigma},
\\
\Phi^1_2&= \frac{-\Gamma^1_1\Gamma^3_2\Omega^2_2-
\Gamma^1_2\Gamma^2_1\Omega^3_2+ \Gamma^3_1\Gamma^1_2\Omega^2_2+
\Gamma^1_1\Omega^3_2\Gamma^2_2+ \Omega^1_2\Gamma^2_1\Gamma^3_2-
\Omega^1_2\Gamma^3_1\Gamma^2_2}{\Delta^2\Sigma},
\\
\Phi^2_2&= -\frac{\Upsilon_3\Gamma^1_2\Omega^2_2-
\Upsilon_3\Omega^1_2\Gamma^2_2- \Omega^2_2\Upsilon_1\Gamma^3_2+
\Upsilon_2\Omega^1_2\Gamma^3_2- \Omega^3_2\Gamma^1_2\Upsilon_2+
\Omega^3_2\Upsilon_1\Gamma^2_2}{\Sigma},
\\
\Phi^3_2&= \frac{\Gamma^2_1\Omega^3_2\Upsilon_1-
\Gamma^2_1\Upsilon_3\Omega^1_2+ \Upsilon_3\Gamma^1_1\Omega^2_2-
\Omega^3_2\Gamma^1_1\Upsilon_2- \Omega^2_2\Gamma^3_1\Upsilon_1+
\Gamma^3_1\Upsilon_2\Omega^1_2}{\Sigma},
\\
\Phi^1_3&= \frac{-\Gamma^1_1\Gamma^3_2\Omega^2_3-
\Gamma^1_2\Gamma^2_1\Omega^3_3+ \Gamma^3_1\Gamma^1_2\Omega^2_3+
\Gamma^1_1\Omega^3_3\Gamma^2_2+ \Omega^1_3\Gamma^2_1\Gamma^3_2-
\Omega^1_3\Gamma^3_1\Gamma^2_2}{\Delta^2\Sigma},
\\
\Phi^2_3&= -\frac{\Upsilon_3\Gamma^1_2\Omega^2_3-
\Upsilon_3\Omega^1_3\Gamma^2_2- \Omega^2_3\Upsilon_1\Gamma^3_2+
\Upsilon_2\Omega^1_3\Gamma^3_2- \Omega^3_3\Gamma^1_2\Upsilon_2+
\Omega^3_3\Upsilon_1\Gamma^2_2}{\Sigma},
\\
\Phi^3_3&= \frac{\Gamma^2_1\Omega^3_3\Upsilon_1-
\Gamma^2_1\Upsilon_3\Omega^1_3+ \Upsilon_3\Gamma^1_1\Omega^2_3-
\Omega^3_3\Gamma^1_1\Upsilon_2- \Omega^2_3\Gamma^3_1\Upsilon_1+
\Gamma^3_1\Upsilon_2\Omega^1_3}{\Sigma},
\endaligned
\]
where $\Sigma$ in the denominator is given explicitly by:
\[\footnotesize\aligned
\Sigma=\Delta^2\big(\Gamma^3_1\Gamma^1_2\Upsilon_2
-\Gamma^1_1\Gamma^3_2\Upsilon_2- \Gamma^1_2\Gamma^2_1\Upsilon_3+
\Gamma^1_1\Upsilon_3\Gamma^2_2+ \Upsilon_1\Gamma^2_1\Gamma^3_2-
\Upsilon_1\Gamma^3_1\Gamma^2_2\big),
\endaligned
\]
and where\,\,---\,\,again admitting the notational coincidences
$x\equiv x_1$ and $y\equiv x_2$\,\,---\,\,the functions
$\Omega^\bullet_\bullet$ are:
\begin{equation*}
\footnotesize\aligned
\Omega^1_i&= \Delta^2\Gamma^1_{ix}-
5\Delta\Delta_{x_i}\Gamma^1_1+ \Delta\Lambda^1_1\Gamma^1_{iu_1}-
4\Delta_{u_1}\Lambda^1_1\Gamma^1_{i}-
\Delta\Lambda^1_{1u_1}\Gamma^1_{i}- \Delta\lambda^1_{u_2}\Gamma^2_1+
\\
&+\Delta_{u_2}\Lambda^1_1\Gamma^1_{21}-
\Delta\Lambda^1_{1u_3}\Gamma^3_{i}+
\Delta_{u_3}\Lambda^1_1\Gamma^3_{i}+
\Delta\Lambda^2_1\Gamma^1_{iu_2}-
5\Delta_{u_2}\Lambda^2_1\Gamma^1_{i}+
\Delta\Lambda^3_1\Gamma^1_{iu_3}-
5\Delta_{u_3}\Lambda^3_1\Gamma^1_{i},
\\
\Omega^2_i&= \Delta^2\Gamma^2_{ix}- 5\Delta\Delta_{x_i}\Gamma^2_1+
\Delta\Lambda^1_1\Gamma^2_{iu_1}-
5\Delta_{u_1}\Lambda^1_1\Gamma^2_{i}-
\Delta\Lambda^2_{1u_1}\Gamma^1_{i}+
\Delta_{u_1}\Lambda^2_1\Gamma^1_{i}+
\\
&+ \Delta\Lambda^2_{1}\Gamma^2_{iu_2}-
4\Delta_{u_2}\Lambda^2_1\Gamma^2_{i}-
\Delta\Lambda^2_{1u_2}\Gamma^2_{i}-
\Delta\Lambda^2_{1u_3}\Gamma^3_{i}+
\Delta_{u_3}\Lambda^2_1\Gamma^3_{i}+
\Delta\Lambda^3_1\Gamma^2_{iu_3}-
5\Delta_{u_3}\Lambda^3_1\Gamma^2_{i},
\\
\Omega^3_i&= \Delta^2\Gamma^3_{ix}-
5\Delta\Delta_{x_i}\Gamma^3_1+ \Delta\Lambda^1_1\Gamma^3_{iu_1}-
5\Delta_{u_1}\Lambda^1_1\Gamma^3_{i}+
\Delta\Lambda^2_{1}\Gamma^3_{iu_2}-
5\Delta_{u_2}\Lambda^2_1\Gamma^3_{i}-
\\
&- \Delta\Lambda^3_{1u_1}\Gamma^1_{i}+
\Delta_{u_1}\Lambda^3_1\Gamma^1_{i}-
\Delta\Lambda^3_{1u_2}\Gamma^2_{i}+
\Delta_{u_2}\Lambda^3_{1}\Gamma^2_{i}+
\Delta\Lambda^3_1\Gamma^3_{iu_3}-
4\Delta_{u_3}\Lambda^3_1\Gamma^3_{i}-
\Delta\Lambda^3_{1u_3}\Gamma^3_{i},
\\
\Omega^1_3&= \Delta^2\Gamma^1_{2x}- 5\Delta\Delta_y\Gamma^1_1+
\Delta\Lambda^1_1\Gamma^1_{2u_1}-
4\Delta_{u_1}\Lambda^1_1\Gamma^1_{3}-
\Delta\Lambda^1_{2u_1}\Gamma^1_{3}- \Delta\lambda^1_{u_2}\Gamma^2_1+
\\
&+\Delta_{u_2}\Lambda^1_1\Gamma^1_{21}-
\Delta\Lambda^1_{2u_3}\Gamma^3_{3}+
\Delta_{u_3}\Lambda^1_1\Gamma^3_{3}+
\Delta\Lambda^2_1\Gamma^1_{2u_2}-
5\Delta_{u_2}\Lambda^2_1\Gamma^1_{3}+
\Delta\Lambda^3_1\Gamma^1_{2u_3}-
5\Delta_{u_3}\Lambda^3_1\Gamma^1_{3},
\endaligned
\end{equation*}
\begin{equation*}
\footnotesize\aligned
 \Omega^2_3&= \Delta^2\Gamma^2_{2x}- 5\Delta\Delta_y\Gamma^2_1+
\Delta\Lambda^1_1\Gamma^2_{2u_1}-
5\Delta_{u_1}\Lambda^1_1\Gamma^2_{3}-
\Delta\Lambda^2_{2u_1}\Gamma^1_{3}+
\Delta_{u_1}\Lambda^2_1\Gamma^1_{3}+
\\
&+ \Delta\Lambda^2_{2}\Gamma^2_{2u_2}-
4\Delta_{u_2}\Lambda^2_1\Gamma^2_{3}-
\Delta\Lambda^2_{2u_2}\Gamma^2_{3}-
\Delta\Lambda^2_{2u_3}\Gamma^3_{3}+
\Delta_{u_3}\Lambda^2_1\Gamma^3_{3}+
\Delta\Lambda^3_1\Gamma^2_{2u_3}-
5\Delta_{u_3}\Lambda^3_1\Gamma^2_{3},
\\
 \Omega^3_3&= \Delta^2\Gamma^3_{2x}- 5\Delta\Delta_y\Gamma^3_1+
\Delta\Lambda^1_1\Gamma^3_{2u_1}-
5\Delta_{u_1}\Lambda^1_1\Gamma^3_{3}+
\Delta\Lambda^2_{2}\Gamma^3_{2u_2}-
5\Delta_{u_2}\Lambda^2_1\Gamma^3_{3}-
\\
&- \Delta\Lambda^3_{2u_1}\Gamma^1_{3}+
\Delta_{u_1}\Lambda^3_1\Gamma^1_{3}-
\Delta\Lambda^3_{2u_2}\Gamma^2_{3}+
\Delta_{u_2}\Lambda^3_{2}\Gamma^2_{3}+
\Delta\Lambda^3_1\Gamma^3_{2u_3}-
4\Delta_{u_3}\Lambda^3_1\Gamma^3_{3}-
\Delta\Lambda^3_{2u_3}\Gamma^3_{3}.
\endaligned
\end{equation*}
It yet remains to compute the $3$ among $10$ structure Lie brackets:
\[
\big[\mathcal{T},\,\mathcal{S}\big],
\ \ \ \ \ \ \ \ \ \ \ \ \
\big[\mathcal{T},\,\overline{\mathcal{S}}\big],
\ \ \ \ \ \ \ \ \ \ \ \ \
\big[\mathcal{S},\,\overline{\mathcal{S}}\big].
\]

\subsection{Two structure brackets of length 5}
Next, for the two iterated Lie brackets:
\[
\aligned
\big[\mathcal{T},\,\mathcal{S}\big]
&
=
\Big[
i[\mathcal{L},\overline{\mathcal{L}}],\,
\big[\mathcal{L},\,i[\mathcal{L},\overline{\mathcal{L}}]\big]
\Big],
\\
\big[\mathcal{T},\,\overline{\mathcal{S}}\big]
&
=
\Big[
i[\mathcal{L},\overline{\mathcal{L}}],\,
\big[\overline{\mathcal{L}},\,i[\mathcal{L},\overline{\mathcal{L}}]\big]
\Big],
\endaligned
\]
that are visibly of length $5$. Again the Jacobi identity helps us to
specify their expressions.

\begin{Lemma}
\label{E-F-G}
The coefficients of the two Lie brackets:
\[
\aligned \big[\mathcal{T},\,\mathcal{S}\big] & = E\,\mathcal{T} +
F\,\mathcal{S} + G\,\overline{\mathcal{S}},
\\
\big[\mathcal{T},\,\overline{\mathcal{S}}\big] & =
\overline{E}\,\mathcal{T} + \overline{G}\,\mathcal{S} +
\overline{F}\,\overline{\mathcal{S}},
\endaligned
\]
are three complex-valued functions $E$, $F$, $G$ which can be
expressed as follows in terms of $P$, $Q$, $R$, $A$, $B$ and their
first-order frame derivatives:
\[
\footnotesize
\aligned
E
&
=
-\,i\,\overline{\mathcal L}(P)
-
i\,A\,Q
-
i\,\overline{P}\,R
+
i\,\mathcal L(A)
+
i\,B\,P
+
i\,A\,\overline{B},
\\
F
&
=
-\,i\,\overline{\mathcal L}(Q)
-
i\,R\,\overline{R}
+
i\,A
+
i\,\mathcal{L}(B)
+
i\,B\,\overline{B},
\\
G
&
=
-\,i\,P
-
i\,\overline{B}\,Q
-
i\,R\,\overline{Q}
-
i\,\overline{\mathcal L}(R)
+
i\,B\,R
+
i\,\overline{B}\overline{B}
+
i\,\mathcal L(\overline{B}).
\endaligned
\]
\end{Lemma}

\proof
A glance at the explicit expressions of the three complex fields
$\mathcal T,\mathcal S,\overline{\mathcal S}$ shows that they are
some combinations of the three coordinate fields $\frac{\partial}{\partial
u_1}$, $\frac{\partial}{\partial u_2}$, $\frac{\partial}{\partial
u_3}$, {\em without} any $\frac{\partial}{\partial z}$,
$\frac{\partial}{\partial \overline{z}}$. Therefore, the two
considered Lie brackets must merely be some combinations of these three
coordinate fields $\frac{\partial}{\partial u_1}$,
$\frac{\partial}{\partial u_2}$, $\frac{\partial}{\partial u_3}$,
too. Then by inversion, they must become linear combinations
of the three fields $\mathcal{ T}$, $\mathcal{ S}$,
$\overline{ \mathcal{ S}}$.
In fact, using~\thetag{\ref{inversion-d-du-T-S}}, one can
expand the following Jacobi identity:
\begin{equation*}
\footnotesize\aligned
-i\,
[\mathcal T,\mathcal S]
&
=
\big[[\mathcal
L,\overline{\mathcal L}],\mathcal S\big]
\\
&
=
\big[\,
\overline{\mathcal L},[\mathcal S,\mathcal L]\big]+\big[\mathcal
L,[\overline{\mathcal L},\mathcal S]\big]
\\
&
=
[\overline{\mathcal
L},-P\,\mathcal T-Q\,\mathcal S-R\,\overline{\mathcal S}]+[{\mathcal
L},A\,\mathcal T+B\,\mathcal S+\overline{B}\,\overline{\mathcal
S}]
\\
&
=
-\overline{\mathcal
L}(P)\,\mathcal T-P\,\overline{\mathcal S}-\overline{\mathcal
L}(Q)\,\mathcal S-Q\,(A\,\mathcal T+\zero{B\,\mathcal
S}+\overline{B}\,\overline{\mathcal S})-\overline{\mathcal
L}(R)\,\overline{\mathcal S}-
\\
& \ \ \ \ \
-
R(\overline{P}\,\mathcal
T+\overline{R}\,\mathcal S+\overline{Q}\,\overline{\mathcal
S})+A\,\mathcal S+\mathcal L(A)\,\mathcal T+
\\
& \ \ \ \ \
+
B(P\,\mathcal T+\zero{Q\,\mathcal
S}+R\,\overline{\mathcal S})+ \mathcal L(B)\,\mathcal
S+\overline{B}\,(A\mathcal T+B\,\mathcal
S+\overline{B}\,\overline{\mathcal S})+\mathcal
L(\overline{B})\,\overline{\mathcal S}.
\endaligned
\end{equation*}
Here the underlined terms vanish by pair. Now, extracting
the coefficients of $\mathcal T$, $\mathcal S$, $\overline{\mathcal
S}$ gives the desired expressions of $E,F$ and $G$,
respectively.
\endproof

\subsection{The last structure bracket of length 6}
Now, the last remaining iterated bracket:
\[
\aligned
\big[\mathcal{S},\,\overline{\mathcal{S}}\big]
&
=
\Big[
\big[\mathcal{L},\,i[\mathcal{L},\overline{\mathcal{L}}]\big],\,\,
\big[\overline{\mathcal{L}},\,i[\mathcal{L},\overline{\mathcal{L}}]\big]
\Big]
\endaligned
\]
is of length $6$, and because Jacobi identities at this level start to
become more complex, we must take care of how to re-express it.  In a
preceding publication ({\em see} equations
$\overset{1}{\,=\,}$ and $\overset{2}{\,=\,}$ page 1818
of~\cite{Merker-Sabzevari-CEJM}), we showed that, in a free Lie
algebra generated by two vectors $h_1$ and $h_2$, the following two
relations hold true:
\[
\!\!\!\!\!\!\!\!\!\!\!\!\!\!\!\!\!\!\!\!
\!\!\!\!\!\!\!\!\!\!\!\!
\aligned
0
&
=
\big[\big[h_1,[h_1,h_2]\big],
\big[h_2,[h_1,h_2]\big]\big]
-
\big[[h_1,h_2],
\big[h_1,\big[h_2,[h_1,h_2]\big]\big]\big]
-
\big[h_2,\big[h_2,\big[h_1,\big[h_1,
[h_1,h_2]\big]\big]\big]\big]
+
\\
&
\ \ \ \ \ \ \ \ \ \ \ \ \ \ \ \ \ \ \ \ \ \ \ \ \ \ \ \ \ \ \ \ \ \ \
\ \ \ \ \ \ \ \ \ \ \ \ \ \ \ \ \ \ \ \ \ \ \ \ \ \ \ \ \ \ \ \ \ \ \
\ \ \ \ \ \ \ \ \ \ \ \ \ \ \ \ \ \ \ \ \ \ \ \ \ \ \ \ \ \ \ \ \
+
\big[h_2,\big[h_1,\big[h_1,\big[h_2,
[h_1,h_2]\big]\big]\big]\big],
\\
0
&
=
\big[\big[h_1,[h_1,h_2]\big],
\big[h_2,[h_1,h_2]\big]\big]
+
\big[[h_1,h_2],
\big[h_1,\big[h_2,[h_1,h_2]\big]\big]\big]
+
\big[h_1,\big[h_2,\big[h_1,\big[h_2,
[h_1,h_2]\big]\big]\big]\big]
-
\\
&
\ \ \ \ \ \ \ \ \ \ \ \ \ \ \ \ \ \ \ \ \ \ \ \ \ \ \ \ \ \ \ \ \ \ \
\ \ \ \ \ \ \ \ \ \ \ \ \ \ \ \ \ \ \ \ \ \ \ \ \ \ \ \ \ \ \ \ \ \ \
\ \ \ \ \ \ \ \ \ \ \ \ \ \ \ \ \ \ \ \ \ \ \ \ \ \ \ \ \ \ \ \ \
-
\big[h_1,\big[h_1,\big[h_2,\big[h_2,
[h_1,h_2]\big]\big]\big]\big].
\endaligned
\]
Then by a plain addition, replacing $h_1 := \mathcal{ L}$ and
$h_2 := \overline{ \mathcal{ L}}$, we may express our length-six
Lie bracket in terms of simple-words Lie brackets:
\[
\!\!\!\!\!\!\!\!\!\!\!\!\!\!\!\!\!\!\!\!
\footnotesize
\aligned
-\,2\,
\Big[
\big[\mathcal{L},\,[\mathcal{L},\overline{\mathcal{L}}]\big],\,
\big[\overline{\mathcal{L}},\,[\mathcal{L},\overline{\mathcal{L}}]\big]
\Big]
&
=
-\,
\Big[\overline{\mathcal{L}},\,\Big[\overline{\mathcal{L}},\,
\big[\mathcal{L},\,\big[\mathcal{L},\,
[\mathcal{L},\overline{\mathcal{L}}]
\big]\big]\Big]\Big]
+
\Big[\overline{\mathcal{L}},\,\Big[\mathcal{L},\,
\big[\mathcal{L},\,\big[\overline{\mathcal{L}},
[\mathcal{L},\overline{\mathcal{L}}]
\big]\big]\Big]\Big]
+
\\
&
\ \ \ \ \
+
\Big[\mathcal{L},\,\Big[\overline{\mathcal{L}},\,
\big[\mathcal{L},\,\big[\overline{\mathcal{L}},
[\mathcal{L},\overline{\mathcal{L}}]
\big]\big]\Big]\Big]
-
\Big[\mathcal{L},\,\Big[\mathcal{L},\,
\big[\overline{\mathcal{L}},\,\big[\overline{\mathcal{L}},
[\mathcal{L},\overline{\mathcal{L}}]
\big]\big]\Big]\Big].
\endaligned
\]
Notably, this expression points out the purely imaginary character
{\em of both sides}, hence it is appropriate for
what will follow. By expanding the four
simple-words Lie brackets, we obtain:

\begin{Lemma}
\label{Theta-Lambda}
The coefficients of the last, tenth structure bracket:
\[
\aligned
\big[\mathcal{S},\,\overline{\mathcal{S}}\big]
&
=
i\,J\,\mathcal{T} + K\,\mathcal{S} -
\overline{K}\,\overline{\mathcal{S}},
\endaligned
\]
are one complex-valued function $K$ and one real-valued function
$J$ which can be expressed as follows in terms of $P$, $Q$, $R$,
$A$, $B$ and their frame derivatives up to order $2$:
\[
\!\!\!\!\!\!\!\!\!\!\!\!\!\!\!\!\!\!\!\!
\footnotesize
\aligned
-\,2\,J
&
=
-\,\overline{\mathcal{L}}\big(\overline{\mathcal{L}}(P)\big)
+
\overline{\mathcal{L}}\big(\mathcal{L}(A)\big)
+
\mathcal{L}\big(\overline{\mathcal{L}}(A)\big)
-
\mathcal{L}\big(\mathcal{L}(\overline{P})\big)
\\
&
\ \ \ \ \
-\,Q\,\overline{\mathcal{L}}(A)
-
2\,A\,\overline{\mathcal{L}}(Q)
-
R\,\overline{\mathcal{L}}(\overline{P})
-
2\,\overline{P}\,\overline{\mathcal{L}}(R)
-
2\,AR\overline{R}
-
2\,P\overline{P}
-
\overline{B}\overline{P}Q
-
\overline{P}\overline{Q}R
-
\\
&
\ \ \ \ \
-
\overline{R}\mathcal{L}(P)
-
2\,P\,\mathcal{L}(\overline{R})
-
\overline{Q}\,\mathcal{L}(A)
-
2\,A\mathcal{L}(\overline{Q})
-
PQ\overline{R}
-
BP\overline{Q}
+
\\
&
\ \ \ \ \
+
2\,P\overline{\mathcal{L}}(B)
+
B\overline{\mathcal{L}}(P)
+
2\,A\overline{\mathcal{L}}(\overline{B})
+
\overline{B}\,\overline{\mathcal{L}}(A)
+
2\,A\,\mathcal{L}(B)
+
2\,AA
+
2\,AB\overline{B}
+
2\,\overline{P}\,\mathcal{L}(\overline{B})
+
\\
&
\ \ \ \ \
+
B\overline{P}R
+
\overline{B}\overline{B}\overline{P}
+
B\,\mathcal{L}(A)
+
\overline{B}\,\mathcal{L}(\overline{P})
+
BBP
+
\overline{B}P\overline{R},
\endaligned
\]
\[
\!\!\!\!\!\!\!\!\!\!\!\!\!\!\!\!\!\!\!\!
\footnotesize
\aligned
2i\,K
&
=
-\,\overline{\mathcal{L}}\big(\overline{\mathcal{L}}(Q)\big)
+
\overline{\mathcal{L}}\big(\mathcal{L}(B)\big)
+
\mathcal{L}\big(\overline{\mathcal{L}}(B)\big)
-
\mathcal{L}\big(\mathcal{L}(\overline{R})\big)
-
\\
&
\ \ \ \ \
-
2\,\overline{R}\,\overline{\mathcal{L}}(R)
-
R\,\overline{\mathcal{L}}(\overline{R})
-
B\,\overline{\mathcal{L}}(Q)
-
BR\overline{R}
-
2\,P\overline{R}
-
\overline{Q}R\overline{R}
-
2\,\mathcal{L}(\overline{P})
-
\overline{R}\,\mathcal{L}(Q)
-
\\
&
\ \ \ \ \
-\,2\,Q\,\mathcal{L}(\overline{R})
-
\overline{Q}\,\mathcal{L}(B)
-
2\,B\,\mathcal{L}(\overline{Q})
-
A\overline{Q}
-
\overline{P}Q
-
QQ\overline{R}
-
BQ\overline{Q}
+
\\
&
\ \ \ \ \
+
2\,\overline{\mathcal{L}}(A)
+
\overline{B}\,\overline{\mathcal{L}}(B)
+
2\,B\,\overline{\mathcal{L}}(\overline{B})
+
3\,B\,\mathcal{L}(B)
+
3\,AB
+
BBQ
+
2\,BB\overline{B}
+
2\,\overline{R}\,\mathcal{L}(\overline{B})
+
\\
&
\ \ \ \ \
+
\overline{B}\overline{B}\overline{R}
+
\overline{B}\,\mathcal{L}(\overline{R})
+
\overline{B}\overline{P}
+
Q\,\overline{\mathcal{L}}(B).
\qed
\endaligned
\]
\end{Lemma}

\subsection{Further relations}
\label{length-six}

Moreover, we also showed in \cite{Merker-Sabzevari-CEJM} ({\it see} equations $\overset{9}{\,=\,}, \overset{10}{\,=\,},\overset{11}{\,=\,}$ pp.~1818--1819 of this paper) that exactly {\em three} independent linear relations hold
between simple-word Lie brackets, which, applied to $h_1 := \mathcal{
L}$ and $h_2 := \overline{ \mathcal{ L}}$, are:
\[
\!\!\!\!\!\!\!\!\!\!\!\!\!\!\!\!\!\!\!\!
\!\!\!\!\!\!\!\!\!\!\!
\footnotesize
\aligned
0
&
\overset{1}{=}
\Big[\mathcal{L},\,\Big[\mathcal{L},\,
\big[\mathcal{L},\,\big[\overline{\mathcal{L}},
[\mathcal{L},\overline{\mathcal{L}}]
\big]\big]\Big]\Big]
-
2\,
\Big[\mathcal{L},\,\Big[\overline{\mathcal{L}},\,
\big[\mathcal{L},\,\big[\mathcal{L},
[\mathcal{L},\overline{\mathcal{L}}]
\big]\big]\Big]\Big]
+
\Big[\overline{\mathcal{L}},\,\Big[\mathcal{L},\,
\big[\mathcal{L},\,\big[\mathcal{L},
[\mathcal{L},\overline{\mathcal{L}}]
\big]\big]\Big]\Big],
\\
0 & \overset{2}{=}
\Big[\overline{\mathcal{L}},\,\Big[\overline{\mathcal{L}},\,
\big[{\mathcal{L}},\,\big[\overline{\mathcal{L}},
[\mathcal{L},\overline{\mathcal{L}}] \big]\big]\Big]\Big] - 2\,
\Big[\overline{\mathcal{L}},\,\Big[\mathcal{L},\,
\big[\overline{\mathcal{L}},\,\big[\overline{\mathcal{L}},
[\mathcal{L},\overline{\mathcal{L}}] \big]\big]\Big]\Big] +
\Big[\mathcal{L},\,\Big[\overline{\mathcal{L}},\,
\big[\overline{\mathcal{L}},\,\big[\overline{\mathcal{L}},
[\mathcal{L},\overline{\mathcal{L}}] \big]\big]\Big]\Big],
\\
0 & \overset{3}{=} \Big[\mathcal{L},\,\Big[\mathcal{L},\,
\big[\overline{\mathcal{L}},\,\big[\overline{\mathcal{L}},
[\mathcal{L},\overline{\mathcal{L}}] \big]\big]\Big]\Big] - 3\,
\Big[\mathcal{L},\,\Big[\overline{\mathcal{L}},\,
\big[\mathcal{L},\,\big[\overline{\mathcal{L}},
[\mathcal{L},\overline{\mathcal{L}}] \big]\big]\Big]\Big] + 3\,
\Big[\overline{\mathcal{L}},\,\Big[\mathcal{L},\,
\big[{\mathcal{L}},\,\big[\overline{\mathcal{L}},
[\mathcal{L},\overline{\mathcal{L}}] \big]\big]\Big]\Big] -
\\
&
\ \ \ \ \ \ \ \ \ \ \ \ \ \ \ \ \ \ \ \ \ \ \ \ \ \ \ \ \ \ \ \ \ \ \ \
\ \ \ \ \ \ \ \ \ \ \ \ \ \ \ \ \ \ \ \ \ \ \ \ \ \ \ \ \ \ \ \ \ \ \ \
\ \ \ \ \ \ \ \ \ \ \ \ \ \ \ \ \ \ \ \ \ \ \ \ \ \ \ \ \ \ \ \ \ \ \ \
\ \ \ \ \ \ \ \ \ \ \ \ \
-
\Big[\overline{\mathcal{L}},\,\Big[\overline{\mathcal{L}},\,
\big[\mathcal{L},\,\big[\mathcal{L},
[\mathcal{L},\overline{\mathcal{L}}]
\big]\big]\Big]\Big].
\endaligned
\]

Applying the computations relevant to the above three equations, one
sees that the two equations 1 and 2 bring two same outcomes, since each of
them is the conjugation of each other with a negative sign. Also for
the equation 3, the coefficient of $\overline{\mathcal S}$ is the
conjugation of that of $\mathcal S$, again with a negative sign.
Overall, these equations give the following five independent
relationships between the fundamental functions $A,B,P,Q,R$, their
conjugations and derivations:
\begin{equation}\label{length 6}
\footnotesize
\aligned
0&\overset{1}{\equiv} 2\,\mathcal
L(\overline{\mathcal L}(P))-\mathcal L(\mathcal
L(A))-\overline{\mathcal L}(\mathcal L(P))-
\\
&2\,P\mathcal L(B)-B\mathcal L(P)-2\,A\mathcal L(\overline
B)-\overline B\mathcal L(A)+P\overline{\mathcal L}(Q)+A\mathcal L(Q)+
\\
&+2\,Q\mathcal L(A)-Q\overline{\mathcal L}(P)+A\overline{\mathcal
L}(R)+2\,R\mathcal L(\overline P)+\overline P\mathcal
L(R)-R\overline{\mathcal L}(A)-
\\
&-PB\overline B-A\overline B^2+PBQ+2\,AQ\overline
B-AQ^2-2\,ABR+2\,RP\overline R+2\,AR\overline Q-QR\overline
P-R\overline B\,\overline P, \ \ \ \ \ \ \ \ \ \ \ \ \ \ \ \ \ \ \ \
\ \ \ \ \ \ \ \ \ \ \ \ \ \ \ \ \ \ \ \ \ \ \ \ \ \ \ \ \ \ \ \ \ \ \
\ \ \ \ \ \ \ \ \ \ \ \ \ \ \ \ \ \ \ \ \ \ \
\\
0&\overset{2}{\equiv}2\,\mathcal L(\overline{\mathcal
L}(Q))-\mathcal L(\mathcal L(B))-\overline{\mathcal L}(\mathcal
L(Q))-
\\
&-2\,\mathcal L(A)-2\,B\mathcal L(\overline B)-\overline B\mathcal
L(B)+B\overline{\mathcal L}(R)+2\,R\mathcal L(\overline R)+\overline
R\mathcal L(R)-R\overline{\mathcal L}(B)+\overline{\mathcal L}(P)+
\\
&+2\,R\overline P+BQ\overline B-A\overline B-B\overline
B^2+AQ+QR\overline R+2\,BR\overline Q-2\,B^2R-R\overline B\,\overline
R,
\\
0&\overset{3}{\equiv}2\,\mathcal L(\overline{\mathcal L}(R))-\mathcal
L(\mathcal L(\overline B))-
\overline{\mathcal{L}}\big(\mathcal{L}(R)\big)-
\\
&-3\,\overline B\mathcal L(\overline B)+\overline B\mathcal
L(Q)+2\,Q\mathcal L(\overline B)-2\,R\mathcal L(B)-B\mathcal
L(R)+R\overline{\mathcal L}(Q)+\overline B\,\overline{\mathcal L}(R)+
\\
&+2\,R\mathcal L(\overline Q)+\overline Q\mathcal
L(R)-Q\overline{\mathcal L}(R)-\overline{\mathcal L}\mathcal
L(R)-R\overline{\mathcal L}(\overline B)+\mathcal L(P)+
\\
&+2\,Q\overline B^2-QP-Q^2\overline B-\overline B^3+P\overline
B-2\,AR-2\,BR\overline B+BQR+2\,R^2\overline R+R\overline
B\,\overline Q-QR\overline Q,
\\
0&\overset{4}{\equiv}-3\overline{\mathcal L}(\mathcal L(A))-\mathcal L(\mathcal
L(\overline P))+3\,\mathcal
L(\overline{\mathcal L}(A))+\overline{\mathcal L}(\overline{\mathcal
L}(P))-
\\
&-2\,A\mathcal L(\overline Q)-\overline Q\mathcal L(A)+3\,B\mathcal
L(A)+3\,\overline B\mathcal L(\overline P)-3\,B\overline{\mathcal
L}(P)-3\,\overline B\,\overline{\mathcal L}(A)+2\,A\overline{\mathcal
L}(Q)-
\endaligned
\end{equation}
\[\footnotesize
\aligned
&+Q\overline{\mathcal L}(A)-2\,P\mathcal L(\overline R)-\overline
R\mathcal L(P)+2\,\overline P\,\overline{\mathcal
L}(R)+R\overline{\mathcal L}(\overline P)-
\\
 &-BP\overline Q+3\,B^2P+2\,A\overline B\,\overline
Q-2\,BQA-3\,\overline B^2\overline P+Q\overline B\,\overline
P-PQ\overline R+3\,P\overline B\,\overline R-3\,BR\overline
P+R\overline P\,\overline Q,
\\
0&\overset{5}{\equiv}-3\,\overline{\mathcal L}(\mathcal
L(B))+3\,\mathcal L(\overline{\mathcal L}(B))+\overline{\mathcal
L}(\overline{\mathcal L}(Q))-\mathcal L(\mathcal L(\overline R))+
\\
&+3\,B\mathcal L(B)-3\,\overline B\,\overline{\mathcal
L}(B)+Q\overline{\mathcal L}(B)-B\overline{\mathcal
L}(Q)-2\,Q\mathcal L(\overline R)-\overline R\mathcal L(Q)-
\\
&-2\,B\mathcal L(\overline Q)-\overline Q\mathcal L(B)+3\,\overline
B\mathcal L(\overline R)+2\,\overline R\,\overline{\mathcal
L}(R)+R\overline{\mathcal L}(\overline R)-2\,\mathcal L(\overline P)-
\\
&-Q\overline P-A\overline Q-BQ\overline Q+3\,AB+3\,\overline
B\,\overline P+2\,B\overline B\,\overline Q+B^2Q-Q^2\overline
R+4\,Q\overline B\,\overline R-
\\
&-3\,BR\overline R-3\,\overline B^2\overline R+R\overline
Q\,\overline R.
\endaligned
\]

\section{Passage to a dual coframe
\\
and its Darboux-Cartan structure}
\label{passage-coframe}
\HEAD{\ref{passage-coframe}.~Passage to a dual coframe
and its Darboux-Cartan structure}{
Masoud {\sc Sabzevari} (Shahrekord) and Jo\"el {\sc Merker} (LM-Orsay)}

On the natural agreement that:
\[
\text{\rm the coframe}\ \ \
\big\{du_3,du_2,du_1,\,dz,\,d\overline{z}\big\} \ \ \ \text{\rm is
dual to the frame} \ \ \ \big\{ {\textstyle{\frac{\partial}{\partial
u_3}}},\, {\textstyle{\frac{\partial}{\partial u_2}}},\,
{\textstyle{\frac{\partial}{\partial u_1}}},\,
{\textstyle{\frac{\partial}{\partial z}}},\,
{\textstyle{\frac{\partial}{\partial\overline{z}}}} \big\},
\]
let us introduce the coframe:
\[
\big\{\overline{\sigma_0},\,\sigma_0,\,\rho_0,\,\overline{\zeta_0},\,\zeta_0\big\}
\ \ \ \text{\rm which is dual to the frame} \ \ \
\big\{\overline{\mathcal{S}},\,
\mathcal{S},\,\mathcal{T},\,\overline{\mathcal{L}},\,\mathcal{L}\big\},
\]
that is to say which satisfies by definition:
\[
\begin{array}{ccccc}
\overline{\sigma_0}(\overline{\mathcal{S}})=1 \ \ \ & \ \ \
\overline{\sigma_0}(\mathcal{S})=0 \ \ \ & \ \ \
\overline{\sigma_0}(\mathcal{T})=0 \ \ \ & \ \ \
\overline{\sigma_0}(\overline{\mathcal{L}})=0 \ \ \ & \ \ \
\overline{\sigma_0}\big(\mathcal{L}\big)=0,
\\
\sigma_0(\overline{\mathcal{S}})=0 \ \ \ & \ \ \
\sigma_0(\mathcal{S})=1 \ \ \ & \ \ \ \sigma_0(\mathcal{T})=0 \ \ \ &
\ \ \ \sigma_0(\overline{\mathcal{L}})=0 \ \ \ & \ \ \
\sigma_0\big(\mathcal{L}\big)=0,
\\
\rho_0(\overline{\mathcal{S}})=0 \ \ \ & \ \ \ \rho_0(\mathcal{S})=0
\ \ \ & \ \ \ \rho_0(\mathcal{T})=1 \ \ \ & \ \ \
\rho_0(\overline{\mathcal{L}})=0 \ \ \ & \ \ \ \rho_0(\mathcal{L})=0,
\\
\overline{\zeta_0}(\overline{\mathcal{S}})=0 \ \ \ & \ \ \
\overline{\zeta_0}(\mathcal{S})=0 \ \ \ & \ \ \
\overline{\zeta_0}(\mathcal{T})=0 \ \ \ & \ \ \
\overline{\zeta_0}(\overline{\mathcal{L}})=1 \ \ \ & \ \ \
\overline{\zeta_0}\big(\mathcal{L}\big)=0,
\\
\zeta_0(\overline{\mathcal{S}})=0 \ \ \ & \ \ \
\zeta_0(\mathcal{S})=0 \ \ \ & \ \ \ \zeta_0(\mathcal{T})=0 \ \ \ & \
\ \ \zeta_0(\overline{\mathcal{L}})=0 \ \ \ & \ \ \
\zeta_0\big(\mathcal{L}\big)=1.
\end{array}
\]
Since neither $\mathcal{ T}$, nor $\mathcal{ S}$, nor $\overline{
\mathcal{ S}}$ incorporate any $\frac{ \partial}{ \partial u_j}$, $j
= 1, 2, 3$, we have:
\[
\zeta_0 = dz \ \ \ \ \ \ \ \ \ \ \ \ \ \text{\rm and} \ \ \ \ \ \ \ \
\ \ \ \ \ \overline{\zeta_0} = d\overline{z}.
\]
In order to launch the Cartan algorithm, initially we need the
expressions of the five 2-forms
$d\overline{\sigma}_0,d\sigma_0,d\rho_0,d\overline{\zeta}_0,
d\zeta_0$ in terms of the wedge products of
$\overline{\sigma}_0,\sigma_0,\rho_0,\overline{\zeta}_0, \zeta_0$.

To
find them, we remember that if
a frame $\big\{ \mathcal{ L}_1, \dots, \mathcal{
L}_n\big\}$ on an open subset of $\R^n$ enjoys the Lie structure:
\[
\big[\mathcal{L}_{i_1},\,\mathcal{L}_{i_2}\big] =
\sum_{k=1}^n\,a_{i_1,i_2}^k\,\mathcal{L}_k \ \ \ \ \ \ \ \ \ \ \ \ \
{\scriptstyle{(1\,\leqslant\,i_1\,<\,i_2\,\leqslant\,n)}},
\]
then its dual coframe
$\{ \omega^1, \dots, \omega^n \}$ enjoys the
Darboux-Cartan structure:
\[
d\omega^k = - \sum_{1\leqslant i_1<i_2\leqslant n}\,
a_{i_1,i_2}^k\,\omega^{i_1}\wedge\omega^{i_2} \ \ \ \ \ \ \ \ \ \ \ \
\ {\scriptstyle{(k\,=\,1\,\cdots\,n)}}.
\]

Granted this reminder,
it is convenient to rewrite the ten Lie brackets under the form of a
convenient auxiliary array:
\[
\footnotesize
\begin{array}{cccccccccccc}
& & \overline{\mathcal{S}} & & \mathcal{S} & & \mathcal{T} & &
\overline{\mathcal{L}} & & \mathcal{L}
\\
& & \boxed{d\overline{\sigma_0}} & & \boxed{d\sigma_0} & &
\boxed{d\rho_0} & & \boxed{d\overline{\zeta_0}} & & \boxed{d\zeta_0}
\\
\big[\overline{\mathcal{S}},\,\mathcal{S}\big] & = &
\overline{K}\cdot\overline{\mathcal{S}} & + &
-\,K\cdot\mathcal{S} & + & -i\,J\cdot\mathcal{T} & + & 0 &
+ & 0 & \boxed{\overline{\sigma_0}\wedge\sigma_0}
\\
\big[\overline{\mathcal{S}},\,\mathcal{T}\big] & = &
-\,\overline{F}\cdot\overline{\mathcal{S}} & + &
-\,\overline{G}\cdot\mathcal{S} & + &
-\,\overline{E}\cdot\mathcal{T} & + & 0 & + & 0 &
\boxed{\overline{\sigma_0}\wedge\rho_0}
\\
\big[\overline{\mathcal{S}},\,\overline{\mathcal{L}}\big] & = &
-\,\overline{Q}\cdot\overline{\mathcal{S}} & + &
-\,\overline{R}\cdot\mathcal{S} & + & -\,\overline{P}\cdot\mathcal{T}
& + & 0 & + & 0 & \boxed{\overline{\sigma_0}\wedge\overline{\zeta_0}}
\\
\big[\overline{\mathcal{S}},\,\mathcal{L}\big] & = &
-\,\overline{B}\cdot\overline{\mathcal{S}} & + &
-\,{B}\cdot\mathcal{S} & + & -A\cdot\mathcal{T} & + & 0 & + & 0 &
\boxed{\overline{\sigma_0}\wedge\zeta_0}
\\
\big[\mathcal{S},\,\mathcal{T}\big] & = &
-\,G\cdot\overline{\mathcal{S}} & + & -\,F\cdot \mathcal{S} & +
& -\,E\cdot\mathcal{T} & + & 0 & + & 0 &
\boxed{\sigma_0\wedge\rho_0}
\\
\big[\mathcal{S},\,\overline{\mathcal{L}}\big] & = &
-\,\overline{B}\cdot\overline{\mathcal{S}} & + & -\,B\cdot
\mathcal{S} & + & -\,A\cdot\mathcal{T} & + & 0 & + & 0 &
\boxed{\sigma_0\wedge\overline{\zeta_0}}
\\
\big[\mathcal{S},\,\mathcal{L}\big] & = &
-\,R\cdot\overline{\mathcal{S}} & + & -\,Q\cdot\mathcal{S} & + &
-\,P\cdot\mathcal{T} & + & 0 & + & 0 & \boxed{\sigma_0\wedge\zeta_0}
\\
\big[\mathcal{T},\,\overline{\mathcal{L}}\big] & = &
-\,\overline{\mathcal{S}} & + & 0 & + & 0 & + & 0 & + & 0 &
\boxed{\rho_0\wedge\overline{\zeta_0}}
\\
\big[\mathcal{T},\,\mathcal{L}\big] & = & 0 & + & -\,\mathcal{S} & +
& 0 & + & 0 & + & 0 & \boxed{\rho_0\wedge\zeta_0}
\\
\big[\overline{\mathcal{L}},\,\mathcal{L}\big] & = & 0 & + & 0 & + &
i\,\mathcal{T} & + & 0 & + & 0 &
\boxed{\overline{\zeta_0}\wedge\zeta_0}
\end{array}
\]
Thank to this array, we can {\em vertically read} the expressions of
the $10$ forms of degree $2$ that provides the associated
Darboux-Cartan structure, putting an overall minus sign:
\begin{equation}
\label{d0-general} \boxed{ \aligned d\overline{\sigma}_0 & =
-\,\overline{K} \cdot \overline{\sigma}_0\wedge\sigma_0 +
\overline{F}\cdot \overline{\sigma}_0\wedge\rho_0 +
\overline{Q}\cdot \overline{\sigma}_0\wedge\overline{\zeta}_0 +
\overline{B}\cdot \overline{\sigma}_0\wedge\zeta_0 +
\\
& \ \ \ \ \ + G\cdot \sigma_0\wedge\rho_0 + \overline{B}\cdot
\sigma_0\wedge\overline{\zeta}_0 + R\cdot \sigma_0\wedge\zeta_0 +
\rho_0\wedge\overline{\zeta}_0,
\\
d\sigma_0 & = K \cdot \overline{\sigma}_0\wedge\sigma_0 +
\overline{G}\cdot \overline{\sigma}_0\wedge\rho_0 +
\overline{R}\cdot \overline{\sigma}_0\wedge\overline{\zeta}_0 +
{B}\cdot \overline{\sigma}_0\wedge\zeta_0 +
\\
& \ \ \ \ \ + F\cdot \sigma_0\wedge\rho_0 + B\cdot
\sigma_0\wedge\overline{\zeta}_0 + Q\cdot \sigma_0\wedge\zeta_0 +
\rho_0\wedge\zeta_0,
\\
d\rho_0 & = i\,J\cdot \overline{\sigma}_0\wedge\sigma_0 +
\overline{E}\cdot \overline{\sigma}_0\wedge\rho_0 +
\overline{P}\cdot \overline{\sigma}_0\wedge\overline{\zeta}_0 +
A\cdot \overline{\sigma}_0\wedge\zeta_0 +
\\
& \ \ \ \ \ + E\cdot \sigma_0\wedge\rho_0 + A\cdot
\sigma_0\wedge\overline{\zeta}_0 + P\cdot \sigma_0\wedge\zeta_0 -
i\,\overline{\zeta}_0\wedge\zeta_0,
\\
d\overline{\zeta}_0 & = 0,
\\
d\zeta_0 & = 0.
\endaligned}
\end{equation}

\subsection{Ambiguity matrix}
Consider now a local biholomorphic equivalence:
\[
h\colon\ \ \ (z,w) \longmapsto \big(f(z,w),\,g(z,w)\big) =: (z',w')
\]
between any two real analytic local CR-generic maximally minimal real
submanifolds:
\[
M^5 \subset \C_{(z,w)}^4
\ \ \ \ \ \ \ \ \ \
\text{\rm and}
\ \ \ \ \ \ \ \ \ \ {M'}^5\subset\C_{(z',w')}^4.
\]
As we saw in what precedes, the assumption that both $M^5$ and $M'^5$
are maximally minimal means that the two sets of five vector fields:
\[
\Big\{ \mathcal{L},\,\,\overline{\mathcal{L}},\,\, \mathcal T, \,\,
\mathcal S, \,\, \overline{\mathcal S} \Big\} \ \ \ \ \ \ \ \ \ \
\text{\rm and} \ \ \ \ \ \ \ \ \ \ \Big\{
\mathcal{L}',\,\,\overline{\mathcal{L}}',\,\, \mathcal T', \,\,
\mathcal S', \,\, \overline{\mathcal S}' \Big\}
\]
make up frames for $TM^5 \otimes_\R \C$ and for
$T {M'}^5 \otimes_R \C$, respectively.

By a quick inspection of the proof of Proposition~\ref{10-ambiguity-matrix},
one easily convinces oneself that for it to hold
true, we in fact did not at all use the assumption that
$M^5$ was the model cubic $M_{\sf c}^5$. So in
the general case, we also have:

\begin{Proposition}
The initial ambiguity matrix associated to the equivalence problem
under local biholomorphic transformations for maximally minimal
CR-generic $3$-codimensional submanifolds $M^5 \subset \C^4$ is of the
general form:
\[
\left(\!\!
\begin{array}{ccccc}
{\sf a}\overline{\sf a}\overline{\sf a} & 0 & \overline{\sf c} &
\overline{\sf e} & \overline{\sf d}
\\
0 & {\sf a}{\sf a}\overline{\sf a} & {\sf c} & {\sf d} & {\sf e}
\\
0 & 0 & {\sf a}\overline{\sf a} & \overline{\sf b} & {\sf b}
\\
0 & 0 & 0 & \overline{\sf a} & 0
\\
0 & 0 & 0 & 0 & {\sf a}
\end{array}
\!\!\right),
\]
where ${\sf a}$, ${\sf b}$, ${\sf c}$, ${\sf e}$, ${\sf d}$ are
complex numbers. Moreover, the collection of all these matrices makes
up a real 10-dimensional matrix Lie subgroup of ${\sf GL}_5 ( \C)$.
\qed
\end{Proposition}

\subsection{Setting up the equivalence problem}

Quite similarly as in the case where $M^5$ was the cubic
model $M_{\sf c}^5$\,\,---\,\,though the levels of complexity will
rapidly diverge\,\,---, with the
initial coframe $\big\{ \overline{ \sigma}_0,
\sigma_0, \rho_0, \overline{ \zeta}_0, \zeta_0 \big\}$
dual to the explicitly computed frame
$\big\{ \overline{ \mathcal{ S}}, \mathcal{ S},
\mathcal{ T}, \overline{ \mathcal{ L}}, \mathcal{ L}
\big\}$,
the lifted coframe is then (one must transpose the
ambiguity matrix):
\begin{equation}
\label{lift-coframe-general}
\aligned \left(\!\!
\begin{array}{c}
\overline{\sigma}
\\
\sigma
\\
\rho
\\
\overline{\zeta}
\\
\zeta
\end{array}
\!\!\right) := \underbrace{\left(\!\!
\begin{array}{ccccc}
{\sf a}\overline{\sf a}\overline{\sf a} & 0 & 0 & 0 & 0
\\
0 & {\sf a}{\sf a}\overline{\sf a} & 0 & 0 & 0
\\
\overline{\sf c} & {\sf c} & {\sf a}\overline{\sf a} & 0 & 0
\\
\overline{\sf e} & {\sf d} & \overline{\sf b} & \overline{\sf a} & 0
\\
\overline{\sf d} & {\sf e} & {\sf b} & 0 & {\sf a}
\end{array}
\!\!\right)}_{=:g} \left(\!\!
\begin{array}{c}
\overline{\sigma_0}
\\
\sigma_0
\\
\rho_0
\\
\overline{\zeta_0}
\\
\zeta_0
\end{array}
\!\!\right),
\endaligned
\end{equation}
that is to say:
\[
\aligned \overline{\sigma} & = {\sf a}\overline{\sf a}\overline{\sf
a}\, \overline{\sigma_0},
\\
\sigma & = {\sf a}{\sf a}\overline{\sf a}\,\sigma_0,
\\
\rho & = \overline{\sf c}\,\overline{\sigma_0} + {\sf c}\,\sigma_0 +
{\sf a}\overline{\sf a}\,\rho_0,
\\
\overline{\zeta} & = \overline{\sf e}\,\overline{\sigma_0} + {\sf
d}\,\sigma_0 + \overline{b}\,\rho_0 + \overline{\sf
a}\,\overline{\zeta_0},
\\
\zeta & = \overline{\sf d}\,\overline{\sigma_0} + {\sf e}\,\sigma_0 +
{\sf b}\,\rho_0 + {\sf a}\,\zeta_0.
\endaligned
\]
Again, the $1$-form $\rho$ is real and the $1$-forms
$\overline{\sigma}$ and $\overline{ \zeta}$ are the conjugate of
$\sigma$ and $\zeta$.

Now, our objective is to perform the equivalence method with these
general data, taking advantage of what has been already finalized
in the simpler case of the model in Section \ref{equivalence-model}.

\section{Absorption and normalization}
\label{abs-norm}
\HEAD{\ref{abs-norm}.~Absorption and normalization}{
Masoud {\sc Sabzevari} (Shahrekord) and Jo\"el {\sc Merker} (LM-Orsay)}

Proceeding exactly as in the beginning of Subsection~\ref{g-minus-1},
differentiating both sides of~\thetag{\ref{lift-coframe-general}}
yields in matrix notation:
\[
\footnotesize
\aligned
d\,\left(\!
\begin{array}{c}
\overline{\sigma} \\
\sigma \\
\rho \\
\overline{\zeta} \\
\zeta \\
\end{array}
\!\right)
&
=
dg\wedge\,\left(\!
\begin{array}{c}
\overline{\sigma}_0 \\
\sigma_0 \\
\rho_0 \\
\overline{\zeta}_0 \\
\zeta_0
\end{array}
\!\right)
+g\cdot\left(%
\begin{array}{c}
d\overline{\sigma}_0 \\
d\sigma_0 \\
d\rho_0 \\
d\overline{\zeta}_0 \\
d\zeta_0 \\
\end{array}%
\right)
\\
&
=
\underbrace{\left(
\begin{array}{ccccc}
2\,\overline{\alpha}_1+\alpha_1 & 0 & 0 & 0 & 0 \\
0 & 2\,\alpha_1+\overline{\alpha}_1 & 0 & 0 & 0 \\
\overline{\alpha}_2 & {\alpha}_2 & \alpha_1+\overline{\alpha}_1 & 0 & 0 \\
\overline{\alpha}_3 & \overline{\alpha}_4 & \overline{\alpha}_5 &
\overline{\alpha}_1 & 0 \\
{\alpha}_4 & {\alpha}_3 & {\alpha}_5 & 0 & {\alpha}_1 \\
\end{array}%
\right)}_{\omega_{\sf MC}:=dg\cdot g^{-1}}\,\,\,
\wedge\!\!\!\!\!
\underbrace{\left(%
\begin{array}{c}
\overline{\sigma} \\
\sigma \\
\rho \\
\overline{\zeta} \\
\zeta \\
\end{array}%
\right)}_{
g\cdot(\overline{\sigma}_0,\sigma_0,\rho_0,\overline{\zeta}_0,\zeta_0)^t}
+
\\
&
\ \ \ \ \
+
\left(\!\!
\begin{array}{c}
{\sf a}\overline{\sf a}\overline{\sf a}\,d\overline{\sigma}_0
\\
{\sf a}{\sf a}\overline{\sf a}\,d\sigma_0
\\
\overline{\sf c}\,d\overline{\sigma}_0 + {\sf c}\,d\sigma_0 + {\sf
a}\overline{\sf a}\,d\rho_0
\\
\overline{\sf e}\,d\overline{\sigma}_0+ {\sf
d}\,d\sigma_0+\overline{\sf b}\,d\rho_0+ \overline{\sf
a}\,d\overline{\zeta}_0
\\
\overline{\sf d}\,d\overline{\sigma}_0+{\sf e}\,d\sigma_0+ {\sf
b}\,d\rho_0+ {\sf a}\,d\zeta_0
\end{array}
\!\!\right).
\endaligned
\]
Of course here, the Maurer-Cartan forms are the same as the ones
computed before
when dealing with the model:
\[
\aligned {\alpha_1} & := \frac{d{\sf a}}{\sf a},
\\
{\alpha_2} &:= \frac{d{\sf c}}{{\sf a}^2\overline{\sf a}} -
\frac{{\sf c}\,d{\sf a}}{{\sf a}^3\overline{\sf a}} - \frac{{\sf
c}\,d\overline{\sf a}}{{\sf a}^2 \overline{\sf a}^2},
\\
{\alpha_3} & := -\,\frac{{\sf c}\,d{\sf b}}{{\sf a}^3\overline{\sf
a}^2} + \bigg( \frac{{\sf b}{\sf c}}{{\sf a}^4\overline{\sf a}^2} -
\frac{{\sf e}}{{\sf a}^3\overline{\sf a}} \bigg)\,d{\sf a} +
\frac{1}{{\sf a}^2\overline{\sf a}}\,d{\sf e},
\\
{\alpha_4} & := \frac{d\overline{\sf d}}{{\sf
a}\overline{\sf a}^2} - \frac{\overline{\sf c}\,d{\sf b}}{{\sf
a}^2\overline{\sf a}^3} + \bigg( \frac{{\sf b}\overline{\sf c}}{{\sf
a}^3\overline{\sf a}^3} - \frac{\overline{\sf d}}{{\sf
a}^2\overline{\sf a}^2} \bigg)\, d{\sf a},
\\
{\alpha_5} & := \frac{d{\sf b}}{{\sf a}\overline{\sf a}} - \frac{{\sf
b}\,d{\sf a}}{{\sf a}^2\overline{\sf a}}.
\endaligned
\]
Also, from Subsection~\ref{g-minus-1} again, we know that
the inverse of the general matrix $g$ of
our ambiguity group is:
\begin{eqnarray*}
g^{-1}=\left(%
\begin{array}{ccccc}
\frac{1}{{\sf a}\overline{\sf a}^2} & 0 & 0 & 0 & 0 \\
0 & \frac{1}{{\sf a}^2\overline{\sf a}} & 0 & 0 & 0 \\
-\frac{\overline{\sf c}}{{\sf a}^2\overline{\sf a}^3} & -\frac{\sf c}{{\sf
a}^3\overline{\sf a}^2} &
\frac{1}{{\sf a}\overline{\sf a}} & 0 & 0 \\
\frac{\overline{\sf b}\,\overline{\sf c}-\overline{\sf e}{\sf a}\overline{\sf
a}}{{\sf a}^2\overline{\sf a}^4} &
\frac{\overline{\sf b}\sf c-{\sf a}\overline{a}\sf d}{{\sf a}^3\overline{\sf a}^3}
& -\frac{\overline{\sf b}}{{\sf a}\overline{\sf a}^2} &
\frac{1}{\overline{\sf a}} & 0 \\
\frac{{\sf b}\overline{\sf c}-{\sf a}\overline{\sf a}\overline{\sf d}}{{\sf
a}^3\overline{\sf a}^3} &
\frac{{\sf bc}-{\sf ea}\overline{\sf a}}{{\sf a}^4\overline{\sf a}^2} & -\frac{\sf
b}{{\sf a}^2\overline{\sf a}} &
 0 & \frac{1}{{\sf a}} \\
\end{array}%
\right),
\end{eqnarray*}
hence we have the same inversion formulas as we had
when dealing with the model:
\begin{equation}
\label{0-frame-general} \footnotesize\aligned \sigma_0 & = \frac{1}{{\sf
a}^2\overline{\sf a}}\,\sigma,
\\
\rho_0 & = -\,\frac{\overline{\sf c}}{{\sf a}^2\overline{\sf a}^3}\,
\overline{\sigma} - \frac{\sf c}{{\sf a}^3\overline{\sf a}^2}\,\sigma
+ \frac{1}{{\sf a}\overline{\sf a}}\,\rho,
\\
\zeta_0 & = \frac{{\sf b}\overline{\sf c}-{\sf a}\overline{\sf
a}\overline{\sf d}}{ {\sf a}^3\overline{\sf a}^3}\,\overline{\sigma}
+ \frac{{\sf b}{\sf c}-{\sf a}\overline{\sf a}{\sf e}}{ {\sf
a}^4\overline{\sf a}^2}\,\sigma - \frac{{\sf b}}{{\sf
a}^2\overline{\sf a}}\,\rho + \frac{1}{{\sf a}}\,\zeta.
\endaligned
\end{equation}

\medskip

{\em However, at this precise point, computations start to be differ and to
become substantially harder.}

\medskip

Indeed, leaving aside the trivial:
\[
d\zeta_0
=
0,
\]
coming back to the two initial structure
equations~\thetag{\ref{d0-general}}
that were set up in the preceding section\,\,---, we
modify appropriately the order of appearance of terms\,\,---:
\begin{equation}
\label{struct-0-rew}
\aligned
d\sigma_0
=
-\,K\cdot
\sigma_0\wedge\overline{\sigma}_0
+
F\cdot\sigma_0\wedge\rho_0
+
Q\cdot\sigma_0\wedge\zeta_0
+
B\cdot\sigma_0\wedge\overline{\zeta}_0
+
\\
+
\overline{G}\cdot\overline{\sigma}_0\wedge\rho_0
+
B\cdot\overline{\sigma}_0\wedge\zeta_0
+
\overline{R}\cdot\overline{\sigma}_0\wedge\overline{\zeta}_0
+
\\
+
\rho_0\wedge\zeta_0,
\endaligned
\end{equation}
\[
\aligned
d\rho_0
=
-\,i\,J\cdot\sigma_0\wedge\overline{\sigma}_0
+
E\cdot\sigma_0\wedge\rho_0
+
P\cdot\sigma_0\wedge\zeta_0
+
A\cdot\sigma_0\wedge\overline{\zeta}_0
+
\\
+
\overline{E}\cdot\overline{\sigma}_0\wedge\rho_0
+
\overline{P}\cdot\overline{\sigma}_0\wedge\overline{\zeta}_0
+
A\cdot\overline{\sigma}_0\wedge\zeta_0
+
\\
+
i\,\zeta_0\wedge\overline{\zeta}_0,
\endaligned
\]
we must replace the so obtained values of $\sigma_0$,
$\overline{\sigma}_0$, $\rho_0$, $\zeta_0$, $\overline{\zeta}_0$ in
terms of $\sigma$, $\overline{\sigma}$, $\rho$, $\zeta$,
$\overline{\zeta}$. In the process, we organize the computations so
that group variables appear first, so that monomials ${\sf a}^\mu
\overline{\sf a}^\nu$ always land at denominator place, so that the
alphabetical order is respected in all numerators, and so that the
initial structure functions always appear after group variables.
Of course, we remember that the $5$ coefficient-functions:
\[
E,\ \ \ \ \ \ \
F,\ \ \ \ \ \ \
G,\ \ \ \ \ \ \
J,\ \ \ \ \ \ \
K
\]
express in fact in terms of the $5$ functions:
\[
P,\ \ \ \ \ \ \
Q,\ \ \ \ \ \ \
R,\ \ \ \ \ \ \
A,\ \ \ \ \ \ \
B
\]
and their frame derivatives up to order $2$, but
we prevent from inserting these expressions
at this stage, planning to perform the replacements later.

\medskip

We obtain:
\[
\footnotesize
\aligned
d\sigma_0
&
=
\sigma\wedge\overline{\sigma}\,
\bigg[
-\,
\frac{1}{{\sf a}^3\overline{\sf a}^3}\,K
-
\frac{\overline{\sf c}}{
{\sf a}^4\overline{\sf a}^4}\,
F
+
\frac{{\sf b}\overline{\sf c}}{
{\sf a}^5\overline{\sf a}^4}\,
Q
-
\frac{\overline{\sf d}}{
{\sf a}^4\overline{\sf a}^3}\,Q
+
\frac{\overline{\sf b}\overline{\sf c}}{
{\sf a}^4\overline{\sf a}^5}\,
B
-
\frac{\overline{\sf e}}{{\sf a}^3\overline{\sf a}^4}\,
B
+
\frac{{\sf c}}{{\sf a}^4\overline{\sf a}^4}\,
\overline{G}
-
\\
&
\ \ \ \ \ \ \ \ \ \ \ \ \ \ \ \ \
-
\frac{{\sf b}{\sf c}}{
{\sf a}^5\overline{\sf a}^4}\,
B
+
\frac{{\sf e}}{
{\sf a}^4\overline{\sf a}^3}\,
B
-
\frac{\overline{\sf b}{\sf c}}{
{\sf a}^4\overline{\sf a}^5}\,
\overline{R}
+
\frac{{\sf d}}{
{\sf a}^3\overline{\sf a}^4}\,
\overline{R}
+
\frac{{\sf c}\overline{\sf d}}{
{\sf a}^5\overline{\sf a}^4}
-
\frac{{\sf e}\overline{\sf c}}{
{\sf a}^5\overline{\sf a}^4}
\bigg]
+
\\
&
+
\sigma\wedge\rho\,
\bigg[
\frac{1}{{\sf a}^3\overline{\sf a}^2}\,F
-
\frac{{\sf b}}{
{\sf a}^4\overline{\sf a}^2}\,
Q
-
\frac{\overline{\sf b}}{
{\sf a}^3\overline{\sf a}^3}\,
B
+
\frac{{\sf e}}{
{\sf a}^4\overline{\sf a}^4}
\bigg]
+
\\
&
+
\sigma\wedge\zeta\,
\bigg[
\frac{1}{{\sf a}^3\overline{\sf a}}\,Q
-
\frac{{\sf c}}{{\sf a}^4\overline{\sf a}^2}
\bigg]
+
\sigma\wedge\overline{\zeta}\,
\bigg[
\frac{1}{{\sf a}^2\overline{\sf a}^2}\,B
\bigg]
+
\\
&
+
\overline{\sigma}\wedge\rho\,
\bigg[
\frac{1}{{\sf a}^2\overline{\sf a}^3}\,
\overline{G}
-
\frac{{\sf b}}{{\sf a}^3\overline{\sf a}^3}\,B
-
\frac{\overline{\sf b}}{{\sf a}^2\overline{\sf a}^4}\,
\overline{R}
+
\frac{\overline{\sf d}}{{\sf a}^3\overline{\sf a}^3}
\bigg]
+
\\
&
+
\overline{\sigma}\wedge\zeta\,
\bigg[
\frac{1}{{\sf a}^2\overline{\sf a}^2}\,B
-
\frac{\overline{\sf c}}{{\sf a}^3\overline{\sf a}^3}
\bigg]
+
\overline{\sigma}\wedge\overline{\zeta}\,
\bigg[
\frac{1}{{\sf a}\overline{\sf a}^3}\,\overline{R}
\bigg]
+
\rho\wedge\zeta\,
\bigg[
\frac{1}{{\sf a}^2\overline{\sf a}}
\bigg],
\endaligned
\]
and:
\[
\footnotesize
\aligned
d\rho_0
&
=
\sigma\wedge\overline{\sigma}\,
\bigg[
-\,i\,
\frac{1}{{\sf a}^3\overline{\sf a}^3}\,J
-
\frac{\overline{\sf c}}{{\sf a}^4\overline{\sf a}^4}\,E
+
\frac{{\sf b}\overline{\sf c}}{{\sf a}^5\overline{\sf a}^4}\,P
-
\frac{\overline{\sf d}}{{\sf a}^4\overline{\sf a}^3}\,P
+
\frac{\overline{\sf b}\overline{\sf c}}{{\sf a}^4\overline{\sf a}^5}\,A
-
\frac{\overline{\sf e}}{{\sf a}^3\overline{\sf a}^4}\,A
+
\\
&
\ \ \ \ \ \ \ \ \ \ \ \ \ \ \ \ \ \
+
\frac{{\sf c}}{{\sf a}^4\overline{\sf a}^4}\,
\overline{E}
-
\frac{\overline{\sf b}{\sf c}}{
{\sf a}^4\overline{\sf a}^5}\,
\overline{P}
+
\frac{{\sf d}}{{\sf a}^3\overline{\sf a}^4}\,
\overline{P}
-
\frac{{\sf b}{\sf c}}{{\sf a}^5\overline{\sf a}^4}\,A
+
\frac{{\sf e}}{{\sf a}^4\overline{\sf a}^3}\,A
-
i\,
\frac{{\sf b}{\sf c}\overline{\sf e}}{
{\sf a}^5\overline{\sf a}^5}
\,-
\\
&
\ \ \ \ \ \ \ \ \ \ \ \ \ \ \ \ \ \
-\,
i\,\frac{\overline{\sf b}\overline{\sf c}{\sf e}}{
{\sf a}^5\overline{\sf a}^5}
+
i\,
\frac{{\sf e}\overline{\sf e}}{
{\sf a}^4\overline{\sf a}^4}
+
i\,
\frac{{\sf b}\overline{\sf c}{\sf d}}{
{\sf a}^5\overline{\sf a}^5}
+
i\,
\frac{\overline{\sf b}{\sf c}\overline{\sf d}}{
{\sf a}^5\overline{\sf a}^5}
-
i\,
\frac{{\sf d}\overline{\sf d}}{
{\sf a}^4\overline{\sf a}^4}
\bigg]
+
\\
&
+
\sigma\wedge\rho
\bigg[
\frac{1}{{\sf a}^3\overline{\sf a}^2}\,E
-
\frac{{\sf b}}{{\sf a}^4\overline{\sf a}^2}\,P
-
\frac{\overline{\sf b}}{{\sf a}^3\overline{\sf a}^3}\,A
+
i\,
\frac{\overline{\sf b}{\sf e}}{{\sf a}^4\overline{\sf a}^3}
-
i\,
\frac{{\sf b}{\sf d}}{{\sf a}^4\overline{\sf a}^3}
\bigg]
+
\\
&
+
\sigma\wedge\zeta
\bigg[
\frac{1}{{\sf a}^3\overline{\sf a}}\,P
-
i\,
\frac{\overline{\sf b}{\sf c}}{{\sf a}^4\overline{\sf a}^3}
+
i\,
\frac{{\sf d}}{{\sf a}^3\overline{\sf a}^2}
\bigg]
+
\\
&
+
\sigma\wedge\overline{\zeta}\,
\bigg[
\frac{1}{{\sf a}^2\overline{\sf a}^2}\,A
+
i\,\frac{{\sf b}{\sf c}}{
{\sf a}^4\overline{\sf a}^3}
-
i\,\frac{{\sf e}}{{\sf a}^3\overline{\sf a}^2}
\bigg]
+
\\
&
+
\overline{\sigma}\wedge\rho\,
\bigg[
\frac{1}{{\sf a}^2\overline{\sf a}^3}\overline{E}
-
\frac{\overline{\sf b}}{{\sf a}^2\overline{\sf a}^4}\,
\overline{P}
-
\frac{{\sf b}}{{\sf a}^3\overline{\sf a}^3}\,\overline{A}
-
i\,
\frac{{\sf b}\overline{\sf e}}{
{\sf a}^3\overline{\sf a}^4}
+
i\,
\frac{\overline{\sf b}\overline{\sf d}}{
{\sf a}^3\overline{\sf a}^4}
\bigg]
+
\\
&
+
\overline{\sigma}\wedge\zeta\,
\bigg[
\frac{1}{{\sf a}^2\overline{\sf a}^2}\,A
-
i\,\frac{\overline{\sf b}\overline{\sf c}}{
{\sf a}^3\overline{\sf a}^4}
+
i\,
\frac{\overline{\sf e}}{
{\sf a}^2\overline{\sf a}^3}
\bigg]
+
\\
&
+
\overline{\sigma}\wedge\overline{\zeta}\,
\bigg[
\frac{1}{{\sf a}\overline{\sf a}^3}\,
\overline{P}
+
i\,
\frac{{\sf b}\overline{\sf c}}{{\sf a}^3\overline{\sf a}^4}
-
i\,
\frac{\overline{\sf d}}{{\sf a}^2\overline{\sf a}^3}
\bigg]
+
\\
&
+
\rho\wedge\zeta\,\bigg[
i\,\frac{\overline{\sf b}}{{\sf a}^2\overline{\sf a}^2}
\bigg]
+
\rho\wedge\overline{\zeta}\,
\bigg[
-\,i\,
\frac{{\sf b}}{{\sf a}^2\overline{\sf a}^2}
\bigg]
+
\frac{i}{{\sf a}\overline{\sf a}}\,
\zeta\wedge\overline{\sf a}.
\endaligned
\]

It yet remains to compute the three last terms:
\[
\aligned
&
{\sf a}^2\overline{\sf a}\,d\sigma_0,
\\
&
{\sf c}\,d\sigma_0+\overline{\sf c}\,
d\overline{\sigma}_0+{\sf a}\overline{\sf a}\,d\rho_0,
\\
&
{\sf e}\,d\sigma_0+\overline{\sf d}\,d\overline{\sigma}_0
+{\sf b}\,d\sigma_0+\zero{{\sf a}\,d\zeta_0},
\endaligned
\]
which happen to be complicated and to give
rise to torsion coefficients. We do this, and
this provides the complete structure equations:
\begin{equation}
\label{structure-firts-absorbtion}
 \aligned d\sigma & =
\big(2\,\alpha_1+\overline{\alpha}_1\big) \wedge \sigma+
\\
& \ \ \ \ \
+
U_1\,\sigma\wedge\overline{\sigma} +
U_2\,\sigma\wedge\rho + U_3\,\sigma\wedge\zeta +
U_4\,\sigma\wedge\overline{\zeta}
+
\\
&
\ \ \ \ \ \ \ \ \ \ \ \ \ \ \ \ \ \ \ \ \ \ \ \
+ U_5\,\overline{\sigma}\wedge\rho +
U_6\,\overline{\sigma}\wedge\zeta
+
U_7\,\overline{\sigma}\wedge\overline{\zeta}
+
\\
&
\ \ \ \ \ \ \ \ \ \ \ \ \ \ \ \ \ \ \ \ \ \
\ \ \ \ \ \ \ \ \ \ \ \ \ \ \ \ \ \ \ \ \
+ {\rho\wedge\zeta},
\endaligned
\end{equation}
\[
\aligned
d\rho & = \alpha_2\wedge\sigma +
\overline{\alpha}_2\wedge\overline{\sigma} + \alpha_1\wedge\rho +
\overline{\alpha}_1\wedge\rho +
\\
& \ \ \ \ \ +
V_1\,\sigma\wedge\overline{\sigma} +
V_2\,\sigma\wedge\rho + V_3\,\sigma\wedge\zeta +
V_4\,\sigma\wedge\overline{\zeta} +
\\
& \ \ \ \ \ \ \ \ \ \ \ \ \ \ \ \ \ \ \ \ \ \ \ \ +
\overline{V_2}\,\overline{\sigma}\wedge\rho +
\overline{V_4}\,\overline{\sigma}\wedge\zeta +
\overline{V_3}\,\overline{\sigma}\wedge\overline{\zeta} +
\\
& \ \ \ \ \ \ \ \ \ \ \ \ \ \ \ \ \ \ \ \ \ \ \ \ \ \ \ \ \ \ \ \ \ \
\ \ \ \ \ \ \ \ \ + V_8\,\rho\wedge\zeta +
\overline{V_8}\,\rho\wedge\overline{\zeta} +
\\
& \ \ \ \ \ \ \ \ \ \ \ \ \ \ \ \ \ \ \ \ \ \ \ \ \ \ \ \ \ \ \ \ \ \
\ \ \ \ \ \ \ \ \ \ \ \ \ \ \ \ \ \ \ \ \ \ \ \ \ \ \ \ \ \ +
i\,\zeta\wedge\overline{\zeta},
\endaligned
\]
\[
\aligned
d\zeta & = \alpha_3\wedge\sigma
+ \alpha_4\wedge\overline{\sigma} + \alpha_5\wedge\rho +
\alpha_1\wedge\zeta +
\\
& \ \ \ \ \ + W_1\,\sigma\wedge\overline{\sigma} +
W_2\,\sigma\wedge\rho + W_3\,\sigma\wedge\zeta +
W_4\,\sigma\wedge\overline{\zeta} +
\\
& \ \ \ \ \ \ \ \ \ \ \ \ \ \ \ \ \ \ \ \ \ \ \ \ +
W_5\,\overline{\sigma}\wedge\rho + W_6\,\overline{\sigma}\wedge\zeta
+ W_7\,\overline{\sigma}\wedge\overline{\zeta} +
\\
& \ \ \ \ \ \ \ \ \ \ \ \ \ \ \ \ \ \ \ \ \ \ \ \ \ \ \ \ \ \ \ \ \ \
\ \ \ \ \ \ \ \ \ + W_8\,\rho\wedge\zeta +
W_9\,\rho\wedge\overline{\zeta} +
\\
& \ \ \ \ \ \ \ \ \ \ \ \ \ \ \ \ \ \ \ \ \ \ \ \ \ \ \ \ \ \ \ \ \ \
\ \ \ \ \ \ \ \ \ \ \ \ \ \ \ \ \ \ \ \ \ \ \ \ \ \ \ \ \ +
W_{10}\,\zeta\wedge\overline{\zeta},
\endaligned
\]
expressed in terms of the lifted coframe,
where the appearing torsion coefficients are as follows.
For $d\sigma$:
\[
\footnotesize
\aligned
U_1
&
=
-\,\frac{1}{{\sf a}\overline{\sf a}^2}\,
K
-
\frac{\overline{\sf c}}{{\sf a}^2\overline{\sf a}^3}\,F
+
\frac{{\sf b}\overline{\sf c}}{{\sf a}^3\overline{\sf a}^3}\,Q
-
\frac{\overline{\sf d}}{{\sf a}^2\overline{\sf a}^2}\,Q
+
\frac{\overline{\sf b}\overline{\sf c}}{
{\sf a}^2\overline{\sf a}^4}\,B
-
\frac{\overline{\sf e}}{{\sf a}\overline{\sf a}^3}\,B
+
\frac{{\sf c}}{{\sf a}^2\overline{\sf a}^3}\,
\overline{G}
-
\\
&
\ \ \ \ \
-
\frac{{\sf b}{\sf c}}{{\sf a}^3\overline{\sf a}^3}\,B
+
\frac{{\sf e}}{{\sf a}^2\overline{\sf a}^2}\,B
-
\frac{\overline{\sf b}{\sf c}}{{\sf a}^2\overline{\sf a}^4}\,
\overline{R}
+
\frac{{\sf d}}{{\sf a}\overline{\sf a}^3}\,
\overline{R}
+
\frac{{\sf c}\overline{\sf d}}{{\sf a}^3\overline{\sf a}^3}
-
\frac{\overline{\sf c}{\sf e}}{{\sf a}^3\overline{\sf a}^3},
\\
U_2
&
=
\frac{1}{{\sf a}\overline{\sf a}}\,F
-
\frac{{\sf b}}{{\sf a}^2\overline{\sf a}}\,Q
-
\frac{\overline{\sf b}}{{\sf a}\overline{\sf a}^2}\,B
+
\frac{{\sf e}}{{\sf a}^2\overline{\sf a}},
\\
U_3
&
=
\frac{1}{{\sf a}}\,Q
-
\frac{{\sf c}}{{\sf a}^2\overline{\sf a}},
\\
U_4
&
=
\frac{1}{\overline{\sf a}}\,B,
\\
U_5
&
=
\frac{1}{\overline{\sf a}^2}\,\overline{G}
-
\frac{{\sf b}}{{\sf a}\overline{\sf a}^2}\,B
-
\frac{\overline{\sf b}}{\overline{\sf a}^3}\,
\overline{R}
+
\frac{\overline{\sf d}}{{\sf a}\overline{\sf a}^2},
\endaligned
\]
\[
\footnotesize
\aligned
U_6
&
=
\frac{1}{\overline{\sf a}}\,B
-
\frac{\overline{\sf c}}{{\sf a}\overline{\sf a}^2}, \ \ \ \ \ \ \ \ \ \ \ \ \ \ \ \ \ \ \ \ \ \ \ \ \ \ \ \ \ \ \ \ \ \ \ \ \ \ \ \ \ \ \ \ \ \ \ \ \ \ \ \ \ \ \ \ \ \ \ \ \ \ \ \ \ \ \ \ \ \ \ \ \ \ \ \ \ \ \ \ \ \ \ \ \ \ \ \ \ \ \ \ \ \ \ \ \ \ \ \ \ \ \ \ \ \ \ \ \ \ \ \
\\
U_7
&
=
\frac{{\sf a}}{\overline{\sf a}^2}\,\overline{R}.
\endaligned
\]

For $d\rho$:
\[
\footnotesize
\aligned
V_1
&
=
-\,\frac{{\sf c}}{{\sf a}^3\overline{\sf a}^3}\,K
-
\frac{{\sf c}\overline{\sf c}}{{\sf a}^4\overline{\sf a}^4}\,F
+
\frac{{\sf b}{\sf c}\overline{\sf c}}{
{\sf a}^5\overline{\sf a}^4}\,Q
-
\frac{{\sf c}\overline{\sf d}}{{\sf a}^4\overline{\sf a}^3}\,Q
+
\frac{\overline{\sf b}{\sf c}\overline{\sf c}}{
{\sf a}^4\overline{\sf a}^5}\,B
-
\frac{{\sf c}\overline{\sf e}}{{\sf a}^3\overline{\sf a}^4}\,B
+
\frac{{\sf c}{\sf c}}{{\sf a}^4\overline{\sf a}^4}\,\overline{G}
-
\\
&
\ \ \ \
-\,
\frac{{\sf b}{\sf c}{\sf c}}{{\sf a}^5\overline{\sf a}^4}\,B
+
\frac{{\sf c}{\sf e}}{{\sf a}^4\overline{\sf a}^3}\,B
-
\frac{\overline{\sf b}{\sf c}{\sf c}}{
{\sf a}^4\overline{\sf a}^5}\,\overline{R}
+
\frac{{\sf c}{\sf d}}{{\sf a}^3\overline{\sf a}^4}\,\overline{R}
+
\frac{{\sf c}{\sf c}\overline{\sf d}}{
{\sf a}^5\overline{\sf a}^4}
-
\frac{{\sf e}{\sf c}\overline{\sf c}}{
{\sf a}^5\overline{\sf a}^4}
+
\\
&
\ \ \ \ \
+
\frac{\overline{\sf c}}{{\sf a}^3\overline{\sf a}^3}\,
\overline{K}
+
\frac{{\sf c}\overline{\sf c}}{{\sf a}^4\overline{\sf a}^4}\,
\overline{F}
-
\frac{\overline{\sf b}{\sf c}\overline{\sf c}}{
{\sf a}^4\overline{\sf a}^5}\,\overline{Q}
+
\frac{\overline{\sf c}{\sf d}}{{\sf a}^3\overline{\sf a}^4}\,
\overline{Q}
-
\frac{{\sf b}{\sf c}\overline{\sf c}}{{\sf a}^5\overline{\sf a}^4}\,
\overline{B}
+
\frac{\overline{\sf c}{\sf e}}{{\sf a}^4\overline{\sf a}^3}\,
\overline{B}
-
\frac{\overline{\sf c}\overline{\sf c}}{
{\sf a}^4\overline{\sf a}^4}\,G
+
\\
&
\ \ \ \ \
+
\frac{\overline{\sf b}\overline{\sf c}\overline{\sf c}}{
{\sf a}^4\overline{\sf a}^5}\,
\overline{B}
-
\frac{\overline{\sf c}\overline{\sf e}}{
{\sf a}^3\overline{\sf a}^4}\,\overline{B}
+
\frac{{\sf b}\overline{\sf c}\overline{\sf c}}{
{\sf a}^5\overline{\sf a}^4}\,R
-
\frac{\overline{\sf c}\overline{\sf d}}{
{\sf a}^4\overline{\sf a}^3}\,R
-
\frac{\overline{\sf c}\overline{\sf c}{\sf d}}{
{\sf a}^4\overline{\sf a}^5}
+
\frac{{\sf c}\overline{\sf c}\overline{\sf e}}{
{\sf a}^4\overline{\sf a}^5}
-
\\
&
\ \ \ \ \
-\,i\,
\frac{1}{{\sf a}^2\overline{\sf a}^2}\,J
-
\frac{\overline{\sf c}}{{\sf a}^3\overline{\sf a}^3}\,E
+
\frac{{\sf b}\overline{\sf c}}{{\sf a}^4\overline{\sf a}^3}\,P
-
\frac{\overline{\sf d}}{{\sf a}^3\overline{\sf a}^2}\,P
+
\frac{\overline{\sf b}\overline{\sf c}}{{\sf a}^3\overline{\sf a}^4}\,A
-
\frac{\overline{\sf e}}{{\sf a}^2\overline{\sf a}^3}\,A
+
\frac{{\sf c}}{{\sf a}^3\overline{\sf a}^3}\,\overline{E}
-
\frac{\overline{\sf b}{\sf c}}{{\sf a}^3\overline{\sf a}^4}\,
\overline{P}
+
\frac{{\sf d}}{{\sf a}^2\overline{\sf a}^3}\,\overline{P}
-
\\
&
\ \ \ \ \
-
\frac{{\sf b}{\sf c}}{{\sf a}^4\overline{\sf a}^3}\,A
+
\frac{{\sf e}}{{\sf a}^3\overline{\sf a}^2}\,A
-
i\,
\frac{{\sf b}{\sf c}\overline{\sf e}}{
{\sf a}^4\overline{\sf a}^4}
-
i\,\frac{\overline{\sf b}\overline{\sf c}{\sf e}}{
{\sf a}^4\overline{\sf a}^4}
+
i\,
\frac{{\sf e}\overline{\sf e}}{{\sf a}^3\overline{\sf a}^3}
+
i\,\frac{{\sf b}\overline{\sf c}{\sf d}}{
{\sf a}^4\overline{\sf a}^4}
+
i\,
\frac{\overline{\sf b}{\sf c}\overline{\sf d}}{
{\sf a}^4\overline{\sf a}^4}
-
i\,
\frac{{\sf d}\overline{\sf d}}{
{\sf a}^3\overline{\sf a}^3},
\endaligned
\]
\[
\footnotesize
\aligned
V_2
&
=
\frac{{\sf c}}{{\sf a}^3\overline{\sf a}^2}\,F
-
\frac{{\sf b}{\sf c}}{{\sf a}^4\overline{\sf a}^2}\,Q
-
\frac{\overline{\sf b}{\sf c}}{{\sf a}^3\overline{\sf a}^3}\,B
+
\frac{{\sf c}{\sf e}}{{\sf a}^4\overline{\sf a}^2}
+
\frac{\overline{\sf c}}{{\sf a}^3\overline{\sf a}^2}\,G
-
\frac{\overline{\sf b}\overline{\sf c}}{
{\sf a}^3\overline{\sf a}^3}\,\overline{B}
-
\frac{{\sf b}\overline{\sf c}}{{\sf a}^4\overline{\sf a}^2}\,R
+
\frac{\overline{\sf c}{\sf d}}{
{\sf a}^3\overline{\sf a}^3}
+
\\
&
\ \ \ \ \
+
\frac{1}{{\sf a}^2\overline{\sf a}}\,E
-
\frac{{\sf b}}{{\sf a}^3\overline{\sf a}}\,P
-
\frac{\overline{\sf b}}{{\sf a}^2\overline{\sf a}^2}\,A
+
i\,\frac{\overline{\sf b}{\sf e}}{
{\sf a}^3\overline{\sf a}^2}
-
i\,\frac{{\sf b}{\sf d}}{
{\sf a}^3\overline{\sf a}^2},
\\
V_3
&
=
\frac{{\sf c}}{{\sf a}^3\overline{\sf a}}\,Q
-
\frac{{\sf c}{\sf c}}{{\sf a}^4\overline{\sf a}^2}
+
\frac{\overline{\sf c}}{{\sf a}^3\overline{\sf a}}\,R
+
\frac{1}{{\sf a}^2}\,P
-
i\,
\frac{\overline{\sf b}{\sf c}}{
{\sf a}^3\overline{\sf a}^2}
+
i\,\frac{{\sf d}}{{\sf a}^2\overline{\sf a}},
\\
V_4
&
=
\frac{{\sf c}}{{\sf a}^2\overline{\sf a}^2}\,B
+
\frac{\overline{\sf c}}{{\sf a}^2\overline{\sf a}^2}\,
\overline{B}
-
\frac{{\sf c}\overline{\sf c}}{
{\sf a}^3\overline{\sf a}^3}
+
\frac{1}{{\sf a}\overline{\sf a}}\,A
+
i\,
\frac{{\sf b}{\sf c}}{{\sf a}^3\overline{\sf a}^2}
-
i\,\frac{{\sf e}}{{\sf a}^2\overline{\sf a}},
\\
V_8
&
=
\frac{{\sf c}}{{\sf a}^2\overline{\sf a}}
+
i\,\frac{\overline{\sf b}}{{\sf a}\overline{\sf a}}.
\endaligned
\]
Lastly, for $d\zeta$:
\[
\footnotesize
\aligned
W_1
&
=
-\,
\frac{{\sf e}}{{\sf a}^3\overline{\sf a}^3}\,K
-
\frac{\overline{\sf c}{\sf e}}{{\sf a}^4\overline{\sf a}^4}\,
F
+
\frac{{\sf b}\overline{\sf c}{\sf e}}{{\sf a}^5\overline{\sf a}^4}\,Q
-
\frac{\overline{\sf d}{\sf e}}{{\sf a}^4\overline{\sf a}^3}\,Q
+
\frac{\overline{\sf b}\overline{\sf c}{\sf e}}{
{\sf a}^4\overline{\sf a}^5}\,B
-
\frac{{\sf e}\overline{\sf e}}{{\sf a}^3\overline{\sf a}^4}\,B
+
\frac{{\sf c}{\sf e}}{{\sf a}^4\overline{\sf a}^4}\,
\overline{G}
-
\\
&
\ \ \ \ \
-\,
\frac{{\sf b}{\sf c}{\sf e}}{{\sf a}^5\overline{\sf a}^4}\,B
+
\frac{{\sf e}{\sf e}}{{\sf a}^4\overline{\sf a}^3}\,B
-
\frac{\overline{\sf b}{\sf c}{\sf e}}{{\sf a}^4\overline{\sf a}^5}\,
\overline{R}
+
\frac{{\sf d}{\sf e}}{{\sf a}^3\overline{\sf a}^4}\,\overline{R}
+
\frac{{\sf c}\overline{\sf d}{\sf e}}{{\sf a}^5\overline{\sf a}^4}
-
\frac{\overline{\sf c}{\sf e}{\sf e}}{
{\sf a}^5\overline{\sf a}^4}
+
\\
&
\ \ \ \ \
+
\frac{\overline{\sf d}}{{\sf a}^3\overline{\sf a}^3}\,
\overline{K}
+
\frac{{\sf c}\overline{\sf d}}{{\sf a}^4\overline{\sf a}^4}\,
\overline{F}
-
\frac{\overline{\sf b}{\sf c}\overline{\sf d}}{
{\sf a}^4\overline{\sf a}^5}\,
\overline{Q}
+
\frac{{\sf d}\overline{\sf d}}{{\sf a}^3\overline{\sf a}^4}\,
\overline{Q}
-
\frac{{\sf b}{\sf c}\overline{\sf d}}{{\sf a}^5
\overline{\sf a}^4}\,\overline{B}
+
\frac{\overline{\sf d}{\sf e}}{{\sf a}^4\overline{\sf a}^3}\,
\overline{B}
-
\frac{\overline{\sf c}\overline{\sf d}}{{\sf a}^4\overline{\sf a}^4}\,G
+
\\
&
\ \ \ \ \
+
\frac{\overline{\sf b}\overline{\sf c}\overline{\sf d}}{
{\sf a}^4\overline{\sf a}^5}\,\overline{B}
-
\frac{\overline{\sf d}\overline{\sf e}}{{\sf a}^3\overline{\sf a}^4}\,
\overline{B}
+
\frac{{\sf b}\overline{\sf c}\overline{\sf d}}{
{\sf a}^5\overline{\sf a}^4}\,R
-
\frac{\overline{\sf d}\overline{\sf d}}{{\sf a}^4\overline{\sf a}^3}\,
R
-
\frac{\overline{\sf c}{\sf d}\overline{\sf d}}{
{\sf a}^4\overline{\sf a}^5}
+
\frac{{\sf c}\overline{\sf d}\overline{\sf e}}{
{\sf a}^4\overline{\sf a}^5}
-
\\
&
\ \ \ \ \
-\,i\,
\frac{{\sf b}}{{\sf a}^3\overline{\sf a}^3}\,J
-
\frac{{\sf b}\overline{\sf c}}{{\sf a}^4\overline{\sf a}^4}\,E
+
\frac{{\sf b}{\sf b}\overline{\sf c}}{{\sf a}^5\overline{\sf a}^4}\,P
-
\frac{{\sf b}\overline{\sf d}}{{\sf a}^4\overline{\sf a}^3}\,P
+
\frac{{\sf b}\overline{\sf b}\overline{\sf c}}{{\sf a}^4
\overline{\sf a}^5}\,A
-
\frac{{\sf b}\overline{\sf e}}{{\sf a}^3\overline{\sf a}^4}\,A
+
\frac{{\sf b}{\sf c}}{{\sf a}^4\overline{\sf a}^4}\,\overline{E}
-
\frac{{\sf b}\overline{\sf b}{\sf c}}{
{\sf a}^4\overline{\sf a}^5}\,\overline{P}
+
\frac{{\sf b}{\sf d}}{{\sf a}^3\overline{\sf a}^4}\,
\overline{P}
-
\\
&
\ \ \ \ \
-\,
\frac{{\sf b}{\sf b}{\sf c}}{{\sf a}^5\overline{\sf a}^4}\,A
+
\frac{{\sf b}{\sf e}}{{\sf a}^4\overline{\sf a}^3}\,A
-
i\,
\frac{{\sf b}{\sf b}{\sf c}\overline{\sf e}}{
{\sf a}^5\overline{\sf a}^5}
-
i\,
\frac{{\sf b}\overline{\sf b}\overline{\sf c}{\sf e}}{
{\sf a}^5\overline{\sf a}^5}
+
i\,
\frac{{\sf b}{\sf e}\overline{\sf e}}{
{\sf a}^4\overline{\sf a}^4}
+
i\,
\frac{{\sf b}{\sf b}\overline{\sf c}{\sf d}}{
{\sf a}^5\overline{\sf a}^5}
+
i\,
\frac{{\sf b}\overline{\sf b}{\sf c}\overline{\sf d}}{
{\sf a}^5\overline{\sf a}^5}
-
i\,\frac{{\sf b}{\sf d}\overline{\sf d}}{
{\sf a}^4\overline{\sf a}^4},
\\
W_2
&
=
\frac{{\sf e}}{{\sf a}^3\overline{\sf a}^2}\,F
-
\frac{{\sf b}{\sf e}}{{\sf a}^4\overline{\sf a}^2}\,Q
-
\frac{\overline{\sf b}{\sf e}}{{\sf a}^3\overline{\sf a}^3}\,B
+
\frac{{\sf e}{\sf e}}{{\sf a}^4\overline{\sf a}^2}
+
\frac{\overline{\sf d}}{{\sf a}^3\overline{\sf a}^2}\,G
-
\frac{\overline{\sf b}\overline{\sf d}}{
{\sf a}^3\overline{\sf a}^3}\,\overline{B}
-
\frac{{\sf b}\overline{\sf d}}{
{\sf a}^4\overline{\sf a}^2}\,R
+
\frac{{\sf d}\overline{\sf d}}{
{\sf a}^3\overline{\sf a}^3}
+
\\
&
\ \ \ \ \
+
\frac{{\sf b}}{{\sf a}^3\overline{\sf a}^2}\,E
-
\frac{{\sf b}{\sf b}}{{\sf a}^4\overline{\sf a}^2}\,P
-
\frac{{\sf b}\overline{\sf b}}{{\sf a}^3\overline{\sf a}^3}\,A
+
i\,
\frac{{\sf b}\overline{\sf b}{\sf e}}{
{\sf a}^4\overline{\sf a}^3}
-
i\,
\frac{{\sf b}{\sf b}{\sf d}}{
{\sf a}^4\overline{\sf a}^3},
\endaligned
\]
\[
\footnotesize
\aligned
W_3
&
=
\frac{{\sf e}}{{\sf a}^3\overline{\sf a}}\,Q
-
\frac{{\sf c}{\sf e}}{{\sf a}^4\overline{\sf a}^2}
+
\frac{\overline{\sf d}}{{\sf a}^3\overline{\sf a}}\,R
+
\frac{{\sf b}}{{\sf a}^3\overline{\sf a}}\,P
-
i\,\frac{{\sf b}\overline{\sf b}{\sf c}}{
{\sf a}^4\overline{\sf a}^3}
+
i\,
\frac{{\sf b}{\sf d}}{{\sf a}^3\overline{\sf a}^2},
\\
W_4
&
=
\frac{{\sf e}}{{\sf a}^2\overline{\sf a}^2}\,B
+
\frac{\overline{\sf d}}{{\sf a}^2\overline{\sf a}^2}\,
\overline{B}
-
\frac{{\sf c}\overline{\sf d}}{
{\sf a}^3\overline{\sf a}^3}
+
\frac{{\sf b}}{{\sf a}^2\overline{\sf a}^2}\,A
+
i\,\frac{{\sf b}{\sf b}{\sf c}}{
{\sf a}^4\overline{\sf a}^3}
-
i\,
\frac{{\sf b}{\sf e}}{
{\sf a}^3\overline{\sf a}^2},
\\
 W_5 & = \frac{{\sf e}}{{\sf a}^2\overline{\sf
a}^3}\, \overline{G} - \frac{{\sf b}{\sf e}}{ {\sf a}^3\overline{\sf
a}^3}\,B - \frac{\overline{\sf b}{\sf e}}{ {\sf a}^2\overline{\sf
a}^4}\, \overline{R} + \frac{\overline{\sf d}{\sf e}}{ {\sf
a}^3\overline{\sf a}^3} + \frac{\overline{\sf d}}{{\sf
a}^2\overline{\sf a}^3}\, \overline{F} - \frac{\overline{\sf
b}\overline{\sf d}}{ {\sf a}^2\overline{\sf a}^4}\, \overline{Q} -
\frac{{\sf b}\overline{\sf d}}{ {\sf a}^3\overline{\sf a}^3}\,
\overline{B} + \frac{\overline{\sf e}\overline{\sf d}}{ {\sf
a}^2\overline{\sf a}^4} +
\\
& \ \ \ \ \ + \frac{{\sf b}}{{\sf a}^2\overline{\sf a}^3}\,
\overline{E} - \frac{{\sf b}\overline{\sf b}}{{\sf a}^2\overline{\sf
a}^4}\, \overline{P} - \frac{{\sf b}{\sf b}}{ {\sf a}^3\overline{\sf
a}^3}\,A - i\, \frac{{\sf b}{\sf b}\overline{\sf e}}{{\sf a}^3
\overline{\sf a}^4} + i\, \frac{{\sf b}\overline{\sf b}\overline{\sf
d}}{ {\sf a}^3\overline{\sf a}^4},
\\
W_6
&
=
\frac{{\sf e}}{{\sf a}^2\overline{\sf a}^2}\,B
-
\frac{\overline{\sf c}{\sf e}}{{\sf a}^3\overline{\sf a}^3}
+
\frac{\overline{\sf d}}{{\sf a}^2\overline{\sf a}^2}\,
\overline{B}
+
\frac{{\sf b}}{{\sf a}^2\overline{\sf a}^2}\,A
-
i\,
\frac{{\sf b}\overline{\sf b}\overline{\sf c}}{
{\sf a}^3\overline{\sf a}^4}
+
i\,
\frac{{\sf b}\overline{\sf e}}{
{\sf a}^2\overline{\sf a}^3},
\\
W_7
&
=
\frac{{\sf e}}{{\sf a}\overline{\sf a}^3}\,
\overline{R}
+
\frac{\overline{\sf d}}{{\sf a}\overline{\sf a}^3}\,
\overline{Q}
-
\frac{\overline{\sf c}\overline{\sf d}}{
{\sf a}^2\overline{\sf a}^4}
+
\frac{{\sf b}}{{\sf a}\overline{\sf a}^3}\,
\overline{P}
+
i\,\frac{{\sf b}{\sf b}\overline{\sf c}}{
{\sf a}^3\overline{\sf a}^4}
-
i\,
\frac{{\sf b}\overline{\sf d}}{
{\sf a}^2\overline{\sf a}^3},
\\
W_8
&
=
\frac{{\sf e}}{{\sf a}^2\overline{\sf a}}
+
i\,
\frac{{\sf b}\overline{\sf b}}{{\sf a}^2\overline{\sf a}^2},
\\
W_9
&
=
\frac{\overline{\sf d}}{{\sf a}\overline{\sf a}^2}
-
i\,
\frac{{\sf b}{\sf b}}{
{\sf a}^2\overline{\sf a}^2},
\\
W_{10}
&
=
i\,\frac{{\sf b}}{{\sf a}\overline{\sf a}}.
\endaligned
\]

\subsection{First loop absorbtion}

Similarly as when
we treated the model
in Subsection~\ref{first-loop-model},
according to Proposition~\ref{prop-changes}, we
must modify the five $1$-forms $\alpha_1$, $\alpha_2$, $\alpha_3$,
$\alpha_4$, $\alpha_5$ by adding to them general
linear combinations of the
$1$-forms $\sigma$, $\overline{ \sigma}$, $\rho$, $\zeta$,
$\overline{ \zeta}$:
\[
\aligned
\alpha_1
&
\longmapsto
\alpha_1
+
p_1\,\sigma
+
q_1\,\overline{\sigma}
+
r_1\,\rho
+
s_1\,\zeta
+
t_1\,\overline{\zeta},
\\
\alpha_2
&
\longmapsto
\alpha_2
+
p_2\,\sigma
+
q_2\,\overline{\sigma}
+
r_2\,\rho
+
s_2\,\zeta
+
t_2\,\overline{\zeta},
\\
\alpha_3
&
\longmapsto
\alpha_3
+
p_3\,\sigma
+
q_3\,\overline{\sigma}
+
r_3\,\rho
+
s_3\,\zeta
+
t_3\,\overline{\zeta},
\\
\alpha_4
&
\longmapsto
\alpha_4
+
p_4\,\sigma
+
q_4\,\overline{\sigma}
+
r_4\,\rho
+
s_4\,\zeta
+
t_4\,\overline{\zeta},
\\
\alpha_5
&
\longmapsto
\alpha_5
+
p_5\,\sigma
+
q_5\,\overline{\sigma}
+
r_5\,\rho
+
s_5\,\zeta
+
t_5\,\overline{\zeta},
\endaligned
\]
with $25$ arbitrary real analytic functions $p_i$, $q_i$, $r_i$,
$s_i$, $t_i$. Then the expressions of $d\sigma$, $d\rho$, $d\zeta$
become\,\,---\,\,mind now that, contrary to the model case, the two torsion
coefficients $U_4$ and $U_7$ do not vanish\,\,---:
\[
\!\!\!\!\!\!\!\!\!\!\!\!\!\!\!\!\!\!\!\!\!\!\!\!\!\!\!\!\!\!\!
\small
\aligned
d\sigma
&
=
\big(2\alpha_1+\overline{\alpha}_1\big)
\wedge
\sigma
+
\\
&
+
\sigma\wedge\overline{\sigma}
\big[U_1-2q_1-\overline{p}_1\big]
+
\sigma\wedge\rho
\big[
U_2-2r_1-\overline{r}_1
\big]
+
\sigma\wedge\zeta
\big[
U_3-2s_1-\overline{t}_1
\big]
+
\sigma\wedge\overline{\zeta}
\big[
U_4-2t_1-\overline{s}_1
\big]
\\
&
\ \ \ \ \ \ \ \ \ \ \ \ \ \ \ \ \ \ \ \ \ \ \ \ \ \ \ \ \ \ \ \ \ \ \
\ \ \ \ \ \ \,
+
\overline{\sigma}\wedge\rho
\big[U_5\big]
+
\overline{\sigma}\wedge\zeta
\big[U_6\big]
+
\overline{\sigma}\wedge\overline{\zeta}
\big[U_7\big]
+
\\
&
\ \ \ \ \ \ \ \ \ \ \ \ \ \ \ \ \ \ \ \ \ \ \ \ \ \ \ \ \ \ \ \ \ \ \
\ \ \ \ \ \ \ \ \ \ \ \ \ \ \ \ \ \ \ \ \ \ \ \ \ \ \ \ \,
+
\rho\wedge\zeta,
\\
d\rho
&
=
\alpha_2\wedge\sigma
+
\overline{\alpha}_2\wedge\overline{\sigma}
+
\alpha_1\wedge\rho
+
\overline{\alpha}_1\wedge\rho
+
\\
&
\ \ \ \ \
+
\sigma\wedge\overline{\sigma}
\big[
V_1-q_2-\overline{q}_2
\big]
+
\sigma\wedge\rho
\big[
V_2-r_2+p_1+\overline{q}_1
\big]
+
\sigma\wedge\zeta
\big[
V_3-s_2
\big]
+
\sigma\wedge\overline{\zeta}
\big[
V_4-t_2
\big]
\\
&
\ \ \ \ \ \ \ \ \ \ \ \ \ \ \ \ \ \ \ \ \ \ \ \ \ \ \ \ \ \ \ \ \ \ \
+
\overline{\sigma}\wedge\rho
\big[
\overline{V}_2-\overline{r}_2+q_1+\overline{p}_1
\big]
+
\overline{\sigma}\wedge\zeta
\big[
\overline{V}_4-\overline{t}_2
\big]
+
\overline{\sigma}\wedge\overline{\zeta}
\big[
\overline{V}_3-\overline{s}_2
\big]
+
\endaligned
\]
\[
\!\!\!\!\!\!\!\!\!\!\!\!\!\!\!\!\!\!\!\!\!\!\!\!\!\!\!\!\!\!\!
\small
\aligned
&
+
\rho\wedge\zeta
\big[
V_8-s_1-\overline{t}_1
\big]
+
\rho\wedge\overline{\zeta}
\big[
\overline{V}_8-t_1-\overline{s}_1
\big]
+ \ \ \ \ \ \ \ \ \ \ \ \ \ \ \ \ \ \ \ \ \ \ \ \ \ \ \ \ \ \ \ \ \ \ \ \ \ \ \ \ \ \ \ \ \ \ \ \ \ \ \ \
\\
&
+
i\,\zeta\wedge\overline{\zeta},
\endaligned
\]
\[
\!\!\!\!\!\!\!\!\!\!\!\!\!\!\!\!\!\!\!\!\!\!\!\!\!\!\!\!\!\!\!
\small
\aligned
d\zeta
&
=
\alpha_3\wedge\sigma
+
\alpha_4\wedge\overline{\sigma}
+
\alpha_5\wedge\rho
+
\alpha_1\wedge\zeta
+
\\
&
\ \ \ \ \
+
\sigma\wedge\overline{\sigma}
\big[
W_1-q_3+p_4
\big]
+
\sigma\wedge\rho
\big[
W_2-r_3+p_5
\big]
+
\sigma\wedge\zeta
\big[
W_3-s_3+p_1
\big]
+
\sigma\wedge\overline{\zeta}
\big[
W_4-t_3
\big]
\\
&
\ \ \ \ \ \ \ \ \ \ \ \ \ \ \ \ \ \ \ \ \ \ \ \ \ \ \ \ \ \ \ \ \ \ \
+
\overline{\sigma}\wedge\rho
\big[
W_5-r_4+q_5
\big]
+
\overline{\sigma}\wedge\zeta
\big[
W_6-s_4+q_1
\big]
+
\overline{\sigma}\wedge\overline{\zeta}
\big[
W_7-t_4
\big]
+
\\
&
\ \ \ \ \ \ \ \ \ \ \ \ \ \ \ \ \ \ \ \ \ \ \ \ \ \ \ \ \ \ \ \ \ \ \
\ \ \ \ \ \ \ \ \ \ \ \ \ \ \ \ \ \ \ \ \ \ \ \ \ \ \ \ \ \ \ \ \ \ \
+
\rho\wedge\zeta
\big[
W_8-s_5+r_1
\big]
+
\rho\wedge\overline{\zeta}
\big[
W_9-t_5
\big]
+
\\
&
\ \ \ \ \ \ \ \ \ \ \ \ \ \ \ \ \ \ \ \ \ \ \ \ \ \ \ \ \ \ \ \ \ \ \
\ \ \ \ \ \ \ \ \ \ \ \ \ \ \ \ \ \ \ \ \ \ \ \ \ \ \ \ \ \ \ \ \ \ \
\ \ \ \ \ \ \ \ \ \ \ \ \ \ \ \ \ \ \ \ \ \ \ \ \ \ \ \ \ \ \ \ \ \ \
\ \ \
+
\zeta\wedge\overline{\zeta}
\big[
W_{10}-t_1
\big].
\endaligned
\]

\subsection{First loop normalization}
\label{First loop normalization}

 In order to know what are the
precise linear combinations of the $22$ torsion coefficients:
\[
\aligned
&
U_1,\ \ \ \ \ U_2,\ \ \ \ \ U_3,\ \ \ \ \ U_4,\ \ \ \ \ U_5,
\ \ \ \ \ U_6,\ \ \ \ \ U_7,
\\
&
V_1,\ \ \ \ \ V_2,\ \ \ \ \ V_3,\ \ \ \ \ V_4,\ \ \ \ \ V_8,
\\
&
W_1,\ \ \ \ \ W_2,\ \ \ \ \ W_3,\ \ \ \ \ W_4,\ \ \ \ \
W_5,\ \ \ \ \ W_6,\ \ \ \ \ W_7,\ \ \ \ \ W_8, \ \ \ \ \
W_9,\ \ \ \ \ W_{10}
\endaligned
\]
that are {\em necessarily normalizable},
one must determine all possible linear
combinations of the following $7 + 5 + 10 = 22$
equations\,\,---\,\,including their (unwritten) conjugates\,\,---:
\[\footnotesize\aligned
\left[
\aligned
U_1
&
=
2q_1+\overline{p}_1,
\\
U_2
&
=
2r_1+\overline{r}_1,
\\
U_3
&
=
2s_1+\overline{t}_1,
\\
U_4
&
=
2t_1+\overline{s}_1,
\\
U_5
&
=
0,
\\
U_6
&
=
0,
\\
U_7
&
=
0,
\endaligned\right.
\ \ \ \ \ \ \ \ \ \ \ \ \ \ \ \ \ \ \ \ \
\left[
\aligned
V_1
&
=
q_2-\overline{q}_2,
\\
V_2
&
=
r_2-p_1-\overline{q}_1,
\\
V_3
&
=
s_2,
\\
V_4
&
=
t_2,
\\
V_8
&
=
s_1+\overline{t}_1,
\endaligned\right.
\ \ \ \ \ \ \ \ \ \ \ \ \ \ \ \ \ \ \ \ \
\left[
\aligned
W_1
&
=
q_3-p_4,
\\
W_2
&
=
r_3-p_5,
\\
W_3
&
=
s_3-p_1,
\\
W_4
&
=
t_3,
\\
W_5
&
=
r_4-q_5,
\\
W_6
&
=
s_4-q_1,
\\
W_7
&
=
t_4,
\\
W_8
&
=
s_5-r_1,
\\
W_9
&
=
t_5,
\\
W_{10}
&
=
t_1,
\endaligned\right.
\endaligned
\]
so as to obtain null right-hand sides. Visually, some complete
appropriate linear combinations are:
\begin{equation}
\label{U-5-6-7}
\aligned
0
&
=
U_5,
\\
0
&
=
U_6,
\\
0
&
=
U_7,
\\
0
&
=
U_3+\overline{U}_4-3\,V_8,
\\
0
&
=
\overline{U}_4-V_8-\overline{W}_{10},
\endaligned
\end{equation}
and that is all (exercise).
Now, if one just replaces the appearing torsion coefficients,
one plainly obtains {\em five} normalizable linear
combinations:
\begin{equation}
\label{normalizable-first loop}
\footnotesize\boxed{\aligned
U_5
&
=
\frac{1}{\overline{\sf a}^2}\,\overline{G}
-
\frac{\sf b}{{\sf a}\overline{\sf a}^2}\,B
-
\frac{\overline{\sf b}}{\overline{\sf a}^3}\,
\overline{R}
+
\frac{\overline{\sf d}}{{\sf a}\overline{\sf a}^2},
\\
U_6
&
=
\frac{1}{\overline{\sf a}}\,B
-
\frac{\overline{\sf c}}{{\sf a}\overline{\sf a}^2},
\\
U_7
&
=
\frac{{\sf a}}{\overline{\sf a}^2}\,\overline{R},
\\
U_3+\overline{U}_4-3\,V_8
&
=
\frac{1}{\sf a}\,Q
-
4\,\frac{\sf c}{{\sf a}^2\overline{\sf a}}
+
\frac{1}{\sf a}\,\overline{B}
-
3i\,
\frac{\overline{\sf b}}{{\sf a}\overline{\sf a}},
\\
\overline{U}_4-V_8-\overline{W}_{10}
&
=
\frac{1}{\sf a}\,\overline{B}
-
\frac{\sf c}{{\sf a}^2\overline{\sf a}}.
\endaligned}
\end{equation}
Visibly, the second and the fifth combinations are conjugate
of each other, {\em hence the fifth can be removed}.
A confirmation of computational correctness is yielded by the
following obvious

\begin{Observation}
One recovers all equations along with the equivalence process
applied to the model in Section~\ref{equivalence-model} just by
assigning the value zero to all of the functions:
\[
\aligned
P,\ \ \ \ \
Q,\ \ \ \ \
R,
\\
A,\ \ \ \ \
B,
\\
E,\ \ \ \ \
F,\ \ \ \ \
G,
\\
J,\ \ \ \ \
K,
\endaligned
\]
appearing in the structure equations~\thetag{
\ref{struct-0-rew}} of the initial
coframe $\big\{ \sigma_0, \overline{ \sigma}_0,
\rho_0, \zeta_0, \overline{\zeta}_0 \big\}$.
\qed
\end{Observation}

But then, in the model case, instead of the third normalizable
expression above:
\[
\frac{\sf a}{\overline{\sf a}^2}\,
\overline{R},
\]
we had only the trivial combination `$0$', just because $R = 0$ in the
model case. It is thus necessary and unavoidable to {\em distinguish
two branches} in the future issue of the equivalence procedure:

\medskip\noindent$\square$
either $R \equiv 0$ as a function on our
CR-generic maximally minimal submanifold $M^5 \subset \C^4$;

\medskip\noindent$\square$
or else, $R \not \equiv 0$, so that after relocalizing
if necessary the consideration at another Zariski-generic central point
(shifting the origin of the coordinate system), we
may assume that $R \not \equiv 0$ vanishes at no point.

\medskip
Of course, in this second case, we will be able to normalize
the complex parameter ${\sf a}$, but let us
continue at first the procedure of equivalence under the
assumption that $R \equiv 0$, a `branch' which is closer to the model case.

\bigskip

\section{The branch $R \equiv 0$}
\label{branch-R-0}
\HEAD{\ref{branch-R-0}.~The branch $R \equiv 0$}{
Masoud {\sc Sabzevari} (Shahrekord) and Jo\"el {\sc Merker} (LM-Orsay)}

\subsection{Normalization of three group parameters}
\label{Normalization-1-r=0}
Assuming therefore that $R \equiv 0$
vanishes identically, the above five normalizable expressions
\thetag{\ref{normalizable-first loop}} reduce to exactly three:
\begin{equation}
\label{normalizable-first loop-R=0}
\aligned
U_5
&
=
\frac{1}{\overline{\sf a}^2}\,\overline{G}
-
\frac{\sf b}{{\sf a}\overline{\sf a}^2}\,B
-
\zero{\frac{\overline{\sf b}}{\overline{\sf a}^3}\,
\overline{R}}
+
\frac{\overline{\sf d}}{{\sf a}\overline{\sf a}^2},
\\
U_6
&
=
\frac{1}{\overline{\sf a}}\,B
-
\frac{\overline{\sf c}}{{\sf a}\overline{\sf a}^2},
\\
U_3+\overline{U}_4-3\,V_8
&
=
\frac{1}{\sf a}\,Q
-
4\,\frac{\sf c}{{\sf a}^2\overline{\sf a}}
+
\frac{1}{\sf a}\,\overline{B}
-
3i\,
\frac{\overline{\sf b}}{{\sf a}\overline{\sf a}}.
\endaligned
\end{equation}
Equating the second expression to zero then specifies the
expression of the group parameter $\sf c$ as:
\[
\boxed{{\sf c}
:=
{\sf a}\overline{\sf a}\,\overline{B}.}
\]
This changes the third expression into a simplified form:
\[
\frac{1}{\sf a}\,Q
-
3\,\frac{1}{\sf a}\,\overline{B}
-
3i\,\frac{\overline{\sf b}}{{\sf a}\overline{\sf a}},
\]
and then, equating it similarly to zero normalizes the expression of $\sf
b$ as:
\[
\boxed{{\sf b}
:=
{\sf a}\,
\big(
-i\,B
+
{\textstyle{\frac{i}{3}}}\,\overline{Q}
\big).
}
\]
Lastly, putting this into the first expression
of~\thetag{\ref{normalizable-first loop-R=0}} and equating it to zero,
we also determine the expression of $\sf d$:
\[
{\sf d}
=
\overline{\sf a}\,
\big(
-\,G
+
i\,\overline{B}\,\overline{B}
-
{\textstyle{\frac{i}{3}}}\,\overline{B}\,Q
\big).
\]
But if we remember at this stage
from Lemma~\ref{E-F-G}
that the function $G$ has an expression
in terms of $P$, $Q$, $R$, $A$, $B$,
then still under the assumption
that $R = 0$, the normalization
of ${\sf d}$ is in fact more precisely:
\[
\boxed{
{\sf d}
=
\overline{\sf a}\,
\big(
-\,i\,\mathcal{L}\big(\overline{B}\big)
+
i\,P
+
{\textstyle{\frac{2i}{3}}}\,\overline{B}\,Q
\big).
}
\]

\subsection{Changing up the initial coframe}

Now, if we insert these normalized values
of ${\sf b}$, ${\sf c}$, ${\sf d}$
into the
expressions of the $1$-forms of our lifted coframe:
\[
\aligned
\sigma
&
=
{\sf a}^2\overline{\sf a}\,\sigma_0,
\\
\rho
&
=
{\sf c}\cdot\sigma_0
+
\overline{\sf c}\cdot\overline{\sigma}_0
+
{\sf a}\overline{\sf a}\cdot\rho_0,
\\
\zeta
&
=
{\sf e}\cdot\sigma_0
+
\overline{\sf d}\cdot\overline{\sigma}_0
+
{\sf b}\cdot\rho_0
+
{\sf a}\cdot\zeta_0,
\endaligned
\]
we realize after reorganization:
\[
\aligned
\sigma
&
=
{\sf a}^2\overline{\sf a}
\cdot
\sigma_0,
\\
\rho
&
=
{\sf a}\overline{\sf a}
\cdot
\big(
\underbrace{\overline{B}\,\sigma_0
+
B\,\overline{\sigma}_0
+
\rho_0}_{=:\,\rho_0^\sim}
\big),
\\
\zeta
&
=
{\sf e}\cdot\sigma_0
+
{\sf a}
\cdot
\big[
\underbrace{
\big(-\overline{G}-i\,BB
+
{\textstyle{\frac{i}{3}}}\,B\overline{Q}
\big)\,\overline{\sigma_0}
+
\big(-i\,B+
{\textstyle{\frac{i}{3}}}\,\overline{Q}\big)\cdot\rho_0
+
\zeta_0}_{=:\,\zeta_0^\sim}
\big]
\endaligned
\]
that it is natural to introduce the two new {\em modified initial
$1$-forms}:
\[
\rho_0^\sim
\ \ \ \ \ \ \ \ \ \ \ \ \ \ \ \
\text{\rm and}
\ \ \ \ \ \ \ \ \ \ \ \ \ \ \ \
\zeta_0^\sim
\]
in terms of which the {\em reduced} lifted coframe rewrites:
\[
\left[
\aligned
\sigma
&
=
{\sf a}^2\overline{\sf a}
\cdot
\sigma_0,
\\
\rho
&
=
{\sf a}\overline{\sf a}
\cdot
\rho_0^\sim,
\\
\zeta
&
=
{\sf e}\cdot\sigma_0
+
{\sf a}\cdot\zeta_0^\sim,
\endaligned\right.
\]
so that the corresponding {\em reduced} matrix group becomes:
\[
G^\sim
:=
\left\{
g^\sim
=
\left(\!\!
\begin{array}{ccccc}
{\sf a}\overline{\sf a}^2 & 0 & 0 & 0 & 0
\\
0 & {\sf a}^2\overline{\sf a} & 0 & 0 & 0
\\
0 & 0 & {\sf a}\overline{\sf a} & 0 & 0
\\
\overline{\sf e} & 0 & 0 & \overline{\sf a} & 0
\\
0 & {\sf e} & 0 & 0 & {\sf a}
\end{array}
\!\!\right)
\colon\
{\sf a},\,{\sf e}\in\C
\right\};
\]
of course, it is immediate that $(\overline\sigma_0, \sigma_0,
\rho_0^\sim, \overline\zeta_0^\sim,  \zeta_0^\sim)$ still constitutes
a coframe on our generic submanifold $M$. Furthermore, a simple
computation shows that:
\[
dg^\sim
\cdot
{g^\sim}^{-1}
=
\left(\!\!
\begin{array}{ccccc}
2\overline{\beta}_1+\beta_1 & 0 & 0 & 0 & 0
\\
0 & 2\beta_1+\overline{\beta}_1 & 0 & 0 & 0
\\
0 & 0 & \beta_1 + \overline{\beta}_1 & 0 & 0
\\
\overline{\beta}_2 & 0 & 0 & \overline{\beta}_1 & 0
\\
0 & \beta_2 & 0 & 0 & \beta_1
\end{array}
\!\!\right),
\]
where:
\[
\beta_1 := \frac{d{\sf a}}{\sf a} \ \ \ \ \ \ \ \ \ \ \ \ \ \ \ \ \ \
\text{\rm and} \ \ \ \ \ \ \ \ \ \ \ \ \ \ \ \ \ \ \beta_2 :=
\frac{d{\sf e}}{{\sf a}^2\overline{\sf a}}-\frac{\sf e}{{\sf
a}^3\overline{\sf a}}\,d{\sf a}.
\]

\subsection{Second-loop absorbtion and normalization}
\label{Second-loop-R=0}

Now, let us pursue the computations with our initial coframe $\big\{
\sigma_0, \overline{ \sigma}_0, \rho_0, \zeta_0, \overline{ \zeta}_0
\big\}$.  We have to replace the normalized values of ${\sf b}$,
${\sf c}$, ${\sf d}$ obtained above in the definition of the lifted
coframe.  Let us abbreviate these three normalizations as:
\[
\aligned
{\sf b}
&
:=
{\sf a}\,{\bf B}_0,
\ \ \  \ \ \ \ \ \ \ \ \ \ \ \ \ \ \ \ \
\text{\rm where:}\ \ \ \ \
\boxed{\,{\bf B}_0
:=
-\,i\,B
+
{\textstyle{\frac{i}{3}}}\,\overline{Q},\,}
\\
{\sf c}
&
:=
{\sf a}\overline{\sf a}\,{\bf C}_0,
\ \ \  \ \ \ \ \ \ \ \ \ \ \ \ \ \ \
\text{\rm where:}\ \ \ \ \
\boxed{\,{\bf C}_0
:=
\overline{B},\,}
\\
{\sf d}
&
:=
\overline{\sf a}\,{\bf D}_0,
\ \ \  \ \ \ \ \ \ \ \ \ \ \ \ \ \ \ \ \
\text{\rm where:}\ \ \ \ \
\boxed{\,{\bf D}_0
:=
-\,i\,
\mathcal{L}\big(\overline{B}\big)
+
i\,P
+
{\textstyle{\frac{2i}{3}}}\,
\overline{B}\,Q,\,}
\endaligned
\]
in terms of three new functions ${\bf B}_0$, ${\bf C}_0$,
${\bf D}_0$ defined on the base manifold $M^5 \subset \C^4$.
Then the new lifted coframe becomes:
\[
\aligned
\sigma
&
=
{\sf a}^2\overline{\sf a}\cdot\sigma_0,
\\
\rho
&
=
{\sf a}\overline{\sf a}\cdot
\big(
{\bf C}_0\,\sigma_0+\overline{\bf C}_0\,\overline{\sigma}_0+\rho_0
\big),
\\
\zeta
&
=
{\sf e}\cdot\sigma_0
+
{\sf a}\cdot
\big(
\overline{\bf D}_0\,\overline{\sigma}_0
+
{\bf B}_0\,\rho_0
+
\zeta_0
\big).
\endaligned
\]
Now, we apply the exterior differentiation operator and we obtain, in
terms of the Maurer-Cartan forms $\beta_1$ and $\beta_2$:
\begin{equation}
\label{set-up-1}
 \aligned d\sigma & =
\big(2\beta_1+\overline{\beta}_1\big) \wedge \sigma +
\\
&
\ \ \ \ \
+
{\sf a}^2\overline{\sf a}
\cdot
d\sigma_0,
\\
d\rho
&
=
\big(\beta_1+\overline{\beta}_1\big)
\wedge
\rho
+
\\
&
\ \ \ \ \
+
{\sf a}\overline{\sf a}\cdot
\big(
{\bf C}_0\,d\sigma_0
+
\overline{\bf C}_0\,d\overline{\sigma}_0
+
d\rho_0
+
\\
&
\ \ \ \ \ \ \ \ \ \ \ \ \ \ \ \ \ \
+
d{\bf C}_0\wedge\sigma_0
+
d\overline{\bf C}_0\wedge\overline{\sigma}_0\big),
\\
d\zeta
&
=
\beta_2\wedge\sigma
+
\beta_1\wedge\zeta
+
\\
&
\ \ \ \ \
+
{\sf e}\cdot d\sigma_0
+
{\sf a}\cdot
\big(
\overline{\bf D}_0\,d\overline{\sigma}_0
+
{\bf B}_0\,d\rho_0
+
d\zeta_0
+
\\
&
\ \ \ \ \ \ \ \ \ \ \ \ \ \ \ \ \ \ \ \ \ \ \ \ \ \ \ \ \ \ \
+
d\overline{\bf D}_0\wedge\overline{\sigma}_0
+
d{\bf B}_0\wedge\rho_0\big).
\endaligned
\end{equation}
One readily observes there is no change in the structure equation of
$d\sigma$, hence the torsion coefficients $U_\bullet$ will be
essentially unchanged. Intentionally here, each right-hand side is organized in
three lines having distinct meanings. On the first line, only
Maurer-Cartan forms appear. On the second line, only exterior
derivatives of the base $1$-forms appear. One easily convinces
oneself that the contribution of these second-line terms to the new
torsion coefficients $U_i'$, $V_j'$, $W_k'$:
\[
\aligned d\sigma& = \big(2\,\beta_1+\overline{\beta}_1\big) \wedge
\sigma+
\\
& \ \ \ \ \ + U_1'\,\sigma\wedge\overline{\sigma} +
U_2'\,\sigma\wedge\rho + U_3'\,\sigma\wedge\zeta +
U_4'\,\sigma\wedge\overline{\zeta} +
\\
& \ \ \ \ \ \ \ \ \ \ \ \ \ \ \ \ \ \ \ \ \ \ \ \ +
U_5'\,\overline{\sigma}\wedge\rho + U_6'\,\overline{\sigma}\wedge\zeta
+
U_7'\,\overline{\sigma}\wedge\overline{\zeta} +
\\
&
\ \ \ \ \ \ \ \ \ \ \ \ \ \ \ \ \ \ \ \ \ \ \ \ \ \ \ \ \ \
\ \ \ \ \ \ \ \ \ \ \ \ \ \ \ \ \ \ \ \ \ \ \ \ \
+ {\rho\wedge\zeta},
\endaligned
\]
\[
\aligned d\rho & = \big(\beta_1+\overline{\beta}_1\big)\wedge\rho +
\\
& \ \ \ \ \ + V_1'\,\sigma\wedge\overline{\sigma} +
V_2'\,\sigma\wedge\rho + V_3'\,\sigma\wedge\zeta +
V_4'\,\sigma\wedge\overline{\zeta} +
\\
& \ \ \ \ \ \ \ \ \ \ \ \ \ \ \ \ \ \ \ \ \ \ \ \ +
\overline{V_2}'\,\overline{\sigma}\wedge\rho +
\overline{V_4}'\,\overline{\sigma}\wedge\zeta +
\overline{V_3}'\,\overline{\sigma}\wedge\overline{\zeta} +
\\
& \ \ \ \ \ \ \ \ \ \ \ \ \ \ \ \ \ \ \ \ \ \ \ \ \ \ \ \ \ \ \ \ \ \
\ \ \ \ \ \ \ \ \ + V_8'\,\rho\wedge\zeta +
\overline{V_8}'\,\rho\wedge\overline{\zeta} +
\\
& \ \ \ \ \ \ \ \ \ \ \ \ \ \ \ \ \ \ \ \ \ \ \ \ \ \ \ \ \ \ \ \ \ \
\ \ \ \ \ \ \ \ \ \ \ \ \ \ \ \ \ \ \ \ \ \ \ \ \ \ \ \ \ \ \ \ +
i\,\zeta\wedge\overline{\zeta},
\endaligned
\]
\[
\aligned d\zeta & = \beta_2\wedge\sigma + \beta_1\wedge\zeta +
\\
& \ \ \ \ \ + W_1'\,\sigma\wedge\overline{\sigma} +
W_2'\,\sigma\wedge\rho +
W_3'\,\sigma\wedge\zeta +
W_4'\,\sigma\wedge\overline{\zeta} +
\\
& \ \ \ \ \ \ \ \ \ \ \ \ \ \ \ \ \ \ \ \ \ \ \ \ +
W_5'\,\overline{\sigma}\wedge\rho +
W_6'\,\overline{\sigma}\wedge\zeta+
W_7'\,\overline{\sigma}\wedge\overline{\zeta} +
\\
& \ \ \ \ \ \ \ \ \ \ \ \ \ \ \ \ \ \ \ \ \ \ \ \ \ \ \ \ \ \ \ \ \ \
\ \ \ \ \ \ \ \ \ \ \ \ \ +
W_8'\,\rho\wedge\zeta +
W_9'\,\rho\wedge\overline{\zeta} +
\\
& \ \ \ \ \ \ \ \ \ \ \ \ \ \ \ \ \ \ \ \ \ \ \ \ \ \ \ \ \ \ \ \ \ \
\ \ \ \ \ \ \ \ \ \ \ \ \ \ \ \ \ \ \ \ \ \ \ \ \ \ \ \ \ \ \ \ \ \ \
\ +
W_{10}'\,\zeta\wedge\overline{\zeta},
\endaligned
\]
just consists in taking back the previous torsion coefficients
$U_i$, $V_j$, $W_k$ shown above and replacing the
values of ${\sf b}$, ${\sf c}$, ${\sf d}$ by their
normalized values:
\[
\aligned U_i'\big\vert_{{\sf replace}({\sf b},{\sf c},{\sf d})} & =
\mathmotsf{second line of}\,\big(d\sigma) \ \ \ \ \ \ \ \ \ \ \ \ \
{\scriptstyle{(i\,=\,1,\,2,\,3,\,4,\,5,\,6,\,7)}},
\\
V_j\big\vert_{{\sf replace}({\sf b},{\sf c},{\sf d})} & =
\mathmotsf{second line of}\,\big(d\rho) \ \ \ \ \ \ \ \ \ \ \ \ \
{\scriptstyle{(j\,=\,1,\,2,\,3,\,4,\,8)}},
\\
W_k\big\vert_{{\sf replace}({\sf b},{\sf c},{\sf d})} & =
\mathmotsf{second line of}\,\big(d\zeta) \ \ \ \ \ \ \ \ \ \ \ \ \
{\scriptstyle{(k\,=\,1,\,2,\,3,\,4,\,5,\,6,\,7,\,8,\,9,\,10)}}.
\endaligned
\]

Now, to compute the third lines in the new torsion coefficients, we
apply a lemma, the proof of which is a direct consequence of a formula
shown above for the inverse $g^{ -1}$ of a general element $g$ in our
ambiguity matrix group.

\begin{Lemma}
\label{Lemma-G-0}
The exterior differential:
\[
d{\bf G}_0
=
\mathcal{S}({\bf G}_0)\cdot\sigma_0
+
\overline{\mathcal{S}}({\bf G}_0)\cdot\overline{\sigma}_0
+
\mathcal{T}({\bf G}_0)\cdot\rho_0
+
\mathcal{L}({\bf G}_0)\cdot\zeta_0
+
\overline{\mathcal{L}}({\bf G}_0)\cdot\overline{\zeta}_0
\]
of any function ${\bf G}_0$ on the base
manifold $M^5 \subset \C^4$ re-expresses,
in terms of the lifted coframe, as:
\[
\!\!\!\!\!\!\!\!\!\!\!\!\!\!\!\!\!\!\!\!
\footnotesize
\aligned
d{\bf G}_0
&
=
\sigma
\cdot
\bigg(
\frac{1}{{\sf a}^2\overline{\sf a}}\,
\mathcal{S}({\bf G}_0)
-
\frac{{\sf c}}{{\sf a}^3\overline{\sf a}^2}\,
\mathcal{T}({\bf G}_0)
+
\frac{{\sf b}{\sf c}}{{\sf a}^4\overline{\sf a}^2}\,
\mathcal{L}({\bf G}_0)
-
\frac{{\sf e}}{{\sf a}^3\overline{\sf a}}\,
\mathcal{L}({\bf G}_0)
+
\frac{\overline{\sf b}{\sf c}}{{\sf a}^3\overline{\sf a}^3}\,
\overline{\mathcal{L}}({\bf G}_0)
-
\frac{{\sf d}}{{\sf a}^2\overline{\sf a}^2}\,
\overline{\mathcal{L}}({\bf G}_0)
\bigg)
+
\\
&
\ \ \ \ \
+
\overline{\sigma}
\cdot
\bigg(
\frac{1}{{\sf a}\overline{\sf a}^2}\,
\overline{\mathcal{S}}({\bf G}_0)
-
\frac{\overline{\sf c}}{{\sf a}^2\overline{\sf a}^3}\,
\mathcal{T}({\bf G}_0)
+
\frac{{\sf b}\overline{\sf c}}{{\sf a}^3\overline{\sf a}^3}\,
\mathcal{L}({\bf G}_0)
-
\frac{\overline{\sf d}}{{\sf a}^2\overline{\sf a}^2}\,
\mathcal{L}({\bf G}_0)
+
\frac{\overline{\sf b}\overline{\sf c}}{{\sf a}^2\overline{\sf a}^4}\,
\overline{\mathcal{L}}({\bf G}_0)
-
\frac{\overline{\sf e}}{{\sf a}\overline{\sf a}^3}\,
\overline{\mathcal{L}}({\bf G}_0)
\bigg)
+
\\
&
\ \ \ \ \
+
\rho\cdot
\bigg(
\frac{1}{{\sf a}\overline{\sf a}}\,
\mathcal{T}({\bf G}_0)
-
\frac{{\sf b}}{{\sf a}^2\overline{\sf a}}\,
\mathcal{L}({\bf G}_0)
-
\frac{\overline{\sf b}}{{\sf a}\overline{\sf a}^2}\,
\overline{\mathcal{L}}({\bf G}_0)
\bigg)
+
\\
&
\ \ \ \ \
+
\zeta\cdot
\bigg(
\frac{1}{{\sf a}}\,
\mathcal{L}({\bf G}_0)
\bigg)
+
\\
&
\ \ \ \ \
+
\overline{\zeta}\cdot
\bigg(
\frac{1}{\overline{\sf a}}\,
\overline{\mathcal{L}}({\bf G}_0)
\bigg),
\endaligned
\]
before any normalization of the coefficients.
\qed
\end{Lemma}

Then all the new torsion coefficients may be computed
completely. Although the computations start to become quite
substantial, an accessible focused computation, left to the reader,
yields the expression of the new fourth torsion coefficient in
$d\rho$:
\[
V_4' = -\,\frac{1}{{\sf a}\overline{\sf a}}\,
\overline{\mathcal{L}}(\overline{B}) + \frac{1}{{\sf a}\overline{\sf
a}}\,A + \frac{2}{{\sf a}\overline{\sf a}}\, B\overline{B} -
\frac{1}{3}\,\frac{1}{{\sf a}\overline{\sf a}}\,
\overline{B}\,\overline{Q} - i\, \frac{{\sf e}}{{\sf
a}^2\overline{\sf a}}.
\]
In just a while, we will show that this torsion coefficient is
essential, hence normalizable.  Assigning to it the value $0$, we
deduce that the group parameter ${\sf e}$ can always be normalized as:
\[
\boxed{
{\sf e}
:=
{\sf a}\cdot
\big(
i\,\overline{\mathcal{L}}(\overline{B})
-
i\,A
-
2i\,B\overline{B}
+
{\textstyle{\frac{i}{3}}}\,
\overline{B}\,\overline{Q}
\big).
}
\]
We may also abbreviate the appearing auxiliary-normalizing
function defined on the base manifold $M^5 \subset \C^4$ as:
\[
{\sf e}
:=
{\sf a}\cdot{\bf E}_0
\ \ \  \ \ \ \ \ \ \ \ \ \ \ \ \ \ \ \ \
\text{\rm where:}\ \ \ \ \
\boxed{\,{\bf E}_0
:=
i\,\overline{\mathcal{L}}(\overline{B})
-
i\,A
-
2i\,B\overline{B}
+
{\textstyle{\frac{i}{3}}}\,
\overline{B}\,\overline{Q}.\,}
\]

\subsection{New torsion coefficients.}
In general, adding second lines to third lines by means of
Lemma~\ref{Lemma-G-0} just above, we obtain temporary (unclosed,
unfinished) formulas for the new torsion coefficients.

At first, in the $U_i$, plain replacements of the values of ${\sf b}$,
${\sf c}$, ${\sf d}$ must be done:
\[
U_i' := U_i \big\vert_{{\sf replace}({\sf b},{\sf c},{\sf d})} \ \ \
\ \ \ \ \ \ \ \ \ \
{\scriptstyle{(i\,=\,1,\,2,\,3,\,4,\,5,\,6,\,7)}}.
\]
But for the $V_i$ and for the $W_k$, supplementary terms coming from
the exterior differentiations of the three auxiliary-normalizing
functions ${\bf B}_0$, ${\bf C}_0$, ${\bf D}_0$ appear
in the concerned third lines. Applying the
lemma, we obtain formulas in which one must replace ${\sf b}$, ${\sf
c}$, ${\sf d}$ afterwards:
\begin{equation}
\!\!\!\!\!\!\!\!\!\!\!\!\!\!\!\!\!\!\!\!\!\!\!\!\!\!\!\!\!
\label{new} \footnotesize \aligned V_1' := V_1 & - \frac{1}{{\sf
a}^2\overline{\sf a}^2}\, \overline{\mathcal{S}}\big({\bf C}_0\big) +
\frac{\overline{\sf c}}{{\sf a}^3\overline{\sf a}^3}\,
\mathcal{T}\big({\bf C}_0\big) - \frac{{\sf b}\overline{\sf c}}{{\sf
a}^4\overline{\sf a}^3}\, \mathcal{L}\big({\bf C}_0\big) +
\frac{\overline{\sf d}}{{\sf a}^3\overline{\sf a}^2}\,
\mathcal{L}\big({\bf C}_0\big) - \frac{\overline{\sf b}\overline{\sf
c}}{{\sf a}^3\overline{\sf a}^4} \overline{\mathcal{L}}\big({\bf
C}_0\big) + \frac{\overline{\sf e}}{{\sf a}^2\overline{\sf a}^3}\,
\overline{\mathcal{L}}\big({\bf C}_0\big) +
\\
&
+
\frac{1}{{\sf a}^2\overline{\sf a}^2}\,
\mathcal{S}\big(\overline{\bf C}_0\big)
-
\frac{{\sf c}}{{\sf a}^3\overline{\sf a}^3}\,
\mathcal{T}\big(\overline{\bf C}_0\big)
+
\frac{{\sf b}{\sf c}}{{\sf a}^4\overline{\sf a}^3}\,
\mathcal{L}\big(\overline{\bf C}_0\big)
-
\frac{{\sf e}}{{\sf a}^3\overline{\sf a}^2}\,
\mathcal{L}\big(\overline{\bf C}_0\big)
+
\frac{\overline{\sf b}{\sf c}}{{\sf a}^3\overline{\sf a}^4}\,
\overline{\mathcal{L}}\big(\overline{\bf C}_0\big)
-
\frac{{\sf d}}{{\sf a}^2\overline{\sf a}^3}\,
\overline{\mathcal{L}}\big(\overline{\bf C}_0\big),
\endaligned
\end{equation}
\[
\footnotesize \aligned V_2' & := V_2 - \frac{1}{{\sf
a}^2\overline{\sf a}}\, \mathcal{T}\big({\bf C}_0\big) + \frac{{\sf
b}}{{\sf a}^3\overline{\sf a}}\, \mathcal{L}\big({\bf C}_0\big) +
\frac{\overline{\sf b}}{{\sf a}^2\overline{\sf a}^2}\,
\overline{\mathcal{L}}\big({\bf C}_0\big),
\\
V_3' & := V_3 - \frac{1}{{\sf a}^2}\, \mathcal{L}\big({\bf
C}_0\big),
\\
V_4' & := V_4 - \frac{1}{{\sf a}\overline{\sf a}}\,
\overline{\mathcal{L}}\big({\bf C}_0\big),
\\
V_8' & := V_8,
\endaligned
\]
\[
\!\!\!\!\!\!\!\!\!\!\!\!\!\!\!\!\!\!\!\!\!\!\!\!\!\!\!\!\!
\footnotesize \aligned W_1' & := W_1 + \frac{1}{{\sf
a}^2\overline{\sf a}^3}\, \mathcal{S}\big(\overline{\bf D}_0\big) -
\frac{{\sf c}}{{\sf a}^3\overline{\sf a}^4}\,
\mathcal{T}\big(\overline{\bf D}_0\big) + \frac{{\sf b}{\sf c}}{{\sf
a}^4\overline{\sf a}^4}\, \mathcal{L}\big(\overline{\bf D}_0\big) -
\frac{{\sf e}}{{\sf a}^3\overline{\sf a}^3}\,
\mathcal{L}\big(\overline{\bf D}_0\big) + \frac{\overline{\sf b}{\sf
c}}{{\sf a}^3\overline{\sf a}^5}\,
\overline{\mathcal{L}}\big(\overline{\bf D}_0\big) - \frac{{\sf
d}}{{\sf a}^2\overline{\sf a}^4}\,
\overline{\mathcal{L}}\big(\overline{\bf D}_0\big) -
\\
&
\ \ \ \ \ \ \ \ \ \
-
\frac{\overline{\sf c}}{{\sf a}^3\overline{\sf a}^4}\,
\mathcal{S}\big({\bf B}_0\big)
+
\frac{{\sf c}}{{\sf a}^3\overline{\sf a}^4}\,
\overline{\mathcal{S}}\big({\bf B}_0\big)
+
\frac{\overline{\sf c}{\sf e}}{{\sf a}^4\overline{\sf a}^4}\,
\mathcal{L}\big({\bf B}_0\big)
-
\frac{{\sf c}\overline{\sf d}}{{\sf a}^4\overline{\sf a}^4}\,
\mathcal{L}\big({\bf B}_0\big)
+
\frac{\overline{\sf c}{\sf d}}{{\sf a}^3\overline{\sf a}^5}\,
\overline{\mathcal{L}}\big({\bf B}_0\big)
-
\frac{{\sf c}\overline{\sf e}}{{\sf a}^3\overline{\sf a}^5}\,
\overline{\mathcal{L}}\big({\bf B}_0\big),
\endaligned
\]
\[
\footnotesize \aligned W_2' & := W_2 + \frac{1}{{\sf
a}^2\overline{\sf a}^2}\, \mathcal{S}\big({\bf B}_0\big) - \frac{{\sf
e}}{{\sf a}^3\overline{\sf a}^2}\, \mathcal{L}\big({\bf B}_0\big) -
\frac{{\sf d}}{{\sf a}^2\overline{\sf a}^3}\,
\overline{\mathcal{L}}\big({\bf B}_0\big),
\\
W_3' & := W_3 + \frac{{\sf c}}{{\sf a}^3\overline{\sf a}^2}\,
\mathcal{L}\big({\bf B}_0\big),
\\
W_4' & := W_4 + \frac{{\sf c}}{{\sf a}^2\overline{\sf a}^3}\,
\overline{\mathcal{L}}\big({\bf B}_0\big),
\\
 W_5' & := W_5 - \frac{1}{{\sf a}\overline{\sf
a}^3}\, \mathcal{T}\big(\overline{\bf D}_0\big) + \frac{{\sf b}}{{\sf
a}^2\overline{\sf a}^3}\, \mathcal{L}\big(\overline{\bf D}_0\big) +
\frac{\overline{\sf b}}{{\sf a}\overline{\sf a}^4}\,
\overline{\mathcal{L}}\big(\overline{\bf D}_0\big) + \frac{1}{{\sf
a}\overline{\sf a}^3}\, \overline{\mathcal{S}}\big({\bf B}_0\big) -
\frac{\overline{\sf d}}{{\sf a}^2\overline{\sf a}^3}\,
\mathcal{L}\big({\bf B}_0\big) - \frac{\overline{\sf e}}{{\sf
a}\overline{\sf a}^4}\, \overline{\mathcal{L}}\big({\bf B}_0\big),
\\
W_6' & := W_6 - \frac{1}{{\sf a}\overline{\sf a}^2}\,
\mathcal{L}\big(\overline{\bf D}_0\big) + \frac{\overline{\sf
c}}{{\sf a}^2\overline{\sf a}^3}\, \mathcal{L}\big({\bf B}_0\big),
\\
W_7' & := W_7 - \frac{1}{\overline{\sf a}^3}\,
\overline{\mathcal{L}}\big(\overline{\bf D}_0\big) +
\frac{\overline{\sf c}}{{\sf a}\overline{\sf a}^4}\,
\overline{\mathcal{L}}\big({\bf B}_0\big),
\\
 W_8' & := W_8 - \frac{1}{{\sf a}\overline{\sf
a}}\, \mathcal{L}\big({\bf B}_0\big),
\\
W_9' & := W_9 - \frac{1}{\overline{\sf a}^2}\,
\overline{\mathcal{L}}\big({\bf B}_0\big),
\\
W_{10}' & := W_{10}.
\endaligned
\]

Now, a direct computation shows that among the above torsion
coefficients, $V'_3$ vanishes identically as soon as we substitute
the obtained expressions of $\bf{B}_0, \bf{C}_0$ and $\bf{D}_0$.

\begin{Lemma}
After determining the group parameters $\sf b,c,d$, the torsion
coefficient $V'_3$ vanishes identically.
\end{Lemma}

\proof Putting the obtained expressions of $\sf b,c,d$ in the
expression of $V_3$ gives:
\[
V_3=\frac{1}{{\sf a}^2}\mathcal L(\overline B).
\]
Now, it remains only to subtract $\frac{1}{{\sf a}^2}\mathcal L({\bf
C}_0)$ and easily see the vanishing of $V'_3$.
\endproof

\subsection{Normalizable essential torsion combinations}
As is known from the general Cartan equivalence
procedure, we must next modify the two remaining Maurer-Cartan
$1$-forms $\beta_1 = \alpha_1$ and $\beta_2 = \alpha_3$ by adding to
them general linear combinations of the $1$-forms $\sigma$,
$\overline{ \sigma}$, $\rho$, $\zeta$, $\overline{ \zeta}$:
\[
\aligned
\beta_1
&
\longmapsto
\beta_1
+
p_1\,\sigma
+
q_1\,\overline{\sigma}
+
r_1\,\rho
+
s_1\,\zeta
+
t_1\,\overline{\zeta},
\\
\beta_2
&
\longmapsto
\beta_2
+
p_2\,\sigma
+
q_2\,\overline{\sigma}
+
r_2\,\rho
+
s_2\,\zeta
+
t_2\,\overline{\zeta},
\endaligned
\]
which then modifies the torsion coefficients in some way.
Of course, from the first loop we already know
(still in the branch $R = 0$) from~\thetag{ \ref{U-5-6-7}} that:
\[
\aligned
U_5' & = 0, \ \ \ \ \ \ \ \ \ \ \ \ \ \ \ \ \ \ \ U_6'=0, \ \
\ \ \ \ \ \ \ \ \ \ \ \ \ \ \ \ \ U_7'=0,
\\
V_8' & = {\textstyle{\frac{1}{3}}}\,U_3' +
{\textstyle{\frac{1}{3}}}\,\overline{U}_4', \ \ \ \ \ \ \ \ \ \ \ \ \
\ \ W_{10}' = {\textstyle{\frac{2}{3}}}\,U_4' -
{\textstyle{\frac{1}{3}}}\,\overline{U}_3',
\endaligned
\]
and we must take account of these normalizations.
In fact, we may skip computations
and just use the previously obtained formulas~\thetag{ \ref{1-norm-eqs}},
setting simply in them:
\[
\aligned
&
p_2=q_2=r_2=s_2=t_2=0,
\ \ \ \ \
p_3=p_2,\ \ q_3=q_2,\ \ r_3=r_2,\ \ s_3=s_2,\ \ t_3=t_2,
\\
&
p_4=q_4=r_4=s_4=t_4=0,
\ \ \ \ \
p_5=q_5=r_5=s_5=t_5=0,
\endaligned
\]
which yields the following equations\,\,---\,\,including their
(unwritten) conjugates\,\,---\,\,:
\[\footnotesize\aligned
\left[
\aligned
U_1' & = 2q_1+\overline{p}_1,
\\
U_2' & = 2r_1+\overline{r}_1,
\\
U_3' & = 2s_1+\overline{t}_1,
\\
U_4' & = 2t_1+\overline{s}_1,
\\
0
&
=
0,
\\
0
&
=
0,
\\
0
&
=
0,
\endaligned\right.
\ \ \ \ \ \ \ \ \ \ \
\left[
\aligned V_1' & =
0,
\\
V_2' & = -p_1-\overline{q}_1,
\\
V_3' & = 0,
\\
V_4' & = 0,
\\
{\textstyle{\frac{1}{3}}}\,U_3'
+
{\textstyle{\frac{1}{3}}}\,\overline{U}_4'
&
=
s_1+\overline{t}_1,
\endaligned\right.
\ \ \ \ \ \ \ \ \ \ \
\left[
\aligned
W_1' & =
q_2,
\\
W_2' & = r_2,
\\
W_3' & = s_2-p_1,
\\
W_4' & = t_2,
\\
W_5' & = 0,
\\
W_6' & = -q_1,
\\
W_7' & = 0,
\\
W_8' & = -r_1,
\\
W_9' & = 0,
\\
{\textstyle{\frac{2}{3}}}\,U_4' -
{\textstyle{\frac{1}{3}}}\,\overline{U}_3' & = t_1,
\endaligned\right.
\endaligned
\]
so as to obtain null right-hand sides.
Then visually, we realize that $8$ new appropriate
linear combinations potentially provide normalizations:
\begin{equation}
\label{torsion-second-loop} \aligned 0 & = V_1',
\\
0 & = V_3',
\\
0 & = V_4',
\\
0 & = W_5',
\\
0 & = W_7',
\\
0 & = W_9',
\\
0 & = U_2'+2\,W_8'+\overline{W}_8',
\\
0 & = \overline{U}_1'+V_2'+\overline{W}_6'.
\endaligned
\end{equation}
As we already pointed out $V_4'$ is indeed normalizable.

Abbreviate the last two essential torsions as:
\[
\aligned
X_1'
&
:=
U_2'+2\,W_8'+\overline{W}_8',
\\
X_2'
&
:=
\overline{U}_1'+V_2'+\overline{W}_6'.
\endaligned
\]

\subsection{Principles for the systematic computation of
new torsion coefficients.} Now, we remember that our initial
frame:
\[
\big\{
\mathcal{L},\ \
\overline{\mathcal{L}},\ \
\mathcal{T},\ \
\mathcal{S},\ \
\overline{\mathcal{S}}
\big\}
\]
on the base manifold $M^5 \subset \C^4$ had its last three elements
given by:
\[
\aligned
\mathcal{T}
&
=
i\,\big[\mathcal{L},\,\overline{\mathcal{L}}\big],
\\
\mathcal{S}
&
=
\big[\mathcal{L},\,\mathcal{T}\big],
\\
\overline{\mathcal{S}}
&
=
\big[\overline{\mathcal{L}},\,\mathcal{T}\big],
\endaligned
\]
complemented by the following
Lie bracket structure:
\[
\aligned
\big[
\mathcal{L},\,\mathcal{S}
\big]
&
=
P\cdot\mathcal{T}
+
Q\cdot\mathcal{S}
+
R\cdot\overline{\mathcal{S}},
\\
\big[\overline{\mathcal{L}},\,\mathcal{S}\big]
&
=
A\cdot\mathcal{T}
+
B\cdot\mathcal{S}
+
\overline{B}\cdot\overline{\mathcal{S}},
\\
\big[\mathcal{L},\,\overline{\mathcal{S}}\big]
&
=
A\cdot\mathcal{T}
+
B\cdot\mathcal{S}
+
\overline{B}\cdot\overline{\mathcal{S}},
\\
\big[\overline{\mathcal{L}},\,\overline{\mathcal{S}}\big]
&
=
\overline{P}\cdot\mathcal{T}
+
\overline{R}\cdot\mathcal{S}
+
\overline{Q}\cdot\overline{\mathcal{S}},
\\
\big[\mathcal{T},\,\mathcal{S}\big]
&
=
E\cdot\mathcal{T}
+
F\cdot\mathcal{S}
+
G\cdot\overline{\mathcal{S}},
\\
\big[\mathcal{T},\,\overline{\mathcal{S}}\big]
&
=
\overline{E}\cdot\mathcal{T}
+
\overline{G}\cdot\mathcal{S}
+
\overline{F}\cdot\overline{\mathcal{S}},
\\
\big[\mathcal{S},\,\overline{\mathcal{S}}\big]
&
=
i\,J\cdot\mathcal{T}
+
K\cdot\mathcal{S}
-
\overline{K}\cdot\overline{\mathcal{S}}.
\endaligned
\]
In fact, the $5$ functions $E$, $F$, $G$, $J$, $K$ express themselves
in terms of the $5$ functions $P$, $Q$, $R$, $A$, $B$ and
their coframes derivatives.
Consequently, an important observation
is in order:

\medskip\noindent{\bf Principle.}
{\em All the subsequent computations
must necessarily be achieved only in terms of
the $5$ functions $P$, $Q$, $R$, $A$, $B$ that are really fundamental and
independent.\qed}\medskip

We must therefore replace in the torsion coefficients
(after setting $R = 0$):
\[
\footnotesize
\aligned
E
&
=
-\,i\,\overline{\mathcal L}(P)
-
i\,A\,Q
-
i\,\overline{P}\,R
+
i\,\mathcal L(A)
+
i\,B\,P
+
i\,A\,\overline{B},
\\
F
&
=
-\,i\,\overline{\mathcal L}(Q)
-
i\,R\,\overline{R}
+
i\,A
+
i\,\mathcal{L}(B)
+
i\,B\,\overline{B},
\\
G
&
=
-\,i\,P
-
i\,\overline{B}\,Q
-
i\,R\,\overline{Q}
-
i\,\overline{\mathcal L}(R)
+
i\,B\,R
+
i\,\overline{B}\overline{B}
+
i\,\mathcal L(\overline{B}).
\endaligned
\]

\[
\footnotesize
\aligned
-\,2\,J
&
=
-\,\overline{\mathcal{L}}\big(\overline{\mathcal{L}}(P)\big)
+
\overline{\mathcal{L}}\big(\mathcal{L}(A)\big)
+
\mathcal{L}\big(\overline{\mathcal{L}}(A)\big)
-
\mathcal{L}\big(\mathcal{L}(\overline{P})\big)
\\
&
\ \ \ \ \
-\,Q\,\overline{\mathcal{L}}(A)
-
2\,A\,\overline{\mathcal{L}}(Q)
-
R\,\overline{\mathcal{L}}(\overline{P})
-
2\,\overline{P}\,\overline{\mathcal{L}}(R)
-
2\,AR\overline{R}
-
2\,P\overline{P}
-
\overline{B}\overline{P}Q
-
\overline{P}\overline{Q}R
-
\\
&
\ \ \ \ \
-
\overline{R}\mathcal{L}(P)
-
2\,P\,\mathcal{L}(\overline{R})
-
\overline{Q}\,\mathcal{L}(A)
-
2\,A\mathcal{L}(\overline{Q})
-
PQ\overline{R}
-
BP\overline{Q}
+
\\
&
\ \ \ \ \
+
2\,P\overline{\mathcal{L}}(B)
+
B\overline{\mathcal{L}}(P)
+
2\,A\overline{\mathcal{L}}(\overline{B})
+
\overline{B}\,\overline{\mathcal{L}}(A)
+
2\,A\,\mathcal{L}(B)
+
2\,AA
+
2\,AB\overline{B}
+
2\,\overline{P}\,\mathcal{L}(\overline{B})
+
\\
&
\ \ \ \ \
+
B\overline{P}R
+
\overline{B}\overline{B}\overline{P}
+
B\,\mathcal{L}(A)
+
\overline{B}\,\mathcal{L}(\overline{P})
+
BBP
+
\overline{B}P\overline{R},
\endaligned
\]
\[
\footnotesize
\aligned
2i\,K
&
=
-\,\overline{\mathcal{L}}\big(\overline{\mathcal{L}}(Q)\big)
+
\overline{\mathcal{L}}\big(\mathcal{L}(B)\big)
+
\mathcal{L}\big(\overline{\mathcal{L}}(B)\big)
-
\mathcal{L}\big(\mathcal{L}(\overline{R})\big)
-
\\
&
\ \ \ \ \
-
2\,\overline{R}\,\overline{\mathcal{L}}(R)
-
R\,\overline{\mathcal{L}}(\overline{R})
-
B\,\overline{\mathcal{L}}(Q)
-
BR\overline{R}
-
2\,P\overline{R}
-
\overline{Q}R\overline{R}
-
2\,\mathcal{L}(\overline{P})
-
\overline{R}\,\mathcal{L}(Q)
-
\\
&
\ \ \ \ \
-\,2\,Q\,\mathcal{L}(\overline{R})
-
\overline{Q}\,\mathcal{L}(B)
-
2\,B\,\mathcal{L}(\overline{Q})
-
A\overline{Q}
-
\overline{P}Q
-
QQ\overline{R}
-
BQ\overline{Q}
+
\\
&
\ \ \ \ \
+
2\,\overline{\mathcal{L}}(A)
+
\overline{B}\,\overline{\mathcal{L}}(B)
+
2\,B\,\overline{\mathcal{L}}(\overline{B})
+
3\,B\,\mathcal{L}(B)
+
3\,AB
+
BBQ
+
2\,BB\overline{B}
+
2\,\overline{R}\,\mathcal{L}(\overline{B})
+
\\
&
\ \ \ \ \
+
\overline{B}\overline{B}\overline{R}
+
\overline{B}\,\mathcal{L}(\overline{R})
+
\overline{B}\overline{P}
+
Q\,\overline{\mathcal{L}}(B).
\endaligned
\]
Moreover, for any real analytic function ${\bf G}_0$
on the base manifold $M^5 \subset \C^4$, we must
expand its last three frame derivatives precisely as:
\[\footnotesize
\aligned
\mathcal{T}\big({\bf G}_0\big)
&
=
i\,\mathcal{L}\big(\overline{\mathcal{L}}\big({\bf G}_0\big)\big)
-
i\,\overline{\mathcal{L}}\big(\mathcal{L}\big({\bf G}_0\big)\big),
\\
\mathcal{S}\big({\bf G}_0\big)
&
=
i\,\mathcal{L}\big(\mathcal{L}\big(\overline{\mathcal{L}}\big(
{\bf G}_0\big)\big)\big)
-
2i\,\mathcal{L}\big(\overline{\mathcal{L}}\big(\mathcal{L}\big({\bf G}_0
\big)\big)\big)
+
i\,\overline{\mathcal{L}}\big(\mathcal{L}\big(\mathcal{L}\big(
{\bf G}_0\big)\big)\big),
\\
\overline{\mathcal{S}}\big({\bf G}_0\big)
&
=
-\,i\,\overline{\mathcal{L}}\big(\overline{\mathcal{L}}\big(\mathcal{L}\big(
{\bf G}_0\big)\big)\big)
+
2i\,\overline{\mathcal{L}}\big(\mathcal{L}\big(\overline{\mathcal{L}}\big(
{\bf G}_0\big)\big)\big)
-
i\,\mathcal{L}\big(\overline{\mathcal{L}}\big(\overline{\mathcal{L}}\big(
{\bf G}_0\big)\big)\big).
\endaligned
\]
Doing so, we obtain the following formulas that are useful to expand
in an systematic way the new torsion coefficients on a computer
machine:
\[\footnotesize
\aligned
\mathcal{S}\big(\overline{\bf D}_0\big)
&
=
-\,\mathcal{L}\big(\mathcal{L}\big(\overline{\mathcal{L}}
\big(\overline{\mathcal{L}}(B)
\big)\big)\big)
+
2\,\mathcal{L}\big(\overline{\mathcal{L}}\big(
\mathcal{L}\big(\overline{\mathcal{L}}(B)
\big)\big)\big)
-
\overline{\mathcal{L}}\big(\mathcal{L}\big(\mathcal{L}
\big(\overline{\mathcal{L}}(B)
\big)\big)\big)
+
\\
&
\ \ \ \ \
+
\mathcal{L}\big(\mathcal{L}\big(\overline{\mathcal{L}}
\big(\overline{P}\big)
\big)\big)
-2\,
\mathcal{L}\big(\overline{\mathcal{L}}\big(\mathcal{L}
\big(\overline{P}\big)
\big)\big)
+
\overline{\mathcal{L}}\big(\mathcal{L}\big(\mathcal{L}
\big(\overline{P}\big)
\big)\big)
\\
&
\ \ \ \ \
+
{\textstyle{\frac{2}{3}}}\,\overline{Q}\,
\mathcal{L}\big(\mathcal{L}\big(\overline{\mathcal{L}}
(B)
\big)\big)
-
{\textstyle{\frac{4}{3}}}\,\overline{Q}\,
\mathcal{L}\big(\overline{\mathcal{L}}\big(\mathcal{L}
(B)
\big)\big)
+
{\textstyle{\frac{2}{3}}}\,\overline{Q}\,
\overline{\mathcal{L}}\big(\mathcal{L}\big(\mathcal{L}
(B)
\big)\big)
\\
&
\ \ \ \ \
+
{\textstyle{\frac{2}{3}}}\,B\,
\mathcal{L}\big(\mathcal{L}\big(\overline{\mathcal{L}}
(\overline{Q})
\big)\big)
-
{\textstyle{\frac{4}{3}}}\,B\,
\mathcal{L}\big(\overline{\mathcal{L}}\big(\mathcal{L}
(\overline{Q})
\big)\big)
+
{\textstyle{\frac{2}{3}}}\,B\,
\overline{\mathcal{L}}\big(\mathcal{L}\big(\mathcal{L}
(\overline{Q})
\big)\big),
\endaligned
\]
\[\footnotesize
\aligned
\overline{\mathcal{S}}\big(\overline{\bf D}_0\big)
&
=
\overline{\mathcal{L}}\big(\overline{\mathcal{L}}\big(
\mathcal{L}\big(\overline{\mathcal{L}}(B)
\big)\big)\big)
-
2\,
\overline{\mathcal{L}}\big(\mathcal{L}\big(
\overline{\mathcal{L}}\big(\overline{\mathcal{L}}(B)
\big)\big)\big)
+
\mathcal{L}\big(\overline{\mathcal{L}}\big(
\overline{\mathcal{L}}\big(\overline{\mathcal{L}}(B)
\big)\big)\big)
-
\\
&
\ \ \ \ \
-\,
\overline{\mathcal{L}}\big(\overline{\mathcal{L}}\big(
\mathcal{L}\big(
\overline{P}
\big)\big)\big)
+
2\,
\overline{\mathcal{L}}\big(\mathcal{L}\big(
\overline{\mathcal{L}}\big(
\overline{P}
\big)\big)\big)
-
\mathcal{L}\big(\overline{\mathcal{L}}\big(
\overline{\mathcal{L}}\big(
\overline{P}
\big)\big)\big)
-
\\
&
\ \ \ \ \
-\,{\textstyle{\frac{2}{3}}}\,\overline{Q}\,
\overline{\mathcal{L}}\big(\overline{\mathcal{L}}\big(
\mathcal{L}
(B)\big)\big)
+
{\textstyle{\frac{4}{3}}}\,\overline{Q}\,
\overline{\mathcal{L}}\big(\mathcal{L}\big(
\overline{\mathcal{L}}
(B)\big)\big)
-
{\textstyle{\frac{2}{3}}}\,\overline{Q}\,
\mathcal{L}\big(\overline{\mathcal{L}}\big(
\overline{\mathcal{L}}
(B)\big)\big)
-
\\
&
\ \ \ \ \
-\,{\textstyle{\frac{2}{3}}}\,B\,
\overline{\mathcal{L}}\big(\overline{\mathcal{L}}\big(
\mathcal{L}
\big(\overline{Q}\big)
\big)\big)
+
{\textstyle{\frac{4}{3}}}\,B\,
\overline{\mathcal{L}}\big(\mathcal{L}\big(
\overline{\mathcal{L}}
\big(\overline{Q}\big)
\big)\big)
-
{\textstyle{\frac{2}{3}}}\,B\,
\mathcal{L}\big(\overline{\mathcal{L}}\big(
\overline{\mathcal{L}}
\big(\overline{Q}\big)
\big)\big),
\endaligned
\]
\[
\footnotesize\aligned \mathcal{T}\big(\overline{\bf D}_0\big) & = -\,
\mathcal{L}\big(\overline{\mathcal{L}}\big( \overline{\mathcal{L}}
(B) \big)\big) + \overline{\mathcal{L}}\big(\mathcal{L}\big(
\overline{\mathcal{L}} (B) \big)\big) +
\mathcal{L}\big(\overline{\mathcal{L}} \big(\overline{P}\big) \big) -
\overline{\mathcal{L}}\big(\mathcal{L} \big(\overline{P}\big) \big) +
\\
&
\ \ \ \ \
+
{\textstyle{\frac{2}{3}}}\,\overline{Q}\,
\mathcal{L}\big(\overline{\mathcal{L}}
(B)
\big)
-
{\textstyle{\frac{2}{3}}}\,\overline{Q}\,
\overline{\mathcal{L}}\big(\mathcal{L}
(B)
\big)
+
{\textstyle{\frac{2}{3}}}\,B\,
\mathcal{L}\big(\overline{\mathcal{L}}
\big(\overline{Q}\big)
\big)
-
{\textstyle{\frac{2}{3}}}\,B\,
\overline{\mathcal{L}}\big(\mathcal{L}
\big(\overline{Q}\big)
\big),
\endaligned
\]
with quite similar formulas for ${\bf B}_0$ and ${\bf E}_0$.

In subsection \ref{length-six}, we obtained five relations
between the fundamental functions $A,B,P,Q,R$, extracted from the
iterated Lie brackets of the length six. By the assumption $R\equiv 0$,
they become:
\begin{equation}\label{length6-R=0}
\!\!\!\!\!\!\!\!\!\!\!\!\!\!\!\!\!\!\!\!
\footnotesize \aligned 0&\overset{1}{\equiv} 2\,\mathcal
L(\overline{\mathcal L}(P))-\mathcal L(\mathcal
L(A))-\overline{\mathcal L}(\mathcal L(P))- 2\,P\mathcal
L(B)-B\mathcal L(P)-2\,A\mathcal L(\overline B)-\overline B\mathcal
L(A)+P\overline{\mathcal L}(Q)+
\\
&+A\mathcal L(Q)+2\,Q\mathcal L(A)-Q\overline{\mathcal
L}(P)-PB\overline B-A\overline
B^2+PBQ+2\,AQ\overline B-AQ^2,
\\
0&\overset{2}{\equiv}2\,\mathcal L(\overline{\mathcal L}(Q))-\mathcal
L(\mathcal L(B))-\overline{\mathcal L}(\mathcal L(Q))-
\\
&-2\,\mathcal L(A)-2\,B\mathcal L(\overline B)-\overline B\mathcal
L(B)+\overline{\mathcal L}(P)+BQ\overline B-A\overline B-B\overline
B^2+AQ,
\\
0&\overset{3}{\equiv}-\mathcal
L(\mathcal L(\overline B))-
\\
&-3\,\overline B\mathcal L(\overline B)+\overline B\mathcal
L(Q)+2\,Q\mathcal L(\overline B)+\mathcal L(P)+2\,Q\overline
B^2-QP-Q^2\overline B-\overline B^3+P\overline B,
\\
 0&\overset{4}{\equiv}-3\overline{\mathcal L}(\mathcal L(A))-\mathcal L(\mathcal
L(\overline P))+3\,\mathcal
L(\overline{\mathcal L}(A))+\overline{\mathcal L}(\overline{\mathcal
L}(P))-
\\
&-2\,A\mathcal L(\overline Q)-\overline Q\mathcal L(A)+3\,B\mathcal
L(A)+3\,\overline B\mathcal L(\overline P)-3\,B\overline{\mathcal
L}(P)-3\,\overline B\,\overline{\mathcal L}(A)+2\,A\overline{\mathcal
L}(Q)+Q\overline{\mathcal L}(A)-
\\
&-BP\overline Q+3\,B^2P+2\,A\overline B\,\overline
Q-2\,BQA-3\,\overline B^2\overline P+Q\overline B\,\overline P,
\\
0&\overset{5}{\equiv}-3\,\overline{\mathcal L}(\mathcal
L(B))+3\,\mathcal L(\overline{\mathcal L}(B))+\overline{\mathcal
L}(\overline{\mathcal L}(Q))+
\\
&+3\,B\mathcal L(B)-3\,\overline B\,\overline{\mathcal
L}(B)+Q\overline{\mathcal L}(B)-B\overline{\mathcal
L}(Q)-2\,B\mathcal L(\overline Q)-\overline Q\mathcal
L(B)-2\,\mathcal L(\overline P)-
\\
&-Q\overline P-A\overline Q-BQ\overline Q+3\,AB+3\,\overline
B\,\overline P+2\,B\overline B\,\overline Q+B^2Q.
\endaligned
\end{equation}

These five equations enable us to simplify the results
obtained during the computations. In particular, we have:

\begin{Lemma}
After determining the group parameter $\sf e$, the two (normalizable)
torsion coefficients $W'_7$ and $X_2'$ vanish identically.
\end{Lemma}

\proof
After determining $\sf e$, the expressions of $W'_7$ and
$X_2'$ take the forms:
\[\footnotesize\aligned
W'_7&=-\frac{i}{\overline{\sf a}^3}\Big(\overline{\mathcal
L}(\overline{\mathcal L}(B))-\overline{\mathcal L}(\overline
P)-B\overline{\mathcal L}(\overline Q)+3\,B\overline{\mathcal
L}(B)-2\,\overline{\mathcal L}(B)\overline Q+
\\
&\ \ \ \ \ \ \ \ \ \ +\overline Q\,\overline P+B\overline
Q^2-2\,B^2\overline Q+B^3-B\overline P\Big),
\\
X_2'&=-\frac{i}{6{\sf a}\overline{\sf a}^2}\Big(\overline{\mathcal
{L}}(\overline{\mathcal{L}}(Q))-3\,\overline{\mathcal L}({\mathcal
L}(B))+3\,{\mathcal L}(\overline{\mathcal L}(B))-2\,{\mathcal
L}(\overline P)-B\overline{\mathcal L}(Q)+
\\
&+3\,AB+3\,B{\mathcal L}(B) -\overline Q{\mathcal L}(B)-A\overline
Q-Q\overline P-2\,B{\mathcal L}(\overline Q)-3\,\overline B\overline
{\mathcal L}(B)+B^2Q+
\\
&+Q\overline{\mathcal L}(B)+3\,\overline B\,\overline P-BQ\overline
Q+2\overline BB\overline Q\Big).
\endaligned
\]
Now, it suffices to use the already presented equation
$\overset{5}{=}$ and the conjugation of $\overset{3}{=}$ to extract
the expressions of $\overline{\mathcal L}(\mathcal L(B))$ and of
$\mathcal L(\mathcal L(\overline B))$ in $W'_7$, respectively and
subsequently to insert them in $X_2'$ and in $W'_7$.
\endproof

\section{Four group parameter general normalizations}
\label{four-normalizations}
\HEAD{\ref{branch-R-neq-0}.~Four group parameter general normalizations}{
Masoud {\sc Sabzevari} (Shahrekord) and Jo\"el {\sc Merker} (LM-Orsay)}

\subsection{Setting up the torsion coefficients}
After normalizing the group parameter $\sf e$ and inserting it into
the expressions of the lifted coframe, one obtains:
\[
\aligned \sigma&={\sf a}^2\overline{\sf a}\cdot\,\sigma_0,
\\
\rho&={\sf a}\overline{\sf a}\cdot\,(\underbrace{\mathbf{
C}_0\,\sigma_0+{\mathbf{\overline C}}_0\overline\sigma_0
+\rho_0}_{\rho_0^{\sim}}),
\\
\zeta&={\sf
a}\cdot\big[\underbrace{\mathbf{E}_0\,\sigma_0+\underbrace{\mathbf{\overline
D}_0\,\overline\sigma_0+\mathbf{
B}_0\cdot\rho_0+\zeta_0}_{\zeta_0^\sim}}_{=:\zeta_0^{\approx}}\big].
\endaligned
 \]
In other words, the lifted coframe converts into the form:
\begin{eqnarray*}
\left[
\begin{array}{l}
\sigma={\sf a}^2\overline{\sf a}\cdot\sigma_0,
\\
\rho={\sf a}\overline{\sf a}.\rho_0^\sim,
\\
\zeta={\sf a}\cdot\zeta_0^\approx,
\end{array}
\right.
\end{eqnarray*}
with the reduced matrix group:
\begin{equation*}\footnotesize\aligned
G^\approx:=\left\{g^\approx:=
\left(%
\begin{array}{ccccc}
  {\sf a}\overline{\sf a}^2 & 0 & 0 & 0 & 0 \\
  0 & {\sf a}^2\overline{\sf a} & 0 & 0 & 0 \\
  0 & 0 & {\sf a}\overline{\sf a} & 0 & 0 \\
  0 & 0 & 0 & \overline{\sf a} & 0 \\
  0 & 0 & 0 & 0 & {\sf a} \\
\end{array}%
\right) : {\sf a}\in\mathbb C\right\},
\endaligned
\end{equation*}
and with the modified Maurer-Cartan matrix:
\begin{equation*}
\footnotesize\aligned
 dg^\approx\cdot {g^\approx}^{-1}=
\left(%
\begin{array}{ccccc}
  2\,\overline\beta+\beta & 0 & 0 & 0 & 0 \\
  0 & 2\,\beta+\overline\beta & 0 & 0 & 0 \\
  0 & 0 & \beta+\overline\beta & 0 & 0 \\
  0 & 0 & 0 & \overline\beta & 0 \\
  0 & 0 & 0 & 0 & \beta \\
\end{array}%
\right).
\endaligned
\end{equation*}
Here, there remains just one Maurer-Cartan 1-form:
\[
\beta:=\beta_1=\frac{d\sf a}{\sf a}.
\]

Now, applying
the exterior differentiation operator on the already obtained lifted
coframe gives:
\[
\aligned d\sigma&=(2\beta+\overline\beta)\wedge\sigma+
\\ &
\ \ \ \ +{\sf a}^2\overline{\sf a}\cdot d\sigma_0,
\\
d\rho&=(\beta+\overline\beta)\wedge\rho+
\\
& \ \ \ +{\sf a}\overline{\sf a}\cdot(\mathbf{
C}_0\,d\sigma_0+{\mathbf{\overline C}}_0\,d\overline\sigma_0
+d\rho_0+
\\
& \ \ \ \ \ \ \ \ +d\mathbf{ C}_0\wedge\sigma_0+d{\mathbf{\overline
C}}_0\wedge\overline\sigma_0),
\\
d\zeta&=\beta\wedge\zeta+
\\
&\ \ \ \ \ +{\sf a}\cdot(\mathbf{E}_0\,d\sigma_0+\mathbf{\overline
D}_0\,d\overline\sigma_0+\mathbf{ B}_0d\rho_0+d\zeta_0+
\\
& \ \ \ \ \ \ \ \ \ \ \
+d\mathbf{E}_0\wedge\sigma_0+d\mathbf{\overline
D}_0\wedge\overline\sigma_0+d\mathbf{ B}_0\wedge\rho_0).
\endaligned
\]
Comparing with \thetag{\ref{set-up-1}}, one easily verifies that the
expressions of $d\sigma$ and $d\rho$ are unchanged. Hence
except possible replacements of ${\sf e} = {\sf a} \,
{\bf E}_0$, we have no
essential change in the expressions of the
new torsion coefficients $U_i''$ and
$V_j''$. Nevertheless, some of the torsion coefficients $W_k''$
change after determining $\sf e$. More precisely, our computations
show that just the four torsion coefficients $W_1''$, $W_2''$, $W_3''$,
$W_4''$ convert into the following modified forms:
\begin{equation}
\label{new-new} \aligned {W''_1}&:= W_1'-\frac{1}{{\sf
a}^2\overline{\sf a}^3}\,\overline{\mathcal
S}(\mathbf{E}_0)+\frac{\overline{\sf c}}{{\sf a}^3\overline{\sf
a}^4}\,{\mathcal T}(\mathbf{E}_0)-\frac{{\sf b}\overline{\sf c}}{{\sf
a}^4\overline{\sf a}^4}\,\mathcal L(\mathbf{E}_0)+\frac{\overline{\sf
d}}{{\sf a}^3\overline{\sf a}^3}\,{\mathcal L}(\mathbf{E}_0)-
\\
& \ \ \ \ \ \ \ \ \ \ \ \ -\frac{\overline{\sf b}\overline{\sf
c}}{{\sf a}^3\overline{\sf a}^5}\,\overline{\mathcal
L}(\mathbf{E}_0)+\frac{\overline{\sf e}}{{\sf a}^2\overline{\sf
a}^4}\,\overline{\mathcal L}(\mathbf{E}_0),
\\
{W_2''}&:= W_2'-\frac{1}{{\sf a}^2\overline{\sf a}^2}\,\mathcal
T(\mathbf{E}_0)+\frac{\sf b}{{\sf a}^3\overline{\sf a}^2}\,\mathcal
L(\mathbf{E}_0)+\frac{\overline{\sf b}}{{\sf a}^2\overline{\sf
a}^3}\,\overline{\mathcal L}(\mathbf{E}_0),
\\
{W_3''}&:= W_3'-\frac{1}{{\sf a}^2\overline{\sf a}}\,\mathcal
L(\mathbf{E}_0),
\\
{W_4''}&:= W_4'-\frac{1}{{\sf a}\overline{\sf
a}^2}\,\overline{\mathcal L}(\mathbf{E}_0).
\endaligned
\end{equation}

\subsection{The remaining normalizable expressions in the second loop}

At the second loop absorbtion, we found eight normalizable expressions
\thetag{\ref{torsion-second-loop}}. Among them, $V'_3$ vanished
automatically and also $V'_4 \equiv 0$ after determining $\sf
e$. Subsequently, $W'_7$ and $X_2'$ vanished identically, as soon as
we took account of the above relations $\overset{1}{=}$,
$\overset{2}{=}$, $\overset{3}{=}$, $\overset{4}{=}$, $\overset{5}{=}$
coming from a study of Jacobi identities between iterated Lie brackets
of length six.  Then, for the moment, a computer-assisted calculation
provides the following expressions for all the essential
torsions~\thetag{\ref{torsion-second-loop}}:
\begin{equation}
\!\!\!\!\!\!\!\!\!\!\!\!\!\!\!\!\!\!\!\!
\label{expression-2nd-loop}
\footnotesize\aligned
V''_1&:=-\frac{i}{3\,{\sf a}^2\overline{\sf a}^2}\bigg(3\,\mathcal
L(\mathcal L(\overline{\mathcal L}(B)))+3\,\mathcal L(\mathcal
L(\overline{\mathcal L}(\overline Q)))-6\,\mathcal
L(\overline{\mathcal L}(\mathcal L(\overline
Q)))+2\,\overline{\mathcal L}(\mathcal L(\mathcal L(\overline
Q)))+
\\
&
+3\,\overline{\mathcal L}(\overline{\mathcal L}(\mathcal
L(\overline B)))+3\,\overline{\mathcal L}(\overline{\mathcal L}(\mathcal
L(Q)))-6\,\overline{\mathcal L}(\mathcal L(\overline{\mathcal
L}(Q)))+2\,\mathcal L(\overline{\mathcal L}(\overline{\mathcal
L}(Q))) -
\\
&-4\,\overline{\mathcal L}(\overline{\mathcal L}(P))+9\,\mathcal
L(\overline{\mathcal L}(A))-7\,\mathcal L(\mathcal L(\overline
P))+2\,\overline Q\overline{\mathcal L}(\mathcal L(\overline
B))-6\,\overline B\mathcal L(\overline{\mathcal L}(\overline
Q))+10\,\overline B\overline{\mathcal L}(\mathcal L(\overline Q))+
\\
&+2\,Q\mathcal L(\overline{\mathcal L}(B))+10\,B\mathcal
L(\overline{\mathcal L}(Q))- 6\,B\overline{\mathcal L}(\mathcal
L(Q))-4\,\overline Q\mathcal L(\overline{\mathcal L}(Q))+2\,\overline
Q\overline{\mathcal L}(\mathcal L(Q))-3\,B\mathcal L(\mathcal
L(\overline Q))-
\\
&-3\,\overline B\overline{\mathcal L}(\overline{\mathcal L}(Q))-
4\,Q\overline{\mathcal L}(\mathcal L(\overline Q))+2\,Q\mathcal
L(\overline{\mathcal L}(\overline Q)) -2\,\overline P \mathcal
L(Q)+6\,\overline P \mathcal L(\overline B)-2\,P\overline{\mathcal
L}(\overline Q)+6\,P\overline{\mathcal L}(B)+
\\
&+3\,A\mathcal L(B)-5\,A\mathcal L(\overline Q)+
3\,A\overline{\mathcal L}(\overline B)+A\overline{\mathcal
L}(Q)+2\,Q\overline{\mathcal L}(A)-\overline Q\mathcal
L(A)+7\,B\overline{\mathcal L}(P)-9\,\overline B\overline{\mathcal
L}(A)+
\\
&+16\,\overline B\mathcal L(\overline P)- 2\,\overline
BQ\overline{\mathcal L}(\overline Q)- 2\,\mathcal
L(B)\overline{\mathcal L}(Q)+2\,B^2\mathcal L(Q)+2\,\mathcal
L(Q)\overline{\mathcal L}(B)+6\,\mathcal L(B)^2+
9\,\overline{\mathcal L}(\overline B)\mathcal L(B)-
\\
&-6\,B^2\overline B^2-9\,B\overline B\overline{\mathcal L}(\overline
B)+5\,BQ\overline{\mathcal L}(\overline B)-9\,B\overline B\mathcal
L(B)+5\,B^2\overline BQ+ 5\,\overline B\,\overline Q\mathcal L(B)
+6\,\overline BQ\overline{\mathcal L}(B)+
\\
&+4\,B\overline B\overline{\mathcal L}(Q)+4\,B\overline B\mathcal
L(\overline Q)+ 6\,B\overline Q\mathcal L(\overline B)-9\,\overline
B^2\overline{\mathcal L}(B)-9\,B^2\mathcal L(\overline B)-3\,\mathcal
L(\overline B)\overline{\mathcal L}(B)-6\,\overline{\mathcal
L}(\overline B)\overline{\mathcal L}(Q)-
\\
&-2\,\overline B\,\overline Q\overline{\mathcal
L}(Q)+6\,\overline{\mathcal L}(\overline B)^2+4\,\overline
B\,\overline Q\overline{\mathcal L}(\overline B) -2\,Q\overline
Q\overline{\mathcal L}(\overline B)-4\,\overline Q\overline{\mathcal
L}(P)-2\,\overline{\mathcal L}(\overline B)\mathcal L(\overline
Q)+2\,\mathcal L(\overline B)\overline{\mathcal L}(\overline Q)-
\\
&-4\,Q\mathcal L(\overline P)-6\,\mathcal L(B)\mathcal L(\overline
Q)+ 2\,\overline B^2\overline{\mathcal L}(\overline Q)-2\,B\mathcal
L(Q)\overline Q-2\,Q\overline Q\mathcal L(B)+4\,BQ\mathcal L(B)-
\\
&-2\,BQ\mathcal L(\overline Q) -6\,\overline B^2\,\overline
P+10\,\overline B\,\overline QA-BP\overline Q+2\,\overline
B\,\overline PQ+4\,ABQ-9\,AB\overline B -4\,QA\overline Q-
\\
&-3\,BQ\overline B\,\overline Q+5\,\overline B^2\overline
QB+3\,B^2P\bigg)
\\
V''_3&=0,
\\
V''_4&=0,
\\
W''_5&=-\frac{1}{3\,{\sf a}\overline{\sf
a}^3}\bigg(2\,\overline{\mathcal L}(\mathcal L(\overline{\mathcal
L}(\overline Q)))-3\,\overline{\mathcal L}(\mathcal
L(\overline{\mathcal L}(B)))+3\,\overline{\mathcal
L}(\overline{\mathcal L}(\mathcal L(B)))-\overline{\mathcal
L}(\overline{\mathcal L}(\mathcal L(\overline Q)))-\mathcal
L(\overline{\mathcal L}(\overline{\mathcal L}(\overline Q)))+
\\
&+3\mathcal L(\overline{\mathcal L}(\overline
P))-3\,\overline{\mathcal L}(\mathcal L(\overline P))+ 3\,\overline
Q\mathcal L(\overline{\mathcal L}(B))-2\,\overline
Q\overline{\mathcal L}(\mathcal L(B))+2\,B\mathcal
L(\overline{\mathcal L}(\overline Q))-2\,B\overline{\mathcal
L}(\mathcal L(\overline Q))-
\\
&-3\,B\mathcal L(\overline{\mathcal L}(B))
-Q\overline{\mathcal L}(\overline{\mathcal L}(B))+3\,\overline
B\overline{\mathcal L}(\overline{\mathcal L}(B))+
+QB\overline{\mathcal L}(\overline Q) -4\,B\overline
B\overline{\mathcal L}(\overline Q)+3\,\overline
P\,\overline{\mathcal L}(\overline B)
\\
&
-3\,\overline{\mathcal
L}(\overline B)\overline{\mathcal L}(B)+3\,A\overline{\mathcal
L}(B)+2\,B\overline Q\overline{\mathcal L}(\overline B)
+15\,B\overline B\overline{\mathcal L}(B)-7\,\overline B\,\overline
Q\,\overline{\mathcal L}(B)-\overline Q^2\mathcal L(B)-
\\
&
-3\,B^2\mathcal
L(B)-2\,B\overline Q\mathcal L(\overline Q)+2\,Q\overline
Q\overline{\mathcal L}(B)-3\,BQ\overline{\mathcal L}(B)+
2\,\overline{\mathcal L}(B)\mathcal L(\overline Q)- 2\,\overline
P\mathcal L(\overline Q)-
\\
&
-3\,\mathcal L(B)\overline{\mathcal
L}(B)+3\,B\overline{\mathcal L}(A)-\overline Q\overline{\mathcal
L}(A)-A\overline{\mathcal L}(\overline Q)+2\,B^2\mathcal L(\overline
Q)-3\,\overline B\,\overline{\mathcal L}(\overline P)+
\\
&+4\,B\overline Q\mathcal L(B)+\mathcal L(B)\overline{\mathcal
L}(\overline Q)+Q\overline{\mathcal L}(\overline P) +3\,AB^2 +
6\,B^3\overline B -4\,AB\overline Q-3\,B\overline B\overline
P-10\,B^2\overline B\overline Q+
\\
&+3\,\overline B\,\overline Q\,\overline P+4\,B\overline B\,\overline
Q^2-B^3Q+A\overline Q^2 + BQ\overline P-BQ\overline Q^2-Q\overline
Q\overline P+2\,B^2Q\overline Q
\bigg),
\\
W''_7&=0,
\endaligned
\end{equation}
\[
\!\!\!\!\!\!\!\!\!\!\!\!\!\!\!\!\!\!\!\!
\footnotesize
\aligned
W''_9&=\frac{i}{9\overline{\sf a}^2}\bigg(18\,\overline{\mathcal
L}(B)-3\,\overline{\mathcal L}(\overline Q)-9\,\overline
P-12\,B\overline Q+9\,B^2+\overline Q^2\bigg),
\\
X_1''&=-\frac{i}{9{\sf a}\overline{\sf a}}\bigg(6\,\overline{\mathcal
L}(Q)+6\,\mathcal L(\overline Q)-18\,\mathcal
L(B)-18\,\overline{\mathcal L}(\overline B)+27\,B\overline
B-6\,BQ-6\,\overline B\,\overline Q+2\,Q\overline Q+9\,A\bigg),
\\
X_2''&=0.
\endaligned
\]

\begin{Theorem}
\label{th-loop-2}
When at least one among the above four
independent essential torsions $V''_1$, $W''_5$,
$W''_9$, $X_1''$ does not vanish identically, one can normalize
either ${\sf a}$ or ${\sf a} \overline{ \sf a}$
after relocalization to a neighborhood of a generic point,
and more precisely:

\begin{itemize}

\smallskip\item[{\bf (i)}]
when $X_1'' \not\equiv 0$, setting $X_1'' := - \frac{ i}{ 9}$
(noticing that $X_1''$ is a
purely imaginary valued function), one
normalizes ${\sf a}\overline{\sf a}$, and in this case, the three
remaining expressions $W''_9$, $W''_5$, $V''_1$ are the invariants of
the equivalence problem;

\smallskip\item[{\bf (ii)}]
when $W''_9\not\equiv 0$ but $X_1'' \equiv 0$, setting $W_9'' :=
\frac{ i}{ 9}$, one normalizes ${\sf a}$, and in this case the two
remaining expressions $W''_5$ and $V''_1$ are the invariants of the
equivalence problem;

\smallskip\item[{\bf (iii)}]
when $W''_5\neq 0$ but $X_1'' \equiv W_9'' \equiv 0$, setting
$W_5'' := \frac{1}{3}$, one normalizes ${\sf a}$, and
in this case,
$V''_1$ is the single invariant of the problem.

\smallskip\item[{\bf (iv)}]
when, if $V''_1\neq 0$ but  $X_1'' \equiv W_9'' \equiv W_5'' \equiv 0$,
setting $V_1'' := -\, \frac{ i}{ 3}$
(noticing that $V_1''$ is a
purely imaginary valued function),
one normalizes ${\sf a} \overline{\sf a}$.

\smallskip\item[{\bf (v)}]
when $X_1'' \equiv W_9'' \equiv W_5'' \equiv V_1'' \equiv 0$,
one has to start
the third loop of the Cartan equivalence procedure under the assumption:
\begin{equation}
\label{assum-second-loop}
V''_1
=
W''_5
=
W''_9
=
X''_1
\equiv 0.
\end{equation}
\end{itemize}
\end{Theorem}

\proof
The reason of choosing first the expressions with the lowest degree of
derivations\,\,---\,\,namely $X_1''$, $W_9''$ of $\big\{ \mathcal{L},
\overline{\mathcal{L}} \big\}$-derivation order $1$, and afterwards
$W_5''$, $V_1''$ of $\big\{ \mathcal{L}, \overline{\mathcal{L}}
\big\}$-derivation order $3$\,\,---\,\,is that the vanishing of the
less complex ones may possibly cause the vanishing of the ones having
greater complexity.

To check the
independency of these expressions, we use the {\sc Maple} command
{\sc IdealMembership}. The procedure will be described in the next
subsection \ref{third-loop}.
\endproof

\subsection{Third loop normalizable essential torsion combinations}
\label{third-loop}

Now, the assumptions \thetag{\ref{assum-second-loop}} lead us to
start the third-loop normalization. At this time, our structure
equations have the following form:
\begin{equation}
\label{struc-equ-third-loop}
\footnotesize \aligned d\sigma & =
\big(2\,\beta+\overline{\beta}\big) \wedge \sigma+
\\
& \ \ \ \ \ + U_1'''\,\sigma\wedge\overline{\sigma} +
U_2'''\,\sigma\wedge\rho + U_3'''\,\sigma\wedge\zeta +
U_4'''\,\sigma\wedge\overline{\zeta} +
\\
& \ \ \ \ \ \ \ \ \ \ \ \ \ \ \ \ \ \ \ \ \ \ \ \ +
\zero{U_5'''}\,\overline{\sigma}\wedge\rho +
\zero{U_6'''}\,\overline{\sigma}\wedge\zeta +
\zero{U_7'''}\,\overline{\sigma}\wedge\overline{\zeta} +
\\
& \ \ \ \ \ \ \ \ \ \ \ \ \ \ \ \ \ \ \ \ \ \ \ \ \ \ \ \ \ \ \ \ \ \
\ \ \ \ \ \ \ \ \ \ \ \ \ \ \ \ \ \ \ \ \ + {\rho\wedge\zeta},
\\
 d\rho & = \big(\beta+\overline{\beta}\big)\wedge\rho +
\\
& \ \ \ \ \ + \zero{V_1'''}\,\sigma\wedge\overline{\sigma} +
V_2'''\,\sigma\wedge\rho + \zero{V_3'''}\,\sigma\wedge\zeta +
\zero{V_4'''}\,\sigma\wedge\overline{\zeta} +
\\
& \ \ \ \ \ \ \ \ \ \ \ \ \ \ \ \ \ \ \ \ \ \ \ \ +
\overline{V_2}'''\,\overline{\sigma}\wedge\rho +
\zero{\overline{V_4}'''}\,\overline{\sigma}\wedge\zeta +
\zero{\overline{V_3}'''}\,\overline{\sigma}\wedge\overline{\zeta} +
\\
& \ \ \ \ \ \ \ \ \ \ \ \ \ \ \ \ \ \ \ \ \ \ \ \ \ \ \ \ \ \ \ \ \ \
\ \ \ \ \ \ \ \ \ + V_8'''\,\rho\wedge\zeta +
\overline{V_8}'''\,\rho\wedge\overline{\zeta} +
\\
& \ \ \ \ \ \ \ \ \ \ \ \ \ \ \ \ \ \ \ \ \ \ \ \ \ \ \ \ \ \ \ \ \ \
\ \ \ \ \ \ \ \ \ \ \ \ \ \ \ \ \ \ \ \ \ \ \ \ \ \ \ \ \ \ \ \ +
i\,\zeta\wedge\overline{\zeta},
\\
d\zeta & = \beta\wedge\zeta +
\\
& \ \ \ \ \ + W_1'''\,\sigma\wedge\overline{\sigma} +
W_2'''\,\sigma\wedge\rho + W_3'''\,\sigma\wedge\zeta +
W_4'''\,\sigma\wedge\overline{\zeta} +
\\
& \ \ \ \ \ \ \ \ \ \ \ \ \ \ \ \ \ \ \ \ \ \ \ \ +
\zero{W_5}'''\,\overline{\sigma}\wedge\rho +
W_6'''\,\overline{\sigma}\wedge\zeta+
\zero{W_7'''}\,\overline{\sigma}\wedge\overline{\zeta} +
\\
& \ \ \ \ \ \ \ \ \ \ \ \ \ \ \ \ \ \ \ \ \ \ \ \ \ \ \ \ \ \ \ \ \ \
\ \ \ \ \ \ \ \ \ \ \ \ \ + W_8'''\,\rho\wedge\zeta +
\zero{W_9'''}\,\rho\wedge\overline{\zeta} +
\\
& \ \ \ \ \ \ \ \ \ \ \ \ \ \ \ \ \ \ \ \ \ \ \ \ \ \ \ \ \ \ \ \ \ \
\ \ \ \ \ \ \ \ \ \ \ \ \ \ \ \ \ \ \ \ \ \ \ \ \ \ \ \ \ \ \ \ \ \ \
\ + W_{10}'''\,\zeta\wedge\overline{\zeta},
\endaligned
\end{equation}
where we underline the torsions that vanish identically.

Starting this step of normalization, one has to replace the single remaining
Maurer-Cartan 1-form $\beta$ as:
\[
\beta\mapsto\beta+p\,\sigma+q\,\overline\sigma+r\,\rho+s\,\zeta+t\,\overline\zeta
\]
and to proceed with
the same line of computations as in the former steps. Then,
by anticipation, one
obtains $4$ new potentially normalizable expressions:
\[
\left[ \aligned 0 & = W_1''',
\\
0 & = W_2''',
\\
0 & = W_4''',
\\
0 & = V_2''' - W_3''' - \overline{W}_6'''=:Y'''.
\endaligned\right.
\]

Putting the lastly obtained expressions of $V'''_2$, $W'''_3$, $W'''_6$
in $Y'''$ immediately implies that:

\begin{Lemma}
The last normalizable expression $Y'''$ vanishes identically.
\qed
\end{Lemma}

Similarly to the second step of normalization, if at least one of the
above expressions does not vanish, then it can be employed to
determine the last group parameter $\sf a$ and next, the remaining
nonzero normalizable expressions will be the invariants of the
equivalence problem. Otherwise, one has to start the prolongation
procedure.

Before proceeding, let us check whether one of the above
normalizable expressions can be expressed as a combination of the
remaining ones or those in \thetag{\ref{assum-second-loop}}. For this
aim, first we extract the zero-order terms of each expression. Since
these terms do not admit any derivations of $\mathcal L$ or
$\overline{\mathcal L}$, then it is not required to consider also the
derivations of the mentioned expressions. Now, we check the
independency of these extracted zero-order expressions. For this aim,
we use the {\sc Maple} command {\sc IdealMembership}\footnote{It is
in fact a time and memory consuming command of {\sc Maple} and that
it is why that we use first just the zero-order terms of the
expressions instead of their corresponding lengthy ones.} and realize
that the first order terms of $W'''_4$ can be eliminated by the
first-order terms of $W'''_9$ and $X_1'''$. Then, we {\it surmise} that
$W'''_4$ may be expressed as a combination of\,\,---\,\,assumed to be
vanishing\,\,---\,\,$W'''_9$ and $X_1'''$. The maximum $\big\{
\mathcal{L}, \overline{ \mathcal{L}} \big\}$-derivation order in
$W'''_4$ equals $2$, while it equals $1$ in the expressions of
$W'''_9$ and $X_1'''$.  Hence, we guess that $W'''_4$ can be expressed
as a combination of $W'''_9$, $X_1'''$ and
$\mathcal{L}(X_1'''),\overline{\mathcal
L}(X_1'''),\mathcal{L}(W'''_9),\ldots,\overline{\mathcal
L}(\overline W'''_9)$. After somehow tremendous computations, we
realized that:
\[
W'''_4
=
-\frac{6\,\overline B-2\,Q}{5\,\sf a}\,W'''_9+\frac{\overline
Q}{5\overline{\sf a}} X_1'''+\frac{3}{5\,\sf a}\mathcal
L(W'''_9)-\frac{3}{5\overline{\sf a}}\overline{\mathcal L}(X_1''').
\]

\begin{Lemma}
Under the assumptions
\thetag{\ref{assum-second-loop}},
the normalizable expression $W'''_4$ vanishes identically.\qed
\end{Lemma}

Our inspections show that it is not possible to express one of the
two remaining expressions $W'''_1,W'''_2$ in terms of the other one
or in terms of the expressions $V''_1$, $W'''_5$, $W'''_9$, $X_1'''$. Hence,
in this step we encounter two essential torsion coefficients having the
following expressions\footnote{The full expression of $W'''_1$ is much too
long\,\,---\,\,it involves
$\big\{ \mathcal{L},
\overline{ \mathcal{L}} \big\}$-derivations
of order $4$\,\,---, hence we do
not typeset its expanded expression. Nevertheless, it is available in
the {\sc Maple} worksheet \cite{Merker-Sabzevari-Maple}}:
\begin{equation}
\!\!\!\!\!\!\!\!\!\!\!\!\!\!\!\!\!\!\!\!
\label{expressions-loop-three}
\footnotesize \aligned
W'''_1
&
=
W_1''
\Big\vert_{\text{\sf replace}\,\,{\sf e}={\sf a}{\bf E}_0}
-
\frac{1}{{\sf a}^2\overline{\sf
a}^3}\,\overline{\mathcal S}(\mathbf{E}_0)+\frac{\overline{\sf
c}}{{\sf a}^3\overline{\sf a}^4}\,{\mathcal
T}(\mathbf{E}_0)-\frac{{\sf b}\overline{\sf c}}{{\sf
a}^4\overline{\sf a}^4}\,\mathcal \overline{\mathcal
L}(\mathbf{E}_0)+\frac{\overline{\sf d}}{{\sf a}^3\overline{\sf
a}^3}\,{\mathcal \overline{\mathcal L}(}(\mathbf{E}_0)-
\\
& \ \ \ \ \ \ \ \ \ \ \ \ -\frac{\overline{\sf b}\overline{\sf
c}}{{\sf a}^3\overline{\sf a}^5}\,\overline{\mathcal
\overline{\mathcal L}}(\mathbf{E}_0)+\frac{\overline{\sf e}}{{\sf
a}^2\overline{\sf a}^4}\,\overline{\mathcal \overline{\mathcal
L}}(\mathbf{E}_0)
\\
&
=
\frac{1}{9\,{\sf a}^2\overline{\sf a}^3}\,
\Big(
\text{\footnotesize\sf long expression}
\Big),
\\
W'''_2&=-\frac{1}{405{\sf a}^2\overline{\sf
a}^2}\bigg[1485\,\overline{\mathcal L}(\overline{\mathcal L}(\mathcal
L(Q)))-2610\,\overline{\mathcal L}(\mathcal L(\overline{\mathcal
L}(Q)))+ 270\,\mathcal L(\overline{\mathcal L}(\overline{\mathcal
L}(Q)))-1080\,\mathcal L(\overline{\mathcal L}(\mathcal L(\overline
Q)))+
\\
&+270\,\overline{\mathcal L}(\mathcal L(\mathcal L(\overline
Q)))+540\,\mathcal L(\mathcal L(\overline{\mathcal L}(\overline
Q)))+405\,\mathcal L(\mathcal L(\overline{\mathcal L}(B)))-
675\,\overline{\mathcal L}(\overline{\mathcal L}(\mathcal L(\overline
B)))+4320\,\mathcal L(\overline{\mathcal L}(A))+
\\
&+111\,P\overline Q^2+ 3150\,B\mathcal L(\overline{\mathcal
L}(Q))+156\,Q^2\overline P-198Q\mathcal L(\overline{\mathcal
L}(B))-1980\,\overline B\overline{\mathcal L}(\mathcal L(\overline
Q))+999\,\overline B\mathcal L(\overline{\mathcal L}(\overline Q))-
\\
&- 918\,\overline Q\overline{\mathcal L}(\mathcal L(\overline
B))+6345\,B\overline{\mathcal L}(\mathcal L(\overline
B))+2241\,\overline B\mathcal L(\overline{\mathcal L}(B))-
1980\,\mathcal L(\mathcal L(\overline P))-2250\,B\overline{\mathcal
L}(\mathcal L(Q))-
\endaligned
\end{equation}
\[
\!\!\!\!\!\!\!\!\!\!\!\!\!\!\!\!\!\!\!\!
\scriptsize
\aligned
&-1320\,\overline Q\mathcal L(\overline{\mathcal L}(Q))- 270B\mathcal
L(\mathcal L(\overline Q))+1215\,\overline B\overline{\mathcal
L}(\overline{\mathcal L}(Q))+ 18\,Q\mathcal L(\overline{\mathcal
L}(\overline Q))+120\,Q\overline{\mathcal L}(\mathcal L(\overline
Q))+
\\
&+120\,\mathcal L(Q)\overline{\mathcal L}(\overline
Q)+858\,\overline{\mathcal L}(\mathcal L(Q))\overline Q-
120\,Q\overline{\mathcal L}(\overline{\mathcal
L}(Q))+1305\,A\overline{\mathcal L}(Q)-2610\,A\mathcal L(\overline
Q)+
\\
 &+1350\,A\overline{\mathcal L}(\overline
B)-270\,P\overline{\mathcal L}(B)+90\,P\overline{\mathcal
L}(\overline Q)+ 810\,\overline P\mathcal L(\overline
B)-270\,\mathcal L(Q)\overline P+52\,\mathcal L(Q)\overline Q^2+
\\
&+ 1827\,\overline B\mathcal L(\overline P)-2745B\overline{\mathcal
L}(P)+2808\,\overline B^2\overline{\mathcal
L}(B)-3780\,\big(\overline{\mathcal L}(\overline
B)\big)^2-120\,\big(\overline{\mathcal
L}(Q)^2\big)-180\,\big(\mathcal L(\overline Q)^2\big)-
\\
&-68\,Q\overline Q\mathcal L(\overline Q)+ 156\,BQ\mathcal
L(\overline Q)+11394\,\overline{\mathcal L}(\overline B)B\overline
B-3042\,\overline{\mathcal L}(\overline B)BQ+456\,\overline
B\overline Q\mathcal L(\overline Q)- 2052\,B\overline B\mathcal
L(\overline Q)-
\\
&-3978\,B\overline B\overline{\mathcal L}(Q)+3366\,B\overline
Q\mathcal L(\overline B)-540\,\overline BQ\overline{\mathcal
L}(B)+150\,\overline BQ\overline{\mathcal L}(\overline Q)+
894\,\overline B\overline Q\overline{\mathcal L}(Q)-5670\,B^2\mathcal
L(\overline B)-
\\
&-210\,B^2Q^2+4860\mathcal L(\overline B)\overline{\mathcal
L}(B)+1800\,\overline{\mathcal L}(\overline B)\overline{\mathcal
L}(Q)-336\,\overline Q\overline{\mathcal
L}(P)+1440\,\overline{\mathcal L}(\overline B)\mathcal L(\overline
Q)-
\\
&-450\,\mathcal L(\overline B)\overline{\mathcal L}(\overline Q)-
378\,\overline B^2\overline{\mathcal L}(\overline Q)-267\,\mathcal
L(\overline B)\overline Q^2+24\,Q\mathcal L(\overline P)-90\,\mathcal
L(Q)\overline{\mathcal L}(B)+630\,B^2\mathcal L(Q)-
\\
&-192\,\overline{\mathcal L}(B)Q^2+12\,Q^2\overline{\mathcal
L}(\overline Q)-180\,\mathcal L(\overline Q)\overline{\mathcal
L}(Q)-2934\,\overline{\mathcal L}(\overline B)\overline B\overline
Q+780\,\overline QQ\overline{\mathcal L}(\overline
B)+294\,QB\overline{\mathcal L}(Q)-
\\
&-546\,B\mathcal L(Q)\overline Q-88\,Q\overline{\mathcal
L}(Q)\overline Q+3780\,B^2P-1674\,\overline B^2\overline
P-432\,\overline B^2\overline Q^2+3753\,AB\overline
B-444\,QA\overline Q-
\\
&-1584\,ABQ+873\,\overline B\overline QA- 1728\,BP\overline
Q+90\,\overline B\overline PQ+405\,A^2+2583\,\overline B^2\,\overline
QB+3096\,B^2\overline BQ-
\\
&-20\,Q^2\overline Q^2-1242\,B^2\overline B^2+278\,\overline
Q^2Q\overline B+234\,\overline QBQ^2-2868\,\overline B\,\overline Q
BQ\bigg].
\endaligned
\]

Similarly to what we did for the second loop, non-vanishing of either
of the above two essential torsion coefficients enables us to
determine either ${\sf a}$ or ${\sf a}\overline{\sf a}$, and to yet
consider the remaining ones as the invariants of the problem. Hence,
pursuant to Theorem \ref{th-loop-2}, we have:

\begin{Theorem}
\label{th-loop-3} Assume that the four expressions $V'''_1$,
$W'''_5$, $W'''_9$, $X_1'''$
vanish identically.

\begin{itemize}

\smallskip\item[{\bf (vi)}]
When $W'''_2\neq 0$, one can normalize ${\sf a}\overline{\sf a}$ and
$W'''_1$ is the single invariant of the equivalence problem.

\smallskip\item[{\bf (vii)}]
When
$W'''_1\neq 0$ but $W_2 ''' \equiv 0$,
one can normalize ${\sf a}$ from ${\sf a}^2 \overline{\sf a}^3$ and
the equivalence problem
has no invariant.

\smallskip\item[{\bf (viii)}]
Otherwise, when both these two
normalizable expressions vanish identically, one has to start
the prolongation procedure under the assumptions:
\begin{equation}
\label{assum-third-loop}\aligned
0
&
\equiv
V'''_1
\equiv
W'''_5
\equiv
W'''_9
\equiv
X_1''',
\\
0
&
\equiv
W'''_1
\equiv
W'''_2
\equiv 0.
\endaligned
\end{equation}
\end{itemize}
\end{Theorem}

\subsection{Prolongation}
After determining the four group parameters $\sf b,c,d,e$ and in the
case that the equations \thetag{\ref{assum-third-loop}} hold, we
have to prolong the
equivalence problem with one undetermined group parameter $\sf
a$ and with one Maurer-Cartan 1-form:
\[
\beta
=
\frac{d\sf a}{\sf a}.
\]
Here the structure equations have the form ({\it cf.
\thetag{\ref{struc-equ-third-loop}}}):
\begin{equation}
\label{Struc-eq-after-third-loop}
 \aligned d\sigma & = \big(2\,\beta+\overline{\beta}\big) \wedge
\sigma+
\\
& \ \ \ \ \ + U_1'''\,\sigma\wedge\overline{\sigma} +
U_2'''\,\sigma\wedge\rho + U_3'''\,\sigma\wedge\zeta +
U_4'''\,\sigma\wedge\overline{\zeta} + {\rho\wedge\zeta}, \ \ \ \
\endaligned
\end{equation}
\[
\aligned d\rho & = \big(\beta+\overline{\beta}\big)\wedge\rho +
\\
& \ \ \ \ \ +  V_2'''\,\sigma\wedge\rho +
\overline{V_2}'''\,\overline{\sigma}\wedge\rho +
 V_8'''\,\rho\wedge\zeta +
\overline{V_8}'''\,\rho\wedge\overline{\zeta} +
i\,\zeta\wedge\overline{\zeta},
\endaligned
\]
\[
\aligned d\zeta & = \beta\wedge\zeta +
\\
& \ \ \ \ \ + W_3'''\,\sigma\wedge\zeta +
W_6'''\,\overline{\sigma}\wedge\zeta+ W_8'''\,\rho\wedge\zeta +
W_{10}'''\,\zeta\wedge\overline{\zeta}, \ \ \ \ \ \ \ \ \ \ \ \ \
\endaligned
\]
which, after the replacement:
\[
\beta\longmapsto\beta+p\,\sigma+q\,\overline\sigma
+
r\,\rho+s\,\zeta+t\,\overline\zeta
\]
take the form:
\begin{equation*}
\!\!\!\!\!\!\!\!\!\!\!\!\!\!\!\!\!\!\!\!
\!\!\!\!\!\!\!\!\!\!\!\!\!\!\!\!\!\!\!\!
\small
\aligned d\sigma & = \big(2\,\beta+\overline{\beta}\big) \wedge
\sigma+
\\
& \ \ \ \ \ + \big(U_1'''-\overline
p-2\,q\big)\,\sigma\wedge\overline{\sigma} + \big(U_2'''-2\,r-\overline
 r\big)\,\sigma\wedge\rho + \big(U_3'''-2\,s-\overline
 t\big)\,\sigma\wedge\zeta +
\big(U_4'''-2\,t-\overline  s\big)\,\sigma\wedge\overline{\zeta} +
{\rho\wedge\zeta},
\endaligned
\end{equation*}
\[
\!\!\!\!\!\!\!\!\!\!\!\!\!\!\!\!\!\!\!\!
\!\!\!\!\!\!\!\!\!\!\!\!\!\!\!\!\!\!\!\!
\small
\aligned d\rho & = \big(\beta+\overline{\beta}\big)\wedge\rho +
\\
& \ \ \ \ \ +  \big(V_2'''+p+\overline q\big)\,\sigma\wedge\rho +
\big(\overline{V_2}'''+\overline p+q\big)\,\overline{\sigma}\wedge\rho
+
 \big(V_8'''-s-\overline
t\big)\,\rho\wedge\zeta + \big(\overline{V_8}'''-\overline
s-t\big)\,\rho\wedge\overline{\zeta} +
i\,\zeta\wedge\overline{\zeta}, \ \ \ \ \ \
\endaligned
\]
\[
\!\!\!\!\!\!\!\!\!\!\!\!\!\!\!\!\!\!\!\!
\!\!\!\!\!\!\!\!\!\!\!\!\!\!\!\!\!\!\!\!
\small
\aligned d\zeta & = \beta\wedge\zeta +
\\
& \ \ \ \ \ + \big(W_3'''+p\big)\,\sigma\wedge\zeta +
\big(W_6'''+q\big)\,\overline{\sigma}\wedge\zeta+
\big(W_8'''+r\big)\,\rho\wedge\zeta +
\big(W_{10}'''-t\big)\,\zeta\wedge\overline{\zeta}. \ \ \ \ \ \ \ \ \
\ \ \ \ \ \ \ \ \ \ \ \ \ \ \ \ \ \ \ \ \ \ \ \ \ \ \ \ \ \ \ \ \ \ \
\
\endaligned
\]
Taking account of the vanishing expressions of the previous subsections,
one verifies that it is possible to annihilate all the above coefficients
by determining:
\[
p:=-W'''_3, \ \ \ \ q:=-W'''_6, \ \ \ \ r:=-W'''_8, \ \ \ \
s:=V'''_8-\overline W'''_{10}, \ \ \ \ t:=W'''_{10}.
\]
Then, putting:
\begin{equation}
\label{beta} \beta:=\frac{d\sf a}{\sf
a}+W'''_3\,\sigma+W'''_6\,\overline\sigma+W'''_8\,\rho+\big(\overline
W'''_{10}-V'''_8\big)\,\zeta-W'''_{10}\,\overline\zeta,
\end{equation}
the structure equations convert into the form:
\begin{equation}
\label{prolongation-1}
\aligned
d\sigma
&
=
\big(2\,\beta+\overline\beta\big)\wedge\sigma+\rho\wedge\zeta,
\\
d\rho&=\big(\beta+\overline\beta\big)\wedge\rho+i\,\zeta\wedge\overline\zeta,
\\
d\zeta&=\beta\wedge\zeta.
\endaligned
\end{equation}

Before starting the prolongation procedure similarly to the
procedure we performed in subsection~\ref{involutive}, one
should observe the non-involutiveness of the above structure
equations. Since all of the coefficients $p,q,r,s,t$ were
determined, there remains no free variable, hence
the above modified $\beta$ is the {\it unique} 1-form enjoying the
structure equations. According to Proposition
\ref{lifted-prop}, we can therefore transform the $G^\approx$-structure
equivalence problem of the 5-dimensional base manifolds
$M^5\subset\mathbb C^4$ to the $\{e\}$-{\it structure} problem on the
7-dimensional prolonged spaces $M^5\times G^\approx\subset\mathbb
C^5:=\mathbb C\{z,w_1,w_2,w_3,{\sf a}\}$ with $G^\approx$ as follows:
\begin{equation*}\footnotesize\aligned
G^\approx:=\left\{g^\approx:=
\left(%
\begin{array}{ccccc}
{\sf a}\overline{\sf a}^2 & 0 & 0 & 0 & 0 \\
0 & {\sf a}^2\overline{\sf a} & 0 & 0 & 0 \\
0 & 0 & {\sf a}\overline{\sf a} & 0 & 0 \\
0 & 0 & 0 & \overline{\sf a} & 0 \\
0 & 0 & 0 & 0 & {\sf a} \\
\end{array}%
\right) : {\sf a}\in\mathbb C\right\}.
\endaligned
\end{equation*}
In this case, the lifted coframe
$\{\sigma,\overline\sigma,\rho,\zeta,\overline\zeta\}$ will be
extended by the two 1-forms $\beta$ and $\overline\beta$. Hence, to find
the structure equations of the new equivalence problem, it suffices
to extend the structure equations \thetag{\ref{prolongation-1}} by
computing the exterior derivation $d\beta$ of $\beta$ in \thetag{\ref{beta}} in terms of
these 7 lifted 1-forms.

To do this, first let us compute $d\beta$
directly, taking account of Lemma \ref{Lemma-G-0}:
\begin{equation}
\!\!\!\!\!\!\!\!\!\!\!\!\!\!\!\!\!\!\!\!
\!\!\!\!\!\!\!\!\!\!\!\!\!\!\!\!\!\!\!\!
\footnotesize\aligned
\label{dbeta-1} d\beta&=\zero{d\big(\frac{d\sf a}{\sf
a}\big)}+W'''_3\,d\sigma+W'''_6\,d\overline\sigma+W'''_8\,d\rho+\big(\overline
W'''_{10}-V'''_8\big)\,d\zeta-W'''_{10}\,d\overline\zeta+
\\
&+\bigg(\mathcal S(W'''_3)\,\sigma_0+\overline{\mathcal
S}(W'''_3)\,\overline\sigma_0+\mathcal T(W'''_3)\,\rho_0+\mathcal
L(W'''_3)\,\zeta_0+\overline{\mathcal
L}(W'''_3)\,\overline\zeta_0\bigg)\wedge\sigma+
\\
&+\bigg(\mathcal S(W'''_6)\,\sigma_0+\overline{\mathcal
S}(W'''_6)\,\overline\sigma_0+\mathcal T(W'''_6)\,\rho_0+\mathcal
L(W'''_6)\,\zeta_0+\overline{\mathcal
L}(W'''_6)\,\overline\zeta_0\bigg)\wedge\overline\sigma+
\\
&+\bigg(\mathcal S(W'''_8)\,\sigma_0+\overline{\mathcal
S}(W'''_8)\,\overline\sigma_0+\mathcal T(W'''_8)\,\rho_0+\mathcal
L(W'''_8)\,\zeta_0+\overline{\mathcal
L}(W'''_8)\,\overline\zeta_0\bigg)\wedge\rho+
\\
&+\bigg(\mathcal S(\overline
W'''_{10}-V'''_8)\,\sigma_0+\overline{\mathcal S}(\overline
W'''_{10}-V'''_8)\,\overline\sigma_0+\mathcal T(\overline
W'''_{10}-V'''_8)\,\rho_0+\mathcal L(\overline
W'''_{10}-V'''_8)\,\zeta_0+\overline{\mathcal L}(\overline
W'''_{10}-V'''_8)\,\overline\zeta_0\bigg)\wedge\zeta-
\\
&-\bigg(\mathcal S(W'''_{10})\,\sigma_0+\overline{\mathcal
S}(W'''_{10})\,\overline\sigma_0+\mathcal T(W'''_{10})\,\rho_0+\mathcal
L(W'''_{10})\,\zeta_0+\overline{\mathcal
L}(W'''_{10})\,\overline\zeta_0\bigg)\wedge\overline\zeta.
\endaligned
\end{equation}
We need the following useful lemma:

\begin{Lemma}
\label{Lemma-dbeta}
 The exterior derivation $d\beta$ of the 1-form
$\beta$ has the form:
\begin{equation*}
d\beta:=T_1\,\sigma\wedge\zeta+T_2\,\overline\sigma\wedge\zeta+T_3\,\rho\wedge\zeta+T_4\,\zeta\wedge\overline\zeta,
\end{equation*}
for some certain functions $T_1,\ldots,T_4$.
\end{Lemma}

\proof According to \thetag{\ref{dbeta-1}}, the expression of
$d\beta$ admits no variable $\sf a$ or 1-form $d\sf a$ in its
expression. Hence, it will be independent of any wedge product of the
form $\bullet\wedge\beta$ and $\bullet\wedge\overline\beta$. On the
other hand, differentiating the expression of $d\zeta$ in
\thetag{\ref{prolongation-1}} gives:
\[
0\equiv d\beta\wedge\zeta-\beta\wedge
d\zeta=d\beta\wedge\zeta-\zero{\beta\wedge\beta\wedge\zeta},
\]
which, according to Cartan's Lemma \ref{Cartan-lemma}, implies that:
\[
d\beta:=\mathcal F\wedge\zeta
\]
for some certain 1-form $\mathcal F$. Now, since seven 1-forms
$\sigma,\overline\sigma,\rho,\zeta,\overline\zeta,\beta,\overline\beta$
constitute a basis for the set of all 1-forms on the prolonged space,
$d\beta$ may be expressed as follows:
\[
d\beta:=T_1\,\sigma\wedge\zeta+T_2\,\overline\sigma\wedge\zeta+T_3\,\rho\wedge\zeta+T_4\,\zeta\wedge\overline\zeta,
\]
for some certain functions $T_1,\ldots,T_4$, as claimed.
\endproof

According to this lemma, to compute the exterior derivation
$d\beta$, it suffices to compute only four coefficients
$T_1,\ldots,T_4$ instead of computing all 21 coefficients of the
possible wedge products between seven 1-forms
$\sigma,\overline\sigma,\rho,\zeta,\overline\zeta,\beta,\overline\beta$
in \thetag{\ref{dbeta-1}}. Extracting the coefficients of
$\sigma\wedge\zeta,\overline{\sigma}
\wedge\zeta,\rho\wedge\zeta,\overline\zeta\wedge\zeta$
in the expression \thetag{\ref{dbeta-1}} gives:
\begin{equation}
\!\!\!\!\!\!\!\!\!\!\!\!\!\!\!\!\!\!\!\!
\!\!\!\!\!\!\!\!\!\!\!\!\!\!\!\!\!\!\!\!\!\!\!\!\!
\label{T} \footnotesize\aligned T_1&=-\frac{1}{\sf a}\mathcal
L(W'''_3)+\frac{1}{{\sf a}^2\overline{\sf a}}\mathcal S(\overline
W'''_{10}-V'''_8)-\frac{\sf c}{{\sf a}^3\overline{\sf a}^2}\mathcal
T(\overline W'''_{10}-V'''_8)+\frac{{\sf bc}-{\sf a}\overline{\sf a}\sf
e}{{\sf a}^4\overline{\sf a}^2}\,\mathcal L(\overline
W'''_{10}-V'''_8)+\frac{\overline{\sf b}{\sf c}-{\sf a}\overline{\sf
a}\sf d}{{\sf a}^3\overline{\sf a}^3}\,\overline{\mathcal
L}(\overline W'''_{10}-V'''_8),
\\
T_2&=-\frac{1}{\sf a}\mathcal L(W'''_6)+\frac{1}{{\sf a}\overline{\sf
a}^2}\overline{\mathcal S}(\overline
W'''_{10}-V'''_8)-\frac{\overline{\sf c}}{{\sf a}^2\overline{\sf
a}^3}\mathcal T(\overline W'''_{10}-V'''_8)+ \frac{{\sf b}\overline{\sf
c}-{\sf a}\overline{\sf a}\overline{\sf d}}{{\sf a}^3\overline{\sf
a}^3}\,\mathcal L(\overline W'''_{10}-V'''_8)+\frac{\overline{\sf
b}\overline{\sf c}-{\sf a}\overline{\sf a}\,\overline{\sf e}}{{\sf
a}^2\overline{\sf a}^4}\,\overline{\mathcal L}(\overline
W'''_{10}-V'''_8),
\\
T_3&=W'''_3-\frac{1}{\sf a}\,\mathcal L(W'''_8)+\frac{1}{{\sf
a}\overline{\sf a}}\,\mathcal T(\overline W'''_{10}-V'''_8)-\frac{\sf
b}{{\sf a}^2\overline{\sf a}}\,\mathcal L(\overline
W'''_{10}-V'''_8)-\frac{\overline{\sf b}}{{\sf a}\overline{\sf
a}^2}\,\overline{\mathcal L}(\overline W'''_{10}-V'''_8),
\\
T_4&=i\,W'''_8-\frac{1}{\overline{\sf a}}\,\overline{\mathcal
L}(\overline W'''_{10}-V'''_8)-\frac{1}{\sf a}\,\mathcal L(W'''_{10}).
\endaligned
\end{equation}

Therefore, the $\{e\}$-structure equivalence problem on the prolonged
space $M\times G^\approx$ enjoys the structure equations of the form:
\begin{equation}
\label{struct-eq-after prolongation} \aligned
d\sigma&=\big(2\,\beta+\overline\beta\big)\wedge\sigma+\rho\wedge\zeta,
\\
d\rho&=\big(\beta+\overline\beta\big)\wedge\rho+i\,\zeta\wedge\overline\zeta,
\\
d\zeta&=\beta\wedge\zeta,
\\
d\beta&=T_1\,\sigma\wedge\zeta+T_2\,\overline\sigma\wedge\zeta+T_3\,\rho\wedge\zeta+T_4\,\zeta\wedge\overline\zeta,
\endaligned
\end{equation}
with four essential invariants $T_1,T_2,T_3,T_4$ as above.

\begin{Theorem}
In the case that six functions $V'''_1$, $W'''_5$,
$W'''_9$, $X_1'''$, $W'''_1$ and
$W'''_2$ vanish identically (see \thetag{\ref{new}} and
\thetag{\ref{new-new}} for their expressions), then the under
consideration equivalence problem has the four essential invariants
$T_1,T_2,T_3$ and $T_4$ as \thetag{\ref{T}}, in terms of the five
complex variables $z,w_1,w_2,w_3,{\sf a}$.
\end{Theorem}

In Section \ref{equivalence-model}, we considered the equivalence problem of an arbitrary 5-dimensional CR-manifold $M^5$ to the cubic model $M^5_{\sf c}$. In Theorem \ref{theorem-cubic}, we observed that such equivalency holds if and only if the structure equations associated to the lifted coframe $\{\sigma,\overline\sigma,\rho,\zeta,\overline\zeta,\alpha,\overline\alpha\}$ of $M^5$ takes the structure equations as \thetag{\ref{d4}}. Now, by a careful comparison between two structure equations \thetag{\ref{d4}} and \thetag{\ref{struct-eq-after prolongation}} and thanks to the above theorem, one realizes that these two structure equations take the same form whenever the appearing invariants vanish, identically. Then we have:

\begin{Corollary}
\label{complete-model}
An arbitrary 5-dimensional CR-manifold $M^5$ is equivalence, through some biholomorphism, to the cubic model $M^5_{\sf c}$ if and only if we have:
\[
\aligned
0&\equiv V'''_1=W'''_5=W'''_9=X_1'''=W'''_1,
\\
0&\equiv T_1=T_2=T_3=T_4.
\endaligned
\]
\qed
\end{Corollary}

This corollary completes the procedure at the end of Section \ref{equivalence-model}. On the other hand, in \cite[Proposition 12]{Beloshapka-2004}, Beloshapka proved that the cubic model $M^5_{\sf c}$ is the most symmetric nondegenerate surface, {\it i.e.} the dimension of the symmetry group of each $M^5$ is not greater than that of $M^5_{\sf c}$. Now, granted the above result one finds out that\,\,---\,\,{\it see} also Corollary \ref{corollary-cubic}\,\,---\,\,:

\begin{Corollary}
If an arbitrary CR-manifold $M^5$ satisfies the assumptions of the above Corollary \ref{complete-model}, then its symmetry group has the maximum dimension, equal to the dimension of the symmetry group of $M^5_{\sf c}$.
\qed
\end{Corollary}

\section{The branch $R \not\equiv 0$}
\label{branch-R-neq-0}
\HEAD{\ref{branch-R-neq-0}.~The branch $R \neq 0$}{
Masoud {\sc Sabzevari} (Shahrekord) and Jo\"el {\sc Merker} (LM-Orsay)}

\subsection{Normalization of four group parameters}

After performing the long and complicated computations of the under
consideration equivalence problem with the assumption $R\equiv 0$ in
Sections \ref{branch-R-0} and \ref{four-normalizations}, now we have somehow simpler
computations to conclude the memoir by inspecting the equivalence
problem in the case $R\not\equiv 0$. For this aim, we have to return
to the subsection \ref{First loop normalization} of the first loop
normalization of the structure equations
\thetag{\ref{structure-firts-absorbtion}}:
\begin{equation}
\label{structure-firts-absorbtion-R-neq-0} \aligned
d\sigma & = \big(2\,\alpha_1+\overline{\alpha}_1\big) \wedge \sigma+
\\
& \ \ \ \ \ + U_1\,\sigma\wedge\overline{\sigma} +
U_2\,\sigma\wedge\rho + U_3\,\sigma\wedge\zeta +
U_4\,\sigma\wedge\overline{\zeta} +
\\
& \ \ \ \ \ \ \ \ \ \ \ \ \ \ \ \ \ \ \ \ \ \ \ \ +
U_5\,\overline{\sigma}\wedge\rho + U_6\,\overline{\sigma}\wedge\zeta
+ U_7\,\overline{\sigma}\wedge\overline{\zeta} +
\\
& \ \ \ \ \ \ \ \ \ \ \ \ \ \ \ \ \ \ \ \ \ \ \ \ \ \ \ \ \ \ \ \ \ \
\ \ \ \ \ \ \ \ \ + {\rho\wedge\zeta},
\endaligned
\end{equation}
\[
\aligned d\rho & = \alpha_2\wedge\sigma +
\overline{\alpha}_2\wedge\overline{\sigma} + \alpha_1\wedge\rho +
\overline{\alpha}_1\wedge\rho +
\\
& \ \ \ \ \ + V_1\,\sigma\wedge\overline{\sigma} +
V_2\,\sigma\wedge\rho + V_3\,\sigma\wedge\zeta +
V_4\,\sigma\wedge\overline{\zeta} +
\\
& \ \ \ \ \ \ \ \ \ \ \ \ \ \ \ \ \ \ \ \ \ \ \ \ +
\overline{V_2}\,\overline{\sigma}\wedge\rho +
\overline{V_4}\,\overline{\sigma}\wedge\zeta +
\overline{V_3}\,\overline{\sigma}\wedge\overline{\zeta} +
\\
& \ \ \ \ \ \ \ \ \ \ \ \ \ \ \ \ \ \ \ \ \ \ \ \ \ \ \ \ \ \ \ \ \ \
\ \ \ \ \ \ \ \ \ + V_8\,\rho\wedge\zeta +
\overline{V_8}\,\rho\wedge\overline{\zeta} +
\\
& \ \ \ \ \ \ \ \ \ \ \ \ \ \ \ \ \ \ \ \ \ \ \ \ \ \ \ \ \ \ \ \ \ \
\ \ \ \ \ \ \ \ \ \ \ \ \ \ \ \ \ \ \ \ \ \ \ \ \ \ \ \ \ \ +
i\,\zeta\wedge\overline{\zeta},
\endaligned
\]
\[
\aligned d\zeta & = \alpha_3\wedge\sigma +
\alpha_4\wedge\overline{\sigma} + \alpha_5\wedge\rho +
\alpha_1\wedge\zeta +
\\
& \ \ \ \ \ + W_1\,\sigma\wedge\overline{\sigma} +
W_2\,\sigma\wedge\rho + W_3\,\sigma\wedge\zeta +
W_4\,\sigma\wedge\overline{\zeta} +
\\
& \ \ \ \ \ \ \ \ \ \ \ \ \ \ \ \ \ \ \ \ \ \ \ \ +
W_5\,\overline{\sigma}\wedge\rho + W_6\,\overline{\sigma}\wedge\zeta
+ W_7\,\overline{\sigma}\wedge\overline{\zeta} +
\\
& \ \ \ \ \ \ \ \ \ \ \ \ \ \ \ \ \ \ \ \ \ \ \ \ \ \ \ \ \ \ \ \ \ \
\ \ \ \ \ \ \ \ \ + W_8\,\rho\wedge\zeta +
W_9\,\rho\wedge\overline{\zeta} +
\\
& \ \ \ \ \ \ \ \ \ \ \ \ \ \ \ \ \ \ \ \ \ \ \ \ \ \ \ \ \ \ \ \ \ \
\ \ \ \ \ \ \ \ \ \ \ \ \ \ \ \ \ \ \ \ \ \ \ \ \ \ \ \ \ +
W_{10}\,\zeta\wedge\overline{\zeta}.
\endaligned
\]
In this early
normalization step, we found the following four\footnote{In fact they
were five but the second and the fifth ones were conjugate of each
other.} normalizable expressions in \thetag{\ref{normalizable-first
loop}}:

\begin{equation}
\label{normalizable-first loop-2}\aligned U_5 & =
\frac{1}{\overline{\sf a}^2}\,\overline{G} - \frac{\sf b}{{\sf
a}\overline{\sf a}^2}\,B - \frac{\overline{\sf b}}{\overline{\sf
a}^3}\, \overline{R} + \frac{\overline{\sf d}}{{\sf a}\overline{\sf
a}^2},
\\
U_6 & = \frac{1}{\overline{\sf a}}\,B - \frac{\overline{\sf c}}{{\sf
a}\overline{\sf a}^2},
\\
U_7 & = \frac{{\sf a}}{\overline{\sf a}^2}\,\overline{R},
\\
U_3+\overline{U}_4-3\,V_8 & = \frac{1}{\sf a}\,Q - 4\,\frac{\sf
c}{{\sf a}^2\overline{\sf a}} + \frac{1}{\sf a}\,\overline{B} - 3i\,
\frac{\overline{\sf b}}{{\sf a}\overline{\sf a}}.
\endaligned
\end{equation}
Normalizing these expressions, this time with the assumption
$R\not\equiv 0$, enables us to determine the first group parameter
$\sf a$ besides $\sf b,c$ and $\sf d$ ({\it cf.} subsection
\ref{Normalization-1-r=0}). Normalizing the third expression to $1$
and the remaining ones to $0$, and taking account of the expression
of $G$ in Lemma \ref{E-F-G}, one obtains:
\begin{equation}
\boxed{\aligned \label{a-b-c-d} {\sf a}&:={\bf A}_0,
\\
{\sf b}&:={\bf A}_0\,\big(-i\,B+\textstyle{\frac{i}{3}}\,\overline
Q\big),
\\
{\sf c}&:={\bf A}_0\,{\bf \overline A}_0\,B,
\\
{\sf d}&:=\overline{\bf A}_0\,\Big(i\,\overline{\mathcal
L}(R)-i\,\mathcal L(\overline B)+\textstyle{\frac{4i}{3}}\,\overline
QR+iP+\textstyle{\frac{2i}{3}}\,\overline BQ-2i\,BR\Big),
\endaligned}
\end{equation}
where ${\bf A}_0$ is a nonzero complex function satisfying $\frac{{\bf
A}_0^2}{\overline{\bf A}_0}=R$.

\subsection{Changing up the initial coframe}

After determining four of the five group parameters in the previous
section, now the equivalence problem takes the form:
\begin{eqnarray*}
\left[
\begin{array}{l}
\sigma={\bf A}_0^2{\bf \overline A}\cdot\sigma_0,
\\
\rho=\rho_0^\dag,
\\
\zeta={\sf e}\cdot\sigma_0+\zeta_0^\dag,
\end{array}
\right.
\end{eqnarray*}
where the two new modified initial 1-forms are:
\[
\!\!\!\!\!\!\!\!\!\!\!\!\!\!\!\!\!\!\!\!
\aligned
\rho_0^\dag&={\bf A}_0{\bf\overline A}_0\big(\overline
B\sigma_0+B\overline\sigma_0+\rho_0\big),
\\
\zeta_0^\dag&={\bf A}_0\,
\Big[\big(-i\,{\mathcal L}(\overline
R)+i\,\overline{\mathcal L}(B)-\textstyle{\frac{4i}{3}}\,Q\overline
R-i\overline P-\textstyle{\frac{2i}{3}}\,B\overline Q+2i\,\overline
B\,\overline
R\big)\overline\sigma_0+
\\
&
\ \ \ \ \ \ \ \ \ \
+
\big(-i\,B+\textstyle{\frac{i}{3}}\,\overline
Q\big)\sigma_0+\zeta_0\Big].
\endaligned
\]
Under this new initial coframe
$(\overline\sigma_0,\sigma_0,\rho_0^\dag,
\overline\zeta_0^\dag,\zeta_0^\dag)^t$,
the ambiguity group reduces to the form:
\[
\footnotesize\aligned
G^\dag := \left\{ g^\dag = \left(\!\!
\begin{array}{ccccc}
{\bf A}_0\overline{\bf A}_0^2 & 0 & 0 & 0 & 0
\\
0 & {\bf A}_0^2\overline{\bf A}_0 & 0 & 0 & 0
\\
0 & 0 & 1 & 0 & 0
\\
\overline{\sf e} & 0 & 0 & 1 & 0
\\
0 & {\sf e} & 0 & 0 & 1
\end{array}
\!\!\right) \colon {\sf e}\in\C \right\},
\endaligned
\]
and the new Maurer-Cartan form is:
\[
\aligned
dg^\dag\cdot {g^\dag}^{-1}=\left(%
\begin{array}{ccccc}
0 & 0 & 0 & 0 & 0 \\
0 & 0 & 0 & 0 & 0 \\
0 & 0 & 0 & 0 & 0 \\
\overline\gamma & 0 & 0 & 0 & 0 \\
0 & \gamma & 0 & 0 & 0 \\
\end{array}%
\right),
\endaligned
\]
where:
\[
\gamma:=\frac{d\sf e}{{\bf A}_0^2\overline{\bf A}_0}.
\]

\subsection{Second-loop absorbtion and normalization}
Similar to what we did in subsection \ref{Second-loop-R=0}, let us
introduce the notations:
\[
\aligned {\sf a} & := {\bf A}_0,\ \ \ \  \ \ \ \ \ \ \ \ \ \ \ \ \ \
\ \ \ \ \  \ \ \  \ \text{\rm where:}\ \ \ \ \ \frac{{\bf
A}_0^2}{\overline{\bf A}_0}=R,
\\
{\sf b} & := {\bf A}_0\,{\bf B}'_0, \ \ \  \ \ \ \ \ \ \ \ \ \ \ \ \
\ \ \ \ \ \text{\rm where:}\ \ \ \ \ {\bf B}'_0 := -\,i\,B +
{\textstyle{\frac{i}{3}}}\,\overline{Q},
\\
{\sf c} & := {\bf A}_0\overline{\bf A}_0\,{\bf C}'_0, \ \ \  \ \ \ \
\ \ \ \ \ \ \ \ \ \text{\rm where:}\ \ \ \ \ {\bf C}'_0 :=
\overline{B},
\\
{\sf d} & := \overline{\bf A}_0\,{\bf D}'_0, \ \  \ \ \ \ \ \ \ \ \ \
\ \ \ \ \ \  \ \ \ \text{\rm where:}\ \ \ \ \ {\bf D}'_0 :=
i\,\overline{\mathcal L}\big(R\big)-\,i\,
\mathcal{L}\big(\overline{B}\big) +
\\
&
\ \ \ \ \ \ \ \ \ \ \ \ \ \ \ \ \ \ \ \ \ \ \ \ \ \ \ \ \ \ \ \ \ \ \ \
\ \ \ \ \ \ \ \ \ \ \ \ \ \ \ \ \ \ \ \ \ \ \ \ \ \ \ \ \
+
\textstyle{\frac{4i}{3}}\,\overline Q R+ i\,P +
{\textstyle{\frac{2i}{3}}}\, \overline{B}\,Q-2i\,BR.
\endaligned
\]
In this case, the new lifted coframe takes the form:
\[
\aligned \sigma & = {\bf A}_0^2\overline{\bf A}_0\cdot\sigma_0,
\\
\rho & = {\bf A}_0\overline{\bf A}_0\cdot \big( {\bf
C}'_0\,\sigma_0+\overline{\bf C}'_0\,\overline{\sigma}_0+\rho_0
\big),
\\
\zeta & = {\sf e}\cdot\sigma_0 + {\bf A}_0\cdot \big( \overline{\bf
D}'_0\,\overline{\sigma}_0 + {\bf B}'_0\,\rho_0 + \zeta_0 \big).
\endaligned
\]
To proceed with the Cartan's method, one has to consider the exterior
differentiation of these expressions:
\[ \aligned d\sigma&={\bf A}_0^2\overline{\bf A}_0\cdot
d\overline\sigma_0+
\\
&+2\,{\bf A}_0\overline{\bf A}_0\,d{\bf A}_0\wedge\sigma_0+{\bf
A}_0^2\,d\overline{\bf A}_0\wedge\sigma_0,
\\
d\rho&={\bf A}_0\overline{\bf A}_0\cdot \big( {\bf
C}'_0\,d\sigma_0+\overline{\bf C}'_0\,d\overline{\sigma}_0+d\rho_0
\big)+
\\
&+ \big({\bf A}_0\overline{\bf A}_0\,d{\bf C}'_0+{\bf A}_0{\bf
C}'_0\,d\overline{\bf A}_0+\overline{\bf A}_0{\bf C}'_0\,d{\bf
A}_0\big)\wedge\sigma_0+
\\
&+\big({\bf A}_0\overline{\bf A}_0\,d\overline{\bf C}'_0+{\bf
A}_0\overline{\bf C}'_0\,d\overline{\bf A}_0+\overline{\bf
A}_0\,\overline{\bf C}'_0\,d{\bf A}_0\big)\wedge\overline\sigma_0+
\\
&+\big(\overline{\bf A}_0\,d{\bf A}_0+{\bf A}_0\,d\overline{\bf
A}_0\big)\wedge\rho_0,
\\
d\zeta&=\gamma\wedge\sigma+
\\
&+{\sf e}\cdot d\sigma_0 + {\bf A}_0\cdot \big( \overline{\bf
D}'_0\,d\overline{\sigma}_0 + {\bf B}'_0\,d\rho_0 + d\zeta_0 \big)+
\\
&+\big({\bf A}_0\,d\overline{\bf D}'_0+\overline{\bf D}'_0\,d{\bf
A}_0\big)\wedge\overline\sigma_0+\big({\bf A}_0\,d{\bf B}'_0+{\bf
B}'_0\,d{\bf A}_0\big)\wedge\rho_0+d{\bf A}_0\wedge\zeta_0.
\endaligned
\]
Reading these new structure equations in terms of the lifted coframe
({\it cf.} subsection \ref{Second-loop-R=0}) gives the following
modified structure equations:
\begin{equation}
\label{structure-firts-absorbtion-R-neq-0-NEW} \aligned d\sigma & =
U^{\sf new}_1\,\sigma\wedge\overline{\sigma} + U^{\sf
new}_2\,\sigma\wedge\rho + U^{\sf new}_3\,\sigma\wedge\zeta + U^{\sf
new}_4\,\sigma\wedge\overline{\zeta}+\overline\sigma\wedge\overline\zeta
+ {\rho\wedge\zeta},
\endaligned
\end{equation}
\[
\aligned d\rho & =V^{\sf new}_1\,\sigma\wedge\overline{\sigma} +
V^{\sf new}_2\,\sigma\wedge\rho + V^{\sf new}_3\,\sigma\wedge\zeta +
V^{\sf new}_4\,\sigma\wedge\overline{\zeta} + \ \ \ \ \ \ \ \ \ \ \ \
\ \ \ \ \ \ \ \ \ \ \ \ \ \ \ \
\\
& \ \ \ \ \ \ \ \ \ \ \ \ + \overline V^{\sf
new}_2\,\overline{\sigma}\wedge\rho + \overline V^{\sf
new}_4\,\overline{\sigma}\wedge\zeta + \overline V^{\sf
new}_3\,\overline{\sigma}\wedge\overline{\zeta} +
\\
& \ \ \ \ \ \ \ \ \ \  \ \ \ \ \ \ \ \ + V^{\sf
new}_8\,\rho\wedge\zeta + \overline V^{\sf
new}_8\,\rho\wedge\overline{\zeta} +
\\
& \ \ \ \ \ \ \ \ \ \  \ \ \ \ \ \ \ \ \ \ \ \ \ \ +
i\,\zeta\wedge\overline{\zeta},
\endaligned
\]
\[
\aligned d\zeta & = \gamma\wedge\sigma + \ \ \ \ \ \ \ \ \ \ \ \ \ \
\ \ \ \ \ \ \ \ \ \ \ \ \ \ \ \ \ \ \ \ \ \ \ \ \ \ \ \ \ \ \ \ \ \ \
\ \ \ \ \ \ \ \ \ \ \ \ \ \ \ \ \ \ \ \ \ \ \ \ \ \ \ \ \ \ \ \ \ \ \
 \ \ \ \ \ \ \ \ \ \  \ \ \ \ \ \ \ \ \ \ \ \ \ \ \
\\
& \ \ \ \ \ + W^{\sf new}_1\,\sigma\wedge\overline{\sigma} + W^{\sf
new}_2\,\sigma\wedge\rho + W^{\sf new}_3\,\sigma\wedge\zeta + W^{\sf
new}_4\,\sigma\wedge\overline{\zeta} +
\\
& \ \ \ \ \ \ \ \ \ \ \ \ \ \ \ \ \ \ \ \ \ \ \ \ + W^{\sf
new}_5\,\overline{\sigma}\wedge\rho + W^{\sf
new}_6\,\overline{\sigma}\wedge\zeta + W^{\sf
new}_7\,\overline{\sigma}\wedge\overline{\zeta} +
\\
& \ \ \ \ \ \ \ \ \ \ \ \ \ \ \ \ \ \ \ \ \ \ \ \ \ \ \ \ \ \ \ \ \ \
\ \ \ \ \ \ \ \ \ + W^{\sf new}_8\,\rho\wedge\zeta + W^{\sf
new}_9\,\rho\wedge\overline{\zeta} +
\\
& \ \ \ \ \ \ \ \ \ \ \ \ \ \ \ \ \ \ \ \ \ \ \ \ \ \ \ \ \ \ \ \ \ \
\ \ \ \ \ \ \ \ \ \ \ \ \ \ \ \ \ \ \ \ \ \ \ \ \ \ \ \ \ + W^{\sf
new}_{10}\,\zeta\wedge\overline{\zeta},
\endaligned
\]
with the new extended torsion coefficients which can be presented as
follows after a large amount of simplification\footnote{Of course
here we have ${\sf a}={\bf A}_0$, but in these expressions we still
keep the notation $\sf a$ to better emphasize their appearance in
denominators.}:
\[
\!\!\!\!\!\!\!\!\!\!\!\!\!\!\!\!\!\!\!\!\!\!\!\!\!\!\!\!\!\!\!\!\!
\footnotesize\aligned U^{\sf new}_1&:=U_1-2\,\frac{{\bf
A}_0\overline{\bf A}_0}{{\sf a}^2\overline{\sf
a}}\,\bigg(\frac{1}{{\sf a}\overline{\sf a}^2}\,\overline{\mathcal
S}({\bf A}_0)-\frac{\overline{\sf c}}{{\sf a}^2\overline{\sf
a}^3}\,\mathcal T({\bf A}_0)+\frac{{\sf b}\overline{\sf c}-{\sf
a}\overline{\sf a}\overline{\sf d}}{{\sf a}^3\overline{\sf
a}^3}\,\mathcal L({\bf A}_0)+\frac{\overline{\sf b}\,\overline{\sf
c}-{\sf a}\overline{\sf a}\,\overline{\sf e}}{{\sf a}^2\overline{\sf
a}^4}\,\overline{\mathcal L}({\bf A}_0)\bigg)-
\\
&\ \ \ \ \ \ \ \ \ \ \ -\frac{{\bf A}_0^2}{{\sf a}^2\overline{\sf
a}}\,\bigg(\frac{1}{{\sf a}\overline{\sf a}^2}\,\overline{\mathcal
S}(\overline{\bf A}_0)-\frac{\overline{\sf c}}{{\sf a}^2\overline{\sf
a}^3}\,\mathcal T(\overline{\bf A}_0)+\frac{{\sf b}\overline{\sf
c}-{\sf a}\overline{\sf a}\overline{\sf d}}{{\sf a}^3\overline{\sf
a}^3}\,\mathcal L(\overline{\bf A}_0)+\frac{\overline{\sf
b}\,\overline{\sf c}-{\sf a}\overline{\sf a}\,\overline{\sf e}}{{\sf
a}^2\overline{\sf a}^4}\,\overline{\mathcal L}(\overline{\bf
A}_0)\bigg),
\endaligned
\]
\[
\footnotesize\aligned
U^{\sf new}_2&:=U_2-2\,\frac{{\bf A}_0\overline{\bf A}_0}{{\sf
a}^2\overline{\sf a}}\,\bigg(\frac{1}{{\sf a}\overline{\sf
a}}\,\mathcal T({\bf A}_0)-\frac{\sf b}{{\sf a}^2\overline{\sf
a}}\,\mathcal L({\bf A}_0)-\frac{\overline{\sf b}}{{\sf
a}\overline{\sf a}^2}\,\overline{\mathcal L}({\bf A}_0)\bigg)-
\\
&-
\frac{{\bf A}_0^2}{{\sf a}^2\overline{\sf a}}\,\bigg(\frac{1}{{\sf
a}\overline{\sf a}}\,\mathcal T(\overline{\bf A}_0)-\frac{\sf b}{{\sf
a}^2\overline{\sf a}}\,\mathcal L(\overline{\bf
A}_0)-\frac{\overline{\sf b}}{{\sf a}\overline{\sf
a}^2}\,\overline{\mathcal L}(\overline{\bf A}_0)\bigg),
\\
U^{\sf new}_3&:=U_3-2\,\frac{{\bf A}_0\overline{\bf A}_0}{{\sf
a}^3\overline{\sf a}}\,\mathcal L({\bf A}_0)-\frac{{\bf A}_0^2}{{\sf
a}^3\overline{\sf a}}\,\mathcal L(\overline{\bf A}_0),
\\
U^{\sf new}_4&:=U_4-2\,\frac{{\bf A}_0\overline{\bf A}_0}{{\sf
a}^2\overline{\sf a}^2}\,\overline{\mathcal L}({\bf A}_0)-\frac{{\bf
A}_0^2}{{\sf a}^2\overline{\sf a}^2}\,\overline{\mathcal
L}(\overline{\bf A}_0),
\\
V^{\sf new}_1&:=V_1-\frac{{\bf
A}_0\overline{\bf A}_0}{{\sf a}^2\overline{\sf
a}}\,\bigg(\frac{1}{{\sf a}\overline{\sf a}^2}\,\overline{\mathcal
S}({\bf C}'_0)-\frac{\overline{\sf c}}{{\sf a}^2\overline{\sf
a}^3}\,{\mathcal T}({\bf C}'_0)+\frac{{\sf b}\overline{\sf c}-{\sf
a}\overline{\sf a}\,\overline{\sf d}}{{\sf a}^3\overline{\sf
a}^3}\,{\mathcal L}({\bf C}'_0)+\frac{\overline{\sf b}\,\overline{\sf
c}-{\sf a}\overline{\sf a}\,\overline{\sf e}}{{\sf a}^2\overline{\sf
a}^4}\,\overline{\mathcal L}({\bf C}'_0)\bigg)+
\\
&+\frac{{\bf A}_0\overline{\bf A}_0}{{\sf a}\overline{\sf
a}^2}\,\bigg(\frac{1}{{\sf a}^2\overline{\sf a}}\,{\mathcal
S}(\overline{\bf C}'_0)-\frac{{\sf c}}{{\sf a}^3\overline{\sf
a}^2}\,{\mathcal T}(\overline{\bf C}'_0)+\frac{{\sf b}{\sf c}-{\sf
a}\overline{\sf a}\,{\sf e}}{{\sf a}^4\overline{\sf a}^2}\,{\mathcal
L}(\overline{\bf C}'_0)+\frac{\overline{\sf b}{\sf c}-{\sf
a}\overline{\sf a}{\sf d}}{{\sf a}^3\overline{\sf
a}^3}\,\overline{\mathcal L}(\overline{\bf C}'_0)\bigg),
\\
V^{\sf new}_2&:=V_2-\frac{{\bf A}_0\overline{\bf A}_0}{{\sf
a}^2\overline{\sf a}}\,\bigg(\frac{1}{{\sf a}\overline{\sf
a}}\,\mathcal T({\bf C}'_0)-\frac{\sf b}{{\sf a}^2\overline{\sf
a}}\,\mathcal L({\bf C}'_0)-\frac{\overline{\sf b}}{{\sf
a}\overline{\sf a}^2}\,\overline{\mathcal L}({\bf C}'_0)\bigg)+
\\
&+\frac{\overline{\bf A}_0}{{\sf a}\overline{\sf
a}}\,\bigg(\frac{1}{{\sf a}^2\overline{\sf a}}\,\mathcal S({\bf
A}_0)-\frac{\sf c}{{\sf a}^3\overline{\sf a}^2}\,\mathcal T({\bf
A}_0)+\frac{{\sf bc}-{\sf a}\overline{\sf a}\sf e}{{\sf
a}^4\overline{\sf a}^2}\,\mathcal L({\bf A}_0)+\frac{\overline{\sf
b}{\sf c}-{\sf a}\overline{\sf a}\sf d}{{\sf a}^3\overline{\sf
a}^3}\,\overline{\mathcal L}({\bf A}_0)\bigg)+
\\
&+\frac{{\bf A}_0}{{\sf a}\overline{\sf a}}\,\bigg(\frac{1}{{\sf
a}^2\overline{\sf a}}\,\mathcal S(\overline{\bf A}_0)-\frac{\sf
c}{{\sf a}^3\overline{\sf a}^2}\,\mathcal T(\overline{\bf
A}_0)+\frac{{\sf bc}-{\sf a}\overline{\sf a}\sf e}{{\sf
a}^4\overline{\sf a}^2}\,\mathcal L(\overline{\bf
A}_0)+\frac{\overline{\sf b}{\sf c}-{\sf a}\overline{\sf a}\sf
d}{{\sf a}^3\overline{\sf a}^3}\,\overline{\mathcal L}(\overline{\bf
A}_0)\bigg),
\\
V^{\sf new}_3&:=V_3-\frac{{\bf A}_0\overline{\bf A}_0}{{\sf
a}^3\overline{\sf a}}\,\mathcal L({\bf C}'_0),
\\
V^{\sf new}_4&:=V_4-\frac{{\bf A}_0\overline{\bf A}_0}{{\sf
a}^2\overline{\sf a}^2}\,\overline{\mathcal L}({\bf C}'_0),
\\
W^{\sf new}_1&:=W_1+\frac{{\bf A}_0}{{\sf a}\overline{\sf
a}^2}\,\bigg(\frac{1}{{\sf a}^2\overline{\sf a}}\,{\mathcal
S}(\overline{\bf D}'_0)-\frac{{\sf c}}{{\sf a}^3\overline{\sf
a}^2}\,{\mathcal T}(\overline{\bf D}'_0)+\frac{{\sf b}{\sf c}-{\sf
a}\overline{\sf a}\sf e}{{\sf a}^4\overline{\sf a}^2}\,{\mathcal
L}(\overline{\bf D}'_0)+\frac{\overline{\sf b}{\sf c}-{\sf
a}\overline{\sf a}{\sf d}}{{\sf a}^3\overline{\sf
a}^3}\,\overline{\mathcal L}(\overline{\bf D}'_0)\bigg)-
\\
&-\frac{{\bf A}_0\overline{\sf c}}{{\sf a}^2\overline{\sf
a}^3}\,\bigg(\frac{1}{{\sf a}^2\overline{\sf a}}\,{\mathcal S}({\bf
B}'_0)-\frac{{\sf c}}{{\sf a}^3\overline{\sf a}^2}\,{\mathcal T}({\bf
B}'_0)+\frac{{\sf b}{\sf c}-{\sf a}\overline{\sf a}{\sf e}}{{\sf
a}^4\overline{\sf a}^2}\,{\mathcal L}({\bf B}'_0)+\frac{\overline{\sf
b}{\sf c}-{\sf a}\overline{\sf a}{\sf d}}{{\sf a}^3\overline{\sf
a}^3}\,\overline{\mathcal L}({\bf B}'_0)\bigg)+
\\
&+\frac{{\bf A}_0{\sf c}}{{\sf a}^3\overline{\sf
a}^2}\,\bigg(\frac{1}{{\sf a}\overline{\sf a}^2}\,\overline{\mathcal
S}({\bf B}'_0)-\frac{\overline{\sf c}}{{\sf a}^2\overline{\sf
a}^3}\,{\mathcal T}({\bf B}'_0)+\frac{{\sf b}\overline{\sf c}-{\sf
a}\overline{\sf a}\,\overline{\sf d}}{{\sf a}^3\overline{\sf
a}^3}\,{\mathcal L}({\bf B}'_0)+\frac{\overline{\sf b}\,\overline{\sf
c}-{\sf a}\overline{\sf a}\,\overline{\sf e}}{{\sf a}^2\overline{\sf
a}^4}\,\overline{\mathcal L}({\bf B}'_0)\bigg)+
\\
&+\frac{{\sf e}}{{\sf a}^3\overline{\sf a}}\,\bigg(\frac{1}{{\sf
a}\overline{\sf a}^2}\,\overline{\mathcal S}({\bf
A}_0)-\frac{\overline{\sf c}}{{\sf a}^2\overline{\sf a}^3}\,{\mathcal
T}({\bf A}_0)+\frac{{\sf b}\overline{\sf c}-{\sf a}\overline{\sf
a}\,\overline{\sf d}}{{\sf a}^3\overline{\sf a}^3}\,{\mathcal L}({\bf
A}_0)+\frac{\overline{\sf b}\,\overline{\sf c}-{\sf a}\overline{\sf
a}\,\overline{\sf e}}{{\sf a}^2\overline{\sf
a}^4}\,\overline{\mathcal L}({\bf A}_0)\bigg),
\\
W^{\sf new}_2&:=W_2+\frac{{\bf A}_0}{{\sf a}\overline{\sf
a}}\,\bigg(\frac{1}{{\sf a}^2\overline{\sf a}}\,\mathcal S({\bf
B}'_0)-\frac{\sf e}{{\sf a}^3\overline{\sf a}}\,\mathcal L({\bf
B}'_0)-\frac{\sf d}{{\sf a}^2\overline{\sf a}^2}\,\overline{\mathcal
L}({\bf B}'_0)\bigg)+
\\
&+\frac{\sf e}{{\sf a}^3\overline{\sf
a}}\,\bigg(\frac{1}{{\sf a}\overline{\sf a}}\,\mathcal T({\bf
A}_0)-\frac{\sf b}{{\sf a}^2\overline{\sf a}}\,\mathcal L({\bf
A}_0)-\frac{\overline{\sf b}}{{\sf a}\overline{\sf
a}^2}\,\overline{\mathcal L}({\bf A}_0)\bigg),
\\
W^{\sf new}_3&:=W_3+\frac{{\bf A}_0\sf c}{{\sf a}^4\overline{\sf
a}^2}\,\mathcal L({\bf B}'_0)+\frac{1}{{\sf a}}\,\bigg(\frac{1}{{\sf
a}^2\overline{\sf a}}\,\mathcal S({\bf A}_0)-\frac{\sf c}{{\sf
a}^3\overline{\sf a}^2}\,\mathcal T({\bf A}_0)+\frac{{\sf bc}}{{\sf
a}^4\overline{\sf a}^2}\,\mathcal L({\bf A}_0)+\frac{\overline{\sf
b}{\sf c}-{\sf a}\overline{\sf a}\sf d}{{\sf a}^3\overline{\sf
a}^3}\,\overline{\mathcal L}({\bf A}_0)\bigg),
\\
W^{\sf new}_4&:=W_4+\frac{{\bf A}_0\sf c}{{\sf a}^3\overline{\sf
a}^3}\,\overline{\mathcal L}({\bf B}'_0)+\frac{\sf e}{{\sf
a}^3\overline{\sf a}^2}\,\overline{\mathcal L}({\bf A}_0),
\\
W^{\sf new}_5&:=W_5+\frac{{\bf A}_0}{{\sf a}\overline{\sf
a}}\,\bigg(\frac{1}{{\sf a}\overline{\sf a}^2}\,\overline{\mathcal
S}({\bf B}'_0)-\frac{\overline{\sf c}}{{\sf a}^2\overline{\sf
a}^3}\,\mathcal T({\bf B}'_0)+\frac{{\sf b}\overline{\sf c}-{\sf
a}\overline{\sf a}\overline{\sf d}}{{\sf a}^3\overline{\sf
a}^3}\,\mathcal L({\bf B}'_0)+\frac{\overline{\sf b}\,\overline{\sf
c}-{\sf a}\overline{\sf a}\,\overline{\sf e}}{{\sf a}^2\overline{\sf
a}^4}\,\overline{\mathcal L}({\bf B}'_0)\bigg)-
\\
&-\frac{{\bf A}_0}{{\sf a}\overline{\sf a}^2}\,\bigg(\frac{1}{{\sf
a}\overline{\sf a}}\,\mathcal T(\overline{\bf D}'_0)-\frac{\sf
b}{{\sf a}^2\overline{\sf a}}\,\mathcal L(\overline{\bf
D}'_0)-\frac{\overline{\sf b}}{{\sf a}\overline{\sf
a}^2}\,\overline{\mathcal L}(\overline{\bf D}'_0)\bigg)+\frac{{\bf
A}_0\overline{\sf c}}{{\sf a}^2\overline{\sf
a}^3}\,\bigg(\frac{1}{{\sf a}\overline{\sf a}}\,\mathcal T({\bf
B}'_0)-\frac{\sf b}{{\sf a}^2\overline{\sf a}}\,\mathcal L({\bf
B}'_0)-\frac{\overline{\sf b}}{{\sf a}\overline{\sf
a}^2}\,\overline{\mathcal L}({\bf B}'_0)\bigg),
\\
W^{\sf new}_6&:=W_6+\frac{1}{{\sf a}}\,\bigg(\frac{1}{{\sf
a}\overline{\sf a}^2}\,\overline{\mathcal S}({\bf
A}_0)-\frac{\overline{\sf c}}{{\sf a}^2\overline{\sf a}^3}\,\mathcal
T({\bf A}_0)+\frac{{\sf b}\overline{\sf c}-{\sf a}\overline{\sf
a}\overline{\sf d}}{{\sf a}^3\overline{\sf a}^3}\,\mathcal L({\bf
A}_0)+\frac{\overline{\sf b}\,\overline{\sf c}-{\sf a}\overline{\sf
a}\,\overline{\sf e}}{{\sf a}^2\overline{\sf
a}^4}\,\overline{\mathcal L}({\bf A}_0)\bigg)-
\\
&-\frac{{\bf A}_0}{{\sf
a}^2\overline{\sf a}^2}\,\mathcal L(\overline{\bf D}'_0)+\frac{{\bf
A}_0\overline{\sf c}}{{\sf a}^3\overline{\sf a}^3}\,\mathcal L({\bf
B}'_0),
\endaligned
\]
\[
\footnotesize\aligned
W^{\sf new}_7&:=W_7-\frac{{\bf A}_0}{{\sf a}\overline{\sf
a}^3}\,\overline{\mathcal L}(\overline{\bf D}'_0)-\frac{\overline{\bf
D}'_0}{{\sf a}\overline{\sf a}^3}\,\overline{\mathcal L}({\bf
A}_0)+\frac{{\bf A}_0\overline{\sf c}}{{\sf a}^2\overline{\sf
a}^4}\,\overline{\mathcal L}({\bf B}'_0),
\\
W^{\sf new}_8&:=W_8-\frac{{\bf A}_0}{{\sf a}^2\overline{\sf
a}}\,\mathcal L({\bf B}'_0)+\frac{1}{\sf a}\,\bigg(\frac{1}{{\sf
a}\overline{\sf a}}\,\mathcal T({\bf A}_0)-\frac{\sf b}{{\sf
a}^2\overline{\sf a}}\,\mathcal L({\bf A}_0)-\frac{\overline{\sf
b}}{{\sf a}\overline{\sf a}^2}\,\overline{\mathcal L}({\bf
A}_0)\bigg),
\\
W^{\sf new}_9&:=W_9-\frac{{\bf A}_0}{{\sf a}\overline{\sf
a}^2}\,\overline{\mathcal L}({\bf B}'_0),
\\
W^{\sf new}_{10}&:=W_{10}-\frac{1}{{\sf a}\overline{\sf
a}}\,\overline{\mathcal L}({\bf A}_0).
\endaligned
\]
After computing the above preparatory equations, we are now ready to
apply the second normalization-absorbtion step. Replacing the single
Maurer-Cartan form $\gamma$ by:
\begin{equation}
\label{gamma} \gamma\longmapsto\gamma+
p\,\sigma+q\,\overline\sigma+r\,\rho+s\,\zeta+t\,\overline\zeta
\end{equation}
has no effect on the two first expressions $d\sigma$ and $d\rho$ of
\thetag{\ref{structure-firts-absorbtion-R-neq-0-NEW}}, but it changes
the third one as follows:
\[
\footnotesize\aligned d\zeta & = \gamma\wedge\sigma + \ \ \ \ \ \ \ \
\ \ \ \ \ \ \ \ \ \ \ \ \ \ \ \ \ \ \ \ \ \ \ \ \ \ \ \ \ \ \ \ \ \ \
\ \ \ \ \ \ \ \ \ \ \ \ \ \ \ \ \ \ \ \ \ \ \ \ \ \ \ \ \ \ \ \ \ \ \
\ \ \ \ \ \
 \ \ \ \ \ \ \ \ \ \
\\
& \ \ \ \ \ + \big(W^{\sf
new}_1-q\big)\,\sigma\wedge\overline{\sigma} + \big(W^{\sf
new}_2-r\big)\,\sigma\wedge\rho + \big(W^{\sf
new}_3-s\big)\,\sigma\wedge\zeta + \big(W^{\sf
new}_4-t\big)\,\sigma\wedge\overline{\zeta} +
\\
& \ \ \ \ \ \ \ \ \ \ \ \ \ \ \ \ \ \ \ \ \ \ \ \ + W^{\sf
new}_5\,\overline{\sigma}\wedge\rho + W^{\sf
new}_6\,\overline{\sigma}\wedge\zeta + W^{\sf
new}_7\,\overline{\sigma}\wedge\overline{\zeta} +
\\
& \ \ \ \ \ \ \ \ \ \ \ \ \ \ \ \ \ \ \ \ \ \ \ \ \ \ \ \ \ \ \ \ \ \
\ \ \ \ \ \ \ \ \ + W^{\sf new}_8\,\rho\wedge\zeta + W^{\sf
new}_9\,\rho\wedge\overline{\zeta} +
\\
& \ \ \ \ \ \ \ \ \ \ \ \ \ \ \ \ \ \ \ \ \ \ \ \ \ \ \ \ \ \ \ \ \ \
\ \ \ \ \ \ \ \ \ \ \ \ \ \ \ \ \ \ \ \ \ \ \ \ \ \ \ \ \ + W^{\sf
new}_{10}\,\zeta\wedge\overline{\zeta}.
\endaligned
\]
Visibly, one can annihilate all the coefficients appearing at the second
line of this expression by putting:
\[
p=0, \ \ \ q=W^{\sf new}_1, \ \ \ r=W^{\sf new}_2, \ \ \ s=W^{\sf
new}_3, \ \ \ t=W^{\sf new}_4.
\]
In other words, we can annihilate the second line of the above expression
by modifying the single Maurer-Cartan 1-form into the form:
\[
\gamma\mapsto \gamma-W^{\sf new}_1\,\overline\sigma-W^{\sf
new}_2\,\rho-W^{\sf new}_3\,\zeta-W^{\sf new}_4\,\overline\zeta.
\]
Therefore, we find 15 normalizable expressions:
\begin{equation}
\label{normalizable-second-loop-R0} \aligned \left[
\begin{array}{l}0\equiv U^{\sf new}_1=U^{\sf new}_2=U^{\sf new}_3=U^{\sf new}_4,
\\
0\equiv V^{\sf new}_1=V^{\sf new}_2=V^{\sf new}_3=V^{\sf
new}_4=V^{\sf new}_8,
\\
0\equiv W^{\sf new}_5=W^{\sf new}_6=W^{\sf new}_7=W^{\sf
new}_8=W^{\sf new}_9=W^{\sf new}_{10}.
\end{array}
\right.
\endaligned
\end{equation}
Many of these torsion coefficients include the as yet undetermined
parameter $\sf e$, but as in the branch $R= 0$, let us plainly employ
$V^{\sf new}_4$ to normalize $\sf e$. Of course, we quite similarly
have:
\[
V^{\sf new}_4=-i\,\frac{\sf e}{{\sf a}^2\overline{\sf
a}}+\frac{1}{3{\sf a}\overline{\sf a}}\Big(6\,B\overline
B+3\,A-\overline B\,\overline Q-3\,\overline{\mathcal L}(\overline
B)\Big).
\]
Normalizing this expression to $0$ determines the expression of the
last group parameter $\sf e$ as:
\[\boxed{\aligned
 {\sf e}=-\textstyle{\frac{i}{3}}\,{\bf A}_0\,\big(6\,B\overline
B-3\,\overline{\mathcal L}(\overline B)+3\,A-\overline B\,\overline
Q\big).
\endaligned}
\]

After determining the last group parameter ${\sf e}$ in such a way,
the remaining fourteen normalizable expressions change into the form:
\[
\!\!\!\!\!\!\!\!\!\!\!\!\!\!\!\!\!\!\!\!\!\!\!\!\!\!\!\!\!\!\!\!\!\!\!\!\!
\!\!\!\!\!\!\!\!\!\!\!\!\!\!\!\!\!\!
\footnotesize\aligned U^{\sf new}_1&= \textstyle\frac{i}{6\,{\bf
A}_0^2\overline{\bf A}^3} \big(12\,\overline{\bf
A}_0\overline{\mathcal L}(\overline{\mathcal L}(\mathcal L({\bf
A}_0)))-24\,\overline{\bf A}_0\overline{\mathcal L}(\mathcal
L(\overline{\mathcal L}({\bf A}_0)))+12\,\overline{\bf A}_0\mathcal
L(\overline{\mathcal L}(\overline{\mathcal L}({\bf A}_0)))+ 6\,{\bf
A}_0\overline{\mathcal L}(\overline{\mathcal L}(\mathcal
L(\overline{\bf A}_0)))-
\\
&-12\,{\bf A}_0\overline{\mathcal L}(\mathcal L(\overline{\mathcal
L}(\overline{\bf A}_0)))+6\,{\bf A}_0\mathcal L(\overline{\mathcal
L}(\overline{\mathcal L}(\overline{\bf A}_0)))-3\,{\bf
A}_0\overline{\bf A}_0\mathcal L(\mathcal L(\overline R))-2\,{\bf
A}_0\overline{\bf A}_0\overline{\mathcal L}(\overline{\mathcal L}(Q))
+ 6\,{\bf A}_0\overline{\bf A}_0\mathcal L(\overline{\mathcal L}(B))+
\\
&+12\,\overline{\bf A}_0B\mathcal L(\overline{\mathcal L}({\bf
A}_0))-12\overline{\bf A}_0B\overline{\mathcal L}(\mathcal L({\bf
A}_0))+ 6\,{\bf A}_0B\mathcal L(\overline{\mathcal L}(\overline{\bf
A}_0))-6\,{\bf A}_0B\overline{\mathcal L}(\mathcal L(\overline{\bf
A}_0))-8a\overline{\bf A}_0\mathcal L(\overline P)+6\,{\bf
A}_0\overline{\bf A}_0\overline{\mathcal L}(A)-
\\
&-3\,{\bf A}_0\overline{\bf A}_0R\overline{\mathcal L}(\overline
R)+2\,{\bf A}_0\overline{\bf A}_0B\overline{\mathcal L}(Q)-3\,{\bf
A}_0\overline{\bf A}_0\overline R\mathcal L(Q)-4\,{\bf
A}_0\overline{\bf A}_0\overline Q\mathcal L(B)+ 3\,{\bf
A}_0\overline{\bf A}_0\overline B\mathcal L(\overline R)-2\,{\bf
A}_0\overline{\bf A}_0Q\overline{\mathcal L}(B)-
\\
&-8\,{\bf A}_0\overline{\bf A}_0B\mathcal L(\overline Q)+6\,{\bf
A}_0\overline{\bf A}_0B\overline{\mathcal L}(\overline B)+ 12\,{\bf
A}_0\overline{\bf A}_0B\mathcal L(B)-12\overline{\bf A}_0B\overline
Q\mathcal L({\bf A})_0-16\overline{\bf A}_0Q\overline R\mathcal
L({\bf A}_0)+24\,\overline{\bf A}_0\overline B\overline R\mathcal
L({\bf A}_0)+
\\
&+12\,B\overline B\,\overline{\bf A}_0\overline{\mathcal L}({\bf
A}_0)-6\,{\bf A}_0B\overline Q\mathcal L(\overline{\bf A}_0)-8\,{\bf
A}_0Q\overline R\mathcal L(\overline{\bf A}_0)+12\,{\bf A}_0\overline
B\,\overline R\mathcal L(\overline{\bf A}_0)+6\,{\bf A}_0B\overline
B\overline{\mathcal L}(\overline{\bf A}_0)+12\,\overline{\bf
A}_0B^2\mathcal L({\bf A}_0)-
\\
&-12\,\overline{\bf A}_0\mathcal L({\bf A}_0)\mathcal L(\overline
R)+12\,\overline{\bf A}_0\mathcal L({\bf A}_0)\overline{\mathcal
L}(B)-12\,\overline{\bf A}_0\overline P\mathcal L({\bf
A}_0)-12\,\overline{\bf A}_0\overline{\mathcal L}({\bf A}_0)\mathcal
L(B)+12\,\overline{\bf A}_0A\overline{\mathcal L}({\bf A}_0)+6\,{\bf
A}_0B^2\mathcal L(\overline{\bf A}_0)-
\\
&-6\,{\bf A}_0\mathcal L(\overline{\bf A}_0)\mathcal L(\overline
R)+6\,{\bf A}_0\mathcal L(\overline{\bf A}_0)\overline{\mathcal
L}(B)-6\,{\bf A}_0\overline P\mathcal L(\overline{\bf A}_0)- 6\,{\bf
A}_0\overline{\mathcal L}(\overline{\bf A}_0)\mathcal L(B)+6\,{\bf
A}_0A\overline{\mathcal L}(\overline{\bf A}_0) +2\,{\bf
A}_0\overline{\bf A}_0B\overline B\overline Q-
\\
&-4\,{\bf A}_0\overline{\bf A}_0A\overline Q+2\,{\bf
A}_0\overline{\bf A}_0Q\overline P+5\,{\bf A}_0\overline{\bf
A}_0Q^2\overline R-2\,{\bf A}_0\overline{\bf A}_0B^2Q+3\,{\bf
A}_0\overline{\bf A}_0\overline B^2\overline R+6\,{\bf
A}_0\overline{\bf A}_0\overline B\overline P-6\,{\bf
A}_0\overline{\bf A}_0B^2\overline B-
\\
&-9\,{\bf A}_0\overline{\bf A}_0BR\overline R+5\,{\bf
A}_0\overline{\bf A}_0R\overline R\overline Q+2\,{\bf
A}_0\overline{\bf A}_0BQ\overline Q-8\,{\bf A}_0\overline{\bf
A}_0Q\overline B\,\overline R \big),
\\
 U^{\sf new}_2&=\textstyle\frac{i}{3\,{\bf
A}_0^2\overline{\bf A}_0^2}\,\big(-6\,\overline{\bf A}_0\mathcal
L(\overline{\mathcal L}({\bf A}_0))+6\,\overline{\bf
A}_0\overline{\mathcal L}(\mathcal L({\bf A}_0))-3\,{\bf A}_0\mathcal
L(\overline{\mathcal L}(\overline{\bf A}_0))+ 3\,{\bf
A}_0\overline{\mathcal L}(\mathcal L(\overline{\bf A}_0))-3\,{\bf
A}_0\overline{\bf A}_0\overline{\mathcal L}(Q)-
\\
&-3\,{\bf A}_0\overline{\bf A}_0R\overline R+3\,{\bf
A}_0\overline{\bf A}_0\mathcal L(B)+ 3\,{\bf A}_0\overline{\bf
A}_0\overline{\mathcal L}(\overline B)- 6\,\overline{\bf
A}_0B\mathcal L({\bf A}_0)+2\,\overline{\bf A}_0\overline Q\mathcal
L({\bf A}_0)+6\,\overline{\bf A}_0\overline B\,\overline{\mathcal
L}({\bf A}_0)-2\,\overline{\bf A}_0Q\overline{\mathcal L}({\bf A}_0)-
\\
&-3\,{\bf A}_0B\mathcal L(\overline{\bf A}_0)+{\bf A}_0\overline
Q\mathcal L(\overline{\bf A}_0)+3\,{\bf A}_0\overline
B\overline{\mathcal L}(\overline{\bf A}_0)-{\bf
A}_0Q\overline{\mathcal L}(\overline{\bf A}_0) +4\,{\bf
A}_0\overline{\bf A}_0BQ-{\bf A}_0\overline{\bf A}_0Q\overline
Q-6\,{\bf A}_0\overline{\bf A}_0B\overline B+{\bf A}_0\overline{\bf
A}_0\,\overline B\,\overline Q \big),
\\
 U^{\sf new}_3&=2\,\overline U^{\sf
new}_4-3\,\overline W^{\sf new}_{10},
\\
U^{\sf new}_4&=\textstyle\frac{1}{{\bf A}_0\overline{\bf
A}_0^2}\,\big(-2\,\overline{\bf A}_0\overline{\mathcal L}({\bf
A}_0)-{\bf A}_0\overline{\mathcal L}(\overline{\bf A}_0)+{\bf
A}_0\overline{\bf A}_0B\big),
\\
V^{\sf new}_1&=-\textstyle\frac{i}{18\,{\bf
A}_0^2\overline{\bf A}_0^2}\big(-18\,\overline{\mathcal L}(\mathcal
L(\overline{\mathcal L}(B)))+18\,\mathcal L(\mathcal
L(\overline{\mathcal L}(\overline Q)))- 36\,\mathcal
L(\overline{\mathcal L}(\mathcal L(\overline
Q)))+12\,\overline{\mathcal L}(\mathcal L(\mathcal L(\overline
Q)))+18\,\overline{\mathcal L}(\overline{\mathcal L}(\mathcal
L(\overline B)))+
\\
&+6\,\overline{\mathcal L}(\overline{\mathcal L}(\overline{\mathcal
L}(Q)))-18\,\overline{\mathcal L}(\overline{\mathcal
L}(P))-6\,\mathcal L(\mathcal L(\overline P))-12\,B\overline{\mathcal
L}(\mathcal L(\overline Q))+9R\overline{\mathcal L}(\overline
P)+18\,\mathcal L(\overline{\mathcal L}(A))+18\,\mathcal
L(\overline{\mathcal L}(\overline P))-
\\
&-12\,\overline{\mathcal L}(\mathcal L(\overline P))+30\,\overline
Q\mathcal L(\overline{\mathcal L}(B))+18\,B\mathcal
L(\overline{\mathcal L}(\overline Q))+27\,B^3R+9\,\overline B\mathcal
L(\mathcal L(\overline R))-24\,\overline B\,\overline{\mathcal
L}(\overline{\mathcal L}(Q))+ 18\,\overline B\mathcal
L(\overline{\mathcal L}(B))+
\\
&+6B\mathcal L(\mathcal L(\overline Q))+9\,B\overline{\mathcal
L}(\overline{\mathcal L}(R))-36\,B\overline{\mathcal L}(\mathcal
L(\overline B))-2\,\overline Q\overline{\mathcal
L}(\overline{\mathcal L}(Q))-36\,\overline B\mathcal
L(\overline{\mathcal L}(\overline Q))-24\,Q\overline{\mathcal
L}(\mathcal L(\overline Q))+
\\
&+12\,Q\mathcal L(\overline{\mathcal L}(\overline Q))+24\,B^2\mathcal
L(\overline Q)-36\,B\mathcal L(\overline{\mathcal
L}(B))+60\,\overline B\overline{\mathcal L}(\mathcal L(\overline
Q))+12\,\overline Q\,\overline{\mathcal L}(\mathcal L(\overline
B))-27\,\overline BP\overline R+24\,\overline{\mathcal L}(B)\mathcal
L(\overline Q)+
\\
&+18\,A\overline{\mathcal L}(B)-6\,B\overline B\,\overline{\mathcal
L}(\overline Q)+6\,B\mathcal L(\overline P)+9\overline R\mathcal
L(P)+18\,\overline P\,\overline{\mathcal L}(\overline
B)+18\,AR\overline R+ 6\,A\mathcal L(\overline
Q)+42\,B\overline{\mathcal L}(P)+18\,Q\overline{\mathcal L}(A)-
\\
&-6\,\overline Q\mathcal L(A)-36\,\overline B\overline{\mathcal
L}(A)+36\,\overline B\mathcal L(\overline P)- 18\,\overline
B^2\,\overline P-18\,B^2P+36\,P\overline{\mathcal L}(B)+18\,\overline
P\,\overline{\mathcal L}(R)+18\,P\mathcal L(\overline
R)-18\,A\mathcal L(B)-
\\
&-18\,\overline P\mathcal L(\overline Q)-14\,\overline Q\mathcal
L(\overline P)-54\,\overline B^2\,\overline{\mathcal
L}(B)-6\,\overline Q\,\overline{\mathcal L}(A)+18\,B^2\mathcal
L(\overline B)+ 18\,B\overline{\mathcal L}(A)-27\,\overline
B^2\mathcal L(\overline R)-27\,B^2\overline{\mathcal L}(R)-
\\
&-12\,B\overline Q\mathcal L(\overline B)+48\,\overline
BQ\overline{\mathcal L}(B)+ 12\,B\overline B\mathcal L(\overline
Q)-90\,B\overline B\overline{\mathcal L}(\overline B)+48\,B\overline
B\,\overline{\mathcal L}(Q)-36\,B\overline B\mathcal
L(B)+9\,\overline BR\overline{\mathcal L}(\overline R)+ 9\,\overline
B\,\overline R\mathcal L(Q)+
\\
&+9\,BR\overline{\mathcal L}(\overline Q)+48\,BQ\overline{\mathcal
L}(\overline B)+9\,B\overline R\mathcal L(R)+ 36\,\overline{\mathcal
L}(\overline B)\mathcal L(B)+24\,Q\overline R\overline{\mathcal
L}(R)+18\,B\overline Q\,\overline{\mathcal L}(R)-36\,\overline
B\,\overline R\,\overline{\mathcal L}(R)+
\\
&+24\,\overline QR\mathcal L(\overline R)+ 18\,\overline BQ\mathcal
L(\overline R)-36\,BR\mathcal L(\overline R)+18\,\overline{\mathcal
L}(R)\mathcal L(\overline R)-18\,\mathcal L(\overline
B)\overline{\mathcal L}(B)-18\,\overline{\mathcal
L}(B)\overline{\mathcal L}(\overline B)+12B\overline
Q\overline{\mathcal L}(\overline B)-
\\
&-32\,B\overline Q\mathcal L(\overline Q)-18\,\overline Q\,\overline
B\,\overline{\mathcal L}(B)+4\overline QQ\overline{\mathcal
L}(B)-24\,B^2\overline B\,\overline Q+ 36\,B\overline
B\,\overline{\mathcal L}(B)+60\,B\overline Q\mathcal
L(B)-36\,\mathcal L(B)\overline{\mathcal L}(B)+18\,\overline
P\mathcal L(B)-
\\
&-36\,B^2\mathcal L(B)-16\,\overline Q^2\mathcal L(B)+12\,\mathcal
L(B)\overline{\mathcal L}(\overline Q)-6\,A\overline{\mathcal
L}(\overline Q)-12\,Q\mathcal L(\overline P)- 6\,\overline
P\overline{\mathcal L}(Q)-4B\overline Q\,\overline{\mathcal
L}(Q)+24\,\overline Q\,\overline B\,\overline{\mathcal L}(\overline
B)-
\\
&-12\,Q\overline Q\,\overline{\mathcal L}(\overline B)- 12\,\overline
B\,\overline Q\,\overline{\mathcal L}(Q)- 12\,\overline
BQ\overline{\mathcal L}(\overline Q)-12\,P\overline{\mathcal
L}(\overline Q)-12\,\overline Q\,\overline{\mathcal
L}(P)-12\,\overline{\mathcal L}(\overline B)\mathcal L(\overline Q)+
12\,\overline{\mathcal L}(\overline Q)\mathcal L(\overline B)+
\\
&+12\,\overline B^2\overline{\mathcal L}(\overline
Q)-36\,\overline{\mathcal L}(\overline B)\overline{\mathcal
L}(Q)+36\,\overline{\mathcal L}(\overline B)^2 +9\,\overline
P\,\overline Q R+ 9\,PQ\overline R+6\,BP\overline Q-27\,BR\overline
P+6\,\overline B\,\overline PQ-36\,AB\overline B+27\,\overline
B^3\,\overline R+
\\
&+2\,\overline Q^2A-36\,B^2R\overline Q-12\,\overline
BB^2Q-36\,\overline B^2Q\overline R+ 6\,\overline B^2B\overline
Q+9\,B\overline Q^2R+12\,ABQ+9\,\overline BQ^2\overline
R+30\,\overline BA\overline Q- 63\,\overline B\overline R\,\overline
QR+
\endaligned
\]
\[
\!\!\!\!\!\!\!\!\!\!\!\!\!\!\!\!\!\!\!\!\!\!\!\!\!\!\!\!\!\!\!\!\!\!\!\!\!
\!\!\!\!\!\!\!\!\!\!\!\!\!\!\!\!\!\!
\footnotesize\aligned
&+126\,\overline B\,\overline RBR-63\,BRQ\overline R+32\,\overline
QRQ\overline R+18\,B^3\overline B+ 18AB^2+12\,\overline Q\,\overline
B\,\overline P+ 14\,\overline BB\overline Q^2-4\,\overline
QQ\overline P-
\\
&-4\,QB\overline Q^2-12\,AB\overline Q+4\,B^2Q\overline
Q-12\,QA\overline Q\big),
\endaligned
\]
\[
\!\!\!\!\!\!\!\!\!\!\!\!\!\!\!\!\!\!\!\!\!\!\!\!\!\!\!\!\!\!\!\!\!\!\!\!\!
\!\!\!\!\!\!\!\!\!\!\!\!\!\!\!\!\!\!
\footnotesize\aligned
V^{\sf new}_2&=\textstyle\frac{i}{9\,{\bf A}_0^3\overline{\bf
A}_0^2}\,\big(- 9\,\overline{\bf A}_0\,\overline{\mathcal L}(\mathcal
L(\mathcal L({\bf A}_0)))-9\,{\bf A}_0\mathcal L(\mathcal
L(\overline{\mathcal L}(\overline{\bf A}_0)))+18\,{\bf A}_0\mathcal
L(\overline{\mathcal L}(\mathcal L(\overline{\bf A}_0)))-9\,{\bf
A}_0\overline{\mathcal L}(\mathcal L(\mathcal L(\overline{\bf
A}_0)))+
\\
&+ 9\,\overline{\bf A}_0\,\overline B\mathcal L(\overline{\mathcal
L}({\bf A}_0))-9\,\overline{\bf A}_0\,\overline B\,\overline{\mathcal
L}(\mathcal L({\bf A}_0))+9\,{\bf A}_0\overline B\mathcal
L(\overline{\mathcal L}(\overline{\bf A}_0))-9\,\overline{\bf
A}_0\mathcal L(\mathcal L(\overline{\mathcal L}({\bf A}_0))- 9\,{\bf
A}_0\overline B\,\overline{\mathcal L}(\mathcal L(\overline{\bf
A}_0))+
\\
&+3\,{\bf A}_0\overline{\bf A}_0\mathcal L(\mathcal L(\overline
Q))+18\,\overline{\bf A}_0\mathcal L(\overline{\mathcal L}(\mathcal
L({\bf A}_0))+9\,{\bf A}_0P\overline{\mathcal L}(\overline{\bf
A}_0)+9\,\overline{\bf A}_0\mathcal L({\bf A}_0)\overline{\mathcal
L}(\overline B)+9{\bf A}_0^2\overline B\,\overline{\mathcal
L}(\overline B)-9\,{\bf A}_0\overline{\bf A}_0\mathcal L(A)-
\\
&-9\,\overline{\bf A}_0\,\overline{\mathcal L}({\bf A}_0)\mathcal
L(\overline B)+3\,{\bf A}_0\overline{\bf A}_0\,\overline{\mathcal
L}(P)-9\,{\bf A}_0A\mathcal L(\overline{\bf A}_0)+ 9\,\overline{\bf
A}_0\,\overline{\mathcal L}({\bf A}_0)\overline{\mathcal
L}(R)-9\,\overline{\bf A}_0\,\overline B^2\,\overline{\mathcal
L}({\bf A}_0)- 3\,{\bf A}_0^2Q\overline{\mathcal L}(\overline B)-
\\
&-9\,{\bf A}_0\overline{\mathcal L}(\overline{\bf A})_0\mathcal
L(\overline B)-9\,{\bf A}_0\overline{\mathcal L}(\overline{\bf
A}_0)\overline B^2+9\,{\bf A}_0\overline{\mathcal L}(\overline{\bf
A}_0)\overline{\mathcal L}(R)+9\,P\overline{\bf
A}_0\,\overline{\mathcal L}({\bf A}_0)+9\,{\bf A}_0\mathcal
L(\overline{\bf A}_0)\overline{\mathcal L}(\overline
B)-9\,\overline{\bf A}_0A\mathcal L({\bf A}_0)-
\\
&-18\,{\bf A}_0BR\overline{\mathcal L}(\overline{\bf A}_0)+12\,{\bf
A}_0\overline QR\overline{\mathcal L}(\overline{\bf A}_0)+3\,{\bf
A}_0\overline{\bf A}_0\,\overline B\,\overline{\mathcal L}(Q)-9\,{\bf
A}_0\overline{\bf A}_0\,\overline B\mathcal L(B)- 9\,{\bf
A}_0\overline B\,\overline{\bf A}_0\,\overline{\mathcal L}(\overline
B)-4\,{\bf A}_0\overline{\bf A}_0\,\overline Q^2R-
\\
&-3\,{\bf A}_0\overline{\bf A}_0\,\overline Q\,\overline{\mathcal
L}(R)+ 9\,{\bf A}_0\overline{\bf A}_0B\overline{\mathcal
L}(R)-9\,B\overline B\,\overline{\bf A}_0\mathcal L({\bf
A}_0)+9\,\overline BQ\overline{\bf A}_0\,\overline{\mathcal L}({\bf
A}_0)+12\,\overline QR\overline{\bf A}_0\,\overline{\mathcal L}({\bf
A}_0)- 18\,BR\overline{\bf A}_0\,\overline{\mathcal L}({\bf A}_0)-
\\
&-9\,B\overline B{\bf A}_0\mathcal L(\overline{\bf A}_0)+9\,\overline
BQ{\bf A}_0\overline{\mathcal L}(\overline{\bf A}_0)+3\,{\bf
A}_0\overline{\bf A}_0\,\overline Q\mathcal L(\overline B)-9\,{\bf
A}_0\overline{\bf A}_0B\mathcal L(\overline B)-3\,{\bf
A}_0\overline{\bf A}_0\,\overline B\mathcal L(\overline Q)+3\,{\bf
A}_0^2AQ+3\,{\bf A}_0^2\overline B^2\overline Q-
\\
&-9{\bf A}_0^2A\overline B-18\,{\bf A}_0^2B\overline B^2+18\,{\bf
A}_0\overline B^2\overline{\bf A}_0B-{\bf A}_0^2\overline BQ\overline
Q+ 6\,{\bf A}_0^2\overline BBQ+3\,{\bf A}_0\overline{\bf
A}_0AQ+9\,{\bf A}_0\overline{\bf A}_0\,\overline PR-18\,{\bf
A}_0\overline{\bf A}_0B^2R+
\\
&+9\,{\bf A}_0\overline{\bf A}_0\overline BR\overline R- 2\,{\bf
A}_0\overline{\bf A}_0\,\overline BQ\overline Q+18\,{\bf
A}_0\overline{\bf A}_0B\overline QR- 3\,{\bf A}_0\overline{\bf
A}_0\,\overline QP+9\,{\bf A}_0\overline{\bf A}_0A\overline B \big),
\\
V^{\sf new}_3&=\textstyle\frac{1}{3\,{\bf
A}_0^2}\,\big(-3\,\overline{\mathcal L}(R)+9\,B\,R-4\,\overline
QR\big),
\\
V^{\sf new}_8&=\overline U^{\sf new}_4-\overline W^{\sf new}_{10},
\\
W^{\sf new}_5&=-\textstyle\frac{1}{27\,{\bf
A}_0\overline{\bf A}_0^3}\,\big(27\,\mathcal L(\overline{\mathcal
L}(\mathcal L(\overline R)))-27\,\overline{\mathcal L}(\mathcal
L(\mathcal L(\overline R))) +18\,\overline{\mathcal L}(\mathcal
L(\overline{\mathcal L}(\overline Q)))- 9\,\overline{\mathcal
L}(\overline{\mathcal L}(\mathcal L(\overline Q)))-9\,\mathcal
L(\overline{\mathcal L}(\overline{\mathcal L}(\overline Q)))+
\\
&+9\,\overline{\mathcal L}(\overline{\mathcal L}(\overline{\mathcal
L}(Q)))+ 27\,\mathcal L(\overline{\mathcal L}(\overline
P))-45\,\overline{\mathcal L}(\mathcal L(\overline P))+18\,B\mathcal
L(\overline{\mathcal L}(\overline Q))-36\,B\overline{\mathcal
L}(\mathcal L(\overline Q))-36\,\overline R\overline{\mathcal
L}(\mathcal L(Q))-
\endaligned
\]
\[
\!\!\!\!\!\!\!\!\!\!\!\!\!\!\!\!\!\!\!\!\!\!\!\!\!\!\!\!\!\!\!\!\!\!\!\!\!
\!\!\!\!\!\!\!\!\!\!\!\!\!\!
\footnotesize\aligned
&-27\,Q\overline{\mathcal L}(\mathcal L(\overline R))+27\,\overline
B\,\overline{\mathcal L}(\mathcal L(\overline R))+36\,\overline
R\mathcal L(\overline{\mathcal L}(Q))+36\,Q\mathcal
L(\overline{\mathcal L}(\overline R))-54\,\overline B\mathcal
L(\overline{\mathcal L}(\overline R))+ 24\,BQ^2\overline
R-180\,B\overline BQ\overline R-
\\
&-9\,\overline Q\,\overline{\mathcal L}(\overline{\mathcal
L}(Q))+27\,B\mathcal L(\mathcal L(\overline R))-9\,\overline
Q\mathcal L(\mathcal L(\overline R))-18\,\overline R\mathcal
L(\mathcal L(\overline Q))54\,A\overline{\mathcal
L}(B)-18\,B\overline B\,\overline{\mathcal L}(\overline
Q)-18\,B\mathcal L(\overline P)-
\\
&-12\,\overline Q\,\overline R\mathcal L(Q)+54\,\overline
P\,\overline{\mathcal L}(\overline B)-18\,\overline P\mathcal
L(\overline Q)+18\,\overline Q\mathcal L(\overline P)-18\,\overline
Q\,\overline{\mathcal L}(A)+54\,B\overline{\mathcal
L}(A)+36\,B\overline R\mathcal L(Q)+12\,\overline
RQ\overline{\mathcal L}(Q)-
\\
&-18\,\mathcal L(\overline R)\mathcal L(\overline Q)-54\,\overline
BQ\overline{\mathcal L}(\overline R)+27\,\mathcal L(\overline
R)\overline{\mathcal L}(\overline B)-54\,\overline{\mathcal
L}(B)\overline{\mathcal L}(\overline B)+36\,B\overline
Q\,\overline{\mathcal L}(\overline B)-27\,\mathcal L(\overline
R)R\overline R-24\,Q\overline R\mathcal L(\overline Q)+
\\
&+36\,Q\overline R\,\overline{\mathcal L}(\overline B)+27\,R\overline
R\overline{\mathcal L}(B)+54\,\overline B\,\overline R\mathcal
L(\overline Q)- 54\,\overline B\,\overline R\,\overline{\mathcal
L}(\overline B)+63\,\overline Q\,\overline B\mathcal L(\overline
R)-18\,\overline Q\,\overline B\,\overline{\mathcal
L}(B)-24\,\overline QQ\mathcal L(\overline R)-
\\
&-135\,B\overline B\mathcal L(\overline R)+108\,B\overline
B\,\overline{\mathcal L}(B)-72\,AB\overline Q+54\,BQ\mathcal
L(\overline R)- 18\,A\overline{\mathcal L}(\overline
Q)+12\,Q^2\overline{\mathcal L}(\overline R)+54\,\overline
B^2\,\overline{\mathcal L}(\overline R)+ 36\,\overline
R\,\overline{\mathcal L}(P)-
\\
&-9\,\overline P\overline{\mathcal L}(Q)+ 18\,\overline
Q^2A+54\,B^3\overline B+ 54\,AB^2+54\,\overline B\,\overline R^2R-
36\,Q\overline R^2R-27\,\overline PR\overline R+18\,\overline
BB\overline Q^2- 108\,\overline Q\,\overline B^2\,\overline R-
\\
&-16\,\overline QQ^2\overline R-81\,BP\overline R+ 27\,\overline
QP\overline R-18\,B\overline QR\overline R+102\,\overline
Q\,\overline BQ\overline R-72\,B^2\overline B\,\overline
Q+162\,B\overline B^2\overline R+18\,Q\overline RA-54\,\overline
B\,\overline RA \big),
\endaligned
\]
\[
\!\!\!\!\!\!\!\!\!\!\!\!\!\!\!\!\!\!\!\!\!\!\!\!\!\!\!\!\!\!\!\!\!\!\!\!\!
\!\!\!\!\!\!\!\!\!\!\!\!\!\!
\footnotesize\aligned W^{\sf new}_6&=-\textstyle\frac{i}{3\,{\bf
A}_0^2\overline{\bf A}_0^2}\,\big(3\,\mathcal L(\overline{\mathcal
L}(\overline{\mathcal L}({\bf A}_0)))+3\,\overline{\mathcal
L}(\overline{\mathcal L}(\mathcal L({\bf A}_0)))-
6\,\overline{\mathcal L}(\mathcal L(\overline{\mathcal L}({\bf
A}_0)))+3\,B\mathcal L(\overline{\mathcal L}({\bf A}_0))-
3\,B\overline{\mathcal L}(\mathcal L({\bf A}_0))+
\\
&+3\,{\bf A}_0\mathcal L(\overline{\mathcal L}(B))-3\,{\bf
A}_0\mathcal L(\mathcal L(\overline R))-4\,{\bf A}_0Q\mathcal
L(\overline R)-3\,{\bf A}_0\overline B\,\overline{\mathcal
L}(B)+3\,B\overline B\,\overline{\mathcal L}({\bf A}_0)+9\,{\bf
A}_0\overline B\mathcal L(\overline R)+
\\
&+6\,\overline B\,\overline R\mathcal L({\bf A}_0)-3\,B\overline
Q\mathcal L({\bf A}_0)-4\,Q\overline R\mathcal L({\bf
A}_0)+3\,B^2\mathcal L({\bf A}_0)-3\,\mathcal L({\bf A}_0)\mathcal
L(\overline R)+3\,\mathcal L({\bf A}_0)\overline{\mathcal
L}(B)-3\,\overline P\mathcal L({\bf A}_0)-
\\
&-3\,\overline{\mathcal L}({\bf A}_0)\mathcal
L(B)+3\,A\overline{\mathcal L}({\bf A}_0)+ 6\,{\bf A}_0\overline
R\mathcal L(\overline B)-3\,{\bf A}_0\overline Q\mathcal L(B)-4\,{\bf
A}_0\overline R\mathcal L(Q)+6\,{\bf A}_0B\mathcal L(B)-
\\
&-3\,{\bf A}_0B\mathcal L(\overline Q)-3\,{\bf A}_0\mathcal
L(\overline P)+3\,{\bf A}_0\overline B\,\overline P-3\,{\bf
A}_0B^2\overline B- 6\,{\bf A}_0\overline B^2\,\overline R+4\,{\bf
A}_0Q\overline B\,\overline R+3\,{\bf A}_0\overline BB\overline Q
\big),
\\
W^{\sf new}_7&=-\textstyle\frac{i}{9\,{\bf
A}_0\overline{\bf A}_0^3}\,\big(- 9\,{\bf A}_0\overline{\mathcal
L}(\mathcal L(\overline R))-12\,Q\overline R\overline{\mathcal
L}({\bf A}_0)+18\,{\bf A}_0\overline B\,\overline{\mathcal
L}(\overline R)-12\,\overline R{\bf A}_0\overline{\mathcal
L}(Q)-6\,B\overline Q\overline{\mathcal L}({\bf A}_0)-
\\
&-12\,{\bf A}_0Q\overline{\mathcal L}(\overline R)+ 9\,\overline
R{\bf A}_0\overline{\mathcal L}(\overline B)-18\,{\bf A}_0B\mathcal
L(\overline R)+12\,{\bf A}_0\overline Q\mathcal L(\overline
R)+18\,\overline B\,\overline R\,\overline{\mathcal L}({\bf
A}_0)-9\,\overline{\mathcal L}({\bf A}_0)\mathcal L(\overline R)+
\\
&+9\,\overline{\mathcal L}({\bf A}_0)\overline{\mathcal
L}(B)-9\,\overline P\overline{\mathcal L}({\bf A}_0) + 54\,{\bf
A}_0B\overline B\,\overline R-24\,{\bf A}_0BQ\overline R-27\,{\bf
A}_0\overline Q\,\overline B\,\overline R+16\,{\bf A}_0\overline
QQ\overline R+9\,{\bf A}_0A\overline R \big),
\endaligned
\]
\[
\footnotesize\aligned
 W^{\sf new}_8&=-\textstyle\frac{i}{9\,{\bf
A}_2\overline{\bf A}_0}\,\big(- 9\,\mathcal L(\overline{\mathcal
L}({\bf A}_0))+9\,\overline{\mathcal L}(\mathcal L({\bf
A}_0))-9\,{\bf A}_0\overline{\mathcal L}(\overline B)-9\,{\bf
A}_0\mathcal L(B)+3\,{\bf A}_0\mathcal L(\overline Q)- \ \ \ \ \ \ \
\ \ \ \ \ \ \ \ \ \ \ \ \ \ \ \ \ \ \ \ \ \ \ \ \ \ \ \ \ \ \ \ \ \ \
\ \ \ \ \ \ \
\\
&-9\,B\mathcal L({\bf A}_0)+3\,\overline Q\mathcal L({\bf
A}_0)+9\,\overline B\overline{\mathcal L}({\bf
A}_0)-3\,Q\overline{\mathcal L}({\bf A}_0)+9\,{\bf A}_0B\overline
B+3\,{\bf A}_0BQ-{\bf A}_0Q\overline Q+9\,{\bf A}_0A \big),
\\
 W^{\sf
new}_9&=\textstyle\frac{i}{9\,\overline{\bf
A}_0^2}\,\big(18\,\overline{\mathcal L}(B)-9\,\mathcal L(\overline
R)-3\overline{\mathcal L}(\overline Q)-12\,Q\overline R-9\,\overline
P-12\,B\overline Q+18\,\overline B\,\overline R+9\,B^2+\overline Q^2
\big), \ \ \ \ \ \ \ \ \ \ \ \ \ \ \ \ \ \ \ \ \ \ \ \ \ \ \ \ \ \ \
\ \ \ \ \ \ \ \ \ \ \ \ \ \ \ \ \ \ \  \
\\
W^{\sf new}_{10}&=\textstyle\frac{1}{3\,{\bf A}_0\overline{\bf
A}_0}\,\big(3\,{\bf A}_0B-{\bf A}_0\overline Q-3\,\overline{\mathcal
L}({\bf A}_0)\big).
\endaligned
\]

After determining the last group parameter $\sf e$ and also modifying
the single Maurer-Cartan form as \thetag{\ref{gamma}}, now the last
structure equations
\thetag{\ref{structure-firts-absorbtion-R-neq-0-NEW}} are converted
into the {\it final form}:

\begin{equation}
\label{structure-firts-absorbtion-R-neq-0-Final} \aligned d\sigma & =
U^{\sf new}_1\,\sigma\wedge\overline{\sigma} + U^{\sf
new}_2\,\sigma\wedge\rho + \big(2\,\overline U^{\sf
new}_4-3\,\overline W^{\sf new}_{10}\big)\,\sigma\wedge\zeta + U^{\sf
new}_4\,\sigma\wedge\overline{\zeta}+\overline\sigma\wedge\overline\zeta
+ {\rho\wedge\zeta},
\endaligned
\end{equation}
\[
\aligned d\rho & =V^{\sf new}_1\,\sigma\wedge\overline{\sigma} +
V^{\sf new}_2\,\sigma\wedge\rho + V^{\sf new}_3\,\sigma\wedge\zeta +
 \overline V^{\sf new}_2\,\overline{\sigma}\wedge\rho + \overline
V^{\sf new}_3\,\overline{\sigma}\wedge\overline{\zeta} +
\\
&\ \ \ + \big(\overline U^{\sf new}_4-\overline W^{\sf
new}_{10}\big)\,\rho\wedge\zeta + \big(U^{\sf new}_4-W^{\sf
new}_{10}\big)\,\rho\wedge\overline{\zeta} +
i\,\zeta\wedge\overline{\zeta},\ \ \ \ \ \ \ \ \ \ \ \ \ \ \ \ \ \ \
\ \ \ \ \ \ \ \ \ \ \ \ \ \ \ \ \ \ \
\endaligned
\]
\[
\aligned d\zeta & = \gamma\wedge\sigma + \ \ \ \ \ \ \ \ \ \ \ \ \ \
\ \ \ \ \ \ \ \ \ \ \ \ \ \ \ \ \ \ \ \ \ \ \ \ \ \ \ \ \ \ \ \ \ \ \
\ \ \ \ \ \ \ \ \ \ \ \ \ \ \ \ \ \ \ \ \ \ \ \ \ \ \ \ \ \ \ \ \ \ \
 \ \ \ \ \ \ \ \ \ \  \ \ \ \ \ \ \ \ \ \ \ \ \ \ \
\\
& \ \ \ \ \ + W^{\sf new}_1\,\sigma\wedge\overline{\sigma} + W^{\sf
new}_2\,\sigma\wedge\rho + W^{\sf new}_3\,\sigma\wedge\zeta + W^{\sf
new}_4\,\sigma\wedge\overline{\zeta} +
\\
& \ \ \ \ \ \ \ \ \ \ \ \ \ \ \ \ \ \ \ \ \ \ \ \ + W^{\sf
new}_5\,\overline{\sigma}\wedge\rho + W^{\sf
new}_6\,\overline{\sigma}\wedge\zeta + W^{\sf
new}_7\,\overline{\sigma}\wedge\overline{\zeta} +
\\
& \ \ \ \ \ \ \ \ \ \ \ \ \ \ \ \ \ \ \ \ \ \ \ \ \ \ \ \ \ \ \ \ \ \
\ \ \ \ \ \ \ \ \ + W^{\sf new}_8\,\rho\wedge\zeta + W^{\sf
new}_9\,\rho\wedge\overline{\zeta} +
\\
& \ \ \ \ \ \ \ \ \ \ \ \ \ \ \ \ \ \ \ \ \ \ \ \ \ \ \ \ \ \ \ \ \ \
\ \ \ \ \ \ \ \ \ \ \ \ \ \ \ \ \ \ \ \ \ \ \ \ \ \ \ \ \ + W^{\sf
new}_{10}\,\zeta\wedge\overline{\zeta},
\endaligned
\]
with twelve essential invariants:
\begin{equation}
\label{Essential-Invariants} \aligned \left[
\begin{array}{l}0\equiv U^{\sf new}_1=U^{\sf new}_2=U^{\sf new}_4,
\\
0\equiv V^{\sf new}_1=V^{\sf new}_2=V^{\sf new}_3,
\\
0\equiv W^{\sf new}_5=W^{\sf new}_6=W^{\sf new}_7=W^{\sf
new}_8=W^{\sf new}_9=W^{\sf new}_{10}.
\end{array}
\right.
\endaligned
\end{equation}

\begin{Theorem}
Two real analytic generic CR-submanifolds of $\mathbb C^4$
represented as the graph of three homogeneous defining equations of
the form:
\[\aligned
w_1-\overline w_1&=2i\,z\overline z+{\rm O}(3),
\\
w_2-\overline w_2&=2i\,z\overline z\,(z+\overline z)+{\rm O}(4),
\\
w_3-\overline w_3&=2\,z\overline z\,(z-\overline z)+{\rm O}(4),
\endaligned
\]
with nonzero corresponding essential functions $R$ are equivalent if
and only if the essential invariants
\thetag{\ref{Essential-Invariants}} associated to them are
identically equal.
\end{Theorem}

\end{document}